\documentclass[reqno]{gsm-l}

\usepackage[margin=1.5in]{geometry}
\usepackage[colorinlistoftodos]{todonotes}

\newcounter{todocounter}

\usepackage[T1]{fontenc}
\usepackage{tgpagella}
\usepackage{mathpazo}

\usepackage[pagebackref]{hyperref} \usepackage{xcolor}
\usepackage{amsmath,amsthm}
\usepackage{graphicx}
\usepackage{tensor}
\usepackage{color}
\usepackage[mathscr]{eucal}
\usepackage{MnSymbol}
\usepackage[frame,ps,matrix,arrow,curve,rotate,all,2cell,tips,color]{xy}
\usepackage{epic,eepic}
\usepackage{enumerate}

\newtheorem{theorem}[subsection]{Theorem}
\newtheorem{lemma}[subsection]{Lemma}

\newtheorem{proposition}[subsection]{Proposition}
\newtheorem{corollary}[subsection]{Corollary}

\theoremstyle{definition}
\newtheorem{definition}[subsection]{Definition}
\newtheorem{example}[subsection]{Example}
\newtheorem{remark}[subsection]{Remark}
\newtheorem{assumption}[subsection]{Assumption}
\newtheorem{convention}[subsection]{Convention}
\newtheorem{notation}[subsection]{Notation}
\newtheorem{conjecture}[subsection]{Conjecture}

\newtheorem{motivation}[subsection]{Motivation}

\newtheorem{problem}[section]{Problem}

\numberwithin{section}{chapter}
\numberwithin{equation}{subsection}

\usepackage{tikz}
\usetikzlibrary{matrix,arrows,decorations.pathmorphing}
\usetikzlibrary{backgrounds,positioning,shapes}
\usetikzlibrary{fit,petri,shapes.misc}

\tikzset{auto}

\tikzset{empty/.style={circle,inner sep=0pt,minimum size=6mm}}
\tikzset{emptyvt/.style={circle,inner sep=0pt,minimum size=0mm}}

\tikzset{plain/.style={circle,draw,very thick,
inner sep=0pt,minimum size=6mm}}

\tikzset{fatplain/.style={rounded rectangle,draw,very thick,minimum size=6mm}}

\tikzset{bigplain/.style={rounded rectangle,draw,very thick,minimum size=.8cm}}

\tikzset{yellowvt/.style={circle,draw,fill=yellow,very thick,inner sep=0pt,minimum size=6mm}}

\newcommand{\cable}[1]{\draw [thick, fill=lightgray] {#1} circle [radius=.175];}
\newcommand{\midcable}[1]{\draw [thick, fill=lightgray] {#1} circle [radius=.2];}
\newcommand{\smallcable}[1]{\draw [thick, fill=lightgray] {#1}circle [radius=.1];}

\tikzset{arrow/.style={->,thick}}
\tikzset{dashedarrow/.style={->,dashed,thick}}
\tikzset{dottedarrow/.style={->,dotted,thick}}
\tikzset{mapto/.style={|->,thick}}

\tikzset{implies/.style={thick,double,double equal sign distance,-implies}}

\tikzset{line/.style={thick}}
\tikzset{dottedline/.style={dotted,thick}}
\tikzset{dashedline/.style={dashed,thick}}

\tikzset{inputleg/.style={<-,thick}}
\tikzset{outputleg/.style={->,thick}}
\tikzset{dottedinput/.style={<-,dotted,thick}}

\tikzstyle{background}=[rectangle,fill=gray!40,inner sep=0.1cm,rounded corners=3mm]

\newcommand{\adjoint}{
\nicearrow\xymatrix{ \ar@<2pt>[r] & \ar@<2pt>[l]}}
\renewcommand{\hookrightarrow}{\nicexy{\ar@{^{(}->}[r] &}}
\newcommand{\nicearrow}{\SelectTips{cm}{10}}
\newcommand{\nicexy}{\nicearrow\xymatrix@C+5pt}

\newcommand{\pushout}{\ar@{}[dr]|-{\mathrm{pushout}}}
\newcommand{\drrpushout}{\ar@{}[drr]|(0.90){\Searrow}}

\renewcommand{\to}{\hspace{-.1cm}\nicearrow\xymatrix@C-.3cm{\ar[r]&}\hspace{-.1cm}}

\newcommand{\sigmasubstar}{\sigma_\bullet}
\newcommand{\sigmasupstar}{\sigma^\bullet}
\newcommand{\sigmasubdot}{\sigma_{\bullet}}
\newcommand{\sigmasupdot}{\sigma^{\bullet}}

\newcommand{\dT}{\partial_T}
\newcommand{\dU}{\partial_U}

\newcommand{\comp}{\circ}
\newcommand{\defn}{\overset{\mathrm{def}}{=\joinrel=}\,}

\newcommand{\finv}{f^{-1}}
\newcommand{\fphiinv}{f^{-1}_{\varphi}}
\newcommand{\ginv}{g^{-1}}

\newcommand{\compone}{\comp_1}
\newcommand{\comptwo}{\comp_2}
\newcommand{\compthree}{\comp_3}
\newcommand{\compi}{\comp_i}
\newcommand{\compj}{\comp_j}
\newcommand{\compk}{\comp_k}
\newcommand{\compjonel}{\comp_{j-1+l}}
\newcommand{\compionej}{\comp_{i-1+j}}

\newcommand{\xcompiw}{\uX \compi \uW}

\newcommand{\cphi}{C_{\varphi}}
\newcommand{\cphithree}{\cphi^{\geq 3}}
\newcommand{\cphizerozero}{\cphi^{(0,0)}}
\newcommand{\cphizeroatleastzero}{\cphi^{(0,\geq 0)}}
\newcommand{\cphizerozeroatleastzero}{C_{\varphi_0}^{(0,\geq 0)}}
\newcommand{\cphizerozeropm}{C_{\varphi, \pm}^{(0,0)}}
\newcommand{\cphizeroone}{\cphi^{(0,1)}}
\newcommand{\cphionezero}{\cphi^{(1,0)}}
\newcommand{\cphioneone}{\cphi^{(1,1)}}
\newcommand{\cphitwozero}{\cphi^{(2,0)}}
\newcommand{\cphizeroatleastone}{\cphi^{(0,\geq 1)}}
\newcommand{\cphizeroatleasttwo}{\cphi^{(0,\geq 2)}}

\newcommand{\cphiatleastthreezero}{\cphi^{(\geq 3, 0)}}
\newcommand{\cphioneatleastone}{\cphi^{(1,\geq 1)}}

\newcommand{\cphiatleastoneatleastzero}{\cphi^{(\geq 1, \geq 0)}}
\newcommand{\cphiatleasttwoatleastone}{\cphi^{(\geq 2,\geq 1)}}

\newcommand{\cpsi}{C_{\psi}}
\newcommand{\cpsione}{C_{\psione}}

\newcommand{\cpsithree}{\cpsi^{\geq 3}}
\newcommand{\cpsizerozero}{\cpsi^{(0,0)}}
\newcommand{\cpsizeroatleastzero}{\cpsi^{(0,\geq 0)}}
\newcommand{\cpsionezeroatleastzero}{\cpsione^{(0,\geq 0)}}
\newcommand{\cpsizerozeropm}{C_{\psi, \pm}^{(0,0)}}

\newcommand{\cpsionezero}{\cpsi^{(1,0)}}

\newcommand{\cpsizeroatleastone}{\cpsi^{(0,\geq 1)}}

\newcommand{\czeta}{C_{\zeta}}
\newcommand{\fphi}{f_{\varphi}}
\newcommand{\fpsi}{f_{\psi}}
\newcommand{\frho}{f_{\rho}}
\newcommand{\grho}{g_{\rho}}
\newcommand{\fzeta}{f_{\zeta}}
\newcommand{\gphi}{g_{\varphi}}
\newcommand{\gphiinv}{\gphi^{-1}}
\newcommand{\gpsi}{g_{\psi}}
\newcommand{\gzeta}{g_{\zeta}}

\newcommand{\bC}{\mathbb{C}}

\newcommand{\bN}{\mathbb{N}}

\newcommand{\bR}{\mathbb{R}}

\newcommand{\tensorunit}{\mathbb{1}}
\newcommand{\unit}{\mathbb{1}}

\newcommand{\As}{\mathsf{As}}

\newcommand{\Com}{\mathsf{Com}}

\newcommand{\DS}{\mathsf{DS}}

\newcommand{\Fin}{\mathsf{Fin}}
\newcommand{\Fins}{\Fin_S}
\newcommand{\Fint}{\Fin_T}
\newcommand{\Finset}{\Fin_{\set}}

\newcommand{\G}{\mathsf{G}}

\newcommand{\GrOp}{\mathsf{GrOp}}

\newcommand{\I}{\mathsf{I}}

\newcommand{\im}{\mathsf{Im}}

\newcommand{\Map}{\mathsf{Map}}

\newcommand{\End}{\mathsf{End}}

\newcommand{\Hist}{\mathsf{Hist}}
\newcommand{\Histone}{\Hist^1}

\newcommand{\sO}{\mathsf{O}}
\renewcommand{\O}{\mathsf{O}}

\renewcommand{\P}{\mathsf{P}}
\newcommand{\sP}{\mathsf{P}}

\newcommand{\prof}{\mathsf{Prof}}

\newcommand{\proffins}{\prof\left(\Fins\right)}
\newcommand{\profs}{\prof(S)}

\newcommand{\profns}{\prof^{\geq n}(S)}
\newcommand{\proft}{\prof(T)}
\newcommand{\profonet}{\prof^{\geq 1}(T)}
\newcommand{\profu}{\prof(U)}

\newcommand{\Q}{\mathsf{Q}}

\renewcommand{\Box}{\mathsf{Box}}
\newcommand{\boxs}{\Box_S}
\newcommand{\boxt}{\Box_T}

\newcommand{\profboxs}{\prof\left(\boxs\right)}
\newcommand{\profboxsboxs}{\profboxs \times \boxs}

\newcommand{\relation}{\mathsf{Rel}}
\newcommand{\rela}{\relation_A}
\newcommand{\relb}{\relation_B}

\newcommand{\T}{\mathsf{T}}

\newcommand{\betaone}{\beta_1}
\newcommand{\betatwo}{\beta_2}
\newcommand{\betathree}{\beta_3}

\newcommand{\phione}{\varphi_1}
\newcommand{\phitwo}{\varphi_2}
\newcommand{\phicompipsi}{\varphi \compi \psi}
\newcommand{\phicomponepsi}{\varphi \compone \psi}
\newcommand{\pione}{\pi_1}
\newcommand{\pitwo}{\pi_2}
\newcommand{\psione}{\psi_1}
\newcommand{\psitwo}{\psi_2}

\newcommand{\dm}{\mathsf{Dm}}
\newcommand{\supply}{\mathsf{Sp}}
\newcommand{\DN}{\mathsf{DN}}
\newcommand{\dn}{\DN}

\newcommand{\dnalpha}{\dn_{\alpha}}

\newcommand{\dnphi}{\dn_{\varphi}}

\newcommand{\dnpsi}{\dn_{\psi}}
\newcommand{\dnzeta}{\dn_{\zeta}}

\newcommand{\ewbeta}{\beta^{\mathsf{w}}_-}
\newcommand{\ewpi}{\pi^{\mathsf{w}}_-}
\newcommand{\ewpitwo}{\pi^{\mathsf{w}}_{2-}}
\newcommand{\ewphi}{\varphi^{\mathsf{w}}_-}
\newcommand{\ewpsi}{\psi^{\mathsf{w}}_-}

\newcommand{\iwpi}{\pi^{\mathsf{w}}_+}
\newcommand{\iwpitwo}{\pi^{\mathsf{w}}_{2+}}
\newcommand{\iwphi}{\varphi^{\mathsf{w}}_+}
\newcommand{\iwpsi}{\psi^{\mathsf{w}}_+}

\newcommand{\lp}{\mathsf{lp}}

\newcommand{\looppi}{\pi^{\lp}}

\newcommand{\looppitwo}{\pi^{\lp}_2}

\newcommand{\piin}{\pi^{\mathsf{sp}}}
\newcommand{\piini}{\piin_+}
\newcommand{\pitwoini}{\piin_{2+}}
\newcommand{\piine}{\piin_-}

\newcommand{\valpha}{v_{\alpha}}

\newcommand{\vphi}{v_{\varphi}}
\newcommand{\vpsi}{v_{\psi}}

\newcommand{\salpha}{s_{\alpha}}
\newcommand{\sbeta}{s_{\beta}}
\newcommand{\sbetaone}{s_{\beta_1}}
\newcommand{\sbetatwo}{s_{\beta_2}}
\newcommand{\sbetathree}{s_{\beta_3}}

\newcommand{\sphi}{s_{\varphi}}
\newcommand{\spsi}{s_{\psi}}
\newcommand{\spi}{s_{\pi}}
\newcommand{\spione}{s_{\pi_1}}
\newcommand{\spitwo}{s_{\pi_2}}
\newcommand{\ssigma}{s_{\sigma}}
\newcommand{\ssigmaone}{s_{\sigma_1}}
\newcommand{\ssigmatwo}{s_{\sigma_2}}
\newcommand{\ssigmaj}{s_{\sigma_j}}
\newcommand{\ssigmaupone}{s_{\sigma^1}}
\newcommand{\ssigmauptwo}{s_{\sigma^2}}
\newcommand{\szeta}{s_{\zeta}}

\newcommand{\dmalpha}{\dm_{\alpha}}

\newcommand{\dmbetaone}{\dm_{\beta_1}}
\newcommand{\dmbetatwo}{\dm_{\beta_2}}
\newcommand{\dmbetathree}{\dm_{\beta_3}}

\newcommand{\dmphi}{\dm_{\varphi}}
\newcommand{\dmpsi}{\dm_{\psi}}

\newcommand{\dmpione}{\dm_{\pi_1}}
\newcommand{\dmpitwo}{\dm_{\pi_2}}
\newcommand{\dmsigma}{\dm_{\sigma}}
\newcommand{\dmsigmaone}{\dm_{\sigma_1}}
\newcommand{\dmsigmatwo}{\dm_{\sigma_2}}
\newcommand{\dmsigmaj}{\dm_{\sigma_j}}
\newcommand{\dmsigmap}{\dm_{\sigma_p}}
\newcommand{\dmsigmaupone}{\dm_{\sigma^1}}
\newcommand{\dmsigmauptwo}{\dm_{\sigma^2}}
\newcommand{\dmsigmaupj}{\dm_{\sigma^j}}
\newcommand{\dmsigmaupr}{\dm_{\sigma^r}}

\newcommand{\dmzeta}{\dm_{\zeta}}

\newcommand{\supplyalpha}{\supply_{\alpha}}

\newcommand{\supplybetaone}{\supply_{\beta_1}}
\newcommand{\supplybetatwo}{\supply_{\beta_2}}
\newcommand{\supplybetathree}{\supply_{\beta_3}}

\newcommand{\supplyphi}{\supply_{\varphi}}
\newcommand{\supplypsi}{\supply_{\psi}}

\newcommand{\supplypione}{\supply_{\pi_1}}
\newcommand{\supplypitwo}{\supply_{\pi_2}}
\newcommand{\supplysigma}{\supply_{\sigma}}
\newcommand{\supplysigmaone}{\supply_{\sigma_1}}
\newcommand{\supplysigmatwo}{\supply_{\sigma_2}}
\newcommand{\supplysigmaupone}{\supply_{\sigma^1}}
\newcommand{\supplysigmauptwo}{\supply_{\sigma^2}}
\newcommand{\supplysigmaj}{\supply_{\sigma_j}}
\newcommand{\supplysigmaupj}{\supply_{\sigma^j}}
\newcommand{\supplyzeta}{\supply_{\zeta}}

\newcommand{\dnphicompipsi}{\dn_{\phicompipsi}}
\newcommand{\dmphicompipsi}{\dm_{\phicompipsi}}
\newcommand{\supplyphicompipsi}{\supply_{\phicompipsi}}
\newcommand{\sphicompipsi}{s_{\phicompipsi}}
\newcommand{\vphicompipsi}{v_{\phicompipsi}}

\newcommand{\WD}{\mathsf{WD}}
\newcommand{\wddot}{\WD_\bullet}

\newcommand{\wdzero}{\WD_0}
\newcommand{\UWD}{\mathsf{UWD}}
\newcommand{\uwd}{\UWD}
\newcommand{\uwds}{\UWD^S}

\newcommand{\uwdset}{\UWD^{\set}}

\newcommand{\ua}{\underline{a}}
\newcommand{\ub}{\underline{b}}
\newcommand{\uc}{\underline{c}}

\newcommand{\um}{\underline{m}}

\newcommand{\us}{\underline{s}}
\newcommand{\ut}{\underline{t}}

\newcommand{\uw}{\underline{w}}
\newcommand{\ux}{\underline{x}}
\newcommand{\uy}{\underline{y}}
\newcommand{\uz}{\underline{z}}

\newcommand{\uW}{\underline{W}}
\newcommand{\uX}{\underline{X}}
\newcommand{\uY}{\underline{Y}}
\newcommand{\uZ}{\underline{Z}}

\newcommand{\chizero}{\chi^0}
\newcommand{\udelta}{\underline{\delta}}
\newcommand{\uepsilon}{\underline{\epsilon}}
\newcommand{\ulambda}{\underline{\lambda}}
\newcommand{\uphi}{\underline{\phi}}

\newcommand{\uomega}{\underline{\omega}}
\newcommand{\upsi}{\underline{\psi}}
\newcommand{\usigma}{\underline{\sigma}}
\newcommand{\utheta}{\underline{\theta}}

\newcommand{\calc}{\mathcal{C}}
\newcommand{\C}{\calc}

\newcommand{\fC}{\mathfrak{C}}

\newcommand{\set}{\mathsf{Set}}

\newcommand{\smallprof}[1]
{\raisebox{.05cm}{\scalebox{0.8}{#1}}}

\newcommand{\smallbinom}[2]
{\raisebox{.05cm}{\scalebox{0.8}{$\binom{#1}{#2}$}}}

\newcommand{\apm}{a_{\pm}}
\newcommand{\aminusplusminusa}
{\smallprof{$\binom{A \setminus \apm}{A}$}}

\newcommand{\capba}
{\smallprof{$\binom{B}{A}$}}

\newcommand{\bjiuaji}
{\smallprof{$\binom{b^j_i}{\ua^j_i}$}}
\newcommand{\bjua}
{\smallprof{$\binom{b_j}{\ua}$}}

\newcommand{\bw}
{\smallprof{$\binom{B}{W}$}}

\newcommand{\ciua}
{\smallprof{$\binom{c_i}{\ua}$}}

\newcommand{\cub}
{\smallprof{$\binom{c}{\ub}$}}

\newcommand{\ciubi}
{\smallprof{$\binom{c_i}{\ub_i}$}}
\newcommand{\cibcompja}
{\smallprof{$\binom{c_i}{\ub \compj \ua}$}}

\newcommand{\cjnothing}
{\smallprof{$\binom{c_j}{}$}}
\newcommand{\cjub}
{\smallprof{$\binom{c_j}{\ub}$}}
\newcommand{\cjubj}
{\smallprof{$\binom{c_j}{\ub_j}$}}

\newcommand{\ciub}
{\smallprof{$\binom{c_i}{\ub}$}}
\newcommand{\fcifub}
{\smallprof{$\binom{fc_i}{f\ub}$}}

\newcommand{\cjasupj}
{\smallprof{$\binom{c_j}{\ua_j}$}}

\newcommand{\ccsingle}
{\smallprof{$\binom{c}{c}$}}
\newcommand{\fcsingle}
{\smallprof{$\binom{fc}{fc}$}}
\newcommand{\cici}
{\smallprof{$\binom{c_i}{c_i}$}}
\newcommand{\cjcj}
{\smallprof{$\binom{c_j}{c_j}$}}
\newcommand{\cjubjtauj}
{\smallprof{$\binom{c_j}{\ub_j\tau_j}$}}

\newcommand{\csigmajubsigmaj}
{\smallprof{$\binom{c_{\sigma(j)}}{\ub_{\sigma(j)}}$}}

\newcommand{\csigmaiub}
{\smallprof{$\binom{c_{\sigma(i)}}{\ub}$}}
\newcommand{\csigmaiubtau}
{\smallprof{$\binom{c_{\sigma(i)}}{\ub\tau}$}}

\newcommand{\dplusginvbdplusb}
{\smallprof{$\binom{D \amalg \ginv b}{D \amalg b}$}}

\newcommand{\dnothing}
{\smallprof{$\binom{d}{}$}}
\newcommand{\ddnothing}
{\smallprof{$\binom{\{d, d\}}{}$}}

\newcommand{\dinothing}
{\smallprof{$\binom{d_i}{}$}}
\newcommand{\djnothing}
{\smallprof{$\binom{d_j}{}$}}

\newcommand{\dua}
{\smallprof{$\binom{d}{\ua}$}}

\newcommand{\dub}
{\smallprof{$\binom{d}{\ub}$}}
\newcommand{\duboneubn}
{\smallprof{$\binom{d}{\ub_1,\ldots,\ub_n}$}}
\newcommand{\dubsigmaoneubsigman}
{\smallprof{$\binom{d}{\ub_{\sigma(1)},\ldots,\ub_{\sigma(n)}}$}}
\newcommand{\dubonetauoneubntaun}
{\smallprof{$\binom{d}{\ub_1\tau_1, \ldots , \ub_n\tau_n}$}}

\newcommand{\dconecn}
{\smallprof{$\binom{d}{c_1, \ldots, c_n}$}}
\newcommand{\fdfconefcn}
{\smallprof{$\binom{fd}{fc_1, \ldots, fc_n}$}}
\newcommand{\ducsigma}
{\smallprof{$\binom{d}{\uc\sigma}$}}

\newcommand{\dccompib}
{\smallprof{$\binom{d}{\uc \compi \ub}$}}
\newcommand{\fdfccompib}
{\smallprof{$\binom{fd}{f(\uc \compi \ub)}$}}

\newcommand{\dccompia}
{\smallprof{$\binom{d}{\uc \compi \ua}$}}

\newcommand{\dccompjb}
{\smallprof{$\binom{d}{\uc \compj \ub}$}}
\newcommand{\dccompiacompjminusonepluslb}
{\smallprof{$\binom{d}{(\uc \compi \ua) \comp_{j-1+l} \ub}$}}

\newcommand{\dccompjbcompia}
{\smallprof{$\binom{d}{(\uc \compj \ub) \comp_{i} \ua}$}}
\newcommand{\dccompibcompja}
{\smallprof{$\binom{d}{\uc \compi (\ub \compj \ua)}$}}
\newcommand{\dccompibcompioneja}
{\smallprof{$\binom{d}{(\uc \compi \ub) \comp_{i-1+j} \ua}$}}

\newcommand{\dccompsigmaib}
{\smallprof{$\binom{d}{\uc \comp_{\sigma(i)} \ub}$}}
\newcommand{\dcsigmacompibtau}
{\smallprof{$\binom{d}{(\uc\sigma) \compi (\ub\tau)}$}}
\newcommand{\dccompsigmaibsigmatau}
{\smallprof{$\binom{d}{(\uc \comp_{\sigma(i)} \ub)(\sigma \compi \tau)}$}}

\newcommand{\dd}
{\smallprof{$\binom{d}{d}$}}
\newcommand{\diplusonedi}
{\smallprof{$\binom{D_{i+1}}{D_i}$}}

\newcommand{\dcsingle}
{\smallprof{$\binom{d}{c}$}}
\newcommand{\duc}
{\smallprof{$\binom{d}{\uc}$}}
\newcommand{\fdfuc}
{\smallprof{$\binom{fd}{f\uc}$}}

\newcommand{\dcsigma}
{\smallprof{$\binom{d}{\uc\sigma}$}}
\newcommand{\fdfcsigma}
{\smallprof{$\binom{fd}{f\uc\sigma}$}}

\newcommand{\UZ}
{\smallprof{$\binom{U}{Z}$}}
\newcommand{\vvzero}
{\smallprof{$\binom{V}{V_0}$}}
\newcommand{\vuoneutwo}
{\smallprof{$\binom{V}{U_1,U_2}$}}
\newcommand{\vuoneun}
{\smallprof{$\binom{V}{U_1,\ldots,U_N}$}}
\newcommand{\vzerouoneun}
{\smallprof{$\binom{V_0}{U_1,\ldots,U_N}$}}
\newcommand{\vjminusonevj}
{\smallprof{$\binom{V_{j-1}}{V_j}$}}
\newcommand{\vzerovr}
{\smallprof{$\binom{V_{0}}{V_r}$}}
\newcommand{\vx}
{\smallprof{$\binom{V}{X}$}}
\newcommand{\vz}
{\smallprof{$\binom{V}{Z}$}}

\newcommand{\wa}
{\smallprof{$\binom{W}{A}$}}
\newcommand{\wv}
{\smallprof{$\binom{W}{V}$}}
\newcommand{\wjminusonewj}
{\smallprof{$\binom{W_{j-1}}{W_j}$}}
\newcommand{\wminusww}
{\smallprof{$\binom{W \setminus w}{W}$}}
\newcommand{\wx}
{\smallprof{$\binom{W}{X}$}}
\newcommand{\wy}
{\smallprof{$\binom{W}{Y}$}}
\newcommand{\wz}
{\smallprof{$\binom{W}{Z}$}}

\newcommand{\xw}
{\smallprof{$\binom{X}{W}$}}
\newcommand{\xiuw}
{\smallprof{$\binom{X_i}{\uW}$}}
\newcommand{\xiwonewm}
{\smallprof{$\binom{X_i}{W_1,\ldots,W_m}$}}
\newcommand{\xiuwi}
{\smallprof{$\binom{X_i}{\uW_i}$}}
\newcommand{\xjuw}
{\smallprof{$\binom{X_j}{\uW}$}}

\newcommand{\xx}
{\smallprof{$\binom{X}{X}$}}
\newcommand{\xonetwox}
{\smallprof{$\binom{X_{12}}{X}$}}
\newcommand{\xtwothreex}
{\smallprof{$\binom{X_{23}}{X}$}}
\newcommand{\xthreefourx}
{\smallprof{$\binom{X_{34}}{X}$}}
\newcommand{\xixi}
{\smallprof{$\binom{X_i}{X_i}$}}

\newcommand{\xkxkxkone}
{\smallprof{$\binom{X^{\geq k}}{X_k,\, X^{\geq k+1}}$}}

\newcommand{\xpm}{x_{\pm}}
\newcommand{\xsubpm}{x_{\pm}}
\newcommand{\xminusplusminus}{X \setminus x_{\pm}}
\newcommand{\xminusxpm}
{X \setminus x_{\pm}}
\newcommand{\xminusxpmx}
{\smallprof{$\binom{\xminusxpm}{X}$}}

\newcommand{\xminusplusminusx}
{\smallprof{$\binom{\xminusplusminus}{X}$}}
\newcommand{\xminusxx}
{\smallprof{$\binom{X \setminus x}{X}$}}
\newcommand{\xminusxxbar}
{\smallprof{$\binom{\xbar \setminus \xsubpm}{\xbar}$}}

\newcommand{\xminusxxprime}
{\smallprof{$\binom{X \setminus x}{X'}$}}
\newcommand{\xminusxxprimeminusx}
{\smallprof{$\binom{X \setminus x}{X' \setminus x}$}}

\newcommand{\xprimex}
{\smallprof{$\binom{X'}{X}$}}
\newcommand{\xprimeux}
{\smallprof{$\binom{X'}{\uX}$}}
\newcommand{\xstarx}
{\smallprof{$\binom{X^*}{X}$}}
\newcommand{\xstarxprime}
{\smallprof{$\binom{X^*}{X'}$}}
\newcommand{\xstaryminusxone}
{\smallprof{$\binom{X^*}{Y \setminus x(1)}$}}

\newcommand{\xxprime}
{\smallprof{$\binom{X}{X'}$}}
\newcommand{\xxempty}
{\smallprof{$\binom{X}{X, \varnothing}$}}
\newcommand{\xprimeminusxx}
{\smallprof{$\binom{X' \setminus x}{X}$}}
\newcommand{\xprimeminusxxprime}
{\smallprof{$\binom{X' \setminus x}{X'}$}}
\newcommand{\xprimeminusxxminusx}
{\smallprof{$\binom{X' \setminus x}{X \setminus x}$}}

\newcommand{\xminusxonex}
{\smallprof{$\binom{X \setminus x^1}{X}$}}
\newcommand{\xminusxtwox}
{\smallprof{$\binom{X \setminus x^2}{X}$}}
\newcommand{\xminusxxminusxone}
{\smallprof{$\binom{X \setminus x}{X \setminus x^1}$}}
\newcommand{\xminusxxminusxtwo}
{\smallprof{$\binom{X \setminus x}{X \setminus x^2}$}}

\newcommand{\xxonexn}
{\smallprof{$\binom{X}{X_1, \ldots, X_n}$}}

\newcommand{\xchoosey}
{\smallprof{$\binom{X}{Y}$}}
\newcommand{\xplusyy}
{\smallprof{$\binom{X \amalg Y}{Y}$}}
\newcommand{\xplusyprimexy}
{\smallprof{$\binom{X \amalg Y'}{X, Y}$}}
\newcommand{\xplusyprimexyprime}
{\smallprof{$\binom{X \amalg Y'}{X, Y'}$}}
\newcommand{\xplusyprimexplusy}
{\smallprof{$\binom{X \amalg Y'}{X \amalg Y}$}}
\newcommand{\xplusyxplusyprime}
{\smallprof{$\binom{X \amalg Y}{X \amalg Y'}$}}
\newcommand{\xplusyxy}
{\smallprof{$\binom{X \amalg Y}{X, Y}$}}
\newcommand{\xbarplusyxbary}
{\smallprof{$\binom{\xbar \amalg y}{\xbar, y}$}}
\newcommand{\xplusyxybar}
{\smallprof{$\binom{\xbar \amalg \ybar}{\xbar, \ybar}$}}
\newcommand{\xplusyxprimey}
{\smallprof{$\binom{X \amalg Y}{X', Y}$}}
\newcommand{\xplusychoosexyprime}
{\smallprof{$\binom{X \amalg Y}{X, Y'}$}}
\newcommand{\xprimeplusyxprimey}
{\smallprof{$\binom{X' \amalg Y}{X', Y}$}}
\newcommand{\xprimeplusyxy}
{\smallprof{$\binom{X' \amalg Y}{X, Y}$}}
\newcommand{\xprimeplusyxplusy}
{\smallprof{$\binom{X' \amalg Y}{X \amalg Y}$}}
\newcommand{\xplusyxyprime}
{\smallprof{$\binom{X' \amalg Y'}{X', Y'}$}}
\newcommand{\xprimeplusyprimexy}
{\smallprof{$\binom{X' \amalg Y'}{X, Y}$}}
\newcommand{\xprimeplusyprimexplusy}
{\smallprof{$\binom{X' \amalg Y'}{X \amalg Y}$}}
\newcommand{\xplusyxprimeplusy}
{\smallprof{$\binom{X \amalg Y}{X' \amalg Y}$}}

\newcommand{\xplusyminusxxy}
{\smallprof{$\binom{(X \amalg Y) \setminus \{x\}}{X, Y}$}}
\newcommand{\xplusyminusxxminusxy}
{\smallprof{$\binom{(X \amalg Y) \setminus \{x\}}{X \setminus x, Y}$}}
\newcommand{\xplusyminusxxplusy}
{\smallprof{$\binom{(X \amalg Y) \setminus \{x\}}{X \amalg Y}$}}

\newcommand{\xplusyz}
{\smallprof{$\binom{X \amalg Y \amalg Z}{X \amalg Y, Z}$}}
\newcommand{\xyplusz}
{\smallprof{$\binom{X \amalg Y \amalg Z}{X, Y \amalg Z}$}}
\newcommand{\xplusyplusz}
{\smallprof{$\binom{X \amalg Y \amalg Z}{X, Y, Z}$}}

\newcommand{\xz}
{\smallprof{$\binom{X}{Z}$}}

\newcommand{\ypm}{y_{\pm}}
\newcommand{\yuw}
{\smallprof{$\binom{Y}{\uW}$}}
\newcommand{\yiuw}
{\smallprof{$\binom{Y_i}{\uW}$}}
\newcommand{\yux}
{\smallprof{$\binom{Y}{\uX}$}}
\newcommand{\yuxsigma}
{\smallprof{$\binom{Y}{\uX\sigma}$}}

\newcommand{\yw}
{\smallprof{$\binom{Y}{W}$}}
\newcommand{\yx}
{\smallprof{$\binom{Y}{X}$}}
\newcommand{\yxbar}
{\smallprof{$\binom{\ybar}{\xbar}$}}
\newcommand{\yxy}
{\smallprof{$\binom{Y}{X, y}$}}
\newcommand{\yprimex}
{\smallprof{$\binom{Y'}{X}$}}
\newcommand{\yxprime}
{\smallprof{$\binom{Y}{X'}$}}
\newcommand{\yprimexprime}
{\smallprof{$\binom{Y'}{X'}$}}
\newcommand{\yxplusy}
{\smallprof{$\binom{Y}{X \amalg Y}$}}

\newcommand{\yonex}
{\smallprof{$\binom{Y_1}{X}$}}
\newcommand{\ytwox}
{\smallprof{$\binom{Y_2}{X}$}}
\newcommand{\yonexone}
{\smallprof{$\binom{Y_1}{X_1}$}}
\newcommand{\ytwoxtwo}
{\smallprof{$\binom{Y_2}{X_2}$}}
\newcommand{\yoneplusytwoyoneytwo}
{\smallprof{$\binom{Y_1 \amalg Y_2}{Y_1, Y_2}$}}
\newcommand{\xoneplusxtwoxonextwo}
{\smallprof{$\binom{X_1 \amalg X_2}{X_1, X_2}$}}
\newcommand{\yoneplusytwoxonextwo}
{\smallprof{$\binom{Y_1 \amalg Y_2}{X_1, X_2}$}}
\newcommand{\yoneplusytwoxoneplusxtwo}
{\smallprof{$\binom{Y_1 \amalg Y_2}{X_1 \amalg X_2}$}}

\newcommand{\yxonextwo}
{\smallprof{$\binom{Y}{X_1, X_2}$}}

\newcommand{\yxonextwobar}
{\smallprof{$\binom{\ybar}{\xbar_1, \xbar_2}$}}
\newcommand{\yxonetwo}
{\smallprof{$\binom{Y}{X_{12}}$}}
\newcommand{\yxtwothree}
{\smallprof{$\binom{Y}{X_{23}}$}}
\newcommand{\yxthreefour}
{\smallprof{$\binom{Y}{X_{34}}$}}

\newcommand{\yuponetwox}
{\smallprof{$\binom{Y^{12}}{X}$}}
\newcommand{\yuptwothreex}
{\smallprof{$\binom{Y^{23}}{X}$}}
\newcommand{\yupthreefourx}
{\smallprof{$\binom{Y^{34}}{X}$}}

\newcommand{\yiux}
{\smallprof{$\binom{Y_i}{\uX}$}}
\newcommand{\yjux}
{\smallprof{$\binom{Y_j}{\uX}$}}
\newcommand{\yxonexn}
{\smallprof{$\binom{Y}{X_1, \ldots, X_n}$}}
\newcommand{\yxonexbign}
{\smallprof{$\binom{Y}{X_1, \ldots, X_N}$}}
\newcommand{\yxonexnbar}
{\smallprof{$\binom{\ybar}{\xbar_1, \ldots, \xbar_n}$}}
\newcommand{\yxonexbignbar}
{\smallprof{$\binom{\ybar}{\xbar_1, \ldots, \xbar_N}$}}
\newcommand{\yxcompiw}
{\smallprof{$\binom{Y}{\uX \compi \uW}$}}
\newcommand{\yzbar}
{\smallprof{$\binom{\ybar}{\zbar}$}}
\newcommand{\yminusyx}
{\smallprof{$\binom{Y \setminus y}{X}$}}

\newcommand{\yplusxyx}
{\smallprof{$\binom{Y \amalg X}{Y, X}$}}
\newcommand{\yminusxoney}
{\smallprof{$\binom{Y \setminus x(1)}{Y}$}}

\newcommand{\yminusypm}
{Y \setminus y_{\pm}}
\newcommand{\yminusypmx}
{\smallprof{$\binom{\yminusypm}{X}$}}
\newcommand{\yminusypmxminusxpm}
{\smallprof{$\binom{\yminusypm}{\xminusxpm}$}}
\newcommand{\yminusypmy}
{\smallprof{$\binom{\yminusypm}{Y}$}}

\newcommand{\yy}
{\smallprof{$\binom{Y}{Y}$}}
\newcommand{\yyzero}
{\smallprof{$\binom{Y}{Y_0}$}}
\newcommand{\yjminusoneyj}
{\smallprof{$\binom{Y_{j-1}}{Y_j}$}}

\newcommand{\yjyjplusone}
{\smallprof{$\binom{Y_j}{Y_{j+1}}$}}
\newcommand{\yzeroyp}
{\smallprof{$\binom{Y_0}{Y_p}$}}
\newcommand{\yzeroxonexbign}
{\smallprof{$\binom{Y_0}{X_1,\ldots,X_N}$}}

\newcommand{\yprimey}
{\smallprof{$\binom{Y'}{Y}$}}
\newcommand{\yyprime}
{\smallprof{$\binom{Y}{Y'}$}}

\newcommand{\yyone}
{\smallprof{$\binom{Y}{Y_1}$}}
\newcommand{\yytwo}
{\smallprof{$\binom{Y}{Y_2}$}}
\newcommand{\yyuponetwo}
{\smallprof{$\binom{Y}{Y^{12}}$}}
\newcommand{\yyuptwothree}
{\smallprof{$\binom{Y}{Y^{23}}$}}
\newcommand{\yyupthreefour}
{\smallprof{$\binom{Y}{Y^{34}}$}}

\newcommand{\yminusyy}
{\smallprof{$\binom{Y \setminus y}{Y}$}}
\newcommand{\ypluszyz}
{\smallprof{$\binom{Y \amalg Z}{Y, Z}$}}
\newcommand{\yz}
{\smallprof{$\binom{Y}{Z}$}}
\newcommand{\ybarz}
{\smallprof{$\binom{\ybar}{Z}$}}

\newcommand{\zw}
{\smallprof{$\binom{Z}{W}$}}
\newcommand{\zminuswz}
{\smallprof{$\binom{Z\setminus w}{Z}$}}

\newcommand{\zx}
{\smallprof{$\binom{Z}{X}$}}
\newcommand{\zxbar}
{\smallprof{$\binom{Z}{\xbar}$}}
\newcommand{\zbarxbar}
{\smallprof{$\binom{\zbar}{\xbar}$}}
\newcommand{\zxonextwo}
{\smallprof{$\binom{Z}{X_1,X_2}$}}

\newcommand{\zxonexnconecp}
{\smallprof{$\binom{Z}{X_1,\ldots,X_N,c_1,\ldots,c_p}$}}
\newcommand{\zux}
{\smallprof{$\binom{Z}{\uX}$}}
\newcommand{\zxprime}
{\smallprof{$\binom{Z}{X'}$}}

\newcommand{\yminusyxminusx}
{\smallprof{$\binom{Y \setminus y}{X \setminus x}$}}

\newcommand{\zy}
{\smallprof{$\binom{Z}{Y}$}}

\newcommand{\zuy}
{\smallprof{$\binom{Z}{\uY}$}}
\newcommand{\zycompjxcompiw}
{\smallprof{$\binom{Z}{(\uY \compj \uX) \compi \uW}$}}
\newcommand{\zycompixcompionejw}
{\smallprof{$\binom{Z}{(\uY \compi \uX) \compionej \uW}$}}

\newcommand{\zminuszz}
{\smallprof{$\binom{Z \setminus z}{Z}$}}
\newcommand{\zzerozq}
{\smallprof{$\binom{Z_{0}}{Z_q}$}}
\newcommand{\zjminusonezj}
{\smallprof{$\binom{Z_{j-1}}{Z_j}$}}

\newcommand{\emptyprof}
{\smallprof{$\binom{\varnothing}{\varnothing}$}}
\newcommand{\emptyx}
{\smallprof{$\binom{\varnothing}{X}$}}
\newcommand{\xempty}
{\smallprof{$\binom{X}{\varnothing}$}}

\newcommand{\emptynothing}
{\smallprof{$\binom{\varnothing}{}$}}
\newcommand{\smallxnothing}
{\smallprof{$\binom{x}{}$}}
\newcommand{\xnothing}
{\smallprof{$\binom{X}{}$}}
\newcommand{\xnminusonexn}
{\smallprof{$\binom{X_{n-1} \amalg X_n}{X_{n-1}, X_n}$}}
\newcommand{\ynothing}
{\smallprof{$\binom{Y}{}$}}
\newcommand{\ybarnothing}
{\smallprof{$\binom{\ybar}{}$}}
\newcommand{\smallynothing}
{\smallprof{$\binom{y}{}$}}
\newcommand{\smallyempty}
{\smallprof{$\binom{y}{\varnothing}$}}
\newcommand{\smallyemptyy}
{\smallprof{$\binom{y}{\varnothing, y}$}}

\newcommand{\starnothing}
{\smallprof{$\binom{*}{}$}}

\newcommand{\inp}{\mathsf{in}}
\newcommand{\out}{\mathsf{out}}
\newcommand{\rd}{\mathsf{rd}}
\newcommand{\up}{\mathsf{up}}

\newcommand{\fin}{f^{\inp}}
\newcommand{\fout}{f^{\out}}
\newcommand{\frd}{f^{\rd}}
\newcommand{\fup}{f^{\up}}
\newcommand{\frdminusx}{\frd_{\minusx}}
\newcommand{\fupminusx}{\fup_{\minusx}}

\newcommand{\gin}{g^{\inp}}
\newcommand{\gout}{g^{\out}}
\newcommand{\grd}{g^{\rd}}
\newcommand{\gup}{g^{\up}}

\newcommand{\smallxin}{x^{\inp}}
\newcommand{\smallxout}{x^{\out}}

\newcommand{\uin}{U^{\inp}}
\newcommand{\uout}{U^{\out}}
\newcommand{\vuin}{\uin_v}
\newcommand{\vuout}{\uout_v}
\newcommand{\vin}{V^{\inp}}
\newcommand{\vout}{V^{\out}}
\newcommand{\vvin}{\vin_v}
\newcommand{\vvout}{\vout_v}
\newcommand{\win}{W^{\inp}}
\newcommand{\wout}{W^{\out}}
\newcommand{\vwin}{\win_v}
\newcommand{\vwout}{\wout_v}
\newcommand{\xin}{X^{\inp}}
\newcommand{\xout}{X^{\out}}
\newcommand{\vxin}{X^{\inp}_v}
\newcommand{\vxout}{X^{\out}_v}
\newcommand{\vxyin}{(X \amalg Y)^{\inp}_v}
\newcommand{\vxyout}{(X \amalg Y)^{\out}_v}
\newcommand{\xprimein}{X'^{\inp}}
\newcommand{\xprimeout}{X'^{\out}}
\newcommand{\vxstarin}{X^{*\inp}_v}
\newcommand{\vxstarout}{X^{*\out}_v}

\newcommand{\fbar}{\overline{f}}

\newcommand{\wbar}{\overline{W}}
\newcommand{\xbar}{\overline{X}}
\newcommand{\ybar}{\overline{Y}}
\newcommand{\zbar}{\overline{Z}}

\newcommand{\varphibar}{\overline{\varphi}}

\newcommand{\psibar}{\overline{\psi}}

\newcommand{\uone}{U^1}
\newcommand{\utwo}{U^2}
\newcommand{\utwoplus}{\utwo_+}
\newcommand{\utwominus}{\utwo_-}

\newcommand{\wminuswin}{(W \setminus w)^{\inp}}
\newcommand{\wminuswout}{(W \setminus w)^{\out}}
\newcommand{\vwminuswin}{(W \setminus w)^{\inp}_v}
\newcommand{\vwminuswout}{(W \setminus w)^{\out}_v}

\newcommand{\xminusxin}{(X \setminus x)^{\inp}}
\newcommand{\xminusxout}{(X \setminus x)^{\out}}
\newcommand{\vxminusxin}{(X \setminus x)^{\inp}_v}
\newcommand{\vxminusxout}{(X \setminus x)^{\out}_v}

\newcommand{\vxminusxtwoin}{(X \setminus x^2)^{\inp}_v}

\newcommand{\vxminusxiout}{(X \setminus x^i)^{\out}_v}

\newcommand{\vxplus}{v(x_+)}
\newcommand{\vxoneplus}{v(\xoneplus)}
\newcommand{\vxtwoplus}{v(\xtwoplus)}
\newcommand{\vxiplus}{v(\xiplus)}

\newcommand{\vxminus}{v(x_-)}
\newcommand{\vxoneminus}{v(\xoneminus)}
\newcommand{\vxtwominus}{v(\xtwominus)}

\newcommand{\xplusy}{X \amalg Y}
\newcommand{\vxplusyin}{(\xplusy)^{\inp}_v}
\newcommand{\vxplusyout}{(\xplusy)^{\out}_v}

\newcommand{\yin}{Y^{\inp}}
\newcommand{\yout}{Y^{\out}}
\newcommand{\vyin}{Y^{\inp}_v}
\newcommand{\vyout}{Y^{\out}_v}

\newcommand{\zin}{Z^{\inp}}
\newcommand{\zout}{Z^{\out}}
\newcommand{\vzin}{\zin_v}
\newcommand{\vzout}{\zout_v}
\newcommand{\zminuszin}{(Z \setminus z)^{\inp}}
\newcommand{\zminuszout}{(Z \setminus z)^{\out}}
\newcommand{\vzminuswin}{(Z \setminus w)^{\inp}_v}
\newcommand{\vzminuswout}{(Z \setminus w)^{\out}_v}
\newcommand{\vzminuszin}{(Z \setminus z)^{\inp}_v}
\newcommand{\vzminuszout}{(Z \setminus z)^{\out}_v}
\newcommand{\zminuswin}{(Z \setminus w)^{\inp}}
\newcommand{\zminuswout}{(Z \setminus w)^{\out}}

\newcommand{\xplus}{x_+}
\newcommand{\xoneplus}{x^1_+}
\newcommand{\xtwoplus}{x^2_+}
\newcommand{\xiplus}{x^i_+}

\newcommand{\xminus}{x_-}
\newcommand{\xoneminus}{\xminus^1}
\newcommand{\xtwominus}{\xminus^2}

\newcommand{\minusxsuponetwo}{\setminus x^{12}}
\newcommand{\minusxonextwo}{\setminus \{x^1,x^2\}}

\newcommand{\minusx}{\setminus x}

\newcommand{\fminusx}{f_{\minusx}}
\newcommand{\finminusx}{\fin_{\minusx}}
\newcommand{\foutminusx}{\fout_{\minusx}}
\newcommand{\xminusx}{X \setminus x}
\newcommand{\minusxplus}{\setminus \xplus}
\newcommand{\minusxoneplus}{\setminus \xoneplus}
\newcommand{\minusxtwoplus}{\setminus \xtwoplus}
\newcommand{\minusxiplus}{\setminus \xiplus}

\newcommand{\minusy}{\setminus y}

\DeclareMathOperator{\id}{id}
\DeclareMathOperator{\Id}{Id}

\DeclareMathOperator{\Ob}{Ob}

\newcommand{\andspace}{\qquad\text{and}\qquad}

\newcommand{\forspace}{\quad\text{for}\quad}

\newcommand{\impliesspace}{\quad\text{implies}\quad}
\newcommand{\orspace}{\qquad\text{or}\qquad}

\newcommand{\minushspace}{\hspace{-.1cm}}

\makeindex

\begin{document}

\frontmatter

\title{Operads of Wiring Diagrams}

\author{Donald Yau}
\address{The Ohio State University at Newark, Newark, OH, USA}
\email{yau.22@osu.edu}

\subjclass[2000]{}


\date{\today}

\begin{abstract}
Wiring diagrams and undirected wiring diagrams are graphical languages for describing interconnected processes and their compositions.  These objects have enormous potentials for applications in many different disciplines, including computer science, cognitive neuroscience, dynamical systems, network theory, and circuit diagrams.  It is known that the collection of wiring diagrams is an operad and likewise for undirected wiring diagrams.  

This monograph is a comprehensive study of the combinatorial structure of various operads of wiring diagrams and undirected wiring diagrams.  Our first main objective is to prove a finite presentation theorem for each operad of wiring diagrams, describing each one in terms of just a few operadic generators and a small number of generating relations.  For example, the operad of wiring diagrams has $8$ generators and $28$ generating relations, while the operad of undirected wiring diagrams has $6$ generators and $17$ generating relations.

Our second main objective is to prove a corresponding finite presentation theorem for algebras over each operad of wiring diagrams.  As applications we provide finite presentations for the propagator algebra, the algebra of discrete systems, the algebra of open dynamical systems, and the (typed) relational algebra.  We also provide a partial verification of Spivak's conjecture regarding the quotient-freeness of the relational algebra.

Our third main objective is to construct explicit operad maps among the several operads of wiring diagrams.  In particular, there is a surjective operad map from the operad of all wiring diagrams, including delay nodes, to the operad of undirected wiring diagrams.

This monograph is intended for graduate students, mathematicians, scientists, and engineers interested in operads and wiring diagrams.  Assuming no prior knowledge of categories, operads, and wiring diagrams, this monograph is self-contained and can be used as a supplement in a graduate course and for independent study.  There are over $100$ graphical illustrations and a chapter with a list of problems.
\end{abstract}

\maketitle

\cleardoublepage
\thispagestyle{empty}
\cleardoublepage

\setcounter{page}{7}

\tableofcontents


\newcommand{\where}[1]{\> \pageref{#1} \>}
\newcommand{\blob}{ \> \> \> \hspace{1em}}

\chapter*{List of Notations}

\begin{tabbing}
\hspace{1.5cm}
\=\textbf{Notation}\hspace{1.3cm}
\= \textbf{Page}\hspace{.25cm}
\= \textbf{Description} \\

\textbf{Ch. \ref{ch02-wiring-diagrams}}
\> $S$ \where{notation:s-class} a class\\

\> $\profs$ \where{notation:profs} the collection of $S$-profiles\\

\> $\us = (s_1,\ldots,s_n)$ \where{notation:us} an $S$-profile of length $n$\\

\> $\profns$ \where{notation:profns} $S$-profiles of length at least $n$\\

\> $\Sigma_n$ \where{notation:sigma-n} the symmetric group on $n$ letters\\

\> $(\O,\tensorunit,\gamma)$ \where{notation:colored-operad-gamma} an $S$-colored operad\\

\> $\O\duc$ \where{notation:oduc} the $\duc$-entry of $\O$\\

\> $\uc\sigma$ \where{notation:c-sigma} the right permutation of $\uc$ by $\sigma$\\

\> $\unit_c$ \where{notation:unit-c} the $c$-colored unit\\

\> $\gamma$ \where{notation:operad-composition-gamma} operadic composition\\

\> $\sigma\langle k_1, \ldots , k_n \rangle$ \where{notation:block-permutation} block permutation\\

\> $\tau_1 \oplus \cdots \oplus \tau_n$ \where{notation:block-sum} block sum\\

\> $\End(X)$ \where{ex:endomorphism} endomorphism operad of $X$\\

\> $\As$ \where{notation:As} associative operad\\

\> $\Com$ \where{notation:Com} commutative operad\\

\> $\GrOp$ \where{notation:grop} operad of graph operations\\

\> $(\O,\tensorunit,\comp)$ \where{notation:colored-operad-comp} an $S$-colored operad\\

\> $\compi$ \where{notation:operadic-compi} $\compi$-composition in an operad\\

\> $\sigma \compi \tau$ \where{notation:sigma-compi-tau} $\compi$-composition of permutations\\

\> $\Fin$ \where{notation:fin} the category of finite sets and functions\\

\> $\Fins$ \where{notation:fins} the category of $S$-finite sets\\

\> $(\xin,\xout)$ \where{notation:sbox} an $S$-box\\

\> $\xin$ \where{notation:xin} the set of inputs of an $S$-box $X$\\

\> $\xout$ \where{notation:xout} the set of outputs of an $S$-box $X$\\

\> $\boxs$ \where{notation:xout} the collection of $S$-boxes\\

\> $v$ \where{notation:value-assignment} value assignment\\

\> $\varnothing$ \where{notation:emptybox} the empty $S$-box\\

\> $\dn$ \where{notation:delay-nodes} the set of delay nodes\\

\> $\dm$ \where{notation:demand} the set of demand wires\\

\> $\supply$ \where{notation:supply} the set of supply wires\\

\> $s : \dm \to \supply$ \where{notation:supplier} the supplier assignment\\

\> $\ewphi$ \where{notation:external-wasted} the set of external wasted wires of $\varphi$\\

\> $\iwphi$ \where{notation:internal-wasted} the set of internal wasted wires of $\varphi$\\

\> $\WD\yux$ \where{notation:wd-yux} the $\yux$-entry of $\WD$\\

\> $\WD$ \where{notation:wd} the collection of wiring diagrams\\

\> $\tensorunit_Y$ \where{wd-unit} the $Y$-colored unit in $\WD$\\

\blob\\

\textbf{Ch. \ref{ch03-generating-wd}}
\> $\epsilon$ \where{notation:empty-wd} the empty wiring diagram\\

\> $\delta_d$ \where{notation:dn-wd} a $1$-delay node\\

\> $\tau_{X,Y}$ \where{notation:name-change-wd} a name change\\

\> $\theta_{X,Y}$ \where{notation:2-cell} a $2$-cell\\

\> $\lambda_{X,x}$ \where{notation:1-loop} a $1$-loop\\

\> $\sigma_{X,x_1,x_2}$ \where{notation:insplit} an in-split\\

\> $\sigma^{Y,y_1,y_2}$ \where{notation:outsplit} an out-split\\

\> $\omega_{Y,y}$ \where{notation:wasted-wire} a $1$-wasted wire\\

\> $\omega^{X,x}$ \where{notation:external-wasted-wire} a $1$-internal wasted wire\\

\> $\varphi_1 \comp \cdots \comp \varphi_k$ \where{iterated-compone} (iterated) $\compone$\\

\> $|T|$ \where{notation:cardinality} the cardinality of a finite set\\

\blob\\

\textbf{Ch. \ref{ch04-decomposition}}
\> $\looppi$ \where{notation:loop-element} the set of  loop elements in $\pi$\\

\> $\piini$ \where{notation:int-supplied} the set of internally supplied elements in $\pi$\\

\> $\piine$ \where{notation:ext-supplied} the set of externally supplied elements in $\pi$\\

\blob\\

\textbf{Ch. \ref{ch05-stratified-presentation}}
\> $|\Psi|$ \where{notation:simplex-composition} the composition of a simplex $\Psi$\\

\> $(\psi_1, \ldots, \psi_n)$ \where{notation:simplex} an $n$-simplex in $\WD$\\

\> $\wddot$ \where{notation:normal-wd} the collection of normal wiring diagrams\\

\> $\wdzero$ \where{notation:strict-wd} the collection of strict wiring diagrams\\

\blob\\

\textbf{Ch. \ref{ch06-wd-algebras}}
\> $A_c$ \where{notation:a-sub-c} the $c$-colored entry of an algebra $A$\\

\> $\mu$ \where{notation:algebra-structure-map} the structure map of an algebra\\

\> $\mu_{\zeta} \compi \mu_{\xi}$ \where{notation:compi-algebra-map} $\compi$-composition of structure maps\\

\> $\partial$ \where{notation:truncation} truncation\\

\> $\Hist^k$ \where{notation:hist} the set of $k$-historical propagators\\

\> $D_k$ \where{ex:moment-delay} $k$-moment delay function\\

\> $\sP_X$ \where{notation:p-sub-x} the set of $1$-historical propagators of type $X$\\

\> $(-)_{\xplus}$ \where{notation:sub-xplus} take only the $v(x_+)$-entry\\

\> $(-)_{\minusxplus}$ \where{notation:sub-minus-xplus} remove the $v(x_+)$-entry\\

\> $\DS(X)$ \where{notation:ds} the collection of $X$-discrete systems\\

\> $\ux_I$ \where{xminusn} take only the entries indexed by $I$\\

\> $\ux_{\setminus I}$ \where{xminusn} remove the entries indexed by $I$\\

\> $\G_X$ \where{g-sub-x} the $X$-colored entry of the algebra\\
\blob of open dynamical systems\\

\textbf{Ch. \ref{ch07-undirected-wiring-diagrams}}
\> $\nicexy@C-.5cm{\cdot \ar[r] & \cdot & \cdot \ar[l]}$ \where{uwd-cospan} a cospan\\

\> $\uwd$ \where{uwdyux} the collection of undirected wiring diagrams\\

\> $X_{[i,j]}$ \where{x-sub-ij} $X_i \amalg \cdots \amalg X_j$ or $\varnothing$\\

\> $\cphi$ \where{def:compi-uwd} the set of cables in $\varphi$\\

\blob\\

\textbf{Ch. \ref{ch08-generating-uwd}}
\> $\epsilon$ \where{def:uwd-gen-empty} the empty cell\\

\> $\omega_*$ \where{def:uwd-gen-1output} a $1$-output wire\\

\> $\tau_f$ \where{def:uwd-gen-namechange} an undirected name change\\

\> $\theta_{(X,Y)}$ \where{def:uwd-gen-2cell} an undirected $2$-cell\\

\> $\lambda_{(X,\xsubpm)}$ \where{uwd-loop} a loop\\

\> $\sigma^{(X,x_1,x_2)}$ \where{uwd-split} a split\\

\blob\\

\textbf{Ch. \ref{ch09-decomp-uwd}}
\> $\cpsi^{(m,n)}$ \where{notation:cable-subsets} the set of $(m,n)$-cables\\

\> $\cpsizerozero$ \where{notation:cable-subsets} the set of wasted cables\\

\> $\cpsi^{(\geq m,n)}$ \where{notation:cable-subsets} the set of $(\geq m, n)$-cables\\

\> $\cpsi^{(m, \geq n)}$ \where{notation:cable-subsets} the set of $(m, \geq n)$-cables\\

\> $\cpsi^{(\geq m, \geq n)}$ \where{notation:cable-subsets} the set of $(\geq m, \geq n)$-cables\\

\blob\\

\textbf{Ch. \ref{ch10-stratified-uwd}}
\> $(\psi_1,\ldots,\psi_n)$  \where{rk:uwd-simplex-notation} an $n$-simplex in $\uwd$\\

\blob\\

\textbf{Ch. \ref{ch11-uwd-algebras}}
\> $A^X$ \where{def:relational-construction} the set of functions $X \to A$\\

\> $\powerset(X)$ \where{def:relational-construction}  the power set of $X$\\

\> $\rela(X)$ \where{def:relational-construction} $\powerset\left(A^X\right)$\\

\> $\rela$ \where{def:relational-algebra-biased} the relational algebra of $A$\\

\> $\relation$ \where{def:typed-relational-algebra} the typed relational algebra\\

\blob\\

\textbf{Ch. \ref{ch12-maps}}
\> $\ybar$ \where{wddot-to-uwd-colors} $\yin \amalg \yout$\\

\> $\chi$ \where{wddot-uwd-chi} the operad map $\wddot \to \uwd$\\

\> $\chizero$ \where{chizero} the operad map $\wdzero \to \uwd$\\

\blob\\

\textbf{Ch. \ref{ch13-wd-uwd}}
\> $\rho$ \where{wd-uwd-rho} the operad map $\WD \to \uwd$\\

\end{tabbing}

\mainmatter

\chapter{Introduction}
\label{ch01-introduction}


\section{What are Wiring Diagrams?}

Wiring diagrams form a kind of graphical language that describes operations or processes with multiple inputs and multiple outputs and how such operations are wired together to form a larger and more complicated operation.  Some visual examples of wiring diagrams are in \eqref{wd-first-example}, Chapter \ref{ch03-generating-wd}, and Example \ref{ex:factoring-pi}.  The first type of wiring diagrams that we are going to study in this monograph was first introduced in \cite{rupel-spivak}, with variants studied in \cite{spivak15,spivak15b,vsl}.  In \cite{spivak15,spivak15b} wiring diagrams without delay nodes (Def. \ref{def:wd-without-dn}), which we call \emph{normal} wiring diagrams, were used to study mode-dependent networks \index{mode-dependent networks}, discrete systems, and dynamical systems.\index{dynamical systems}  In \cite{vsl} wiring diagrams without delay nodes and whose supplier assignments are bijections (Def. \ref{def:strict-wd}), which we call \emph{strict} wiring diagrams, were used to study open dynamical systems.  

Wiring diagrams are by nature directed, in the sense that every operation has a finite set of inputs and a finite set of outputs, each element of which is allowed to carry a value of some kind.  There is also an undirected version of wiring diagrams \cite{spivak13,spivak14}.  Unlike a wiring diagram, in an undirected wiring diagram, each operation is a finite set, each element of which is again allowed to carry a value.  Some visual examples of undirected wiring diagrams are in \eqref{uwd-first-picture}, Example \ref{ex:uwd-jointly-surjective}, and Chapters \ref{ch08-generating-uwd} and \ref{ch09-decomp-uwd}.  For those familiar with operad theory, the distinction between wiring diagrams and undirected wiring diagrams is similar to that between operads and cyclic operads.  Just as cyclic operads are not simpler than operads, undirected wiring diagrams are not really simpler than wiring diagrams and have their own subtlety.

The main reason that wiring diagrams and undirected wiring diagrams are important is that they have enormous potential for applications in many different disciplines.  Wiring diagrams and undirected wiring diagrams allow one to consider a finite collection of related operations, wired together in some way, as an operation itself.  Such an operation can then be considered as a single operation within a yet larger collection of operations, and so forth. For instance, a finite collection of related operations may be a group of neurons in a certain region of the brain, a collection of codes within a large computer program, or a few related agents within a large supply-chain.  In fact, the authors of \cite{rupel-spivak} cited both computer science and cognitive neuroscience as potential applications of wiring diagrams.  Furthermore, in \cite{spivak15b,vsl} wiring diagrams were used to study dynamical systems and to model certain kinds of differential equations.  Many potential fields of applications are mentioned in the introduction of \cite{vsl}.  In \cite{spivak13} undirected wiring diagrams were used to study database relational queries, plug-and-play devices, recursion, and circuit diagrams.  

The substitution process involving wiring diagrams and undirected wiring diagrams described in the previous paragraph can be captured precisely using the notion of colored operads.  A \emph{colored operad}, or just an operad, is a mathematical object that describes operations with multiple inputs and one output and their compositions.  A colored operad in which there are only unary operations is exactly a category.  If one restricts even further to just the $1$-colored case in which the unary operations form a set, then one has exactly a monoid, such as the set of integers under addition.  Therefore, a colored operad is a multiple-input generalization of a category, and in fact colored operads are also called symmetric multicategories.  Multicategories without symmetric group actions were introduced by Lambek \cite{lambek}.  One-colored operads, together with the name \emph{operad}, were introduced by May \cite{may72} in the topological setting.  See \cite{may97} for the definition of a one-colored operad in a symmetric monoidal category.  The book \cite{yau-operad} is an elementary introduction to colored operads in symmetric monoidal categories.  The  book \cite{jy2} has more in-depth discussion of colored operads and related objects.

In \cite{rupel-spivak} Rupel and Spivak observed that the collection of wiring diagrams is a colored operad $\WD$, in which the operadic composition corresponds precisely to the substitution process described above.  Each colored operad has associated algebras, on which the colored operad acts.  The operadic action is similar to the action of an associative algebra on a left module.  The authors of \cite{rupel-spivak} defined a $\WD$-algebra, called the \emph{propagator algebra}, that models a certain kind of input-output process.  

Closely related colored operads of sub-classes of wiring diagrams were introduced in \cite{spivak15,spivak15b,vsl}.  We will denote by $\wddot$ (Def. \ref{def:wd-without-dn}) the operad of \emph{normal} wiring diagrams, meaning those without delay nodes.  In \cite{spivak15b} Spivak introduced a $\wddot$-algebra, called the \emph{algebra of discrete systems}, that is closely related to a Moore machine, also known as a discrete state machine.  Also we will write $\wdzero$ (Def. \ref{def:strict-wd}) for the operad of \emph{strict} wiring diagrams, meaning those without delay nodes and whose supplier assignments are bijections.  In \cite{vsl} a $\wdzero$-algebra, called the \emph{algebra of open dynamical systems}, was defined that models a certain kind of differential equations.  Likewise, in \cite{spivak13} Spivak constructed the colored operad $\uwd$ of undirected wiring diagrams.  Spivak also defined a $\uwd$-algebra called the \emph{relational algebra}, which was used to model relational queries in database.

\section{Purposes of this Monograph}

This monograph is a comprehensive study of the combinatorial structure of the operads $\WD$, $\wddot$, $\wdzero$, and $\uwd$ of (normal/strict/undirected) wiring diagrams, their algebras, and the relationships between these operads.  Specifically, our main results are of the following three types.
\begin{description}
\item[Finite Presentation for Operads]
For each of the operads $\WD$, $\wddot$, $\wdzero$, and $\uwd$, we prove a finite presentation theorem that describes the operad in terms of just a few operadic generators and a small number of generating relations.  For the operad of wiring diagrams $\WD$, there are $8$ generating wiring diagrams and $28$  generating relations.  For the smaller operads $\wddot$ and $\wdzero$ of normal and strict wiring diagrams, the numbers of operadic generators and of generating relations are $(7,28)$ and $(4,8)$, respectively.  For the operad of undirected wiring diagrams $\uwd$, there are $6$ operadic generators and $17$ generating relations.
\item[Finite Presentation for Algebras]
For each of the operads $\WD$, $\wddot$, $\wdzero$, and $\uwd$, we prove a corresponding finite presentation theorem for their algebras. To be more precise, we describe $\WD$-algebras using $8$ generating structure maps and $28$ generating axioms.  So finite presentation refers to the $\WD$-algebra structure maps, not the elements in the underlying set.  Similar finite presentations are also obtained for the algebras over the operads $\wddot$ of normal wiring diagrams, $\wdzero$ of strict wiring diagrams, and $\uwd$ of undirected wiring diagrams.  As applications we provide finite presentations for the propagator algebra over $\WD$, the algebra of discrete systems over $\wddot$, the algebra of open dynamical systems over $\wdzero$, and the (typed) relational algebra over $\uwd$.  Along the way, we provide a partial verification of Spivak's conjecture \cite{spivak13} regarding the quotient-freeness of the relational algebra.
\item[Maps Between Operads]
We construct a commutative diagram
\[\nicexy{\wdzero \ar[dr]_-{\chizero} \ar[r] & \wddot \ar[d]^-{\chi} \ar[r] & \WD \ar[dl]^-{\rho}\\
& \uwd &}\]
of operad maps, in which the horizontal maps are operad inclusions.  For each of the operad maps $\chizero$, $\chi$, and $\rho$, we compute precisely the image in $\uwd$.  In particular, the operad map $\rho : \WD \to \uwd$ is surjective.  The existence of the operad map $\rho$ answers a question raised in \cite{rupel-spivak}.
\end{description}
The end of this chapter contains several tables that summarize the main results and contain some key references.

The finite presentation theorems for the operads $\WD$, $\wddot$, $\wdzero$, and $\uwd$ reduce the structure of these operads and their algebras to just a few generators and relations.  For example, our finite presentation theorem for the operad $\WD$ reduces the understanding of this operad to just $8$ simple wiring diagrams and $28$ simple relations among them.  Likewise, the structure map of a general $\WD$-algebra can be quite involved, as can be seen in the propagator algebra \cite{rupel-spivak}.  Our finite presentation theorem for $\WD$-algebras reduces the definition and understanding of $\WD$-algebras to just $8$  simple generating structure maps and $28$ generating axioms, almost all of which are trivial to check in practice.  We will further illuminate this point when we discuss the finite presentations for the propagator algebra over $\WD$, the algebra of discrete systems over $\wddot$, the algebra of open dynamical systems over $\wdzero$, and the (typed) relational algebra over $\uwd$.

To give our finite presentation theorems for the operads $\WD$, $\wddot$, $\wdzero$, and $\uwd$, and for their algebras an even more familiar feel, let us recall a few other places where various kinds of finite presentations occur.  Each type of finite presentation is a way to reduce a large, usually infinite, collection of conditions to a finite, or at least smaller, collection of conditions, thereby making the relevant structure more manageable.
\begin{enumerate}
\item
In basic group theory and commutative ring theory \cite{atiyah-macdonald,rotman}, it is quite common to consider finite presentation of groups and modules.  For example, over a commutative Noetherian ring, every finitely generated module is also finitely presented \cite{rotman} (Prop. 7.59).
\item
A cornerstone in category theory, Mac Lane's Coherence\index{coherence} Theorem  \cite{maclane63,maclane} can be regarded as a finite presentation theorem for monoidal categories.  Roughly speaking, this theorem says that, in any \index{monoidal category}monoidal category, the infinite collection of commutative formal diagrams has a finite presentation.  The generators are the associativity isomorphism, the left and the right units, and their inverses.  The generating relations are the Pentagon Axiom\index{Pentagon Axiom} for $4$-fold iterated tensor products and two unity axioms.  
\item
In the linear setting, the operads for associative algebras, commutative algebras, Lie algebras, Leibniz algebras, Poisson algebras, and many others, are finitely presented \cite{gk}.  For example, the associative operad has one generator, which in its algebras corresponds to the usual multiplication $A \otimes A \to A$ of an associative algebra.  The associative operad has one generating relation, which in its algebras corresponds to the usual associativity axiom, $(ab)c = a(bc)$, of an associative algebra.
\item
In \cite{baez-erbele,bsz} a finite presentation is given for the symmetric monoidal category of signal-flow\index{signal-flow graphs} graphs.  In applications signal-flow graphs form another kind of graphical language that describes processes with inputs and outputs and relations between them.  One main difference between (undirected) wiring diagrams and signal-flow graphs is that the composition of signal-flow graphs is done by \emph{grafting}.  This means that the outputs of one signal-flow graph are connected to the inputs of another signal-flow graph.  This is similar to the situation in \index{string diagrams}string diagrams \cite{jsv,ssr15}  On the other hand, the operadic composition of (undirected) wiring diagrams is done by \emph{substitution}, which is pictorially depicted in \eqref{wd-compi-picture} for wiring diagrams and in \eqref{uwd-compi-picture} for undirected wiring diagrams.  In more conceptual terms, the collection of signal-flow graphs is a prop, hence an algebra over the operad for props \cite{jy2} (Theorem 14.1).  On the other hand, the collection of (undirected) wiring diagrams is an operad.
\item
Closer to the topic of this monograph is \cite{jy2} (Ch.7), where finite presentations are given for various \index{graph groupoids}graph groupoids including those for colored operads, \index{props}colored props, and \index{wheeled prop}colored wheeled props.  In fact, the way we phrase and verify our finite presentation theorems for the operads of (undirected) wiring diagrams and for their algebras is conceptually similar to the presentation in \cite{jy2} (Ch.7).  One way to explain this conceptual similarity is that, for both (undirected)  wiring diagrams and the graphs for, say, colored wheeled props, the composition is done by substitution.  However, wiring diagrams are in several ways more complicated than the graphs in \cite{jy2}.  In fact, the graphs in \cite{jy2} do not have delay nodes, internal and external wasted wires, and split wires, all of which can happen in a wiring diagram.  See, for example, the wiring diagram in \eqref{wd-first-example}.
\end{enumerate}

\section{Audience and Features}

The main results in this monograph--namely, the finite presentation theorems for the various operads of (undirected) wiring diagrams and for their algebras as well as operad maps between them--are new.  So this monograph should be useful to mathematicians with an interest in operads and (undirected) wiring diagrams.  Furthermore, due to the wide variety of potential applications, we also intend to make this monograph and this subject accessible to scientists and engineers.  

With such a large audience in mind, the prerequisite for this monograph has been kept to an absolute  minimum.  In particular, we assume the reader  is comfortable with basic concepts of sets, functions, and mathematical induction.  No prior knowledge of categories, operads, and (undirected) wiring diagrams is assumed.  The presentation of the material proceeds at a fairly leisurely pace and is roughly at the advanced undergraduate to beginning graduate level.  To motivate various constructions and concepts, we have many examples and a lot of discussion that explains the intuition behind the scene.  Furthermore, there are over $100$ pictures throughout this monograph that help the reader visualize (undirected) wiring diagrams.  Finally, to solidify one's understanding of the subject, the reader may work through the problems in Chapter \ref{ch-problems}.  There are enough problems there to keep one busy for a few days.

\section{Chapter Summaries}

This monograph is divided into three parts.
\begin{description}
\item[Part 1] Wiring Diagrams (Chapters \ref{ch02-wiring-diagrams}--\ref{ch06-wd-algebras})

This part contains the finite presentation theorems for the operad $\WD$ of wiring diagrams, the operad $\wddot$ of normal wiring diagrams, the operad $\wdzero$ of strict wiring diagrams, and their algebras.

\item[Part 2] Undirected Wiring Diagrams (Chapters \ref{ch07-undirected-wiring-diagrams}--\ref{ch11-uwd-algebras})

This part contains the finite presentation theorems for the operad $\uwd$ of undirected wiring diagrams and for their algebras.

\item[Part 3] Maps Between Operads of Wiring Diagrams (Chapters \ref{ch12-maps}--\ref{ch:further-reading})

This part contains the construction and description of various operad maps between the operads $\WD$, $\wddot$, $\wdzero$, and $\uwd$.  It also contains a chapter with a list of problems and a chapter with references for further reading.
\end{description}
Each part begins with a brief introduction and a reading guide.  Below is a brief description of each chapter.

In Chapter \ref{ch02-wiring-diagrams}, to keep this document self-contained, we first recall two equivalent definitions of a colored operad.  The first definition is in terms of May's operad structure map $\gamma$, and the other one is in terms of the $\compi$-compositions.  After recalling the definition of a wiring diagram, we provide a proof of the fact from \cite{rupel-spivak} that the collection of wiring diagrams $\WD$ is a colored operad.  The main difference here is that we use the definition of a colored operad based on the $\compi$-compositions.  In this monograph, we prefer to work with the $\compi$-compositions rather than May's operad structure map $\gamma$ because the $\compi$-compositions are more convenient in phrasing and verifying our finite presentation theorems.

In Chapter \ref{ch03-generating-wd} we introduce $8$ \emph{generating wiring diagrams} and $28$ \emph{elementary relations} among them.  On the one hand, one may regard this chapter as a long list of concrete examples of wiring diagrams and their operadic compositions.  On the other hand, in later chapters we will see that these finite collections of generating wiring diagrams and elementary relations are sufficient to describe the operad $\WD$ of wiring diagrams, its variants $\wddot$ and $\wdzero$, and their algebras.

For the finite presentation theorems for the operad $\WD$ of wiring diagrams and its variants $\wddot$ and $\wdzero$, we will need to be able to decompose every wiring diagram in terms of the $8$ generating wiring diagrams in a highly structured way.  The purpose of Chapter \ref{ch04-decomposition} is to supply all the steps needed to establish such a decomposition.

The finite presentation theorems for the operad $\WD$ of wiring diagrams as well as its two variants $\wddot$ and $\wdzero$ are given in Chapter \ref{ch05-stratified-presentation}; see Theorems \ref{thm:wd-generator-relation}, \ref{thm:without-dn-coherence}, and \ref{thm:strict-wd-coherence}.  Since we are not working in the linear setting (e.g., of modules) where we can take quotients, we need to be extra careful in phrasing our finite presentations for the operads $\WD$, $\wddot$, and $\wdzero$.  For this purpose, a crucial concept is that of a \emph{stratified presentation}, which is the highly structured decomposition mentioned in the previous paragraph.  The results in Chapter \ref{ch04-decomposition} guarantees the existence of a stratified presentation for each wiring diagram.  This implies the finite generation parts of our finite presentation theorems for $\WD$, $\wddot$, and $\wdzero$.  The relation parts of the finite presentation theorems are phrased in terms of our concept of an \emph{elementary equivalence}.  Roughly speaking, an elementary equivalence means replacing one side of either (i) an elementary relation in Chapter \ref{ch03-generating-wd} or (ii) an operad associativity/unity axiom for the generating wiring diagrams, by the other side.

In Chapter \ref{ch06-wd-algebras}  we discuss finite presentations for algebras over the operads $\WD$, $\wddot$, and $\wdzero$.  In each case, the finite presentation for algebras is a consequence of the finite presentation theorem for the corresponding operad of wiring diagrams.  To illustrate the finite presentation for $\WD$-algebras, we will describe the propagator algebra in terms of $8$ generating structure maps and $28$ generating axioms.  To illustrate our finite presentation for $\wddot$-algebras, we will describe the algebra of discrete systems in terms of $7$ generating structure maps and $28$ generating axioms.   To illustrate our finite presentation for $\wdzero$-algebras, we will similarly describe the algebra of open dynamical systems in terms of $4$ generating structure maps and $8$ generating axioms.  This finishes Part 1 on wiring diagrams.

Part 2 begins with Chapter \ref{ch07-undirected-wiring-diagrams}, where we first recall the notion of an undirected wiring diagram.  Then we give a proof of the fact that the collection of undirected wiring diagrams forms an operad $\uwd$.  As in Chapter \ref{ch02-wiring-diagrams}, the operad structure on $\uwd$ as well as its proof are both given in terms of the $\compi$-compositions because the finite presentation theorems are easier to phrase this way.  One subtlety about the operad $\uwd$ is that undirected wiring diagrams may have wasted cables (Def. \ref{def:pre-uwd}), which are cables that are not soldered to any wires.  As opposed to what was stated in \cite{spivak14} (Example 7.4.2.10),  wasted cables \emph{cannot} be excluded from the definition of undirected wiring diagrams.  In fact, wasted cables can actually be created from operadic composition of undirected wiring diagrams without wasted cables.  We will make this point precise in Example \ref{ex:uwd-jointly-surjective}.

Chapter \ref{ch08-generating-uwd} is the undirected analogue of Chapter \ref{ch03-generating-wd}.  In this chapter, we describe $6$ \emph{generating undirected wiring diagrams} and $17$ \emph{elementary relations} among them.  On the one hand, one may regard this chapter as a long list of concrete examples of undirected wiring diagrams and their operadic compositions.  On the other hand, in later chapters we will see that these finite collections of generating undirected wiring diagrams and elementary relations are sufficient to describe the operad $\uwd$ of undirected wiring diagrams.

Chapter \ref{ch09-decomp-uwd} is the undirected analogue of Chapter \ref{ch04-decomposition}.  In this chapter, we show that each undirected wiring diagram can be decomposed in terms of the generating undirected wiring diagrams in a highly structured way.  Such a decomposition is needed to establish the finite presentation theorem for the operad $\uwd$.

Chapter \ref{ch10-stratified-uwd} is the undirected analogue of Chapter \ref{ch05-stratified-presentation}.  In this chapter, we establish the finite presentation theorem for the operad $\uwd$ of undirected wiring diagrams; see Theorem \ref{thm:uwd-generator-relation}.  This result is phrased in terms of the generating undirected wiring diagrams and an undirected version of an elementary equivalence.

Chapter \ref{ch11-uwd-algebras} contains the finite presentation theorem for $\uwd$-algebras.  This result is a consequence of the finite presentation theorem for the operad $\uwd$.  It describes each $\uwd$-algebra in terms of $6$ generating structure maps and $17$ generating axioms, almost all of which are trivial to check in practice.  We will illustrate this point with the relation algebra and the typed relational algebra from \cite{spivak13}.  We will also provide a partial verification of Spivak's conjecture regarding the quotient-freeness of the relational algebra.  This finishes Part 2 on undirected wiring diagrams.

Part 3 begins with Chapter \ref{ch12-maps}, in which we first construct the operad inclusions $\wdzero \to \wddot \to \WD$.  Recall that $\wddot$ is the operad of normal wiring diagrams--i.e., those without delay nodes--and that $\wdzero$ is the operad of strict wiring diagrams--i.e., those without delay nodes and whose supplier assignments are bijections.  Then we construct an operad map $\chi : \wddot \to \uwd$, essentially by forgetting directions, and its restriction $\chizero : \wdzero \to \uwd$.  For each of the operad maps $\chi$ and $\chizero$, we compute precisely the image in $\uwd$.  In the terminology of Notation \ref{notation:cable-subsets}, the image of the operad map $\chi$ consists of precisely the undirected wiring diagrams \emph{without} wasted cables and $(0, \geq 2)$-cables.  The image of the operad map $\chizero$ consists of precisely the undirected wiring diagrams with only $(1,1)$-cables and $(2,0)$-cables.

In Chapter \ref{ch13-wd-uwd} we extend the operad map $\chi : \wddot \to \uwd$ to an operad map $\rho : \WD \to \uwd$ that is defined for all wiring diagrams.  We prove that the operad map $\rho$ is surjective, so every undirected wiring diagram is the image of some wiring diagram.  The operad map $\rho$ is slightly subtle because wiring diagrams may have delay nodes, while undirected wiring diagrams do not seem to have an exact analogue of delay nodes.  In fact, for this reason Rupel and Spivak \cite{rupel-spivak} (Section 4.1) expressed doubt about the possibility that there be an operad map $\WD \to \uwd$.  We will see that delay nodes, far from being an obstruction, play a critical role in the surjectivity of the operad map $\rho$.

Chapter \ref{ch-problems} contains some problems about operads and (undirected) wiring diagrams arising from the earlier chapters.  Chapter \ref{ch:further-reading} contains some relevant references on categories, operads, props, and their applications in the sciences.  This finishes Part 3.

\section{References for the Main Results and Examples}

The following table summaries the key references for the various operads of (undirected) wiring diagrams and their finite presentation theorems.

\begin{center}
\begin{small}
\begin{tabular}{|c|c|c|c|c|}\hline
& Operad of & \shortstack{Finite\\ Presentation} & Generators & Relations \\ \hline
$\WD$ & \shortstack{wiring diagrams\\ (Theorem \ref{wd-operad}}) & Theorem \ref{thm:wd-generator-relation} & \shortstack{8\\(Def. \ref{def:generating-wiring-diagrams})} & \shortstack{28\\(Def. \ref{def:elementary-relations})}\\ \hline
$\wddot$ & \shortstack{normal wiring\\ diagrams (Prop. \ref{prop:without-dn-operads})} & Theorem \ref{thm:without-dn-coherence} & \shortstack{7\\(Def. \ref{def:generating-without-dn}(1))} & \shortstack{28\\(Def. \ref{def:generating-without-dn}(4))}\\ \hline
$\wdzero$ &  \shortstack{strict wiring\\ diagrams (Prop. \ref{prop:strict-operads})} & Theorem \ref{thm:strict-wd-coherence} & \shortstack{4\\(Def. \ref{def:generating-strict}(1))} & \shortstack{8\\(Def. \ref{def:generating-strict}(5))}\\ \hline
$\uwd$ & \shortstack{undirected wiring\\ diagrams (Theorem \ref{uwd-operad})} & Theorem \ref{thm:uwd-generator-relation} & \shortstack{6\\(Def. \ref{def:generating-uwd})} & \shortstack{17\\(Def. \ref{def:elementary-relation-uwd})}\\ \hline
\end{tabular}
\end{small}
\end{center}

The following table summarizes the key references for the finite presentation theorems for algebras over the various operads of (undirected) wiring diagrams.

\begin{center}
\begin{small}
\begin{tabular}{|c|c|c|c|}\hline
\shortstack{Algebras\\ over} & \shortstack{Finite\\ Presentation} & (Generators, Relations) & \shortstack{Key Example} \\ \hline
$\WD$ & Theorem \ref{thm:wd-algebra} & \shortstack{(8, 28)\\ (Def. \ref{def:wd-algebra})} & \shortstack{propagator algebra\\ (Theorem \ref{prop:propagator-algebra-is-algebra})}\\ \hline
$\wddot$ & Theorem \ref{thm:normal-wd-algebra} & \shortstack{(7, 28)\\ (Def. \ref{def:normal-algebra})} & \shortstack{discrete systems\\ (Theorem \ref{thm:algebra-ds})}\\ \hline
$\wdzero$ & Theorem \ref{thm:strict-wd-algebra} & \shortstack{(4, 8)\\ (Def. \ref{def:strict-algebra})} & \shortstack{open dynamical systems\\ (Theorem \ref{prop:ods-algebra-is-algebra})}\\ \hline
$\uwd$ & Theorem \ref{thm:uwd-algebra} & \shortstack{(6, 17)\\ (Def. \ref{def:uwd-algebra-biased})} & \shortstack{(typed) relational algebra\\ (Theorems \ref{relational-algebra-is-algebra} and \ref{typed-relational-algebra})}\\ \hline
\end{tabular}
\end{small}
\end{center}

The following table summarizes the key references for the operad maps between the various operads of (undirected) wiring diagrams.

\begin{center}
\begin{small}
\begin{tabular}{|c|c|c|c|c|}\hline
Operad Map & Reference & Note/Image\\ \hline
$\wdzero \to \wddot \to \WD$ & Prop. \ref{prop:strict-normal-wd} & inclusions\\ \hline
$\nicexy{\wddot \ar[r]^-{\chi} & \uwd}$ & Theorem \ref{wddot-uwd-operad-map} & \shortstack{no $(0,0)$- and $(0, \geq 2)$-cables\\ (Theorem \ref{thm:chi-image})}\\ \hline 
$\nicexy{\wdzero \ar[r]^-{\chizero} & \uwd}$ & Theorem \ref{chizero-wdzero-uwd} & only $(1,1)$- and $(2,0)$-cables\\ \hline
$\nicexy{\WD \ar[r]^-{\rho} & \uwd}$ & Theorem \ref{wd-uwd-operad-map} & \shortstack{surjective\\ (Theorem \ref{thm:wd-uwd-surjective})}\\ \hline
\end{tabular}
\end{small}
\end{center}

\part{Wiring Diagrams}

The main purpose of this part is to describe the combinatorial structure of 
\begin{enumerate}
\item the operad $\WD$ of wiring diagrams,
\item the operad $\wddot$ of normal wiring diagrams, and 
\item the operad $\wdzero$ of strict wiring diagrams.  
\end{enumerate}
A normal wiring diagram is a wiring diagram without delay nodes.  A strict wiring diagram is a normal wiring diagram  whose supplier assignment is a bijection.  For each of these three operads, we prove a finite presentation theorem that describes the operad in terms of a few operadic generators and a small number of generating relations.

Operads and wiring diagrams are recalled in Chapter \ref{ch02-wiring-diagrams}.  Operadic generators and generating relations for the operad $\WD$ of wiring diagrams are presented in Chapter \ref{ch03-generating-wd}.  Various decompositions of wiring diagrams are given in Chapter \ref{ch04-decomposition}.  Using the results in Chapters \ref{ch03-generating-wd} and \ref{ch04-decomposition}, the finite presentation theorems for the operads $\WD$, $\wddot$, and $\wdzero$ are proved in Chapter \ref{ch05-stratified-presentation}.  In Chapter \ref{ch06-wd-algebras} we prove the corresponding finite presentation theorems for $\WD$-algebras, $\wddot$-algebras, and $\wdzero$-algebras and discuss the main examples of the propagator algebra, the algebra of discrete systems, and the algebra of open dynamical systems.  Each finite presentation theorem for algebras describes the algebras in terms of finitely many generating structure maps and relations among these maps.

\textbf{Reading Guide}: In this reading guide we describe what can be safely skipped in Part 1 during the first reading.  The purpose is to help the reader get to the main results and examples quicker without getting bogged down by all the technical details.

In Chapter \ref{ch02-wiring-diagrams}, the reader who already knows about colored operads and categories may skip Section \ref{sec:colored-operads} and start reading about wiring diagrams from Def. \ref{def:Fins}.  In Section \ref{sec:operad-structure} about the operad structure on $\WD$, the reader may wish to skip the proofs of the Lemmas and just study the pictures.  Section \ref{sec:internal-wasted} about internal wasted wires may be skipped during the first reading.

In Chapter \ref{ch04-decomposition} the various decompositions of wiring diagrams are outlined in the introduction.  The reader may read that introduction, followed by Motivation \ref{mot:psi-alpha-phi} and \ref{mot:phi-2cells} and Examples \ref{ex:factoring-pi} and \ref{ex:factor-pitwo}, which provide pictures that illustrate the decompositions.

In Section \ref{sec:coherence-wd}, after the initial definitions and examples, the reader may wish to skip the proofs of Lemmas \ref{lemma:simplex-to-stratified}, \ref{lemma:stratified-type1}, and \ref{lemma:stratified-type2} and go straight to Theorem \ref{thm:wd-generator-relation}, the finite presentation theorem for wiring diagrams.

The reader who already knows about operad algebras may skip Section \ref{sec:operad-algebras}.

\chapter{Wiring Diagrams}
\label{ch02-wiring-diagrams}

This purpose of this chapter is to recall the definitions of colored operads and wiring diagrams.  In Section \ref{sec:colored-operads} we recall two equivalent definitions of a colored operad, one in terms of the structure map $\gamma$ \eqref{operadic-composition} and the other in terms of the $\compi$-compositions \eqref{operadic-compi}.  

Wiring diagrams are defined in Section \ref{sec:wiring-diagrams}.  The main difference between our definition of a wiring diagram and the original one in \cite{rupel-spivak} is that we allow the wires to carry values in an arbitrary class $S$ instead of just the class of pointed sets.  This added flexibility will be important in later chapters when we discuss operad algebras.  Indeed, in Section \ref{sec:propagator} when we discuss the propagator algebra, we will take $S$ to be the class of pointed sets.  On the other hand, in Section \ref{sec:algebra-ods} when we discuss the algebra of open dynamical systems, we will take $S$ to be a set of representatives of isomorphism classes of second-countable smooth manifolds.

In Section \ref{sec:operad-structure} we define the operad structure on wiring diagrams in terms of $\compi$-compositions.  Although we could also have defined this operad structure in terms of $\gamma$ as in \cite{rupel-spivak}, we chose to use $\compi$-compositions because the finite presentation theorems in Chapter \ref{ch05-stratified-presentation} can be phrased and proved more easily using the latter.

\section{Colored Operads}
\label{sec:colored-operads}

For brief discussion about \index{class}\emph{classes} in the set-theoretic sense, the reader is referred to \cite{halmos,pinter}.  In this monograph, the reader can safely take the word \emph{class} to just mean a collection of objects, such as sets, pointed sets, and real functions. First we introduce some notations for the colors in a colored operad.

\begin{definition}
\label{def:profile}
Suppose $S$\label{notation:s-class} is a class.  
\begin{enumerate}
\item
Denote by \index{profile}$\profs$\label{notation:profs} the class of finite ordered sequences of elements in $S$.  An element in $\profs$ is called an \emph{$S$-profile} or just a \emph{profile} if $S$ is clear from the context.  
\item
A typical $S$-profile of length\index{length of a profile} $n$ is denoted by $\us = (s_1, \ldots, s_n)$\label{notation:us} with $|\us|$ denoting its length.  
\item
The empty $S$-profile\index{empty profile} is denoted by $\varnothing$. 
\item
For $n \geq 0$ denote by $\profns \subseteq \profs$\label{notation:profns} the sub-class of $S$-profiles of length at least $n$.
\end{enumerate}
\end{definition}

\begin{motivation}\label{mot:operad}
Before we define an operad, let us first motivate its definition with a simple but important example.  Suppose $X$ is a set and $\Map(X^n,X)$ is the set of functions from $X^n = X \times \cdots \times X$, with $n \geq 0$ factors of $X$, to $X$.  If $f \in \Map(X^n,X)$ with $n \geq 1$ and $g_i \in \Map(X^{m_i},X)$ for each $1 \leq i \leq n$, then one can form the new function
\[f \circ (g_1,\ldots, g_n) \in \Map\bigl(X^{m_1 + \cdots + m_n},X\bigr)\]
by first applying the $g_i's$ simultaneously and then applying $f$.  Moreover, we may even allow the inputs and the output of each function to be from different sets, i.e., functions $X_{c_1} \times \cdots \times X_{c_n} \to X_d$.  In this case, the above composition is defined if and only if the outputs of the $g_i$'s match with the inputs of $f$.  

A function $f \in \Map(X_{c_1} \times \cdots \times X_{c_n},X_d)$ may be depicted as follows.
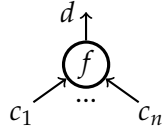
\begin{center}
\begin{tikzpicture}
\matrix[row sep=1cm, column sep=1cm]{
\node [plain, label=below:$...$] (f) {$f$};\\ };
\draw [outputleg] (f) to node[at end]{$d$} +(0,.7cm);
\draw [inputleg] (f) to node[below left=.1cm]{$c_1$} +(-.7cm,-.5cm);
\draw [inputleg] (f) to node[below right=.1cm]{$c_n$} +(.7cm,-.5cm);
\end{tikzpicture}
\end{center}
When $n \geq 1$, the composition $f \circ (g_1, \ldots, g_n)$ corresponds to the $2$-level tree:
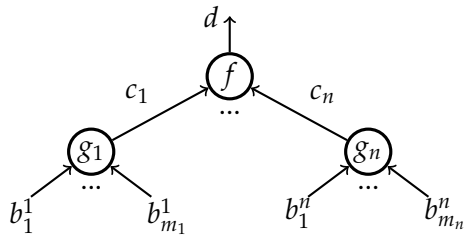
\begin{center}
\begin{tikzpicture}
\matrix[row sep=.05cm, column sep=1.2cm]{
& \node [plain, label=below:$...$] (f) {$f$}; &\\
\node [plain, label=below:$...$] (g1) {$g_1$}; &&
\node [plain, label=below:$...$] (gn) {$g_n$};\\ };
\draw [outputleg] (f) to node[at end]{$d$} +(0,.8cm);
\draw [arrow] (g1) to node{$c_1$} (f);
\draw [arrow] (gn)  to node[swap]{$c_n$} (f);
\draw [inputleg] (g1) to node[below left=.1cm]{$b^1_1$} +(-.8cm,-.6cm);
\draw [inputleg] (g1) to node[below right=.1cm]{$b^1_{m_1}$} +(.8cm,-.6cm);
\draw [inputleg] (gn) to node[below left=.1cm]{$b^n_1$} +(-.8cm,-.6cm);
\draw [inputleg] (gn) to node[below right=.1cm]{$b^n_{m_n}$} +(.8cm,-.6cm);
\end{tikzpicture}
\end{center}
Here each $g_i \in \Map\bigl(X_{b^i_1} \times \cdots \times X_{b^i_{m_i}}, X_{c_i}\bigr)$, $n \geq 1$, and each $m_i \geq 0$.  Together with permutations of the inputs, the above collection of functions satisfies some associativity, unity, and equivariance conditions.  An operad is an abstraction of this example that allows one to encode operations with multiple, possibly zero, inputs and one output and their compositions.
\end{motivation}

With the above motivation in mind, next we define colored operads.  See, for example, \cite{yau-operad} for more in-depth discussion of colored operads.  For each integer $n \geq 0$, the symmetric group on $n$ letters is denoted by \index{symmetric group}$\Sigma_n$.\label{notation:sigma-n}

\begin{definition}
\label{def:colored-operad}
Suppose $S$ is a class.  An \emph{$S$-colored operad} \index{colored operad} \index{operad}$(\O, \unit, \gamma)$\label{notation:colored-operad-gamma} consists of the following data.
\begin{enumerate}
\item
For any $d \in S$ and $\uc \in \profs$ with length $n \geq 0$, $\O$ is equipped with a class\label{notation:oduc}
\[\O\duc = \O\smallbinom{d}{c_1,\ldots,c_n}\]
called the \emph{entry of $\O$}\index{entry of an operad} with  \emph{input profile}\index{input profile} $\uc$ and \emph{output color}\index{output color} $d$. An element in $\O\duc$ is called an \emph{$n$-ary element}\index{n-ary element@$n$-ary element} in $\O$.
\item
For $\duc \in \profs \times S$ as above and a permutation $\sigma \in \Sigma_n$, $\O$ is equipped with a bijection
\begin{equation}
\label{operad-right-equivariance}
\nicexy{\O\duc \ar[r]^-{\sigma}_-{\cong} 
& \O\dcsigma}
\end{equation}
called the \emph{right action}\index{right action} or the \index{symmetric group action}\emph{symmetric group action}, in which\label{notation:c-sigma}
\[\uc\sigma = (c_{\sigma(1)}, \ldots, c_{\sigma(n)})\]
is the right permutation\index{right permutation} of $\uc$ by $\sigma$.
\item
For each $c \in S$, $\O$ is equipped with a specific element \label{notation:unit-c}$\unit_c \in \O\ccsingle$, called the \index{colored unit}\emph{$c$-colored unit}.
\item
For $\duc \in \profs \times S$ as above with $n \geq 1$, suppose $\ub_1, \ldots, \ub_n \in \profs$ and $\ub = (\ub_1,\ldots,\ub_n) \in \profs$ is their \index{concatenation}concatenation.  Then $\O$ is equipped with a map\label{notation:operad-composition-gamma}
\begin{equation}
\label{operadic-composition}
\nicexy{
\O\duc \times \prod\limits_{i=1}^n \O\ciubi \ar[r]^-{\gamma}
& \O\dub}
\end{equation}
called the \index{operadic composition}\emph{operadic composition}.  For $y \in \O\duc$ and $x_i \in \O\ciubi$ for $1 \leq i \leq n$, the image of the operadic composition is written as
\[\gamma\bigl(y; x_1, \ldots, x_n\bigr) \in \O\dub.\]
\end{enumerate}
This data is required to satisfy the following associativity, unity, and equivariance axioms.
\begin{description}
\item[Associativity Axiom]
Suppose that:\index{associativity of an operad}
\begin{itemize}
\item
in \eqref{operadic-composition}
\[\ub_j = \left(b^j_1, \ldots , b^j_{k_j}\right) \in \profs\]
has length $k_j \geq 0$ for each $1 \leq j \leq n$ such that at least one $k_j > 0$;
\item
$\ua^j_i \in \profs$ for each $1 \leq j \leq n$ and $1 \leq i \leq k_j$;
\item
for each $1 \leq j \leq n$, 
\begin{equation}
\label{ua-sub-j}
\ua_j = 
\begin{cases}
\left(\ua^j_1, \ldots , \ua^j_{k_j}\right)
& \text{if $k_j > 0$},\\
\varnothing & \text{if $k_j = 0$};
\end{cases}
\end{equation}
\item
$\ua = (\ua_1,\ldots , \ua_n)$ is their concatenation.
\end{itemize}
Then the \emph{associativity diagram}
\begin{equation}
\label{operad-associativity}
\nicexy{\O\duc \times 
\left[\prod\limits_{j=1}^n \O\cjubj\right] 
\times \prod\limits_{j=1}^n 
\left[\prod\limits_{i=1}^{k_j} \O\bjiuaji\right] 
\ar[r]^-{(\gamma, \Id)} 
\ar[d]_{\text{permute}}^-{\cong} 
&
\O\dub \times \prod\limits_{j=1}^{n}
\left[\prod\limits_{i=1}^{k_j} \O\bjiuaji\right] 
\ar[dd]^{\gamma}
\\
\O\duc \times 
\prod\limits_{j=1}^n \left[\O\cjubj 
\times \prod\limits_{i=1}^{k_j} \sO\bjiuaji\right] 
\ar[d]_{(\Id, \smallprod_j \gamma)} 
&
\\
\O\duc \times 
\prod\limits_{j=1}^n \O\cjasupj 
\ar[r]^-{\gamma} 
& \O\dua}
\end{equation}
is commutative.
\item[Unity Axioms]
Suppose $d \in S$.\index{unity of an operad}
\begin{enumerate}
\item
If $\uc = (c_1,\ldots,c_n) \in \profs$ has length $n \geq 1$, then the \emph{right unity diagram}\index{right unity}
\begin{equation}
\label{right-unity}
\nicexy{
\O\duc \times \{*\}^{n}
\ar[d]_-{(\Id, \smallprod \tensorunit_{c_j})} 
\ar[r]^-{\cong}
& 
\O\duc \ar[d]^-{=}
\\
\O\duc \times \prod\limits_{j=1}^n \O\cjcj
\ar[r]^-{\gamma}
&
\O\duc}
\end{equation}
is commutative.  Here $\{*\}$ is the one-point set, and $\{*\}^n$ is its $n$-fold product.
\item
If $\ub \in \profs$, then the \index{left unity}\emph{left unity diagram}
\begin{equation}
\label{left-unity}
\nicexy{
\{*\} \times \O\dub
\ar[d]_-{(\tensorunit_d, \Id)} 
\ar[r]^-{\cong}
& 
\O\dub \ar[d]^-{=}
\\
\O\dd \times \O\dub
\ar[r]^-{\gamma}
&
\O\dub}
\end{equation}
is commutative.
\end{enumerate}
\item[Equivariance Axioms]
Suppose that in \eqref{operadic-composition} $|\ub_j| = k_j \geq 0$.\index{equivariance of an operad}
\begin{enumerate}
\item
For each permutation $\sigma \in \Sigma_n$, the \index{top equivariance}\emph{top equivariance diagram} 
\begin{equation}
\label{operadic-eq-1}
\nicexy@C+.3cm{\O\duc \times \prod\limits_{j=1}^n \O\cjubj 
\ar[d]_-{\gamma} \ar[r]^-{(\sigma, \sigma^{-1})}
& \sO\ducsigma \times 
\prod\limits_{j=1}^n \O\csigmajubsigmaj \ar[d]^-{\gamma}
\\
\sO\duboneubn \ar[r]^-{\sigma\langle k_1, \ldots , k_n\rangle}
& \sO\dubsigmaoneubsigman}
\end{equation}
is commutative.  Here $\sigma\langle k_1, \ldots , k_n \rangle \in \Sigma_{k_1+\cdots+k_n}$\label{notation:block-permutation} is the block permutation\index{block permutation} induced by $\sigma$ that permutes the $n$ consecutive blocks of lengths $k_1, \ldots, k_n$, leaving the relative order within each block unchanged.
\item
Given permutations $\tau_j \in \Sigma_{k_j}$ for $1 \leq j \leq n$, the \index{bottom equivariance}\emph{bottom equivariance diagram}
\begin{equation}
\label{operadic-eq-2}
\nicexy@C+.3cm{\O\duc \times \prod\limits_{j=1}^n \O\cjubj
\ar[d]_-{\gamma} \ar[r]^-{(\Id, \smallprod \tau_j)}
& 
\O\duc \times \prod\limits_{j=1}^n \O\cjubjtauj
\ar[d]^-{\gamma}
\\
\O\duboneubn \ar[r]^-{\tau_1 \oplus \cdots \oplus \tau_n}
& \sO\dubonetauoneubntaun}
\end{equation}
is commutative.  Here the block sum\index{block sum} $\tau_1 \oplus \cdots \oplus \tau_n \in \Sigma_{k_1+\cdots+k_n}$\label{notation:block-sum} is the image of $(\tau_1, \ldots, \tau_n)$ under the inclusion $\Sigma_{k_1} \times \cdots \times \Sigma_{k_n} \to \Sigma_{k_1 + \cdots + k_n}$.
\end{enumerate}
\end{description}
A \emph{colored operad} is a $C$-colored operad for some class $C$.  We will also call a colored operad simply as an \emph{operad}.
\end{definition}

\begin{remark}
 A $1$-colored operad \cite{kelly05,mss,may72,may97}, with $S = \{*\}$, is usually called an \emph{operad} or a \index{symmetric operad}\emph{symmetric operad}.  The \index{underline notation}underline notation for $\uc \in \profs$ and the \index{vertical notation}vertical notation for $\duc \in \profs \times S$ originated in \cite{jy1}.  See \cite{yau-operad} (section 9.6) for more discussion of these notations.  In a $1$-colored operad $\sO$, the entry of $\sO$ whose input profile has length $n$ is denoted $\sO(n)$.
\end{remark}

\begin{example}[Endomorphism Operads]\label{ex:endomorphism}
\begin{enumerate}
\item Each set $X$ yields a $1$-colored operad $\End(X)$, called the \index{endomorphism operad} \emph{endomorphism operad}, whose $n$th entry is the set $\Map(X^n,X)$ of functions with the operad structure in Motivation \ref{mot:operad}.
\item More generally, for a non-empty class $\fC$, suppose $X = \{X_c\}_{c\in \fC}$ is a $\fC$-indexed class of sets.  Then there is a $\fC$-colored endomorphism operad whose $\duc$-entry is the set of functions $\Map(X_{c_1} \times \cdots \times X_{c_n}, X_d)$ with the operad structure in Motivation \ref{mot:operad}. 
\end{enumerate}
\end{example}

\begin{example}[Monoids as Operads]\label{ex:monoid}
Recall that a \index{monoid}\emph{monoid} is a triple $(A,\mu,1)$ consisting of a set $A$, a binary operation $\mu : A \times A \to A$, and an element $1 \in A$ such that the following two conditions are satisfied.
\begin{description}
\item[Associativity] $\mu\bigl(\mu(a,b),c\bigr) = \mu\bigl(a,\mu(b,c)\bigr)$ for all $a,b,c \in A$.
\item[Unity] $\mu(1,a) = a = \mu(a,1)$.
\end{description}
A monoid is uniquely determined by the operad $\sO$ with $\sO(1) = A$, $\sO(n) = \varnothing$ for all $n \not= 1$, $\mu$ as the only non-trivial operadic composition, and $1 \in A$ as the unit. 
\end{example}

\begin{example}[Associative and Commutative Operads]\label{ex:as-com}
\begin{enumerate}
\item There is a $1$-colored operad \label{notation:As} $\As$, called the \index{associative operad}\emph{associative operad}, whose $n$th entry is the symmetric group $\Sigma_n$ for each $n \geq 0$ with group multiplication as the symmetric group action and the identity $e \in \Sigma_1$ as the unit.  Given permutations $\sigma \in \Sigma_n$ with $n \geq 1$ and $\tau_i \in \Sigma_{k_i}$ for each $1 \leq i \leq n$, the operadic composition is given by
\[\gamma\bigl(\sigma; \tau_1, \ldots, \tau_n\bigr) = \sigma\langle k_1,\ldots, k_n\rangle \comp (\tau_1 \oplus \cdots \oplus \tau_n) \in \Sigma_{k_1 + \cdots + k_n}.\]
Here $\sigma\langle k_1, \ldots , k_n \rangle \in \Sigma_{k_1+\cdots+k_n}$ is the block permutation in \eqref{operadic-eq-1}, and $\tau_1 \oplus \cdots \oplus \tau_n \in \Sigma_{k_1+\cdots+k_n}$ is the block sum in \eqref{operadic-eq-2}.  The algebras of the operad $\As$, in the sense of Def. \ref{def1:operad-algebra}, are precisely monoids.
\item There is a $1$-colored operad \label{notation:Com} $\Com$, called the \index{commutative operad}\emph{commutative operad}, whose $n$th entry is a single point $*$ for each $n \geq 0$.  Its operad structure maps are all trivial maps, and its algebras, in the sense of Def. \ref{def1:operad-algebra}, are precisely commutative monoids, i.e., monoids whose multiplication maps are commutative.
\end{enumerate}
\end{example}

\begin{example}[Operad of Graph Operations]\label{ex:graph-op}
We now discuss an operad from non-commutative probability theory \cite{male}.  By a \emph{finite directed graph}, or just a \index{graph}\emph{graph}, we mean a pair of finite sets $(V,E)$ with $V$ non-empty such that an element in $E$ is an ordered pair $(u,v) \in V^2$, where each such ordered pair may appear in $E$ more than once.  Elements in $V$ and $E$ are called vertices and edges, respectively, and an edge $e = (u,v)$ is said to have initial vertex $u$ and terminal vertex $v$, denoted $e : u \to v$.  An edge of the form $(v,v)$ is called a \index{loop}loop at $v$.  An edge of the form $(u,v)$ or $(v,u)$ is said to \emph{connect} $u$ and $v$.  We say that a graph $(V,E)$ is \index{conncted graph}\emph{connected} if for each pair of distinct vertices $u$ and $v$, there exist edges $e_i$ for $1 \leq i \leq n$ for some $n \geq 1$ such that each $e_i$ connects $v_{i-1}$ and $v_i$ with $v_0 = u$ and $v_n = v$.

A \index{graph operation}\emph{graph operation} is a connected graph $(V,E)$ equipped with
\begin{enumerate}
\item an ordering $\sigma$ of the set $E$ of edges and
\item two possibly equal vertices $\inp$ and $\out$, called the \emph{input} and the \emph{output}.  
\end{enumerate}
An  isomorphism of graph operations is a pair of bijections $(V,E) \to (V',E')$ on vertices and edges that preserves the initial and the terminal vertices of each edge, the ordering on edges, and the input and the output.  We only consider graph operations up to isomorphisms.  That is, if there is an isomorphism
\[(V,E,\sigma,\inp,\out) \to (V',E',\sigma',\inp',\out')\]
of graph operations, then we consider them to be the same.  For each $n \geq 0$, denote by \label{notation:grop} $\GrOp_n$ the set of graph operations with $n$ edges.  So $\GrOp_0$ contains only the graph with one vertex, which is both the input and the output, and no edges.  Here is a graph operation with two vertices and four edges, two of which are loops:
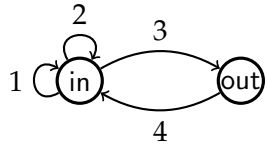
\begin{center}
\begin{tikzpicture}
\matrix[row sep=1cm, column sep=1.5cm]{
\node [plain] (in) {$\inp$}; & \node [plain] (out) {$\out$};\\};
\draw [arrow, out=210, in=150, looseness=4] (in) to node{$1$} (in);
\draw [arrow, out=120, in=60, looseness=4] (in) to node{$2$} (in);
\draw [arrow, bend left] (in) to node{$3$} (out);
\draw [arrow, bend left] (out) to node{$4$} (in);
\end{tikzpicture}
\end{center}

There is an operad structure on graph operations given by \index{edge substitution}\emph{edge substitution} as follows.  Suppose $G \in \GrOp_n$ with $n \geq 1$ and $G_i \in \GrOp_{m_i}$ for $1 \leq i \leq n$.  Then the operadic composition
\[G(G_1,\ldots,G_n) \in \GrOp_{m_1 + \cdots + m_n}\]
is obtained from $G$ by replacing the $i$th edge $e_i$ in $G$ by $G_i$ and by identifying the initial (resp., terminal) vertex of $e_i$ with the input (resp., output) of $G_i$.  The edge ordering of the operadic composition is induced by those of $G$ and of the $G_i$'s.  The input and the output are inherited from $G$.  The symmetric group action on $\GrOp_n$ is given by permutation of the edge ordering.  The operadic unit is the graph operation $\inp \to \out$ with two vertices and one edge from the input to the output.

For example, suppose $G$, $H$, and $K$ are the following graph operations in $\GrOp_2$:
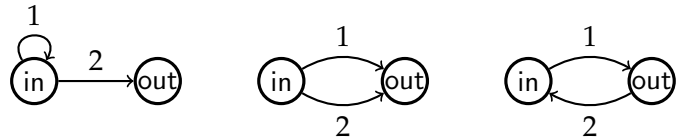
\begin{center}
\begin{tikzpicture}
\matrix[row sep=1cm, column sep=1cm]{
\node [plain] (ing) {$\inp$}; & \node [plain] (outg) {$\out$};
& \node [plain] (inh) {$\inp$}; & \node [plain] (outh) {$\out$};
& \node [plain] (ink) {$\inp$}; & \node [plain] (outk) {$\out$};\\};
\draw [arrow, out=120, in=60, looseness=4] (ing) to node{$1$} (ing);
\draw [arrow] (ing) to node{$2$} (outg);
\draw [arrow, bend left] (inh) to node{$1$} (outh);
\draw [arrow, bend right] (inh) to node[swap]{$2$} (outh);
\draw [arrow, bend left] (ink) to node{$1$} (outk);
\draw [arrow, bend left] (outk) to node{$2$} (ink);
\end{tikzpicture}
\end{center}
Then the operadic composition $G(H,K)$ is the graph operation with four edges above.  The algebras of the operad $\GrOp$ are closely related to traffic spaces and non-commutative probability, as we will discuss in Example \ref{ex:traffic} below.
\end{example}

Due to the presence of the colored units, a colored operad can also be defined in terms of a binary product, called the $\compi$-composition, which leads to axioms that are sometimes easier to check in practice and that we will use in most of the rest of this monograph.  In the one-colored linear setting, this alternative formulation of an operad was first made explicit in \cite{markl96}.  

\begin{motivation}
Using functions as in Motivation \ref{mot:operad}, the operadic $\compi$-compositions can be visualized as follows.  Suppose
\[f \in \Map(X_{c_1} \times \cdots \times X_{c_n}, X_d) \andspace 
g \in \Map(X_{b_1} \times \cdots \times X_{b_k}, X_{c_i})\]
are functions for some $1 \leq i \leq n$, pictorially depicted as follows.
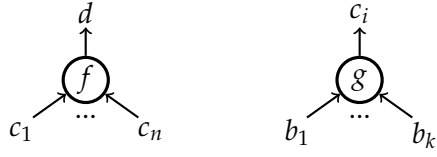
\begin{center}
\begin{tikzpicture}
\matrix[row sep=1cm, column sep=3cm]{
\node [plain, label=below:$...$] (f) {$f$}; 
& \node [plain, label=below:$...$] (g) {$g$};\\ };
\draw [outputleg] (f) to node[above=.1cm]{$d$} +(0,.7cm);
\draw [inputleg] (f) to node[below left=.1cm]{$c_1$} +(-.7cm,-.5cm);
\draw [inputleg] (f) to node[below right=.1cm]{$c_n$} +(.7cm,-.5cm);
\draw [outputleg] (g) to node[above=.1cm]{$c_i$} +(0,.7cm);
\draw [inputleg] (g) to node[below left=.05cm]{$b_1$} +(-.7cm,-.5cm);
\draw [inputleg] (g) to node[below right=.1cm]{$b_k$} +(.7cm,-.5cm);
\end{tikzpicture}
\end{center}
Then their $\compi$-composition
\[f \compi g = f \comp \bigl(\Id^{i-1}, g, \Id^{n-i}\bigr)\]
is the picture
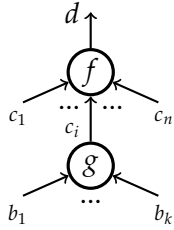
\begin{center}
\label{f-compi-g-pic}
\begin{tikzpicture}
\matrix[row sep=.6cm, column sep=1cm]{
\node [plain] (f) {$f$};\\
\node [plain, label=below:$...$] (g) {$g$}; \\};
\draw [arrow] (g) to node[near start]{\scriptsize{$c_i$}} 
node[near end]{...} node[near end, swap]{...} (f);
\draw [outputleg] (f) to node[at end]{$d$} +(0,.8cm);
\draw [inputleg] (f) to node[below left=.15cm]{\scriptsize{$c_1$}} +(-.9cm,-.4cm);
\draw [inputleg] (f) to node[below right=.15cm]{\scriptsize{$c_n$}} +(.9cm,-.4cm);
\draw [inputleg] (g) to node[below left=.15cm]{\scriptsize{$b_1$}} +(-.9cm,-.4cm);
\draw [inputleg] (g) to node[below right=.15cm]{\scriptsize{$b_k$}} +(.9cm,-.4cm);
\end{tikzpicture}
\end{center}
in which the output of $g$ is used as the $i$th input of $f$.  The operadic $\compi$-composition is an abstraction of this $f \compi g$ of functions.
\end{motivation}

With the above motivation in mind, we now recall this alternative formulation of an operad.  The elementary relations in Section \ref{sec:elementary-relations} are almost all stated in terms of the $\compi$-compositions in the following definition.

\begin{definition}
\label{def:pseudo-operad}
Suppose $S$ is a class.  An \emph{$S$-colored operad}\index{colored operad} \index{operad} $\left(\O, \tensorunit, \comp\right)$\label{notation:colored-operad-comp} consists of the following data.
\begin{enumerate}
\item
It has the same data as in (1)--(3) in Def. \ref{def:colored-operad}.
\item
For each  $d \in S$, $\uc = (c_1, \ldots , c_n) \in \profs$ with length $n \geq 1$, $\ub \in \profs$, and $1 \leq i \leq n$, it is equipped with a map\label{notation:operadic-compi}
\begin{equation}
\label{operadic-compi}
\nicexy{\O\duc \times \O\ciub
\ar[r]^-{\compi} 
& \O\dccompib}
\end{equation}
called the \index{compi-composition@$\compi$-composition} \index{operadic composition} \index{composition} \emph{$\compi$-composition}, where
\begin{equation}
\label{compi-profile}
\uc \compi \ub 
= \bigl(\underbrace{c_1, \ldots, c_{i-1}}_{\emptyset ~\text{if}~ i=1}, \ub, \underbrace{c_{i+1}, \ldots, c_n}_{\emptyset ~\text{if}~ i=n}\bigr).
\end{equation}
\end{enumerate}
This data is required to satisfy the following associativity, unity, and equivariance axioms.  Suppose $d \in S$, $\uc = (c_1, \ldots , c_n) \in \profs$, $\ub \in \profs$ with length $|\ub| = m$, and $\ua \in \profs$ with length $|\ua| = l$.
\begin{description}
\item[Associativity Axioms]
There are two associativity axioms.\index{associativity of an operad}
\begin{enumerate}
\item
Suppose $n \geq 2$ and $1 \leq i < j \leq n$.  Then the \emph{horizontal associativity diagram}\index{horizontal associativity}
\begin{equation}
\label{compi-associativity}
\nicexy{
\O\duc \times \O\ciua \times \O\cjub
\ar[d]_-{\mathrm{permute}}^-{\cong} 
\ar[r]^-{(\compi, \Id)}
& 
\O\dccompia \times \O\cjub
\ar[dd]^-{\comp_{j-1+l}}
\\
\O\duc \times \O\cjub \times \O\ciua
\ar[d]_-{(\compj,\Id)} &
\\
\O\dccompjb \times \sO\ciua
\ar[r]^-{\compi}
& \O\dccompjbcompia =
\O\dccompiacompjminusonepluslb
}
\end{equation}
is commutative.
\item
Suppose $n,m \geq 1$, $1 \leq i \leq n$, and $1 \leq j \leq m$.  Then the  \emph{vertical associativity diagram}\index{vertical associativity}
\begin{equation}
\label{compi-associativity-two}
\nicexy{\O\duc \times \O\ciub \times \O\bjua
\ar[d]_-{(\compi, \Id)} \ar[r]^-{(\Id, \compj)}
&
\O\duc \times \O\cibcompja \ar[d]^-{\compi}
\\
\O\dccompib \times \O\bjua
\ar[r]^-{\comp_{i-1+j}}
& 
\O\dccompibcompioneja
=
\O\dccompibcompja}
\end{equation}
is commutative.
\end{enumerate}
\item[Unity Axioms]
There are two unity axioms.\index{unity of an operad}
\begin{enumerate}
\item
The \emph{left unity diagram}\index{left unity}
\begin{equation}
\label{compi-left-unity}
\nicexy{\{*\} \times \O\duc 
\ar[dr]_-{\cong} \ar[r]^-{(\tensorunit_d, \Id)}
&
\O\dd \times \O\duc \ar[d]^-{\compone}
\\
&
\O\duc}
\end{equation}
is commutative.
\item
If $n \geq 1$ and $1 \leq i \leq n$, then the \emph{right unity diagram}\index{right unity}
\begin{equation}
\label{compi-right-unity}
\nicexy{\O\duc \times \{*\}
\ar[dr]_-{\cong} \ar[r]^-{(\Id, \tensorunit_{c_i})}
& 
\O\duc \times \O\cici
\ar[d]^-{\compi}
\\
& \O\duc}
\end{equation}
is commutative.
\end{enumerate}
\item[Equivariance Axiom]
Suppose $|\uc| = n \geq 1$, $1 \leq i \leq n$, $\sigma \in \Sigma_n$, and $\tau \in \Sigma_m$.  Then the \emph{equivariance diagram}\index{equivariance of an operad}
\begin{equation}
\label{compi-eq}
\nicexy{
\O\duc \times \O\csigmaiub
\ar[d]_-{(\sigma,\tau)} \ar[r]^-{\comp_{\sigma(i)}}
& 
\O\dccompsigmaib
\ar[d]^-{\sigma \compi \tau}
\\
\O\ducsigma \times \O\csigmaiubtau
\ar[r]^-{\compi}
& 
\O\dcsigmacompibtau
=
\O\dccompsigmaibsigmatau}
\end{equation}
is commutative, where\label{notation:sigma-compi-tau}
\[
\sigma \compi \tau
= \underbrace{\sigma\langle1^{i-1},m,1^{n-i}\rangle}_{\scriptsize\text{block permutation}} 
\underbrace{\bigl(\id^{i-1} \oplus \tau \oplus \id^{n-i}\bigr)}_{\scriptsize\text{block sum}} 
\in \Sigma_{n+m-1}.
\]
On the right side, the block sum permutes the interval $[i,i+m-1]$ via $\tau$.  The block permutation induced by $\sigma$ regards the interval $[i,i+m-1]$ as a single block, within which the relative order is unchanged.
\end{description}
Each $\compi$-composition is also called an \emph{operadic composition} or just a \emph{composition}.
\end{definition}

\begin{remark}
Without the symmetric group action \eqref{operad-right-equivariance} and the equivariance axiom \eqref{compi-eq}, a non-symmetric colored operad as in Def. \ref{def:pseudo-operad} is exactly a \emph{multicategory} \index{multicategory} as defined by Lambek \cite{lambek} (p.103-105).  In \cite{lambek}:
\begin{enumerate}
\item
The elements in an entry $\O\duc$ of a colored operad $\O$ are called \emph{multimaps}.\index{multimaps}
\item
The $\compi$-composition \eqref{operadic-compi} is called a \index{substitution}\emph{substitution} or a \emph{cut}.\index{cut}
\item
The horizontal associativity axiom \eqref{compi-associativity} is called the \emph{commutative law}.\index{commutative law}
\item
The vertical associativity axiom \eqref{compi-associativity-two} is called the \emph{associative law}.\index{associative law}
\end{enumerate}
The reader is cautioned that there are several typographical errors in \cite{lambek} (p.104-105) in the definition of a multicategory. 
\end{remark}

Due to the presence of the colored units, the two definitions of a colored operad are in fact equivalent.

\begin{proposition}
\label{prop:operad-def-equiv}
Definitions \ref{def:colored-operad} and \ref{def:pseudo-operad} of an $S$-colored operad are equivalent.
\end{proposition}

\begin{proof}
In the $1$-colored case, a proof can be found in \cite{markl08} (Prop. 13).  For the general colored case, the proof is similar and can be found in, e.g., \cite{yau-operad} (section 16.4).  The correspondence goes as follows.  Given the operadic composition $\gamma$ \eqref{operadic-composition}, $y \in \O\duc$ with $|\uc| = n \geq 1$, and $x \in \O\ciub$ with $1 \leq i \leq n$, one defines the $\compi$-composition as
\[
y \compi x = \gamma\Bigl(y; \underbrace{\tensorunit_{c_1}, \ldots, \tensorunit_{c_{i-1}}}_{\text{$\varnothing$ if $i=1$}}, x, \underbrace{\tensorunit_{c_{i+1}}, \ldots, \tensorunit_{c_n}}_{\text{$\varnothing$ if $i=n$}}\Bigr).
\]
Conversely, given the $\compi$-compositions \eqref{operadic-compi}, $y \in \O\duc$ with $|\uc| = n \geq 1$, and $x_i \in \O\ciubi$ for $1 \leq i \leq n$ with $k_i = |\ub_i|$, one defines the operadic composition $\gamma$ as
\begin{equation}
\label{gamma-in-comps}
\gamma\bigl(y;x_1, \ldots, x_n\bigr) =
\Bigl(\bigl((y \compone x_1) \comp_{k_1+1} x_2\bigr) \cdots\Bigr) \comp_{k_1+\cdots+k_{n-1}+1} x_n.
\end{equation}
On the right side, every pair of parentheses starts on the left.
\end{proof}

\section{Defining Wiring Diagrams}
\label{sec:wiring-diagrams}

The operad of wiring diagrams $\WD$ has $\boxs$ as its class of colors.  So before we define wiring diagrams, we first define $\boxs$.  We begin by recalling the basic definition of a \emph{category}.  The reader may consult \cite{awodey,leinster,maclane} for more in-depth discussion of category theory.  In this monograph, we do not need anything fancy from category theory.  All that the reader needs to know is that a category consists of a collection of objects, such as sets or finite sets, and maps between them that can be composed and that satisfy some natural unity and associativity axioms with respect to compositions.

\begin{motivation}\label{mot:category}
To motivate the definition of a category, consider the collection of groups and group homomorphisms.  Given group homomorphisms $f : G_1 \to G_2$ and $g : G_2 \to G_3$, there is a composition $g \comp f : G_1 \to G_3$.  The identity map of each group serves as the identity for composition.  Furthermore, composition of group homomorphisms is associative, in the sense that given a group homomorphism $h : G_3 \to G_4$, there is an equality 
\[(h \comp g) \comp f=h \comp (g \comp f)\]
of group homomorphisms.  One can think of a category as an abstraction of the collection of groups, group homomorphisms, their composition, and the unity and the associativity axioms governing composition.
\end{motivation}

\begin{definition}
\label{def:category}
A \emph{category} $\C$\index{category} consists of the following data.
\begin{description}
\item[Objects] It is equipped with a collection $\Ob(\C)$ of \emph{objects}.\index{object}  For an object $a$ in $\C$, we will write either $a \in \Ob(\C)$ or simply $a \in \C$.
\item[Morphisms] For any objects $a,b \in \Ob(\C)$, it is equipped with a collection $\C(a,b)$ of \emph{morphisms}\index{morphism} with domain $a$ and codomain $b$.  A morphism $f$ in $\C(a,b)$ is usually denoted by $f : a \to b$ and is also called a \emph{map} from $a$ to $b$.
\item[Composition] For any objects $a,b,c \in \Ob(\C)$ and morphisms $f : a \to b$ and $g : b \to c$, it is equipped with a morphism $g \comp f : a \to c$, called the \emph{composition} of $f$ and $g$. 
\item[Identities] Each object $a \in \Ob(\C)$ is equipped with a morphism $\tensorunit_a : a \to a$, called the \emph{identity morphism of $a$}. 
\end{description}
The above data is required to satisfy the following axioms.
\begin{description}
\item[Unity] For any objects $a,b,c \in \Ob(\C)$ and morphisms $g : a \to b$ and $h : c \to a$, there are equalities of morphisms\index{unity of a category}
\[g \comp \tensorunit_a = g \in \C(a,b) 
\andspace \tensorunit_a \comp h = h \in \C(c,a).\]
\item[Associativity]
For any objects $a,b,c,d \in \Ob(\C)$ and morphisms $f : a \to b$, $g : b \to c$, and $h : c \to d$, there is an equality of morphisms\index{associativity of a category}
\[(h \comp g) \comp f = h \comp (g \comp f)\]
in $\C(a,d)$.
\end{description}
\end{definition}

\begin{example}
Here are some basic examples of categories.  In each case, the identity morphisms and the composition are the obvious ones.
\begin{enumerate}
\item There is an empty category with no objects and no morphisms.
\item There is a category $*$ with only one object $*$ and only the identity morphism of $*$.
\item There is a category $\set$ whose objects are sets and whose morphisms are functions between sets.
\item There is a category $\Fin$ whose objects are finite sets and whose morphisms are functions between finite sets.  Given any disjoint finite sets $X_1,\ldots, X_n$, their \emph{coproduct} $\coprod_{i=1}^n X_i$ is the finite set given by their disjoint union.  If the $X_i$'s are not disjoint, we can still define their coproduct, but we must first replace each $X_i$ (or just the ones with $i \geq 2$) by an $X_i'$ equipped with a bijection to $X_i$ such that the resulting $X_1', \ldots, X_n'$ are disjoint.  Then the coproduct $\coprod_{i=1}^n X_i$ is defined as the disjoint union of the $X_i'$ for $1 \leq i \leq n$, and it is well-defined up to isomorphism.  If $n=0$ then the coproduct is defined as the empty set $\varnothing$.
\end{enumerate}
In what follows, if the identity morphisms and the composition are obvious, then we will omit mentioning them. 
\end{example}

\begin{example}
A monoid $(A,\mu,1)$ (Example \ref{ex:monoid}) is a category with one object, whose only morphism set is $A$.  Composition and identity are those of $A$, i.e., the multiplication $\mu$ and the unit element $1$.  Therefore, one can think of a category as a \index{monoid with multiple objects}\emph{monoid with multiple objects}. 
\end{example}

\begin{example}
Suppose $\sO$ is an $S$-colored operad (Def. \ref{def:colored-operad}). Then $\sO$ determines a category $\C$ whose objects are the elements in $S$.  For $c,d \in S$ the morphism object $\C(c,d)$ is $\sO\dcsingle$, and the identity morphism of $c$ is the $c$-colored unit of $\sO$.  Composition in $\C$ is the restriction of the operadic composition in $\sO$.  Therefore, one can think of a colored operad as a generalization of a category in which the domain of each morphism is a finite sequence of objects. 
\end{example}

For wiring diagrams, we will usually consider finite sets in which each element is allowed to carry a value of some kind.  The precise notion is given in the following definition.

\begin{definition}
\label{def:Fins}
Suppose $S$ is a non-empty class, and $\Fin$\label{notation:fin}\index{Fin@$\Fin$} is the category of finite sets\index{category of finite sets} and functions between them.  
\begin{enumerate}
\item
Denote by $\Fins$\label{notation:fins}\index{Fin@$\Fins$} the category in which:
\begin{itemize}
\item
an object is a pair $(X,v)$ with $X \in \Fin$ and $v : X \to S$ a function;
\item
a map $(X,v_X) \to (Y,v_Y)$ is a function $X \to Y$ such that the diagram
\begin{equation}
\label{maps-fins}
\nicexy{
X \ar[r] \ar[dr]_-{v_X} & Y \ar[d]^-{v_Y}
\\
& S
}
\end{equation}
is commutative.
\end{itemize}
\item
An object $(X,v) \in \Fins$ is called an \emph{$S$-finite set}.\index{S-finite set@$S$-finite set}
\item
For $(X,v) \in \Fins$, we call $v$ the \emph{value assignment}\index{value assignment} for $X$.  For each $x \in X$, $v(x) \in S$ is called the \emph{value of $x$}.
\item
If $(X_i, v_i) \in \Fins$ for $1 \leq i \leq n$, then their \emph{coproduct}\index{coproduct of finite sets} $X = \coprod_{i=1}^n X_i \in \Fins$ has value assignment $\coprod_{i=1}^n v_i$.
\item
The empty $S$-finite set is denoted by $\varnothing$.
\end{enumerate}
\end{definition}

\begin{definition}
\label{def:s-box}
Suppose $S$ is a non-empty class.  
\begin{enumerate}
\item
An \emph{$S$-box}\index{S-box@$S$-box}\index{box} is a pair \label{notation:sbox}$X = (\xin, \xout) \in \Fins \times \Fins$.  If $S$ is clear from the context, then we will drop $S$ and call $X$ a \emph{box}. 
\begin{itemize}
\item
An element of $\xin$\label{notation:xin} is called an \emph{input}\index{input} of $X$.
\item
An element of $\xout$\label{notation:xout} is called an \emph{output}\index{output} of $X$.
\item
We write $v = v_X : \xin \amalg \xout \to S$\label{notation:value-assignment} for the value assignment for $X$.
\end{itemize} 
\item
The class of $S$-boxes is denoted by $\boxs$.
\item
The empty $S$-box, denoted by \label{notation:emptybox}$\emptyset$, is the $S$-box with $\emptyset^{\inp} = \emptyset^{\out}$ the empty set.
\end{enumerate}
\end{definition}

\begin{convention}
From now on, whenever $\Fins$ or $\boxs$ is used, we always assume that the class $S$ is non-empty.
\end{convention}

\begin{remark}
In \cite{spivak15} (Def. 3.1) an $S$-box is called a \emph{signed finite set}.\index{signed finite set}  It is a slight generalization of what appears in \cite{rupel-spivak,vsl}.  In \cite{rupel-spivak} $S$ is the class of pointed sets, where an $S$-box is called a \emph{black box}.\index{black box}  In \cite{vsl} $S$ is a set of representatives of isomorphism classes of second-countable smooth manifolds and smooth maps, or more generally the class of objects in a category with finite products.  
\end{remark}

\begin{convention}
\label{conv:box}
For the purpose of visualization, an $S$-box $X$ will be drawn as follows.
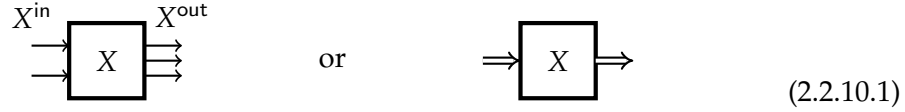
\begin{equation}
\label{box-picture}
\begin{tikzpicture}
\draw [ultra thick] (2,1) rectangle (3,2);
\node at (2.5,1.5) {$X$};
\draw [arrow] (1.5,1.7) -- (2,1.7);
\draw [arrow] (1.5,1.3) -- (2,1.3);
\draw [arrow] (3,1.7) -- (3.5,1.7);
\draw [arrow] (3,1.5) -- (3.5,1.5);
\draw [arrow] (3,1.3) -- (3.5,1.3);
\node at (1.5,2.1) {$\xin$};
\node at (3.5,2.1) {$\xout$};
\node at (5.5,1.5) {or};
\draw [ultra thick] (8,1) rectangle (9,2);
\draw [implies] (7.5,1.5) -- (8,1.5);
\draw [implies] (9,1.5) -- (9.5,1.5);
\node at (8.5,1.5) {$X$};
\end{tikzpicture}
\end{equation}
On the left, the inputs of $X$ are depicted as arrows going into the box, and the outputs of $X$ are depicted as arrows leaving the box.  Alternatively, if we do not need to specify the sizes of $\xin$ and $\xout$, then we depict them using a generic arrow $\Rightarrow$, as in the picture on the right.  The value of each $x \in \xin \amalg \xout$ is either not depicted in the picture for simplicity, or it is written near the corresponding arrow.
\end{convention}

\begin{definition}
Suppose $X_1, \ldots, X_n$ are $S$-boxes for some $n \geq 0$.  Define the $S$-box $X = \coprod_{i=1}^n X_i$, called the \emph{coproduct},\index{coproduct of boxes} as follows.
\begin{enumerate}
\item
If $n=0$, then $X = \emptyset$, the empty $S$-box.
\item
If $n \geq 1$, then:
\[\xin = \coprod_{i=1}^n \xin_i,\quad
\xout = \coprod_{i=1}^n \xout_i, \andspace
v_X = \coprod\limits_{i=1}^n v_{X_i}.\]
\end{enumerate}
\end{definition}

\begin{motivation}
Before we define wiring diagrams, let us first provide a motivating example.  A typical wiring diagram looks like this:
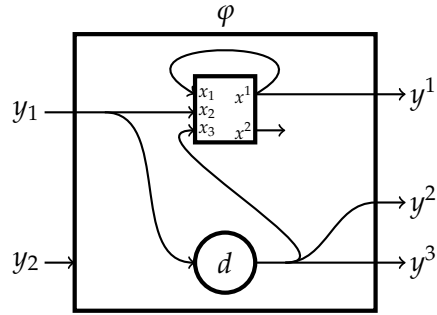
\begin{center}
\begin{tikzpicture}[scale=0.8]
\draw [ultra thick] (1,0.2) rectangle (6,4.8);
\node at (3.5,5.1) {$\varphi$};
\draw [ultra thick] (3,3) rectangle (4,4.1);
\node at (3.2,3.8) {\tiny{$x_1$}};
\node at (3.2,3.5) {\tiny{$x_2$}};
\node at (3.2,3.2) {\tiny{$x_3$}};
\node at (3.8,3.8) {\tiny{$x^1$}};
\node at (3.8,3.2) {\tiny{$x^2$}};
\draw [ultra thick] (3.5,1) circle [radius=0.5];
\node at (3.5,1) {$d$};
\draw [thick] (0.5,3.5) -- (1,3.5);
\node at (0.2,3.5) {$y_1$};
\draw [arrow] (0.5,1) -- (1,1);
\node at (0.2,1) {$y_2$};
\draw [arrow] (6,3.8) -- (6.5,3.8);
\node at (6.8,3.8) {$y^1$};
\draw [arrow] (6,2) -- (6.5,2);
\node at (6.8,2) {$y^2$};
\draw [arrow] (6,1) -- (6.5,1);
\node at (6.8,1) {$y^3$};
\draw [arrow] (1,3.5) -- (3,3.5);
\draw [arrow, thick] (1.5,3.5) to [out=0, in=180] (3,1);
\draw [thick] (4,1) -- (6,1);
\draw [thick] (4.5,1) to [out=0, in=180] (6,2);
\draw [arrow, looseness=1.2] (4.5,1) to [out=0, in=180] (3,3.2);
\draw [arrow] (4,3.2) -- (4.5,3.2);
\draw [thick] (4,3.8) -- (6,3.8);
\draw [arrow, looseness=5] (4,3.8) to [out=30, in=150] (3,3.8);
\end{tikzpicture}
\end{center}
There is an output box $Y$ (the outermost box in the picture), a finite number of input boxes $X$ (which the above picture only has one), and a finite number of delay nodes (which the above picture again only has one).  To make sense of such a picture, first of all, for each output $y^i$ of $Y$ we need to specify where the arrow is coming from.  In the example above, $y^1$ comes from $x^1$, and both $y^2$ and $y^3$ come from the delay node $d$.  We will say that each $y^i$ is a demand wire, and $y^1$ (resp., $y^2$ and $y^3$) is supplied by the supply wire $x^1$ (resp., $d$).

Similarly, for each input $x_i$ of $X$ and the delay node $d$, we again need to specify where the arrow is coming from.  So $d$ is also a demand wire, and by tracing the arrow ending at $d$ backward, we see that it is supplied by the supply wire $y_1$.  Starting at the input $x_1$ (resp., $x_2$ and $x_3$) and tracing the arrow backward, we see that it is supplied by $x^1$ (resp., $y_1$ and $d$).  We will come back to this example precisely in Example \ref{wd-first-example} below.

The above exercise tells us that the outputs of $Y$, the inputs of $X$, and the delay nodes are demand wires, in the sense that each of their elements demands a supply wire.  The supply wires consist of the inputs of $Y$, the outputs of $X$, and the delay nodes.  Each supply wire may supply multiple demand wires or none at all.  For instance, in the above example, the supply wire $x^1$ supplies both the demand wires $x_1$ and $y^1$.  On the other hand, the supply wires $y_2$ and $x^2$ do not supply to any demand wires at all.  We will call them wasted wires.

To avoid pathological situations, one requirement of a wiring diagram is that a demand wire in the outputs of $Y$ should \emph{not} be supplied by a supply wire in the inputs of $Y$.  In other words, we exclude pictures like this
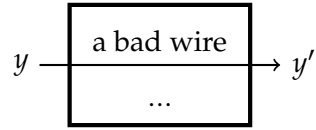
\begin{center}
\begin{tikzpicture}[scale=0.8]
\draw [ultra thick] (.5,0) rectangle (3.5,2);
\node at (2,1.4) {a bad wire};
\draw [arrow] (0,1) -- (4,1);
\node at (-.3,1) {$y$};
\node at (4.4,1) {$y'$};
\node at (2,.3) {$\cdots$};
\end{tikzpicture}
\end{center}
in which the demand wire $y'$ is directly supplied by the supply wire $y$.  So we insist that an output of $Y$ be supplied by either an output of an input box $X$ or a delay node $d$.  We will call this the non-instantaneity requirement.
\end{motivation}

Wiring diagrams will be defined as equivalence classes of prewiring diagrams, as defined below.  

\begin{definition}
\label{def:wiring-diagram}
Suppose $S$ is a class.  An \emph{$S$-prewiring diagram} \index{prewiring diagram} is a tuple
\begin{equation}
\label{phi-wd}
\varphi = \bigl(\uX, Y, \dn, v, s\bigr)
\end{equation}
consisting of the following data.
\begin{enumerate}
\item
$Y = (\yin, \yout) \in \boxs$, called the \emph{output box}\index{output box} of $\varphi$.  
\begin{itemize}
\item
An element in $\yin$ is called a \emph{global input}\index{global input} for $\varphi$.
\item
An element in $\yout$ is called a \emph{global output}\index{global output} for $\varphi$.
\end{itemize}
\item
$\uX = (X_1, \ldots, X_n)$ is a $\boxs$-profile for some $n \geq 0$; i.e., each $X_i \in \boxs$.  
\begin{itemize}
\item
We call $X_i$ the \emph{$i$th input box}\index{input box} of $\varphi$.  
\item
Denote by $X = \coprod_{i=1}^n X_i \in \boxs$ the coproduct.
\end{itemize}
\item
$(\dn,v) \in \Fins$\label{notation:delay-nodes}\index{DN@$\dn$} is an $S$-finite set.  An element of $\dn$ is called a \emph{delay node}.\index{delay node}  Define:
\begin{itemize}
\item
$\dm = \yout \amalg \xin \amalg \dn \in \Fins$.\label{notation:demand}\index{Dm@$\dm$}\index{demand}  An element of $\dm$ is called a \emph{demand wire} in $\varphi$.
\item
$\supply = \yin \amalg \xout \amalg \dn \in \Fins$.\label{notation:supply}\index{Sp@$\supply$}\index{supply}  An element of $\supply$ is called a \emph{supply wire} in $\varphi$.
\end{itemize}
Furthermore:
\begin{itemize}
\item
When $\dn$ is regarded as a subset of $\dm$, an element in $\xin \amalg \dn$ is called an \emph{internal input}\index{internal input} for $\varphi$.
\item
When $\dn$ is regarded as a subset of $\supply$, an element in $\xout \amalg \dn$ is called an \emph{internal output}\index{internal output} for $\varphi$.
\end{itemize} 
\item
With a slight abuse of notation, we write
\[
\nicexy{
\dm \coprod\limits_{\dn} \supply = \yin \amalg \yout \amalg \xin \amalg \xout \amalg \dn \ar[r]^-{v} &  S
}\]
for the coproduct of the value assignments for $X$, $Y$, and $\dn$.
\item
$s : \dm \to \supply \in \Fins$\label{notation:supplier} is a map, called the \emph{supplier assignment}\index{supplier assignment} for $\varphi$, such that
\begin{equation}
\label{non-instant}
y \in \yout \subseteq \dm \impliesspace
s(y) \in \supply \setminus \yin = \xout \amalg \dn.
\end{equation}
\begin{itemize}
\item
The condition \eqref{non-instant} is called the \emph{non-instantaneity requirement}.\index{non-instantaneity requirement}
\item
For $w \in \dm$, we call $s(w) \in \supply$ the \emph{supplier}\index{supplier} or the \emph{supply wire}\index{supply wire} of $w$.  So non-instantaneity says that the supply wire of a global output cannot be a global input.
\item
A supply wire $w \in \yin$ that does not belong to the image of the supplier assignment $s$ is called an \emph{external wasted wire}.\index{external wasted wire}  The set of external wasted wires in $\varphi$ is denoted by \label{notation:external-wasted}$\ewphi$.
\item
A supply wire $w \in \xout \amalg \dn$ that does not belong to the image of the supplier assignment $s$ is called an \emph{internal wasted wire}.\index{internal wasted wire}  The set of internal wasted wires in $\varphi$ is denoted by \label{notation:internal-wasted}$\iwphi$.
\end{itemize}
\end{enumerate}
If $S$ is clear from the context, then we will drop $S$ and call $\varphi$ a  \emph{prewiring diagram}.  If we need to emphasize $\varphi$, then we will use subscripts such as $\dm_{\varphi}$, $\supplyphi$, and $\sphi$.
\end{definition}

\begin{remark}
In the constructions that follow, the compatibility of the value assignments with the various supplier assignments are usually immediate because, in each prewiring diagram, the supplier assignment $s : \dm \to \supply$ is a map in $\Fins$.  Therefore, we will omit checking such compatibility.  
\end{remark}

\begin{definition}
\label{wd-equivalence}
Suppose $S$ is a class, $\varphi = (\uX,Y,\dn,v,s)$ is an $S$-prewiring diagram as in \eqref{phi-wd}, and $\varphi' = (\uX,Y,\dn',v',s')$ is another $S$-prewiring diagram with the same input boxes $\uX$ and output box $Y$.  
\begin{enumerate}
\item
An \emph{equivalence}\index{equivalence of prewiring diagrams} $f : \varphi \to \varphi'$ is an isomorphism $f_0 : (\dn,v) \to (\dn',v') \in \Fins$ such that the diagram
\[
\nicexy@C+.7cm{
\yout \amalg \xin \amalg \dn 
= \dm_{\varphi} \ar[d]_-{s_{\varphi}}
\ar[r]^-{\Id \amalg \Id \amalg f_0}
&
\dm_{\varphi'} = \yout \amalg \xin \amalg \dn' \ar[d]^-{s_{\varphi'}}
\\
\yin \amalg \xout \amalg \dn  = \supply_{\varphi}
\ar[r]^-{\Id \amalg \Id \amalg f_0}
&
\supply_{\varphi'} = \yin \amalg \xout \amalg \dn'
}\]
in $\Fins$ is commutative.
\item
We say that $\varphi$ and $\varphi'$ are \emph{equivalent} if there exists an equivalence $\varphi \to \varphi'$.  
\item
An \emph{$S$-wiring diagram} is an equivalence class of $S$-prewiring diagrams.  If $S$ is clear from the context, we will drop $S$ and just say \emph{wiring diagram}.\index{wiring diagram}
\item
The class of $S$-wiring diagrams with input boxes $\uX = (X_1, \ldots, X_n)$ and output box $Y$ is denoted by\label{notation:wd-yux}
\begin{equation}
\label{wdyux}
\WD\yux \orspace \WD\yxonexn.
\end{equation}
The class of all $S$-wiring diagrams is denoted by \index{WD@$\WD$}\label{notation:wd}$\WD$.  If we want to emphasize $S$, then we will write $\WD^S$.
\end{enumerate}
\end{definition}

\begin{convention}
\label{conv:prewiring}
To simplify the presentation, we usually suppress the difference between a prewiring diagram and a wiring diagram.  In any given wiring diagram, $\dn$ is a finite set in which each element $d$ is equipped with a value $v(d) \in S$.  We will suppress the difference between each delay node $d$ and its value $v(d)$. In other words, each wiring diagram has a unique representative in which:
\begin{itemize}
\item
each delay node is an element in $S$;
\item
$v : \dn \to S$ sends each delay node to itself.  
\end{itemize}
Unless otherwise specified, in the rest of this book, we will always use this representative of a wiring diagram.  For a wiring diagram $\varphi \in \WD\yux$, we will  sometimes draw it as
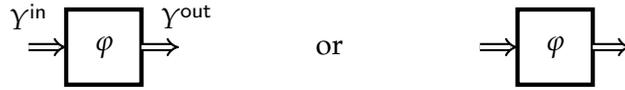
\begin{center}
\begin{tikzpicture}
\draw [ultra thick] (2,1) rectangle (3,2);
\draw [implies] (1.5,1.5) -- (2,1.5);
\draw [implies] (3,1.5) -- (3.5,1.5);
\node at (1.5,1.9) {$\yin$};
\node at (3.6,1.9) {$\yout$};
\node at (2.5,1.5) {$\varphi$};
\node at (5.5,1.5) {or};
\draw [ultra thick] (8,1) rectangle (9,2);
\draw [implies] (7.5,1.5) -- (8,1.5);
\draw [implies] (9,1.5) -- (9.5,1.5);
\node at (8.5,1.5) {$\varphi$};
\end{tikzpicture}
\end{center}
without drawing the input boxes $\uX$, the delay nodes, and the supplier assignment.
\end{convention}

\begin{remark}
Def. \ref{wd-equivalence} of an $S$-wiring diagram is a slight generalization of the one given in \cite{rupel-spivak}, where $S$ was taken to be the class of pointed sets.  Note that when $S$ is a proper class (e.g., the class of pointed sets), the collection $\WD\yux$ \eqref{wdyux} is also a proper class, not a set.  This is the reason why in Def. \ref{def:colored-operad} we allow $\O\duc$ to be a class, in contrast to what was stated in \cite{rupel-spivak} (Def. 2.1.2B).
\end{remark}

\begin{example}
Suppose $S$ is a non-empty class and:
\begin{itemize}
\item
$\xin = \{x_1,x_1,x_3\}$, $\xout = \{x^1,x^2\}$, $\yin = \{y_1,y_2\}$, $\yout = \{y^1,y^2, y^3\}$, and $\dn = \{d\}$;
\item
$v(x_2) = v(x_3) = v(y_1) = v(y^2) = v(y^3) = v(d) \in S$;
\item
$v(x_1) = v(x^1) = v(y^1) \in S$, and $v(x^2)$ and $v(y_2)$ in $S$ are arbitrary;
\item
$s(x_1) = s(y^1) = x^1$, $s(x_2) = s(d) = y_1$, and $s(x_3) = s(y^2) = s(y^3) = d$.
\end{itemize}
Then the above data defines a wiring diagram $\varphi \in \WD\yx$ with one input box $X$ and output box $Y$, which can be depicted as follows:
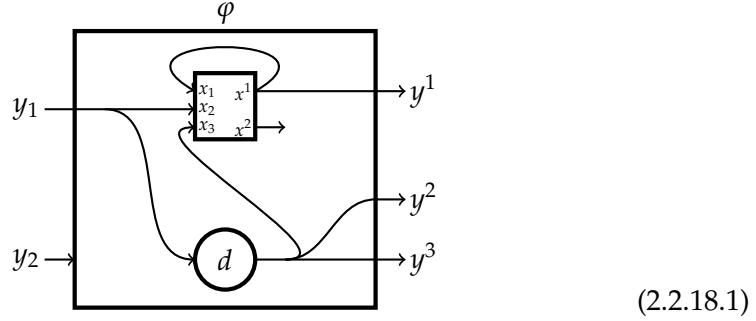
\begin{equation}
\label{wd-first-example}
\begin{tikzpicture}[scale=0.8]
\draw [ultra thick] (1,0.2) rectangle (6,4.8);
\node at (3.5,5.1) {$\varphi$};
\draw [ultra thick] (3,3) rectangle (4,4.1);
\node at (3.2,3.8) {\tiny{$x_1$}};
\node at (3.2,3.5) {\tiny{$x_2$}};
\node at (3.2,3.2) {\tiny{$x_3$}};
\node at (3.8,3.8) {\tiny{$x^1$}};
\node at (3.8,3.2) {\tiny{$x^2$}};
\draw [ultra thick] (3.5,1) circle [radius=0.5];
\node at (3.5,1) {$d$};
\draw [thick] (0.5,3.5) -- (1,3.5);
\node at (0.2,3.5) {$y_1$};
\draw [arrow] (0.5,1) -- (1,1);
\node at (0.2,1) {$y_2$};
\draw [arrow] (6,3.8) -- (6.5,3.8);
\node at (6.8,3.8) {$y^1$};
\draw [arrow] (6,2) -- (6.5,2);
\node at (6.8,2) {$y^2$};
\draw [arrow] (6,1) -- (6.5,1);
\node at (6.8,1) {$y^3$};
\draw [arrow] (1,3.5) -- (3,3.5);
\draw [arrow, thick] (1.5,3.5) to [out=0, in=180] (3,1);
\draw [thick] (4,1) -- (6,1);
\draw [thick] (4.5,1) to [out=0, in=180] (6,2);
\draw [arrow, looseness=1.2] (4.5,1) to [out=0, in=180] (3,3.2);
\draw [arrow] (4,3.2) -- (4.5,3.2);
\draw [thick] (4,3.8) -- (6,3.8);
\draw [arrow, looseness=5] (4,3.8) to [out=30, in=150] (3,3.8);
\end{tikzpicture}
\end{equation}
The output box $Y$ is drawn as the outermost box.  The single input box $X$ is the smaller box.  The delay node $d$ is drawn as a circle, which will be our convention from now on.  The supply wires $y_2 \in \yin$ and $x^2 \in \xout$ are not in the image of the supplier assignment $s : \dm \to \supply$, so $y_2$ (resp., $x^2$) is an external (resp., internal) wasted wire.
\end{example}

\section{Operad Structure}
\label{sec:operad-structure}

Throughout this section, $S$ denotes a class.  In this section, using Def. \ref{def:pseudo-operad} of a colored operad, we define the colored operad structure on the collection $\WD$ of $S$-wiring diagrams (Def. \ref{wd-equivalence}).  

The equivariant structure is given by permuting the labels of the input boxes.

\begin{definition}[Equivariance in $\WD$]
\label{wd-equivariance}
Suppose $Y \in \boxs$, $\uX$ is a $\boxs$-profile of length $n$, and $\sigma \in \Sigma_n$ is a permutation.  Define the function\index{equivariance in $\WD$}
\begin{equation}
\label{wd-permutation}
\nicexy{\WD\yxonexn = \WD\yux \ar[r]^-{\sigma}_-{\cong} 
& \WD\yuxsigma = \WD\smallbinom{Y}{X_{\sigma(1)}, \ldots, X_{\sigma(n)}}}
\end{equation}
by sending $\varphi = (\dn,v,s) \in \WD\yux$ to $\varphi = (\dn,v,s) \in \WD\yuxsigma$, using the fact that $\coprod_{i=1}^n X_i = \coprod_{i=1}^n X_{\sigma(i)} \in \boxs$.
\end{definition}

Next we define the colored units.  For a box $Y$, the $Y$-colored unit can be depicted as follows.
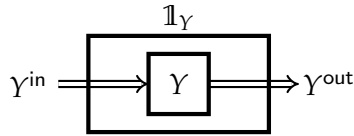
\begin{center}
\begin{tikzpicture}[scale=.8]
\draw [implies] (0.5,1.5) -- (2,1.5);
\draw [implies] (3,1.5) -- (4.5,1.5);
\draw [ultra thick] (2,1) rectangle (3,2);
\draw [ultra thick] (1,0.7) rectangle (4,2.3);
\node at (2.5,1.5) {$Y$};
\node at (0,1.5) {$\yin$};
\node at (5,1.5) {$\yout$};
\node at (2.5,2.6) {$\tensorunit_Y$};
\end{tikzpicture}
\end{center}

\begin{definition}[Units in $\WD$]
\label{wd-units}
For each $Y \in \boxs$, the \emph{$Y$-colored unit} \index{units in $\WD$} is defined as the wiring diagram
\begin{equation}
\label{wd-unit}
\tensorunit_Y = \Bigl(\dn = \varnothing, s = \Id\Bigr) \in \WD\yy
\end{equation}
with:
\begin{itemize}
\item
no delay nodes and trivial $v : \dn \to S$;
\item
supplier assignment
\[\nicexy{\dm = \yout \amalg \yin \ar[r]^-{\Id}
& \yin \amalg \yout = \supply}\]
the identity function.
\end{itemize}
\end{definition}

\begin{motivation}
Before we define the $\compi$ composition  \eqref{operadic-compi} in $\WD$, let us first describe the intuition behind the definition.  Pictorially, the $\compi$ composition $\varphi \compi \psi$ in $\WD$ replaces the $i$th input box $X_i$ in $\varphi$ with the wiring diagram $\psi$, using $\xin_i$ and $\xout_i$ for the necessary wire connections.  Here is a picture to keep in mind for $\phicompipsi$.
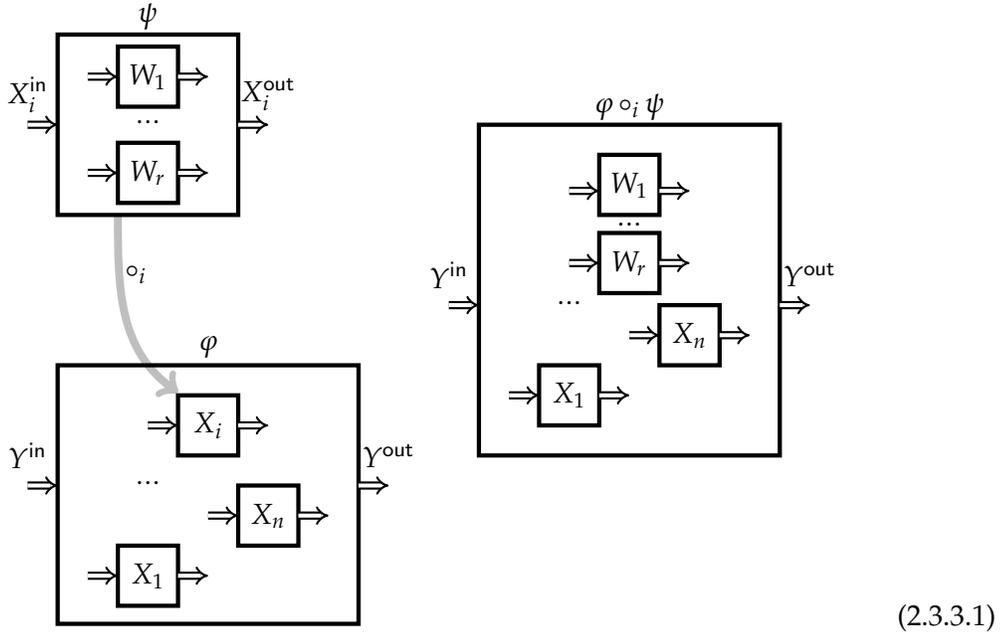
\begin{equation}
\label{wd-compi-picture}
\begin{tikzpicture}[scale=.8]
\draw [lightgray, ->, line width=3pt] (1,6) to [out=-90, in=145] (2,3.05);
\node at (1.3,5) {$\compi$};
\draw [ultra thick] (0,6) rectangle (3,9);
\draw [ultra thick] (1,7.8) rectangle (2,8.8);
\draw [ultra thick] (1,6.2) rectangle (2,7.2);
\draw [implies] (-0.5,7.5) -- (0,7.5);
\draw [implies] (3,7.5) -- (3.5,7.5);
\draw [implies] (0.5,8.3) -- (1,8.3);
\draw [implies] (2,8.3) -- (2.5,8.3);
\draw [implies] (0.5,6.7) -- (1,6.7);
\draw [implies] (2,6.7) -- (2.5,6.7);
\node at (1.5,9.3) {$\psi$};
\node at (1.5,8.3) {$W_1$};
\node at (1.5,6.7) {$W_r$};
\node at (-0.5,8) {$\xin_i$};
\node at (3.5,8) {$\xout_i$};
\node at (1.5,7.5) {$\cdots$};
\draw [ultra thick] (0,-.8) rectangle (5,3.5);
\draw [ultra thick] (2,2) rectangle (3,3);
\draw [ultra thick] (1,-.5) rectangle (2,.5);
\draw [ultra thick] (3,.5) rectangle (4,1.5);
\draw [implies] (-0.5,1.5) -- (0,1.5);
\draw [implies] (5,1.5) -- (5.5,1.5);
\node at (-0.5,2) {$\yin$};
\node at (5.5,2) {$\yout$};
\node at (2.5,3.8) {$\varphi$};
\draw [implies] (0.5,0) -- (1,0);
\draw [implies] (2,0) -- (2.5,0);
\node at (1.5,0) {$X_1$};
\draw [implies] (1.5,2.5) -- (2,2.5);
\draw [implies] (3,2.5) -- (3.5,2.5);
\node at (2.5,2.5) {$X_i$};
\draw [implies] (2.5,1) -- (3,1);
\draw [implies] (4,1) -- (4.5,1);
\node at (3.5,1) {$X_n$};
\node at (1.5,1.5) {$\cdots$};
\node at (9.5,7.8) {$\phicompipsi$};
\draw [ultra thick] (7,2) rectangle (12,7.5);
\node at (6.5,5) {$\yin$};
\node at (12.5,5) {$\yout$};
\draw [implies] (6.5,4.5) -- (7,4.5);
\draw [implies] (12,4.5) -- (12.5, 4.5);
\node at (8.5,4.5) {$\cdots$};
\draw [ultra thick] (8,2.5) rectangle (9,3.5);
\node at (8.5,3) {$X_1$};
\draw [implies] (7.5,3) -- (8,3);
\draw [implies] (9,3) -- (9.5,3);
\draw [ultra thick] (10,3.5) rectangle (11,4.5);
\node at (10.5,4) {$X_n$};
\draw [implies] (9.5,4) -- (10,4);
\draw [implies] (11,4) -- (11.5,4);
\draw [ultra thick] (9,6) rectangle (10,7);
\node at (9.5,6.5) {$W_1$};
\draw [implies] (8.5,6.4) -- (9,6.4);
\draw [implies] (10,6.4) -- (10.5,6.4);
\draw [ultra thick] (9,4.7) rectangle (10,5.7);
\node at (9.5,5.2) {$W_r$};
\draw [implies] (8.5,5.2) -- (9,5.2);
\draw [implies] (10,5.2) -- (10.5,5.2);
\node at (9.5,5.8) {$\cdots$};
\end{tikzpicture}
\end{equation}
For simplicity we did not draw any delay nodes.  In the actual definition below, we will take the coproduct of the sets of delay nodes in $\varphi$ and $\psi$.
\end{motivation}

\begin{definition}[$\compi$-Composition in $\WD$]
\label{def:compi-wd}
Suppose:\index{compi-composition in WD@$\compi$-composition in $\WD$}
\begin{itemize}
\item
$\varphi = (\dnphi, \vphi, \sphi) \in \WD\yux$ with  $\uX = (X_1, \ldots, X_n)$, $n \geq 1$, and $1 \leq i \leq n$;
\item
$\psi = (\dnpsi, \vpsi, \spsi) \in \WD\xiuw$ with $|\uW| = r \geq 0$.  
\end{itemize}
Define the wiring diagram
\[
\phicompipsi \in \WD\yxcompiw
\]
as follows, where the $\boxs$-profile
\[
\xcompiw = \Bigl(X_1, \ldots, X_{i-1}, \uW, X_{i+1}, \ldots , X_n\Bigr)
\]
is as in \eqref{compi-profile}.
\begin{enumerate}
\item
$\dnphicompipsi = \dnphi \amalg \dnpsi \in \Fins$, so $\vphicompipsi = \vphi \amalg \vpsi$.
\item
The supplier assignment for $\phicompipsi$,
\begin{equation}
\label{compi-supply}
\nicexy{\dmphicompipsi = \yout \amalg \coprod\limits_{j\not= i} \xin_j \amalg
\coprod\limits_{k=1}^r \win_k \amalg \dnphi \amalg \dnpsi
\ar[d]_-{\sphicompipsi}\\
\supplyphicompipsi = \yin \amalg \coprod\limits_{j\not= i} \xout_j \amalg
\coprod\limits_{k=1}^r \wout_k \amalg \dnphi \amalg \dnpsi}
\end{equation}
in which the coproduct $\coprod_{j\not= i}$ is indexed by all $j \in \{1,\ldots,i-1,i+1,\ldots, n\}$, is defined as follows.  Suppose $z \in \dmphicompipsi$.
\begin{enumerate}
\item
If $z \in \yout \amalg \coprod_{j\not= i} \xin_j \amalg \dnphi \subseteq \dmphi$, then
\begin{equation}
\label{compi-supply1}
\sphicompipsi(z) = 
\begin{cases}
\sphi(z) & \text{ if $\sphi(z) \not\in \xout_i$};\\
\spsi\sphi(z) & \text{ if $\sphi(z) \in \xout_i$}.
\end{cases}
\end{equation}
\item
If $z \in \coprod_{k=1}^r \win_k \amalg \dnpsi \subseteq \dmpsi$, then
\begin{equation}
\label{compi-supply2}
\sphicompipsi(z) = 
\begin{cases}
\spsi(z) & \text{ if $\spsi(z) \not\in \xin_i$};\\
\sphi\spsi(z) & \text{ if $\spsi(z) \in \xin_i$ and $\sphi\spsi(z) \not\in \xout_i$};\\
\spsi\sphi\spsi(z) & \text{ if $\spsi(z) \in \xin_i$ and $\sphi\spsi(z) \in \xout_i$}.
\end{cases}
\end{equation}
\end{enumerate}
\end{enumerate}
This finishes the definition of $\phicompipsi$.
\end{definition}

\begin{lemma}
\label{compi-well-defined}
Def. \ref{def:compi-wd} indeed defines a wiring diagram $\phicompipsi$ in $\WD\yxcompiw$.
\end{lemma}

\begin{proof}
We need to check that the supplier assignment for $\phicompipsi$ satisfies the non-instantaneity requirement \eqref{non-instant}.  So suppose $y \in \yout$.  We must show that $\sphicompipsi(y) \not\in \yin$.  By \eqref{compi-supply1} we have
\[
\sphicompipsi(y) = 
\begin{cases}
\sphi(y) \in \coprod\limits_{j\not= i} \xout_j \amalg \dnphi 
& \text{ if $\sphi(y) \not\in \xout_i$};\\
\spsi\sphi(y) \in \coprod\limits_{k=1}^r \wout_k \amalg \dnpsi
& \text{ if $\sphi(y) \in \xout_i$}.
\end{cases}
\]
Here we have used the non-instantaneity requirement for both $\varphi$ and $\psi$. So in either case we have that  $\sphicompipsi(y) \not\in \yin$.
\end{proof}

Many examples of the $\compi$-composition in $\WD$ will be given in Chapter \ref{ch03-generating-wd}.  We now prove that, equipped with the structure above, $\WD$ is a colored operad.

\begin{lemma}
\label{compi-unity}
The $\compi$-composition in Def. \ref{def:compi-wd} satisfies the left unity axiom \eqref{compi-left-unity}, the right unity axiom \eqref{compi-right-unity}, and the equivariance axiom \eqref{compi-eq}.\index{unity in $\WD$}\index{equivariance in $\WD$}
\end{lemma}

\begin{proof}
This follows from a direct inspection of the definitions of the equivariant structure \eqref{wd-permutation} and the colored units in $\WD$ \eqref{wd-unit}.
\end{proof}

Next we show that $\WD$ satisfies the associativity axioms \eqref{compi-associativity} and \eqref{compi-associativity-two}.  The reader may wish to skip the proofs of the following two Lemmas and simply look at the pictures during the first reading.  

\begin{motivation}
For the horizontal associativity axiom \eqref{compi-associativity}, one should keep in mind the following picture for the iterated operadic composition $(\varphi \compj \zeta) \compi \psi$.
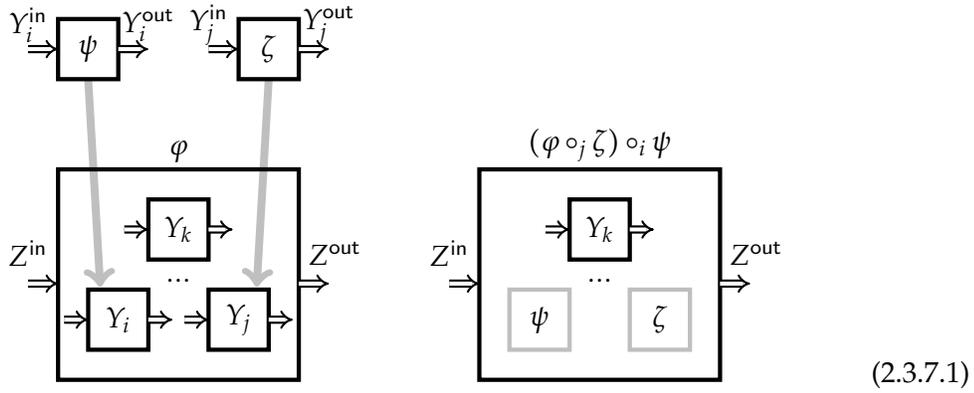
\begin{equation}
\label{wd-hass-picture}
\begin{tikzpicture}[scale=.8]
\draw [lightgray, ->, line width=3pt] (1.5,5) to (1.7,1.55);
\draw [lightgray, ->, line width=3pt] (4.5,5) to (4.3,1.55);
\draw [ultra thick] (1,5) rectangle (2,6);
\node at (0.5,5.9) {$\yin_i$};
\node at (2.5,5.9) {$\yout_i$};
\node at (1.5,5.5) {$\psi$};
\draw [implies] (.5,5.5) -- (1,5.5);
\draw [implies] (2,5.5) -- (2.5,5.5);
\draw [ultra thick] (4,5) rectangle (5,6);
\node at (3.5,5.9) {$\yin_j$};
\node at (5.5,5.9) {$\yout_j$};
\node at (4.5,5.5) {$\zeta$};
\draw [implies] (3.5,5.5) -- (4,5.5);
\draw [implies] (5,5.5) -- (5.5,5.5);
\draw [ultra thick] (1,0) rectangle (5,3.5);
\draw [implies] (.5,1.6) -- (1,1.6);
\draw [implies] (5,1.6) -- (5.5,1.6);
\node at (0.5,2.1) {$\zin$};
\node at (5.6,2.1) {$\zout$};
\node at (3,3.8) {$\varphi$};
\node at (3,1.6) {$\cdots$};
\draw [ultra thick] (1.5,.5) rectangle (2.5,1.5);
\draw [implies] (1.1,1) -- (1.5,1);
\draw [implies] (2.5,1) -- (2.9,1);
\node at (2,1) {$Y_i$};
\draw [ultra thick] (3.5,.5) rectangle (4.5,1.5);
\draw [implies] (3.1,1) -- (3.5,1);
\draw [implies] (4.5,1) -- (4.9,1);
\node at (4,1) {$Y_j$};
\draw [ultra thick] (2.5,2) rectangle (3.5,3);
\draw [implies] (2.1,2.5) -- (2.5,2.5);
\draw [implies] (3.5,2.5) -- (3.9,2.5);
\node at (3,2.5) {$Y_k$};
\draw [ultra thick] (8,0) rectangle (12,3.5);
\draw [implies] (7.5,1.6) -- (8,1.6);
\draw [implies] (12,1.6) -- (12.5,1.6);
\node at (7.5,2.1) {$\zin$};
\node at (12.6,2.1) {$\zout$};
\node at (10,3.9) {$(\varphi \compj \zeta) \compi \psi$};
\node at (10,1.6) {$\cdots$};
\draw [ultra thick, lightgray] (8.5,.5) rectangle (9.5,1.5);
\node at (9,1) {$\psi$};
\draw [ultra thick, lightgray]  (10.5,.5) rectangle (11.5,1.5);
\node at (11,1) {$\zeta$};
\draw [ultra thick] (9.5,2) rectangle (10.5,3);
\draw [implies] (9.1,2.5) -- (9.5,2.5);
\draw [implies] (10.5,2.5) -- (10.9,2.5);
\node at (10,2.5) {$Y_k$};
\end{tikzpicture}
\end{equation}
Note that on the right side, $\psi$ and $\zeta$ are depicted as gray boxes because their output boxes--namely $Y_i$ and $Y_j$--are no longer input boxes in $(\varphi \compj \zeta) \compi \psi$.  Furthermore, for simplicity the delay nodes are not drawn.
\end{motivation}

\begin{lemma}
\label{compi-horizontal-associative}
The $\compi$-composition in Def. \ref{def:compi-wd} satisfies the horizontal associativity axiom \eqref{compi-associativity}.\index{horizontal associativity in $\WD$}
\end{lemma}

\begin{proof}
Suppose:
\begin{itemize}
\item
$\varphi \in \WD\zuy$ with $|\uY| = n \geq 2$ and $1 \leq i < j \leq n$;
\item
$\psi \in \WD\yiuw$ with $|\uW| = l$;
\item
$\zeta \in \WD\yjux$ with $|\uX| = m$.
\end{itemize}
We must show that
\begin{equation}
\label{horizontal-associativity}
\left( \varphi \compj \zeta\right) \compi \psi
= 
\left(\varphi \compi \psi\right) \compjonel \zeta \in \WD\zycompjxcompiw.
\end{equation}
By Lemma \ref{compi-well-defined} we already know that both sides are well-defined wiring diagrams in the indicated entry of $\WD$.  Moreover, both sides have $\dnphi \amalg \dnpsi \amalg \dnzeta \in \Fins$ as the set of delay nodes.  So it remains to show that their supplier assignments are equal.

Note that both sides in \eqref{horizontal-associativity} have demand wires
\[\dm = \zout \amalg \coprod_{p\not= i,j} \yin_p \amalg \coprod_{q=1}^l \win_q \amalg \coprod_{r=1}^m \xin_r \amalg \dnphi \amalg \dnpsi \amalg \dnzeta\]
in which the coproduct $\coprod_{p\not= i,j}$ is indexed by all $1 \leq p \leq n$ such that $p\not= i,j$.  Similarly, both sides in  \eqref{horizontal-associativity} have supply wires
\[\supply = \zin \amalg \coprod_{p\not=i, j} \yout_p \amalg \coprod_{q=1}^l \wout_q \amalg \coprod_{r=1}^m \xout_r \amalg \dnphi \amalg \dnpsi \amalg \dnzeta.\]

Using the definitions \eqref{compi-supply1} and \eqref{compi-supply2}, it follows from direct inspection that both sides in \eqref{horizontal-associativity} have the following supplier assignment $s : \dm \to \supply$.
\begin{enumerate}
\item
If $v \in \zout \amalg \coprod_{p\not=i,j} \yin_p \amalg \dnphi \subseteq \dmphi$, then 
\[s(v) = \begin{cases}
\sphi(v) & \text{ if $\sphi(v) \not\in \yout_i \amalg \yout_j$};\\
\spsi\sphi(v) & \text{ if $\sphi(v) \in \yout_i$};\\
\szeta\sphi(v) & \text{ if $\sphi(v) \in \yout_j$}.
\end{cases}\]
\item
If $v \in \coprod_{q=1}^l \win_q \amalg \dnpsi \subseteq \dmpsi$, then
\[s(v) = \begin{cases}
\spsi(v) & \text{ if $\spsi(v) \not\in \yin_i$};\\
\sphi\spsi(v) & \text{ if $\spsi(v) \in \yin_i$ and $\sphi\spsi(v) \not\in \yout_i \amalg \yout_j$};\\
\spsi\sphi\spsi(v) & \text{ if $\spsi(v) \in \yin_i$ and $\sphi\spsi(v) \in \yout_i$};\\
\szeta\sphi\spsi(v) & \text{ if $\spsi(v) \in \yin_i$ and $\sphi\spsi(v) \in \yout_j$}.
\end{cases}\]
\item
Finally, if $v \in \coprod_{r=1}^m \xin_r \amalg \dnzeta \subseteq \dmzeta$, then
\[s(v) = \begin{cases}
\szeta(v) & \text{ if $\szeta(v) \not\in \yin_j$};\\
\sphi\szeta(v) & \text{ if $\szeta(v) \in \yin_j$ and $\sphi\szeta(v) \not\in \yout_i \amalg \yout_j$};\\
\spsi\sphi\szeta(v) & \text{ if $\szeta(v) \in \yin_j$ and $\sphi\szeta(v) \in \yout_i$};\\
\szeta\sphi\szeta(v) & \text{ if $\szeta(v) \in \yin_j$ and $\sphi\szeta(v) \in \yout_j$}.
\end{cases}\]
\end{enumerate}
This finishes the proof of the desired equality \eqref{horizontal-associativity}.
\end{proof}

\begin{motivation}
For the vertical associativity axiom, one should keep the following picture of $\varphi \compi (\psi \compj \zeta)$ in mind.
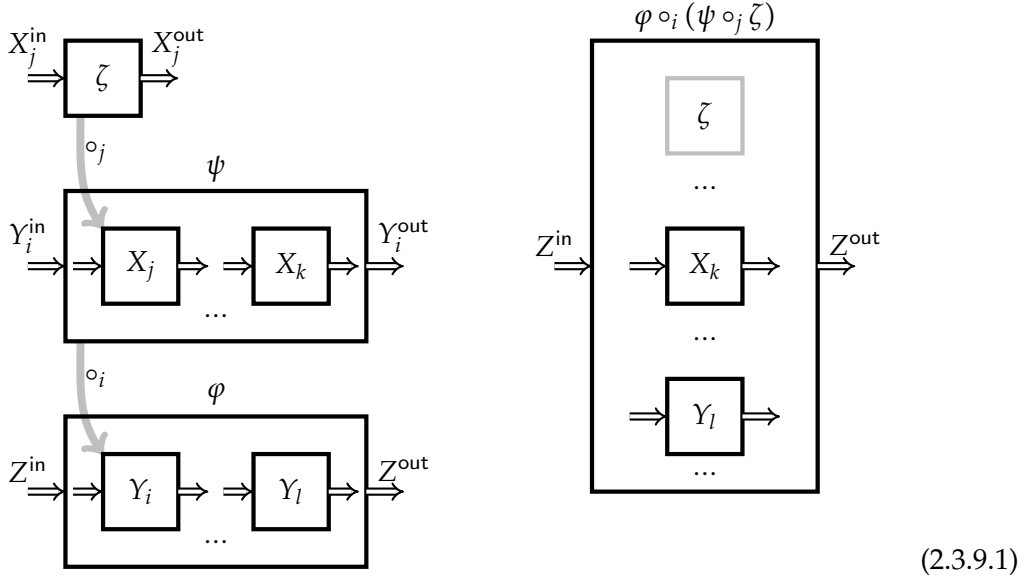
\begin{equation}
\label{wd-vass-picture}
\begin{tikzpicture}[scale=1]
\draw [lightgray, ->, line width=3pt] (1.2,5) to [out=-90, in=125] (1.5,3.5);
\node at (1.4,4.5) {$\compj$};
\draw [lightgray, ->, line width=3pt] (1.2,2) to [out=-90, in=125]  (1.5,.5);
\node at (1.4,1.5) {$\compi$};
\draw [ultra thick] (1,5) rectangle (2,6);
\node at (0.5,5.9) {$\xin_j$};
\node at (2.5,5.9) {$\xout_j$};
\node at (1.5,5.5) {$\zeta$};
\draw [implies] (.5,5.5) -- (1,5.5);
\draw [implies] (2,5.5) -- (2.5,5.5);
\draw [ultra thick] (1,2) rectangle (5,4);
\draw [implies] (.5,3) -- (1,3);
\draw [implies] (5,3) -- (5.5,3);
\node at (0.5,3.4) {$\yin_i$};
\node at (5.5,3.4) {$\yout_i$};
\node at (3,4.3) {$\psi$};
\node at (3,2.3) {$\cdots$};
\draw [ultra thick] (1.5,2.5) rectangle (2.5,3.5);
\draw [implies] (1.1,3) -- (1.5,3);
\draw [implies] (2.5,3) -- (2.9,3);
\node at (2,3) {$X_j$};
\draw [ultra thick] (3.5,2.5) rectangle (4.5,3.5);
\draw [implies] (3.1,3) -- (3.5,3);
\draw [implies] (4.5,3) -- (4.9,3);
\node at (4,3) {$X_k$};
\draw [ultra thick] (1,-1) rectangle (5,1);
\draw [implies] (.5,0) -- (1,0);
\draw [implies] (5,0) -- (5.5,0);
\node at (0.5,.3) {$\zin$};
\node at (5.5,.3) {$\zout$};
\node at (3,1.3) {$\varphi$};
\node at (3,-.7) {$\cdots$};
\draw [ultra thick] (1.5,-.5) rectangle (2.5,.5);
\draw [implies] (1.1,0) -- (1.5,0);
\draw [implies] (2.5,0) -- (2.9,0);
\node at (2,0) {$Y_i$};
\draw [ultra thick] (3.5,-.5) rectangle (4.5,0.5);
\draw [implies] (3.1,0) -- (3.5,0);
\draw [implies] (4.5,0) -- (4.9,0);
\node at (4,0) {$Y_l$};
\draw [ultra thick] (8,0) rectangle (11,6);
\node at (9.5,6.3) {$\varphi \compi (\psi \compj \zeta)$};
\draw [implies] (7.5,3) -- (8,3);
\draw [implies] (11,3) -- (11.5,3);
\node at (7.5,3.3) {$\zin$};
\node at (11.5,3.3) {$\zout$};
\node at (9.5,.2) {$\cdots$};
\node at (9.5,2) {$\cdots$};
\node at (9.5,4) {$\cdots$};
\draw [ultra thick, gray, semitransparent]  (9,4.5) rectangle (10,5.5);
\node at (9.5,5) {$\zeta$};
\draw [ultra thick] (9,2.5) rectangle (10,3.5);
\draw [implies] (8.5,3) -- (9,3);
\draw [implies] (10,3) -- (10.5,3);
\node at (9.5,3) {$X_k$};
\draw [ultra thick] (9,.5) rectangle (10,1.5);
\draw [implies] (8.5,1) -- (9,1);
\draw [implies] (10,1) -- (10.5,1);
\node at (9.5,1) {$Y_l$};
\end{tikzpicture}
\end{equation}
Once again on the right side, $\zeta$ is depicted as a gray box because its output box $X_j$ is no longer an input box in $\varphi \compi (\psi \compj \zeta)$.  Furthermore, for simplicity the delay nodes are not drawn.
\end{motivation}

\begin{lemma}
\label{compi-vertical-associative}
The $\compi$-composition in Def. \ref{def:compi-wd} satisfies the vertical associativity axiom \eqref{compi-associativity-two}.\index{vertical associativity in $\WD$}
\end{lemma}

\begin{proof}
Suppose:
\begin{itemize}
\item
$\varphi \in \WD \zuy$ with $|\uY| = n \geq 1$ and $1 \leq i \leq n$;
\item
$\psi \in \WD\yiux$ with $|\uX| = m \geq 1$ and $1 \leq j \leq m$;
\item
$\zeta \in \WD\xjuw$ with $|\uW| = l$.
\end{itemize}
We must show that
\begin{equation}
\label{vertical-assoc}
\left(\varphi \compi \psi\right) \compionej \zeta 
=
\varphi \compi \left(\psi \compj \zeta\right)
\in \WD\zycompixcompionejw.
\end{equation}
By Lemma \ref{compi-well-defined} we already know that both sides are well-defined wiring diagrams in the indicated entry of $\WD$.  Moreover, both sides have $\dnphi \amalg \dnpsi \amalg \dnzeta \in \Fins$ as the set of delay nodes.  So it remains to show that their supplier assignments are equal.

Note that both sides in \eqref{vertical-assoc} have demand wires
\[\dm = \zout \amalg \coprod_{p\not= i} \yin_p \amalg \coprod_{q\not=j} \xin_q  \amalg \coprod_{r=1}^l \win_r \amalg \dnphi \amalg \dnpsi \amalg \dnzeta\]
in which $\coprod_{p\not= i}$ is indexed by all $1 \leq p \leq n$ such that $p\not= i$, and $\coprod_{q\not= j}$ is indexed by all $1 \leq q \leq m$ such that $q \not= j$.  Similarly, both sides in  \eqref{vertical-assoc} have supply wires
\[\supply = \zin \amalg \coprod_{p\not=i} \yout_p \amalg \coprod_{q\not=j} \xout_q \amalg \coprod_{r=1}^l \wout_r  \amalg \dnphi \amalg \dnpsi \amalg \dnzeta.\]

Using the definitions \eqref{compi-supply1} and \eqref{compi-supply2}, it follows from direct inspection that both sides in \eqref{vertical-assoc} have the following supplier assignment $s : \dm \to \supply$.
\begin{enumerate}
\item
If $v \in \zout \amalg \coprod_{p\not= i} \yin_p \amalg \dnphi \subseteq \dmphi$, then
\[s(v) = \begin{cases}
\sphi(v) & \text{ if $\sphi(v) \not\in \yout_i$};\\
\spsi\sphi(v) & \text{ if $\sphi(v) \in \yout_i$ and $\spsi\sphi(v) \not\in \xout_j$};\\
\szeta\spsi\sphi(v) & \text{ if $\sphi(v) \in \yout_i$ and $\spsi\sphi(v) \in \xout_j$}.
\end{cases}\]
\item
If $v \in \coprod_{q\not=j} \xin_q \amalg \dnpsi \subseteq \dmpsi$, then
\[s(v) = \begin{cases}
\spsi(v) & \text{ if $\spsi(v) \not\in \xout_j \amalg \yin_i$};\\
\szeta\spsi(v) & \text{ if $\spsi(v) \in \xout_j$};\\
\sphi\spsi(v) & \text{ if $\spsi(v) \in \yin_i$ and $\sphi\spsi(v) \not\in \yout_i$};\\
\spsi\sphi\spsi(v) & \text{ if $\spsi(v) \in \yin_i$, $\sphi\spsi(v) \in \yout_i$, and $\spsi\sphi\spsi(v) \not\in \xout_j$};\\
\szeta\spsi\sphi\spsi(v) & \text{ if  $\spsi(v) \in \yin_i$, $\sphi\spsi(v) \in \yout_i$, and $\spsi\sphi\spsi(v) \in \xout_j$}.
\end{cases}\]
\item
If $v \in \coprod_{r=1}^l \win_r \amalg \dnzeta \subseteq \dmzeta$, then
\begin{footnotesize}
\[s(v) = \begin{cases}
\szeta(v) & \text{ if $\szeta(v) \not\in \xin_j$};\\
\spsi\szeta(v) & \text{ if $\szeta(v) \in \xin_j$ and $\spsi\szeta(v) \not\in \xout_j \amalg \yin_i$};\\
\szeta\spsi\szeta(v) & \text{ if $\szeta(v) \in \xin_j$ and $\spsi\szeta(v) \in \xout_j$};\\
\sphi\spsi\szeta(v) & \text{ if $\szeta(v) \in \xin_j$, $\spsi\szeta(v) \in \yin_i$, and $\sphi\spsi\szeta(v) \not\in \yout_i$};\\
\spsi\sphi\spsi\szeta(v) & \text{ if  $\szeta(v) \in \xin_j$, $\spsi\szeta(v) \in \yin_i$, $\sphi\spsi\szeta(v) \in \yout_i$, and $\spsi\sphi\spsi\szeta(v) \not\in \xout_j$};\\
\szeta\spsi\sphi\spsi\szeta(v) & \text{ if  $\szeta(v) \in \xin_j$, $\spsi\szeta(v) \in \yin_i$, $\sphi\spsi\szeta(v) \in \yout_i$, and $\spsi\sphi\spsi\szeta(v) \in \xout_j$}.
\end{cases}\]
\end{footnotesize}
\end{enumerate}
This finishes the proof of the desired equality \eqref{vertical-assoc}.
\end{proof}

\begin{theorem}
\label{wd-operad}
For any class $S$, when equipped with the structure in Def. \ref{wd-equivariance}--\ref{def:compi-wd}, $\WD$ in Def. \ref{wd-equivalence} is a $\boxs$-colored operad, called the operad of wiring diagrams. \index{WD is an operad@$\WD$ is an operad}
\end{theorem}

\begin{proof}
In view of Def. \ref{def:pseudo-operad}, this follows from Lemmas \ref{compi-unity}, \ref{compi-horizontal-associative}, and \ref{compi-vertical-associative}.
\end{proof}

\begin{remark}
When $S$ is the class of pointed sets, Theorem \ref{wd-operad} is proved in \cite{rupel-spivak} (section 2).  The proof in \cite{rupel-spivak} is somewhat different than ours because it uses Def. \ref{def:colored-operad} of a colored operad instead of Def. \ref{def:pseudo-operad}.
\end{remark}

\section{Summary of Chapter \ref{ch02-wiring-diagrams}}

\begin{enumerate}
\item An $S$-colored operad consists of a class $\sO\dconecn$ for each $d, c_1, \ldots, c_n \in S$ and $n \geq 0$ together with symmetric group actions, an operadic composition, and colored units that satisfy the associativity, unity, and equivariance axioms.
\item Every operad is determined by the $\compi$-compositions.
\item An $S$-wiring diagram has a finite number of input boxes, an output box, an $S$-finite set of delay nodes, and a supplier assignment that satisfies the non-instantaneity requirement.
\item For each class $S$, the collection of $S$-wiring diagrams $\WD$ is a $\boxs$-colored operad.
\end{enumerate}

\chapter{Generators and Relations}
\label{ch03-generating-wd}

Fix a class $S$.  The purpose of this chapter is to describe a finite number of wiring diagrams that we will later show to be sufficient to describe the entire operad $\WD$ of wiring diagrams (Theorems \ref{thm:wd-generator-relation}) as well as its variants $\wddot$ (Theorem \ref{thm:without-dn-coherence}) and $\wdzero$ (Theorem \ref{thm:strict-wd-coherence}).  One may also regard this chapter as consisting of a long list of examples of wiring diagrams.

In Section \ref{sec:generating-wd} we describe eight wiring diagrams, called the \emph{generating wiring diagrams}.  In Theorem \ref{stratified-presentation-exists} we will show that they generate the operad $\WD$ of wiring diagrams.  This means that every wiring diagram can be obtained from finitely many generating wiring diagrams via iterated operadic compositions.  For now one may think of the generating wiring diagrams as examples of wiring diagrams.

In Section \ref{sec:internal-wasted} we explain why a wiring diagram with an internal wasted wire is not among the generating wiring diagrams.  More concretely, we will observe in Prop. \ref{prop:internal-wasted-wire} that an internal wasted wire can be generated using two generating wiring diagrams.

In Section \ref{sec:elementary-relations} we describe $28$ \emph{elementary relations} among the generating wiring diagrams.  In Theorem \ref{thm:wd-generator-relation} we will show that these elementary relations together with the operad associativity and unity axioms--\eqref{compi-associativity}, \eqref{compi-associativity-two}, \eqref{compi-left-unity}, and \eqref{compi-right-unity}--for the  generating wiring diagrams generate \emph{all} the relations in the operad $\WD$ of wiring diagrams.  In other words, suppose an arbitrary wiring diagram can be built in two ways using the generating wiring diagrams.  Then there exists a finite sequence of steps connecting them in which each step is given by one of the $28$ elementary relations or an operad associativity/unity axiom for the generating wiring diagrams.   For now one may think of the elementary relations as examples of the operadic composition in the operad $\WD$.

\section{Generating Wiring Diagrams}
\label{sec:generating-wd}

Recall the definition of a wiring diagram (Def. \ref{wd-equivalence}).  In this section, we introduce $8$ wiring diagrams, called the generating wiring diagrams.  They will be used in later chapters to give a finite presentation for the operad $\WD$ of wiring diagrams.

\begin{definition}
\label{def:empty-wd}
Define the \emph{empty wiring diagram} \label{notation:empty-wd}$\epsilon \in \WD\emptynothing$ with:\index{empty wiring diagram} \index{epsilon@$\epsilon$}
\begin{enumerate}
\item
no input boxes;
\item
the empty box $\varnothing$ (Def. \ref{def:s-box}) as the output box;
\item
no delay nodes (i.e., $\dn = \varnothing$);
\item
supplier assignment $s : \dm = \varnothing \to \varnothing = \supply$ the trivial function.
\end{enumerate}
\end{definition}

The next wiring diagram has a delay node as depicted in the following picture, where we use the convention that delay nodes are drawn as circles as in \eqref{wd-first-example}.
\begin{center}
\begin{tikzpicture}[scale=0.7]
\draw [ultra thick] (0,1.2) rectangle (3,2.8);
\node at (1.5,3.2) {$\delta_d$};
\draw [ultra thick] (1.5,2) circle [radius=0.5];
\node at (1.5,2) {$d$};
\draw [arrow] (-.5,2) -- (1,2);
\draw [arrow] (2,2) -- (3.5,2);
\end{tikzpicture}
\end{center}

\begin{definition}
\label{def:one-dn}
Suppose $d \in S$.  Denote also by $d \in \boxs$ the box with one input and one output, both also denoted by $d$ and have values $d \in S$. Define the \emph{$1$-delay node} \label{notation:dn-wd}$\delta_d \in \WD\dnothing$ as the wiring diagram with:\index{deltad@$\delta_d$}\index{$1$-delay node}
\begin{enumerate}
\item
no input boxes;
\item
the output box $d = (\{d\}, \{d\})$, in which both $d$'s have values $d \in S$;
\item
$\dn = \{d\}$, in which $d$ has value $d \in S$;
\item
supplier assignment 
\[
\nicexy{
\dm = \yout \amalg \dn = \{d\} \amalg \{d\} \ar[r]^-{s}
& \{d\} \amalg \{d\} = \yin \amalg \dn = \supply
}\]
the identity function that takes $d \in \yout$ to $d \in \dn$ and $d \in \dn$ to $d \in \yin$.
\end{enumerate}
\end{definition}

Next we define the wiring diagram:
\begin{center}
\begin{tikzpicture}[scale=.7]
\draw [ultra thick] (2,1) rectangle (3,2);
\draw [ultra thick] (1,0.7) rectangle (4,2.5);
\node at (2.5,2.8) {$\tau_{X,Y}$};
\node at (2.5,1.5) {$X$};
\node at (3.6,2.2) {$Y$};
\draw [arrow] (0.5,1.7) -- (2,1.7);
\draw [arrow] (0.5,1.3) -- (2,1.3); 
\draw [arrow] (3,1.8) -- (4.5,1.8);
\draw [arrow] (3,1.5) -- (4.5,1.5); 
\draw [arrow] (3,1.2) -- (4.5,1.2); 
\end{tikzpicture}
\end{center}

\begin{definition}
\label{def:name-change}
Suppose $X, Y \in \boxs$ together with isomorphisms $\fin : \xin \cong \yin$ and $\fout : \yout \cong \xout$ in $\Fins$.  Define the wiring diagram $\tau_{f} \in \WD\yx$ with:\index{tauxy@$\tau_{X,Y}$}\index{name change}
\begin{enumerate}
\item
one input box $X$ and output box $Y$;
\item
no delay nodes;
\item
supplier assignment
\[\nicexy@C+1cm{
\dm = \yout \amalg \xin \ar[r]^-{s \,=\, \fout \amalg \fin}
& \yin \amalg \xout = \supply}\]
the coproduct of the given isomorphisms.
\end{enumerate}
We will often suppress the given isomorphisms and simply write \label{notation:name-change-wd}$\tau_{X,Y}$ or even $\tau$, which will be called a \emph{name change}.
\end{definition}

Next we define the wiring diagram:
\begin{center}
\begin{tikzpicture}[scale=.7]
\draw [ultra thick] (1,2) rectangle (2,3);
\draw [implies] (-.5,2.5) to (1,2.5);
\draw [implies] (2,2.5) to (3.5,2.5);
\node at (1.5,2.5) {$X$};
\draw [ultra thick] (1,.5) rectangle (2,1.5);
\draw [implies] (-.5,1) to (1,1);
\draw [implies] (2,1) to (3.5,1);
\node at (1.5,1) {$Y$};
\draw [ultra thick] (0,0.2) rectangle (3,3.3);
\node at (1.5,3.6) {$\theta_{X,Y}$};
\end{tikzpicture}
\end{center}

\begin{definition}
\label{def:theta-wd}
Suppose $X,Y \in \boxs$.  Define the wiring diagram \label{notation:2-cell}$\theta_{X,Y} \in \WD\xplusyxy$ with:\index{thetaxy@$\theta_{X,Y}$} \index{2-cell@$2$-cell}
\begin{enumerate}
\item
two input boxes $(X,Y)$ and output box $X \amalg Y$;
\item
no delay nodes;
\item
supplier assignment
\[
\nicexy{
\dm = \left[\xout \amalg \yout\right] \amalg \left[\xin \amalg \yin\right]
\ar[r]^-{s} &
\left[\xin \amalg \yin\right] \amalg \left[\xout \amalg \yout\right] = \supply
}\]
the identity map.
\end{enumerate}
We will call $\theta_{X,Y}$ a \emph{$2$-cell}.
\end{definition}

Next we define the wiring diagram:
\begin{center}
\begin{tikzpicture}[scale=.8]
\draw [ultra thick] (2,1) rectangle (3,2);
\node at (2.5,1.4) {$X$};
\node at (2.2,1.8) {\tiny{$x_-$}};
\node at (2.8,1.8) {\tiny{$x_+$}};
\draw [implies] (.5,1.3) to (2,1.3);
\node at (-.1,1.8) {$\xminusxin$};
\draw [implies] (3,1.3) to (4.5,1.3);
\node at (5.2,1.8) {$\xminusxout$};
\draw [ultra thick] (1,0.7) rectangle (4,2.5);
\node at (2.5,2.8) {$\lambda_{X, x}$};
\draw [arrow, looseness=3.5] (3,1.8) to [out=30, in=150] (2,1.8);
\end{tikzpicture}
\end{center}

\begin{definition}
\label{def:loop-wd}
Suppose:
\begin{itemize}
\item
$X \in \boxs$, and $(x_+, x_-) \in \xout \times \xin$ such that $v(x_+) = v(x_-) \in S$.
\item
$X \setminus x \in \boxs$ is obtained from $X$ by removing $x_{\pm}$, so $(X \setminus x)^{\inp} = \xin \setminus \{x_-\}$ and $(X \setminus x)^{\out} = \xout \setminus \{x_+\}$.  
\end{itemize}
Define the wiring diagram \label{notation:1-loop}$\lambda_{X,x} \in \WD\xminusxx$ with:\index{lambdaxx@$\lambda_{X,x}$} \index{$1$-loop}
\begin{enumerate}
\item
one input box $X$ and output box $X \setminus x$;
\item
no delay nodes;
\item
supplier assignment
\[
\nicexy@R-.5cm{
\dm = 
\overbrace{\left[\xout \setminus \{x_+\}\right]}^{\xminusxout} 
\amalg 
\overbrace{\left[\xin \setminus \{x_-\}\right] \amalg \{x_-\}}^{\xin}
\ar[d]_-{s}
\\
\supply = 
\underbrace{\left[\xin \setminus \{x_-\}\right]}_{\xminusxin} 
\amalg
\underbrace{\left[\xout \setminus \{x_+\}\right] \amalg \{x_+\}}_{\xout}
}\]
given by $s(x_-) = x_+$ and the identity function everywhere else.
\end{enumerate}
We will call the wiring diagram $\lambda_{X,x}$ a \emph{$1$-loop}.
\end{definition}

Next we define the wiring diagram:
\begin{center}
\begin{tikzpicture}[scale=.8]
\draw [ultra thick] (2,1) rectangle (3,2.1);
\node at (2.7,1.5) {$X$};
\node at (2.2,1.9) {\tiny{$x_1$}};
\node at (2.2,1.5) {\tiny{$x_2$}};
\draw [implies] (.5,1.2) to (2,1.2);
\node at (0,1.6) {$\yin$};
\node at (.7,1.9) {\tiny{$x_{12}$}};
\draw [thick] (.5,1.7) -- (1.2,1.7);
\draw [arrow] (1.2,1.7) to [out=0, in=180] (2,1.9);
\draw [arrow] (1.2,1.7) to [out=0, in=180] (2,1.5);
\draw [implies] (3,1.5) to (4.5,1.5);
\node at (5,1.6) {$\yout$};
\draw [ultra thick] (1,0.7) rectangle (4,2.4);
\node at (2.5,2.7) {$\sigma_{X, x_1, x_2}$};
\end{tikzpicture}
\end{center}

\begin{definition}
\label{def:insplit-wd}
Suppose:
\begin{itemize}
\item
$X \in \boxs$, and $x_1,x_2 \in \xin$ are distinct elements such that $v(x_1) = v(x_2) \in S$.
\item
$Y = X/(x_1=x_2) \in \boxs$ is obtained from $X$ by identifying $x_1$ and $x_2$, so $\yin = \xin/(x_1=x_2)$ and $\yout = \xout$.  The identified element of $x_1$ and $x_2$ in $\yin$ will be denoted by $x_{12}$.
\end{itemize}
Define the wiring diagram \label{notation:insplit}$\sigma_{X,x_1,x_2} \in \WD\yx$ with:\index{in-split}\index{sigmasubx@$\sigma_{X,x_1,x_2}$}
\begin{enumerate}
\item
one input box $X$ and output box $Y$;
\item
no delay nodes;
\item
supplier assignment
\[
\nicexy{
\dm = \yout \amalg \xin = \xout \amalg \xin
\ar[d]_-{s}
& \ni & x_1,\, x_2 \ar@{|->}[d]
\\
\supply = \yin \amalg \xout = \frac{\xin}{(x_1\,=\,x_2)} \amalg \xout
& \ni & x_{12} 
}\]
that sends both $x_1, x_2 \in \xin$ to $x_{12} \in \yin$ and is the identity function everywhere else.
\end{enumerate}  
We will call the wiring diagram $\sigma_{X,x_1,x_2}$ an \emph{in-split}.
\end{definition}

Next we define the wiring diagram:
\begin{center}
\begin{tikzpicture}[scale=.8]
\draw [ultra thick] (2,1) rectangle (3,2);
\node at (2.5,1.3) {$X$};
\node at (2.75,1.8) {\tiny{$y_{12}$}};
\draw [thick] (3,1.8) -- (3.3,1.8);
\draw [arrow] (3.3,1.8) to [out=0, in=180] (4.5,2);
\draw [arrow] (3.3,1.8) to [out=0, in=180] (4.5,1.6);
\draw [implies] (.5,1.5) to (2,1.5);
\node at (0,1.5) {$\yin$};
\draw [implies] (3,1.2) to (4.5,1.2);
\node at (5.5,1.5) {$\yout$};
\node at (4.7,2) {\tiny{$y_{1}$}};
\node at (4.7,1.6) {\tiny{$y_{2}$}};
\draw [ultra thick] (1,0.7) rectangle (4,2.3);
\node at (2.5,2.6) {$\sigma^{Y, y_1, y_2}$};
\end{tikzpicture}
\end{center}

\begin{definition}
\label{def:out-split}
Suppose:
\begin{itemize}
\item
$Y \in \boxs$, and $y_1, y_2 \in \yout$ are distinct elements such that $v(y_1) = v(y_2) \in S$.
\item
$X = Y/(y_1=y_2) \in \boxs$ is obtained from $Y$ by identifying $y_1$ and $y_2$, so $\xin = \yin$ and $\xout = \yout/(y_1 = y_2)$.  The identified element of $y_1$ and $y_2$ in $\xout$ will be denoted by $y_{12}$.
\end{itemize}
Define the wiring diagram \label{notation:outsplit}$\sigma^{Y,y_1,y_2} \in \WD\yx$ with:\index{out-split}\index{sigmasupy@$\sigma^{Y,y_1,y_2}$}
\begin{enumerate}
\item
one input box $X$ and output box $Y$;
\item
no delay nodes;
\item
supplier assignment
\[
\nicexy{
\dm = \yout \amalg \xin = \yout \amalg \yin
\ar[d]_-{s}
& \ni & y_1,\,y_2 \ar@{|->}[d]
\\
\supply = \yin \amalg \xout = \yin \amalg \frac{\yout}{(y_1\, =\,y_2)}
& \ni & y_{12} 
}\]
that sends both $y_1, y_2 \in \yout$ to $y_{12} \in \xout$ and is the identity function everywhere else.
\end{enumerate}  
We will call the wiring diagram $\sigma^{Y,y_1,y_2}$ an \emph{out-split}.
\end{definition}

Next we define the following wiring diagram with an external wasted wire:
\begin{center}
\begin{tikzpicture}[scale=.8]
\draw [ultra thick] (2,1) rectangle (3,2);
\node at (2.5,1.5) {$X$};
\draw [implies] (.5,1.3) to (2,1.3);
\node at (0,1.5) {$\yin$};
\draw [arrow] (.5,1.7) -- (1,1.7);
\node at (.7,1.9) {\tiny{$y$}};
\draw [implies] (3,1.5) to (4.5,1.5);
\node at (5.,1.5) {$\yout$};
\draw [ultra thick] (1,0.7) rectangle (4,2.3);
\node at (2.5,2.6) {$\omega_{Y,y}$};
\end{tikzpicture}
\end{center}

\begin{definition}
\label{def:wasted-wire-wd}
Suppose:
\begin{itemize}
\item
$Y \in \boxs$, and $y \in \yin$.
\item
$X \in \boxs$ is obtained from $Y$ by removing $y$, so $\xin = \yin \setminus \{y\}$ and $\xout = \yout$.
\end{itemize}
Define the wiring diagram \label{notation:wasted-wire}$\omega_{Y,y} \in \WD\yx$ with:\index{1-wasted wire@$1$-wasted wire}\index{omegay@$\omega_{Y,y}$}
\begin{enumerate}
\item
one input box $X$ and output box $Y$;
\item
no delay nodes;
\item
supplier assignment 
\[
\nicexy{
\dm = \yout \amalg \xin = \yout \amalg \left[\yin \setminus \{y\}\right]
\ar[d]_-{s}\\
\supply = \yin \amalg \xout = \yin \amalg \yout
}\]
the inclusion.
\end{enumerate}
We will call the wiring diagram $\omega_{Y,y}$ a \emph{$1$-wasted wire}.
\end{definition}

\begin{definition}
\label{def:generating-wiring-diagrams}
The eight wiring diagrams in Def. \ref{def:empty-wd}--\ref{def:wasted-wire-wd} will be referred to as \emph{generating wiring diagrams}.\index{generating wiring diagrams}
\end{definition}

\begin{remark}
\label{rk:generators-arity}
Among the generating wiring diagrams:
\begin{enumerate}
\item
A $1$-delay node $\delta_d$ (Def. \ref{def:one-dn}) is the only wiring diagram that has a delay node.
\item
A $1$-wasted wire $\omega_{Y,y}$ (Def. \ref{def:wasted-wire-wd}) is the only wiring diagram that has an external wasted wire, namely $y \in \yin$.    
\item
None has an internal wasted wire (Def. \ref{def:wiring-diagram}).  As we will see in Prop. \ref{prop:internal-wasted-wire} below, an internal wasted wire can be generated using a $1$-loop and a $1$-wasted wire, hence is not needed as a generator.
\item
The empty wiring diagram $\epsilon$ (Def. \ref{def:empty-wd}) and a $1$-delay node $\delta_d$ are $0$-ary elements in $\WD$.
\item
A name change $\tau$ (Def. \ref{def:name-change}), a $1$-loop $\lambda_{X,x}$ (Def. \ref{def:loop-wd}), an in-split $\sigma_{X,x_1,x_2}$ (Def. \ref{def:insplit-wd}), an out-split $\sigma^{Y,y_1,y_2}$ (Def. \ref{def:out-split}), and a $1$-wasted wire $\omega_{Y,y}$ are unary elements in $\WD$.
\item
A $2$-cell $\theta_{X,Y}$ (Def. \ref{def:theta-wd}) is a binary element in $\WD$.
\end{enumerate}
\end{remark}

\section{Internal Wasted Wires}
\label{sec:internal-wasted}

Recall from Def. \ref{def:wiring-diagram} that an \emph{internal wasted wire} is an internal output, hence a supply wire, that does not belong to the image of the supplier assignment.  The purpose of this section is to explain why the wiring diagram 
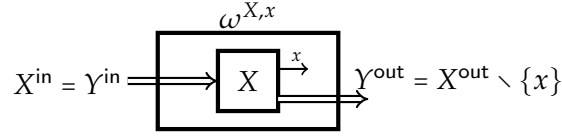
\begin{center}
\begin{tikzpicture}[scale=.8]
\draw [ultra thick] (2,1) rectangle (3,2);
\node at (2.5,1.5) {$X$};
\draw [implies] (.5,1.5) to (2,1.5);
\node at (-.5,1.5) {$\xin = \yin$};
\draw [arrow] (3,1.7) -- (3.5,1.7);
\node at (3.3,1.9) {\tiny{$x$}};
\draw [implies] (3,1.2) to (4.5,1.2);
\node at (6,1.5) {$\yout = \xout \setminus \{x\}$};
\draw [ultra thick] (1,0.7) rectangle (4,2.3);
\node at (2.5,2.6) {$\omega^{X,x}$};
\end{tikzpicture}
\end{center}
that has an internal wasted wire $x$ is not needed as a generating wiring diagram.  First we define this wiring diagram.

\begin{definition}
\label{def:internal-wasted-wire}
Suppose:
\begin{itemize}
\item
$X \in \boxs$, and $x \in \xout$.
\item
$Y = X \setminus \{x\} \in \boxs$ is obtained from $X$ by removing $x$.
\end{itemize}
Define the wiring diagram \label{notation:external-wasted-wire}$\omega^{X,x} \in \WD\yx$ with:\index{$1$-internal wasted wire}\index{omegax@$\omega^{X,x}$}
\begin{enumerate}
\item
one input box $X$ and output box $Y$;
\item
no delay nodes;
\item
supplier assignment
\[
\nicexy{
\dm = \yout \amalg \xin = [\xout \setminus \{x\}] \amalg \xin
\ar[d]_-{s}
\\
\supply = \yin \amalg \xout = \xin \amalg \xout
}\]
the inclusion.
\end{enumerate}
We will call the wiring diagram $\omega^{X,x}$ a \emph{$1$-internal wasted wire}.
\end{definition}

\begin{motivation}
The following observation says that a $1$-internal wasted wire can be obtained as the substitution of a $1$-wasted wire into a $1$-loop, both of which are generating wiring diagrams.  This is expressed in the following picture
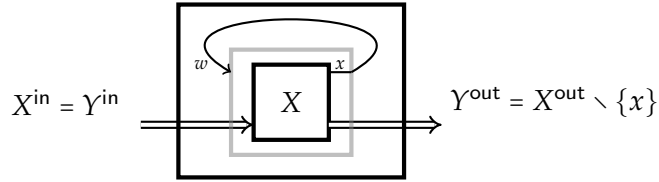
\begin{center}
\begin{tikzpicture}
\draw [ultra thick] (2,1) rectangle (3,2);
\node at (2.5,1.5) {$X$};
\draw [implies] (.5,1.2) to (2,1.2);
\node at (-.5,1.5) {$\xin = \yin$};
\draw [implies] (3,1.2) to (4.5,1.2);
\node at (6,1.5) {$\yout = \xout \setminus \{x\}$};
\draw [ultra thick, gray, semitransparent] (1.7,.8) rectangle (3.3,2.2);
\draw [thick] (3,1.9) -- (3.3,1.9);
\node at (3.15,2) {\tiny{$x$}};
\draw [arrow, looseness=3] (3.3,1.9) to [out=30, in=150] (1.7,1.9);
\node at (1.3,2) {\tiny{$w$}};
\draw [ultra thick] (1,0.5) rectangle (4,2.8);
\end{tikzpicture}
\end{center}
in which the intermediate gray box will be called $W$ below.
\end{motivation}

\begin{proposition}
\label{prop:internal-wasted-wire}
Suppose:
\begin{itemize}
\item
$\omega^{X,x} \in \WD\yx$ is a $1$-internal wasted wire (Def. \ref{def:internal-wasted-wire}).
\item
$W = X \amalg \{w\} \in \boxs$ such that $w \in \win$ satisfies $v(w) = v(x) \in S$.
\item
$\lambda_{W,\{w,x\}} \in \WD\yw$ is the $1$-loop (Def. \ref{def:loop-wd}) in which $x \in \xout=\wout$ is the supply wire of $w \in \win$.
\item
$\omega_{W,w} \in \WD\wx$ is a $1$-wasted wire (Def. \ref{def:wasted-wire-wd}).
\end{itemize}
Then
\begin{equation}
\label{internal-wasted-wire}
\left(\lambda_{W,\{w,x\}}\right) \compone
\left(\omega_{W,w}\right)  
=
\omega^{X,x} \in \WD\yx.
\end{equation}
\end{proposition}

\begin{proof}
By definition both sides of \eqref{internal-wasted-wire} belong to $\WD\yx$ and have no delay nodes.  It remains to check that their supplier assignments are equal.  By the definitions of $\compone$ (Def. \ref{def:compi-wd}), $1$-loop, and $1$-wasted wire, the supplier assignment of the left side $\left(\lambda_{W,\{w,x\}}\right) \compone
\left(\omega_{X,x}\right)$, namely
\[
\nicexy{
\dm = \yout \amalg \xin = [\xout \setminus \{x\}] \amalg \xin
\ar[d]_-{s}
\\
\supply = \yin \amalg \xout = \xin \amalg \xout
}\]
is the inclusion.  By Def. \ref{def:internal-wasted-wire} this is also the supplier assignment of the $1$-internal wasted wire $\omega^{X,x}$.
\end{proof}

As a consequence of \eqref{internal-wasted-wire}, the $1$-internal wasted wire $\omega^{X,x}$ is \emph{not} needed as a generating wiring diagram.

\section{Elementary Relations}
\label{sec:elementary-relations}

The purpose of this section is to introduce $28$ elementary relations among the generating wiring diagrams (Def. \ref{def:generating-wiring-diagrams}).  Each elementary relation is proved by a simple inspection of the relevant definitions of the generating wiring diagrams and operadic compositions, similar to the proofs of Lemma \ref{compi-horizontal-associative}, Lemma \ref{compi-vertical-associative}, and Prop. \ref{prop:internal-wasted-wire} above.  Therefore, we will prove only the first one and omit the proofs for the rest, providing a picture instead in most cases.  We will frequently use the $\compi$-composition \eqref{operadic-compi} in describing these elementary relations.

\begin{notation}
\label{comp-is-compone}
Suppose $\O$ is an $S$-colored operad (Def. \ref{def:pseudo-operad}), and $T$ is a set. 
\begin{enumerate}
\item
If $\varphi \in \O\dcsingle$ where the input profile $(c)$ has length $1$ and if $\phi \in \O\cub$, then we write
\begin{equation}
\label{comp-means-compone}
\varphi \comp \phi = \varphi \compone \phi.
\end{equation}
\item
Suppose $\varphi_1, \ldots, \varphi_k \in \sO$ such that each of $\varphi_1, \ldots, \varphi_{k-1}$ belongs to an entry of $\O$ whose input profile has length $1$.  Then we write
\begin{equation}
\label{iterated-compone}
\varphi_1 \comp \cdots \comp \varphi_k
=
\Bigl( \cdots(\varphi_1 \compone \varphi_2) \compone \cdots \Bigr) \compone \varphi_k
\end{equation}
whenever the right side is defined, in which each pair of parentheses starts on the left.  For example, we have
\[
\begin{split}
\varphi_1 \comp \varphi_2 \comp \varphi_3 
&= (\varphi_1 \compone \varphi_2) \compone \varphi_3,\\
\varphi_1 \comp \varphi_2 \comp \varphi_3 \comp \varphi_4 
&= \Bigl((\varphi_1 \compone \varphi_2) \compone \varphi_3\Bigr) \compone \varphi_4.
\end{split}\]
\item
Write \label{notation:cardinality}$|T|$ for the cardinality of $T$.
\end{enumerate}
\end{notation}

The first six relations are about the name change wiring diagrams (Def. \ref{def:name-change}).  The first relation says that two consecutive name changes can be composed into one name change.

\begin{proposition}
\label{prop:move:a1}
Suppose:
\begin{itemize}
\item
$\tau_{Y,Z} \in \WD\zy$ and $\tau_{X,Y} \in \WD\yx$ are name changes.
\item
$\tau_{X,Z} \in \WD\zx$ is the name change given by composing the isomorphisms that define $\tau_{Y,Z}$ and $\tau_{X,Y}$.
\end{itemize}
Then 
\begin{equation}
\label{move:a1}
\left(\tau_{Y,Z}\right) \comp \left(\tau_{X,Y}\right)
=
\tau_{X,Z} \in \WD\zx.
\end{equation}
\end{proposition}

\begin{proof}
We are given isomorphisms $\fin : \xin \cong \yin$ and $\fout : \yout \cong \xout$ for $\tau_{X,Y}$ and isomorphisms $\gin : \yin \cong \zin$ and $\gout : \zout \cong \yout$ for $\tau_{Y,Z}$.  The name change $\tau_{X,Z}$ given by composing these isomorphisms is defined by the isomorphisms
\[\gin \circ \fin : \xin \cong \zin \andspace \fout \circ \gout : \zout \cong \xout.\]

On the other hand, by Def. \ref{def:compi-wd} the composition $\tau_{Y,Z} \circ \tau_{X,Y}$ has no delay nodes, since neither $\tau_{X,Y}$ nor $\tau_{Y,Z}$ has a delay node.  Its supplier assignment is the function
\[\nicexy{\dm = \zout \coprod \xin \ar[r]^-{s} & \zin \coprod \xout = \supply}\]
given by
\[\begin{split}
s(z) &= s_{\tau_{X,Y}} s_{\tau_{Y,Z}}(z) = \fout \gout(z) \forspace z \in \zout;\\
s(x) &= s_{\tau_{Y,Z}} s_{\tau_{X,Y}}(x) = \gin \fin(x) \forspace x \in \xin.\end{split}\]
This is the same supplier assignment as that of the name change $\tau_{X,Z}$ above.
\end{proof}

The next relation says that name changes inside a $2$-cell (Def. \ref{def:theta-wd}) can be rewritten as a name change outside of a $2$-cell.

\begin{proposition}
\label{prop:move:a2}
Suppose:
\begin{itemize}
\item
$\tau_{X,X'} \in \WD\xprimex$ and $\tau_{Y,Y'}  \in \WD\yprimey$ are name changes.
\item
$\tau_{X \amalg Y, X' \amalg Y'} \in \WD\xprimeplusyprimexplusy$ is the name change induced by $\tau_{X,X'}$ and $\tau_{Y,Y'}$.
\item
$\theta_{X',Y'} \in \WD\xplusyxyprime$ and $\theta_{X,Y} \in \WD\xplusyxy$ are $2$-cells.
\end{itemize}
Then
\begin{equation}
\label{move:a2}
\Bigl(\theta_{X',Y'} \compone \tau_{X,X'}\Bigr) \comptwo \tau_{Y,Y'}
=
\left(\tau_{X \amalg Y, X' \amalg Y'}\right) 
\comp \left(\theta_{X,Y}\right)
\in \WD\xprimeplusyprimexy.
\end{equation}
\end{proposition}

The next relation says that a name change inside a $1$-loop (Def. \ref{def:loop-wd}) can be rewritten as a name change of a $1$-loop.

\begin{proposition}
\label{prop:move:a3}
Suppose:
\begin{itemize}
\item
$X \in \boxs$, $(x_+, x_-) \in \xout \times \xin$ such that $v(x_+) = v(x_-) \in S$.
\item
$X \setminus x \in \boxs$ is obtained from $X$ by removing $x_{\pm}$.
\item
$\lambda_{X,x} \in \WD\xminusxx$ is the corresponding $1$-loop.
\item
$\tau_{X,Y} \in \WD\yx$ is a name change such that $(y_+,y_-) \in \yout\times\yin$ corresponds to $(x_+, x_-)$.
\item
$\lambda_{Y,y} \in \WD\yminusyy$ is the corresponding $1$-loop.
\item
$\tau_{X \setminus x,Y \setminus y} \in \WD\yminusyxminusx$ is the name change induced by $\tau_{X,Y}$.
\end{itemize}
Then
\begin{equation}
\label{move:a3}
\left(\lambda_{Y,y}\right) 
\comp \left(\tau_{X,Y}\right)
=
\left(\tau_{X \setminus x, Y \setminus y}\right) 
\comp \left(\lambda_{X,x}\right) 
\in \WD\yminusyx.
\end{equation}
\end{proposition}

The next relation says that a name change inside an in-split (Def. \ref{def:insplit-wd}) can be rewritten as a name change of an in-split.

\begin{proposition}
\label{prop:move:a4}
Suppose:
\begin{itemize}
\item
$X \in \boxs$, and $x_1,x_2 \in \xin$ are distinct elements such that $v(x_1) = v(x_2) \in S$.
\item
$X' = X/(x_1=x_2) \in \boxs$ is obtained from $X$ by identifying $x_1$ and $x_2$.
\item
$\sigma_{X,x_1,x_2} \in \WD\xprimex$ is the corresponding in-split.
\item
$\tau_{X,Y} \in \WD\yx$ is a name change with $y_1,y_2 \in \yin$ corresponding to $x_1,x_2 \in \xin$.
\item
$Y' = Y/(y_1=y_2) \in \boxs$ is obtained from $Y$ by identifying $y_1$ and $y_2$.
\item
$\sigma_{Y,y_1,y_2} \in \WD\yprimey$ is the corresponding in-split.
\item
$\tau_{X',Y'} \in \WD\yprimexprime$ is the name change induced by $\tau_{X,Y}$.
\end{itemize}
Then
\begin{equation}
\label{move:a4}
\left(\sigma_{Y,y_1,y_2}\right) \comp \left(\tau_{X,Y}\right)
=
\left(\tau_{X',Y'}\right) 
\comp \left(\sigma_{X,x_1,x_2}\right) 
\in\WD\yprimex.
\end{equation}
\end{proposition}

The next relation is the out-split (Def. \ref{def:out-split}) analogue of \eqref{move:a4}.

\begin{proposition}
\label{prop:move:a5}
Suppose:
\begin{itemize}
\item
$Y \in \boxs$, and $y_1, y_2 \in \yout$ are distinct elements such that $v(y_1) = v(y_2) \in S$.
\item
$Y' = Y/(y_1=y_2) \in \boxs$ is obtained from $Y$ by identifying $y_1$ and $y_2$.
\item
$\sigma^{Y,y_1,y_2} \in \WD\yyprime$ is an out-split.
\item
$\tau_{X,Y} \in \WD\yx$ is a name change with $x_1,x_2 \in \xout$ corresponding to $y_1,y_2 \in \yout$.
\item
$X' = X/(x_1=x_2) \in \boxs$ is obtained from $X$ by identifying $x_1$ and $x_2$.
\item
$\sigma^{X,x_1,x_2} \in \WD\xxprime$ is the corresponding out-split.
\item
$\tau_{X',Y'} \in \WD\yprimexprime$ is the name change induced by $\tau_{X,Y}$.
\end{itemize}
Then
\begin{equation}
\label{move:a5}
\left(\sigma^{Y,y_1,y_2}\right) 
\comp \left(\tau_{X',Y'}\right)
=
\left(\tau_{X,Y}\right) 
\comp \left(\sigma^{X,x_1,x_2}\right) 
\in\WD\yxprime.
\end{equation}
\end{proposition}

The next relation says that a name change inside a $1$-wasted wire 
(Def. \ref{def:wasted-wire-wd}) can be rewritten as a name change of a $1$-wasted wire.

\begin{proposition}
\label{prop:move:a6}
Suppose:
\begin{itemize}
\item
$Y \in \boxs$, $y \in \yin$, and $Y' = Y \setminus \{y\} \in \boxs$ is obtained from $Y$ by removing $y$.
\item
$\omega_{Y,y} \in \WD\yyprime$ is the corresponding $1$-wasted wire.
\item
$\tau_{X,Y} \in \WD\yx$ is a name change with $x \in \xin$ corresponding to $y \in \yin$.
\item
$X' \in \boxs$ is obtained from $X$ by removing $x$.
\item
$\omega_{X,x} \in \WD\xxprime$ is the corresponding $1$-wasted wire.
\item
$\tau_{X',Y'} \in \WD\yprimexprime$ is the name change induced by $\tau_{X,Y}$.
\end{itemize}
Then
\begin{equation}
\label{move:a6}
\left(\omega_{Y,y}\right) 
\comp \left(\tau_{X',Y'}\right)
=
\left(\tau_{X,Y}\right) 
\comp \left(\omega_{X,x}\right)
\in \WD\yxprime.
\end{equation}
\end{proposition}

The next seven relations are about $2$-cells (Def. \ref{def:theta-wd}).  The following relation says that substituting the empty wiring diagram (Def. \ref{def:empty-wd}) into a $2$-cell yields a colored unit  \eqref{wd-unit}.\index{unity of $2$-cells}

\begin{proposition}
\label{prop:move:b0}
Suppose:
\begin{itemize}
\item
$X \in \boxs$ with $X$-colored unit $\tensorunit_{X} \in \WD\xx$.
\item
$\epsilon \in \WD\emptynothing$ is the empty wiring diagram.
\item
$\theta_{X,\varnothing} \in \WD\xxempty$ is the $2$-cell with input boxes $(X,\varnothing)$ and output box $X$.
\end{itemize}
Then
\begin{equation}
\label{move:b0}
\theta_{X,\varnothing} \comptwo \epsilon 
= 
\tensorunit_{X}
\in \WD\xx.
\end{equation}
\end{proposition}

The next  relation is the associativity property of $2$-cells.\index{associativity of $2$-cells}  It says that, in the picture below, the wiring diagram in the middle can be constructed using two $2$-cells, either as the operadic composition on the left or the one on the right.
\begin{center}
\begin{tikzpicture}[scale=.7]
\draw [ultra thick] (1,2) rectangle (2,3);
\draw [implies] (-.5,2.5) to (1,2.5);
\draw [implies] (2,2.5) to (3.5,2.5);
\node at (1.5,2.5) {$X$};
\draw [ultra thick] (1,.5) rectangle (2,1.5);
\draw [implies] (-.5,1) to (1,1);
\draw [implies] (2,1) to (3.5,1);
\node at (1.5,1) {$Y$};
\draw [ultra thick] (1,-1) rectangle (2,0);
\draw [implies] (-.5,-.5) to (1,-.5);
\draw [implies] (2,-.5) to (3.5,-.5);
\node at (1.5,-.5) {$Z$};
\draw [ultra thick] (0,-1.3) rectangle (3,3.3);
\node at (-2,1) {$=$};
\draw [ultra thick] (-6,2) rectangle (-5,3);
\draw [implies] (-7.5,2.5) to (-6,2.5);
\draw [implies] (-5,2.5) to (-3.5,2.5);
\node at (-5.5,2.5) {$X$};
\draw [ultra thick] (-6,.5) rectangle (-5,1.5);
\draw [implies] (-7.5,1) to (-6,1);
\draw [implies] (-5,1) to (-3.5,1);
\node at (-5.5,1) {$Y$};
\draw [ultra thick] (-6,-1) rectangle (-5,0);
\draw [implies] (-7.5,-.5) to (-6,-.5);
\draw [implies] (-5,-.5) to (-3.5,-.5);
\node at (-5.5,-.5) {$Z$};
\draw [ultra thick] (-7,-1.3) rectangle (-4,3.3);
\draw [ultra thick, gray, semitransparent] (-6.5,.25) rectangle (-4.5,3.15);
\node at (5,1) {$=$};
\draw [ultra thick] (8,2) rectangle (9,3);
\draw [implies] (6.5,2.5) to (8,2.5);
\draw [implies] (9,2.5) to (10.5,2.5);
\node at (8.5,2.5) {$X$};
\draw [ultra thick] (8,.5) rectangle (9,1.5);
\draw [implies] (6.5,1) to (8,1);
\draw [implies] (9,1) to (10.5,1);
\node at (8.5,1) {$Y$};
\draw [ultra thick] (8,-1) rectangle (9,0);
\draw [implies] (6.5,-.5) to (8,-.5);
\draw [implies] (9,-.5) to (10.5,-.5);
\node at (8.5,-.5) {$Z$};
\draw [ultra thick] (7,-1.3) rectangle (10,3.3);
\draw [ultra thick, gray, semitransparent] (7.5,-1.15) rectangle (9.5,1.75);
\end{tikzpicture}
\end{center}

\begin{proposition}
\label{prop:move:b1}
Suppose:
\begin{itemize}
\item
$\theta_{X \amalg Y, Z} \in \WD\xplusyz$ and $\theta_{X,Y} \in \WD\xplusyxy$ are $2$-cells.
\item
$\theta_{X, Y \amalg Z} \in \WD\xyplusz$ and $\theta_{Y,Z} \in \WD\ypluszyz$ are $2$-cells.  
\end{itemize}
Then
\begin{equation}
\label{move:b1}
\left(\theta_{X \amalg Y, Z}\right) 
\compone \left(\theta_{X,Y}\right)
=
\left(\theta_{X, Y \amalg Z}\right) 
\comptwo \left(\theta_{Y,Z}\right)
\in \WD\xplusyplusz.
\end{equation}
\end{proposition}

The next relation is the commutativity property of $2$-cells\index{commutativity of $2$-cells} and uses the equivariant structure in $\WD$ \eqref{wd-permutation}.

\begin{proposition}
\label{prop:move:b2}
Suppose:
\begin{itemize}
\item
$\theta_{X,Y} \in \WD\xplusyxy$ is a $2$-cell.
\item
$(1~2) \in \Sigma_2$ is the non-trivial permutation.
\end{itemize}
Then
\begin{equation}
\label{move:b2}
\theta_{X,Y}(1~2) = \theta_{Y,X} \in \WD\yplusxyx.
\end{equation}
\end{proposition}

The next relation says that substituting a $1$-loop inside a $2$-cell can be rewritten as substituting a $2$-cell inside a $1$-loop.  It gives two different ways to construct the wiring diagram in the middle in the picture below using a $1$-loop and a $2$-cell, either as the operadic composition on the left or the one on the right.
\begin{center}
\begin{tikzpicture}[scale=.8]
\draw [ultra thick] (1,2) rectangle (2,3);
\draw [implies] (-.5,2.5) to (1,2.5);
\draw [implies] (2,2.5) to (3.5,2.5);
\node at (1.25,2.8) {\tiny{$x_-$}};
\node at (1.8,2.8) {\tiny{$x_+$}};
\node at (1.5,2.3) {$X$};
\node at (-1,2.8) {\tiny{$\xin \setminus \{x_-\}$}};
\node at (4,2.8) {\tiny{$\xout \setminus \{x_+\}$}};
\draw [ultra thick] (1,.5) rectangle (2,1.5);
\draw [implies] (-.5,1) to (1,1);
\draw [implies] (2,1) to (3.5,1);
\node at (1.5,1) {$Y$};
\draw [ultra thick] (0,0.2) rectangle (3,3.7);
\draw [arrow, looseness=3.5] (2,2.8) to [out=30, in=150] (1,2.8);
\node at (-1.5,1.75) {$=$};
\draw [ultra thick] (-5,2) rectangle (-4,3);
\draw [implies] (-6.5,2.5) to (-5,2.5);
\draw [implies] (-4,2.5) to (-2.5,2.5);
\node at (-4.5,2.5) {$X$};
\draw [ultra thick] (-5,.5) rectangle (-4,1.5);
\draw [implies] (-6.5,1) to (-5,1);
\draw [implies] (-4,1) to (-2.5,1);
\node at (-4.5,1) {$Y$};
\draw [ultra thick] (-6,0.2) rectangle (-3,3.7);
\draw [arrow, looseness=3.5] (-4,2.8) to [out=30, in=150] (-5,2.8);
\draw [ultra thick, gray, semitransparent] (-5.5,1.75) rectangle (-3.5,3.5);
\draw [ultra thick] (7,2) rectangle (8,3);
\draw [implies] (5.5,2.5) to (7,2.5);
\draw [implies] (8,2.5) to (9.5,2.5);
\node at (7.5,2.5) {$X$};
\draw [ultra thick] (7,.5) rectangle (8,1.5);
\draw [implies] (5.5,1) to (7,1);
\draw [implies] (8,1) to (9.5,1);
\node at (7.5,1) {$Y$};
\draw [ultra thick] (6,0.2) rectangle (9,3.7);
\draw [arrow, looseness=4.5] (8,2.8) to [out=30, in=150] (7,2.8);
\draw [ultra thick, gray, semitransparent] (6.8,.35) rectangle (8.2,3.15);
\node at (4.5,1.75) {$=$};
\end{tikzpicture}
\end{center}

\begin{proposition}
\label{prop:move:b3}
Suppose:
\begin{itemize}
\item
$X \in \boxs$, and $(x_+,x_-) \in \xout \times \xin$ such that $v(x_+) = v(x_-) \in S$.
\item
$X \setminus x \in \boxs$ is obtained from $X$ by removing $x_{\pm}$.
\item
$\theta_{X \setminus x, Y} \in \WD\xplusyminusxxminusxy$ and $\theta_{X,Y} \in \WD\xplusyxy$ are $2$-cells.
\item
$\lambda_{X,x} \in \WD\xminusxx$ and $\lambda_{X\amalg Y, x} \in \WD\xplusyminusxxplusy$ are the corresponding $1$-loops.
\end{itemize}
Then
\begin{equation}
\label{move:b3}
\left(\theta_{X \setminus x, Y}\right) 
\compone \left(\lambda_{X,x}\right)
= 
\left(\lambda_{X \amalg Y,x}\right) 
\comp \left(\theta_{X,Y}\right)
\in \WD\xplusyminusxxy.
\end{equation}
\end{proposition}

All the relations in the rest of this section can be illustrated with pictures similar to the two previous pictures, each one showing how a wiring diagram can be built in two different ways using operadic compositions.  So we will mostly just draw the picture of the wiring diagram being built without the accompanying pictures of the operadic compositions.

The next relation says that substituting an in-split inside a $2$-cell can be rewritten as substituting a $2$-cell inside an in-split.  It gives two different ways to construct the following wiring diagram using an in-split and a $2$-cell:
\begin{center}
\begin{tikzpicture}[scale=.8]
\draw [ultra thick] (1,2) rectangle (2,3);
\draw [implies] (-.5,2.1) to (1,2.1);
\draw [implies] (2,2.5) to (3.5,2.5);
\node at (1.2,2.8) {\tiny{$x_1$}};
\node at (1.2,2.5) {\tiny{$x_2$}};
\node at (1.7,2.5) {$X$};
\node at (-1.3,2.5) {$\frac{\xin}{(x_1\,=\,x_2)}$};
\draw [thick] (-.5,2.7) -- (.2,2.7);
\draw [arrow] (.2,2.7) to [out=0, in=180] (1,2.9);
\draw [arrow] (.2,2.7) to [out=0, in=180] (1,2.5);
\draw [ultra thick] (1,.5) rectangle (2,1.5);
\draw [implies] (-.5,1) to (1,1);
\draw [implies] (2,1) to (3.5,1);
\node at (1.5,1) {$Y$};
\draw [ultra thick] (0,0.2) rectangle (3,3.3);
\end{tikzpicture}
\end{center}

\begin{proposition}
\label{prop:move:b4}
Suppose:
\begin{itemize}
\item
$X \in \boxs$, and $x_1,x_2 \in \xin$ are distinct elements such that $v(x_1) = v(x_2) \in S$.
\item
$X' = X/(x_1=x_2) \in \boxs$ is obtained from $X$ by identifying $x_1$ and $x_2$.
\item
$\theta_{X',Y} \in \WD\xprimeplusyxprimey$ and $\theta_{X,Y} \in \WD\xplusyxy$ are $2$-cells.
\item
$\sigma_{X,x_1,x_2} \in \WD\xprimex$ and $\sigma_{X\amalg Y, x_1, x_2} \in \WD\xprimeplusyxplusy$ are in-splits.
\end{itemize}
Then
\begin{equation}
\label{move:b4}
\left(\theta_{X',Y}\right) 
\compone \left(\sigma_{X,x_1,x_2}\right)
= 
\left(\sigma_{X\amalg Y, x_1,x_2}\right) 
\comp \left(\theta_{X,Y}\right)
\in 
\WD\xprimeplusyxy.
\end{equation}
\end{proposition}

The next relation says that substituting an out-split inside a $2$-cell can be rewritten as substituting a $2$-cell inside an out-split.  It gives two different ways to construct the following wiring diagram using an out-split and a $2$-cell:
\begin{center}
\begin{tikzpicture}[scale=.8]
\draw [ultra thick] (1,2) rectangle (2,3);
\draw [implies] (-.5,2.5) to (1,2.5);
\draw [implies] (2,2.2) to (3.5,2.2);
\node at (1.7,2.7) {\tiny{$x_{12}$}};
\node at (1.5,2.3) {$X'$};
\draw [thick] (2,2.7) -- (2.3,2.7);
\draw [arrow] (2.3,2.7) to [out=0, in=180] (3.5,2.9);
\draw [arrow] (2.3,2.7) to [out=0, in=180] (3.5,2.5);
\node at (3.7, 2.9) {\tiny{$x_1$}};
\node at (3.7,2.5) {\tiny{$x_2$}};
\draw [ultra thick] (1,.5) rectangle (2,1.5);
\draw [implies] (-.5,1) to (1,1);
\draw [implies] (2,1) to (3.5,1);
\node at (1.5,1) {$Y$};
\draw [ultra thick] (0,0.2) rectangle (3,3.3);
\node at (4.4,2.5) {$\xout$};
\end{tikzpicture}
\end{center}

\begin{proposition}
\label{prop:move:b5}
Suppose:
\begin{itemize}
\item
$X \in \boxs$, and $x_1, x_2 \in \xout$ are distinct elements such that $v(x_1) = v(x_2) \in S$.
\item
$X' = X/(x_1=x_2) \in \boxs$ is obtained from $X$ by identifying $x_1$ and $x_2$.
\item
$\theta_{X,Y} \in \WD\xplusyxy$ and $\theta_{X',Y} \in \WD\xprimeplusyxprimey$ are $2$-cells.
\item
$\sigma^{X,x_1,x_2} \in \WD\xxprime$ and $\sigma^{X \amalg Y,x_1,x_2} \in \WD\xplusyxprimeplusy$ are out-splits.
\end{itemize}
Then
\begin{equation}
\label{move:b5}
\left(\theta_{X,Y}\right) 
\compone \left(\sigma^{X,x_1,x_2}\right)
=
\left(\sigma^{X\amalg Y,x_1,x_2}\right) 
\comp \left(\theta_{X',Y}\right)
\in \WD\xplusyxprimey.
\end{equation}
\end{proposition}

The next relation says that substituting a $1$-wasted wire inside a $2$-cell can be rewritten as substituting a $2$-cell inside a $1$-wasted wire.  It gives two different ways to construct the following wiring diagram using a $1$-wasted wire and a $2$-cell:
\begin{center}
\begin{tikzpicture}[scale=.8]
\draw [ultra thick] (1,2) rectangle (2,3);
\draw [implies] (-.5,2.3) to (1,2.3);
\draw [implies] (2,2.5) to (3.5,2.5);
\node at (1.5,2.5) {$X'$};
\draw [arrow] (-.5,3) to (0,3);
\node at (-.4,3.2) {\tiny{$x_0$}};
\draw [ultra thick] (1,.5) rectangle (2,1.5);
\draw [implies] (-.5,1) to (1,1);
\draw [implies] (2,1) to (3.5,1);
\node at (1.5,1) {$Y$};
\draw [ultra thick] (0,0.2) rectangle (3,3.3);
\node at (-1,2.5) {$\xin$};
\end{tikzpicture}
\end{center}

\begin{proposition}
\label{prop:move:b6}
Suppose:
\begin{itemize}
\item
$X \in \boxs$, $x_0 \in \xin$, and $X' = X \setminus \{x_0\} \in \boxs$ is obtained from $X$ by removing $x_0$.
\item
$\theta_{X,Y} \in \WD\xplusyxy$ and $\theta_{X',Y} \in \WD\xprimeplusyxprimey$ are $2$-cells.
\item
$\omega_{X,x_0} \in \WD\xxprime$ and $\omega_{X\amalg Y,x_0} \in \WD\xplusyxprimeplusy$ are $1$-wasted wires.
\end{itemize}
Then
\begin{equation}
\label{move:b6}
\left(\theta_{X,Y}\right) 
\compone \left(\omega_{X,x_0}\right)
= 
\left(\omega_{X\amalg Y,x_0}\right) 
\comp \left(\theta_{X',Y}\right)
\in \WD\xplusyxprimey.
\end{equation}
\end{proposition}

The following six relations are about $1$-loops.  The next relation is the commutativity property of $1$-loops.\index{commutativity of $1$-loops}  It gives two different ways to construct the following wiring diagram, which we will call a \emph{double-loop},\index{double-loop} using two $1$-loops:
\begin{center}
\begin{tikzpicture}
\draw [ultra thick] (2,1) rectangle (3,2);
\node at (2.5,1.2) {\small{$X$}};
\node at (2.2,1.8) {\tiny{$x^1_-$}};
\node at (2.8,1.8) {\tiny{$x^1_+$}};
\node at (2.2,1.4) {\tiny{$x^2_-$}};
\node at (2.8,1.4) {\tiny{$x^2_+$}};
\draw [implies] (.5,1.2) to (2,1.2);
\node at (-.5,1.6) {$\xin \setminus \{x^1_-,x^2_-\}$};
\draw [implies] (3,1.2) to (4.5,1.2);
\node at (5.6,1.6) {$\xout \setminus \{x^1_+,x^2_+\}$};
\draw [arrow, looseness=3] (3,1.8) to [out=30, in=150] (2,1.8);
\draw [arrow, looseness=8] (3,1.5) to [out=30, in=150] (2,1.5);
\draw [ultra thick] (1,0.7) rectangle (4,3);
\end{tikzpicture}
\end{center}

\begin{proposition}
\label{prop:move:c1}
Suppose:
\begin{itemize}
\item
$X \in \boxs$, $x^1_- \not= x^2_- \in \xin$, and $x^1_+ \not= x^2_+ \in \xout$ such that $v(x^1_+) = v(x^1_-) \in S$ and $v(x^2_+) = v(x^2_-) \in S$.
\item
$X \setminus x^1$, $X \setminus x^2$, and $X \setminus x \in \boxs$ are obtained from $X$ by removing $x^1_{\pm}$, $x^2_{\pm}$, and $\{x^1_{\pm}, x^2_{\pm}\}$, respectively.
\item
$\lambda_{X \setminus x^1, x^2} \in \WD\xminusxxminusxone$ and 
$\lambda_{X,x^1} \in \WD\xminusxonex$ are $1$-loops.
\item
$\lambda_{X\setminus x^2, x^1} \in \WD\xminusxxminusxtwo$ and $\lambda_{X,x^2} \in \WD\xminusxtwox$ are $1$-loops.
\end{itemize}
Then
\begin{equation}
\label{move:c1}
\left(\lambda_{X \setminus x^1, x^2}\right) 
\comp \left(\lambda_{X,x^1}\right)
=
\left(\lambda_{X \setminus x^2, x^1}\right) 
\comp \left(\lambda_{X,x^2}\right)
\in \WD\xminusxx.
\end{equation}
\end{proposition}

The next relation is the commutativity property between $1$-loops and in-splits.  It gives two different ways to construct the following wiring diagram using one $1$-loop and one in-split:
\begin{center}
\begin{tikzpicture}
\draw [ultra thick] (1,2) rectangle (2,3);
\draw [implies] (-.5,2.1) to (1,2.1);
\draw [implies] (2,2.5) to (3.5,2.5);
\node at (1.2,2.8) {\tiny{$x_-$}};
\node at (1.8,2.8) {\tiny{$x_+$}};
\node at (1.2,2.6) {\tiny{$x_1$}};
\node at (1.2,2.4) {\tiny{$x_2$}};
\node at (1.7,2.2) {$X$};
\node at (-1.3,2.3) {$\frac{\xin \setminus \{x_-\}}{(x_1\,=\,x_2)}$};
\node at (4.5,2.5) {$\xout \setminus \{x_+\}$};
\draw [thick] (-.5,2.5) -- (.2,2.5);
\draw [arrow] (.2,2.5) to [out=0, in=180] (1,2.6);
\draw [arrow] (.2,2.5) to [out=0, in=180] (1,2.4);
\draw [arrow, looseness=3] (2,2.8) to [out=30, in=150] (1,2.8);
\draw [ultra thick] (0,1.7) rectangle (3,3.4);
\end{tikzpicture}
\end{center}

\begin{proposition}
\label{prop:move:c2}
Suppose:
\begin{itemize}
\item
$X \in \boxs$, $x_-, x_1, x_2 \in \xin$ are distinct, and $x_+ \in \xout$ such that $v(x_+) = v(x_-) \in S$ and $v(x_1) = v(x_2) \in S$.
\item
$X \setminus x \in \boxs$ is obtained from $X$ by removing $x_{\pm}$.
\item
$X' = X/(x_1 = x_2) \in \boxs$ is obtained from $X$ by identifying $x_1$ and $x_2$.
\item
$X' \setminus x \in \boxs$ is obtained from $X'$ by removing $x_{\pm}$.
\item
$\lambda_{X',x} \in \WD\xprimeminusxxprime$ and $\lambda_{X,x} \in \WD\xminusxx$ are $1$-loops.
\item
$\sigma_{X,x_1,x_2} \in \WD\xprimex$ and $\sigma_{X\setminus x,x_1,x_2} \in \WD\xprimeminusxxminusx$ are in-splits.

\end{itemize}
Then
\begin{equation}
\label{move:c2}
\left(\lambda_{X',x}\right) 
\comp \left(\sigma_{X,x_1,x_2}\right)
=
\left(\sigma_{X\setminus x,x_1,x_2}\right)
\comp \left(\lambda_{X,x}\right) 
\in \WD\xprimeminusxx.
\end{equation}
\end{proposition}

The next relation is the commutativity property between $1$-loops and out-splits.  It gives two different ways to construct the following wiring diagram using one $1$-loop and one out-split:
\begin{center}
\begin{tikzpicture}
\draw [ultra thick] (1,2) rectangle (2,3);
\draw [implies] (-.5,2.5) to (1,2.5);
\draw [implies] (2,2.1) to (3.5,2.1);
\node at (1.2,2.8) {\tiny{$x_-$}};
\node at (1.8,2.8) {\tiny{$x_+$}};
\node at (1.8,2.5) {\tiny{$x_{12}$}};
\node at (1.5,2.2) {$X'$};
\node at (-1.5,2.5) {$\xin \setminus \{x_-\}$};
\node at (5,2.5) {$\xout \setminus \{x_+\}$};
\draw [thick] (2,2.6) -- (2.3,2.6);
\draw [arrow] (2.3,2.6) to [out=0, in=180] (3.5,2.8);
\draw [arrow] (2.3,2.6) to [out=0, in=180] (3.5,2.4);
\node at (3.7,2.8) {\tiny{$x_1$}};
\node at (3.7,2.4) {\tiny{$x_2$}};
\draw [arrow, looseness=3] (2,2.8) to [out=30, in=150] (1,2.8);
\draw [ultra thick] (0,1.7) rectangle (3,3.4);
\end{tikzpicture}
\end{center}

\begin{proposition}
\label{prop:move:c3}
Suppose:
\begin{itemize}
\item
$X \in \boxs$, and $(x_+,x_-) \in \xout \times \xin$ such that $v(x_+) = v(x_-) \in S$.
\item
$x_1 \not= x_2 \in \xout \setminus \{x_+\}$ such that $v(x_1) = v(x_2) \in S$.
\item
$X' = X/(x_1 = x_2) \in \boxs$ is obtained from $X$ by identifying $x_1$ and $x_2$.
\item
$X \setminus x \in \boxs$ is obtained from $X$ by removing $x_{\pm}$.
\item
$X' \setminus x \in \boxs$ is obtained from $X'$ by removing $x_{\pm}$.
\item
$\lambda_{X',x} \in \WD\xprimeminusxxprime$ and $\lambda_{X,x} \in \WD\xminusxx$ are $1$-loops.
\item
$\sigma^{X\setminus x, x_1,x_2} \in \WD\xminusxxprimeminusx$ and $\sigma^{X,x_1,x_2} \in \WD\xxprime$ are out-splits.
\end{itemize}
Then
\begin{equation}
\label{move:c3}
\left(\sigma^{X\setminus x, x_1,x_2}\right) 
\comp \left(\lambda_{X',x}\right)
=
\left(\lambda_{X,x}\right)  
\comp \left(\sigma^{X,x_1,x_2}\right)
\in 
\WD\xminusxxprime.
\end{equation}
\end{proposition}

The next relation is the commutativity property between $1$-loops and $1$-wasted wires.  It gives two different ways to construct the following wiring diagram using one $1$-loop and one $1$-wasted wire:
\begin{center}
\begin{tikzpicture}
\draw [ultra thick] (1,2) rectangle (2,3);
\draw [implies] (-.5,2.2) to (1,2.2);
\draw [implies] (2,2.5) to (3.5,2.5);
\node at (1.2,2.8) {\tiny{$x_-$}};
\node at (1.8,2.8) {\tiny{$x_+$}};
\node at (1.5,2.2) {$X'$};
\node at (-1.5,2.5) {$\xin \setminus \{x_-\}$};
\node at (4.5,2.5) {$\xout \setminus \{x_+\}$};
\draw [arrow, looseness=3] (2,2.8) to [out=30, in=150] (1,2.8);
\draw [ultra thick] (0,1.7) rectangle (3,3.4);
\draw [arrow] (-.5,2.8) to (0,2.8);
\node at (-.4,3) {\tiny{$x_0$}};
\end{tikzpicture}
\end{center}

\begin{proposition}
\label{prop:move:c4}
Suppose:
\begin{itemize}
\item
$X \in \boxs$, $(x_+,x_-) \in \xout \times \xin$ such that $v(x_+) = v(x_-) \in S$, and $x_0 \in \xin \setminus \{x_-\}$.
\item
$X' = X \setminus \{x_0\} \in \boxs$ is obtained from $X$ by removing $x_0$.
\item
$X' \setminus x \in \boxs$ is obtained from $X'$ by removing $x_{\pm}$.
\item
$X \setminus x \in \boxs$ is obtained from $X$ by removing $x_{\pm}$.
\item
$\lambda_{X',x} \in \WD\xprimeminusxxprime$ and $\lambda_{X,x} \in \WD\xminusxx$ are $1$-loops.
\item
$\omega_{X\setminus x, x_0} \in \WD\xminusxxprimeminusx$ and $\omega_{X,x_0} \in \WD\xxprime$ are $1$-wasted wires.
\end{itemize}
Then
\begin{equation}
\label{move:c4}
\left(\omega_{X\setminus x, x_0}\right) 
\comp \left(\lambda_{X',x}\right)
=
\left(\lambda_{X,x}\right)
\comp \left(\omega_{X,x_0}\right)
\in \WD\xminusxxprime.
\end{equation}
\end{proposition}

The next relation involves $1$-loops, in-splits, and out-splits.  It says that the following two wiring diagrams are equal:
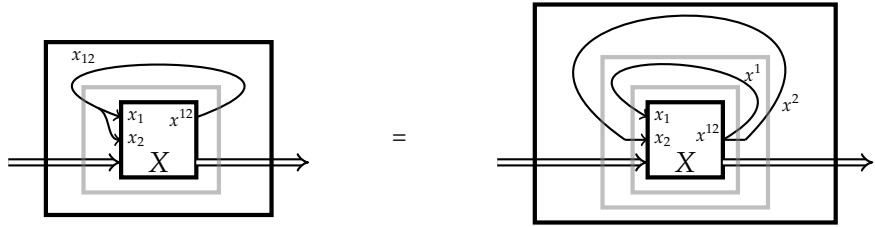
\begin{center}
\begin{tikzpicture}
\draw [ultra thick] (1,2) rectangle (2,3);
\draw [implies] (-.5,2.2) to (1,2.2);
\draw [implies] (2,2.2) to (3.5,2.2);
\node at (1.2,2.8) {\tiny{$x_1$}};
\node at (1.2,2.5) {\tiny{$x_2$}};
\node at (.5,3.6) {\tiny{$x_{12}$}};
\node at (1.8,2.8) {\tiny{$x^{12}$}};
\node at (1.5,2.2) {$X$};
\draw [arrow, looseness=7] (2,2.8) to [out=20, in=160] (1,2.8);
\draw [arrow] (.7,2.92) to [out=-30, in=180] (1,2.5);
\draw [ultra thick, gray, semitransparent] (.5,1.8) rectangle (2.3,3.2);
\draw [ultra thick] (0,1.5) rectangle (3,3.8);
\node at (4.7,2.5) {$=$};
\draw [ultra thick] (8,2) rectangle (9,3);
\draw [implies] (6,2.2) to (8,2.2);
\draw [implies] (9,2.2) to (11,2.2);
\node at (8.2,2.8) {\tiny{$x_1$}};
\node at (8.2,2.5) {\tiny{$x_2$}};
\node at (8.8,2.6) {\tiny{$x^{12}$}};
\node at (9.4,3.4) {\tiny{$x^1$}};
\node at (9.9,3) {\tiny{$x^2$}};
\node at (8.5,2.2) {$X$};
\draw [arrow, looseness=5.5] (9,2.5) to [out=30, in=150] (8,2.8);
\draw [thick] (9,2.5) to (9.3,2.5);
\draw [arrow] (7.7,2.5) to (8,2.5);
\draw [thick, looseness=5.5] (9.3,2.5) to [out=45, in=145] (7.7,2.5);
\draw [ultra thick, gray, semitransparent] (7.8,1.8) rectangle (9.2,3.2);
\draw [ultra thick, gray, semitransparent] (7.4,1.6) rectangle (9.6,3.6);
\draw [ultra thick] (6.5,1.4) rectangle (10.5,4.3);
\end{tikzpicture}
\end{center}
The wiring diagram on the left, in which the gray box is called $X'$ below, is created by substituting  an in-split into a $1$-loop.  The wiring diagram on the right is created by substituting an out-split into a $1$-loop, which is then substituted into another $1$-loop.  The inner gray box is called $Y$, and the outer gray box is called $Y \setminus x(1)$ below.  In both wiring diagrams, the outermost box is called $X^*$.

\begin{proposition}
\label{prop:move:c5}
Suppose:
\begin{itemize}
\item
$Y \in \boxs$, and $x^1 \not= x^2 \in \yout$ such that $v(x^1) = v(x^2) \in S$.
\item
$X = Y/(x^1 = x^2) \in \boxs$ is obtained from $Y$ by identifying $x^1$ and $x^2$, called $x^{12} \in \xout$.
\item
$x_1 \not=x_2 \in \xin = \yin$ such that $v(x^{12}) = v(x_1) = v(x_2) \in S$.
\item
$X' = X/(x_1=x_2) \in \boxs$ is obtained from $X$ by identifying $x_1$ and $x_2$, called $x_{12}$ in $X'$.
\item
$\sigma_{X,x_1,x_2} \in \WD\xprimex$ is an in-split.
\item
$X^* = X \setminus \{x^{12},x_1,x_2\} \in \boxs$ is obtained from $X$ by removing $x^{12},x_1$, and $x_2$.
\item
$\lambda_{X',x} \in \WD\xstarxprime$ is the $1$-loop in which $x^{12}$ is the supply wire of $x_{12}$.
\item
$\sigma^{Y,x^1,x^2} \in \WD\yx$ is an out-split.
\item
$\lambda_{Y,x(1)} \in \WD\yminusxoney$ is a $1$-loop, where $Y \setminus x(1) \in \boxs$ is obtained from $Y$ by removing $\{x_1,x^1\}$.
\item
$\lambda_{Y\setminus x(1), x(2)} \in \WD\xstaryminusxone$ is a $1$-loop, in which $x^2$ is the supply wire of $x_2$.
\end{itemize}
Then
\begin{equation}
\label{move:c5}
\left(\lambda_{X',x}\right)  \comp \left(\sigma_{X,x_1,x_2}\right)
= 
\left(\lambda_{Y\setminus x(1), x(2)}\right)
\comp
\left(\lambda_{Y,x(1)}\right) 
\comp 
\left(\sigma^{Y,x^1,x^2}\right)
\in \WD\xstarx.
\end{equation}
\end{proposition}

The next relation says that  the colored unit of a box can be rewritten as the substitution of an out-split into a $1$-wasted wire and then into a $1$-loop.   This is depicted in the picture
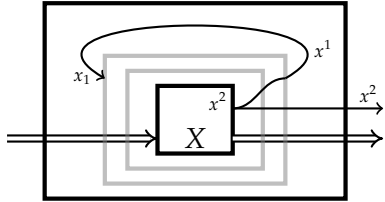
\begin{center}
\begin{tikzpicture}
\draw [ultra thick] (1,2) rectangle (2,2.9);
\draw [implies] (-1,2.2) to (1,2.2);
\draw [implies] (2,2.2) to (4,2.2);
\node at (1.8,2.7) {\tiny{$x^2$}};
\node at (0,3) {\tiny{$x_1$}};
\node at (3.2,3.4) {\tiny{$x^1$}};
\node at (3.8,2.8) {\tiny{$x^2$}};
\node at (1.5,2.2) {$X$};
\draw [ultra thick, gray, semitransparent] (.6,1.8) rectangle (2.4,3.1);
\draw [ultra thick, gray, semitransparent] (.3,1.6) rectangle (2.7,3.3);
\draw [arrow] (2,2.6) -- (4,2.6);
\draw [thick] (2,2.6) to [out=0, in=180] (2.7,3);
\draw [arrow, looseness=2] (2.7,3) to [out=30, in=150] (.3,3);
\draw [ultra thick] (-.5,1.4) rectangle (3.5,4);
\end{tikzpicture}
\end{center}
in which the outer gray box is called $Z$ and the inner gray box is called $Y$ below.

\begin{proposition}
\label{prop:move:c6}
Suppose:
\begin{itemize}
\item
$Z \in \boxs$, and $(x_1,x^1) \in \zin \times \zout$ such that $v(x_1) = v(x^1) \in S$.
\item
$Y = Z \setminus \{x_1\} \in \boxs$ is obtained from $Z$ by removing $x_1 \in \zin$.
\item
$x^1 \not= x^2 \in \yout = \zout$ such that $v(x^1) = v(x^2) \in S$.
\item
$X = Z \setminus \{x_1,x^1\} \in \boxs$ is obtained from $Z$ by removing $\{x_1,x^1\}$.
\item
$\sigma^{Y,x^1,x^2} \in \WD\yx$ is an out-split in which both $x^1, x^2 \in \yout$ have supply wire $x^2 \in \xout$.
\item
$\omega_{Z,x_1} \in \WD\zy$ is a $1$-wasted wire.
\item
$\lambda_{Z,x} \in \WD\xz$ is a $1$-loop in which $x^1 \in \zout$ is the supply wire of $x_1 \in \zin$.
\end{itemize}
Then
\begin{equation}
\label{move:c6}
\left(\lambda_{Z,x}\right) \comp
\left(\omega_{Z,x_1}\right) \comp
\left(\sigma^{Y,x^1,x^2}\right)
= \tensorunit_X
\in \WD\xx.
\end{equation}
\end{proposition}

The following five relations are about in-splits.  The next one is the associativity property of in-splits.\index{associativity of in-splits}  It gives two different ways to construct the following wiring diagram using two in-splits:
\begin{center}
\begin{tikzpicture}
\draw [ultra thick] (2,1) rectangle (3,2);
\node at (2.7,1.5) {$X$};
\node at (2.2,1.8) {\tiny{$x_1$}};
\node at (2.2,1.6) {\tiny{$x_2$}};
\node at (2.2,1.4) {\tiny{$x_3$}};
\draw [implies] (.5,1.1) to (2,1.1);
\node at (0,1.4) {$\yin$};
\draw [arrow] (.5,1.6) -- (2,1.6);
\draw [arrow] (1.2,1.6) to [out=0, in=180] (2,1.8);
\draw [arrow] (1.2,1.6) to [out=0, in=180] (2,1.4);
\draw [implies] (3,1.5) to (4.5,1.5);
\node at (5,1.6) {$\yout$};
\draw [ultra thick] (1,0.7) rectangle (4,2.3);
\end{tikzpicture}
\end{center}

\begin{proposition}
\label{prop:move:d1}
Suppose:
\begin{itemize}
\item
$X \in \boxs$, and $x_1,x_2,x_3 \in \xin$ are distinct elements such that $v(x_1) = v(x_2) = v(x_3) \in S$.
\item
$X_{12} = X/(x_1=x_2) \in \boxs$ is obtained from $X$ by identifying $x_1$ and $x_2$, called $x_{12} \in \xin_{12}$.
\item
$X_{23} = X/(x_2=x_3) \in \boxs$ is obtained from $X$ by identifying $x_2$ and $x_3$, called $x_{23} \in \xin_{23}$.
\item
$Y = X/(x_1=x_2=x_3) \in \boxs$ is  obtained from $X$ by identifying $x_1$, $x_2$, and $x_3$.
\item
$\sigma_{X_{12},x_{12},x_3} \in \WD\yxonetwo$ and $\sigma_{X, x_1,x_2} \in \WD\xonetwox$ are in-splits.
\item
$\sigma_{X_{23},x_1,x_{23}} \in \WD\yxtwothree$ and $\sigma_{X,x_2,x_3} \in \WD\xtwothreex$ are in-splits.
\end{itemize}
Then
\begin{equation}
\label{move:d1}
\left(\sigma_{X_{12},x_{12},x_3}\right) 
\comp \left(\sigma_{X, x_1,x_2}\right)
=
\left(\sigma_{X_{23},x_1,x_{23}}\right) 
\comp \left(\sigma_{X,x_2,x_3}\right)
\in \WD\yx.
\end{equation}
\end{proposition}

The next relation is the commutativity property of in-splits.\index{commutativity of in-splits}  It gives two different ways to construct the following wiring diagram using two in-splits:
\begin{center}
\begin{tikzpicture}
\draw [ultra thick] (2,.8) rectangle (3,2);
\node at (2.7,1.5) {$X$};
\node at (2.2,1.8) {\tiny{$x_1$}};
\node at (2.2,1.6) {\tiny{$x_2$}};
\node at (2.2,1.4) {\tiny{$x_3$}};
\node at (2.2,1.2) {\tiny{$x_4$}};
\draw [implies] (.5,.9) to (2,.9);
\node at (0,1.5) {$\yin$};
\draw [thick] (.5,1.7) -- (1.2,1.7);
\draw [arrow] (1.2,1.7) to [out=0, in=180] (2,1.8);
\draw [arrow] (1.2,1.7) to [out=0, in=180] (2,1.6);
\draw [thick] (.5,1.3) -- (1.2,1.3);
\draw [arrow] (1.2,1.3) to [out=0, in=180] (2,1.4);
\draw [arrow] (1.2,1.3) to [out=0, in=180] (2,1.2);
\draw [implies] (3,1.5) to (4.5,1.5);
\node at (5,1.6) {$\yout$};
\draw [ultra thick] (1,0.5) rectangle (4,2.3);
\end{tikzpicture}
\end{center}

\begin{proposition}
\label{prop:move:d2}
Suppose:
\begin{itemize}
\item
$X \in \boxs$, and $x_1,x_2,x_3,x_4 \in \xin$ are distinct elements such that $v(x_1) = v(x_2)$ and $v(x_3) = v(x_4) \in S$.
\item
$X_{12} = X/(x_1=x_2) \in \boxs$ is obtained from $X$ by identifying $x_1$ and $x_2$, called $x_{12} \in \xin_{12}$.
\item
$X_{34} = X/(x_3=x_4) \in \boxs$ is obtained from $X$ by identifying $x_3$ and $x_4$, called $x_{34} \in \xin_{34}$.
\item
$Y = X/(x_1=x_2;\, x_3=x_4) \in \boxs$ is  obtained from $X$ by (i) identifying $x_1$ and $x_2$ and (ii) identifying $x_3$ and $x_4$.
\item
$\sigma_{X_{12},x_3,x_4} \in \WD\yxonetwo$ and $\sigma_{X, x_1,x_2} \in \WD\xonetwox$ are in-splits.
\item
$\sigma_{X_{34},x_1,x_2} \in \WD\yxthreefour$ and $\sigma_{X,x_3,x_4} \in \WD\xthreefourx$ are in-splits.
\end{itemize}
Then
\begin{equation}
\label{move:d2}
\left(\sigma_{X_{12},x_3,x_4}\right) 
\comp \left(\sigma_{X, x_1,x_2}\right)
=
\left(\sigma_{X_{34},x_1,x_2}\right) 
\comp \left(\sigma_{X,x_3,x_4}\right)
\in \WD\yx.
\end{equation}
\end{proposition}

The next relation is the commutativity property between an in-split and an out-split.  It gives two different ways to construct the following wiring diagram using one in-split and one out-split:
\begin{center}
\begin{tikzpicture}
\draw [ultra thick] (2,1) rectangle (3,2);
\node at (2.7,1.2) {$X$};
\node at (2.2,1.8) {\tiny{$z_1$}};
\node at (2.2,1.4) {\tiny{$z_2$}};
\node at (2.8,1.8) {\tiny{$z^{12}$}};
\draw [thick] (3,1.8) -- (3.3,1.8);
\draw [arrow] (3.3,1.8) to [out=0, in=180] (4.5,2);
\draw [arrow] (3.3,1.8) to [out=0, in=180] (4.5,1.6);
\draw [implies] (.5,1.1) to (2,1.1);
\node at (-.5,1.5) {$\yin = \win$};
\draw [thick] (.5,1.6) -- (1.2,1.6);
\draw [arrow] (1.2,1.6) to [out=0, in=180] (2,1.8);
\draw [arrow] (1.2,1.6) to [out=0, in=180] (2,1.4);
\draw [implies] (3,1.2) to (4.5,1.2);
\node at (6,1.5) {$\yout = \zout$};
\node at (4.7,2) {\tiny{$z^{1}$}};
\node at (4.7,1.6) {\tiny{$z^{2}$}};
\draw [ultra thick] (1,0.7) rectangle (4,2.3);
\end{tikzpicture}
\end{center}

\begin{proposition}
\label{prop:move:d3}
Suppose:
\begin{itemize}
\item
$Z \in \boxs$, and $z^1 \not= z^2 \in \zout$ such that $v(z^1) = v(z^2) \in S$.
\item
$X = Z/(z^1=z^2) \in \boxs$ is obtained from $Z$ by identifying $z^1$ and $z^2$.
\item
$z_1 \not= z_2 \in \zin$ such that $v(z_1) = v(z_2) \in S$.
\item
$Y = Z/(z_1=z_2) \in \boxs$ is obtained from $Z$ by identifying $z_1$ and $z_2$.
\item
$W = Z/(z_1=z_2; \, z^1=z^2) \in \boxs$ is obtained from $Z$ by (i) identifying $z_1$ and $z_2$ and (ii) identifying $z^1$ and $z^2$.
\item
$\sigma^{Y,z^1,z^2} \in \WD\yw$ is an out-split, and $\sigma_{X,z_1,z_2} \in \WD\wx$ is an in-split.
\item
$\sigma_{Z,z_1,z_2} \in \WD\yz$ is an in-split, and $\sigma^{Z,z^1,z^2} \in \WD\zx$ is an out-split.
\end{itemize}
Then
\begin{equation}
\label{move:d3}
\left(\sigma^{Y,z^1,z^2}\right) \comp
\left(\sigma_{X,z_1,z_2}\right)
=
\left(\sigma_{Z,z_1,z_2}\right) \comp
\left(\sigma^{Z,z^1,z^2}\right)
\in \WD\yx.
\end{equation}
\end{proposition}

The next relation is the commutativity property between an in-split and a $1$-wasted wire.  It gives two different ways to construct the following wiring diagram using one in-split and one $1$-wasted wire:
\begin{center}
\begin{tikzpicture}
\draw [ultra thick] (2,1) rectangle (3,2);
\node at (2.7,1.2) {$X$};
\node at (2.2,1.8) {\tiny{$z_1$}};
\node at (2.2,1.4) {\tiny{$z_2$}};
\draw [implies] (.5,1.1) to (2,1.1);
\node at (0,1.5) {$\yin$};
\draw [thick] (.5,1.6) -- (1.2,1.6);
\draw [arrow] (1.2,1.6) to [out=0, in=180] (2,1.8);
\draw [arrow] (1.2,1.6) to [out=0, in=180] (2,1.4);
\draw [arrow] (.5,2) to (1,2);
\node at (.6,2.2) {\tiny{$z$}};
\draw [implies] (3,1.5) to (4.5,1.5);
\draw [ultra thick] (1,0.7) rectangle (4,2.3);
\end{tikzpicture}
\end{center}

\begin{proposition}
\label{prop:move:d4}
Suppose:
\begin{itemize}
\item
$Z \in \boxs$, and $z, z_1, z_2 \in \zin$ are distinct elements such that $v(z_1) = v(z_2) \in S$.
\item
$Y = Z/(z_1=z_2) \in \boxs$ is obtained from $Z$ by identifying $z_1$ and $z_2$.
\item
$X = Z \setminus \{z\} \in \boxs$ is obtained from $Z$ by removing $z \in \zin$.
\item
$W = X/(z_1=z_2) \in \boxs$ is obtained from $X$ by identifying $z_1$ and $z_2$.
\item
$\sigma_{X,z_1,z_2} \in \WD\wx$ and $\sigma_{Z,z_1,z_2} \in \WD\yz$ are in-splits.
\item
$\omega_{Y,z} \in \WD\yw$ and $\omega_{Z,z} \in \WD\zx$ are $1$-wasted wires.
\end{itemize}
Then
\begin{equation}
\label{move:d4}
\left(\omega_{Y,z}\right) \comp \left(\sigma_{X,z_1,z_2}\right)
=
\left(\sigma_{Z,z_1,z_2}\right) \comp \left(\omega_{Z,z}\right)
\in 
\WD\yx.
\end{equation}
\end{proposition}

The next relation says that  the colored unit of a box $X$ can be rewritten as the substitution of a $1$-wasted wire into an in-split.  This is depicted in the picture
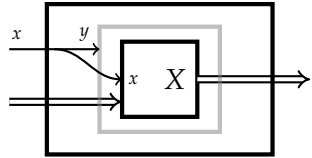
\begin{center}
\begin{tikzpicture}
\draw [ultra thick] (1,2) rectangle (2,3);
\draw [implies] (-.5,2.2) to (1,2.2);
\draw [implies] (2,2.5) to (3.5,2.5);
\node at (1.15,2.5) {\tiny{$x$}};
\node at (.5,3.1) {\tiny{$y$}};
\node at (1.7,2.5) {$X$};
\draw [ultra thick, gray, semitransparent] (.7,1.8) rectangle (2.3,3.2);
\draw [arrow] (-.5,2.9) -- (.7,2.9);
\node at (-.4,3.1) {\tiny{$x$}};
\draw [arrow] (.1,2.9) to [out=0, in=180] (1,2.5);
\draw [ultra thick] (0,1.5) rectangle (3,3.5);
\end{tikzpicture}
\end{center}
in which the intermediate gray box is called $Y$ below.

\begin{proposition}
\label{prop:move:d5}
Suppose:
\begin{itemize}
\item
$Y \in \boxs$, and $x,y \in \yin$ are distinct elements such that $v(x) = v(y) \in S$.
\item
$X = Y/(x=y) \in \boxs$ is obtained from $Y$ by identifying $x$ and $y$.
\item
$\omega_{Y,y} \in \WD\yx$ is a $1$-wasted wire.
\item
$\sigma_{Y,x,y} \in \WD\smallbinom{X}{Y}$ is an in-split.
\end{itemize}
Then
\begin{equation}
\label{move:d5}
\left(\sigma_{Y,x,y}\right) 
\comp \left(\omega_{Y,y}\right)
=
\tensorunit_X
\in \WD\xx.
\end{equation}
\end{proposition}

The following three relations are about out-splits.  The next one is the associativity property of out-splits.\index{associativity of out-splits}  It gives two different ways to construct the following wiring diagram using two out-splits:
\begin{center}
\begin{tikzpicture}
\draw [ultra thick] (2,1) rectangle (3,2);
\node at (2.5,1.5) {$X$};
\draw [implies] (.5,1.5) to (2,1.5);
\node at (0,1.6) {$\yin$};
\draw [arrow] (3,1.8) -- (4.5,1.8);
\draw [arrow] (3.1,1.8) to [out=0, in=180] (4.5,2);
\draw [arrow] (3.1,1.8) to [out=0, in=180] (4.5,1.6);
\node at (4.7,2.1) {\tiny{$y^1$}};
\node at (4.7,1.8) {\tiny{$y^2$}};
\node at (4.7,1.5) {\tiny{$y^3$}};
\draw [implies] (3,1.1) to (4.5,1.1);
\node at (5.5,1.6) {$\yout$};
\draw [ultra thick] (1,0.7) rectangle (4,2.3);
\end{tikzpicture}
\end{center}

\begin{proposition}
\label{prop:move:e1}
Suppose:
\begin{itemize}
\item
$Y \in \boxs$, and $y^1, y^2, y^3 \in \yout$ are distinct elements such that $v(y^1) = v(y^2) = v(y^3) \in S$.
\item
$Y^{12} = Y/(y^1=y^2) \in \boxs$ is obtained from $Y$ by identifying $y^1$ and $y^2$, called $y^{12}$ in $Y^{12}$.
\item
$Y^{23} = Y/(y^2 = y^3) \in \boxs$ is obtained from $Y$ by identifying $y^2$ and $y^3$, called $y^{23}$ in $Y^{23}$.
\item
$X = Y/(y^1 = y^2 = y^3) \in \boxs$  is obtained from $Y$ by identifying $y^1$, $y^2$, and $y^3$.
\item
$\sigma^{Y,y^1,y^2} \in \WD\yyuponetwo$ and $\sigma^{Y^{12},y^{12},y^3} \in \WD\yuponetwox$ are out-splits.
\item
$\sigma^{Y,y^2,y^3} \in \WD\yyuptwothree$ and $\sigma^{Y^{23},y^{1},y^{23}} \in \WD\yuptwothreex$ are out-splits.
\end{itemize}
Then
\begin{equation}
\label{move:e1}
\left(\sigma^{Y,y^1,y^2}\right) 
\comp
\left(\sigma^{Y^{12},y^{12},y^3}\right)
=
\left(\sigma^{Y,y^2,y^3}\right)
\comp
\left(\sigma^{Y^{23},y^{1},y^{23}}\right)
\in \WD\yx.
\end{equation}
\end{proposition}

The next relation is the commutativity property of out-splits.\index{commutativity of out-splits}  It gives two different ways to construct the following wiring diagram using two out-splits:
\begin{center}
\begin{tikzpicture}
\draw [ultra thick] (2,1) rectangle (3,2);
\node at (2.5,1.5) {$X$};
\draw [implies] (.5,1.5) to (2,1.5);
\node at (0,1.6) {$\yin$};
\draw [arrow] (3,1.9) to [out=0, in=180] (4.5,2);
\draw [arrow] (3,1.9) to [out=0, in=180] (4.5,1.8);
\draw [arrow] (3,1.5) to [out=0, in=180] (4.5,1.6);
\draw [arrow] (3,1.5) to [out=0, in=180] (4.5,1.4);
\node at (4.7,2.2) {\tiny{$y^1$}};
\node at (4.7,1.9) {\tiny{$y^2$}};
\node at (4.7,1.6) {\tiny{$y^3$}};
\node at (4.7,1.3) {\tiny{$y^4$}};
\draw [implies] (3,1.1) to (4.5,1.1);
\node at (5.7,1.6) {$\yout$};
\draw [ultra thick] (1,0.7) rectangle (4,2.3);
\end{tikzpicture}
\end{center}

\begin{proposition}
\label{prop:move:e2}
Suppose:
\begin{itemize}
\item
$Y \in \boxs$, and $y^1, y^2, y^3,y^4 \in \yout$ are distinct elements such that $v(y^1) = v(y^2)$ and $v(y^3) = v(y^4) \in S$.
\item
$Y^{12} = Y/(y^1=y^2) \in \boxs$ is obtained from $Y$ by identifying $y^1$ and $y^2$.
\item
$Y^{34} = Y/(y^3 = y^4) \in \boxs$ is obtained from $Y$ by identifying $y^3$ and $y^4$.
\item
$X = Y/(y^1 = y^2;\, y^3 = y^4) \in \boxs$  is obtained from $Y$ by (i) identifying $y^1$ and $y^2$ and (ii) identifying $y^3$ and $y^4$.
\item
$\sigma^{Y,y^1,y^2} \in \WD\yyuponetwo$ and $\sigma^{Y^{12},y^3,y^4} \in \WD\yuponetwox$ are out-splits.
\item
$\sigma^{Y,y^3,y^4} \in \WD\yyupthreefour$ and $\sigma^{Y^{34},y^1,y^2} \in \WD\yupthreefourx$ are out-splits.
\end{itemize}
Then
\begin{equation}
\label{move:e2}
\left(\sigma^{Y,y^1,y^2}\right) \comp
\left(\sigma^{Y^{12},y^3,y^4}\right)
=
\left(\sigma^{Y,y^3,y^4}\right) \comp
\left(\sigma^{Y^{34},y^1,y^2}\right)
\in \WD\yx. 
\end{equation}
\end{proposition}

The next relation is the commutativity property between an out-split and a $1$-wasted wire.  It gives two different ways to construct the following wiring diagram using one out-split and one $1$-wasted wire:
\begin{center}
\begin{tikzpicture}
\draw [ultra thick] (2,1) rectangle (3,2);
\node at (2.5,1.5) {$X$};
\draw [implies] (.5,1.2) to (2,1.2);
\node at (0,1.6) {$\yin$};
\draw [arrow] (.5,1.7) to (1,1.7);
\node at (.7,1.9) {\tiny{$y$}};
\draw [thick] (3,1.7) -- (3.2,1.7);
\draw [arrow] (3.2,1.7) to [out=0, in=180] (4.5,1.9);
\draw [arrow] (3.2,1.7) to [out=0, in=180] (4.5,1.5);
\node at (4.7,2) {\tiny{$y^1$}};
\node at (4.7,1.6) {\tiny{$y^2$}};
\draw [implies] (3,1.2) to (4.5,1.2);
\node at (5.7,1.6) {$\yout$};
\draw [ultra thick] (1,0.7) rectangle (4,2.3);
\end{tikzpicture}
\end{center}

\begin{proposition}
\label{prop:move:e3}
Suppose:
\begin{itemize}
\item
$Y \in \boxs$, $y \in \yin$, and $y^1, y^2 \in \yout$ such that $v(y^1) = v(y^2) \in S$.
\item
$W = Y/(y^1 = y^2) \in \boxs$ is obtained from $Y$ by identifying $y^1$ and $y^2$.
\item
$Z = Y \setminus \{y\} \in \boxs$ is obtained from $Y$ by removing $y \in \yin$.
\item
$X = Z/(y^1=y^2) \in \boxs$ is obtained from $Z$ by identifying $y^1$ and $y^2$.
\item
$\sigma^{Z,y^1,y^2} \in \WD\zx$ and $\sigma^{Y,y^1,y^2} \in \WD\yw$ are out-splits.
\item
$\omega_{Y,y} \in \WD\yz$ and $\omega_{W,y} \in \WD\wx$ are $1$-wasted wires.
\end{itemize}
Then
\begin{equation}
\label{move:e3}
\left(\omega_{Y,y}\right) \comp
\left(\sigma^{Z,y^1,y^2}\right)
=
\left(\sigma^{Y,y^1,y^2}\right) \comp
\left(\omega_{W,y}\right)
\in \WD\yx.
\end{equation}
\end{proposition}

The final relation is the commutativity property of $1$-wasted wires.\index{commutativity of $1$-wasted wires}  It gives two different ways to construct the following wiring diagram using two $1$-wasted wires:
\begin{center}
\begin{tikzpicture}
\draw [ultra thick] (2,1) rectangle (3,2);
\node at (2.5,1.5) {$X$};
\draw [implies] (.5,1.1) to (2,1.1);
\node at (0,1.6) {$\yin$};
\draw [arrow] (.5,2) to (1,2);
\node at (.7,2.2) {\tiny{$y_1$}};
\draw [arrow] (.5,1.5) to (1,1.5);
\node at (.7,1.65) {\tiny{$y_2$}};
\draw [implies] (3,1.5) to (4.5,1.5);
\node at (5,1.6) {$\yout$};
\draw [ultra thick] (1,0.7) rectangle (4,2.3);
\end{tikzpicture}
\end{center}

\begin{proposition}
\label{prop:move:f1}
Suppose:
\begin{itemize}
\item
$Y \in \boxs$, and $y_1,y_2 \in \yin$ are distinct elements.
\item
$Y_1 = Y \setminus \{y_1\} \in \boxs$ is obtained from $Y$ by removing $y_1$.
\item
$Y_2 = Y \setminus \{y_2\} \in \boxs$ is obtained from $Y$ by removing $y_2$.
\item
$X = Y \setminus \{y_1,y_2\} \in \boxs$ is obtained from $Y$ by removing $y_1$ and $y_2$.
\item
$\omega_{Y,y_1} \in \WD\yyone$ and $\omega_{Y_1,y_2} \in \WD\yonex$ are $1$-wasted wires.
\item
$\omega_{Y,y_2} \in \WD\yytwo$ and $\omega_{Y_2,y_1} \in \WD\ytwox$ are $1$-wasted wires.
\end{itemize}
Then
\begin{equation}
\label{move:f1}
\left(\omega_{Y,y_1}\right) \comp \left(\omega_{Y_1,y_2}\right)
=
\left(\omega_{Y,y_2}\right) \comp \left(\omega_{Y_2,y_1}\right)
\in \WD\yx.
\end{equation}
\end{proposition}

\begin{definition}
\label{def:elementary-relations}
The $28$ relations \eqref{move:a1}--\eqref{move:f1} are called \emph{elementary relations}.\index{elementary relations}
\end{definition}

\section{Summary of Chapter \ref{ch03-generating-wd}}

\begin{enumerate}
\item There are eight generating wiring diagrams in $\WD$.
\item Each internal wasted wire can be generated using a $1$-wasted wire and a $1$-loop.
\item There are twenty-eight elementary relations in $\WD$.
\end{enumerate}

\chapter{Decomposition of Wiring Diagrams}
\label{ch04-decomposition}

As part of the finite presentation theorem for the operad $\WD$ of wiring diagrams (Theorem \ref{thm:wd-generator-relation}), in Theorem \ref{stratified-presentation-exists} we will observe that each wiring diagram has a highly structured decomposition into generating wiring diagrams (Def. \ref{def:generating-wiring-diagrams}), called a stratified presentation.  Stratified presentations are also needed to establish the second part of the finite presentation theorem for $\WD$ regarding relations.  The purpose of this chapter is to provide all the steps needed to establish the existence of a stratified presentation for each wiring diagram.  We remind the reader about Notation \ref{comp-is-compone} for (iterated) operadic compositions.

In Section \ref{sec:factoring-wd} we show that each wiring diagram $\psi$ has a specific operadic \index{decomposition} decomposition \eqref{psi-is-alpha-comp-phi}
\[
\psi = \alpha \comp \varphi.\]
An explanation of this decomposition is given just before Def. \ref{def:wd-phi}.  The idea of this decomposition is that we are breaking the complexity of a general wiring diagram into two simpler parts.  On the one hand, the inner wiring diagram $\varphi$ contains all the input boxes and the delay nodes of $\psi$, but its supplier assignment is as simple as possible, namely the identity map. See Lemmas \ref{lemma:phi-small-Nr} and \ref{lemma:phi-presentation}.  On the other hand, the outer wiring diagram $\alpha$ has only one input box and no delay nodes, but its supplier assignment is equal to that of $\psi$.

In Section \ref{sec:unary-wd} we observe that the outer wiring diagram $\alpha$ in the previous decomposition of $\psi$ can be decomposed as \eqref{pi-pione-pitwo} 
\[
\alpha = \pione \comp \pitwo.\]  
Example \ref{ex:factoring-pi} has a concrete wiring diagram that illustrates this decomposition.  The idea of this decomposition is that in a wiring diagram there are usually wires that go backward (i.e., "point to the left"), as in \eqref{wd-first-example}, in a $1$-loop (Def. \ref{def:loop-wd}), and in the pictures just before Prop. \ref{prop:move:b3} and Prop. \ref{prop:move:c1}.  This decomposition breaks the complexity of the wiring diagram $\alpha$ into two simpler parts.  On the one hand, the outer wiring diagram $\pione$ contains all the backward-going wires in $\alpha$ but no wasted wires or split wires (Lemma \ref{lemma:iterated-loops}).  On the other hand, the inner wiring diagram $\pitwo$ contains no backward-going wires, but it has all the wasted wires and split wires in $\alpha$.

In Section \ref{sec:unary-no-loop} we observe that the wiring diagram $\pitwo$ in the previous decomposition of $\alpha$ can be decomposed further as  \eqref{betaonetwothree}
\[
\pitwo = \betaone \comp \betatwo \comp \betathree.\]
Example \ref{ex:factor-pitwo} has a concrete wiring diagram that illustrates this decomposition.  In this decomposition:
\begin{itemize}
\item
The outermost wiring diagram $\betaone$ is an iterated operadic composition of $1$-wasted wires (Lemma \ref{lemma:betaone}).  
\item
The middle wiring diagram $\betatwo$ is an iterated operadic composition of in-splits (Lemma \ref{lemma:betatwo}).  
\item
The innermost wiring diagram $\betathree$ is an iterated operadic composition of out-splits (Lemma \ref{lemma:betathree}).  
\end{itemize}
By convention an empty operadic composition means a colored unit.  In summary, for a wiring diagram $\psi$, we will decompose it as
\[
\psi = \pione \comp  \betaone \comp \betatwo \comp \betathree \comp \varphi.\]

\section{Factoring Wiring Diagrams}
\label{sec:factoring-wd}

\begin{assumption}
\label{assumption:psi}
Throughout this chapter, fix a class $S$.  Suppose 
\begin{equation}
\label{psi-wiring-diagram}
\psi = \left(\dnpsi, \vpsi, \spsi\right) \in \WD\yux
\end{equation}
is a wiring diagram with:
\begin{itemize}
\item
output box $Y \in \boxs$ and input boxes $\uX = (X_1, \ldots, X_N)$ for some $N \geq 0$;
\item
$r$ delay nodes $\dnpsi = \{d_1,\ldots,d_r\}$ for some $r \geq 0$;
\item
value assignment $\vpsi : \yin \amalg \yout \amalg \xin \amalg \xout \amalg \dnpsi \to S$, where $\xin = \coprod_{i=1}^N \xin_i$ and $\xout = \coprod_{i=1}^N \xout_i$;
\item
supplier assignment $\spsi : \dmpsi \to \supplypsi$.
\end{itemize}
Since $N=0$ and $r=0$ are both allowed, $\psi$ is a general wiring diagram.  Furthermore:
\begin{enumerate}
\item
To simplify notation, we will write $\vpsi(d_i) \in S$ simply as $d_i$, so each $\delta_{d_i} \in \WD\dinothing$ is a $1$-delay node (Def. \ref{def:one-dn}). 
\item
$X = \coprod_{i=1}^N X_i \in \boxs$ is the coproduct of the $X_i$'s.
\item
Define $X' \in \boxs$ as
\begin{equation}
\label{xprime-inout}
\xprimein = \xin \amalg \dnpsi \andspace
\xprimeout = \xout \amalg \dnpsi.
\end{equation}
\end{enumerate}
\end{assumption}

\begin{motivation}
The first observation is about the marginal case where $\psi$ has no input boxes and no delay nodes.  So it looks like this
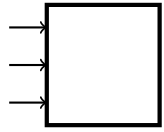
\begin{center}
\begin{tikzpicture}[scale=1]
\draw [ultra thick] (0,0.2) rectangle (1.5,1.8);
\draw [arrow] (-.5,.5) to (0,.5);
\draw [arrow] (-.5,1) to (0,1);
\draw [arrow] (-.5,1.5) to (0,1.5);
\end{tikzpicture}
\end{center}
with finitely many, possibly zero, external wasted wires.  In any case, it can be written in terms of the empty wiring diagram and finitely many $1$-wasted wires.
\end{motivation}

\begin{lemma}
\label{lemma:N=r=0}
Suppose $N = r =0$ in $\psi$; i.e., $\psi \in \WD\ynothing$ has no input boxes and no delay nodes.  Then one of the following two statements is true.
\begin{enumerate}
\item
$\psi = \epsilon \in \WD\emptynothing$, the empty wiring diagram (Def. \ref{def:empty-wd}).
\item
There exist $1$-wasted wires (Def. \ref{def:wasted-wire-wd}) $\omega_1, \ldots, \omega_m$, where $m = |\yin| > 0$, such that
\begin{equation}
\label{only-external-wasted-wires}
\psi = \omega_1 \comp \cdots \comp \omega_m \comp \epsilon.
\end{equation}
\end{enumerate}
\end{lemma}

\begin{proof}
Since $\xin = \xout = \dnpsi = \varnothing$, the supplier assignment of $\psi$ is a function
\[\nicexy{\dmpsi = \yout \ar[r]^-{s}
& \yin = \supplypsi.}\]
The non-instantaneity requirement \eqref{non-instant} then implies $\yout=\varnothing$.  If $m=|\yin| = 0$, then $Y$ is the empty box and $\psi = \epsilon$, the empty wiring diagram, by definition.

If $m > 0$, then every global input $y \in \yin = \{y_1, \ldots, y_m\}$ is an external wasted wire, and the supplier assignment $s : \yout = \varnothing \to \yin$ is the trivial map.  For each $1 \leq j \leq m$, define the box
\[Y_j = \begin{cases}
Y & \text{ if $j=1$};\\
Y \setminus \{y_1, \ldots , y_{j-1}\} & \text{ if $2 \leq j \leq m$};\\
\varnothing & \text{ if $j = m+1$}.
\end{cases}\]
Each $\omega_{Y_j,y_j} \in \WD\yjyjplusone$ is a $1$-wasted wire.  Using the notation \eqref{iterated-compone}, the iterated composition
\[\omega_{Y_1,y_1} \comp \cdots \comp \omega_{Y_m,y_m} \comp \epsilon \in
\WD\ynothing\]
is then a wiring diagram with output box $Y_1 = Y$, no input boxes and no delay nodes, and supplier assignment $\varnothing = \yout \to \yin$ the trivial map.  This is  the same as $\psi$.
\end{proof}

Next, for wiring diagrams $\psi$ not necessarily covered by Lemma \ref{lemma:N=r=0}, we define two relatively simple wiring diagrams that will be shown to provide a decomposition for $\psi$.  Each of these two simpler wiring diagrams will then be analyzed further.  

\begin{motivation}\label{mot:psi-alpha-phi}
This decomposition for $\psi$ is depicted in the following picture:
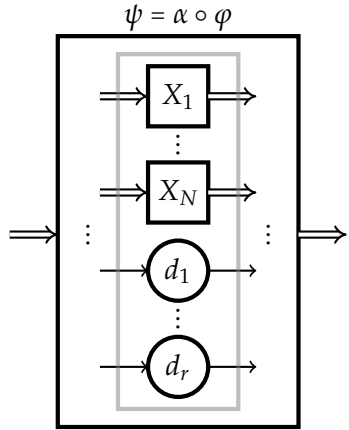
\begin{center}
\begin{tikzpicture}[scale=.8]
\draw [ultra thick] (2,5) rectangle (3,6);
\node at (2.5,5.5) {$X_1$};
\draw [implies] (1.2,5.5) to (2,5.5);
\draw [implies] (3,5.5) to (3.8,5.5);
\node at (2.5,4.7) {$\vdots$};
\draw [ultra thick] (2,3.4) rectangle (3,4.4);
\node at (2.5,3.9) {$X_N$};
\draw [implies] (1.2,3.9) to (2,3.9);
\draw [implies] (3,3.9) to (3.8,3.9);
\draw [ultra thick] (2.5,2.6) circle [radius=0.5];
\node at (2.5,2.6) {$d_1$};
\draw [arrow] (1.2,2.6) to (2,2.6);
\draw [arrow] (3,2.6) to (3.8,2.6);
\node at (2.5,1.8) {$\vdots$};
\draw [ultra thick] (2.5,1) circle [radius=0.5];
\node at (2.5,1) {$d_r$};
\draw [arrow] (1.2,1) to (2,1);
\draw [arrow] (3,1) to (3.8,1);
\draw [ultra thick, gray, semitransparent] (1.5,.3) rectangle (3.5,6.2);
\node at (1,3.2) {$\vdots$};
\node at (4,3.2) {$\vdots$};
\draw [ultra thick] (.5,0) rectangle (4.5,6.5);
\draw [implies] (-.3,3.2) to (.5,3.2);
\draw [implies] (4.5,3.2) to (5.3,3.2);
\node at (2.5,6.8) {$\psi = \alpha \comp \varphi$};
\end{tikzpicture}
\end{center}
Here the intermediate gray box is $X'$  \eqref{xprime-inout}.  In this decomposition, the inside wiring diagram $\varphi$ has all the input boxes and the delay nodes of $\psi$, but its supplier assignment is the identity function.  The outside wiring diagram $\alpha$ has a single input box $X'$ and no delay nodes, but it has the same supplier assignment as $\psi$.
\end{motivation}

\begin{definition}
\label{def:wd-phi}
Suppose $\psi$ is as in Assumption \ref{assumption:psi}.  Define the wiring diagram
\begin{equation}
\label{phi-wiring-diagram}
\varphi = \left(\dnphi, \vphi, \sphi\right) \in \WD\xprimeux
\end{equation}
with:
\begin{itemize}
\item
output box $X'$ \eqref{xprime-inout} and input boxes $\uX = (X_1, \ldots, X_N)$;
\item
delay nodes $\dnphi = \dnpsi = \{d_1, \ldots, d_r\}$;
\item
supplier assignment
\[\nicexy{\dmphi = \xprimeout \amalg \xin \amalg \dnphi
=
\left(\xout \amalg \dnpsi\right) \amalg \left(\xin \amalg \dnpsi\right)
\ar[d]_-{\sphi}
\\
\supplyphi = \xprimein \amalg \xout \amalg \dnphi
=
\left(\xin \amalg \dnpsi\right) \amalg \left(\xout \amalg \dnpsi\right)}\]
the identity function that sends $\xprimeout$ to $(\xout \amalg \dnpsi)$ and $(\xin \amalg \dnpsi)$ to $\xprimein$.
\end{itemize}
\end{definition}

\begin{definition}
\label{def:wd-alpha}
Suppose $\psi$ is as in Assumption \ref{assumption:psi}.  Define the wiring diagram
\begin{equation}
\label{alpha-wiring-diagram}
\alpha = \bigl(\dnalpha, \valpha, \salpha\bigr) \in \WD\yxprime
\end{equation}
with:
\begin{itemize}
\item
one input box $X'$  \eqref{xprime-inout}  and output box $Y$;
\item
no delay nodes;
\item
supplier assignment
\[\nicexy{\dmalpha = \yout \amalg \xprimein = \yout \amalg \xin \amalg \dnpsi = \dmpsi
\ar[d]_-{\salpha}
\\
\supplyalpha = \yin \amalg \xprimeout = \yin \amalg \xout \amalg \dnpsi = \supplypsi}\]
equal to $\spsi$.
\end{itemize}
\end{definition}

\begin{lemma}
\label{lemma:psi-alpha-comp-phi}
Given a wiring diagram $\psi$ \eqref{psi-wiring-diagram}, there is a decomposition
\begin{equation}
\label{psi-is-alpha-comp-phi}
\psi = \alpha \comp \varphi \in \WD\yux
\end{equation}
in which $\alpha \in \WD\yxprime$ is as in \eqref{alpha-wiring-diagram} and $\varphi \in \WD\xprimeux$ is as in \eqref{phi-wiring-diagram}.
\end{lemma}

\begin{proof}
By the definition of $\compone$ (Def. \ref{def:compi-wd}), $\alpha \comp \varphi = \alpha \compone \varphi$ belongs to $\WD\yux$ and has $\dnphi = \dnpsi$ as its set of delay nodes.  It remains to check that its supplier assignment is equal to that of $\psi$.  This follows from a direct inspection because $\sphi$ is the identity function, while $\salpha = \spsi$.
\end{proof}

To obtain the desired stratified presentation of $\psi$, we now begin to analyze the wiring diagram $\varphi$.

\begin{lemma}
\label{lemma:phi-small-Nr}
Consider the wiring diagram $\varphi$ \eqref{phi-wiring-diagram}.
\begin{enumerate}
\item
If $N=r=0$, then $\varphi = \epsilon$, the empty wiring diagram (Def. \ref{def:empty-wd}).
\item
If $(N,r) = (1,0)$, then $\varphi = \tensorunit_{X_1}$, the colored unit of $X_1$ (Def. \ref{wd-units}).
\item
If $(N,r) = (0,1)$, then $\varphi = \delta_{d_1}$, a $1$-delay node (Def. \ref{def:one-dn}).
\end{enumerate}
\end{lemma}

\begin{proof}
All three cases are checked by direct inspection.
\end{proof}

\begin{motivation}\label{mot:phi-2cells}
Next we observe that, for higher values of $N+r$, the wiring diagram $\varphi$ is generated by $2$-cells (Def. \ref{def:theta-wd}) and $1$-delay nodes via iterated operadic compositions, as depicted in the following picture.
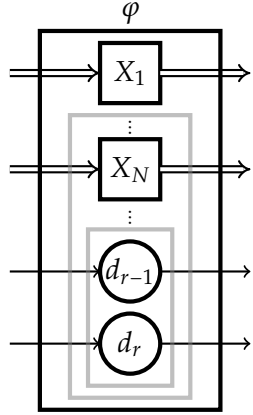
\begin{center}
\begin{tikzpicture}[scale=.8]
\draw [ultra thick] (2,5) rectangle (3,6);
\node at (2.5,5.5) {$X_1$};
\draw [implies] (.5,5.5) to (2,5.5);
\draw [implies] (3,5.5) to (4.5,5.5);
\node at (2.5,4.6) {\tiny{$\vdots$}};
\draw [ultra thick] (2,3.4) rectangle (3,4.4);
\node at (2.5,3.9) {$X_N$};
\draw [implies] (.5,3.9) to (2,3.9);
\draw [implies] (3,3.9) to (4.5,3.9);
\draw [ultra thick] (2.5,2.2) circle [radius=0.5];
\node at (2.5,2.2) {$d_{r-1}$};
\draw [arrow] (.5,2.2) to (2,2.2);
\draw [arrow] (3,2.2) to (4.5,2.2);
\node at (2.5,3.1) {\tiny{$\vdots$}};
\draw [ultra thick] (2.5,1) circle [radius=0.5];
\node at (2.5,1) {$d_r$};
\draw [arrow] (.5,1) to (2,1);
\draw [arrow] (3,1) to (4.5,1);
\draw [ultra thick, gray, semitransparent] (1.5,.1) rectangle (3.5,4.8);
\draw [ultra thick, gray, semitransparent] (1.8,.3) rectangle (3.2,2.9);
\draw [ultra thick] (1,-.1) rectangle (4,6.2);
\node at (2.5,6.5) {$\varphi$};
\end{tikzpicture}
\end{center}
\end{motivation}

The operadic composition $\gamma$ \eqref{operadic-composition} is used in the following observation.

\begin{lemma}
\label{lemma:phi-presentation}
Suppose $N+r \geq 2$ in the wiring diagram $\varphi \in \WD\xprimeux$ \eqref{phi-wiring-diagram}.  Then it admits a decomposition
\begin{equation}
\label{phi-presentation}
\varphi = 
\gamma\left(\utheta;
\left\{\tensorunit_{X_i}\right\}_{i=1}^N,\, 
\left\{\delta_{d_j}\right\}_{j=1}^r 
\right).
\end{equation}
Here
\begin{equation}
\label{theta-comp-nrtwo}
\utheta = \begin{cases}
\theta_1 & \text{ if $N+r=2$},\\
\theta_1 \comptwo \Bigl( \cdots \comptwo 
\bigl(\theta_{N+r-2} \comptwo \theta_{N+r-1}\bigr)\Bigr)
& \text{ if $N+r > 2$}
\end{cases}
\end{equation}
with each $\theta_k$ a $2$-cell, and each $\delta_{d_j} \in \WD\djnothing$ is a $1$-delay node as in Assumption \ref{assumption:psi}. 
\end{lemma}

\begin{proof}
Recall that $\uX = (X_1, \ldots,X_N)$ and $\dnphi = \dnpsi = \{d_1,\ldots,d_r\}$.  For $N+1 \leq j \leq N+r$ define the box $X_j = d_j \in \boxs$ as
\[\xin_j = \{d_j\} = \xout_j.\]
For $1 \leq i \leq N+r$ define the box
\[X^{\geq i} = \coprod_{p=i}^{N+r} X_p \in \boxs.\]
Note that $X^{\geq 1} = X'$ \eqref{xprime-inout}.  Next, for $1 \leq k \leq N+r-1$, define the $2$-cell
\[\theta_k = \theta_{X_k,X^{\geq k+1}} \in \WD\xkxkxkone.\]
Then we have a wiring diagram
\[\utheta \in \WD\smallbinom{X'}{X_1,\ldots,X_{N+r}}\]
in which $\theta^*$ is defined as in \eqref{theta-comp-nrtwo}.  Since $\tensorunit_{X_i} \in \WD\xixi$ and $\delta_{d_j} \in \WD\djnothing$, the operadic composition on the right side of \eqref{phi-presentation} is defined and belongs to $\WD\xprimeux$.  

Since the two sides of \eqref{phi-presentation} both have delay nodes $\{d_1,\ldots,d_r\}$, it remains to check that the supplier assignment of the right side is equal to $\sphi = \Id$.  This follows from a direct inspection because (i) colored units (Def. \ref{wd-units}), $1$-delay nodes (Def. \ref{def:one-dn}), and $2$-cells (Def. \ref{def:theta-wd}) all have identity supplier assignments and because (ii) $\gamma$ \eqref{gamma-in-comps} is an iteration of various $\compi$ (Def. \ref{def:compi-wd}).
\end{proof}

\section{Unary Wiring Diagrams}
\label{sec:unary-wd}

In this section, we analyze wiring diagrams with  exactly one input box and no delay nodes, such as $\alpha$ \eqref{alpha-wiring-diagram}.  We will show that such a wiring diagram can be generated by the generating wiring diagrams (Def. \ref{def:generating-wiring-diagrams}) of name changes (Def. \ref{def:name-change}), $1$-loops (Def. \ref{def:loop-wd}), $1$-wasted wires (Def. \ref{def:wasted-wire-wd}), in-splits (Def. \ref{def:insplit-wd}), and out-splits (Def. \ref{def:out-split}), in this order.  We remind the reader of Notation \ref{comp-is-compone} for (iterated) $\compone$.

We will need a few definitions and notations.  

\begin{definition}
\label{def:loop-element}
Suppose $\pi \in \WD\yx$ is a wiring diagram with one input box $X$ and no delay nodes.
\begin{enumerate}
\item
A \emph{loop element} \index{loop element} in $\pi$ is an element $x \in \xout$ such that there exists $x' \in \xin$ with $x$ as its supply wire.  The set of loop elements in $\pi$ is denoted by \label{notation:loop-element}$\looppi$.
\item
An element $x' \in \xin$ is said to be \emph{internally supplied} \index{internally supplied element} if $\spi(x') \in \xout$.  The set of such elements in $\pi$ is denoted by \label{notation:int-supplied}$\piini$.
\item
An element $x' \in \xin$ is said to be \emph{externally supplied} \index{externally supplied element} if $\spi(x') \in \yin$.  The set of such elements in $\pi$ is denoted by \label{notation:ext-supplied}$\piine$.
\end{enumerate}
\end{definition}

Recall from Definition \ref{def:wiring-diagram} the concepts of external wasted wires $\ewpi$ and of internal wasted wires $\iwpi$ of a wiring diagram $\pi$.

\begin{lemma}
\label{lemma:loop-element}
Suppose $\pi \in \WD\yx$ is a wiring diagram with one input box $X$ and no delay nodes.  Then:
\begin{enumerate}
\item
$\piini \amalg \piine = \xin$.
\item
$\iwpi \amalg \looppi \subseteq \xout$.
\item
$\spi\left(\piini\right) = \looppi$.
\item
$\spi\left(\piine\right) \amalg \ewpi = \yin$.
\item
$\spi\left(\yout\right) \subseteq \xout \setminus \iwpi$.
\end{enumerate}
\end{lemma}

\begin{proof}
All the statements are immediate from the definitions.
\end{proof}

\begin{example}
\label{ex:factoring-pi}
Consider the wiring diagram $\pi \in \WD\yx$ as depicted in the picture
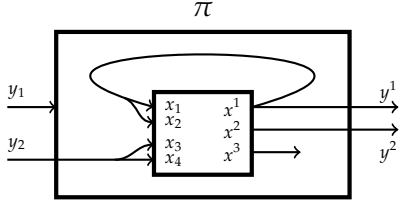
\begin{center}
\begin{tikzpicture}[xscale = 1. 3]
\draw [ultra thick] (1,1.9) rectangle (2,3);
\node at (1.2,2.8) {\tiny{$x_1$}};
\node at (1.2,2.6) {\tiny{$x_2$}};
\node at (1.2,2.3) {\tiny{$x_3$}};
\node at (1.2,2.1) {\tiny{$x_4$}};
\draw [arrow, looseness=7] (2,2.8) to [out=20, in=160] (1,2.8);
\draw [arrow] (.7,2.92) to [out=-30, in=180] (1,2.6);
\draw [arrow] (2,2.8) to (3.5,2.8);
\draw [arrow] (2,2.5) to (3.5,2.5);
\draw [arrow] (2,2.2) to (2.5,2.2);
\node at (1.8,2.8) {\tiny{$x^1$}};
\node at (1.8,2.5) {\tiny{$x^2$}};
\node at (1.8,2.2) {\tiny{$x^3$}};
\draw [ultra thick] (0,1.6) rectangle (3,3.8);
\node at (1.5,4.1) {$\pi$};
\draw [arrow] (-.5,2.8) to (0,2.8);
\draw [arrow] (-.5,2.1) to (1,2.1);
\draw [arrow] (.6,2.1) to [out=0,in=180] (1,2.3);
\node at (-.4,3) {\tiny{$y_1$}};
\node at (-.4,2.3) {\tiny{$y_2$}};
\node at (3.4,3) {\tiny{$y^1$}};
\node at (3.4,2.2) {\tiny{$y^2$}};
\end{tikzpicture}
\end{center}
with
\begin{itemize}
\item
$\xin = \{x_1,x_2,x_3,x_4\}$, $\xout = \{x^1,x^2,x^3\}$, $\yin = \{y_1,y_2\}$, and $\yout = \{y^1, y^2\}$;
\item
$\looppi = \{x^1\}$, $\piini = \{x_1,x_2\}$, $\piine = \{x_3,x_4\}$, $\iwpi = \{x^3\}$, and $\ewpi = \{y_1\}$.
\end{itemize}
Note that we may operadically decompose $\pi$ as follows.
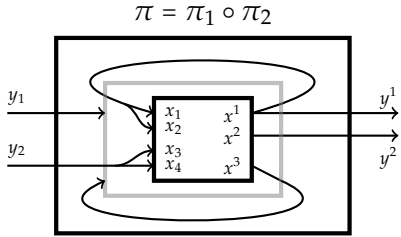
\begin{equation}
\label{pi-decomposition-picture}
\begin{tikzpicture}[xscale = 1. 3]
\draw [ultra thick] (1,1.9) rectangle (2,3);
\node at (1.2,2.8) {\tiny{$x_1$}};
\node at (1.2,2.6) {\tiny{$x_2$}};
\node at (1.2,2.3) {\tiny{$x_3$}};
\node at (1.2,2.1) {\tiny{$x_4$}};
\draw [arrow, looseness=7] (2,2.8) to [out=20, in=160] (1,2.8);
\draw [arrow] (.7,2.92) to [out=-30, in=180] (1,2.6);
\draw [arrow] (2,2.8) to (3.5,2.8);
\draw [arrow] (2,2.5) to (3.5,2.5);
\draw [thick] (2,2.1) to (2.3,1.9);
\node at (1.8,2.8) {\tiny{$x^1$}};
\node at (1.8,2.5) {\tiny{$x^2$}};
\node at (1.8,2.1) {\tiny{$x^3$}};
\draw [ultra thick, gray, semitransparent] (.5,1.7) rectangle (2.3,3.2);
\draw [arrow, looseness=2] (2.3,1.9) to [out=-30, in=210] (.5,1.9);
\draw [ultra thick] (0,1.2) rectangle (3,3.8);
\node at (1.5,4.1) {$\pi = \pi_1 \comp \pi_2$};
\draw [arrow] (-.5,2.8) to (.5,2.8);
\draw [arrow] (-.5,2.1) to (1,2.1);
\draw [arrow] (.6,2.1) to [out=0,in=180] (1,2.3);
\node at (-.4,3) {\tiny{$y_1$}};
\node at (-.4,2.3) {\tiny{$y_2$}};
\node at (3.4,3) {\tiny{$y^1$}};
\node at (3.4,2.2) {\tiny{$y^2$}};
\end{tikzpicture}
\end{equation}
The point of this decomposition is that the inner wiring diagram $\pi_2$ is generated by:
\begin{itemize}
\item
two $1$-wasted wires, one for the external wasted wire $y_1$ and the other for the internal wasted wire $x^3$;
\item
two in-splits, one for $\{x_1,x_2\}$ and the other for $\{x_3,x_4\}$;
\item
one out-split for $x^1$, which is the supply wire of $y^1$, $x_1$, and $x_2$.  
\end{itemize}
At the same time, the outer wiring diagram $\pi_1$ is generated by two $1$-loops, one for the  loop element $x^1$ and the other for the internal wasted wire $x^3$.

With this example as a guide, next we will factor a general wiring diagram with one input box and no delay nodes into two wiring diagrams.  The outer wiring diagram will be generated by name changes and $1$-loops.  The inner wiring diagram will be generated by $1$-wasted wires, in-splits, and out-splits.  The intermediate gray box in \eqref{pi-decomposition-picture} will be called $Z$ below.
\end{example}

\begin{convention}
\label{conv:ignore-name-change}
Using the five elementary relations \eqref{move:a2}, \eqref{move:a3}, \eqref{move:a4}, \eqref{move:a5}, and \eqref{move:a6}, name changes can always be rewritten on the outside (i.e., left side) of an iterated operadic composition in $\WD$.  Moreover, using the elementary relation \eqref{move:a1}, an iteration of name changes can be composed down into just one name change.  To simplify the presentation, in what follows these elementary relations regarding name changes are automatically applied wherever necessary.  With this in mind, we will mostly not mention name changes.
\end{convention}

For a wiring diagram with one input box and no delay nodes, we will decompose it using the wiring diagrams in the next definition.

\begin{definition}
\label{def:pione-pitwo}
Suppose $\pi \in \WD\yx$ is a wiring diagram with one input box $X$ and no delay nodes.
\begin{enumerate}
\item
Define the box $Z \in \boxs$ as
\[
\begin{split}
\zin &= \yin \amalg \iwpi \amalg \looppi;\\
\zout &= \yout \amalg \iwpi \amalg \looppi.
\end{split}
\]
\item
Define the wiring diagram $\pi_1 \in \WD\yz$ with:
\begin{itemize}
\item
one input box $Z$, output box $Y$, and no delay nodes;
\item
supplier assignment
\begin{equation}
\label{pi-one-supplier}
\nicexy{
\dmpione = \yout \amalg \zin 
= \yout \amalg  \left[\yin \amalg \iwpi \amalg \looppi\right]
\ar[d]_-{\spione}
\\
\supplypione = \yin \amalg \zout
= \yin \amalg  \left[\yout \amalg \iwpi \amalg \looppi\right]
}
\end{equation}
the identity function.
\end{itemize}
\item
Define the wiring diagram $\pi_2 \in \WD\zx$ with:
\begin{itemize}
\item
one input box $X$, output $Z$, and no delay nodes;
\item
supplier assignment
\begin{equation}
\label{pi-two-supplier}
\nicexy{
\dmpitwo = \zout \amalg \xin
=  \left[\yout \amalg \iwpi \amalg \looppi\right] 
\amalg \left[\piini \amalg \piine\right]
\ar[d]_-{\spitwo}
\\
\supplypitwo = \zin \amalg \xout
= \left[\yin \amalg \iwpi \amalg \looppi\right] 
\amalg \xout
}
\end{equation}
whose restriction to:
\begin{itemize}
\item
$\yout$ is $\spi : \yout \to \xout$;
\item
$\iwpi \amalg \looppi$ is the subset inclusion into $\xout$;
\item
$\piini$ is $\spi : \piini \to \looppi$;
\item
$\piine$ is $\spi : \piine \to \yin$.
\end{itemize}
This is well-defined by the non-instantaneity requirement \eqref{non-instant} for $\pi$ and Lemma \ref{lemma:loop-element}.
\end{itemize}
\end{enumerate}
\end{definition}

An example of the following decomposition is the picture \eqref{pi-decomposition-picture} above.

\begin{lemma}
\label{lemma:pi-decomposition}
Suppose $\pi \in \WD\yx$ is a wiring diagram with one input box $X$ and no delay nodes.  Then it admits a decomposition\index{decomposition of unary wiring diagrams}
\begin{equation}
\label{pi-pione-pitwo}
\pi = \pi_1 \comp \pi_2
\end{equation}
in which $\pi_1 \in \WD\yz$ and $\pi_2 \in \WD\zx$ are as in Def. \ref{def:pione-pitwo}.
\end{lemma}

\begin{proof}
Both sides of \eqref{pi-pione-pitwo} belong to $\WD\yx$ and have no delay nodes.  So it remains to check that the  supplier assignment  $s$ of $\pi_1 \comp \pi_2$ is equal to $\spi$.  Note that
\[
\dm_{\pi_1 \comp \pi_2} = \yout \amalg \xin = \yout \amalg  \left[\piini \amalg \piine\right].   
\]
By the definitions of $\comp = \compone$ (Def. \ref{def:compi-wd}), $\spione$ \eqref{pi-one-supplier}, and $\spitwo$ \eqref{pi-two-supplier}:
\begin{itemize}
\item
on $\yout$ the supplier assignment $s$ is $\spitwo \spione = \spi \Id = \spi$.
\item
on $\piini \amalg \piine$ the supplier assignment $s$ is $\spitwo\spione = \Id \spi = \spi$.
\end{itemize}
So the supplier assignment of $\pi_1 \comp \pi_2$ is equal to $\spi$.  
\end{proof}

To obtain the desired stratified presentation of $\pi$, next we observe that $\pi_1$ in Def. \ref{def:pione-pitwo} is either a colored unit \eqref{wd-unit} or an iterated operadic composition of $1$-loops (Def. \ref{def:loop-wd}).  An example of $\pi_1$ is the outer wiring diagram in the example \eqref{pi-decomposition-picture}.

\begin{lemma}
\label{lemma:iterated-loops}
Suppose:
\begin{itemize}
\item
$Y,Z \in \boxs$ such that $\zin = \yin \amalg T$ and $\zout = \yout \amalg T$ for some $T \in \Fins$.
\item
$\zeta \in \WD\yz$ is a wiring diagram with no delay nodes and with supplier assignment
\[
\nicexy{
\dmzeta = \yout \amalg \zin = \yout \amalg \left[\yin \amalg T\right]
\ar[d]_-{\szeta} \\
\supplyzeta = \yin \amalg \zout = \yin \amalg \left[\yout \amalg T\right]
}\]
the identity function.
\end{itemize}
Then the following statements hold.
\begin{enumerate}
\item
$\zeta = \tensorunit_{Y}$ if $T = \varnothing$.
\item
If  $p = |T| > 0$, then there exist $1$-loops $\lambda_1,\ldots, \lambda_p$ such that
\begin{equation}
\label{iterated-oneloops}
\zeta = \lambda_1 \comp \cdots \comp \lambda_p.
\end{equation}
\end{enumerate}
\end{lemma}

\begin{proof}
If $T = \varnothing$, then $\zeta \in \WD\yy$ has no delay nodes and has supplier assignment $\szeta = \Id$.  So $\zeta$ is the $Y$-colored unit.

Next suppose $T = \{t_1,\ldots,t_p\}$ with $p > 0$.  For the definitions below, it is convenient to keep in mind the following picture of $\zeta$:
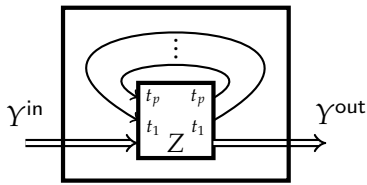
\begin{center}
\begin{tikzpicture}
\draw [ultra thick] (2,1) rectangle (3,2);
\node at (2.5,1.2) {\small{$Z$}};
\node at (2.2,1.8) {\tiny{$t_p$}};
\node at (2.8,1.8) {\tiny{$t_p$}};
\node at (2.2,1.4) {\tiny{$t_1$}};
\node at (2.8,1.4) {\tiny{$t_1$}};
\draw [implies] (.5,1.2) to (2,1.2);
\node at (.5,1.6) {$\yin$};
\draw [implies] (3,1.2) to (4.5,1.2);
\node at (4.7,1.6) {$\yout$};
\draw [arrow, looseness=3] (3,1.8) to [out=30, in=150] (2,1.8);
\draw [arrow, looseness=8] (3,1.5) to [out=30, in=150] (2,1.5);
\node at (2.5,2.45) {$\vdots$};
\draw [ultra thick] (1,0.7) rectangle (4,3);
\end{tikzpicture}
\end{center}
For $0 \leq j \leq p$ define the box $Y_j \in \boxs$ as
\[
Y_j = \begin{cases}
Y & \text{ if $j=0$};\\
Y \amalg \{t_1, \ldots, t_j\} & \text{ if $1 \leq j \leq p$}.
\end{cases}
\]
Here $Y \amalg \{t_1, \ldots, t_j\}$ means a copy of $t_i$ for $1 \leq i \leq j$ is added to each of $\yin$ and $\yout$.  In particular, we have $Y_p = Z$.   For $1 \leq j \leq p$ define the $1$-loop
\[
\lambda_j = \lambda_{Y_j, t_j} \in \WD\yjminusoneyj
\]
in which $t_j \in \yin_j$ has supply wire $t_j \in \yout_j$.  

The iterated operadic composition
\[
\lambda_1 \comp \cdots \comp \lambda_p \in \WD\yzeroyp = \WD\yz
\]
has no delay nodes.  To see that it is equal to $\zeta$, it remains to check that its supplier assignment is equal to $\szeta = \Id$.  This holds because each $1$-loop $\lambda_j$ has identity supplier assignment.
\end{proof}

Observe that Lemma \ref{lemma:iterated-loops} applies to $\pi_1 \in \WD\yz$ in Def.  \ref{def:pione-pitwo} with $T = \iwpi \amalg \looppi$.  So $\pi_1$ is either the $Y$-colored unit or an iterated operadic composition of $1$-loops.

\section{Unary Wiring Diagrams with No Loop Elements}
\label{sec:unary-no-loop}

In order to show that the wiring diagram $\pi_2 \in \WD\zx$ in Def.  \ref{def:pione-pitwo} is generated by $1$-wasted wires, in-splits, and out-splits, first we identify its external wasted wires, internal wasted wires, and loop elements.

\begin{lemma}
\label{lemma:pitwo-wasted-wired}
Consider the wiring diagram $\pi_2 \in \WD\zx$ in Def.  \ref{def:pione-pitwo}.
\begin{enumerate}
\item
The set of external wasted wires in $\pi_2$ is $\ewpitwo = \iwpi \amalg \ewpi$.
\item
The set of internal wasted wires in $\pi_2$ is $\iwpitwo = \varnothing$.
\item
The set of loop elements in $\pi_2$ is $\looppitwo = \varnothing$.
\item
$\spitwo(\zout) = \xout$.
\item
$\zin =  \iwpi \amalg \ewpi \amalg \spitwo(\xin)$.
\end{enumerate}
\end{lemma}

\begin{proof}
\begin{enumerate}
\item
By definition an external wasted wire in $\pi_2$ is an element in $\zin$ that is not in the image of $\spitwo$ \eqref{pi-two-supplier}.  By the definition of $\spitwo$, this is the subset
\[
\iwpi  \amalg [\yin \setminus \spi(\piine)] \subseteq \zin.  
\]
It follows from the non-instantaneity requirement \eqref{non-instant} for $\pi$ that 
\[
\yin \setminus \spi(\piine) = \ewpi.
\]
\item
By definition an internal wasted wire in $\pi_2$ is an element in $\xout$ that is not in the image of $\spitwo$.  An element of $\xout$ is either an internal wasted wire in $\pi$, or else it is the $\spi$-image of an element in $\piini \amalg \yout$.  Since
\[
\begin{split}
\xout 
&= \iwpi \amalg \spi\left(\piini \amalg \yout\right)\\
&= \iwpi \amalg \left[\looppi \cup \spi(\yout)\right],
\end{split}
\]
an inspection of the definition of $\spitwo$ \eqref{pi-two-supplier} reveals that all of $\xout$ is in the image of $\spitwo$.  So $\pi_2$ has no internal wasted wires.
\item
By definition a loop element in $\pi_2$ is an element in $\xout$ that is the supply wire, under $\spitwo$, of some element in $\xin$.  Since $\xin = \piini \amalg \piine$, the definition of $\spitwo$  \eqref{pi-two-supplier}  yields
\[
\spitwo(\piini \amalg \piine)
= \spi(\piini) \amalg \spi(\piine)
\subseteq \looppi \amalg \yin 
\subseteq \zin
= \supplypitwo \setminus \xout.
\]
So $\pi_2$ has no loop elements.
\item
By (2) $\pi_2$ has no internal wasted wires, so $\xout = \spitwo(\zout \amalg \pitwoini)$.  But by (3) $\pi_2$ has no loop elements, so $\pitwoini = \varnothing$ and $\xout = \spitwo(\zout)$.
\item
Since $\pi_2$ has no loop elements by (3), $\spitwo(\xin) \subseteq \zin$.  An element in $\zin$ that is not in $\spitwo(\xin)$ is precisely an external wasted wire in $\pi_2$.  By (1) the set of external wasted wires in $\pi_2$ is $\iwpi \amalg \ewpi$.
\end{enumerate}
\end{proof}

Continuing our analysis of wiring diagrams with one input box and no delay nodes, our next goal is to construct a decomposition for $\pi_2$ (Def.  \ref{def:pione-pitwo}) involving $1$-wasted wires, in-splits, and out-splits.

\begin{example}
\label{ex:factor-pitwo}
Consider the inner wiring diagram $\pi_2 \in \WD\zx$ in the example \eqref{pi-decomposition-picture}, which is depicted in the following picture.
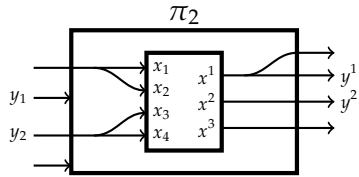
\begin{center}
\begin{tikzpicture}[xscale=1]
\draw [ultra thick] (1,1.8) rectangle (2,3.1);
\node at (1.2,2.9) {\tiny{$x_1$}};
\node at (1.2,2.6) {\tiny{$x_2$}};
\node at (1.2,2.3) {\tiny{$x_3$}};
\node at (1.2,2) {\tiny{$x_4$}};
\draw [arrow] (-.5,2.9) to (1,2.9);
\draw [arrow] (-.5,2) to (1,2);
\draw [arrow] (2,2.8) to (3.5,2.8);
\draw [arrow] (2,2.45) to (3.5,2.45);
\draw [arrow] (2,2.1) to (3.5,2.1);
\node at (1.8,2.8) {\tiny{$x^1$}};
\node at (1.8,2.45) {\tiny{$x^2$}};
\node at (1.8,2.1) {\tiny{$x^3$}};
\draw [ultra thick] (0,1.5) rectangle (3,3.4);
\node at (1.5,3.6) {$\pi_2$};
\node at (3.7,2.8) {\tiny{$y^1$}};
\node at (3.7,2.45) {\tiny{$y^2$}};
\draw [arrow] (-.5,2.5) to (0,2.5);
\draw [arrow] (-.5,1.6) to (0,1.6);
\node at (-.7,2.5) {\tiny{$y_1$}};
\node at (-.7,2) {\tiny{$y_2$}};
\draw [thick] (2.3,2.8) to [out=0, in=180] (3,3.1);
\draw [arrow] (3,3.1) to (3.5,3.1);
\draw [arrow] (.3,2.9) to [out=0, in=180] (1,2.6);
\draw [arrow] (.3,2) to [out=0, in=180] (1,2.3);
\end{tikzpicture}
\end{center}
For this wiring diagram, the desired decomposition is depicted in the picture:
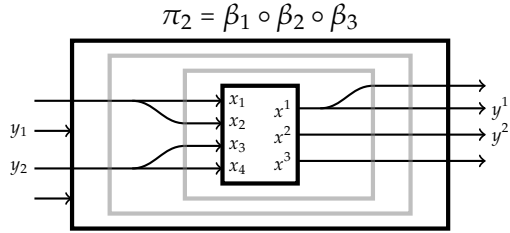
\begin{center}
\begin{tikzpicture}[xscale=1]
\draw [ultra thick] (1,1.8) rectangle (2,3.1);
\node at (1.2,2.9) {\tiny{$x_1$}};
\node at (1.2,2.6) {\tiny{$x_2$}};
\node at (1.2,2.3) {\tiny{$x_3$}};
\node at (1.2,2) {\tiny{$x_4$}};
\draw [arrow] (-1.5,2.9) to (1,2.9);
\draw [arrow] (-1.5,2) to (1,2);
\draw [arrow] (2,2.8) to (4.5,2.8);
\draw [arrow] (2,2.45) to (4.5,2.45);
\draw [arrow] (2,2.1) to (4.5,2.1);
\node at (1.8,2.8) {\tiny{$x^1$}};
\node at (1.8,2.45) {\tiny{$x^2$}};
\node at (1.8,2.1) {\tiny{$x^3$}};
\draw [arrow] (-1.5,2.5) to (-1,2.5);
\draw [arrow] (-1.5,1.6) to (-1,1.6);
\node at (-1.7,2.5) {\tiny{$y_1$}};
\node at (-1.7,2) {\tiny{$y_2$}};
\draw [thick] (2.3,2.8) to [out=0, in=180] (3,3.1);
\draw [arrow] (3,3.1) to (4.5,3.1);
\draw [thick] (-.2,2.9) to [out=0, in=180] (.5,2.6);
\draw [arrow] (.5,2.6) to (1,2.6);
\draw [thick] (-.2,2) to [out=0, in=180] (.5,2.3);
\draw [arrow] (.5,2.3) to (1,2.3);
\draw [ultra thick, gray, semitransparent] (.5,1.6) rectangle (3,3.3);
\draw [ultra thick, gray, semitransparent]  (-.5,1.4) rectangle (3.5,3.5);
\draw [ultra thick] (-1,1.2) rectangle (4,3.7);
\node at (1.5,4) {$\pi_2 = \beta_1 \comp \beta_2 \comp \beta_3$};
\node at (4.7,2.8) {\tiny{$y^1$}};
\node at (4.7,2.45) {\tiny{$y^2$}};
\end{tikzpicture}
\end{center}
The inner gray box will be called $V$, and the outer gray box will be called $W$ below.  Note that:
\begin{itemize}
\item
The outermost wiring diagram $\beta_1 \in \WD\zw$ is generated by two $1$-wasted wires.  
\item
The middle wiring diagram $\beta_2 \in \WD\wv$ is generated by two in-splits.  
\item
The innermost wiring diagram $\beta_3 \in \WD\vx$ is an out-split.  
\end{itemize}
For a general wiring diagram with one input box, no delay nodes, no loop elements, and no internal wasted wires, such a decomposition uses the following definitions.
\end{example}

\begin{definition}
\label{def:beta}
Suppose $X,Z \in \boxs$ and $\beta \in \WD\zx$ is a wiring diagram with no delay nodes and no loop elements.
\begin{enumerate}
\item
Suppose the box $W = Z \setminus \ewbeta \in \boxs$ is obtained from $Z \in \boxs$ by removing the external wasted wires of $\beta$, so $\win = \zin \setminus \ewbeta$ and $\wout = \zout$.
\item
Define the wiring diagram $\beta_1 \in \WD\zw$ as having:
\begin{itemize}
\item
no delay nodes;
\item
supplier assignment
\begin{equation}
\label{supply-betaone}
\nicexy{
\dmbetaone =
\zout \amalg \win = \zout \amalg \left[\zin \setminus  \ewbeta\right]
\ar[d]_-{\sbetaone}
\\
\supplybetaone = \zin \amalg \wout = \zin \amalg \zout
}
\end{equation}
the identity function on $\zout$ and the subset inclusion on $\zin \setminus  \ewbeta \subseteq \zin$.
\end{itemize}
\item
Define the box $V \in \boxs$ as $\vin = \xin$ and $\vout = \zout = \wout$.
\item
Define the wiring diagram $\beta_2 \in \WD\wv$ as having:
\begin{itemize}
\item
no delay nodes;
\item
supplier assignment
\begin{equation}
\label{supply-betatwo}
\nicexy{
\dmbetatwo = \wout \amalg \vin = \wout \amalg \xin
\ar[d]_-{\sbetatwo}
\\
\supplybetatwo = \win \amalg \vout = \left[\zin \setminus \ewbeta\right] \amalg \wout
}
\end{equation}
the coproduct of the identity function on $\wout$ and the restriction of the supplier assignment $\sbeta : \xin \to \zin \setminus \ewbeta$.
\end{itemize}
This is well-defined because $\beta$ has no delay nodes and no loop elements.
\item
Define the wiring diagram $\beta_3 \in \WD\vx$ as having:
\begin{itemize}
\item
no delay nodes;
\item
supplier assignment
\begin{equation}
\label{supply-betathree}
\nicexy{
\dmbetathree = \vout \amalg \xin = \zout \amalg \xin
\ar[d]_-{\sbetathree}
\\
\supplybetathree = \vin \amalg \xout = \xin \amalg \xout
}
\end{equation}
the coproduct of the identity function on $\xin$ and $\sbeta : \zout \to \xout$.
\end{itemize}
This is well-defined because $\beta$ has no delay nodes and because of the non-instantaneity requirement \eqref{non-instant} for $\beta$.
\end{enumerate}
\end{definition}

\begin{lemma}
\label{lemma:sbeta-surjective}
In the context of Def. \ref{def:beta}:
\begin{enumerate}
\item
The map  $\sbeta : \xin \to \zin \setminus \ewbeta$, which is part of $\sbetatwo$, is surjective.
\item
If $\beta$ has no internal wasted wires (such as $\pi_2$  in Def.  \ref{def:pione-pitwo}), then the map  $\sbeta : \zout \to \xout$, which is part of $\sbetathree$, is surjective.
\end{enumerate}
\end{lemma}

\begin{proof}
The first assertion is true because $\beta$ has no delay nodes and because of the non-instantaneity requirement \eqref{non-instant}.  The second assertion is true  because $\beta$ has no delay nodes and no loop elements.
\end{proof}

\begin{lemma}
\label{lemma:beta-factor}
In the context of Def. \ref{def:beta}, there is a decomposition\index{decomposition of unary wiring diagrams}
\begin{equation}
\label{betaonetwothree}
\beta = \beta_1 \comp \beta_2 \comp \beta_3 \in \WD\zx.
\end{equation}
\end{lemma}

\begin{proof}
By construction the iterated operadic composition $\beta_1 \comp \beta_2 \comp \beta_3$ also belongs to $\WD\zx$ and has no delay nodes.  So it remains to check that its supplier assignment $s$ is equal to $\sbeta$.  A direct inspection of \eqref{supply-betaone}, \eqref{supply-betatwo}, and \eqref{supply-betathree} reveals that:
\begin{itemize}
\item
on $\xin \subseteq \dmbetathree$ the supplier assignment $s$ is given by $\Id\sbetatwo\Id_{\xin} = \sbeta$;
\item
on $\zout \subseteq \dmbetaone$ the supplier assignment $s$ is given by $\sbetathree\Id_{\zout}\Id_{\zout} = \sbeta$.
\end{itemize}
So the supplier assignment of  $\beta_1 \comp \beta_2 \comp \beta_3$ is equal to $\sbeta$.
\end{proof}

Note that the decomposition in Lemma \ref{lemma:beta-factor} applies to $\pi_2$ because $\pi_2$ has one input box, no delay nodes, and no loop elements (by Lemma \ref{lemma:pitwo-wasted-wired}).

Next we show that in the decomposition \ref{betaonetwothree}:
\begin{enumerate}
\item
$\beta_1$ is either a colored unit (Def. \ref{wd-units}) or an iterated operadic composition of $1$-wasted wires  (Def. \ref{def:wasted-wire-wd}).  See Lemma \ref{lemma:betaone}.
\item
$\beta_2$ is either a colored unit or an iterated operadic composition of in-splits (Def. \ref{def:insplit-wd}).  See Lemma \ref{lemma:betatwo}.
\item
If $\beta$ has no internal wasted wires, such as $\pi_2$ by Lemma \ref{lemma:pitwo-wasted-wired}, then $\beta_3$ is either a colored unit or an iterated operadic composition of out-splits (Def. \ref{def:out-split}).  See Lemma \ref{lemma:betathree}.
\end{enumerate}
During the first reading, the reader may wish to skip the proofs of the following Lemmas and simply look at the pictures.  The following observation deals with the first statement.

\begin{lemma}
\label{lemma:betaone}
Consider the wiring diagram $\beta_1 \in \WD\zw$ in Def. \ref{def:beta}.
\begin{enumerate}
\item
If $\ewbeta = \varnothing$ (i.e., $\beta$ has no external wasted wires), then $\beta_1 = \tensorunit_Z$, the $Z$-colored unit.
\item
If $q = |\ewbeta| > 0$, then there exist $1$-wasted wires $\omega_1, \ldots, \omega_q$ such that
\begin{equation}
\label{betaone-iterated-omega}
\beta_1 = \omega_1 \comp \cdots \comp \omega_q.
\end{equation}
\end{enumerate}
\end{lemma}

\begin{proof}
Recall that $\beta_1$ has no delay nodes and has supplier assignment \eqref{supply-betaone}
\[
\nicexy{
\dmbetaone =
\zout \amalg \win = \zout \amalg \left[\zin \setminus  \ewbeta\right]
\ar[d]_-{\sbetaone}
\\
\supplybetaone = \zin \amalg \wout = \zin \amalg \zout
}\]
that is the identity function on $\zout$ and the subset inclusion on $\zin \setminus \ewbeta$.  If $\ewbeta = \varnothing$, then $\sbetaone = \Id$ and, therefore, $\beta_1$ is the colored unit.

Next suppose $\ewbeta = \{w_1, \ldots, w_q\} \subseteq \zin$ with $q > 0$.  Recall that $W = Z \setminus \ewbeta$.  For $0 \leq j \leq q$ define the box
\[
Z_j = \begin{cases}
Z & \text{ if $j=0$};\\
Z \setminus \{w_1,\ldots,w_j\} & \text{ if $1 \leq j \leq q$}.
\end{cases}
\]
So in particular $Z_q = W$.  The iterated operadic composition on the right side of \eqref{betaone-iterated-omega} is represented in the following picture.
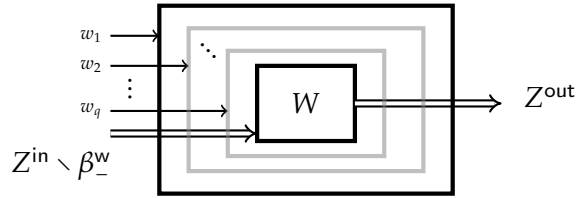
\begin{center}
\begin{tikzpicture}[xscale=1.3]
\draw [ultra thick] (2,1) rectangle (3,2);
\node at (2.5,1.5) {$W$};
\draw [implies] (.5,1.1) to (2,1.1);
\node at (0,.7) {$\zin \setminus \ewbeta$};
\draw [arrow] (.5,2.4) to (1,2.4);
\node at (.3,2.4) {\tiny{$w_1$}};
\draw [arrow] (.5,2) to (1.3,2);
\node at (.3,2) {\tiny{$w_2$}};
\node at (.7,1.7) {$\vdots$};
\draw [arrow] (.5,1.4) to (1.7,1.4);
\node at (.3,1.4) {\tiny{$w_q$}};
\draw [implies] (3,1.5) to (4.5,1.5);
\node at (5,1.6) {$\zout$};
\draw [ultra thick, gray, semitransparent] (1.3,0.6) rectangle (3.7,2.5);
\node at (1.5,2.2) {$\ddots$};
\draw [ultra thick, gray, semitransparent] (1.7,0.8) rectangle (3.3,2.2);
\draw [ultra thick] (1,0.3) rectangle (4,2.8);
\end{tikzpicture}
\end{center}
Here the outermost box is $Z$, the outermost gray box is $Z_1$, and the innermost gray box is $Z_{q-1}$.

For $1 \leq j \leq q$ define the $1$-wasted wires  (Def. \ref{def:wasted-wire-wd})
\[
\omega_j = \omega_{Z_{j-1},w_j} \in \WD\zjminusonezj.
\]
The iterated operadic composition
\[
\omega_1 \comp \cdots \comp \omega_q \in \WD\zzerozq = \WD\zw
\]
has no delay nodes.  So to prove \eqref{betaone-iterated-omega}, it remains to check that its supplier assignment $s$ is equal to $\sbetaone$.
\begin{itemize}
\item
On $\zout \subseteq \dm_{\omega_1}$ the supplier assignment $s$ is the composition of $q$ copies of the identity function, hence is the identity function.  
\item
On $\win \subseteq \dm_{\omega_q}$ the supplier assignment $s$ is the composition of the inclusions $\zin_j \to \zin_{j-1}$ for $1 \leq j \leq q$, which is the inclusion $\win \to \zin$.
\end{itemize}
\end{proof}

\begin{motivation}
To show that $\beta_2$ is either a colored unit or an iterated composition of in-splits, we first need a lemma that says that the following wiring diagram is generated by in-splits.
\begin{center}
\begin{tikzpicture}
\draw [ultra thick] (2,1) rectangle (3,2.2);
\node at (2.7,1.5) {$X$};
\node at (2.2,2) {\tiny{$x_1$}};
\node at (2.2,1.7) {$\vdots$};
\node at (2.2,1.35) {\tiny{$x_k$}};
\draw [implies] (.5,1.1) to (2,1.1);
\node at (0,1.6) {$\yin$};
\node at (.7,1.9) {\tiny{$y$}};
\draw [thick] (.5,1.7) -- (1.2,1.7);
\draw [arrow] (1.2,1.7) to [out=0, in=180] (2,2);
\draw [arrow] (1.2,1.7) to [out=0, in=180] (2,1.4);
\node at (1.8,1.7) {$\vdots$};
\draw [implies] (3,1.5) to (4.5,1.5);
\node at (5,1.6) {$\yout$};
\draw [ultra thick] (1,0.7) rectangle (4,2.4);
\end{tikzpicture}
\end{center}
\end{motivation}

\begin{lemma}
\label{lemma:iterated-insplits}
Suppose:
\begin{itemize}
\item
$X, Y \in \boxs$ such that $\xout = \yout$.
\item
There exist $y \in \yin$ and distinct elements $x_1,\ldots,x_k \in \xin$ with $k \geq 1$ such that
\[
\xin = \left[\yin \setminus \{y\}\right] \amalg \{x_1,\ldots,x_k\}
\]
and $v(y) = v(x_i) \in S$ for all $i$.
\item
$\sigma \in \WD\yx$ is a wiring diagram with no delay nodes and with supplier assignment
\[
\nicexy{
\dmsigma = \yout \amalg \xin = \yout \amalg \left[\yin \setminus \{y\}\right] \amalg \{x_1,\ldots,x_k\}
\ar[d]_-{\ssigma}
\\
\supplysigma = \yin \amalg \xout = \yin \amalg \yout
}\]
given by
\[
\ssigma(z) = \begin{cases}
y & \text{ if $z = x_1,\ldots,x_k$};\\
z & \text{ if $z \in \yout \amalg \left[\yin \setminus \{y\}\right]$}.
\end{cases}
\]
\end{itemize}
Then:
\begin{enumerate}
\item
$\sigma$ is the $Y$-colored unit if $k=1$;
\item
$\sigma$ is an iterated operadic composition of $(k-1)$ in-splits if $k \geq 2$.
\end{enumerate}
\end{lemma}

\begin{proof}
Since $\sigma$ has no delay nodes, if $k=1$, then $\ssigma = \Id$.  So $\sigma$ is a colored unit.

Suppose $k \geq 2$.  We will prove that $\sigma$ is an iterated operadic composition of $(k-1)$ in-splits by induction on $k$.  If $k=2$, then by definition $\sigma$ is the in-split $\sigma_{X,x_1,x_2}$  (Def. \ref{def:insplit-wd}).  

Suppose $k \geq 3$.  We will factor $\sigma$ into two wiring diagrams as depicted in the picture
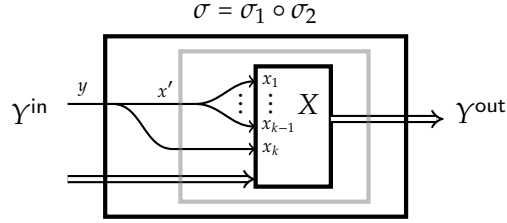
\begin{center}
\begin{tikzpicture}
\draw [ultra thick] (2,.6) rectangle (3,2.2);
\node at (2.7,1.7) {$X$};
\node at (2.2,2) {\tiny{$x_1$}};
\node at (2.2,1.7) {$\vdots$};
\node at (2.3,1.4) {\tiny{$x_{k-1}$}};
\node at (2.2,1.1) {\tiny{$x_{k}$}};
\draw [implies] (-.5,.7) to (2,.7);
\node at (-1,1.6) {$\yin$};
\node at (-.3,1.9) {\tiny{$y$}};
\draw [thick] (-.5,1.7) -- (1.2,1.7);
\draw [thick] (.1,1.7) to [out=0, in=180] (.9,1.1);
\draw [arrow] (.9,1.1) to (2,1.1);
\node at (.8,1.9) {\tiny{$x'$}};
\draw [arrow] (1.2,1.7) to [out=0, in=180] (2,2);
\draw [arrow] (1.2,1.7) to [out=0, in=180] (2,1.4);
\node at (1.8,1.7) {$\vdots$};
\draw [implies] (3,1.5) to (4.5,1.5);
\node at (5,1.6) {$\yout$};
\draw [ultra thick, gray, semitransparent] (1,0.4) rectangle (3.5,2.4);
\draw [ultra thick] (0,.2) rectangle (4,2.6);
\node at (2,2.9) {$\sigma = \sigma_1 \comp \sigma_2$};
\end{tikzpicture}
\end{center}
in which the intermediate gray box will be called $W$ below.  The outer wiring diagram $\sigma_1$ will be an in-split, and the inner wiring diagram $\sigma_2$ will be an iterated operadic composition of $(k-2)$ in-splits.  To define such a decomposition, we will need the following definitions.
\begin{enumerate}
\item
Suppose $W \in \boxs$ such that $\wout = \yout = \xout$ and
\[
\win = \left[\yin \setminus \{y\}\right] \amalg \{x',x_k\}
\]
for some $x_k \not= x'$ such that $v(x') = v(x_k) \in S$.  In particular, we have
\begin{equation}
\label{xin-win}
\xin = \left[\win \setminus \{x'\}\right] \amalg \{x_1,\ldots,x_{k-1}\}.
\end{equation}
\item
Define the wiring diagram $\sigma_1 \in \WD\yw$ with no delay nodes and with supplier assignment
\[
\nicexy{
\dmsigmaone = \yout \amalg \win 
= \yout \amalg  \left[\yin \setminus \{y\}\right] \amalg \{x',x_k\}
\ar[d]_-{\ssigmaone}
\\
\supplysigmaone = \yin \amalg \wout 
= \left[\yin \setminus \{y\}\right] \amalg \{y\} \amalg \yout
}\]
given by
\[
\ssigmaone(z) = \begin{cases}
y & \text{ if $z = x',x_k$};\\
z & \text{ otherwise}.
\end{cases}
\]
\item
Define the wiring diagram $\sigma_2 \in \WD\wx$ with no delay nodes and with supplier assignment
\[
\nicexy{
\dmsigmatwo = \wout \amalg \xin 
= \wout \amalg \left[\yin \setminus \{y\}\right] \amalg \{x_1,\ldots,x_{k-1}\} \amalg \{x_k\} \ar[d]_-{\ssigmatwo}\\
\supplysigmatwo = \win \amalg \xout 
=  \left[\yin \setminus \{y\}\right] \amalg \{x',x_k\} \amalg \wout
}\]
given by
\[
\ssigmatwo(z) = \begin{cases}
x' & \text{ if $z = x_1, \ldots, x_{k-1}$};\\
z & \text{ otherwise}.
\end{cases}
\]
\end{enumerate}
Then $\sigma_1 \comp \sigma_2 \in \WD\yx$ is a wiring diagram with no delay nodes.  To see that it is equal to $\sigma$, it suffices to check that the supplier assignment of $\sigma_1 \comp \sigma_2$ is equal to $\ssigma$.  This is true by a direct inspection of $\ssigmaone$ and $\ssigmatwo$.

By definition $\sigma_1$ is the in-split $\sigma_{W,x',x_k}$.  By \eqref{xin-win} the induction hypothesis applies to $\sigma_2$, which says that it is an iterated operadic composition of $(k-2)$ in-splits.  Combined with the previous paragraph, it follows that $\sigma = \sigma_1 \comp \sigma_2$ is the iterated operadic composition of $(k-1)$ in-splits, finishing the induction.
\end{proof}

Next we consider $\beta_2$.

\begin{lemma}
\label{lemma:betatwo}
The wiring diagram $\beta_2 \in \WD\wv$ in Def. \ref{def:beta} is either a colored unit or an iterated operadic composition of in-splits.
\end{lemma}

\begin{proof}
Recall that $\win = \zin \setminus \ewbeta$, $\vin = \xin$, and $\vout = \zout = \wout$.  The wiring diagram $\beta_2 \in \WD\wv$ has no delay nodes and has supplier assignment \eqref{supply-betatwo}
\[
\nicexy{
\dmbetatwo = \wout \amalg \vin = \zout \amalg \xin
\ar[d]_-{\sbetatwo}
\\
\supplybetatwo = \win \amalg \vout = \left[\zin \setminus \ewbeta\right] \amalg \zout
}
\]
the coproduct of the identity function on $\zout$ and $\sbeta : \xin \to \zin \setminus \ewbeta$.  Write $\win = \zin \setminus \ewbeta = \{z_1, \ldots, z_p\}$, so each $\sbeta^{-1}(z_i)$ is non-empty by Lemma \ref{lemma:sbeta-surjective}.

If $p=0$, then $\zin \setminus \ewbeta = \varnothing = \xin$, and $\beta_2$ is the $W$-colored unit.

Suppose $p > 0$.  We will write $\beta_2$ as an iterated composition as in the following picture.
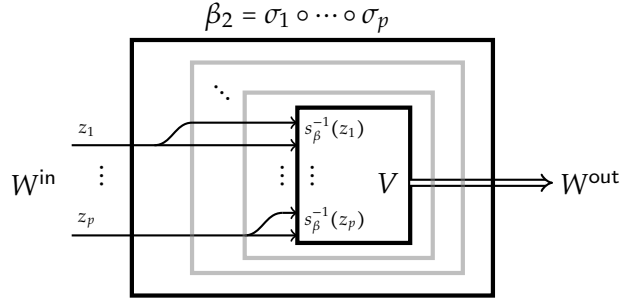
\begin{center}
\begin{tikzpicture}
\draw [ultra thick] (2,.4) rectangle (3.5,2.2);
\node at (3.2,1.2) {$V$};
\node at (2.5,1.9) {\tiny{$\sbeta^{-1}(z_1)$}};
\node at (2.2,1.3) {$\vdots$};
\node at (2.5,.7) {\tiny{$\sbeta^{-1}(z_p)$}};
\draw [arrow] (-1,.5) to (2,.5);
\draw [thick] (1.3,.5) to  [out=0, in=180] (1.8,.8);
\draw [arrow] (1.8,.8) to (2,.8);
\draw [implies] (3.5,1.2) to (5.4,1.2);
\node at (5.9,1.2) {$\wout$};
\draw [ultra thick, gray, semitransparent] (1.3,.2) rectangle (3.8,2.4);
\node at (1.8,1.3) {$\vdots$};
\draw [ultra thick, gray, semitransparent] (.6,0) rectangle (4.2,2.8);
\node at (1,2.4) {$\ddots$};
\draw [arrow] (-1,1.7) to (2,1.7);
\draw [thick] (.1,1.7) to  [out=0, in=180] (.6,2);
\draw [arrow] (.6,2) to (2,2);
\draw [ultra thick] (-.2,-.3) rectangle (4.6,3.1);
\node at (2,3.4) {$\beta_2 = \sigma_1 \comp \cdots \comp \sigma_p$};
\node at (-.8,1.9) {\tiny{$z_1$}};
\node at (-.6,1.3) {$\vdots$};
\node at (-.8,.7) {\tiny{$z_p$}};
\node at (-1.5,1.2) {$\win$};
\end{tikzpicture}
\end{center}
There are $p-1$ gray boxes.  The outermost gray box will be called $W_1$, and the innermost gray box will be called $W_{p-1}$ below.  To define such a decomposition, we will need the following definitions.
\begin{enumerate}
\item
For $0 \leq j \leq p$ define the box $W_j \in \boxs$ with $\wout_j = \wout$ and
\[
\win_j = 
\underbrace{\left[\coprod_{i=1}^j \sbeta^{-1}(z_i)\right]}_{\text{$\varnothing$ if $j=0$}} \amalg \underbrace{\{z_{j+1},\ldots,z_p\}}_{\text{$\varnothing$ if $j=p$}}.
\]
Note that $\win_0 = \win$ by definition, while $\win_p = \xin = \vin$ by Lemma \ref{lemma:sbeta-surjective}.  So $W_0 = W$ and $W_p = V$.
\item
For $1 \leq j \leq p$ define the wiring diagram $\sigma_j \in \WD\wjminusonewj$ with no delay nodes and with supplier assignment
\[
\nicexy{
\dmsigmaj = \wout_{j-1} \amalg \win_j 
= \wout \amalg \left[\coprod\limits_{i=1}^j \sbeta^{-1}(z_i)\right] \amalg \{z_{j+1},\ldots,z_p\}
\ar[d]_-{\ssigmaj}\\
\supplysigmaj = \win_{j-1} \amalg \wout_j
= \left[\coprod\limits_{i=1}^{j-1} \sbeta^{-1}(z_i)\right] \amalg \{z_{j},\ldots,z_p\} 
\amalg \wout
}\]
given by
\[
\ssigmaj(x) = \begin{cases}
z_j & \text{ if $x \in \sbeta^{-1}(z_j)$};\\
x & \text{ otherwise}.
\end{cases}
\]
\end{enumerate}

Since
\[
\win_j = \left[\win_{j-1} \setminus \{z_j\}\right] \amalg \sbeta^{-1}(z_j),
\]
by Lemma \ref{lemma:iterated-insplits} $\sigma_j$ is:
\begin{itemize}
\item
a colored unit if $|\sbeta^{-1}(z_j)| = 1$;
\item
an iterated operadic composition of $\left(|\sbeta^{-1}(z_j)| - 1\right)$ in-splits if $|\sbeta^{-1}(z_j)| \geq 2$.  
\end{itemize}
Therefore, to show that $\beta_2$ is either a colored unit or an iterated operadic composition of in-splits, it is enough to check that there is a decomposition
\[
\beta_2 = \sigma_1 \comp \cdots \comp \sigma_p \in \WD\wv.
\]
Since the iterated operadic composition on the right has no delay nodes, it remains to check that its supplier assignment $s$ is equal to $\sbetatwo$ \eqref{supply-betatwo}.

On $\wout = \zout \subseteq \dmsigmaone$ the supplier assignment $s$ is the composition of $p$ identity functions, hence the identity function.  On
\[
\vin = \xin = \coprod_{i=1}^p \sbeta^{-1}(z_i) \subseteq \dmsigmap
\]
the supplier assignment $s$ sends elements in each $\sbeta^{-1}(z_j)$ to $z_j \in \win$, so it is equal to $\sbetatwo$.
\end{proof}

\begin{motivation}
Next, to show that $\beta_3$ is either a colored unit or an iterated operadic  composition of out-splits, we first need a lemma that says that the following wiring diagram is generated by out-splits.
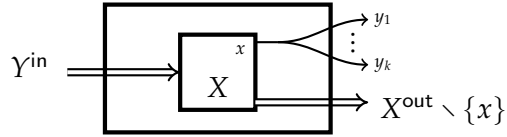
\begin{center}
\begin{tikzpicture}
\draw [ultra thick] (2,1) rectangle (3,2);
\node at (2.5,1.3) {$X$};
\node at (2.8,1.8) {\tiny{$x$}};
\draw [implies] (.5,1.5) to (2,1.5);
\node at (0,1.6) {$\yin$};
\draw [thick] (3,1.9) -- (3.3,1.9);
\draw [arrow] (3.3,1.9) to [out=0, in=180] (4.5,2.2);
\draw [arrow] (3.3,1.9) to [out=0, in=180] (4.5,1.6);
\node at (4.7,2.2) {\tiny{$y_1$}};
\node at (4.7,1.6) {\tiny{$y_k$}};
\node at (4.3,1.9) {$\vdots$};
\draw [implies] (3,1.1) to (4.5,1.1);
\node at (5.5,1) {$\xout \setminus \{x\}$};
\draw [ultra thick] (1,0.7) rectangle (4,2.4);
\end{tikzpicture}
\end{center}
The following observation is the out-split analogue of Lemma \ref{lemma:iterated-insplits}.
\end{motivation}

\begin{lemma}
\label{lemma:iterated-outsplits}
Suppose:
\begin{itemize}
\item
$X, Y \in \boxs$ such that $\xin = \yin$.
\item
There exist $x \in \xout$ and distinct elements $y_1, \ldots, y_k \in \yout$ with $k \geq 1$ such that 
\[
\yout = \left[\xout \setminus \{x\}\right] \amalg \{y_1,\ldots,y_k\}
\]
and $v(x) = v(y_i) \in S$ for $1 \leq i \leq k$.
\item
$\sigma \in \WD\yx$ is a wiring diagram with no delay nodes and with supplier assignment
\[
\nicexy{
\dmsigma = \yout \amalg \xin =
\left[\xout \setminus \{x\}\right] \amalg \{y_1,\ldots,y_k\} \amalg \xin
\ar[d]_-{\ssigma}
\\
\supplysigma = \yin \amalg \xout = \xin \amalg \xout 
}\]
given by
\[
\ssigma(z) = \begin{cases}
x & \text{ if $z = y_1, \ldots, y_k$};\\
z & \text{ if $z \in \left[\xout \setminus \{x\}\right] \amalg \xin$}.
\end{cases}
\]
\end{itemize}
Then:
\begin{enumerate}
\item
$\sigma$ is the $Y$-colored unit if $k=1$;
\item
$\sigma$ is an iterated operadic composition of $(k-1)$ out-splits if $k \geq 2$.
\end{enumerate}
\end{lemma}

\begin{proof}
If $k=1$, then $\ssigma$ is the identity function, so $\sigma$ is a colored unit.

The assertion for $k \geq 2$ is proved by induction.  If $k = 2$, then $\sigma$ is by definition the out-split $\sigma^{Y,y_1,y_2}$ (Def. \ref{def:out-split}).

Suppose $k \geq 3$.  We will factor $\sigma$ into two wiring diagrams as depicted in the picture
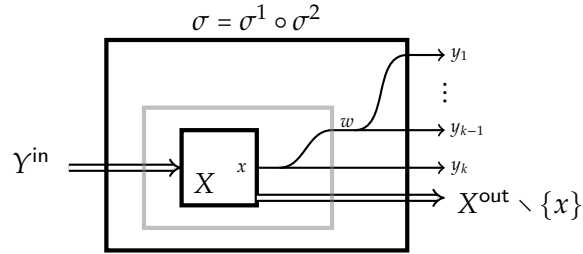
\begin{center}
\begin{tikzpicture}
\draw [ultra thick] (2,1) rectangle (3,2);
\node at (2.3,1.3) {$X$};
\node at (2.8,1.5) {\tiny{$x$}};
\draw [implies] (.5,1.5) to (2,1.5);
\node at (0,1.6) {$\yin$};
\draw [arrow] (3,1.5) -- (5.5,1.5);
\draw [thick] (3.3,1.5) to [out=0, in=180] (4,2);
\node at (5.7,1.5) {\tiny{$y_k$}};
\draw [arrow] (4,2) to (5.5,2);
\draw [thick] (4.3,2) to [out=0, in=180] (5,3);
\draw [arrow] (5,3) to (5.5,3);
\node at (5.7,3) {\tiny{$y_1$}};
\node at (5.5,2.5) {$\vdots$};
\node at (5.8,2) {\tiny{$y_{k-1}$}};
\draw [implies] (3,1.1) to (5.5,1.1);
\node at (6.5,1) {$\xout \setminus \{x\}$};
\draw [ultra thick, gray, semitransparent] (1.5,.7) rectangle (4,2.3);
\node at (4.2,2.1) {\tiny{$w$}};
\draw [ultra thick] (1,0.4) rectangle (5,3.2);
\node at (3,3.5) {$\sigma = \sigma^1 \comp \sigma^2$};
\end{tikzpicture}
\end{center}
in which the intermediate gray box will be called $W$ below.  The inner wiring diagram $\sigma^2$ will be an out-split, and the outer wiring diagram $\sigma^1$ will be an iterated operadic composition of $(k-2)$ out-splits.  To define such a decomposition, we will need the following definitions.
\begin{enumerate}
\item
Define the box $W \in \boxs$ with $\win  = \xin = \yin$ and
\[
\wout = \left[\xout \setminus \{x\}\right] \amalg \{w,y_k\}
\]
for some $w \not= y_k$ such that $v(w) = v(y_k) \in S$.  In particular, we have
\begin{equation}
\label{yout-out}
\yout = \left[\wout \setminus \{w\}\right] \amalg \{y_1, \ldots, y_{k-1}\}. 
\end{equation}
\item
Define the wiring diagram $\sigma^1 \in \WD\yw$ with no delay nodes and with supplier assignment
\[
\nicexy{
\dmsigmaupone = \yout \amalg \win
= \left[\xout \setminus \{x\}\right] \amalg \{y_1,\ldots,y_k\} \amalg \yin
\ar[d]_-{\ssigmaupone}\\
\supplysigmaupone = \yin \amalg \wout
= \yin \amalg  \left[\xout \setminus \{x\}\right] \amalg \{w,y_k\}
}\]
given by
\[
\ssigmaupone(z) = \begin{cases}
w & \text{ if $z = y_1,\ldots,y_{k-1}$};\\
z & \text{ otherwise}.
\end{cases}
\]
\item
Define the wiring diagram $\sigma^2 \in \WD\wx$ with no delay nodes and with supplier assignment
\[
\nicexy{
\dmsigmauptwo = \wout \amalg \xin
=  \left[\xout \setminus \{x\}\right] \amalg \{w,y_k\} \amalg \win
\ar[d]_-{\ssigmauptwo}\\
\supplysigmauptwo = \win \amalg \xout
}\]
given by
\[
\ssigmauptwo(z) = \begin{cases}
x & \text{ if $z = w,y_k$};\\
z & \text{ otherwise}.
\end{cases}
\]
\end{enumerate}
Then $\sigma^1 \comp \sigma^2 \in \WD\yx$ is a wiring diagram with no delay nodes.  To see that it is equal to $\sigma$, it suffices to check that the supplier assignment of $\sigma^1 \comp \sigma^2$ is equal to $\ssigma$.  This is true by a direct inspection of $\ssigmaupone$ and $\ssigmauptwo$.

By \eqref{yout-out} the induction hypothesis applies to $\sigma^1$, which says that it is an iterated operadic composition of $(k-2)$ out-splits.  By definition $\sigma^2$ is the out-split $\sigma^{W,w,y_k}$.    Combined with the previous paragraph, it follows that $\sigma = \sigma^1 \comp \sigma^2$ is the iterated operadic composition of $(k-1)$ out-splits, finishing the induction.
\end{proof}

Finally, we consider $\beta_3$.

\begin{lemma}
\label{lemma:betathree}
In the context of Def. \ref{def:beta}, suppose $\beta$ has no internal wasted wires (such as $\pi_2$ \eqref{pi-two-supplier} by Lemma \ref{lemma:pitwo-wasted-wired}).  Then $\beta_3 \in \WD\vx$ is either a colored unit or an iterated operadic composition of out-splits.
\end{lemma}

\begin{proof}
Recall that the wiring diagram $\beta_3 \in \WD\vx$ has no delay nodes and has supplier assignment
\eqref{supply-betathree}
\[
\nicexy{
\dmbetathree = \vout \amalg \xin = \zout \amalg \xin
\ar[d]_-{\sbetathree}
\\
\supplybetathree = \vin \amalg \xout = \xin \amalg \xout
}
\]
the coproduct of the identity function on $\xin$ and $\sbeta : \zout \to \xout$.  Since $\beta$ has no internal wasted wires, by Lemma \ref{lemma:sbeta-surjective} the map $\sbeta : \zout \to \xout$ is surjective.

Write $\xout = \{x_1,\ldots,x_r\}$, so $\zout = \coprod_{i=1}^r \sbeta^{-1}(x_i)$ with each $\sbeta^{-1}(x_i)$ non-empty.  If $r = 0$, then $\xout = \varnothing = \zout$, and $\beta_3$ is a colored unit.

Suppose $r > 0$.  We will decompose $\beta_3$ as in the picture:
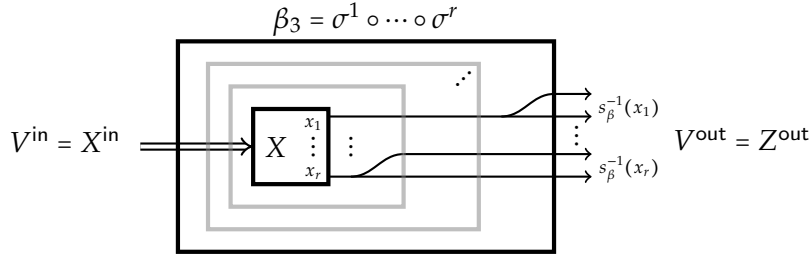
\begin{center}
\begin{tikzpicture}
\draw [ultra thick] (2,1) rectangle (3,2);
\node at (2.3,1.5) {$X$};
\node at (2.8,1.8) {\tiny{$x_1$}};
\node at (2.8,1.5) {$\vdots$};
\node at (2.8,1.15) {\tiny{$x_r$}};
\draw [implies] (.5,1.5) to (2,1.5);
\node at (-.5,1.6) {$\vin = \xin$};
\draw [arrow] (3,1.9) -- (6.5,1.9);
\draw [thick] (5.3,1.9) to [out=0, in=180] (6,2.2);
\draw [arrow] (6,2.2) -- (6.5,2.2);
\node at (7,2) {\tiny{$\sbeta^{-1}(x_1)$}};
\draw [arrow] (3,1.1) -- (6.5,1.1);
\draw [thick] (3.3,1.1) to [out=0, in=180] (4,1.4);
\draw [arrow] (4,1.4) -- (6.5,1.4);
\node at (7,1.2) {\tiny{$\sbeta^{-1}(x_r)$}};
\draw [ultra thick, gray, semitransparent] (1.7,.7) rectangle (4,2.3);
\node at (3.3,1.5) {$\vdots$};
\draw [ultra thick, gray, semitransparent] (1.4,.4) rectangle (5,2.6);
\node at (4.8,2.4) {\reflectbox{$\ddots$}};
\draw [ultra thick] (1,0.1) rectangle (6,2.9);
\node at (6.3,1.65) {$\vdots$};
\node at (8.5,1.6) {$\vout = \zout$};
\node at (3.5,3.2) {$\beta_3 = \sigma^1 \comp \cdots \comp \sigma^r$};
\end{tikzpicture}
\end{center}
The outermost gray box will be called $V_1$, and the innermost gray box will be called $V_{r-1}$ below.  To define such a decomposition, we first need some definitions.
\begin{enumerate}
\item
For $0 \leq j \leq r$ define the box $V_j \in \boxs$ with $\vin_j = \vin = \xin$ and
\[
\vout_j = \underbrace{\{x_1, \ldots, x_j\}}_{\text{$\varnothing$ if $j=0$}} \amalg \underbrace{\coprod_{i=j+1}^r \sbeta^{-1}(x_i)}_{\text{$\varnothing$ if $j=r$}}.
\]
Note that $\vout_0 = \zout = \vout$, so $V_0 = V$.  Also, $\vout_r = \xout$, so $V_r = X$.
\item
For $1 \leq j \leq r$ define the wiring diagram $\sigma^j \in \WD\vjminusonevj$ with no delay nodes and with supplier assignment
\[
\nicexy{
\dmsigmaupj = \vout_{j-1} \amalg \vin_j
= \{x_1, \ldots, x_{j-1}\} \amalg \coprod\limits_{i=j}^r \sbeta^{-1}(x_i) \amalg \xin
\ar[d]_-{\ssigmaj}\\
\supplysigmaupj = \vin_{j-1} \amalg \vout_j = \xin \amalg
 \{x_1, \ldots, x_{j}\} \amalg \coprod\limits_{i=j+1}^r \sbeta^{-1}(x_i)
}\]
given by
\[
\ssigmaj(z) = \begin{cases}
x_j & \text{ if $z \in \sbeta^{-1}(x_j)$};\\
z & \text{ otherwise}.
\end{cases}
\]
\end{enumerate}
The iterated operadic composition
\[
\sigma^1 \comp \cdots \comp \sigma^r \in \WD\vzerovr = \WD\vx
\]
has no delay nodes.  To see that it is equal to $\beta_3$, it remains to check that the former's supplier assignment $s$ is equal to $\sbetathree$.   
\begin{itemize}
\item
On $\xin \subseteq \dmsigmaupr$ the supplier assignment $s$ is the composition of $r$ identity functions, hence the identity function.  
\item
On $\vout = \zout \subseteq \dmsigmaupone$ the supplier assignment $s$ sends elements in each $\sbeta^{-1}(x_i)$ to $x_i$.  
\end{itemize}
So $s$ is equal to $\sbetathree$.

By Lemma \ref{lemma:iterated-outsplits} each $\sigma^j$ for $1 \leq j \leq r$ is either a colored unit or an iterated operadic composition of out-splits.  Therefore, $\beta_3 = \sigma^1 \comp \cdots \comp \sigma^r$ is either a colored unit or an iterated operadic composition of out-splits.
\end{proof}

\section{Summary of Chapter \ref{ch04-decomposition}}

\begin{enumerate}
\item A wiring diagram with no input boxes and no delay nodes is generated by the empty wiring diagram and a finite number of $1$-wasted wires.
\item Every wiring diagram $\psi$ has a decomposition
\[\psi = \pi_1 \comp \beta_1 \comp \beta_2 \comp \beta_3 \comp \varphi\]
in which:
\begin{itemize}
\item $\pi_1$ is generated by $1$-loops;
\item $\beta_1$ is generated by $1$-wasted wires;
\item $\beta_2$ is generated by in-splits;
\item $\beta_3$ is generated by out-splits;
\item $\varphi$ is either the empty wiring diagram or is generated by $2$-cells and $1$-delay nodes.
\end{itemize}
\end{enumerate}

\chapter{Finite Presentation}
\label{ch05-stratified-presentation}

Fix a class $S$, with respect to which the $\boxs$-colored operad $\WD$ of wiring diagrams is defined (Theorem \ref{wd-operad}).  The main purpose of this chapter is to establish finite presentations for the operad $\WD$ of wiring diagrams and its variants $\wddot$ and $\wdzero$.  For the operad $\WD$, our finite presentation means the following two statements.
\begin{enumerate}
\item
The $8$ generating wiring diagrams (Def. \ref{def:generating-wiring-diagrams}) generate the operad $\WD$.  This means that every wiring diagram can be expressed as a finite iterated operadic composition involving only generating wiring diagrams.
\item
If a wiring diagram can be operadically generated by the generating wiring diagrams in two different ways, then there exists a finite sequence of \emph{elementary equivalences} from the first iterated operadic composition to the other one.  An elementary equivalence is induced by either an elementary relation (Def. \ref{def:elementary-relations}) or an operad associativity/unity axiom for the generating wiring diagrams.
\end{enumerate}
In Chapter \ref{ch06-wd-algebras} we will use these finite presentations to describe algebras over the operads $\WD$, $\wddot$, and $\wdzero$ in terms of finitely many generating structure maps and generating axioms corresponding to the generating wiring diagrams and elementary relations.  In Section \ref{sec:propagator} we will use the finite presentation for $\WD$-algebras to study the propagator algebra.  In Section \ref{sec:algebra-ods}  we will use the finite presentation for $\wdzero$-algebras to study the algebra of open dynamical systems.

In Section \ref{sec:stratified} we establish the first part of the finite presentation theorem for $\WD$ by showing that every wiring diagram has a \emph{stratified presentation} (Theorem \ref{stratified-presentation-exists}).  A stratified presentation (Def. \ref{def:stratified-presentation}) is a highly structured iterated operadic composition of the generating wiring diagrams.  The proof of the second part of the finite presentation theorem also requires the use of stratified presentations.

In Section \ref{sec:coherence-wd} we establish the second part of the finite presentation theorem for $\WD$.  We show that any two presentations of the same wiring diagram in terms of generating wiring diagrams are connected by a finite sequence of elementary equivalences (Theorem \ref{thm:wd-generator-relation}).

In Section \ref{sec:coherence-normal} we establish a finite presentation for the operad $\wddot$ of wiring diagrams without delay nodes, which we call \emph{normal} wiring diagrams.  Normal wiring diagrams appeared in Spivak's study of mode-dependent networks and dynamical systems \cite{spivak15,spivak15b}.

In Section \ref{sec:coherence-strict} we restrict further and establish a finite presentation for the operad $\wdzero$ of wiring diagrams without delay nodes and whose supplier assignments are  bijections. We call them \emph{strict} wiring diagrams.  They appeared in \cite{vsl}.  We will use strict wiring diagrams in Section \ref{sec:algebra-ods} to study the algebra of open dynamical systems.

\section{Stratified Presentation}
\label{sec:stratified}

In this section, we define a \emph{stratified presentation} and show that every wiring diagram has a stratified presentation (Theorem \ref{stratified-presentation-exists}).  We also need the concept of a \emph{simplex} to discuss generators and relations in the operad $\WD$ of wiring diagrams.

\begin{motivation}
In plain language, a simplex is a finite parenthesized word whose alphabets are generating wiring diagrams, in which each pair of parentheses has a well defined associated $\compi$-composition.  In particular, a simplex has a well defined operadic composition.  As we have seen in Chapter \ref{ch03-generating-wd}, it is often possible to express a wiring diagram as an operadic composition of generating wiring diagrams in multiple ways.  In other words, a wiring diagram can have many different simplex presentations.  We now start to develop the necessary language to say precisely that any two such simplex presentations of the same wiring diagram are equivalent in some way.
\end{motivation}

\begin{definition}
\label{def:simplex}
Suppose $n \geq 1$.  An \emph{$n$-simplex} \index{simplex} $\Psi$ and its \emph{composition} \index{composition of a simplex} \label{notation:simplex-composition}$|\Psi| \in \WD$ are defined inductively as follows.
\begin{enumerate}
\item
A \emph{$1$-simplex} is a generating wiring diagram  (Def. \ref{def:generating-wiring-diagrams}) $\psi$.  Its \emph{composition} $|\psi|$ is defined as $\psi$ itself.
\item
Suppose $n \geq 2$ and that $k$-simplices for $1 \leq k \leq n-1$ and their compositions in $\WD$  are already defined.  An \emph{$n$-simplex} is a tuple $\Psi = \left(\upsi, i, \uphi\right)$ consisting of
\begin{itemize}
\item
an integer $i \geq 1$,
\item
a $p$-simplex $\upsi$ for some $p \geq 1$, and
\item
a $q$-simplex $\uphi$ for some $q \geq 1$
\end{itemize}
such that:
\begin{enumerate}[(i)]
\item
$p+q=n$;
\item
the operadic composition
\begin{equation}
\label{length-n-simplex}
|\Psi| \defn \left|\upsi\right| \compi \left|\uphi\right|
\end{equation}
is defined in $\WD$ (Def. \ref{def:compi-wd}).  
\end{enumerate}
The wiring diagram $|\Psi|$ in \eqref{length-n-simplex} is the \emph{composition} of $\Psi$.
\end{enumerate}
A \emph{simplex} in $\WD$ is an $m$-simplex in $\WD$ for some $m \geq 1$.  We say that a simplex $\Psi$ is a \emph{presentation} \index{presentation} of the wiring diagram $|\Psi|$.
\end{definition}

\begin{notation}\label{notation:simplex}
To simplify notations, we will sometimes use the right side of \eqref{length-n-simplex} to denote a simplex.  To simplify notations even further, we may even just list the generating wiring diagrams $(\psi_1,\ldots,\psi_n)$ in a simplex in the order in which they appear in the composition \eqref{length-n-simplex}, omitting all the pairs of parentheses and the operadic compositions from the notations.
\end{notation}

\begin{remark}\label{rk:simplex}
In Def. \ref{def:simplex} we could have made the definition for a general operad $\O$ other than $\WD$, using a specified collection of elements in $\O$ in place of the generating wiring diagrams.  Such a definition would be useful in discussing generators and relations in a general operad $\O$.
\end{remark}

In the next three Examples, every $\psi_i$ denotes a generating wiring diagram, and we will use Notation \ref{notation:simplex}.

\begin{example}
A $2$-simplex has the form $\bigl(\psi_1, i, \psi_2\bigr)$, which we abbreviate to  $\psi_1 \compi \psi_2$, for some integer $i \geq 1$.  For instance, suppose $d \in S$, and $X$ is the box with $\xin = \{d\} = \xout$.  Suppose $Y$ is an arbitrary box.  Then there is a $2$-simplex
\[\bigl(\theta_{X,Y}, 1, \delta_d\bigr)\]
in which $\theta_{X,Y}$ is a $2$-cell (Def. \ref{def:theta-wd}) and $\delta_d$ is a $1$-delay node (Def. \ref{def:one-dn}).  Its composition $\theta_{X,Y} \comp_1 \delta_d$ is the wiring diagram
\begin{center}
\begin{tikzpicture}[scale=.7]
\draw [ultra thick] (1.5,2.5) circle [radius=0.5];
\node at (1.5,2.5) {$d$};
\draw [arrow] (-.5,2.5) -- (1,2.5);
\draw [arrow] (2,2.5) -- (3.5,2.5);
\draw [ultra thick] (1,.5) rectangle (2,1.5);
\draw [implies] (-.5,1) to (1,1);
\draw [implies] (2,1) to (3.5,1);
\node at (1.5,1) {$Y$};
\draw [ultra thick] (0,0.2) rectangle (3,3.3);
\node at (1.5,3.7) {$\theta_{X,Y} \comp_1 \delta_d$};
\end{tikzpicture}
\end{center}
in $\WD\xplusyy$ with one delay node.
\end{example}

\begin{example}
A $3$-simplex is an iterated operadic composition in $\WD$ of the form
\[(\psi_1 \compi \psi_2) \compj \psi_3 \orspace
\psi_1 \compi (\psi_2 \compj \psi_3)\]
for some integers $i,j \geq 1$.  Once again these are really abbreviations for the $3$-simplices
\[\bigl(\left(\psi_1, i, \psi_2\right), j, \psi_3\bigr) \orspace
\bigl(\psi_1, i, \left(\psi_2, j, \psi_3\right)\bigr).\]
For instance, continuing the example above, suppose $\lambda_{X \amalg Y, d} \in \WD\yxplusy$ is a $1$-loop (Def. \ref{def:loop-wd}).  Then there is a $3$-simplex
\[\Bigl(\lambda_{X \amalg Y,d}, 1, \bigl(\theta_{X,Y}, 1, \delta_d\bigr)\Bigr)\]
whose composition $\lambda_{X\amalg Y,d} \comp_1 (\theta_{X,Y} \comp_1 \delta_d)$ is the wiring diagram
\begin{center}
\begin{tikzpicture}[scale=.7]
\draw [ultra thick] (1.5,2.5) circle [radius=0.5];
\node at (1.5,2.5) {$d$};
\draw [arrow, looseness=4.5] (2,2.5) to [out=30, in=150] (1,2.5);
\draw [ultra thick] (1,.5) rectangle (2,1.5);
\draw [implies] (-.5,1) to (1,1);
\draw [implies] (2,1) to (3.5,1);
\node at (1.5,1) {$Y$};
\draw [ultra thick] (0,0.2) rectangle (3,3.7);
\end{tikzpicture}
\end{center}
in $\WD\yy$ with one delay node.
\end{example}

\begin{example}
 A $4$-simplex is an iterated operadic composition in $\WD$ of the form
\[\begin{split}
&\bigl((\psi_1 \compi \psi_2) \compj \psi_3\bigr) \compk \psi_4, \quad
\bigl(\psi_1 \compi (\psi_2 \compj \psi_3)\bigr) \compk \psi_4, \quad
(\psi_1 \compi \psi_2) \compj (\psi_3 \compk \psi_4),\\
&\psi_1 \compi \bigl((\psi_2 \compj \psi_3) \compk \psi_4\bigr), \orspace
\psi_1 \compi \bigl(\psi_2 \compj (\psi_3 \compk \psi_4)\bigr)
\end{split}\]
for some integers $i,j,k \geq 1$.  For instance, continuing the previous example, suppose $\yout = \{y\}$ and $Z$ is a box such that $\zout = \{z,z'\}$ with $v(z) = v(z') = v(y)$ and that $Z/(z=z') = Y$.  Suppose $\sigma^{Z,z,z'} \in \WD\zy$ is an out-split (Def. \ref{def:out-split}).  Then there is a $4$-simplex
\[\Bigl(\sigma^{Z,z,z'}, 1, \bigl(\lambda_{X \amalg Y,d}, 1, (\theta_{X,Y}, 1, \delta_d)\bigr) \Bigr)\]
whose composition
\[\sigma^{Z,z,z'} \comp_1 \Bigl(\lambda_{X\amalg Y,d} \comp_1 (\theta_{X,Y} \comp_1 \delta_d)\Bigr)\]
is the wiring diagram
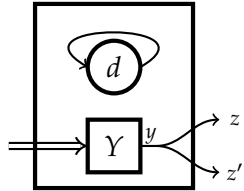
\begin{center}
\begin{tikzpicture}[scale=.7]
\draw [ultra thick] (1.5,2.5) circle [radius=0.5];
\node at (1.5,2.5) {$d$};
\draw [arrow, looseness=4.5] (2,2.5) to [out=30, in=150] (1,2.5);
\draw [ultra thick] (1,.5) rectangle (2,1.5);
\draw [implies] (-.5,1) to (1,1);
\draw [line] (2,1) to (2.3,1);
\node at (2.2,1.2) {\scriptsize{$y$}};
\draw [arrow, out=0, in=180] (2.3,1) to (3.5,1.5);
\node at (3.8,1.5) {\scriptsize{$z$}};
\draw [arrow, out=0, in=180] (2.3,1) to (3.5,.5);
\node at (3.8,.5) {\scriptsize{$z'$}};
\node at (1.5,1) {$Y$};
\draw [ultra thick] (0,0.2) rectangle (3,3.7);
\end{tikzpicture}
\end{center}
in $\WD\zy$ with one delay node.
\end{example}

In Section \ref{sec:coherence-wd} we will show that any two presentations of the same wiring diagram are equivalent in a certain way.  For this purpose, we will need a more structured kind of presentation.

\begin{motivation}
If we think of a simplex as a parenthesized word whose alphabets are generating wiring diagrams, then the stratified simplex in the next definition is a word where the same alphabets must occur in a consecutive string.  For example, all the $1$-loops must occur together as a string $(\lambda_1,\ldots, \lambda_n)$.  Furthermore, we can even insist that these strings for different types of generating wiring diagrams occur in a specific order, with name changes and $1$-loops at the top and with $1$-delay nodes and $2$-cells at the bottom.
\end{motivation}

\begin{definition}
\label{def:stratified-presentation}
A \emph{stratified simplex} \index{stratified simplex} in $\WD$ is a simplex in $\WD$ (Def. \ref{def:simplex}) of one of the following two forms, where Notation \ref{notation:simplex} is used:
\begin{enumerate}
\item
$\bigl(\uomega, \epsilon\bigr)$, where:
\begin{itemize}
\item
$\uomega$ is a possibly empty string of $1$-wasted wires (Def. \ref{def:wasted-wire-wd});
\item
$\epsilon$ is the empty wiring diagram (Def. \ref{def:empty-wd}).
\end{itemize}
\item
$\Bigl(\tau, \ulambda, \uomega, \usigma_*, \usigma^*, \utheta, \udelta\Bigr)$, where:
\begin{itemize}
\item
$\tau$ is a name change (Def. \ref{def:name-change});
\item
$\ulambda$ is a possibly empty string of $1$-loops (Def. \ref{def:loop-wd});
\item
$\uomega$ is a possibly empty string of $1$-wasted wires;
\item
$\usigma_*$ is a possibly empty string of in-splits (Def. \ref{def:insplit-wd});
\item
$\usigma^*$  is a possibly empty string of out-splits (Def. \ref{def:out-split});
\item
$\utheta$ is a possibly empty string of $2$-cells (Def. \ref{def:theta-wd});
\item
$\udelta$ is a possibly empty string of $1$-delay nodes (Def. \ref{def:one-dn}).
\end{itemize}
\end{enumerate}
We call these stratified simplices of type (1) and of type (2), respectively.  If $\Psi$ is a stratified simplex, then we call it a \emph{stratified presentation} \index{stratified presentation} of the wiring diagram $|\Psi|$.
\end{definition}

\begin{remark}
\label{rk:stratified-disjoint}
Stratified simplices of type (1) and of type (2) are mutually exclusive.  Indeed, the composition of a stratified simplex of type (1) has no input boxes and no delay nodes.  On the other hand, the composition of a stratified simplex of type (2) either has at least one input box or at least one delay node or both.
\end{remark}

Using the decompositions in the previous Chapter, we now observe that the generating wiring diagrams generate the operad $\WD$ of wiring diagrams in a highly structured way.

\begin{theorem}
\label{stratified-presentation-exists}
Every wiring diagram has a stratified presentation (Def. \ref{def:stratified-presentation}).\index{stratified presentations exist}
\end{theorem}

\begin{proof}
Suppose $\psi \in \WD\yux$ is a general wiring diagram as in Assumption \ref{assumption:psi} with input boxes $\uX = (X_1,\ldots,X_N)$ and delay nodes $\dnpsi = \{d_1,\ldots,d_r\}$.  Recall Notation \ref{comp-is-compone} for (iterated) $\compone$.

If $N=r=0$, then $\psi$ has a stratified presentation of type (1) by Lemma \ref{lemma:N=r=0}.

Next suppose $N+r \geq 1$.  We use the decomposition $\psi = \alpha \comp \varphi$ \eqref{psi-is-alpha-comp-phi} and show that $\psi$ has a stratified presentation of type (2) using the following observations.
\begin{enumerate}
\item
If $N+r = 1$, then $\varphi$ is either a colored unit, which can be ignored in a simplex by Lemma \ref{compi-unity}, or a $1$-delay node by Lemma \ref{lemma:phi-small-Nr}.  
\item
If $N+r \geq 2$, then $\varphi$ has a stratified presentation $\bigl(\utheta, \udelta\bigr)$ consisting of $2$-cells and $1$-delay nodes by Lemma \ref{lemma:phi-presentation} and by the equivalence between $\gamma$ and the $\compi$-compositions \eqref{gamma-in-comps}.
\item
By definition $\alpha$ \eqref{alpha-wiring-diagram} has one input box and no delay nodes.  There is a decomposition $\alpha = \pi_1 \comp \pi_2$ by Lemma \ref{lemma:pi-decomposition}.  The outer wiring diagram $\pi_1$ is either a colored unit or has a stratified presentation $\left(\ulambda\right)$ consisting of $1$-loops by Lemma \ref{lemma:iterated-loops}.  
\item
Furthermore, by Lemma \ref{lemma:beta-factor} there is a decomposition $\pi_2 = \beta_1 \comp \beta_2 \comp \beta_3$ in which:
\begin{itemize}
\item
$\beta_1$ is a colored unit or has a stratified presentation $\left(\uomega\right)$ consisting of $1$-wasted wires by Lemma \ref{lemma:betaone}.
\item
$\beta_2$ is a colored unit or has a stratified presentation $\left(\usigma_*\right)$ consisting of in-splits by Lemma \ref{lemma:betatwo}.
\item
$\beta_3$  is a colored unit or has a stratified presentation $\left(\usigma^*\right)$ consisting of out-splits by Lemma \ref{lemma:betathree}.
\end{itemize}
\end{enumerate}
Using the decomposition
\[
\psi = \pi_1 \comp \beta_1 \comp \beta_2 \comp \beta_3 \comp \varphi
\]
together with Convention \ref{conv:ignore-name-change}, we obtain the desired stratified presentation of type (2) for $\psi$ when $N+r \geq 1$.
\end{proof}

\section{Finite Presentation for Wiring Diagrams}
\label{sec:coherence-wd}

The purpose of this section is to establish the second part of the finite presentation theorem for the operad $\WD$.  First we define precisely what it means for two presentations of the same wiring diagram to be related to each other.  Recall the $28$ elementary relations (Def. \ref{def:elementary-relations}) and the definition of a colored operad (Def. \ref{def:pseudo-operad}).  In what follows, we will regard each operad associativity or unity axiom as an equality.  We remind the reader of Notation \ref{notation:simplex} regarding simplices in $\WD$.

\begin{motivation}
Recall that a simplex is essentially a parenthesized word whose alphabets are generating wiring diagrams.  In the next definition, we develop the precise concept through which one simplex presentation of a wiring diagram may be replaced by another.  We only allow replacement of strings within a simplex corresponding to either one of the $28$ elementary relations or an operad associativity/unity axiom.  When such a replacement within a simplex is possible, we say that the two simplices are elementarily equivalent.
\end{motivation}

\begin{definition}
\label{def:equivalent-simplices}
Suppose $\Psi$ is an $n$-simplex in $\WD$ as in Def. \ref{def:simplex}.
\begin{enumerate}
\item
A \emph{subsimplex} \index{subsimplex} of $\Psi$ is a simplex in $\WD$ defined inductively as follows.
\begin{itemize}
\item
If $\Psi$ is a $1$-simplex, then a \emph{subsimplex} of $\Psi$ is $\Psi$ itself.
\item
Suppose $n \geq 2$ and $\Psi =  \left(\upsi, i, \uphi\right)$ for some $i \geq 1$, $p$-simplex $\upsi$, and $q$-simplex $\uphi$ with $p+q=n$.  Then a \emph{subsimplex} of $\Psi$ is 
\begin{itemize}
\item
a subsimplex of $\upsi$, 
\item
a subsimplex of $\uphi$, or
\item
$\Psi$ itself.
\end{itemize}
\end{itemize}
If $\Psi'$ is a subsimplex of $\Psi$, then we write $\Psi' \subseteq \Psi$.
\item
An \emph{elementary subsimplex} \index{elementary subsimplex} $\Psi'$ of $\Psi$ is a subsimplex of one of two forms: 
\begin{enumerate}[(i)]
\item
$\Psi'$ is one side (either left or right) of a specified elementary relation (Def. \ref{def:elementary-relations}).
\item
$\Psi'$ is one side (either left or right) of a specified operad associativity or unity axiom--\eqref{compi-associativity}, \eqref{compi-associativity-two}, \eqref{compi-left-unity}, or \eqref{compi-right-unity}--involving only the generating wiring diagrams  (Def. \ref{def:generating-wiring-diagrams}). 
\end{enumerate}
\item
Suppose $\Phi$ is another simplex in $\WD$.  Then $\Psi$ and $\Phi$ are said to be \emph{equivalent} \index{equivalent simplices} if their compositions are equal; i.e., $|\Psi| = |\Phi| \in \WD$.
\item
Suppose:
\begin{itemize}
\item
$\Psi' \subseteq \Psi$ is an elementary subsimplex corresponding to one side of $R$, which is either an elementary relation or an operad associativity/unity axiom for the generating wiring diagrams.
\item
$\Psi"$ is the simplex given by the other side of $R$.
\item
$\Psi^1$ is the simplex obtained from $\Psi$ by replacing the subsimplex $\Psi'$ by $\Psi"$.  
\end{itemize}
We say that $\Psi$ and $\Psi^1$ are \emph{elementarily equivalent}.\index{elementarily equivalent}  Note that elementarily equivalent simplices are also equivalent.
\item
If $\Psi$ and $\Phi$ are elementarily equivalent, we write $\Psi \sim \Phi$ and call this an \emph{elementary equivalence}.\index{elementary equivalence}
\item
Suppose $\Psi_0,\ldots,\Psi_r$ are simplices for some $r \geq 1$ and that there exist elementary equivalences
\[
\Psi_0 \sim \Psi_1 \sim \cdots \sim \Psi_r.
\]
Then we say that $\Psi_0$ and $\Psi_r$ are \emph{connected by a finite sequence of elementary equivalences}.  Note that in this case $\Psi_0$ and $\Psi_r$ are equivalent.
\end{enumerate}
\end{definition}

\begin{remark}
In the definition of an elementary subsimplex and an elementary equivalence, we did \emph{not} use the operad equivariance axiom \eqref{compi-eq}.  The reason is that the associativity and commutativity properties of $2$-cells--namely, the elementary relations \eqref{move:b1} and \eqref{move:b2}--are enough to guarantee the operad equivariance axiom involving only the generating wiring diagrams.
\end{remark}

\begin{example}
In a $3$-simplex $(\psi_1 \compi \psi_2) \compj \psi_3$, both
\[
(\psi_1, \psi_2) = \psi_1 \compi \psi_2
\andspace
(\psi_1,\psi_2,\psi_3) = (\psi_1 \compi \psi_2) \compj \psi_3
\]
are subsimplices.  However, $(\psi_2, \psi_3)$ is not a subsimplex. 
\end{example}

\begin{example}
A given wiring diagram may have many different equivalent presentations.  For  example, suppose $X \in \boxs$.  Then the $1$-simplex consisting of the $X$-colored unit $\tensorunit_X$ \eqref{wd-unit} is elementarily equivalent to:
\begin{enumerate}
\item
the $2$-simplex $\theta_{X,\varnothing} \comptwo \epsilon$  by \eqref{move:b0}; 
\item
the $3$-simplex $\left(\lambda_{Z,x}\right) \comp \left(\omega_{Z,x_1}\right) \comp \left(\sigma^{Y,x^1,x^2}\right)$ by \eqref{move:c6};
\item
the $2$-simplex $\left(\sigma_{Y,x,y}\right) \comp \left(\omega_{Y,y}\right)$ by \eqref{move:d5}.
\end{enumerate}
Any two of these three simplices are connected by a finite sequence of elementary equivalences.  Note that elementarily equivalent simplices may have different lengths.
\end{example}

\begin{example}
\label{ex:hass-3-simplex}
Suppose:
\begin{itemize}
\item
$\theta_{X,Y} \in \WD\xplusyxy$ is a $2$-cell (Def. \ref{def:theta-wd}).
\item
$\theta_{V,W} \in \WD\smallbinom{X}{V,W}$ is a $2$-cell with $X = V \amalg W \in \boxs$.
\item
$\sigma_{T,t_1,t_2} \in \WD\smallbinom{Y}{T}$ is an in-split (Def. \ref{def:insplit-wd}).  
\end{itemize}
Then the $3$-simplices
\[\bigl(\theta_{X,Y} \compone \theta_{V,W}\bigr) \compthree \sigma_{T,t_1,t_2}
\andspace
\bigl(\theta_{X,Y} \comptwo \sigma_{T,t_1,t_2}\bigr) \compone \theta_{V,W}\]
are elementarily equivalent by the horizontal associativity axiom \eqref{compi-associativity}.  This elementary equivalence expresses the fact that the wiring diagram
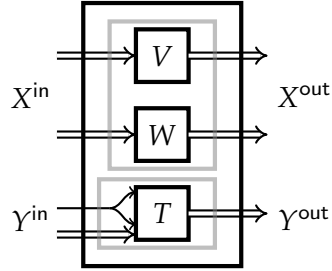
\begin{center}
\begin{tikzpicture}[scale=.7]
\draw [ultra thick] (-6,2) rectangle (-5,3);
\draw [implies] (-7.5,2.5) to (-6,2.5);
\draw [implies] (-5,2.5) to (-3.5,2.5);
\node at (-5.5,2.5) {$V$};
\draw [ultra thick] (-6,.5) rectangle (-5,1.5);
\draw [implies] (-7.5,1) to (-6,1);
\draw [implies] (-5,1) to (-3.5,1);
\node at (-5.5,1) {$W$};
\draw [ultra thick] (-6,-1) rectangle (-5,0);
\draw [line] (-7.5,-.4) to (-6.5,-.4);
\draw [arrow, out=0, in=180] (-6.5,-.4) to (-6,-.1);
\draw [arrow, out=0, in=180] (-6.5,-.4) to (-6,-.7);
\draw [implies] (-7.5,-.9) to (-6,-.9);
\draw [implies] (-5,-.5) to (-3.5,-.5);
\node at (-5.5,-.5) {$T$};
\draw [ultra thick] (-7,-1.5) rectangle (-4,3.5);
\node at (-8,1.75) {$\xin$};
\node at (-2.8,1.75) {$\xout$};
\node at (-8,-.6) {$\yin$};
\node at (-2.8,-.6) {$\yout$};
\draw [ultra thick, gray, semitransparent] (-6.5,.35) rectangle (-4.5,3.15);
\draw [ultra thick, gray, semitransparent] (-6.7,-1.15) rectangle (-4.5,.15);
\end{tikzpicture}
\end{center}
can be created from $\theta_{X,Y}$ by substituting in the two gray boxes in either order.
\end{example}

\begin{convention}
\label{conv:operad-axiom-move}
In what follows, to simplify the presentation, elementary equivalences corresponding to an operad associativity/unity axiom--\eqref{compi-associativity}, \eqref{compi-associativity-two}, \eqref{compi-left-unity}, or \eqref{compi-right-unity}--for the generating wiring diagrams will often be applied tacitly wherever necessary. For instance, an elementary equivalence given by replacing one of the $3$-simplices in Example \ref{ex:hass-3-simplex} by the other one will often \emph{not} be mentioned explicitly.
\end{convention}

Our next goal is to show that any two equivalent simplices are connected by a finite sequence of elementary equivalences.  In other words, with respect to the generating wiring diagrams, the $28$ elementary relations and the operad associativity/unity axioms--\eqref{compi-associativity}, \eqref{compi-associativity-two}, \eqref{compi-left-unity}, and \eqref{compi-right-unity}--for the generating wiring diagrams generate \emph{all} the relations in $\WD$.   During the first reading, the reader may wish to skip the proofs of the following three Lemmas.

The first step is to show that every simplex is connected to a stratified simplex in the following sense.

\begin{lemma}
\label{lemma:simplex-to-stratified}
Every simplex is either a stratified simplex or is connected to an equivalent stratified simplex by a finite sequence of elementary equivalences (Def. \ref{def:equivalent-simplices}).
\end{lemma}

\begin{proof}
Using Notation \ref{notation:simplex} suppose $\Psi = (\psi_1, \ldots, \psi_n)$ is a simplex with composition $|\Psi| = \psi \in \WD\yux$ as in Assumption \ref{assumption:psi}.  So $\psi$ has input boxes $\uX = (X_1,\ldots,X_N)$ and delay nodes $\dnpsi = \{d_1,\ldots,d_r\}$.  Suppose $\Psi$ is not a stratified simplex.  We will show that $\Psi$ is connected to an equivalent stratified simplex by a finite sequence of elementary equivalences.

Using the five elementary relations \eqref{move:a2}--\eqref{move:a6}, first we  move all the name changes (Def. \ref{def:name-change}) in $\Psi$, if there are any, to the left.  Then we use the elementary relation \eqref{move:a1} repeatedly to compose them down into one name change.  Therefore, after a finite sequence of elementary equivalences, we may assume that there is at most one name change in $\Psi$, which is the left-most entry.  If there are further elementary equivalences later that create name changes, we will perform the same procedure without explicitly mentioning it.

The empty wiring diagram $\epsilon \in \WD\emptynothing$ (Def. \ref{def:empty-wd}) and the $1$-delay nodes $\delta_d \in \WD\dnothing$ (Def. \ref{def:one-dn}) have no input boxes, so no operadic composition of the forms $\epsilon \compi -$ or  $\delta_d \compi -$ can be defined.  Therefore, after a finite sequence of elementary equivalences corresponding to the horizontal associativity axiom  \eqref{compi-associativity}, we may assume that $\Psi$ has the form
\[
\Bigl(\tau,\Psi^1,\uepsilon, \udelta\Bigr)
\]
in which:
\begin{itemize}
\item
$\tau$ is a name change;
\item
all the $1$-delay nodes $\udelta$ are at the right-most entries;
\item
all the empty wiring diagrams $\uepsilon$ are just to their left;
\item
$\Psi^1$ is either empty or is a subsimplex involving $2$-cells  (Def. \ref{def:theta-wd}), $1$-loops (Def. \ref{def:loop-wd}), in-splits (Def. \ref{def:insplit-wd}), out-splits  (Def. \ref{def:out-split}), and  $1$-wasted wires (Def. \ref{def:wasted-wire-wd}).
\end{itemize}

Next we use the elementary relations \eqref{move:b3}--\eqref{move:b6} to move all the $2$-cells in $\Psi$ to just the left of $\uepsilon$.  Then we use the elementary relations \eqref{move:c2}--\eqref{move:c4} to move all the remaining $1$-loops to just the right of the name change $\tau$.  After that, we use the elementary relations  \eqref{move:d4} and \eqref{move:e3} to move all the $1$-wasted wires to just the right of the $1$-loops.  Then we use the elementary relation \eqref{move:d3} to move all the in-splits to just the right of the $1$-wasted wires.  So after a finite sequence of elementary equivalences, we may assume that the simplex $\Psi$ has the form
\begin{equation}
\label{almost-stratified}
\Bigl(\tau, \ulambda, \uomega, \usigma_*, \usigma^*, \utheta, \uepsilon, \udelta\Bigr).
\end{equation}
If the string $\uepsilon$ of empty wiring diagrams is empty, then we are done because this is now a stratified simplex of type (2).

So suppose the string $\uepsilon$ in \eqref{almost-stratified} is non-empty.  Using finitely many elementary equivalences corresponding to the elementary relations \eqref{move:b0}--\eqref{move:b2}, we may cancel all the unnecessary empty wiring diagrams in \eqref{almost-stratified}.  If there are no empty wiring diagrams left after the cancellation, then we have a stratified simplex of type (2).  

Suppose that, after the cancellation in the previous paragraph, the resulting string $\uepsilon$ is still non-empty.  Then it must contain a single empty wiring diagram $\epsilon$, and there are no $2$-cells $\utheta$ and no $1$-delay nodes $\udelta$ in the resulting simplex $\Psi$.  Since the output box of $\epsilon$ is the empty box, the current simplex $\Psi$ cannot have any $1$-loops $\ulambda$, in-splits $\usigma_*$, or out-splits $\usigma^*$.  Therefore, in this case the simplex \eqref{almost-stratified} has the form
\begin{equation}
\label{almost-type1-stratified}
\Bigl(\tau, \uomega, \epsilon\Bigr).
\end{equation}

There are now two cases.  First suppose the string $\uomega$ in \eqref{almost-type1-stratified} is empty.  Since the output box of $\epsilon$ is the empty box, in the simplex $\left(\tau,\epsilon\right)$ the name change $\tau$ must be the colored unit of the empty box.  So by the left unity axiom \eqref{compi-left-unity}, the simplex $\left(\tensorunit_{\varnothing}, \epsilon\right)$ is elementarily equivalent to the simplex $(\epsilon)$, which is a stratified simplex of type (1).

Next suppose the string $\uomega$ in \eqref{almost-type1-stratified} is non-empty.  Using finitely many elementary equivalences corresponding to the elementary relation \eqref{move:a6}, the simplex \eqref{almost-type1-stratified} is connected to a simplex of the form
\begin{equation}
\label{almost-type1-stratified-b}
\Bigl(\uomega, \tau, \epsilon\Bigr)
\end{equation}
with $\tau \comp \epsilon$ as one of the operadic compositions.  As in the previous case, the composition $\tau \comp \epsilon$ forces $\tau$ to be the colored unit of the empty box.  So the simplex $\left(\uomega, \tensorunit_{\varnothing}, \epsilon\right)$ in \eqref{almost-type1-stratified-b} is elementarily equivalent to the simplex $\left(\uomega,\epsilon\right)$, which is a stratified simplex of type (1).
\end{proof}

The next step is to show that equivalent stratified simplices are connected.  We begin with stratified simplices of type (1).

\begin{lemma}
\label{lemma:stratified-type1}
Any two equivalent stratified simplices of type (1) are either equal or are connected by a finite sequence of elementary equivalences (Def. \ref{def:equivalent-simplices}).
\end{lemma}

\begin{proof}
Suppose $\Psi^1 = \left(\uomega^1,\epsilon\right)$ and $\Psi^2 = \left(\uomega^2, \epsilon\right)$ are equivalent stratified simplices of type (1) with common composition $\psi$.  Then $\psi$ has no input boxes and no delay nodes, and its output box contains only external wasted wires as in Lemma \ref{lemma:N=r=0}.  Each $1$-wasted wire in each $\Psi^i$ creates one external wasted wire in $\psi$.  So the $1$-wasted wire strings $\uomega^1$ and $\uomega^2$ have the same length.  It follows that the simplices $\Psi^1$ and $\Psi^2$ are connected by a finite sequence of elementary equivalences corresponding to the elementary relation \eqref{move:f1} and the vertical associativity axiom \eqref{compi-associativity-two}.
\end{proof}

\begin{lemma}
\label{lemma:stratified-type2}
Any two equivalent stratified simplices of type (2) are either equal or are connected by a finite sequence of elementary equivalences (Def. \ref{def:equivalent-simplices}).
\end{lemma}

\begin{proof}
The proof consists of a series of reductions.  Suppose $\Psi^1$ and $\Psi^2$ are distinct but equivalent stratified simplices of type (2) with common composition $\psi \in \WD\yux$.  Using elementary equivalences corresponding to
\begin{itemize}
\item
the operad unity axioms \eqref{compi-left-unity} and \eqref{compi-right-unity}, 
\item
the elementary relations \eqref{move:b0},  \eqref{move:c6}, and  \eqref{move:d5} regarding colored units, and 
\item
other elementary relations that move the generating wiring diagrams around the simplices, 
\end{itemize}
we may assume that there are no unnecessary generating wiring diagrams in these stratified simplices.  Here \emph{unnecessary} refers to either a colored unit or generating wiring diagrams whose (iterated) operadic composition is a colored unit.

The name change $\tau^1$ in $\Psi^1$ has output box $Y$ and input box uniquely determined by $\psi$, and the same is true for the name change $\tau^2$ in $\Psi^2$.  It follows that $\tau^1$ is equal to $\tau^2$.  So we may assume that there are no name changes in the two stratified simplices $\Psi^i$.

The string of delay nodes $\udelta^i$ in each simplex $\Psi^i$ represents the set of delay nodes in $\psi$.  Therefore, the two $\Psi^i$ without their strings of delay nodes are also equivalent.  Moreover, if these simplices without delay nodes are connected by a finite sequence of elementary equivalences, then so are the two $\Psi^i$ themselves by the horizontal associativity axiom \eqref{compi-associativity}.  So we may assume that the wiring diagram $\psi$ and the two simplices $\Psi^i$ have no delay nodes.   At this stage, each stratified simplex $\Psi^i$ has the form
\[
\left(\ulambda^i, \uomega^i, \usigma_*^i, \usigma^{*i}, \utheta^i\right).\]

The composition of the string of $2$-cells $|\utheta^i|$ in each simplex $\Psi^i$ has the same input boxes as $\psi$.  So using finitely many elementary equivalences corresponding to the elementary relations \eqref{move:b1} and \eqref{move:b2}, we may assume that the wiring diagram $\psi$ has only one input box and that the simplices $\Psi^i$ have no $2$-cells.  At this stage, each stratified simplex $\Psi^i$ has the form 
\[
\left(\ulambda^i, \uomega^i, \usigma_*^i, \usigma^{*i}\right).\]

Observe that for each $i \in \{1,2\}$, the string of $1$-wasted wires $\uomega^i$ in the simplex $\Psi^i$ corresponds to precisely the set $\ewpsi \amalg \iwpsi$ of external and internal wasted wires in the wiring diagram $\psi$ (Def. \ref{def:wiring-diagram}).  Here an internal wasted wire in $\psi$ is created by applying a $1$-loop to a $1$-wasted wire as in \eqref{internal-wasted-wire}.  Therefore, using finitely many elementary equivalences corresponding to the elementary relations \eqref{move:c1} and \eqref{move:c4}, we may assume that the wiring diagram $\psi$ and the two simplices $\Psi^i$ have no $1$-wasted wires.   At this stage, each stratified simplex $\Psi^i$ has the form
\[
\left(\ulambda^i, \usigma_*^i, \usigma^{*i}\right).\]

Using finitely many elementary equivalences corresponding to the elementary relation \eqref{move:c5}, we may assume that each loop element in the wiring diagram $\psi$ (Def. \ref{def:loop-element}) corresponds to precisely one $1$-loop $\lambda$  in each simplex $\Psi^i$.  At this stage, the $1$-loops in each  simplex $\Psi^i$ are in bijection with the loop elements in $\psi$.  Moreover, the two stratified subsimplices $\left(\usigma_*^1, \usigma^{*1}\right) \subseteq \Psi^1$ and  $\left(\usigma_*^2, \usigma^{*2}\right) \subseteq \Psi^2$ are equivalent.  Therefore, using finitely many elementary equivalences corresponding to the elementary relation \eqref{move:c1}, we may assume that the wiring diagram $\psi$ has no loop elements and that the simplices $\Psi^i$ have no $1$-loops.  So each stratified simplex $\Psi^i$ now has the form
\[
\left(\usigma_*^i, \usigma^{*i}\right).\]

The two stratified subsimplices $\left(\usigma_*^i\right) \subseteq \Psi^i$ of in-splits for $i \in \{1,2\}$ are also equivalent.  They are connected by a finite sequence of elementary equivalences corresponding to the elementary relations \eqref{move:d1} and \eqref{move:d2}.  Likewise, the stratified subsimplices  $\left(\usigma^{*i}\right) \subseteq \Psi^i$ of out-splits for $i \in \{1,2\}$ are connected by a finite sequence of elementary equivalences corresponding to the elementary relations \eqref{move:e1} and \eqref{move:e2}.  So the two simplices $\Psi^i$ are also connected by a finite sequence of elementary equivalences.
\end{proof}

We are now ready for the finite presentation theorem for wiring diagrams.  It describes the wiring diagram operad $\WD$ (Theorem \ref{wd-operad}) in terms of finitely many generators and finitely many relations.

\begin{theorem}
\label{thm:wd-generator-relation}
Consider the operad $\WD$ of wiring diagrams.\index{finite presentation for WD@finite presentation for $\WD$}
\begin{enumerate}
\item
Every wiring diagram can be obtained from finitely many generating wiring diagrams (Def. \ref{def:generating-wiring-diagrams}) via iterated operadic compositions (Def. \ref{def:pseudo-operad}).
\item
Any two equivalent simplices are either equal or are connected by a finite sequence of elementary equivalences  (Def. \ref{def:equivalent-simplices}).
\end{enumerate}
\end{theorem}

\begin{proof}
The first statement is a special case of Theorem \ref{stratified-presentation-exists}.  The second statement is a combination of Remark \ref{rk:stratified-disjoint}, Lemma \ref{lemma:simplex-to-stratified} twice, Lemma \ref{lemma:stratified-type1}, and Lemma \ref{lemma:stratified-type2}.
\end{proof}

\section{Finite Presentation for Normal Wiring Diagrams}
\label{sec:coherence-normal}

In this section, we establish a finite presentation theorem for the operad of wiring diagrams without delay nodes.  Such wiring diagrams are used in \cite{spivak15,spivak15b} to study mode-dependent networks and dynamical systems.  Recall  Def. \ref{def:wiring-diagram}, Def. \ref{wd-equivalence}, and Convention \ref{conv:prewiring} regarding wiring diagrams.

\begin{definition}
\label{def:wd-without-dn}
Fix a class $S$.
\begin{enumerate}
\item
A wiring diagram is said to be \emph{normal} \index{normal wiring diagram} if its set of delay nodes is empty.
\item
The collection of normal wiring diagrams is denoted by \label{notation:normal-wd}$\wddot$.\index{WDdot@$\wddot$}  If we want to emphasize $S$, then we will write $\wddot^S$.
\end{enumerate}
\end{definition}

\begin{example}
\label{ex:generators-without-dn}
Among the $8$ generating wiring diagrams (section \ref{sec:generating-wd}):
\begin{enumerate}
\item
A $1$-delay node $\delta_d$ (Def. \ref{def:one-dn}) is not normal.
\item
The empty wiring diagram $\epsilon$ (Def. \ref{def:empty-wd}), a name change $\tau_{X,Y}$ (Def. \ref{def:name-change}), a $2$-cell $\theta_{X,Y}$ (Def. \ref{def:theta-wd}), a $1$-loop $\lambda_{X,x}$ (Def. \ref{def:loop-wd}), an in-split $\sigma_{X,x_1,x_2}$ (Def. \ref{def:insplit-wd}), an out-split $\sigma^{Y,y_1,y_2}$ (Def. \ref{def:out-split}), and a $1$-wasted wire $\omega_{Y,y}$ (Def. \ref{def:wasted-wire-wd}) are normal.
\end{enumerate}
In particular, there is a proper inclusion $\wddot \subsetneq \WD$.  Furthermore, the $1$-internal wasted wire $\omega^{X,x}$ (Def. \ref{def:internal-wasted-wire}) is normal.
\end{example}

\begin{example}
\label{ex:elementary-relation-without-dn}
All the wiring diagrams that appear in the $28$ elementary relations (section \ref{sec:elementary-relations}) are normal.
\end{example}

\begin{example}
\label{ex:decomp-without-dn}
Among the wiring diagrams in Chapter \ref{ch04-decomposition}:
\begin{enumerate}
\item
$\varphi$ \eqref{phi-wiring-diagram} is not normal, unless $r=0$.
\item
$\psi$ in \eqref{only-external-wasted-wires}, $\alpha$ \eqref{alpha-wiring-diagram}, $\pi_1$ \eqref{pi-one-supplier}, $\pi_2$ \eqref{pi-two-supplier}, $\beta_1$ \eqref{supply-betaone}, $\beta_2$ \eqref{supply-betatwo}, and $\beta_3$ \eqref{supply-betathree} are normal.
\end{enumerate}
\end{example}

\begin{proposition}
\label{prop:without-dn-operads}
With respect to 
\begin{itemize}
\item
the equivariant structure in Def. \ref{wd-equivariance}, 
\item
the colored units in Def. \ref{wd-units}, and 
\item
the $\compi$-compositions in Def. \ref{def:compi-wd}, 
\end{itemize}
$\wddot$ is a $\boxs$-colored operad, called the operad of normal wiring diagrams.\index{wddot is an operad@$\wddot$ is an operad}
\end{proposition}

\begin{proof}
We can reuse the proof of Theorem \ref{wd-operad}--that $\WD$ is a $\boxs$-colored operad--as long as we know that the relevant structure is well-defined in $\wddot$.

The collection $\wddot$ is closed under the equivariant structure map \eqref{wd-permutation}.  Furthermore, each colored unit $\tensorunit_Y$ \eqref{wd-unit} is in  $\wddot$.  

Suppose both $\varphi$ and $\psi$ are normal wiring diagrams such that $\phicompipsi \in \WD$ is defined.  Then $\phicompipsi$ is also normal because
\[
\dnphicompipsi = \dnphi \amalg \dnpsi =\varnothing.
\]
Therefore,  Lemmas \ref{compi-unity}, \ref{compi-horizontal-associative}, and \ref{compi-vertical-associative} all apply to $\wddot$ to show that it is an operad.
\end{proof}

Our next objective is to obtain a version of the finite presentation theorem for $\wddot$.  For this purpose, we will use the following definitions.

\begin{definition}
\label{def:generating-without-dn}
Consider the operad $\wddot$ of normal wiring diagrams.
\begin{enumerate}
\item
A \emph{normal generating wiring diagram} \index{normal generating wiring diagram} is a generating wiring diagram (Def. \ref{def:generating-wiring-diagrams}) except for $1$-delay nodes $\delta_d$ (Def. \ref{def:one-dn}).
\item
A \emph{normal simplex} \index{normal simplex} is defined as in Def. \ref{def:simplex} using normal generating wiring diagrams and $\wddot$ in place of $\WD$.
\item
A \emph{normal stratified simplex} \index{normal stratified simplex} and a \emph{normal stratified presentation} \index{normal stratified presentation} are defined as in Def. \ref{def:stratified-presentation} with $\wddot$ in place of $\WD$, except that a normal stratified simplex of type (2) has the form $\bigl(\tau, \ulambda, \uomega, \usigma_*, \usigma^*, \utheta\bigr)$.
\item
All of Def. \ref{def:equivalent-simplices} is repeated with normal generating wiring diagrams and $\wddot$ in place of $\WD$.
\end{enumerate}
\end{definition}

The following result is the finite presentation theorem for normal wiring diagrams.

\begin{theorem}
\label{thm:without-dn-coherence}
Consider the operad $\wddot$ of normal wiring diagrams.\index{finite presentation for WDdot@finite presentation for $\wddot$}
\begin{enumerate}
\item
Every normal wiring diagram has a normal stratified presentation.
\item
Every normal wiring diagram can be obtained from finitely many normal generating wiring diagrams via iterated operadic compositions (Def. \ref{def:pseudo-operad}).
\item
Any two equivalent normal simplices are connected by a finite sequence of  elementary equivalences  (Def. \ref{def:equivalent-simplices}).
\end{enumerate}
\end{theorem}

\begin{proof}
For statement (1), we reuse the proof of Theorem \ref{stratified-presentation-exists} while assuming $r=0$.  Statement (2) is a special case of statement (1).  

For statement (3) we reuse the proof of Theorem \ref{thm:wd-generator-relation}(2).  In other words, we simply reuse the proofs of Lemma \ref{lemma:simplex-to-stratified}, Lemma \ref{lemma:stratified-type1}, and Lemma \ref{lemma:stratified-type2} while assuming $r=0$ throughout.  The key observation is that, for normal simplices, elementary equivalences as in Def. \ref{def:equivalent-simplices} involve either:
\begin{itemize}
\item
elementary relations (Def. \ref{def:elementary-relations}), none of which involves delay nodes, or 
\item
an operad associativity or unity axiom--\eqref{compi-associativity}, \eqref{compi-associativity-two}, \eqref{compi-left-unity}, or \eqref{compi-right-unity}--for the normal generating wiring diagrams.  
\end{itemize}
So in the $\wddot$ versions of these Lemmas, we simply ignore all the delay nodes in the original proofs.
\end{proof}

\section{Finite Presentation for Strict Wiring Diagrams}
\label{sec:coherence-strict}

In this section, we establish a finite presentation theorem for the operad of strict wiring diagrams.  Such wiring diagrams are used in \cite{vsl} to study open dynamical systems.

\begin{definition}
\label{def:strict-wd}
Fix a class $S$.
\begin{enumerate}
\item
A wiring diagram (Def. \ref{wd-equivalence}) is said to be \emph{strict} \index{strict wiring diagram} if
\begin{enumerate}[(i)]
\item
it is normal (Def. \ref{def:wd-without-dn}) and 
\item
its supplier assignment is a bijection.
\end{enumerate}
\item
The collection of strict wiring diagrams is denoted by \label{notation:strict-wd}$\wdzero$.\index{WDzero@$\wdzero$}  If we want to emphasize $S$, then we will write $\wdzero^S$.
\end{enumerate}
\end{definition}

\begin{remark}
What we call a strict wiring diagram is simply called a wiring diagram in \cite{vsl} (Def. 3.5).  In \cite{vsl} (Remark 2.7) $S$ is a \emph{set} of representatives of isomorphism classes of second-countable smooth manifolds.  The non-instantaneity requirement \eqref{non-instant} in this case is called the \emph{no passing wires} \index{no passing wires} requirement in \cite{vsl}.  As noted in \cite{vsl} (Remark 3.6), strictness implies the non-existence of external wasted wires, internal wasted wires (Def. \ref{def:wiring-diagram}), and split wires, i.e., multiple (at least two) demand wires having the same supply wire.  So strict wiring diagrams are much simpler than a general wiring diagram.  
\end{remark}

\begin{example}
\label{ex:strict-generators}
Among the $8$ generating wiring diagrams (section \ref{sec:generating-wd}):
\begin{enumerate}
\item
A $1$-delay node $\delta_d$ (Def. \ref{def:one-dn}) is not normal (Def. \ref{def:wd-without-dn}), hence also not strict.
\item
The empty wiring diagram $\epsilon$ (Def. \ref{def:empty-wd}), a name change $\tau_{X,Y}$ (Def. \ref{def:name-change}), a $2$-cell $\theta_{X,Y}$ (Def. \ref{def:theta-wd}), and a $1$-loop $\lambda_{X,x}$ (Def. \ref{def:loop-wd}) are strict.
\item
An in-split $\sigma_{X,x_1,x_2}$ (Def. \ref{def:insplit-wd}), an out-split $\sigma^{Y,y_1,y_2}$ (Def. \ref{def:out-split}), and a $1$-wasted wire $\omega_{Y,y}$ (Def. \ref{def:wasted-wire-wd}) are normal but not strict.
\end{enumerate}
In particular, there are proper inclusions $\wdzero \subsetneq \wddot \subsetneq \WD$.   Furthermore, the $1$-internal wasted wire $\omega^{X,x}$ (Def. \ref{def:internal-wasted-wire}) is normal but not strict.
\end{example}

\begin{example}
\label{ex:strict-elementary-relation}
Among the wiring diagrams that appear in the $28$ elementary relations (section \ref{sec:elementary-relations}):
\begin{enumerate}
\item
\eqref{move:a1}, \eqref{move:a2}, \eqref{move:a3}, \eqref{move:b0}, \eqref{move:b1}, \eqref{move:b2}, \eqref{move:b3}, \eqref{move:c1}, \eqref{move:c6}, and \eqref{move:d5} are strict wiring diagrams.
\item
The other $18$ are normal but not strict.
\item
Only \eqref{move:a1}, \eqref{move:a2}, \eqref{move:a3}, \eqref{move:b0}, \eqref{move:b1}, \eqref{move:b2}, \eqref{move:b3}, and \eqref{move:c1} involve only strict wiring diagrams on both sides.  Indeed, both  \eqref{move:c6} and \eqref{move:d5} involve a $1$-wasted wire, which is not strict.
\end{enumerate}
\end{example}

\begin{example}
\label{ex:strict-decomp}
Among the wiring diagrams in Chapter \ref{ch04-decomposition}:
\begin{enumerate}
\item
$\psi$ in \eqref{only-external-wasted-wires} is normal but not strict.
\item
$\alpha$ \eqref{alpha-wiring-diagram}, $\pi_2$ \eqref{pi-two-supplier}, $\beta_1$ \eqref{supply-betaone}, $\beta_2$ \eqref{supply-betatwo}, and $\beta_3$ \eqref{supply-betathree} are normal but not strict in general.
\item
$\pi_1$ \eqref{pi-one-supplier} is a strict wiring diagram by Lemma \ref{lemma:iterated-loops}.
\end{enumerate}
\end{example}

\begin{proposition}
\label{prop:strict-operads}
With respect to 
\begin{itemize}
\item
the equivariant structure in Def. \ref{wd-equivariance}, 
\item
the colored units in Def. \ref{wd-units}, and 
\item
the $\compi$-compositions in Def. \ref{def:compi-wd}, 
\end{itemize}
$\wdzero$ is a $\boxs$-colored operad, called the operad of strict wiring diagrams.\index{wdzero is an operad@$\wdzero$ is an operad}
\end{proposition}

\begin{proof}
The argument is essentially identical to the proof of Prop. \ref{prop:without-dn-operads} with a minor modification.  Suppose both $\varphi$ and $\psi$ are strict wiring diagrams such that $\phicompipsi \in \WD$ is defined.  Then $\phicompipsi$ is also strict.  Indeed, we already know that it is normal.  Next, one can check directly from the definition of the supplier assignment $\sphicompipsi$ \eqref{compi-supply} that it is a bijection because, in all cases, it is defined as a composition of the bijections $\sphi$ and $\spsi$.  Therefore,  Lemmas \ref{compi-unity}, \ref{compi-horizontal-associative}, and \ref{compi-vertical-associative} all apply to $\wdzero$ to show that it is an operad.
\end{proof}

Our next objective is to obtain a version of the finite presentation theorem for $\wdzero$.  For this purpose, we will use the following definitions.

\begin{definition}
\label{def:generating-strict}
Consider the operad $\wdzero$ of strict wiring diagrams.
\begin{enumerate}
\item
A \emph{strict generating wiring diagram} \index{strict generating wiring diagram}  means the empty wiring diagram $\epsilon$ (Def. \ref{def:empty-wd}), a name change $\tau_{X,Y}$ (Def. \ref{def:name-change}), a $2$-cell $\theta_{X,Y}$ (Def. \ref{def:theta-wd}), or a $1$-loop $\lambda_{X,x}$ (Def. \ref{def:loop-wd}).
\item
A \emph{strict simplex} \index{strict simplex} is defined as in Def. \ref{def:simplex} using strict generating wiring diagrams and $\wdzero$ in place of $\WD$.
\item
A \emph{strict stratified simplex} \index{strict stratified simplex} is a stratified simplex (Def. \ref{def:stratified-presentation}) of the form $(\epsilon)$ or $\bigl(\tau, \ulambda, \utheta\bigr)$.
\item
If $\Psi$ is a strict stratified simplex, then we call it a \emph{strict stratified presentation} \index{strict stratified presentation} of the strict wiring diagram $|\Psi|$.
\item
A \emph{strict elementary relation} \index{strict elementary relation} means one of the $8$ elementary relations that involve only strict wiring diagrams on both sides, namely, \eqref{move:a1}, \eqref{move:a2}, \eqref{move:a3}, \eqref{move:b0}, \eqref{move:b1}, \eqref{move:b2}, \eqref{move:b3}, and  \eqref{move:c1}.  See Example \ref{ex:strict-elementary-relation}.
\item
All of Def. \ref{def:equivalent-simplices} is repeated with $\wdzero$ in place of $\WD$ using strict generating wiring diagrams, strict simplices, and strict elementary relations.  The resulting notions are called \emph{strict elementary equivalences},\index{strict elementary equivalences} and so forth.
\end{enumerate}
\end{definition}

The following result is the finite presentation theorem for strict wiring diagrams.

\begin{theorem}
\label{thm:strict-wd-coherence}
Consider the operad $\wdzero$ of strict wiring diagrams.\index{finite presentation for WDzero@finite presentation for $\wdzero$}
\begin{enumerate}
\item
Every strict wiring diagram has a strict stratified presentation.
\item
Every strict wiring diagram can be obtained from finitely many strict generating wiring diagrams via iterated operadic compositions (Def. \ref{def:pseudo-operad}).
\item
Any two equivalent strict simplices are connected by a finite sequence of strict elementary equivalences.
\end{enumerate}
\end{theorem}

\begin{proof}
As in the proof of Theorem \ref{thm:without-dn-coherence}, for statement (1), we reuse the proof of Theorem \ref{stratified-presentation-exists} while assuming the wiring diagram $\psi$ is strict.  In this case, $\psi$ is either the empty wiring diagram $\epsilon$ or has a decomposition (using \eqref{psi-is-alpha-comp-phi} and \eqref{pi-pione-pitwo})
\[
\psi = \pi_1 \comp \pi_2 \comp \varphi
\]
in which $\pi_2$ \eqref{pi-two-supplier} is a name change.  The desired strict stratified presentation follows from Convention \ref{conv:ignore-name-change} and the facts that: 
\begin{itemize}
\item
$\pi_1$ is either a colored unit or has a stratified presentation $(\ulambda)$ by Lemma \ref{lemma:iterated-loops};
\item
$\varphi$ is either a colored unit or has a stratified presentation $(\utheta)$ by Lemma \ref{lemma:phi-presentation}.
\end{itemize}

Statement (2) is a special case of statement (1).  

For statement (3) we use the strict versions of the proofs of Lemma \ref{lemma:simplex-to-stratified}, Lemma \ref{lemma:stratified-type1}, and Lemma \ref{lemma:stratified-type2}.  The key observation is that, in this case, only strict generating wiring diagrams and strict elementary equivalences are used in these proofs.
\end{proof}

\section{Summary of Chapter \ref{ch05-stratified-presentation}}

\begin{enumerate}
\item A simplex in $\WD$ is a finite non-empty parenthesized word of generating wiring diagrams in which each pair of parentheses is equipped with an operadic $\compi$-composition.
\item A stratified simplex in $\WD$ is a simplex of one of the following two forms.
\begin{itemize}
\item
$\bigl(\uomega, \epsilon\bigr)$
\item
$\Bigl(\tau, \ulambda, \uomega, \usigma_*, \usigma^*, \utheta, \udelta\Bigr)$
\end{itemize}
\item Every wiring diagram has a stratified presentation.
\item Two simplices are elementarily equivalent if one can be obtained from the other by replacing a subsimplex $\Psi'$ by an equivalent simplex $\Psi''$ such that $|\Psi'| = |\Psi''|$ is either one of the twenty-eight elementary relations in $\WD$ or an operad associativity/unity axiom involving only the eight generating wiring diagrams.
\item Any two simplex presentations of a given wiring diagram are connected by a finite sequence of elementary equivalences.
\item A normal wiring diagram is a wiring diagram with no delay nodes.
\item The operad $\wddot$ of normal wiring diagrams satisfies a finite presentation theorem involving the seven normal generating wiring diagrams and the twenty-eight elementary relations.
\item A strict wiring diagram is a wiring diagram with no delay nodes and whose supplier assignment is a bijection.
\item The operad $\wdzero$ of strict wiring diagrams satisfies a finite presentation theorem involving the four strict generating wiring diagrams and the eight strict elementary relations.
\end{enumerate}

\chapter{Finite Presentation for Algebras over Wiring Diagrams}
\label{ch06-wd-algebras}

The main purpose of this chapter is to provide finite presentations for algebras over the operads $\WD$ (Theorem \ref{wd-operad}), $\wddot$ (Prop. \ref{prop:without-dn-operads}), and $\wdzero$ (Prop. \ref{prop:strict-operads}).  The advantage of such a finite presentation is that sometimes the general structure map of an operad algebra can be a bit difficult to write down and understand.  On the other hand, our generating structure maps and generating axioms are all fairly easy to write down and understand, as we will illustrate with examples in Sections \ref{sec:propagator}, \ref{sec:discrete-systems}, and \ref{sec:algebra-ods}.

In Section \ref{sec:operad-algebras} we recall the basics of algebras over an operad.

In Section \ref{sec:algebras-wd} we first define a $\WD$-algebra in terms of $8$ generating structure maps and $28$ generating axioms corresponding to the generating wiring diagrams (Def. \ref{def:generating-wiring-diagrams})  and the elementary relations (Def. \ref{def:elementary-relations}), respectively.  Then we observe that this finite presentation for a $\WD$-algebra is in fact equivalent to the general definition of a $\WD$-algebra (Theorem \ref{thm:wd-algebra}).   This is an application of the finite presentation theorem for the operad $\WD$ (Theorem \ref{thm:wd-generator-relation}).

In Section \ref{sec:propagator} we provide a finite presentation for the $\WD$-algebra called the propagator algebra.  In its original form, the propagator algebra was the main example in \cite{rupel-spivak}. 

In Section \ref{sec:algebra-normal-wd} we observe that algebras over the operad $\wddot$ of normal wiring diagrams have a similar finite presentation with $7$ generating structure maps and $28$ generating axioms.  In Section \ref{sec:discrete-systems} we provide a finite presentation for the $\wddot$-algebra called the algebra of discrete systems.  In its original form, this algebra was one of the main examples in \cite{spivak15b}.   

In Section \ref{sec:algebra-strict-wd} we observe that algebras over the operad $\wdzero$ of strict wiring diagrams admit a finite presentation with $4$ generating structure maps and $8$ generating axioms.  In Section \ref{sec:algebra-ods} we provide a finite presentation for the $\wdzero$-algebra called the algebra of open dynamical systems.  In its original form, this algebra was one of the main examples in \cite{vsl}.

\section{Operad Algebras}
\label{sec:operad-algebras}

Let us first recall the definition of an algebra over an operad.  The following definition can be found in \cite{yau-operad} (Def. 13.2.3).  In its original $1$-colored topological form, it was first given in \cite{may72}.  Fix a class $S$ as before.

\begin{motivation}
One can think of an algebra over an operad as a generalization of a module over a ring.  Given a ring $R$, a left $R$-module $M$ is equipped with structure maps $r : M \to M$ for each element $r \in R$ that satisfy some axioms.  In particular, these structure maps are associative in the sense that
\[(r_1r_2)(m) = r_1(r_2m)\]
for $r_1,r_2 \in R$ and $m \in M$.  Furthermore, the multiplicative unit $1_R$ of $R$ acts as the identity map, so $1_R(m) = m$.  For algebras over an operad, there is also an equivariance axiom because operads can model operations with multiple inputs.
\end{motivation}

\begin{definition}
\label{def1:operad-algebra}
Suppose $(\O,\tensorunit,\gamma)$ is an $S$-colored operad as in Def. \ref{def:colored-operad}.  An \emph{$\O$-algebra} \index{operad algebra} is a pair $(A,\mu)$ consisting of the following data.
\begin{enumerate}
\item
For each $c \in S$, $A$ is equipped with a class \label{notation:a-sub-c}$A_c$ called the \emph{$c$-colored entry} of $A$.\index{colored entries of an algebra}
\item
For each $d \in S$, $\uc = (c_1,\ldots,c_n) \in \profs$, and $\zeta \in \O\duc$, $A$ is equipped with a \emph{structure map}\label{notation:algebra-structure-map} \index{structure map of an algebra}
\begin{equation}
\label{operad-structure-map}
\nicexy{
A_{\uc} \defn \prod\limits_{i=1}^n A_{c_i} \ar[r]^-{\mu_{\zeta}} & A_d}
\end{equation}
in which an empty product, for the case $n=0$, means the one-point set $\{*\}$.
\end{enumerate}
This data is required to satisfy the following associativity, unity, and equivariance axioms.
\begin{description}
\item[Associativity]
Suppose $\duc \in \profs \times S$ is as above with $n \geq 1$, $\zeta_0 \in \O\duc$, $\zeta_i \in \O\ciubi$ for each $1 \leq i \leq n$ with $\ub_i \in \profs$, and $\ub = \left(\ub_1,\ldots,\ub_n\right)$ as in \eqref{operadic-composition}.  Write $\zeta = \gamma\bigl(\zeta_0; \zeta_1,\ldots,\zeta_n\bigr) \in \O\dub$.  Then the diagram\index{associativity of an algebra}
\begin{equation}
\label{operad-algebra-associativity}
\nicexy{A_{\ub} \ar[r]^-{\mu_{\zeta}} \ar[d]_-{=} & A_d\\
A_{\ub_1} \times \cdots \times A_{\ub_n} \ar[r]^-{\prod \mu_{\zeta_i}} 
& A_{c_1} \times \cdots \times A_{c_n} = A_{\uc} \ar[u]_-{\mu_{\zeta_0}}}
\end{equation}
is commutative.
\item[Unity]
For each $c \in S$, the structure map\index{unity of an algebra}
\begin{equation}
\label{operad-algebra-unity}
\nicexy{
A_c \ar[r]^-{\mu_{\tensorunit_c}} & A_c}
\end{equation}
is the identity map, where $\tensorunit_c \in \O\ccsingle$ is the $c$-colored unit of $\O$.
\item[Equivariance]
Suppose $\zeta \in \O\duc$ as in \eqref{operad-structure-map}, $\sigma \in \Sigma_n$, and $\zeta^{\sigma} \in \O\ducsigma$ is the image of $\zeta$ under the right action \eqref{operad-right-equivariance}.  Then the diagram\index{equivariance of an algebra}
\begin{equation}
\label{operad-algebra-eq}
\nicexy{
A_{\uc} \ar[d]_-{\mu_\zeta} \ar[r]^-{\sigma^{-1}} & A_{\uc\sigma} \ar[d]^-{\mu_{\zeta^{\sigma}}}\\
A_d \ar[r]^-{=} & A_d}
\end{equation}
is commutative.  Here the top $\sigma^{-1}$ permutes the factors of $A_{\uc}$ from the left, and  $\uc\sigma = \left(c_{\sigma(1)}, \ldots, c_{\sigma(n)}\right)$. 
\end{description}
To simplify the notations, we will sometimes denote an $\O$-algebra by just $A$ and denote the structure map $\mu_{\zeta}$ by $\zeta$.
\end{definition}

Just as operads can be equivalently expressed in terms of the $\compi$-compositions (Prop. \ref{prop:operad-def-equiv}), so can operad algebras.

\begin{definition}
\label{def2:operad-algebra}
Suppose $(\O,\tensorunit,\comp)$ is an $S$-colored operad as in Def. \ref{def:pseudo-operad}.  An \emph{$\O$-algebra} $(A,\mu)$ is defined exactly as in Def. \ref{def1:operad-algebra} except that the associativity axiom \eqref{operad-algebra-associativity} is replaced by the following axiom.
\begin{description}
\item[Associativity]
Suppose:
\begin{itemize}
\item
$d \in S$, $\uc = (c_1, \ldots, c_n) \in \profs$ with $n \geq 1$, and $1 \leq i \leq n$;
\item
$\ub \in \profs$ and $\uc \compi \ub$ as in \eqref{compi-profile};
\item
$\zeta \in \O\duc$, $\xi \in \O\ciub$, and $\zeta \compi \xi \in \O\dccompib$.
\end{itemize}
Then the diagram
\begin{equation}
\label{operad-algebra-associativity2}
\nicexy@C+.4cm{
A_{\uc \compi \ub} \ar[d]_-{=} \ar[r]^-{\mu_{\zeta \compi \xi}} & A_d\\
A_{(c_1,\ldots,c_{i-1})} \times A_{\ub} \times A_{(c_{i+1},\ldots,c_n)} 
\ar[r]^-{(\Id, \mu_\xi, \Id)} & A_{(c_1,\ldots,c_{i-1})} \times A_{c_i} \times A_{(c_{i+1},\ldots,c_n)}  = A_{\uc} \ar[u]_-{\mu_\zeta}}
\end{equation}
is commutative.
\end{description}
\end{definition}

\begin{notation}
\label{notation:algebra-compi}
To simplify the notations, we will sometimes denote the structure map $\mu_{\zeta}$ by $\zeta$.  We will also write the composition in the diagram \eqref{operad-algebra-associativity2} as $\mu_{\zeta} \compi \mu_{\xi}$, called the \emph{$\compi$-composition} of $\mu_{\zeta}$ and $\mu_{\xi}$.   So this associativity axiom states that\label{notation:compi-algebra-map} \index{compi-composition of@$\compi$-composition of structure maps}
\[\mu_{\zeta \compi \xi} = \mu_{\zeta} \compi \mu_{\xi}.\]
In other words, the structure map of the $\compi$-composition $\zeta \compi \xi$ is the $\compi$-composition of the structure maps corresponding to $\zeta$ and $\xi$.
\end{notation}

Using the associativity and the unity axioms in Def. \ref{def:colored-operad} and Def. \ref{def:pseudo-operad}, it is not hard to check that  Def. \ref{def1:operad-algebra} and Def. \ref{def2:operad-algebra} are in fact equivalent.  A proof of this equivalence can be found in \cite{yau-operad} (Prop. 16.7.8 and Exercises 10 and 11 in Chapter 16).

\begin{remark}\label{rk:operad-algebra-in-sets}
In  Def. \ref{def1:operad-algebra} and Def. \ref{def2:operad-algebra}, each entry of an operad algebra has no structure beyond being a class.  We could have just as easily defined operad algebras in a symmetric monoidal category, which is in fact the setting in \cite{yau-operad} (Def. 13.2.3).  One simply replaces the product with the symmetric monoidal product and the one-point set with the monoidal unit.  However, to keep the presentation simple, we chose to define operad algebras whose entries are just classes.  This is sufficient for the main examples in Sections \ref{sec:propagator}, \ref{sec:discrete-systems}, and \ref{sec:algebra-ods}.
\end{remark}

\begin{example}[Monoid Modules]
Suppose $(A,\mu,1)$ is a monoid (Example \ref{ex:monoid}).  Recall that a \index{module over a monoid} \emph{left $A$-module} is a set $M$ equipped with a structure map $a : M \to M$ for each $a \in A$ that is associative and unital.  Associativity means $(ab)m = a(bm)$ for $a,b \in A$ and $m \in M$.  Unity means $1m = m$ for $m \in M$.  If we regard $A$ as a $1$-colored operad $\sO$ concentrated in arity $1$ as in Example \ref{ex:monoid}, then left $A$-modules yield $\sO$-algebras in the sense of Def. \ref{def1:operad-algebra}.  The only slight difference between a left $A$-module and an $\sO$-algebra is that the former has an underlying set, while the latter is allowed to have an underlying class.
\end{example}

\begin{example}[Associative and Commutative Monoids]
This is a continuation of Example \ref{ex:as-com}.
\begin{enumerate}
\item For the associative operad $\As$, an $\As$-algebra with an underlying set is exactly a monoid.
\item For the commutative operad $\Com$, a $\Com$-algebra with an underlying set is exactly a monoid whose multiplication is commutative.
\end{enumerate}
\end{example}

\begin{example}[Traffic Spaces and Probability Spaces]\label{ex:traffic}
This is a continuation of Example \ref{ex:graph-op}, where we discussed the graph operation operad $\GrOp$.  Here we consider $\GrOp$-algebras in the category of complex vector spaces; see Remark \ref{rk:operad-algebra-in-sets}.  In particular, a $\GrOp$-algebra $A$ has an underlying complex vector space, and all the structure maps are linear maps, with tensor products playing the roles of products.  The graph operation $\bullet \in \GrOp_0$, consisting of a single vertex and no edges, yields an element in $A$, also denoted by $\bullet$.  For a general graph operation $G \in \GrOp_n$, the corresponding structure map $A^{\otimes n} \to A$ is also denoted by $G$.  We will write $\delta \in \GrOp_1$ for the graph operation consisting of a single vertex and a loop.

A \index{traffic space}\emph{traffic space} \cite{male} is a pair $(A,\varphi)$ in which 
$A$ is a $\GrOp$-algebra and $\varphi : A \to \bC$ is a linear functional such that the following two conditions are satisfied.
\begin{description}
\item[Unity and Diagonality] $\varphi(\bullet) = 1$ and $\varphi = \varphi \circ \delta$.
\item[Input-Independence] For each graph operation $G \in \GrOp_n$, the graph operation $\delta(G) \in \GrOp_n$ is obtained from $G$ by identifying its input and output.  Suppose $G'$ is a graph operation obtained from $\delta(G)$ by choosing a different vertex as the input/output.  Then $\varphi \circ \delta(G) = \varphi \circ G'$.
\end{description}
For example, suppose $G \in \GrOp_4$ is the graph operation on the left:
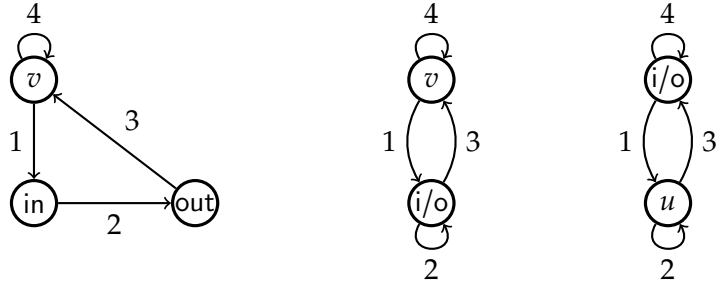
\begin{center}
\begin{tikzpicture}
\matrix[row sep=1cm, column sep=1.5cm]{
\node [plain] (v) {$v$}; &&[1cm]   \node [plain] (v2) {$v$}; &[1cm]  \node [plain] (v3) {$\mathsf{i/o}$};\\
\node [plain] (in) {$\inp$}; &\node [plain] (out) {$\out$}; & \node [plain] (io2) {$\mathsf{i/o}$}; & \node [plain] (io3) {$u$};\\};
\draw [arrow] (v) to node[swap]{$1$} (in);
\draw [arrow] (in) to node[swap]{$2$} (out);
\draw [arrow] (out) to node[swap]{$3$} (v);
\draw [arrow, out=120, in=60, looseness=4] (v) to node{$4$} (v);
\draw [arrow, bend right] (v2) to node[swap]{$1$} (io2);
\draw [arrow, out=240, in=300, looseness=4] (io2) to node[swap]{$2$} (io2);
\draw [arrow, bend right] (io2) to node[swap]{$3$} (v2);
\draw [arrow, out=120, in=60, looseness=4] (v2) to node{$4$} (v2);
\draw [arrow, bend right] (v3) to node[swap]{$1$} (io3);
\draw [arrow, out=240, in=300, looseness=4] (io3) to node[swap]{$2$} (io3);
\draw [arrow, bend right] (io3) to node[swap]{$3$} (v3);
\draw [arrow, out=120, in=60, looseness=4] (v3) to node{$4$} (v3);
\end{tikzpicture}
\end{center}
Then $\delta(G)$ is the graph operation in the middle, and the graph operation on the right is an example of a $G'$.

Traffic spaces play an important role in (non-commutative) probability theory.   Indeed, a \index{non-commutative probability space}\emph{non-commutative probability space}, also known as a quantum probability space, is a pair $(A,\varphi)$ in which:
\begin{enumerate}
\item $A$ is a unital $\bC$-algebra.
\item $\varphi : A \to \bC$ is a unital linear functional.
\end{enumerate}
A \emph{$*$-probability space} is a non-commutative probability space $(A,\varphi)$ in which:
\begin{enumerate}
\item $A$ is equipped with an anti-linear involution $*$ such that $(ab)^* = b^*a^*$ for all $a,b \in A$.
\item $\varphi$ satisfies the positivity condition that $\varphi(a^*a) \geq 0$ for all $a \in A$.
\end{enumerate}
Then a commutative $*$-probability space is an example of a traffic space, since the product of $n$ elements in $A$ is well-defined and is independent of the order of those elements.  It is, furthermore, a $*$-algebra in the following sense.  For each graph operation $G \in \GrOp_n$, its transpose $G^t$ is the graph operation obtained from $G$ by reversing the direction of each edge and swapping the input and the output.  Then 
\[\bigl(G(a_1 \otimes \cdots \otimes a_n)\bigr)^* = G^t(a_1^* \otimes \cdots \otimes a_n^*)\]
for all $a_1, \ldots, a_n \in A$.
\end{example}

\section{Algebras over the Operad of Wiring Diagrams}
\label{sec:algebras-wd}

The purpose of this section is to provide a finite presentation for $\WD$-algebras.  We begin by defining a $\WD$-algebra in terms of the generating wiring diagrams and the elementary relations.  Immediately afterwards we will establish its equivalence with Def. \ref{def2:operad-algebra} when $\O = \WD$.  Recall that $\WD$ is a $\boxs$-colored operad (Theorem \ref{wd-operad}).

\begin{definition}
\label{def:wd-algebra}
A \emph{$\WD$-algebra} \index{WD-algebra@$\WD$-algebra} $A$ consists of the following data.  For each $X \in \boxs$, $A$ is equipped with a class $A_X$ called the \emph{$X$-colored entry} of $A$.  It is equipped with the following $8$ \emph{generating structure maps} corresponding to the generating wiring diagrams (Def. \ref{def:generating-wiring-diagrams}).
\begin{enumerate}
\item
Corresponding to the empty wiring diagram $\epsilon \in \WD\emptynothing$ (Def. \ref{def:empty-wd}), it has a structure map
\begin{equation}
\label{structuremap-empty}
\nicexy{\ast \ar[r]^-{\epsilon} & A_{\varnothing}},
\end{equation}
i.e., a chosen element in $A_{\varnothing}$.
\item
Corresponding to each $1$-delay node $\delta_d \in \WD\dnothing$ (Def. \ref{def:one-dn}), it has a structure map
\begin{equation}
\label{structuremap-delaynode}
\nicexy{\ast \ar[r]^-{\delta_d} & A_d},
\end{equation}
i.e., a chosen element in $A_d$.
\item
Corresponding to each name change $\tau_{X,Y} \in \WD\yx$ (Def. \ref{def:name-change}), it has a structure map
\begin{equation}
\label{structuremap-namechange}
\nicexy{A_X \ar[r]^-{\tau_{X,Y}} & A_Y}
\end{equation}
that is, furthermore, the identity map if $\tau_{X,X}$ is the colored unit $\tensorunit_X$ \eqref{wd-unit}.
\item
Corresponding to each $2$-cell $\theta_{X,Y} \in \WD\xplusyxy$ (Def. \ref{def:theta-wd}), it has a structure map
\begin{equation}
\label{structuremap-2cell}
\nicexy{A_X \times A_Y \ar[r]^-{\theta_{X,Y}} & A_{X \amalg Y}}.
\end{equation}
\item
Corresponding to each $1$-loop $\lambda_{X,x} \in \WD\xminusxx$ (Def. \ref{def:loop-wd}),  it has a structure map
\begin{equation}
\label{structuremap-1loop}
\nicexy{A_X  \ar[r]^-{\lambda_{X,x}} & A_{X \setminus x}}.
\end{equation}
\item
Corresponding to each in-split $\sigma_{X,x_1,x_2} \in \WD\yx$ (Def. \ref{def:insplit-wd}), it has a structure map
\begin{equation}
\label{structuremap-insplit}
\nicexy{A_X  \ar[r]^-{\sigma_{X,x_1,x_2}} & A_Y}.
\end{equation}
\item
Corresponding to each out-split  $\sigma^{Y,y_1,y_2} \in \WD\yx$ (Def. \ref{def:out-split}),  it has a structure map
\begin{equation}
\label{structuremap-outsplit}
\nicexy{A_X  \ar[r]^-{\sigma^{Y,y_1,y_2}} & A_Y}.
\end{equation}
\item
Corresponding to each $1$-wasted wire  $\omega_{Y,y} \in \WD\yx$ (Def. \ref{def:wasted-wire-wd}),  it has a structure map
\begin{equation}
\label{structuremap-wastedwire}
\nicexy{A_X  \ar[r]^-{\omega_{Y,y}} & A_Y}.
\end{equation}
\end{enumerate}
The following $28$ diagrams, called the \emph{generating axioms}, which correspond to the  elementary relations  (Def. \ref{def:elementary-relations}), are required to be commutative.
\begin{enumerate}
\item
In the setting of \eqref{move:a1}, the diagram
\begin{equation}
\label{axiom-a1}
\nicexy{
A_X \ar[r]^-{\tau_{X,Y}} \ar[dr]_-{\tau_{X,Z}} & A_Y \ar[d]^-{\tau_{Y,Z}}\\
& A_Z}
\end{equation}
is commutative.
\item
In the setting of \eqref{move:a2}, the diagram
\begin{equation}
\label{axiom-a2}
\nicexy{
A_X \times A_Y \ar[r]^-{\theta_{X,Y}} \ar[d]_-{(\tau_{X,X'},\tau_{Y,Y'})} 
& A_{X \amalg Y} \ar[d]^-{\tau_{X \amalg Y, X' \amalg Y'}}\\
A_{X'} \times A_{Y'} \ar[r]^-{\theta_{X',Y'}} & A_{X' \amalg Y'}}
\end{equation}
is commutative.
\item
In the setting of \eqref{move:a3}, the diagram
\begin{equation}
\label{axiom-a3}
\nicexy{
A_X  \ar[r]^-{\lambda_{X,x}} \ar[d]_-{\tau_{X,Y}} 
& A_{X \setminus x} \ar[d]^-{\tau_{X \setminus x, Y \setminus y}}\\
A_{Y} \ar[r]^-{\lambda_{Y,y}} & A_{Y \setminus y}}
\end{equation}
is commutative.
\item
In the setting of \eqref{move:a4}, the diagram
\[
\nicexy{
A_X  \ar[r]^-{\sigma_{X,x_1,x_2}} \ar[d]_-{\tau_{X,Y}} 
& A_{X'} \ar[d]^-{\tau_{X',Y'}}\\
A_{Y} \ar[r]^-{\sigma_{Y,y_1,y_2}} & A_{Y'}
}\]
is commutative.
\item
In the setting of \eqref{move:a5}, the diagram
\[
\nicexy{
A_{X'}  \ar[r]^-{\sigma^{X,x_1,x_2}} \ar[d]_-{\tau_{X',Y'}} 
& A_X \ar[d]^-{\tau_{X,Y}}\\
A_{Y'} \ar[r]^-{\sigma^{Y,y_1,y_2}} & A_{Y}
}\]
is commutative.
\item
In the setting of \eqref{move:a6}, the diagram
\[
\nicexy{
A_{X'}  \ar[r]^-{\omega_{X,x}} \ar[d]_-{\tau_{X',Y'}} 
& A_X \ar[d]^-{\tau_{X,Y}}\\
A_{Y'} \ar[r]^-{\omega_{Y,y}} & A_{Y}
}\]
is commutative.
\item
In the setting of \eqref{move:b0}, the diagram
\begin{equation}
\label{axiom-b0}
\nicexy{
A_{X} \times \ast  \ar[r]^-{(\Id,\epsilon)} \ar[d]_-{=} 
& A_X \times A_{\varnothing} \ar[d]^-{\theta_{X,\varnothing}}\\
A_{X} \ar[r]^-{=} & A_{X \amalg \varnothing}}
\end{equation}
is commutative.
\item
In the setting of \eqref{move:b1}, the diagram
\begin{equation}
\label{wd-algebra-2cell-associativity}
\nicexy@C+.4cm{
A_{X} \times A_Y \times A_Z  \ar[r]^-{(\Id,\theta_{Y,Z})} \ar[d]_-{(\theta_{X,Y},\Id)} 
& A_X \times A_{Y \amalg Z} \ar[d]^-{\theta_{X,Y \amalg Z}}\\
A_{X \amalg Y} \times A_Z \ar[r]^-{\theta_{X \amalg Y, Z}} & A_{X \amalg Y \amalg Z}
}
\end{equation}
is commutative.
\item
In the setting of \eqref{move:b2}, the diagram
\begin{equation}
\label{wd-algebra-equivariance}
\nicexy@C+.3cm{
A_Y \times A_X  \ar[r]^-{\text{permute}} \ar[d]_-{\theta_{Y,X}} 
& A_X \times A_Y \ar[d]^-{\theta_{X,Y}} \\
A_{Y \amalg X} \ar[r]^-{=} & A_{X \amalg Y}
}
\end{equation}
is commutative.
\item
In the setting of \eqref{move:b3}, the diagram
\begin{equation}
\label{axiom-b3}
\nicexy{
A_X \times A_Y \ar[r]^-{\theta_{X,Y}} \ar[d]_{(\lambda_{X,x},\Id)} 
& A_{X \amalg Y} \ar[d]^-{\lambda_{X\amalg Y, x}} \\
A_{X \setminus x} \times A_Y \ar[r]^-{\theta_{X\setminus x, Y}} & A_{(X \amalg Y) \setminus \{x\}}}
\end{equation}
is commutative.
\item
In the setting of \eqref{move:b4}, the diagram
\[
\nicexy{
A_X \times A_Y \ar[r]^-{\theta_{X,Y}} \ar[d]_{(\sigma_{X,x_1,x_2},\Id)} 
& A_{X \amalg Y} \ar[d]^-{\sigma_{X\amalg Y, x_1,x_2}} \\
A_{X'} \times A_Y \ar[r]^-{\theta_{X',Y}} & A_{X' \amalg Y}
}\]
is commutative.
\item
In the setting of \eqref{move:b5}, the diagram
\[
\nicexy{
A_{X'} \times A_Y \ar[r]^-{\theta_{X',Y}} \ar[d]_{(\sigma^{X,x_1,x_2},\Id)} 
& A_{X' \amalg Y} \ar[d]^-{\sigma^{X\amalg Y, x_1,x_2}} \\
A_X \times A_Y \ar[r]^-{\theta_{X,Y}} & A_{X \amalg Y}
}\]
is commutative.
\item
In the setting of \eqref{move:b6}, the diagram
\[
\nicexy{
A_{X'} \times A_Y \ar[r]^-{\theta_{X',Y}} \ar[d]_{(\omega_{X,x_0},\Id)} 
& A_{X' \amalg Y} \ar[d]^-{\omega_{X\amalg Y, x_0}} \\
A_X \times A_Y \ar[r]^-{\theta_{X,Y}} & A_{X \amalg Y}
}\]
is commutative.
\item
In the setting of \eqref{move:c1}, the diagram
\begin{equation}
\label{wd-algebra-doubleloop}
\nicexy@C+.4cm{
A_X \ar[r]^-{\lambda_{X,x^2}} \ar[d]_{\lambda_{X,x^1}} 
& A_{X \setminus x^2} \ar[d]^-{\lambda_{X\setminus x^2, x^1}} \\
A_{X \setminus x^1} \ar[r]^-{\lambda_{X \setminus x^1, x^2}} & A_{X \setminus x}
}
\end{equation}
is commutative.
\item
In the setting of \eqref{move:c2}, the diagram
\[
\nicexy{
A_X \ar[r]^-{\lambda_{X,x}} \ar[d]_{\sigma_{X,x_1,x_2}} 
& A_{X \setminus x} \ar[d]^-{\sigma_{X\setminus x, x_1, x_2}} \\
A_{X'} \ar[r]^-{\lambda_{X',x}} & A_{X' \setminus x}
}\]
is commutative.
\item
In the setting of \eqref{move:c3}, the diagram
\[
\nicexy@C+.5cm{
A_{X'} \ar[r]^-{\sigma^{X,x_1,x_2}} \ar[d]_{\lambda_{X',x}} 
& A_X \ar[d]^-{\lambda_{X,x}} \\
A_{X'\setminus x} \ar[r]^-{\sigma^{X\setminus x, x_1, x_2}} & A_{X \setminus x}
}\]
is commutative.
\item
In the setting of \eqref{move:c4}, the diagram
\[
\nicexy@C+.5cm{
A_{X'} \ar[r]^-{\omega_{X,x_0}} \ar[d]_{\lambda_{X',x}} 
& A_X \ar[d]^-{\lambda_{X,x}} \\
A_{X'\setminus x} \ar[r]^-{\omega_{X\setminus x, x_0}} & A_{X \setminus x}
}\]
is commutative.
\item
In the setting of \eqref{move:c5}, the diagram
\begin{equation}
\label{wd-algebra-loopelement}
\nicexy{
A_{X} \ar[r]^-{\sigma^{Y,x^1,x^2}} \ar[d]_{\sigma_{X,x_1,x_2}} 
& A_Y \ar[r]^-{\lambda_{Y,x(1)}}
& A_{Y \setminus x(1)} \ar[d]^-{\lambda_{Y \setminus x(1), x(2)}} \\
A_{X'} \ar[rr]^-{\lambda_{X',x}} && A_{X^*}
}
\end{equation}
is commutative.
\item
In the setting of \eqref{move:c6}, the diagram
\[
\nicexy{
A_{X} \ar[r]^-{=} \ar[d]_{\sigma^{Y,x_1,x_2}} 
& A_X \\
A_Y \ar[r]^-{\omega_{Z,x_1}} & A_{Z} \ar[u]_-{\lambda_{Z,x}}
}\]
is commutative.
\item
In the setting of \eqref{move:d1}, the diagram
\[
\nicexy@C+.5cm{
A_{X} \ar[r]^-{\sigma_{X,x_2,x_3}} \ar[d]_{\sigma_{X,x_1,x_2}} 
& A_{X_{23}} \ar[d]^-{\sigma_{X_{23},x_1,x_{23}}} \\
A_{X_{12}} \ar[r]^-{\sigma_{X_{12},x_{12},x_3}} & A_Y
}\]
is commutative.
\item
In the setting of \eqref{move:d2}, the diagram
\[
\nicexy@C+.5cm{
A_{X} \ar[r]^-{\sigma_{X,x_3,x_4}} \ar[d]_{\sigma_{X,x_1,x_2}} 
& A_{X_{34}} \ar[d]^-{\sigma_{X_{34},x_1,x_2}} \\
A_{X_{12}} \ar[r]^-{\sigma_{X_{12},x_3,x_4}} & A_Y
}\]
is commutative.
\item
In the setting of \eqref{move:d3}, the diagram
\[
\nicexy@C+.5cm{
A_{X} \ar[r]^-{\sigma^{Z,z^1,z^2}} \ar[d]_{\sigma_{X,z_1,z_2}} 
& A_Z \ar[d]^-{\sigma_{Z,z_1,z_2}} \\
A_W \ar[r]^-{\sigma^{Y,z^1,z^2}} & A_Y
}\]
is commutative.
\item
In the setting of \eqref{move:d4}, the diagram
\[
\nicexy{
A_{X} \ar[r]^-{\omega_{Z,z}} \ar[d]_{\sigma_{X,z_1,z_2}} 
& A_Z \ar[d]^-{\sigma_{Z,z_1,z_2}} \\
A_W \ar[r]^-{\omega_{Y,z}} & A_Y
}\]
is commutative.
\item
In the setting of \eqref{move:d5}, the diagram
\[
\nicexy{
A_{X} \ar[r]^-{\omega_{Y,y}} \ar[dr]_{=} 
& A_Y \ar[d]^-{\sigma_{Y,x,y}} \\
 & A_X
}\]
is commutative.
\item
In the setting of \eqref{move:e1}, the diagram
\[
\nicexy@C+.5cm{
A_{X} \ar[r]^-{\sigma^{Y^{23},y^1,y^{23}}} \ar[d]_{\sigma^{Y^{12},y^{12},y^3}} 
& A_{Y^{23}} \ar[d]^-{\sigma^{Y,y^2,y^3}} \\
A_{Y^{12}} \ar[r]^-{\sigma^{Y,y^1,y^2}} & A_Y
}\]
is commutative.
\item
In the setting of \eqref{move:e2}, the diagram
\[
\nicexy@C+.5cm{
A_{X} \ar[r]^-{\sigma^{Y^{34},y^1,y^2}} \ar[d]_{\sigma^{Y^{12},y^3,y^4}} 
& A_{Y^{34}} \ar[d]^-{\sigma^{Y,y^3,y^4}} \\
A_{Y^{12}} \ar[r]^-{\sigma^{Y,y^1,y^2}} & A_Y
}\]
is commutative.
\item
In the setting of \eqref{move:e3}, the diagram
\[
\nicexy{
A_{X} \ar[r]^-{\omega_{W,y}} \ar[d]_{\sigma^{Z,y^1,y^2}} 
& A_W \ar[d]^-{\sigma^{Y,y^1,y^2}} \\
A_Z \ar[r]^-{\omega_{Y,y}} & A_Y
}\]
is commutative.
\item
In the setting of \eqref{move:f1}, the diagram
\[
\nicexy{
A_{X} \ar[r]^-{\omega_{Y_2,y_1}} \ar[d]_{\omega_{Y_1,y_2}} 
& A_{Y_2} \ar[d]^-{\omega_{Y,y_2}} \\
A_{Y_1} \ar[r]^-{\omega_{Y,y_1}} & A_Y
}\]
is commutative.
\end{enumerate}
This finishes the definition of a $\WD$-algebra.
\end{definition}

At this moment we have two definitions of a $\WD$-algebra.  
\begin{enumerate}
\item
On the one hand, in Def. \ref{def2:operad-algebra} with $\O = \WD$, a $\WD$-algebra has a structure map $\mu_{\zeta}$ \eqref{operad-structure-map} for each wiring diagram $\zeta$.  This structure map satisfies the associativity axiom \eqref{operad-algebra-associativity2} for a general operadic composition in $\WD$, together with the unity and the equivariance axioms in Def. \ref{def1:operad-algebra}.  
\item
On the other hand, in Def. \ref{def:wd-algebra} a $\WD$-algebra has $8$ generating structure maps and satisfies $28$ generating axioms.  
\end{enumerate}
We now observe that these two definitions are equivalent, so $\WD$-algebras indeed have a finite presentation as in Def. \ref{def:wd-algebra}.

\begin{theorem}
\label{thm:wd-algebra}
For the operad of wiring diagrams $\WD$ (Theorem \ref{wd-operad}), Def. \ref{def2:operad-algebra} with $\O=\WD$ and Def. \ref{def:wd-algebra} of a $\WD$-algebra are equivalent.\index{finite presentation for WD algebra@finite presentation for $\WD$-algebras}
\end{theorem}

\begin{proof}
First suppose $(A,\mu)$ is a $\WD$-algebra in the sense of Def. \ref{def2:operad-algebra}.  To see that it is also a $\WD$-algebra in the sense of Def. \ref{def:wd-algebra}, first note that the structure map $\mu_?$ \eqref{operad-structure-map} gives the eight generating  structure maps \eqref{structuremap-empty}--\eqref{structuremap-wastedwire}.  Moreover, the generating structure map $\mu_{\tensorunit_X}$ \eqref{structuremap-namechange} is the identity map by the unity axiom \eqref{operad-algebra-unity}.

The generating axiom \eqref{wd-algebra-equivariance} is a special case of the equivariance diagram \eqref{operad-algebra-eq}, so it is commutative.  Each of the other $27$ generating axioms corresponds to an elementary relation that describes two different ways to construct the same wiring diagram as an iterated operadic composition of generating wiring diagrams.  Each such generating axiom asserts that the two corresponding compositions of generating structure maps--defined using the composition in the diagram \eqref{operad-algebra-associativity2}--are equal.  The associativity axiom \eqref{operad-algebra-associativity2} of $(A,\mu)$ applied twice guarantees that two such compositions are indeed equal.

Conversely, suppose $A$ is a $\WD$-algebra in the sense of Def. \ref{def:wd-algebra}, so it has eight generating structure maps that satisfy $28$ generating axioms.  We must show that it is a $\WD$-algebra in the sense of Def. \ref{def2:operad-algebra}  For a wiring diagram $\psi \in \WD$ with a presentation $\Psi$ (Def. \ref{def:simplex}), we define its structure map $\mu_{\psi}$ \eqref{operad-structure-map} inductively as follows.
\begin{enumerate}
\item
If $\Psi$ is a $1$-simplex, then $\Psi = (\psi)$, and $\psi$ is a generating wiring diagram by the definition of a simplex.  In this case, we define $\mu_{\psi}$ as the corresponding generating structure map  \eqref{structuremap-empty}--\eqref{structuremap-wastedwire} of $A$.
\item
Inductively, suppose $\Psi$ is an $n$-simplex for some $n \geq 2$, so $\Psi = (\Phi, i, \Theta)$ for some $i \geq 1$, $p$-simplex $\Phi$, and $q$-simplex $\Theta$ with $p+q = n$.  Since $1 \leq p,q < n$, by the induction hypothesis, the structure maps $\mu_{|\Phi|}$ and $\mu_{|\Theta|}$ are already defined.  Then we define the structure map
\begin{equation}
\label{mu-compi-algebra}
\mu_{\psi} = \mu_{|\Phi|} \compi \mu_{|\Theta|}
\end{equation}
as in Notation \ref{notation:algebra-compi}.
\end{enumerate}

By Theorem \ref{stratified-presentation-exists} every wiring diagram has a stratified presentation, hence a presentation.  To see that the structure map $\mu_{\psi}$ as above is well-defined, we need to show that the map $\mu_{\psi}$ is independent of the choice of a presentation $\Psi$.  Any two presentations of a wiring diagram are by definition equivalent simplices.  By Theorem \ref{thm:wd-generator-relation}(2) ($=$ the relations part of the finite presentation theorem for $\WD$), any two equivalent simplices are either equal or are connected by a finite sequence of elementary equivalences.  Therefore, it suffices to show that every elementary equivalence in $\WD$ yields a commutative diagram involving the generating structure maps of $A$, where $\compi$ is interpreted as in Notation \ref{notation:algebra-compi}.  Recall from Def. \ref{def:equivalent-simplices} that an elementary equivalence comes from either an elementary relation or an operad associativity/unity axiom for the generating wiring diagrams.

It follows from a direct inspection that the operad associativity and unity axioms--\eqref{compi-associativity}, \eqref{compi-associativity-two}, \eqref{compi-left-unity}, and \eqref{compi-right-unity}--for the generating wiring diagrams yield commutative diagrams involving the generating structure maps of $A$.  In fact, the diagrams involving the horizontal and the vertical associativity axioms \eqref{compi-associativity} and \eqref{compi-associativity-two} are commutative because composition of functions is associative.  The diagrams for the two unity axioms \eqref{compi-left-unity} and \eqref{compi-right-unity} are commutative because the generating structure map for a colored unit \eqref{structuremap-namechange} is required to be the identity map.

By definition each of the $28$ generating axioms of $A$ corresponds to an elementary relation (Def. \ref{def:elementary-relations}) and is a commutative diagram.   Therefore, the structure map $\mu_{\psi}$ for each wiring diagram $\psi$ is well-defined.

It remains to check that the structure map $\mu$ satisfies the required unity, equivariance, and associativity axioms.  The unity axiom \eqref{operad-algebra-unity} holds because it is part of the assumption on the generating structure map corresponding to a name change \eqref{structuremap-namechange}.

The associativity axiom \eqref{operad-algebra-associativity2} holds because the structure map $\mu_{\psi}$ is defined above \eqref{mu-compi-algebra} by requiring that the  diagram  \eqref{operad-algebra-associativity2} be commutative.

For the equivariance axiom \eqref{operad-algebra-eq}, first note that it is enough to check it when the wiring diagram in questioned is an iterated operadic composition of $2$-cells.  This is because $2$-cells are the only binary generating wiring diagrams (Remark \ref{rk:generators-arity}).  All other generating wiring diagrams are either $0$-ary or unary, for which equivariance is trivial.

So now suppose $\zeta$ in the equivariance axiom \eqref{operad-algebra-eq} is an iterated operadic composition of $2$-cells.  If $\zeta$ is a $2$-cell and the permutation $\sigma$ is the transposition $(1~2) \in \Sigma_2$, then  the equivariance axiom \eqref{operad-algebra-eq} is true by the generating axiom \eqref{wd-algebra-equivariance}.  The general case now follows from this special case using:
\begin{itemize}
\item
the generating axiom \eqref{wd-algebra-2cell-associativity} corresponding to the associativity property of $2$-cells \eqref{move:b1};
\item
the operad associativity axioms \eqref{compi-associativity} and \eqref{compi-associativity-two} when applied to $2$-cells;
\item 
the fact that the transpositions $(i, i+1)$ for $1 \leq i \leq n-1$ generate the symmetric group $\Sigma_n$.
\end{itemize}
So  $(A,\mu)$ is a $\WD$-algebra in the sense of Def. \ref{def2:operad-algebra}.
\end{proof}

\section{Finite Presentation for the Propagator Algebra}
\label{sec:propagator}

As an illustration of Theorem \ref{thm:wd-algebra}, in this section we provide a finite presentation for the $\WD$-algebra called the \emph{propagator algebra} that was first introduced in \cite{rupel-spivak} (Section 3).  This finite presentation is about the structure maps, not the underlying sets.  As explained in \cite{rupel-spivak} (Section 3.4), the propagator algebra can be used, for example, to provide an iterative description of the Fibonacci sequence.  In contrast to the original definition in \cite{rupel-spivak}, we will define the propagator algebra using finitely many generating structure maps and axioms--$8$ generating structure maps and $28$ generating axioms to be exact.  Since our generating structure maps are rather simple, our description of the propagator algebra is less involved than the original definition and verification in \cite{rupel-spivak}.  The equivalence between the two definitions of the propagator algebra is explained in Remark \ref{rk:propagator-agree}.

\begin{assumption}
Throughout this section, $S$ denotes the class of pointed sets, with respect to which $S$-boxes (Def. \ref{def:s-box}) and the operad $\WD$ are defined (Theorem \ref{wd-operad}).  In a pointed set, the base point is denoted by $*$.  
\end{assumption}

Recall the concept and notations regarding profiles from Def. \ref{def:profile}.  Let us first recall a few definitions from \cite{rupel-spivak} (section 3).  

\begin{definition}
\label{def:historical}
Suppose $T$ and $U$ are pointed sets and $k \geq 0$.
\begin{enumerate}
\item
Define the \emph{truncation} \index{truncation} $\dT : \profonet \to \proft$ as the function
\begin{equation}
\label{truncation}
\dT(t_1,\ldots,t_n) = (t_1, \ldots,t_{n-1}).
\end{equation}
We will often omit the subscript and just write \label{notation:truncation}$\partial$.
\item
A \emph{$k$-historical propagator from $T$ to $U$} \index{historical propagator} is a function $f : \proft \to \profu$ such that:
\begin{enumerate}[(i)]
\item
$|f(\ut)| = |\ut| + k$ for all $\ut \in \proft$;
\item
If $\ut \in \proft$ has length $|\ut| \geq 1$, then 
\begin{equation}
\label{historicity}
\dU f(\ut) = f(\dT \ut).
\end{equation}
The condition \eqref{historicity} is called \emph{historicity}.\index{historicity}
\end{enumerate}
\item
The set of $k$-historical propagators from $T$ to $U$ is denoted by \label{notation:hist}$\Hist^k(T,U)$.
\item
A \emph{historical propagator from $T$ to $U$} is an $m$-historical propagator from $T$ to $U$ for some $m \geq 0$.
\end{enumerate}
\end{definition}

\begin{example}\label{ex:moment-delay}
Given a pointed set $T$ and an integer $k \geq 0$, the function $D_k : \proft \to \proft$ defined as
\[D_k(t_1,\ldots,t_n) = (*, \ldots, *, t_1,\ldots, t_n),\]
in which the right side starts with $k$ entries of the base point $*$, is a $k$-historical propagator, called the \index{moment delay}$k$-moment delay function in \cite{rupel-spivak}.
\end{example}

Before we can define the propagator algebra, we first need to define its entries.

\begin{definition}
\label{propagator-algebra-sets}
Suppose $X = (\xin, \xout) \in \boxs$ (Def. \ref{def:s-box}).
\begin{enumerate}
\item
Define the pointed sets
\begin{equation}
\label{vxin-vxout}
\vxin = \prod_{x \in \xin} v(x) \andspace
\vxout = \prod_{x \in \xout} v(x)
\end{equation}
in which each $v(x)$, a pointed set, is the value of $x$ (Def. \ref{def:Fins}) and an empty product means the one-point set.
\item
Define the set\label{notation:p-sub-x}
\begin{equation}
\label{propagator-x}
\sP_X = \Histone\bigl(\vxin, \vxout\bigr)
\end{equation}
of \emph{$1$-historical propagators of type $X$}.\index{historical propagators of type $X$}
\end{enumerate}
\end{definition}

So a $1$-historical propagator of type $X$ is a function that takes each $\vxin$-profile to an $\vxout$-profile whose length is one higher than before and that satisfies historicity \eqref{historicity}.  

\begin{example}\label{ex:fg}
Suppose $(\bN,1)$ is the pointed set of non-negative integers with base point $1$.  Consider the box $X$ with $\xin$ and $\xout$ both equal to the one-point set $*$ with value $v(*) = (\bN,1)$.
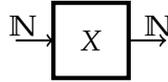
\begin{center}
\begin{tikzpicture}[scale=.5]
\draw [ultra thick] (0,0) rectangle (2,2);
\draw [arrow] (-1,1) -- (0,1);
\draw [arrow] (2,1) -- (3,1);
\node at (1,1) {$X$};
\node at (-.8,1.4) {$\bN$};
\node at (2.8,1.4) {$\bN$};
\end{tikzpicture}
\end{center}
A $1$-historical propagator of type $X$ is a function that takes each finite sequence of non-negative integers to a sequence of non-negative integers whose length is one higher than before and that satisfies historicity.  For example:
\begin{enumerate}
\item The $1$-moment delay function in Example \ref{ex:moment-delay}, given by
\[D_1(m_1,\ldots,m_n) = (1,m_1,\ldots,m_n)\]
for $m_i \in \bN$, is a $1$-historical propagator of type $X$.  For instance, we have $D_1(1,5,6) = (1,1,5,6)$.
\item The function $f : \prof(\bN) \to \prof(\bN)$ defined as
\[f(m_1,\ldots,m_n) = \left(0, m_1, m_1 + m_2, \ldots, \sum_{i=1}^n m_i\right)\]
is a $1$-historical propagator of type $X$, denoted $``\Sigma"$ in \cite{rupel-spivak}.  This $1$-historical propagator takes a sequence of non-negative integers to the sequence whose $i$th entry is the sum of the first $i-1$ entries of the given sequence.  For instance, we have $f(1,5,6) = (0,1,6,12)$.
\item The function $g : \prof(\bN) \to \prof(\bN)$ defined as
\[g(m_1,\ldots,m_n) = \left(1, m_1, m_1m_2, \ldots, \prod_{i=1}^n m_i\right)\]
is a $1$-historical propagator of type $X$. This $1$-historical propagator takes a sequence of non-negative integers to the sequence whose $i$th entry is the product of the first $i-1$ entries of the given sequence.  For instance, we have $g(1,5,6) = (1,1,5,30)$.
\end{enumerate}
\end{example}

\begin{example}\label{ex:hj}
Consider the box $Y$ with $\yin = \{y_1,y_2\}$, $\yout = \{y\}$, and values $v(y_1) = v(y_2) = v(y) = (\bN,1)$. 
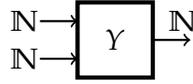
\begin{center}
\begin{tikzpicture}[scale=.5]
\draw [ultra thick] (0,0) rectangle (2,2);
\draw [arrow] (-1,.5) -- (0,.5);
\draw [arrow] (-1,1.5) -- (0,1.5);
\draw [arrow] (2,1) -- (3,1);
\node at (1,1) {$Y$};
\node at (-1.4,.5) {$\bN$};
\node at (-1.4,1.5) {$\bN$};
\node at (2.8,1.5) {$\bN$};
\end{tikzpicture}
\end{center}
A $1$-historical propagator of type $Y$ is a function that takes each finite sequence of ordered pairs of non-negative integers to a sequence of non-negative integers whose length is one higher than before and that satisfies historicity.  For example:
\begin{enumerate}
\item The function $h : \prof(\bN \times \bN) \to \prof(\bN)$ given by
\[h\bigl((m_1,m_1'), \ldots, (m_n,m_n')\bigr) = \left(1,m_1+m_1',\ldots,m_n+m_n'\right)\]
is a $1$-historical propagator of type $Y$, denoted $``+"$ in \cite{rupel-spivak}.  For instance, we have $h\bigl((1,4),(5,2),(6,8)\bigr) = (1,5,7,14)$.
\item The function $j : \prof(\bN \times \bN) \to \prof(\bN)$ given by
\[j\bigl((m_1,m_1'), \ldots, (m_n,m_n')\bigr) = \left(1,m_1m_1',\ldots,m_nm_n'\right)\]
is a $1$-historical propagator of type $Y$.  For instance, we have
\[j\bigl((1,4),(5,2),(6,8)\bigr) = (1,4,10,48).\]
\end{enumerate}
\end{example}

\begin{example}\label{ex:ell}
Consider the box $Z$ with $\zin = \{z_1,z_2,z_3\}$, $\zout = \{z^1,z^2\}$, and all $v(-) = (\bN,1)$.
\begin{center}
\begin{tikzpicture}[scale=.7]
\draw [ultra thick] (0,0) rectangle (2,2);
\draw [arrow] (-1,.2) -- (0,.2);
\draw [arrow] (-1,1) -- (0,1);
\draw [arrow] (-1,1.8) -- (0,1.8);
\draw [arrow] (2,.5) -- (3,.5);
\draw [arrow] (2,1.5) -- (3,1.5);
\node at (1,1) {$Z$};
\node at (-1.4,.2) {$\bN$};
\node at (-1.4,1) {$\bN$};
\node at (-1.4,1.8) {$\bN$};
\node at (3.4,.5) {$\bN$};
\node at (3.4,1.5) {$\bN$};
\node at (.3,.3) {\scriptsize{$z_3$}};
\node at (.3,1) {\scriptsize{$z_2$}};
\node at (.3,1.7) {\scriptsize{$z_1$}};
\node at (1.7,.5) {\scriptsize{$z^2$}};
\node at (1.7,1.5) {\scriptsize{$z^1$}};
\end{tikzpicture}
\end{center}
A $1$-historical propagator of type $Z$ is a function that takes each finite sequence of ordered triples of non-negative integers to a sequence of ordered pairs of non-negative integers whose length is one higher than before and that satisfies historicity.  For example, the function
\[\ell : \prof(\bN\times\bN\times\bN) \to \prof(\bN\times\bN)\]
given by
\[\ell\bigl((m_1,m_1',m_1''), \ldots, (m_n,m_n',m_n'')\bigr) = \bigl((1,1),(m_1,m_1'+m_1''),\ldots,(m_n, m_n'+m_n'')\bigr)\]
is a $1$-historical propagator of type $Z$.  For instance, we have
\[\ell\bigl((3,1,4), (7,2,9), (8,5,10)\bigr) = 
\bigl((1,1), (3,5), (7,11), (8,15)\bigr).\]
We will come back to this example several times below.
\end{example}

The generating structure map of the propagator algebra associated to a $1$-loop requires a few notations in its definition.  So here we define this map first.  The reader should keep in mind that the following definition as well as all the proofs in this section involve simple inductions on the length of some profiles.

\begin{definition}
\label{def:generating-structure-1loop}
Suppose $X \in \boxs$, $x_- \in \xin$, and $x_+ \in \xout$ such that $v(x_-) = v(x_+)$ as pointed sets.  The box $X \setminus x \in \boxs$ is obtained from $X$ by removing $x = \{x_{\pm}\}$.  For $\ut \in \prof\bigl(\vxout\bigr)$, we will write:
\begin{enumerate}[(i)]
\item
$\ut_{\xplus} \in \prof\bigl(v(\xplus)\bigr)$\label{notation:sub-xplus} for the profile obtained from $\ut$ by  taking only the $v(\xplus)$-entry;
\item
$\ut_{\minusxplus} \in \prof\bigl(\vxminusxout\bigr)$\label{notation:sub-minus-xplus}  for the profile obtained from  $\ut$ by removing the $v(\xplus)$-entry.
\end{enumerate}
Suppose $g \in \sP_X$ \eqref{propagator-x}.  Define two functions
\begin{equation}
\label{propagator-1loop-defined}
\nicexy{
\prof\bigl(\vxminusxin\bigr) \ar[r]^-{\lambda g} & \prof\bigl(\vxminusxout\bigr)}
\end{equation}
and
\begin{equation}
\label{propagator-G}
\nicexy{
\prof\bigl(\vxminusxin\bigr) \ar[r]^-{G_g} & \prof\bigl(v(x_+)\bigr)}
\end{equation}
with the properties
\begin{equation}
\label{lambdag-length}
|(\lambda g)(?)| = |?| + 1 = |G_g(?)|
\end{equation}
inductively as follows.
\begin{enumerate}[(i)]
\item
For the empty profile, define
\begin{equation}
\label{lambdag-empty}
\begin{split}
(\lambda g)(\varnothing) &= g(\varnothing)_{\minusxplus} \in \prof\bigl(\vxminusxout\bigr)\\
G_g(\varnothing) &= g(\varnothing)_{\xplus} \in  \prof\bigl(v(x_+)\bigr).
\end{split}
\end{equation}
In each definition in \eqref{lambdag-empty}, the first $\varnothing$ is the empty $\vxminusxin$-profile, and the second $\varnothing$, to which $g$ applies, is the empty $\vxin$-profile.  The profile $g(\varnothing)$ has length $1$ because $g \in \Histone\left(\vxin,\vxout\right)$.  So both $(\lambda g)(\varnothing)$ and $G_g(\varnothing)$ have length $1$.
\item
Inductively, suppose $\uw \in \prof\bigl(\vxminusxin\bigr)$ has length $n \geq 1$.  Define
\begin{equation}
\label{lambdag-general}
\begin{split}
(\lambda g)(\uw) &= g\bigl(\uw, G_g(\partial \uw)\bigr)_{\minusxplus} \in \prof\bigl(\vxminusxout\bigr)\\
G_g(\uw) &= g\bigl(\uw, G_g(\partial \uw)\bigr)_{\xplus} \in  \prof\bigl(v(x_+)\bigr).
\end{split}
\end{equation}
Here $\partial$ is the truncation \eqref{truncation}, so the profile
\[
G_g(\partial \uw) \in \prof\bigl(v(\xplus)\bigr) =  \prof\bigl(v(\xminus)\bigr)
\]
is already defined and has length $n$ by the induction hypothesis.  In each definition in \eqref{lambdag-general}, 
\[
\bigl(\uw, G_g(\partial \uw)\bigr) \in \prof\left(\vxin\right)\]
has length $n$, so its image under $g$ has length $n+1$.  Therefore, both $(\lambda g)(\uw)$ and $G_g(\uw)$ have length $n+1$.
\end{enumerate}
We say that $\lambda g$ and $G_g$ are \emph{defined with respect to $x = \{x_{\pm}\}$}.
\end{definition}

\begin{example}\label{ex2:ell}
This is a continuation of Example \ref{ex:ell}, where the box $Z$ has $\zin = \{z_1,z_2,z_3\}$, $\zout = \{z^1,z^2\}$, and all $v(-) = (\bN,1)$.  For $z_1 \in \zin$ and $z^1 \in \zout$, suppose $Z \setminus z$ is the box obtained from $Z$ by removing $z = \{z_1,z^1\}$.  So we have
\[\begin{split}
\vzminuszin &= v(z_2) \times v(z_3) = \bN \times \bN;\\
\vzminuszout &= v(z^2) = \bN.
\end{split}\]
\begin{center}
\begin{tikzpicture}[scale=.7]
\draw [ultra thick] (0,0) rectangle (2,2);
\draw [arrow] (-2,.2) -- (0,.2);
\draw [arrow] (-2,1) -- (0,1);
\draw [arrow] (2,.5) -- (4,.5);
\draw [arrow, looseness=2] (2,1.8) to [out=30, in=150] (0,1.8);
\node at (1,1) {$Z$};
\node at (.3,.3) {\scriptsize{$z_3$}};
\node at (.3,1) {\scriptsize{$z_2$}};
\node at (.3,1.7) {\scriptsize{$z_1$}};
\node at (1.7,.5) {\scriptsize{$z^2$}};
\node at (1.7,1.7) {\scriptsize{$z^1$}};
\draw [ultra thick] (-1,-.5) rectangle (3,2.7);
\node at (-2.5,1.5) {$\zminuszin$};
\node at (4.5,1.5) {$\zminuszout$};
\node at (1,3.2) {$\lambda_{Z,z} \in \WD\zminuszz$};
\end{tikzpicture}
\end{center}
For the $1$-historical propagator $\ell$ of type $Z$ defined as
\[\ell\bigl((m_1,m_1',m_1''), \ldots, (m_n,m_n',m_n'')\bigr) = \bigl((1,1),(m_1,m_1'+m_1''),\ldots,(m_n, m_n'+m_n'')\bigr),\]
the functions
\[\nicexy{\prof(\bN\times\bN) \ar[r]^-{\lambda\ell} & \prof(\bN)} \andspace 
\nicexy{\prof(\bN\times\bN) \ar[r]^-{G_{\ell}} & \prof(\bN)}\]
in \eqref{propagator-1loop-defined} and \eqref{propagator-G} are given as follows.  Suppose
\[\um = \bigl((m_1,m_1'), \ldots, (m_n,m_n')\bigr) \in \prof(\bN\times\bN)\]
is an $(\bN\times\bN)$-profile of length $n \geq 0$.  Then a simple induction shows that
\[\begin{split}
(\lambda\ell)(\um) &= \bigl(1, m_1+m_1', \ldots, m_n+m_n'\bigr),\\
G_{\ell}(\um) &= (1,1,\ldots,1),
\end{split}\]
in which on the right side of $G_{\ell}(\um)$ there are $n+1$ copies of $1$.  In particular, $\lambda\ell = h$, the $1$-historical propagator of type $Y$ in Example \ref{ex:hj}.  In this example, both $\lambda\ell$ and $G_{\ell}$ are $1$-historical propagators.  This is not an accident, as we show in the next result.
\end{example}

The following observation will use Def. \ref{def:historical} about historical propagators.

\begin{lemma}
\label{lemma:lambdag-historical}
In the context of Def. \ref{def:generating-structure-1loop} with $g \in \sP_X$, the following statements hold.
\begin{enumerate}
\item
$G_g \in \Histone\bigl(\vxminusxin, v(\xplus)\bigr)$.
\item
$\lambda g \in \Histone\bigl(\vxminusxin, \vxminusxout\bigr) = \sP_{X\setminus x}$.
\end{enumerate}
\end{lemma}

\begin{proof}
In this proof, we will abbreviate $G_g$ to $G$.  For statement (1), we are trying to show that $G$ is a $1$-historical propagator from $\vxminusxin$ to $v(\xplus)$.  In view of the property \eqref{lambdag-length}, it remains to check  historicity  \eqref{historicity} for $G$, which we will do by induction.  Suppose $\uw \in \prof\bigl(\vxminusxin\bigr)$ has length $\geq 1$.

If $|\uw|=1$, then $\partial \uw = \varnothing$.  Using the definitions \eqref{lambdag-empty} and  \eqref{lambdag-general} we have:
\[
\begin{split}
\partial G(\uw) &= \partial g\bigl(\uw, G(\varnothing)\bigr)_{\xplus}\\
&= g\bigl(\partial\left(\uw, g(\varnothing)_{\xplus}\right) \bigr)_{\xplus} \quad\text{by historicity of $g$}\\
&= g(\varnothing)_{\xplus} = G(\varnothing).
\end{split}\]
Inductively, suppose $|\uw| \geq 2$.  Using the definition \eqref{lambdag-general}, we have:
\[
\begin{split}
\partial G(\uw) &= \partial g\bigl(\uw, G(\partial \uw)\bigr)_{\xplus}\\
&= g\bigl(\partial\uw, \partial G(\partial\uw) \bigr)_{\xplus} \quad\text{by historicity of $g$}\\
&= g\bigl(\partial \uw, G(\partial \partial \uw)\bigr)_{\xplus} \quad\text{by the induction hypothesis on $G$}\\
&= G(\partial \uw).
\end{split}\]
This finishes the proof of statement (1).

Statement (2) is proved by essentially the same argument as above, except that, in view of  the definitions \eqref{lambdag-empty} and  \eqref{lambdag-general}, the various right-most subscripts $x_+$ are replaced by $\minusxplus$.
\end{proof}

We are now ready to define the propagator algebra in terms of finitely many generating structure maps and generating axioms as in Def. \ref{def:wd-algebra}.  Most of the generating structure maps below are easily seen to be well-defined.  The only exception is the generating structure map associated to a $1$-loop, which we dealt with in Lemma \ref{lemma:lambdag-historical} above.

\begin{definition}
\label{def:propagator-algebra}
Define the \emph{propagator algebra} \index{propagator algebra} $\sP$ as the $\WD$-algebra, in the sense of  Def. \ref{def:wd-algebra}, with $X$-colored entry
\[
\sP_X = \Histone\bigl(\vxin, \vxout\bigr)
\]
as in \eqref{propagator-x} for each $X \in \boxs$.  Its $8$ generating structure maps are defined as follows.
\begin{enumerate}
\item
Corresponding to the empty wiring diagram $\epsilon \in \WD\emptynothing$ (Def. \ref{def:empty-wd}), the structure map
\begin{equation}
\label{propagator-empty}
\nicexy{\{*\} \ar[r]^-{\epsilon} & \sP_{\varnothing} = \Histone\bigl(\{*\}, \{*\}\bigr) = \{*\}}
\end{equation}
sends $*$ to the unique function
\[
\left(\epsilon *\right)\bigl(\underbrace{*,\ldots,*}_{m}\bigr) = 
\bigl(\underbrace{*,\ldots,*}_{m+1}\bigr).
\]
\item
Corresponding to each $1$-delay node $\delta_d \in \WD\dnothing$ (Def. \ref{def:one-dn}) with $d$ a pointed set, the structure map
\begin{equation}
\label{propagator-delaynode}
\nicexy{\{*\} \ar[r]^-{\delta_d} & \sP_d = \Histone\bigl(d, d\bigr)}
\end{equation}
sends $*$ to the function
\[
(\delta_d *)(\ut) = \bigl(*, \ut\bigr)
\]
for each $\ut \in \prof(d)$.  Here the $*$ on the right is the base point in $d$.  In other words, $\delta_d*$ is the $1$-moment delay function in Example \ref{ex:moment-delay}.
\item
Corresponding to each name change $\tau_{X,Y} \in \WD\yx$ (Def. \ref{def:name-change}), the structure map
\begin{equation}
\label{propagator-namechange}
\nicexy{\Histone\bigl(\vxin, \vxout\bigr) = \sP_X \ar[r]^-{\tau_{X,Y}} 
& \sP_Y = \Histone\bigl(\vyin, \vyout\bigr)}
\end{equation}
is the identity map.  Here we are using the fact that, if $x \in \xin \amalg \xout$ and $y \in \yin \amalg \yout$ corresponds to $x$ under $\tau_{X,Y}$, then $v(x) = v(y)$ as pointed sets.
\item
Corresponding to each $2$-cell $\theta_{X,Y} \in \WD\xplusyxy$ (Def. \ref{def:theta-wd}), the structure map
\begin{equation}
\label{propagator-2cell}
\nicexy{\sP_X \times \sP_Y  = \Histone\bigl(\vxin, \vxout\bigr) \times \Histone\bigl(\vyin, \vyout\bigr) \ar[d]^-{\theta_{X,Y}}\\ 
\sP_{X \amalg Y} = \Histone\bigl(\vxplusyin, \vxplusyout\bigr)
= \Histone\bigl(\vxin \times \vyin, \vxout \times \vyout\bigr)}
\end{equation}
is given by
\[
\theta_{X,Y}(f_X,f_Y) = f_X \times f_Y
\]
for $f_X \in \Histone\bigl(\vxin, \vxout\bigr)$ and $f_Y \in \Histone\bigl(\vyin, \vyout\bigr)$.
\item
Corresponding to each $1$-loop $\lambda_{X,x} \in \WD\xminusxx$ (Def. \ref{def:loop-wd}), the structure map
\begin{equation}
\label{propagator-1loop}
\nicexy{\sP_X  \ar[r]^-{\lambda_{X,x}} & \sP_{X \setminus x}}
\end{equation}
sends each $g \in \sP_X$ to $\lambda g \in \sP_{X \setminus x}$ \eqref{propagator-1loop-defined}, which is well-defined by Lemma \ref{lemma:lambdag-historical}.
\item
Corresponding to each in-split $\sigma_{X,x_1,x_2} \in \WD\yx$ (Def. \ref{def:insplit-wd}), the structure map
\begin{equation}
\label{propagator-insplit}
\nicexy@C+.5cm{\sP_X  \ar[r]^-{\sigma_{X,x_1,x_2}} & \sP_Y}
\end{equation}
is given by
\[
\left(\sigma_{X,x_1,x_2} g\right)(\uy) = g(\pi \uy)
\]
for $g \in \sP_X$ and $\uy \in \prof(\vyin)$.  Here $\pi \uy \in \prof(\vxin)$ is the same as $\uy$ except that its $v(x_1)$-entry and $v(x_2)$-entry are both given by the $v(x_{12})$-entry of $\uy$. 
\item
Corresponding to each out-split  $\sigma^{Y,y_1,y_2} \in \WD\yx$ (Def. \ref{def:out-split}), the structure map
\begin{equation}
\label{propagator-outsplit}
\nicexy@C+.5cm{\sP_X  \ar[r]^-{\sigma^{Y,y_1,y_2}} & \sP_Y}
\end{equation}
is given by
\[
\left(\sigma^{Y,y_1,y_2} g\right)(\uy) = \pi g(\uy)
\]
for $g \in \sP_X$ and $\uy \in \prof(\vyin) = \prof(\vxin)$.  Here $\pi g(\uy) \in \prof(\vyout)$ is the same as $g(\uy) \in \prof(\vxout)$ except that its $v(y_1)$-entry and $v(y_2)$-entry are both given by the $v(y_{12})$-entry of $g(\uy)$. 
\item
Corresponding to each $1$-wasted wire  $\omega_{Y,y} \in \WD\yx$ (Def. \ref{def:wasted-wire-wd}), the structure map
\begin{equation}
\label{propagator-wastedwire}
\nicexy{\sP_X  \ar[r]^-{\omega_{Y,y}} & \sP_Y}
\end{equation}
is given by
\[\bigl(\omega_{Y,y} g\bigr)(\ut) = g\bigl(\ut_{\minusy}\bigr)\]
for $g \in \sP_X$ and $\ut \in \prof(\vyin)$.  Here $\ut_{\minusy}  \in \prof(\vxin)$ is obtained from $\ut$ by removing the $v(y)$-entry.
\end{enumerate}
This finishes the definition of the propagator algebra $\sP$.
\end{definition}

The following four examples continue Examples \ref{ex:ell} and \ref{ex2:ell}, where $Z$ is the box with $\zin = \{z_1,z_2,z_3\}$, $\zout = \{z^1,z^2\}$, and all $v(-) = (\bN,1)$.  The $1$-historical propagator $\ell \in \sP_Z = \Histone(\bN^3,\bN^2)$ is defined as 
\[\ell\bigl((m_1,m_1',m_1''), \ldots, (m_n,m_n',m_n'')\bigr) = \bigl((1,1),(m_1,m_1'+m_1''),\ldots,(m_n, m_n'+m_n'')\bigr).\]  
Let us consider the image of $\ell$ under some of the structure maps of $\sP$.

\begin{example}\label{ex:ell-loop}
Suppose $Y$ is the box obtained from $Z$ by removing $z = \{z_1,z^1\}$.  Consider the $1$-loop (Def. \ref{def:loop-wd})
\begin{center}
\begin{tikzpicture}[scale=.9]
\draw [ultra thick] (2,1) rectangle (3,2);
\node at (2.5,1.4) {$Z$};
\node at (2.2,1.75) {\tiny{$z_1$}};
\node at (2.8,1.8) {\tiny{$z^1$}};
\draw [arrow] (.5,1.2) to (2,1.2);
\draw [arrow] (.5,1.5) to (2,1.5);
\node at (-.1,1.4) {$\yin$};
\draw [arrow] (3,1.3) to (4.5,1.3);
\node at (5.2,1.4) {$\yout$};
\draw [ultra thick] (1,0.7) rectangle (4,2.5);
\node at (2.5,2.8) {$\lambda_{Z,z}$};
\draw [arrow, looseness=3.5] (3,1.8) to [out=30, in=150] (2,1.8);
\end{tikzpicture}
\end{center}
in $\WD\yz$, the structure map
\[\nicexy{\sP_Z = \Histone(\bN^3,\bN^2)  \ar[r]^-{\lambda_{Z,z}} & \Histone(\bN^2,\bN) = \sP_{Y}}\]
in \eqref{propagator-1loop}, and $\um = \bigl((m_1,m_1'), \ldots, (m_n,m_n')\bigr) \in \prof(\vyin) = \prof(\bN^2)$. As observed in Example \ref{ex2:ell} we have
\[(\lambda_{Z,z}\ell)(\um) = h(\um) = \left(1,m_1+m_1',\ldots,m_n+m_n'\right)\]
in $\prof(\vyout) = \prof(\bN)$, where $h \in \sP_Y$ is the $1$-historical propagator in Example \ref{ex:hj}.  For instance, we have
\[(\lambda_{Z,z}\ell)\bigl((2,5),(4,9),(3,7)\bigr) = (1,7,13,10).\]
\end{example}

\begin{example}\label{ex:ell-insplit}
Suppose $W$ is the box with $\win = \zin/(z_1=z_2) = \{w,z_3\}$, and $\wout = \zout$.  For the in-split (Def. \ref{def:insplit-wd})
\begin{center}
\begin{tikzpicture}[scale=.9]
\draw [ultra thick] (2,1) rectangle (3,2.1);
\node at (2.7,1.5) {$Z$};
\node at (2.2,1.9) {\tiny{$z_1$}};
\node at (2.2,1.5) {\tiny{$z_2$}};
\node at (2.2,1.2) {\tiny{$z_3$}};
\draw [arrow] (.5,1.2) to (2,1.2);
\node at (-.1,1.6) {$\win$};
\node at (.7,1.9) {\tiny{$w$}};
\draw [thick] (.5,1.7) -- (1.2,1.7);
\draw [arrow] (1.2,1.7) to [out=0, in=180] (2,1.9);
\draw [arrow] (1.2,1.7) to [out=0, in=180] (2,1.5);
\draw [arrow] (3,1.3) to (4.5,1.3);
\draw [arrow] (3,1.8) to (4.5,1.8);
\node at (5.2,1.6) {$\wout$};
\draw [ultra thick] (1,0.7) rectangle (4,2.4);
\node at (2.5,2.7) {$\sigma_{Z, z_1, z_2}$};
\end{tikzpicture}
\end{center}
in $\WD\wz$, the structure map
\[\nicexy@C+.5cm{\sP_Z = \Histone(\bN^3,\bN^2)  \ar[r]^-{\sigma_{Z,z_1,z_2}} & \Histone(\bN^2,\bN^2) = \sP_W}\]
in \eqref{propagator-insplit} applied to the $1$-historical propagator $\ell \in \sP_Z$ is given as follows. For $\um = \bigl((m_1,m_1'),\ldots,(m_n,m_n')\bigr) \in \prof(\vwin) = \prof(\bN^2)$, we have
\[\begin{split}
\bigl(\sigma_{Z,z_1,z_2}\ell\bigr)(\um) &= \ell(\pi\um)\\
&= \ell\bigl((m_1,m_1,m_1'),\ldots,(m_n,m_n,m_n')\bigr)\\
&= \bigl((1,1), (m_1,m_1+m_1'),\ldots,(m_n,m_n+m_n')\bigr)
\end{split}\]
in $\prof(\vwout) = \prof(\bN^2)$.  For instance, we have
\[\bigl(\sigma_{Z,z_1,z_2}\ell\bigr)\bigl((2,5),(4,9),(3,7)\bigr) 
= \bigl((1,1),(2,7),(4,13),(3,10)\bigr).\]
\end{example}

\begin{example}\label{ex:ell-outsplit}
Suppose $V$ is a box with $\vout = \{v,v',z^2\}$ such that $\vout/(v=v') = \zout$, and $\vin = \zin$.  For the out-split (Def. \ref{def:out-split})
\begin{center}
\begin{tikzpicture}[scale=.9]
\draw [ultra thick] (2,1) rectangle (3,2);
\node at (2.3,1.5) {$Z$};
\node at (2.75,1.8) {\tiny{$z^{1}$}};
\node at (2.75,1.2) {\tiny{$z^{2}$}};
\draw [thick] (3,1.8) -- (3.3,1.8);
\draw [arrow] (3.3,1.8) to [out=0, in=180] (4.5,2);
\draw [arrow] (3.3,1.8) to [out=0, in=180] (4.5,1.6);
\draw [arrow] (.5,1.2) to (2,1.2);
\draw [arrow] (.5,1.5) to (2,1.5);
\draw [arrow] (.5,1.8) to (2,1.8);
\node at (0,1.5) {$\vin$};
\draw [arrow] (3,1.2) to (4.5,1.2);
\node at (5.5,1.5) {$\vout$};
\node at (4.7,2) {\tiny{$v$}};
\node at (4.7,1.6) {\tiny{$v'$}};
\draw [ultra thick] (1,0.7) rectangle (4,2.3);
\node at (2.5,2.7) {$\sigma^{V,v,v'}$};
\end{tikzpicture}
\end{center}
in $\WD\vz$, the structure map
\[\nicexy@C+.5cm{\sP_Z = \Histone(\bN^3,\bN^2)  \ar[r]^-{\sigma^{V,v,v'}} & \Histone(\bN^3,\bN^3) = \sP_V}\]
in \eqref{propagator-outsplit} applied to the $1$-historical propagator $\ell \in \sP_Z$ is given as follows. For $\um = \bigl((m_1,m_1',m_1''),\ldots,(m_n,m_n',m_n'')\bigr) \in \prof(\vvin) = \prof(\bN^3)$, we have
\[\begin{split}
\bigl(\sigma^{V,v,v'}\ell\bigr)(\um) &= \pi\ell(\um)\\
&= \pi\bigl((1,1), (m_1,m_1'+m_1''),\ldots,(m_n,m_n'+m_n'')\bigr)\\
&= \bigl((1,1,1), (m_1,m_1,m_1'+m_1''), \ldots, (m_n,m_n,m_n'+m_n'')\bigr)
\end{split}\]
in $\prof(\vvout) = \prof(\bN^3)$.  For instance, we have
\[\bigl(\sigma^{V,v,v'}\ell\bigr)\bigl((2,5,1),(4,9,10),(3,7,6)\bigr) 
= \bigl((1,1,1),(2,2,6),(4,4,19),(3,3,13)\bigr).\]
\end{example}

\begin{example}\label{ex:ell-wasted}
Suppose $U$ is a box such that $U \setminus u = Z$ for some $u \in \uin$ with $v(u) = (\bN,1)$.  For the $1$-wasted wire (Def. \ref{def:wasted-wire-wd})
\begin{center}
\begin{tikzpicture}[scale=.9]
\draw [ultra thick] (2,1) rectangle (3,2);
\node at (2.5,1.5) {$Z$};
\draw [arrow] (.5,1.2) to (2,1.2);
\draw [arrow] (.5,1.5) to (2,1.5);
\draw [arrow] (.5,1.8) to (2,1.8);
\node at (0,1.5) {$\uin$};
\draw [arrow] (.5,2.1) -- (1,2.1);
\node at (.7,2.3) {\tiny{$u$}};
\draw [arrow] (3,1.3) to (4.5,1.3);
\draw [arrow] (3,1.7) to (4.5,1.7);
\node at (5.1,1.5) {$\uout$};
\draw [ultra thick] (1,0.7) rectangle (4,2.3);
\node at (2.5,2.6) {$\omega_{U,u}$};
\end{tikzpicture}
\end{center}
in $\WD\UZ$, the structure map
\[\nicexy{\sP_Z = \Histone(\bN^3,\bN^2)  \ar[r]^-{\omega_{U,u}} & \Histone(\bN^4,\bN^2) = \sP_U}\]
in \eqref{propagator-wastedwire} applied to $\ell \in \sP_Z$ is given as follows.  For
\[\begin{split}
\um &= \left\{(m_i^1,m_i^2,m_i^3,m_i^4)\right\}_{i=1}^n\\
&\in \prof(\vuin) = \prof\bigl(v(u) \times v(z_1) \times v(z_2) \times v(z_3)\bigr) = \prof(\bN^4),
\end{split}\]
we have
\[\begin{split}
\bigl(\omega_{U,u} \ell\bigr)(\um) &= \ell\left(\um_{\setminus u}\right) 
= \ell\left\{(m_i^2,m_i^3,m_i^4)\right\}_{i=1}^n\\
&= \bigl((1,1),(m_1^2,m_1^3+m_1^4),\ldots,(m_n^2,m_n^3+m_n^4)\bigr)
\end{split}\]
in $\prof(\vuout) = \prof(\bN^2)$.  For instance, we have
\[\bigl(\omega_{U,u} \ell\bigr)\bigl((2,5,1,7),(4,9,10,2),(3,7,6,5)\bigr) 
= \bigl((1,1),(5,8),(9,12),(7,11)\bigr).\]
\end{example}

The following observation is the finite presentation theorem for the propagator algebra.

\begin{theorem}
\label{prop:propagator-algebra-is-algebra}
The propagator algebra in Def. \ref{def:propagator-algebra} is actually a $\WD$-algebra in the sense of  Def. \ref{def:wd-algebra},\index{finite presentation for the  propagator algebra} hence also in the sense of Def. \ref{def2:operad-algebra} by Theorem \ref{thm:wd-algebra}.
\end{theorem}

\begin{proof}
We need to check the $28$ generating axioms in Def. \ref{def:wd-algebra}.  The $8$ generating structure maps are all rather simple functions except for $\lambda_{X,x}$ \eqref{propagator-1loop}.  The only generating axioms that are not  obvious are the ones that involve a composition of two such generating structure maps, namely \eqref{wd-algebra-doubleloop} and \eqref{wd-algebra-loopelement}.  These two generating axioms are verified in Lemma \ref{lemma:propagator-doubleloop} and Lemma \ref{lemma:propagator-loopelement} below.
\end{proof}

In preparation for Lemma \ref{lemma:propagator-doubleloop}, we will need a few definitions and notations.  Recall that the generating axiom \eqref{wd-algebra-doubleloop} is really the $\WD$-algebra manifestation of Prop. \ref{prop:move:c1}.  The next definition is essentially the double-loop version of Def. \ref{def:generating-structure-1loop}.
 
\begin{definition}
\label{def:lambdatwo}
Suppose:
\begin{itemize}
\item
$X \in \boxs$, $x^1_- \not= x^2_- \in \xin$, and $x^1_+ \not= x^2_+ \in \xout$ such that $v(x^i_+) = v(x^i_-)$ as pointed sets for each $i \in \{1,2\}$.
\item
$X \setminus x^1$, $X \setminus x^2$, and $X \setminus x \in \boxs$ are obtained from $X$ by removing $x^1 = \{x^1_{\pm}\}$, $x^2 = \{x^2_{\pm}\}$, and $x = \{x^1_{\pm}, x^2_{\pm}\}$, respectively.
\end{itemize}
Suppose $\ut \in \prof\bigl(\vxout\bigr)$.
\begin{enumerate}[(i)]
\item
Write $v(\xplus) = \vxoneplus \times \vxtwoplus$ and $\vxminus = \vxoneminus \times \vxtwominus$.
\item
Write $\ut_{\xplus} \in \prof\bigl(\vxplus\bigr)$ for the profile obtained from $\ut$ by  taking only the $\vxplus$-entries.
\item
For $i\in \{1,2\}$, write $\ut_{\xiplus} \in \prof\bigl(\vxiplus\bigr)$ for the profile obtained from $\ut$ by  taking only the $\vxiplus$-entry.
\item
The profile $\uz_{\xiplus}$ is also defined as long as
\[
\uz \in \prof\left(\prod_{u \in T} v(u)\right)
\]
for some subset $\{\xoneplus, \xtwoplus\} \subseteq T \subseteq \xout$.
\item
Write $\ut_{\minusxplus} \in \prof\bigl(\vxminusxout\bigr)$ for the profile obtained from $\ut$ by removing the $\vxplus$-entries.
\item
For $i\in \{1,2\}$, write $\ut_{\minusxiplus} \in \prof\bigl(\vxminusxiout\bigr)$ for the profile obtained from $\ut$ by removing the $\vxiplus$-entry.
\end{enumerate} 
Suppose $g \in \sP_X$.  Define two functions
\begin{equation}
\label{propagator-2loop-defined}
\nicexy{
\prof\bigl(\vxminusxin\bigr) \ar[r]^-{\lambda^2 g} & \prof\bigl(\vxminusxout\bigr)}
\end{equation}
and
\begin{equation}
\label{propagator-2loop-G}
\nicexy{
\prof\bigl(\vxminusxin\bigr) \ar[r]^-{G^2_g} & \prof\bigl(v(x_+)\bigr)}
\end{equation}
with the properties
\begin{equation}
\label{2loopg-length}
|(\lambda^2 g)(?)| = |?| + 1 = |G^2_g(?)|
\end{equation}
inductively as follows.
\begin{enumerate}[(i)]
\item
For the empty profile, define
\begin{equation}
\label{2loopg-empty}
\begin{split}
(\lambda^2 g)(\varnothing) &= g(\varnothing)_{\minusxplus} \in \prof\bigl(\vxminusxout\bigr)\\
G^2_g(\varnothing) &= g(\varnothing)_{\xplus} \in  \prof\bigl(v(x_+)\bigr).
\end{split}
\end{equation}
In each definition in \eqref{2loopg-empty}, the first $\varnothing$ is the empty $\vxminusxin$-profile, and the second $\varnothing$, to which $g$ applies, is the empty $\vxin$-profile.  The profile $g(\varnothing)$ has length $1$ because $g \in \Histone\left(\vxin,\vxout\right)$.  So both $(\lambda^2 g)(\varnothing)$ and $G^2_g(\varnothing)$ have length $1$.
\item
Inductively, suppose $\uw \in \prof\bigl(\vxminusxin\bigr)$ has length $n \geq 1$.  Define
\begin{equation}
\label{2loopg-general}
\begin{split}
(\lambda^2 g)(\uw) &= g\bigl(\uw, G^2_g(\partial \uw)\bigr)_{\minusxplus} \in \prof\bigl(\vxminusxout\bigr)\\
G^2_g(\uw) &= g\bigl(\uw, G^2_g(\partial \uw)\bigr)_{\xplus} \in  \prof\bigl(v(x_+)\bigr).
\end{split}
\end{equation}
Here $\partial$ is the truncation \eqref{truncation}, so the profile
\[G^2_g(\partial \uw) \in \prof\bigl(v(\xplus)\bigr) =  \prof\bigl(v(\xminus)\bigr)\]
is already defined and has length $n$ by the induction hypothesis.  In each definition in \eqref{2loopg-general},
\[\bigl(\uw, G^2_g(\partial \uw)\bigr) \in \prof\left(\vxin\right)\]
has length $n$, so its image under $g$ has length $n+1$.  Therefore, both $(\lambda^2 g)(\uw)$ and $G^2_g(\uw)$ have length $n+1$.
\end{enumerate}
\end{definition}

\begin{example}
This is a continuation of Examples \ref{ex:ell} and \ref{ex2:ell}, where $Z$ is the box with $\zin = \{z_1,z_2,z_3\}$, $\zout = \{z^1,z^2\}$, and all $v(-) = (\bN,1)$. The $1$-historical propagator $\ell \in \sP_Z = \Histone(\bN^3,\bN^2)$ is defined as 
\[\ell\bigl((m_1,m_1',m_1''), \ldots, (m_n,m_n',m_n'')\bigr) = \bigl((1,1),(m_1,m_1'+m_1''),\ldots,(m_n, m_n'+m_n'')\bigr).\]  
With $w^1 = \{z_1,z^1\}$ and $w^2 = \{z_2,z^2\}$, suppose  $Z \setminus w$ is the box obtained from $Z$ by removing $\{z_1,z_2,z^1,z^2\}$.  Then 
\[\vzminuswin = v(z_3) = \bN \andspace \vzminuswout = *,\]
where $*$ here is the one-point set ($=$ empty product).  
\begin{center}
\begin{tikzpicture}
\draw [ultra thick] (2,1) rectangle (3,2);
\node at (2.2,1.8) {\tiny{$z_1$}};
\node at (2.8,1.8) {\tiny{$z^1$}};
\node at (2.2,1.5) {\tiny{$z_2$}};
\node at (2.8,1.4) {\tiny{$z^2$}};
\node at (2.2,1.15) {\tiny{$z_3$}};
\draw [arrow] (.5,1.2) to (2,1.2);
\node at (-.5,1.6) {$\zminuswin$};
\node at (5.6,1.6) {$\zminuswout$};
\draw [arrow, looseness=3] (3,1.8) to [out=30, in=150] (2,1.8);
\draw [arrow, looseness=8] (3,1.5) to [out=30, in=150] (2,1.5);
\draw [ultra thick] (1,0.7) rectangle (4,3);
\node at (2.5,3.4) {$\lambda^2 \in \WD\zminuswz$};
\end{tikzpicture}
\end{center}
The functions
\[\nicexy{\prof(\bN) \ar[r]^-{\lambda^2\ell} & \prof(*)} \andspace 
\nicexy{\prof(\bN) \ar[r]^-{G^2_{\ell}} & \prof(\bN\times\bN)}\]
in \eqref{propagator-2loop-defined} and \eqref{propagator-2loop-G} are given as follows.  For $\um = (m_1, \ldots , m_n) \in \prof(\bN)$, a simple induction shows that
\[\begin{split}
(\lambda^2\ell)(\um) &= (*,\ldots,*),\\
(G^2_{\ell})(\um) &= \bigl((1,1), (1,1+m_1), (1,1+m_1+m_2),\ldots,(1,1+m_1+\cdots+m_n)\bigr),
\end{split}\]
where on the right side of $(\lambda^2\ell)(\um)$ there are $n+1$ copies of $*$.  For instance, we have
\[(G^2_{\ell})(6,3,2,9) = \bigl((1,1),(1,7),(1,10),(1,12),(1,21)\bigr).\]
\end{example}

The generating axiom \eqref{wd-algebra-doubleloop} for the propagator algebra $\sP$ claims that the diagram
\begin{equation}
\label{propagator-doubleloop-diagram}
\nicexy@C+.4cm{
\sP_X \ar[r]^-{\lambda_{X,x^2}} \ar[d]_{\lambda_{X,x^1}} 
& \sP_{X \setminus x^2} \ar[d]^-{\lambda_{X\setminus x^2, x^1}} \\
\sP_{X \setminus x^1} \ar[r]^-{\lambda_{X \setminus x^1, x^2}} & \sP_{X \setminus x}
}\end{equation}
is commutative.  We will consider the top-right composition, so let us use the abbreviations
\begin{equation}
\label{lambda-abbreviations}
\lambda_2 = \lambda_{X,x^2} \andspace \lambda_1 = \lambda_{X\setminus x^2, x^1}.
\end{equation}
For the proof of Lemma \ref{lemma:propagator-doubleloop}, we will need the following observation.  It provides an explicit formula for the function $G^2_g$ \eqref{propagator-2loop-G} in terms of $G_g$ (defined with respect to $x^2$) and $G_{\lambda_2 g}$ \eqref{propagator-G} (defined with respect to $x^1$) for each  $g \in \sP_X$.

\begin{lemma}
\label{lemma:G2}
In the context of Def. \ref{def:lambdatwo}, suppose $g \in \sP_X$ and $\uw \in \prof\bigl(\vxminusxin\bigr)$ with length at least $1$.  Then the following equalities hold.
\begin{equation}
\label{G2-formula}
\begin{split}
G^2_g(\partial\uw)_{\xoneplus} &= G_{\lambda_2 g}(\partial \uw) \in \prof\bigl(\vxoneplus\bigr)\\
G^2_g(\partial \uw)_{\xtwoplus} &= G_g\bigl(\partial \uw, \partial G_{\lambda_2 g}(\partial \uw)\bigr) \in \prof\bigl(\vxtwoplus\bigr)
\end{split}
\end{equation}
In the above equalities:
\begin{enumerate}
\item
$\partial$ is the truncation \eqref{truncation}.
\item 
$\lambda_2 g \in \sP_{X \setminus x^2}$ by Lemma \ref{lemma:lambdag-historical}(2).
\item 
The function 
\[
\nicexy{\prof\bigl(\vxminusxin\bigr) \ar[r]^-{G_{\lambda_2 g}} & \prof\bigl(\vxoneplus\bigr)}
\]
is defined with respect to $x^1=\{\xoneplus, \xoneminus\}$ \eqref{propagator-G}.
\item 
The function 
\[
\nicexy{\prof\bigl(\vxminusxtwoin\bigr) \ar[r]^-{G_g} & \prof\bigl(\vxtwoplus\bigr)}
\]
is defined with respect to $x^2=\{\xtwoplus,\xtwominus\}$.
\end{enumerate}
\end{lemma}

\begin{proof}
The proof of \eqref{G2-formula} is by induction on $|\uw| \geq 1$.  If $|\uw|=1$, then $\partial \uw = \varnothing$.  So by \eqref{lambdag-empty} and \eqref{2loopg-empty} we have
\[
G_{\lambda_2 g}(\varnothing) 
= (\lambda_2 g)(\varnothing)_{\xoneplus} = \bigl[g(\varnothing)_{\minusxtwoplus}\bigr]_{\xoneplus}
= g(\varnothing)_{\xoneplus} = G^2_g(\varnothing)_{\xoneplus}.
\]
Likewise, we have
\[
G_g(\varnothing) = g(\varnothing)_{\xtwoplus} = G^2_g(\varnothing)_{\xtwoplus}.
\]
This proves the initial case $|\uw| = 1$.

For the induction step, suppose $|\uw| \geq 2$.  For the first equality in \eqref{G2-formula} we have:
\[
\begin{split}
G_{\lambda_2 g}(\partial \uw) 
&= (\lambda_2 g)\Bigl(\partial \uw, G_{\lambda_2 g}(\partial^2 \uw)\Bigr)_{\xoneplus} \quad\text{by \eqref{lambdag-general}}\\
&= g\Bigl(\partial \uw, G_{\lambda_2 g}(\partial^2 \uw), G_g\bigl(\partial^2 \uw, \partial G_{\lambda_2 g}(\partial^2 \uw)\bigr)\Bigr)_{\xoneplus} \quad\text{by \eqref{lambdag-general}}\\
&= g\Bigl(\partial \uw, G^2_g(\partial^2 \uw)_{\xoneplus}, G^2_g(\partial^2 \uw)_{\xtwoplus}\Bigr)_{\xoneplus} \quad\text{by induction hypothesis}\\
&= G^2_g(\partial \uw)_{\xoneplus} \quad \text{by \eqref{2loopg-general}}.
\end{split}
\]
In the second equality above, we used the fact that
\[
\left\{g(\cdots)_{\minusxtwoplus}\right\}_{\xoneplus} = g(\cdots)_{\xoneplus}.
\]

For the second equality in \eqref{G2-formula} we have:
\[
\begin{split}
&G_g\bigl(\partial \uw, \partial G_{\lambda_2 g}(\partial \uw)\bigr)\\
&= G_g\bigl(\partial \uw, G_{\lambda_2 g}(\partial^2 \uw)\bigr) \quad\text{by Lemma \ref{lemma:lambdag-historical}(1)}\\
&= g\Bigl(\partial \uw, G_{\lambda_2 g}(\partial^2 \uw), G_g\bigl(\partial^2 \uw, \partial G_{\lambda_2 g}(\partial^2 \uw)\bigr)\Bigr)_{\xtwoplus} \quad\text{by \eqref{lambdag-general}}\\
&= g\Bigl(\partial \uw, G^2_g(\partial^2 \uw)_{\xoneplus}, G^2_g(\partial^2 \uw)_{\xtwoplus}\Bigr)_{\xtwoplus} \quad\text{by induction hypothesis}\\
&= G^2_g(\partial \uw)_{\xtwoplus} \quad \text{by \eqref{2loopg-general}}.
\end{split}
\]
This finishes the induction.
\end{proof}

\begin{lemma}
\label{lemma:propagator-doubleloop} 
The propagator algebra $\sP$ in Def. \ref{def:propagator-algebra} satisfies the generating axiom \eqref{wd-algebra-doubleloop}; i.e., the diagram \eqref{propagator-doubleloop-diagram} is commutative.
\end{lemma}

\begin{proof}
We will use the abbreviations \eqref{lambda-abbreviations}.  By the symmetry between $x^1 = \{x^1_{\pm}\}$ and $x^2 = \{x^2_{\pm}\}$, it suffices to show that
\begin{equation}
\label{lambda1-lambda2}
\bigl(\lambda_1 \lambda_2 g\bigr)(\uw) = (\lambda^2 g)(\uw)
\end{equation}
for $g \in \sP_X$ and $\uw \in \prof\bigl(\vxminusxin\bigr)$, where $\lambda^2 g$ is defined in \eqref{propagator-2loop-defined}.  We prove \eqref{lambda1-lambda2} by induction on the length $|\uw|$.

If $|\uw| = 0$, then by \eqref{lambdag-empty} and \eqref{2loopg-empty} the left side of \eqref{lambda1-lambda2} is:
\[
\bigl(\lambda_2 g)(\varnothing)_{\minusxoneplus} 
= \left[g(\varnothing)_{\minusxtwoplus}\right]_{\minusxoneplus}
= g(\varnothing)_{\minusxplus} = (\lambda^2 g)(\varnothing).
\]
For the induction step, suppose $|\uw| \geq 1$.  Then the left side of \eqref{lambda1-lambda2} is:
\[
\begin{split}
& (\lambda_2 g)\bigl(\uw, G_{\lambda_2 g}(\partial \uw)\bigr)_{\minusxoneplus} \quad\text{by \eqref{lambdag-general}}\\
&=  g\Bigl(\uw, G_{\lambda_2 g}(\partial \uw), G_g\bigl(\partial \uw, \partial G_{\lambda_2 g}(\partial \uw)\bigr)\Bigr)_{\minusxplus} \quad\text{by \eqref{lambdag-general}}\\
&= g\Bigl(\uw, G^2_g(\partial \uw)_{\xoneplus}, G^2_g(\partial \uw)_{\xtwoplus}\Bigr)_{\minusxplus} \quad \text{by \eqref{G2-formula}}\\
&= (\lambda^2 g)(\uw) \quad\text{by \eqref{2loopg-general}}.
\end{split}
\]
This finishes the induction.
\end{proof}

The proof of the generating axiom \eqref{wd-algebra-loopelement} in Lemma \ref{lemma:propagator-loopelement} below will use the following observation.  Recall that the generating axiom \eqref{wd-algebra-loopelement} is really a $\WD$-algebra manifestation of the generating relation \eqref{move:c5}.

\begin{lemma}
In the context of Prop. \ref{prop:move:c5}, recall that
\[
X = \frac{Y}{(x^1 = x^2)}, \quad X' = \frac{X}{(x_1 = x_2)}, \andspace X^* = X \setminus \{x^{12}, x_1,x_2\}.
\]
Consider the maps:
\[
\nicexy@C+1cm{
\sP_X \ar[d]_-{\sigmasubstar \,=\, \sigma_{X,x_1,x_2}} \ar[r]^-{\sigmasupstar \,=\, \sigma^{Y,x^1,x^2}} & \sP_Y\\
\sP_{X'} &
}\]
Then for $g \in \sP_X$ and $\uw \in \prof\bigl(\vxstarin\bigr)$ with length $\geq 1$, the equality
\begin{equation}
\label{gtwo-sigmag}
G^2_{\sigmasupstar g}(\partial\uw) = \Bigl(G_{\sigmasubstar g}(\partial\uw), G_{\sigmasubstar g}(\partial\uw)\Bigr)
\end{equation}
holds.  Here:
\begin{enumerate}
\item
$\sigmasupstar g \in \sP_Y$ and $\sigmasubstar g \in \sP_{X'}$.
\item
The function
\[
\nicexy{\prof\Bigl(\vxstarin\Bigr) \ar[r]^-{G^2_{\sigmasupstar g}} & \prof\bigl(v(x^1) \times v(x^2)\bigr)
}\]
is defined as in \eqref{propagator-2loop-G}, starting with the box $Y$ and the wires $x^1 \not= x^2 \in \yout$ and $x_1 \not= x_2 \in \yin$.
\item
The function
\[
\nicexy{
\prof\Bigl(\vxstarin\Bigr) \ar[r]^-{G_{\sigmasubstar g}} & \prof\bigl(v(x^{12})\bigr)
}\]
is defined as in \eqref{propagator-G}, starting with the box $X'$ and the wires $x_{12} \in \xprimein$ and $x^{12} \in \xprimeout$.
\end{enumerate}
\end{lemma}

\begin{proof}
The proof is by induction on the length of $\uw$.  If $|\uw| = 1$, the $\partial\uw = \varnothing$.  So we have:
\[
\begin{split}
G^2_{\sigmasupstar g} (\varnothing)
&= (\sigmasupstar g)(\varnothing)_{\{x^1,x^2\}} \quad\text{by \eqref{2loopg-empty}}\\
&= \Bigl(g(\varnothing)_{x^{12}}, g(\varnothing)_{x^{12}}\Bigr) \quad\text{by \eqref{propagator-outsplit}}\\
&= \Bigl((\sigmasubstar g)(\varnothing)_{x^{12}}, (\sigmasubstar g)(\varnothing)_{x^{12}}\Bigr) \quad\text{by \eqref{propagator-insplit}}\\
&= \Bigl(G_{\sigmasubstar g}(\varnothing), G_{\sigmasubstar g}(\varnothing)\Bigr) \quad\text{by \eqref{lambdag-empty}}
\end{split}
\]

For the induction step, suppose $|\uw| \geq 2$.  Then we have:
\[
\begin{split}
G^2_{\sigmasupstar g}(\partial\uw) 
&= (\sigmasupstar g)\Bigl(\partial\uw, G^2_{\sigmasupstar g}(\partial^2 \uw)\Bigr)_{\{x^1,x^2\}} \quad\text{by \eqref{2loopg-general}}\\
&=  (\sigmasupstar g)\Bigl(\partial\uw, G_{\sigmasubstar g}(\partial^2\uw),  G_{\sigmasubstar g}(\partial^2\uw)\Bigr)_{\{x^1,x^2\}} \quad\text{by induction hypothesis}\\
&= \Bigl(g\bigl(\partial\uw, G_{\sigmasubstar g}(\partial^2\uw),  G_{\sigmasubstar g}(\partial^2\uw)\bigr)_{x^{12}},~ \text{same} \Bigr) \quad\text{by \eqref{propagator-outsplit}}\\
&= \Bigl((\sigmasubstar g)\bigl(\partial\uw, G_{\sigmasubstar g}(\partial^2\uw)\bigr)_{x^{12}},~ \text{same}\Bigr) \quad\text{by \eqref{propagator-insplit}}\\
&=  \Bigl(G_{\sigmasubstar g}(\partial\uw), G_{\sigmasubstar g}(\partial\uw)\Bigr) \quad\text{by \eqref{lambdag-general}}
\end{split}
\]
This finishes the induction.
\end{proof}

\begin{lemma}
\label{lemma:propagator-loopelement} 
The propagator algebra $\sP$ in Def. \ref{def:propagator-algebra} satisfies the generating axiom  \eqref{wd-algebra-loopelement}; i.e., the diagram
\[
\nicexy@C+1cm{
\sP_X \ar[d]_-{\sigmasubstar \,=\, \sigma_{X,x_1,x_2}} \ar[r]^-{\sigmasupstar \,=\, \sigma^{Y,x^1,x^2}} & \sP_Y \ar[r]^-{\lambda_{(1)} \,=\, \lambda_{Y,x(1)}} & \sP_{Y \setminus x(1)} \ar[d]^-{\lambda_{(2)} \,=\, \lambda_{Y\setminus x(1), x(2)}}\\
\sP_{X'} \ar[rr]^-{\lambda \,=\, \lambda_{X',x}} && \sP_{X^*}
}\]
is commutative.
\end{lemma}

\begin{proof}
For $g \in \sP_X$ and $\uw \in \prof\bigl(\vxstarin\bigr)$, we will prove the desired equality
\begin{equation}
\label{lambda-sigmastar-g}
(\lambda\sigmasubstar g)(\uw) = \bigl(\lambda_{(2)} \lambda_{(1)} \sigmasupstar g\bigr)(\uw) \in \prof\Bigl(\vxstarout\Bigr)
\end{equation}
by induction on the length of $\uw$.  If $|\uw| = 0$, then both sides of \eqref{lambda-sigmastar-g} are equal to $g(\varnothing)_{\minusxsuponetwo}$.

For the induction step, suppose $|\uw| \geq 1$.  Then we have:
\[
\begin{split}
(\lambda\sigmasubstar g)(\uw) 
&= (\sigmasubstar g)\Bigl(\uw, G_{\sigmasubstar g}(\partial\uw)\Bigr)_{\minusxsuponetwo} \quad\text{by \eqref{lambdag-general}}\\
&= g\Bigl(\uw, G_{\sigmasubstar g}(\partial\uw), G_{\sigmasubstar g}(\partial\uw)\Bigr)_{\minusxsuponetwo} \quad\text{by \eqref{propagator-insplit}}\\
&= g\Bigl(\uw, G^2_{\sigmasupstar g}(\partial\uw)\Bigr)_{\minusxsuponetwo} \quad\text{by \eqref{gtwo-sigmag}}\\
&= (\sigmasupstar g)\Bigl(\uw, G^2_{\sigmasupstar g}(\partial\uw)\Bigr)_{\minusxonextwo} \quad\text{by \eqref{propagator-outsplit}}\\
&= \bigl(\lambda_{(2)} \lambda_{(1)} \sigmasupstar g\bigr)(\uw) \quad\text{by  \eqref{2loopg-general} and \eqref{lambda1-lambda2}}
\end{split}
\]
This finishes the induction.
\end{proof}

\begin{remark}
\label{rk:propagator-agree}
To see that our definition of the propagator algebra $\sP$ in Def. \ref{def:propagator-algebra} agrees with the one in \cite{rupel-spivak} (Section 3), recall that our version of the propagator algebra is based on  Def. \ref{def:wd-algebra}.  On the other hand, 
the propagator algebra in \cite{rupel-spivak} is based on Def. \ref{def1:operad-algebra}, which is equivalent to Def. \ref{def2:operad-algebra}.  A direct inspection of \cite{rupel-spivak} (Announcement 3.3.3 and Eq. (17)) reveals that their structure map of $\sP$, when applied to the generating wiring diagrams (section \ref{sec:generating-wd}), reduces to our $8$ generating structure maps in  Def. \ref{def:propagator-algebra}.  Theorem \ref{thm:wd-algebra} then guarantees that the two definitions are equivalent.
\end{remark}

\section{Algebras over the Operad of Normal Wiring Diagrams}
\label{sec:algebra-normal-wd}

The purpose of this section is to provide a finite presentation for algebras over the $\boxs$-colored operad $\wddot$ of normal wiring diagrams (Prop. \ref{prop:without-dn-operads}).  We begin by defining these algebras in terms of finitely many generators and relations.  Recall from Def. \ref{def:generating-without-dn} that a \emph{normal generating wiring diagram} is a generating wiring diagram that is not a $1$-delay node $\delta_d$.

\begin{definition}
\label{def:normal-algebra}
A $\wddot$-algebra $A$ consists of the following data.\index{wddot-algebra@$\wddot$-algebra} 
\begin{enumerate}
\item
For each $X \in \boxs$, $A$ is equipped with a class $A_X$ called the \emph{$X$-colored entry} of $A$.
\item
It is equipped with the $7$ \emph{generating structure maps} in Def. \ref{def:wd-algebra} corresponding to the normal generating wiring diagrams.
\end{enumerate}
This data is required to satisfy the same $28$ generating axioms in Def. \ref{def:wd-algebra}.
\end{definition}

The next observation is the $\wddot$ version of Theorem \ref{thm:wd-algebra}.   It  guarantees that the two existing definitions of a $\wddot$-algebra are equivalent.  The first one (Def. \ref{def2:operad-algebra}) is in terms of a general structure map satisfying an associativity axiom for a general operadic composition.  The other one (Def. \ref{def:normal-algebra}) is in terms of $7$ generating structure maps and $28$ generating axioms regarding the normal generating wiring diagram.  Therefore, algebras over $\wddot$  have a finite presentation.

\begin{theorem}
\label{thm:normal-wd-algebra}
For the operad $\wddot$ of normal wiring diagrams (Prop.  \ref{prop:without-dn-operads}), Def. \ref{def2:operad-algebra} with $\O=\wddot$ and Def. \ref{def:normal-algebra} of a $\wddot$-algebra are equivalent.\index{finite presentation for WDdot algebra@finite presentation for $\wddot$-algebras}
\end{theorem}

\begin{proof}
Simply restrict the proof of Theorem \ref{thm:wd-algebra} to normal (generating) wiring diagrams.  Instead of Theorem \ref{stratified-presentation-exists}, here we use Theorem \ref{thm:without-dn-coherence} for the existence of a presentation involving only normal generating wiring diagrams.
\end{proof}

\begin{remark}
\label{rk:algebra-ds}
 In \cite{spivak15b} (Def. 4.1--4.4) several closely related $\wddot$-algebras were defined, although they appeared in the language of symmetric monoidal categories.  By Theorem \ref{thm:normal-wd-algebra} each of these $\wddot$-algebras has a finite presentation with $7$ generating structure maps and $28$ generating axioms as in Def. \ref{def:normal-algebra}.  In Section \ref{sec:discrete-systems} we will discuss one of these $\wddot$-algebras and its finite presentation.  In Section \ref{sec:algebra-ods} we will discuss a similar algebra of open dynamical systems over the operad $\wdzero$ of strict wiring diagrams.  Essentially the same formalism applies to the other $\wddot$-algebras in \cite{spivak15b}.
\end{remark}

\section{Finite Presentation for the Algebra of Discrete Systems}
\label{sec:discrete-systems}

The purpose of this section is to provide a finite presentation for the algebra of discrete systems introduced in \cite{spivak15b} (Definition 4.9).  Let us first recall some definitions from \cite{spivak15b} (Sections 2.1 and 4.1).

\begin{assumption}
Throughout this section, $S$ denotes the class of sets.  So $\boxs = \Box_{\set}$, and $\wddot$ is the $\Box_{\set}$-colored operad of normal wiring diagrams (Prop. \ref{prop:without-dn-operads}).  
\end{assumption}

\begin{definition}
\label{def:discrete-systems}
Suppose $A$ and $B$ are sets.  An \emph{$(A,B)$-discrete system}\index{discrete system} is a triple $(T, \frd,\fup)$ consisting of:
\begin{enumerate}
\item a set $T$, called the \emph{state set};\index{state set}
\item a function $\frd : T \to B$, called the \emph{readout function};\index{readout function}
\item a function $\fup : A \times T \to T$, called the \emph{update function}.\index{update function}
\end{enumerate}
\end{definition}

\begin{definition}
\label{def:x-discrete-systems}
Suppose $X = (\xin, \xout) \in \boxs$ is a box.  An \emph{$X$-discrete system} is an $(\vxin, \vxout)$-discrete system, where
\[\vxin = \prod_{x \in \xin} v(x) \andspace 
\vxout = \prod_{x \in \xout} v(x) \in \set\]
as in \eqref{vxin-vxout} (but with sets instead of pointed sets).  In other words, an $X$-discrete system is a triple $(T,\frd,\fup)$ such that $T$ is a set and that
\[\nicexy{T \ar[r]^-{\frd} & \vxout} \andspace 
\nicexy{\vxin \times T \ar[r]^-{\fup} & T}\]
are functions.  The collection of all $X$-discrete systems is denoted by $\DS(X)$.\label{notation:ds}
\end{definition}

\begin{example}
\label{ex:ds-empty}
If $X = \varnothing \in \boxs$ is the empty box, then $\vxin = \vxout = *$ by convention.  A readout function $\frd : T \to *$ gives no information, and $* \times T \cong T$.  So 
\begin{equation}
\label{ds-empty}
\DS(\varnothing) = \Bigl\{(T,\fup) : T \in \set,\, \fup : T \to T \text{ a function}\Bigr\}.
\end{equation}
In particular, the collection $\DS(\varnothing)$ is \emph{not} a set but a proper class.  This example explains why in Def. \ref{def1:operad-algebra} we defined an entry of an operad algebra to be a class and not a set. 
\end{example}

\begin{example}\label{ex:dsx}
Suppose $X$ is the box with $\xin = \{x_1,x_2\}$, $\xout = \{x^1,x^2\}$, and values
\[v(x_1) = \{a_1,b_1\},\quad v(x_2) = \{a_2,b_2\}, \quad v(x^1) = \{a^1,b^1\}, \andspace v(x^2) = \{a^2,b^2\}.\]
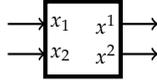
\begin{center}
\begin{tikzpicture}[scale=1]
\draw [ultra thick] (0,0) rectangle (1,1);
\node at (.2,.7) {\scriptsize{$x_1$}};
\node at (.2,.3) {\scriptsize{$x_2$}};
\node at (.8,.7) {\scriptsize{$x^1$}};
\node at (.8,.3) {\scriptsize{$x^2$}};
\draw [arrow] (-.5,.3) to (0,.3);
\draw [arrow] (-.5,.7) to (0,.7);
\draw [arrow] (1,.3) to (1.5,.3);
\draw [arrow] (1,.7) to (1.5,.7);
\end{tikzpicture}
\end{center}
An $X$-discrete system is an $(\vxin,\vxout)$-discrete system, where
\[\begin{split}
\vxin &= v(x_1) \times v(x_2) = \{a_1,b_1\} \times \{a_2,b_2\};\\
\vxout &= v(x^1) \times v(x^2) = \{a^1,b^1\} \times \{a^2,b^2\}.
\end{split}\]
Suppose $T = \{1,2\}$ is the state set.  There is an $X$-discrete system $(T,\frd,\fup)$ with the readout function $\frd : T \to \vxout$ and update function $\fup : \vxin \times T \to T$ defined as follows.
\begin{align*}
\frd(1) &= (a^1,a^2)& \fup\bigl((a_1,a_2),1\bigr) &=1 & \fup\bigl((a_1,a_2),2\bigr) &=1\\
\frd(2) &= (b^1,a^2)& \fup\bigl((a_1,b_2),1\bigr) &=2 & \fup\bigl((a_1,b_2),2\bigr) &=1\\
&& \fup\bigl((b_1,a_2),1\bigr) &=2 & \fup\bigl((b_1,a_2),2\bigr) &=2\\
&& \fup\bigl((b_1,b_2),1\bigr) &=1 & \fup\bigl((b_1,b_2),2\bigr) &=2\\
\end{align*}
Visually it can also be represented by the \index{transition diagram}transition diagram:
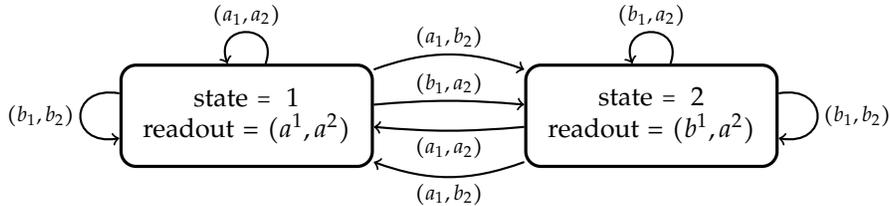
\begin{center}
\begin{tikzpicture}[auto, node distance=2cm,font=\small,
    every node/.style={inner sep=2pt, rectangle, rounded corners=0.2cm, text centered},
    state/.style={draw, very thick, minimum height=3.5em, text width=3.2cm}]
\matrix[row sep=.8cm,column sep=2cm] {
\node [state] (1) {state $= 1$\\ readout = $(a^1,a^2)$}; &
\node [state] (2) {state $= 2$\\ readout = $(b^1,a^2)$};\\};
\draw [arrow, out=70, in=110, looseness=3] (1) to node[swap]{\scriptsize{$(a_1,a_2)$}} (1);
\draw [arrow, out=170, in=190, looseness=3] (1) to node[swap]{\scriptsize{$(b_1,b_2)$}} (1);
\draw [arrow, bend left=20] (1) to node{\scriptsize{$(a_1,b_2)$}} (2);
\draw [arrow, bend left=5] (1) to node{\scriptsize{$(b_1,a_2)$}} (2);
\draw [arrow, out=70, in=110, looseness=3] (2) to node[swap]{\scriptsize{$(b_1,a_2)$}} (2);
\draw [arrow, out=10, in=-10, looseness=3] (2) to node{\scriptsize{$(b_1,b_2)$}} (2);
\draw [arrow, bend left=5] (2) to node{\scriptsize{$(a_1,a_2)$}} (1);
\draw [arrow, bend left=20] (2) to node{\scriptsize{$(a_1,b_2)$}} (1);
\end{tikzpicture}
\end{center}
For example, the arrow labeled $(b_1,a_2)$ from the box for state $1$ to the box for state $2$ represents the value $\fup\bigl((b_1,a_2),1\bigr) = 2$, and likewise for the other arrows.
\end{example}

We now define the algebra of discrete systems in terms of $7$ very simple generating structure maps.

\begin{definition}
\label{def:algebra-discrete-systems}
The \emph{algebra of discrete systems}\index{algebra of discrete systems} is the $\wddot$-algebra $\DS$ in the sense of Def. \ref{def:normal-algebra} defined as follows.  For each $X \in \boxs$, the $X$-colored entry is the class $\DS(X)$ of $X$-discrete systems in Def. \ref{def:x-discrete-systems}

The $7$ generating structure maps--as in Def. \ref{def:wd-algebra} but without $\delta_d$--are defined as follows.
\begin{enumerate}
\item Corresponding to the empty wiring diagram $\epsilon \in \wddot\emptynothing$ (Def. \ref{def:empty-wd}), the chosen element in $\DS(\varnothing)$ \eqref{ds-empty} is the pair $(*, \Id)$ with:
\begin{itemize}
\item the one-point set $*$ as its state set;
\item the identity map as its update function.
\end{itemize}
\item Corresponding to each name change $\tau_{X,Y} \in \wddot\yx$ (Def. \ref{def:name-change}), the structure map
\[\nicexy{\DS(X) \ar[r]^-{\tau_{X,Y}}_-{=} & \DS(Y)}\]
is the identity map, using the fact that $\vxin = \vyin$ and $\vxout = \vyout$.
\item Corresponding to a $2$-cell $\theta_{X,Y} \in \wddot\xplusyxy$ (Def. \ref{def:theta-wd}), it has the structure map
\[\nicexy{\DS(X) \times \DS(Y) \ar[d]_-{\theta_{X,Y}} & \Bigl((T_X,\frd_X, \fup_X), (T_Y,\frd_Y, \fup_Y)\Bigr) \ar@{|->}[d]\\
\DS(\xplusy) & \Bigl(T_X \times T_Y, \frd_X \times \frd_Y, \fup_X \times \fup_Y\Bigr)}\]
using the fact that
\[\vxyin = \vxin \times \vyin \andspace \vxyout = \vxout \times \vyout.\]
\item Corresponding to a $1$-loop $\lambda_{X,x} \in \wddot\xminusxx$  (Def. \ref{def:loop-wd}) with $x = \bigl(x_- , x_+\bigr) \in \xin \times \xout$ (so $v(x_+) = v(x_-)$), it has the structure map
\[\nicexy{\DS(X) \ar[r]^-{\lambda_{X,x}} & \DS(\xminusx)} \quad
\nicexy{\bigl(T, \frd, \fup\bigr) \ar@{|->}[r] & \bigl(T,\frdminusx, \fupminusx\bigr)}\]
in which, for $(\uy,t) \in \vxminusxin \times T$:
\[\begin{split}
\frd(t) &= \Bigl(\frd(t)_{\minusxplus}, \frd(t)_{x_+}\Bigr) \in \vxout = \vxminusxout \times v(x_+);\\
\frdminusx(t) &= \frd(t)_{\minusxplus} \in \vxminusxout;\\  
\fupminusx\bigl(\uy, t\bigr) &= \fup\Bigl(\bigl(\uy, \frd(t)_{x_+}\bigr) , t\Bigr).\\
\end{split}\]
\item Corresponding to an in-split $\sigma_{X,x_1,x_2} \in \wddot\yx$ (Def. \ref{def:insplit-wd}) with $v(x_1) = v(x_2)$ and $Y = X/(x_1 = x_2)$, it has the structure map
\[\nicexy{\DS(X)  \ar[r]^-{\sigma_{X,x_1,x_2}} & \DS(Y)}
\quad \nicexy{\bigl(T, \frd, \fup\bigr) \ar@{|->}[r] & \bigl(T, \frd, \sigmasubdot\fup\bigr)}\]
in which the update function is
\[\bigl(\sigmasubdot\fup\bigr)\bigl(\uy,t\bigr) = \fup\Bigl(\sigmasubdot\uy, t\Bigr)\]
for $(\uy,t) \in \vyin \times T$.  Here $\sigmasubdot\uy \in \vxin$ is obtained from $\uy$ by using the $v(x_{12})$-entry of $\uy$ in both the $v(x_1)$-entry and the $v(x_2)$-entry, where $x_{12} \in \yin$ is the identified element of $x_1$ and $x_2$.
\item Corresponding to an out-split  $\sigma^{Y,y_1,y_2} \in \wddot\yx$ (Def. \ref{def:out-split}) with $v(y_1) = v(y_2)$ and $X = Y/(y_1 = y_2)$,  it has the structure map
\[\nicexy{\DS(X)  \ar[r]^-{\sigma^{Y,y_1,y_2}} & \DS(Y)} \quad
\nicexy{\bigl(T, \frd, \fup\bigr) \ar@{|->}[r] & \bigl(T, \sigmasupdot\frd, \fup\bigr)}\]
in which the readout function is
\[\bigl(\sigmasupdot\frd\bigr)(t) = \sigmasupdot\bigl(\frd(t)\bigr) \in \vyout\]
for $t \in T$.  Here $\sigmasupdot\bigl(\frd(t)\bigr)$ is obtained from $\frd(t) \in \vxout$ by using its $v(y_{12})$-entry in both the $v(y_1)$-entry and the $v(y_2)$-entry, where $y_{12} \in \xout$ is the identified element of $y_1$ and $y_2$.
\item Corresponding to a $1$-wasted wire  $\omega_{Y,y} \in \wddot\yx$ (Def. \ref{def:wasted-wire-wd}) with $y \in \yin$ and $X = Y \setminus y$,  it has the structure map
\[\nicexy{\DS(X)  \ar[r]^-{\omega_{Y,y}} & \DS(Y)} \quad
\nicexy{\bigl(T, \frd, \fup\bigr) \ar@{|->}[r] & \bigl(T, \frd, \omega\fup\bigr)}\]
in which the update function is
\[\bigl(\omega\fup\bigr)(\uz,t) = \fup\bigl(\uz_{\setminus y}, t\bigr)\]
for $(\uz,t) \in \vyin \times T$.  Here $\uz_{\setminus y} \in \vxin$ is obtained from $\uz$ by removing the $v(y)$-entry.
\end{enumerate}
This finishes the definition of the $\wddot$-algebra of discrete systems.
\end{definition}

\begin{remark}
In \cite{spivak15b} (Example 2.7) the image of $\theta_{X,Y}$ is called the \emph{parallel composition}.  The structure map $\lambda_{X,x}$ corresponding to a $1$-loop was discussed in \cite{spivak15b} (Example 2.9).  The structure maps $\sigma_{X,x_1,x_2}$ and $\sigma^{Y,y_1,y_2}$ corresponding to an in-split and an out-split were discussed in \cite{spivak15b} (Example 2.8).
\end{remark}

The following four examples refer to the $X$-discrete system $(T,\frd,\fup)$ in Example \ref{ex:dsx}, where $\xin = \{x_1,x_2\}$, $\xout = \{x^1,x^2\}$, $v(x_i) = \{a_i,b_i\}$, and $v(x^i) = \{a^i,b^i\}$ for $i=1,2$.  Let us consider the effects of some of the structure maps in Def. \ref{def:algebra-discrete-systems} on $(T,\frd,\fup) \in \DS(X)$.

\begin{example}
Suppose $a_1 = a^1$ and $b_1 = b^1$, so $v(x_1) = \{a_1,b_1\} = v(x^1)$.  Suppose $X \setminus x$ is the box obtained from $X$ by removing $x = \{x_1,x^1\}$, so $\vxminusxin = v(x_2)$ and $\vxminusxout = v(x^2)$.  Consider the $1$-loop (Def. \ref{def:loop-wd})
\begin{center}
\begin{tikzpicture}[scale=.9]
\draw [ultra thick] (2,1) rectangle (3,2);
\node at (2.5,1.4) {$X$};
\node at (2.2,1.75) {\tiny{$x_1$}};
\node at (2.8,1.8) {\tiny{$x^1$}};
\draw [arrow] (.5,1.3) to (2,1.3);
\node at (-.1,1.8) {$\xminusxin$};
\draw [arrow] (3,1.3) to (4.5,1.3);
\node at (5.2,1.8) {$\xminusxout$};
\draw [ultra thick] (1,0.7) rectangle (4,2.5);
\node at (2.5,2.8) {$\lambda_{X, x}$};
\draw [arrow, looseness=3.5] (3,1.8) to [out=30, in=150] (2,1.8);
\end{tikzpicture}
\end{center}
in $\wddot\xminusxx$.  Then
\[\lambda_{X,x}(T,\frd,\fup) = (T,\frd_{\setminus x},\fup_{\setminus x}) \in \DS(X \setminus x)\]
is the $(X\setminus x)$-discrete system with update and readout functions
\[\nicexy{\vxminusxin \times T = \{a_2,b_2\} \times T \ar[r]^-{\fup_{\setminus x}} & T \ar[r]^-{\frd_{\setminus x}} & \{a^2,b^2\} = \vxminusxout.}\]
Its transition diagram is:
\begin{center}
\begin{tikzpicture}[auto, node distance=2cm,font=\small,
    every node/.style={inner sep=2pt, rectangle, rounded corners=0.2cm, text centered},
    state/.style={draw, very thick, minimum height=3.5em, text width=2.7cm}]
\matrix[row sep=.8cm,column sep=2cm] {
\node [state] (1) {state $= 1$\\ readout = $a^2$}; &
\node [state] (2) {state $= 2$\\ readout = $a^2$};\\};
\draw [arrow, out=70, in=110, looseness=3] (1) to node[swap]{\scriptsize{$a_2$}} (1);
\draw [arrow] (1) to node{\scriptsize{$b_2$}} (2);
\draw [arrow, out=70, in=110, looseness=3] (2) to node[swap]{\scriptsize{$a_2$}} (2);
\draw [arrow, out=10, in=-10, looseness=3] (2) to node{\scriptsize{$b_2$}} (2);
\end{tikzpicture}
\end{center}
For instance, we have
\[\fup_{\setminus x}(a_2,1) = \fup\bigl((a_2,\frd(1)_{x^1}),1\bigr) = \fup\bigl((a_1,a_2),1\bigr) = 1,\]
which explains the loop at state $1$ labeled $a_2$.  Similar calculation yields the rest of the update function and the readout function.
\end{example}

\begin{example}
Suppose $a_1=a_2$ and $b_1=b_2$, so $v(x_1) = \{a_1,b_1\}=v(x_2)$.  Suppose $Y$ is the box $X/(x_1=x_2)$, so $\vyin = v(x_1)$ and $\vyout = \vxout$.  Consider the in-split (Def. \ref{def:insplit-wd})
\begin{center}
\begin{tikzpicture}[scale=.9]
\draw [ultra thick] (2,1) rectangle (3,2.1);
\node at (2.7,1.5) {$X$};
\node at (2.2,1.8) {\tiny{$x_1$}};
\node at (2.2,1.2) {\tiny{$x_2$}};
\node at (0,1.6) {$\yin$};
\node at (.7,1.7) {\tiny{$y$}};
\draw [thick] (.5,1.5) -- (1.2,1.5);
\draw [arrow] (1.2,1.5) to [out=0, in=180] (2,1.8);
\draw [arrow] (1.2,1.5) to [out=0, in=180] (2,1.2);
\draw [arrow] (3,1.2) to (4.5,1.2);
\draw [arrow] (3,1.8) to (4.5,1.8);
\node at (5,1.6) {$\yout$};
\draw [ultra thick] (1,0.7) rectangle (4,2.4);
\node at (2.5,2.7) {$\sigma_{X, x_1, x_2}$};
\end{tikzpicture}
\end{center}
in $\wddot\yx$.  Then
\[\sigma_{X, x_1, x_2}(T,\frd,\fup) = (T,\frd,\sigmasubdot\fup) \in \DS(Y)\]
is the $Y$-discrete system with update and readout functions
\[\nicexy{\vyin \times T = \{a_1,b_1\} \times T \ar[r]^-{\sigmasubdot\fup} & T \ar[r]^-{\frd} & \vyout = \{a^1,b^1\} \times \{a^2,b^2\}.}\]
Its transition diagram is:
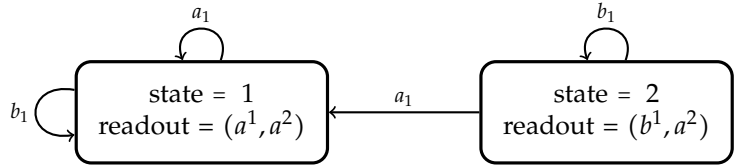
\begin{center}
\begin{tikzpicture}[auto, node distance=2cm,font=\small,
    every node/.style={inner sep=2pt, rectangle, rounded corners=0.2cm, text centered},
    state/.style={draw, very thick, minimum height=3.5em, text width=3.2cm}]
\matrix[row sep=.8cm,column sep=2cm] {
\node [state] (1) {state $= 1$\\ readout = $(a^1,a^2)$}; &
\node [state] (2) {state $= 2$\\ readout = $(b^1,a^2)$};\\};
\draw [arrow, out=70, in=110, looseness=3] (1) to node[swap]{\scriptsize{$a_1$}} (1);
\draw [arrow, out=170, in=190, looseness=3] (1) to node[swap]{\scriptsize{$b_1$}} (1);
\draw [arrow, out=70, in=110, looseness=3] (2) to node[swap]{\scriptsize{$b_1$}} (2);
\draw [arrow] (2) to node[swap]{\scriptsize{$a_1$}} (1);
\end{tikzpicture}
\end{center}
For instance, we have
\[(\sigmasubdot\fup)(a_1,2) = \fup\bigl((a_1,a_2),2\bigr) = 1,\]
which explains the arrow from state $2$ to state $1$ labeled $a_1$.  Similar calculation yields the rest of the update function.
\end{example}

\begin{example}
Suppose $Z$ is a box with $z \not= z'$ in $\zout$ such that $v(z) = v(z') = v(x^1)$ and that $Z/(z=z') = X$.  So $\vzin = \vxin$ and $\vzout = v(x^1) \times v(x^1) \times v(x^2)$.  Consider the out-split (Def. \ref{def:out-split})
\begin{center}
\begin{tikzpicture}[scale=.9]
\draw [ultra thick] (2,1) rectangle (3,2);
\node at (2.5,1.3) {$X$};
\node at (2.75,1.8) {\tiny{$x^1$}};
\draw [thick] (3,1.8) -- (3.3,1.8);
\draw [arrow] (3.3,1.8) to [out=0, in=180] (4.5,2);
\draw [arrow] (3.3,1.8) to [out=0, in=180] (4.5,1.6);
\draw [arrow] (.5,1.2) to (2,1.2);
\draw [arrow] (.5,1.8) to (2,1.8);
\node at (0,1.5) {$\zin$};
\draw [arrow] (3,1.2) to (4.5,1.2);
\node at (5.5,1.5) {$\zout$};
\node at (4.7,2) {\tiny{$z$}};
\node at (4.7,1.6) {\tiny{$z'$}};
\draw [ultra thick] (1,0.7) rectangle (4,2.3);
\node at (2.5,2.6) {$\sigma^{Z,z,z'}$};
\end{tikzpicture}
\end{center}
in $\wddot\zx$.  Then
\[\sigma^{Z,z,z'}(T,\frd,\fup) = (T,\sigmasupdot\frd,\fup) \in \DS(Z)\]
is the $Z$-discrete system with update and readout functions
\[\nicexy{\vzin \times T = \vxin \times T \ar[r]^-{\fup} & T \ar[r]^-{\sigmasupdot\frd} & \vzout.}\]
Its transition diagram is:
\begin{center}
\begin{tikzpicture}[auto, node distance=2cm,font=\small,
    every node/.style={inner sep=2pt, rectangle, rounded corners=0.2cm, text centered},
    state/.style={draw, very thick, minimum height=3.5em, text width=3.6cm}]
\matrix[row sep=.8cm,column sep=2cm] {
\node [state] (1) {state $= 1$\\ readout = $(a^1,a^1,a^2)$}; &
\node [state] (2) {state $= 2$\\ readout = $(b^1,b^1,a^2)$};\\};
\draw [arrow, out=70, in=110, looseness=3] (1) to node[swap]{\scriptsize{$(a_1,a_2)$}} (1);
\draw [arrow, out=170, in=190, looseness=3] (1) to node[swap]{\scriptsize{$(b_1,b_2)$}} (1);
\draw [arrow, bend left=20] (1) to node{\scriptsize{$(a_1,b_2)$}} (2);
\draw [arrow, bend left=5] (1) to node{\scriptsize{$(b_1,a_2)$}} (2);
\draw [arrow, out=70, in=110, looseness=3] (2) to node[swap]{\scriptsize{$(b_1,a_2)$}} (2);
\draw [arrow, out=10, in=-10, looseness=3] (2) to node{\scriptsize{$(b_1,b_2)$}} (2);
\draw [arrow, bend left=5] (2) to node{\scriptsize{$(a_1,a_2)$}} (1);
\draw [arrow, bend left=20] (2) to node{\scriptsize{$(a_1,b_2)$}} (1);
\end{tikzpicture}
\end{center}
This transition diagram is the same as that of $(T,\frd,\fup)$, except for the values of the readout function at states $1$ and $2$.
\end{example}

\begin{example}
Suppose $W$ is a box such that $W \setminus w = X$ for some $w \in \win$.  So $\vwin = v(w) \times \vxin$ and $\vwout = \vxout$.  Consider the $1$-wasted wire (Def. \ref{def:wasted-wire-wd})
\begin{center}
\begin{tikzpicture}[scale=.9]
\draw [ultra thick] (2,1) rectangle (3,2);
\node at (2.5,1.5) {$X$};
\draw [arrow] (.5,1.2) to (2,1.2);
\draw [arrow] (.5,1.6) to (2,1.6);
\node at (0,1.5) {$\win$};
\draw [arrow] (.5,2) -- (1,2);
\node at (.7,2.2) {\tiny{$w$}};
\draw [arrow] (3,1.2) to (4.5,1.2);
\draw [arrow] (3,1.8) to (4.5,1.8);
\node at (5.,1.5) {$\wout$};
\draw [ultra thick] (1,0.7) rectangle (4,2.3);
\node at (2.5,2.6) {$\omega_{W,w}$};
\end{tikzpicture}
\end{center}
in $\wddot\wx$.  Then
\[\omega_{W,w}(T,\frd,\fup) = (T,\frd,\omega\fup) \in \DS(W)\]
is the $W$-discrete system with update and readout functions
\[\nicexy{\vwin \times T = v(w) \times \vxin \times T \ar[r]^-{\omega\fup} & T \ar[r]^-{\frd} & \vwout.}\]
Its transition diagram is:
\begin{center}
\begin{tikzpicture}[auto, node distance=2cm,font=\small,
    every node/.style={inner sep=3pt, rectangle, rounded corners=0.2cm, text centered},
    state/.style={draw, very thick, minimum height=3.5em, text width=3.2cm}]
\matrix[row sep=.8cm,column sep=2cm] {
\node [state] (1) {state $= 1$\\ readout = $(a^1,a^2)$}; &
\node [state] (2) {state $= 2$\\ readout = $(b^1,a^2)$};\\};
\draw [implies, out=70, in=110, looseness=3] (1) to node[swap]{\scriptsize{$(s,a_1,a_2)$}} (1);
\draw [implies, out=170, in=190, looseness=3] (1) to node[swap]{\scriptsize{$(s,b_1,b_2)$}} (1);
\draw [implies, bend left=30] (1) to node{\scriptsize{$(s,a_1,b_2)$}} (2);
\draw [implies, bend left=5] (1) to node{\scriptsize{$(s,b_1,a_2)$}} (2);
\draw [implies, out=70, in=110, looseness=3] (2) to node[swap]{\scriptsize{$(s,b_1,a_2)$}} (2);
\draw [implies, out=10, in=-10, looseness=3] (2) to node{\scriptsize{$(s,b_1,b_2)$}} (2);
\draw [implies, bend left=5] (2) to node{\scriptsize{$(s,a_1,a_2)$}} (1);
\draw [implies, bend left=30] (2) to node{\scriptsize{$(s,a_1,b_2)$}} (1);
\end{tikzpicture}
\end{center}
Here each double arrow $\Rightarrow$ represents the set of arrows as $s$ runs through $v(w)$.  For instance, for each $s \in v(w)$ we have
\[(\omega\fup)\bigl((s,a_1,b_2),1\bigr) = \fup\bigl((a_1,b_2),1\bigr) = 2,\]
which explains the double arrow from state $1$ to state $2$ labeled $(s,a_1,b_2)$.  Similar calculation yields the rest of the update function.
\end{example}

The following observation ensures that $\DS$ is a well-defined $\wddot$-algebra, i.e., that it satisfies the generating axioms.

\begin{theorem}
\label{thm:algebra-ds}
The algebra of discrete systems $\DS$ in Def. \ref{def:algebra-discrete-systems} is actually a $\wddot$-algebra in the sense of Def. \ref{def:normal-algebra}, hence also in the sense of Def. \ref{def2:operad-algebra} by Theorem \ref{thm:normal-wd-algebra}.\index{finite presentation for the algebra of discrete systems}  
\end{theorem}

\begin{proof}
We must check that $\DS$ satisfies the $28$ generating axioms in Def. \ref{def:wd-algebra}, which are all trivial to check.  For example, the generating axiom \eqref{wd-algebra-loopelement} says that, in the setting of \eqref{move:c5} with $X^* = X \setminus \{x^{12}, x_1, x_2\}$ and $v(x^{12}) = v(x_1) = v(x_2)$, the diagram
\[\nicexy{
\DS(X) \ar[r]^-{\sigma^{Y,x^1,x^2}} \ar[d]_{\sigma_{X,x_1,x_2}} 
& \DS(Y) \ar[r]^-{\lambda_{Y,x(1)}}
& \DS(Y \setminus x(1)) \ar[d]^-{\lambda_{Y \setminus x(1), x(2)}} \\
\DS(X') \ar[rr]^-{\lambda_{X',x}} && \DS(X^*)}\]
is commutative.  When applied to a typical element $(T,\frd,\fup) \in \DS(X)$, a simple direct inspection reveals that both compositions in the above diagram yield $(T,\grd,\gup) \in \DS(X^*)$, in which
\[\begin{split}
\grd(t) &= \frd(t)_{\setminus x^{12}};\\
\gup\bigl(\uy, t\bigr) &= \fup\Bigl(\bigl(\uy, \frd(t)_{x^{12}}, \frd(t)_{x^{12}}\bigr), t\Bigr)
\end{split}\]
Here for $(\uy,t) \in \vxstarin \times T$, we have
\[\begin{split}
&\frd(t) = \Bigl(\frd(t)_{\setminus x^{12}}, \frd(t)_{x^{12}}\Bigr) \in \vxout = \vxstarout \times v(x^{12});\\
& \bigl(\uy, \frd(t)_{x^{12}}, \frd(t)_{x^{12}}\bigr) \in \vxstarin \times v(x_1) \times v(x_2) = \vxin.\end{split}\]
The other generating axioms are checked similarly.
\end{proof}

\begin{remark}
\label{rk:ds-algebra-agree}
Our definition of the algebra of discrete systems $\DS$ actually agrees with the one in \cite{spivak15b} (Example 2.7 and Def. 4.9).  To see this, note that Spivak's definition is essentially based on Def. \ref{def1:operad-algebra}, except that it is stated in terms of symmetric monoidal categories.  Spivak's structure map of $\DS$, when applied to the $7$ normal generating wiring diagrams (Def. \ref{def:generating-without-dn}(1)), agrees with ours in Def. \ref{def:algebra-discrete-systems}.  So Theorems  \ref{thm:normal-wd-algebra} and \ref{thm:algebra-ds} imply that the two definitions of $\DS$--namely, the one in \cite{spivak15b} and our Def. \ref{def:algebra-discrete-systems}--are equivalent.
\end{remark}

\section{Algebras over the Operad of Strict Wiring Diagrams}
\label{sec:algebra-strict-wd}

The purpose of this section is to provide a finite presentation for algebras over the $\boxs$-colored operad $\wdzero$ of strict wiring diagrams (Prop. \ref{prop:strict-operads}).  We begin by defining these algebras in terms of finitely many generators and relations.  Recall from Def. \ref{def:generating-strict}) that:
\begin{enumerate}
\item
The \emph{strict generating wiring diagrams} are the empty wiring diagram $\epsilon$, a name change $\tau_{X,Y}$, a $2$-cell $\theta_{X,Y}$, and a $1$-loop $\lambda_{X,x}$.
\item
The \emph{strict elementary relations} are the $8$ elementary relations that involve only strict generating wiring diagrams on both sides.
\end{enumerate}

\begin{definition}
\label{def:strict-algebra}
A $\wdzero$-algebra $A$ consists of the following data.\index{wdzero-algebra@$\wdzero$-algebra}
\begin{enumerate}
\item
For each $X \in \boxs$, $A$ is equipped with a class $A_X$ called the \emph{$X$-colored entry} of $A$.
\item
It is equipped with the $4$ \emph{generating structure maps} in Def. \ref{def:wd-algebra} corresponding to the strict generating wiring diagrams.
\end{enumerate}
This data is required to satisfy the $8$ generating axioms in Def. \ref{def:wd-algebra} corresponding to the strict elementary relations, namely,  \eqref{axiom-a1}, \eqref{axiom-a2}, \eqref{axiom-a3}, \eqref{axiom-b0}, \eqref{wd-algebra-2cell-associativity}, \eqref{wd-algebra-equivariance}, \eqref{axiom-b3}, and \eqref{wd-algebra-doubleloop}.
\end{definition}

The next observation is the $\wdzero$ version of Theorems \ref{thm:wd-algebra} and \ref{thm:normal-wd-algebra}.  It gives a finite presentation for $\wdzero$-algebras.

\begin{theorem}
\label{thm:strict-wd-algebra}
For the operad $\wdzero$ of strict wiring diagrams (Prop. \ref{prop:strict-operads}), Def. \ref{def2:operad-algebra} with $\O=\wdzero$ and Def. \ref{def:strict-algebra} of a $\wdzero$-algebra are equivalent.\index{finite presentation for WDzero algebra@finite presentation for $\wdzero$-algebras}
\end{theorem}

\begin{proof}
Simply restrict the proof of Theorem \ref{thm:wd-algebra} to strict (generating) wiring diagrams.  Instead of Theorem \ref{stratified-presentation-exists}, here we use Theorem \ref{thm:strict-wd-coherence} for the existence of a presentation involving only strict generating wiring diagrams.
\end{proof}

\section{Finite Presentation for the Algebra of Open Dynamical Systems}
\label{sec:algebra-ods}

The purpose of this section is to provide a finite presentation for the algebra of open dynamical systems introduced in \cite{vsl}.  In \cite{vsl} the algebra of open dynamical systems $\G$ was defined and verified using essentially Def. \ref{def1:operad-algebra} but in the form of symmetric monoidal categories and monoidal functors.  Our definition of $\G$ in Def. \ref{def:ods-algebra} is based on Def. \ref{def:strict-algebra}, which involves four relatively simple generating structure maps.  Our verification that $\G$ is actually a $\wdzero$-algebra in Theorem \ref{prop:ods-algebra-is-algebra} boils down to verifying the generating axiom \eqref{wd-algebra-doubleloop} for a double-loop.  This is a simple exercise involving the definition of the generating structure map corresponding to a $1$-loop  \eqref{ods-loop}.   The equivalence between the two definitions of the algebra of open dynamical systems is explained in Remark \ref{rk:ods-algebra-agree}.

Let us first recall the setting of \cite{vsl}.  For the definitions of the basic objects in differential geometry that appear below, the reader may consult, for example, \cite{helgason} (Ch.1).

\begin{assumption}
Throughout this section:
\begin{enumerate}
\item
$S$ is a chosen set of representatives of isomorphism classes of second-countable smooth manifolds, henceforth called \emph{manifolds}.\index{manifolds}
\item
The operad of strict wiring diagrams $\wdzero$ (Prop. \ref{prop:strict-operads}) is defined using this choice of $S$. 
\item
 All the maps between manifolds are smooth maps.  
\item
For a manifold $M$, denote by $\pi : TM \to M$ the projection map of the tangent bundle.\index{tangent bundle}
\end{enumerate}
\end{assumption}

\begin{definition}
Suppose $(M,v) \in \Fins$ is an $S$-finite set (Def. \ref{def:Fins}).
\begin{enumerate}
\item
Define
\begin{equation}
\label{m-sub-v}
M_v = \prod_{m \in M} v(m) \in S.\end{equation}
\item
For a subset $I \subseteq M$ and $\ux = (x_m)_{m \in M} \in M_v$, define
\begin{equation}
\label{xminusn}
\begin{split}
\ux_I &= (x_m)_{m \in I} \in I_v = \prod_{m \in I} v(m)\\
\ux_{\setminus I} &= (x_m)_{m \in M \setminus I} \in (M \setminus I)_v = \prod_{m \in M \setminus I} v(m).
\end{split}
\end{equation}
\end{enumerate}
\end{definition}

\begin{definition}
\label{def:vsl-setting}
An \emph{open dynamical system}, or \emph{ods} for short, is a tuple $\left(M,\uin,\uout,f\right)$ consisting of:\index{open dynamical system}
\begin{enumerate}
\item
manifolds $M$, $\uin$, and $\uout$;
\item
a pair of maps $f = (\fin, \fout)$,
\[\nicexy{
M \times \uin \ar[r]^-{\fin} \ar[dr]_-{\text{project}} & TM \ar[d]^-{\pi} & M \ar[d]^-{\fout}\\
& M & \uout}\]
such that the left diagram is commutative.
\end{enumerate}
\end{definition}

Next is \cite{vsl} (Def. 4.2).

\begin{definition}
\label{def:ods-objects}
For each $X = (\xin,\xout) \in \boxs$, define the class
\begin{equation}
\label{g-sub-x}
\G_X = \Bigl\{(M,f) : M \in \Fins,\, (M_v, \vxin, \vxout, f) \text{ is an ods} \Bigr\}
\end{equation}
in which $M_v$, $\vxin$, and $\vxout$ are as in \eqref{m-sub-v}.
\end{definition}

\begin{example}
\label{ex:ods-empty}
For the empty box $\varnothing = (\varnothing,\varnothing) \in \boxs$ and an $S$-finite set $M$, to say that $(M_v, \varnothing_v = \{*\}, \varnothing_v = \{*\}, f)$ is an open dynamical system means that $f$ is a pair of maps $f = (\fin, \fout)$,
\[
\nicexy{
M_v = M_v \times \{*\} \ar[r]^-{\fin} \ar[dr]_-{\Id} & TM_v \ar[d]^-{\pi} & M_v \ar[d]^-{\fout}\\ 
& M_v & \{*\}}\]
such that the left diagram is commutative.  Since $\fout$ gives no information, $f = \fin$ is equivalent to a vector field \index{vector field} on $M_v$.   So
\begin{equation}
\label{g-sub-empty}
\G_\varnothing = \Bigl\{(M,f) : M \in \Fins, \text{ $f$ is a vector field on $M_v$} \Bigr\}.
\end{equation}
\end{example}

\begin{example}\label{ex:odsx}
Suppose $W$ is the box with $\win = \{w_1,w_2\}$, $\wout = \{w^1,w^2\}$, and all $v(-) = \bR$, so $\vwin = \vwout = \bR^2$.  Suppose $M$ is the one-point set with value $\bR$.  There is an element $(M,f) \in \G_W$ whose structure maps
\[\nicexy{\bR \times \bR^2 \ar[r]^-{\fin} \ar[dr]_-{\pi_1} & TM = \bR^2 \ar[d]^-{\pi_1} & \bR \ar[d]^-{\fout}\\
& \bR & \bR^2}\]
are given by
\[\fin(x,y,z) = (x, ax+by+cz) \andspace \fout(x) = (dx,e^x)\]
for any choice of fixed parameters $a,b,c,d \in \bR$.
\end{example}

We now define the algebra of open dynamical systems in terms of $4$ generating structure maps.

\begin{definition}
\label{def:ods-algebra}
The \emph{algebra of open dynamical systems} \index{algebra of open dynamical systems} is the $\wdzero$-algebra $\G$ in the sense of Def. \ref{def:strict-algebra} defined as follows.  For each box $X \in \boxs$, the $X$-colored entry is $\G_X$ in \eqref{g-sub-x}.

The $4$ generating structure maps (Def. \ref{def:wd-algebra}) are defined as follows.
\begin{enumerate}
\item
Corresponding to the empty wiring diagram $\epsilon \in \wdzero\emptynothing$ (Def. \ref{def:empty-wd}), the structure map
\begin{equation}
\label{ods-empty}
\nicexy{\ast \ar[r]^-{\epsilon} & \G_\varnothing}
\end{equation}
sends $*$ to $(\varnothing, *) \in \G_\varnothing$ \eqref{g-sub-empty}.  Here $\varnothing \in \Fins$ is the empty set, in which $\varnothing_v = \{*\}$, and in the second entry $*$ is the trivial vector field.
\item
Corresponding to a name change $\tau_{X,Y} \in \wdzero\yx$  (Def. \ref{def:name-change}), the structure map
\begin{equation}
\label{ods-name-change}
\nicexy{\G_X \ar[r]^-{\tau_{X,Y}}_-{=} & \G_Y}
\end{equation}
is the identity map, using the fact that $\vxin = \vyin$ and $\vxout = \vyout$.
\item
Corresponding to a $2$-cell $\theta_{X,Y} \in \wdzero\xplusyxy$ (Def. \ref{def:theta-wd}), it has the structure map
\begin{equation}
\label{ods-2cell}
\nicexy{
\G_X \times \G_Y \ar[d]_-{\theta_{X,Y}} & \Bigl((M_X,f_X), (M_Y,f_Y)\Bigr) \ar@{|->}[d]\\
\G_{\xplusy} & \Bigl(M_X \amalg M_Y, f_X \times f_Y\Bigr)}
\end{equation}
in which $M_X \amalg M_Y$ is the coproduct in $\Fins$ (Def. \ref{def:Fins}).
\item
Corresponding to a $1$-loop $\lambda_{X,x} \in \wdzero\xminusxx$  (Def. \ref{def:loop-wd}) with $x = \bigl(x_-, x_+\bigr) \in \xin \times \xout$, it has the structure map
\begin{equation}
\label{ods-loop}
\nicexy{\G_X \ar[r]^-{\lambda_{X,x}} & \G_{\xminusx}}\quad  
\nicexy{(M,f) \ar@{|->}[r] & \bigl(M,\fminusx\bigr).}
\end{equation}
The maps $\fminusx = \left(\finminusx, \foutminusx\right)$ are defined as
\[\begin{split}
\finminusx \bigl(\um, \uy\bigr) &= \fin\Bigl(\um, \bigl(\fout(\um)_{\xplus}, \uy\bigr)\Bigr) \in TM_v\\
\foutminusx(\um) &= \fout(\um)_{\minusxplus} \in \vxminusxout
\end{split}\]
for $\um \in M_v$ and $\uy \in \vxminusxin$.  Recalling that $\fout(\um) \in \vxout$, the elements
\[\fout(\um)_{\xplus} \in v(x_+) = v(x_-) \andspace \fout(\um)_{\minusxplus} \in \vxminusxout\]
are as in \eqref{xminusn}.
\end{enumerate}
This finishes the definition of the $\wdzero$-algebra of open dynamical systems.
\end{definition}

\begin{example}
This is a continuation of Example \ref{ex:odsx}, where $W$ is the box with $\win = \{w_1,w_2\}$ and $\wout = \{w^1,w^2\}$ and $M$ is the one-point set, all with $v(-) = \bR$.  Suppose $W\setminus w$ is the box obtained from $W$ by removing $w = \{w_1,w^1\}$, so 
\[\vwminuswin = \vwminuswout = \bR.\]  
Consider the $1$-loop (Def. \ref{def:loop-wd})
\begin{center}
\begin{tikzpicture}[scale=.9]
\draw [ultra thick] (2,1) rectangle (3,2);
\node at (2.5,1.4) {$W$};
\node at (2.2,1.75) {\tiny{$w_1$}};
\node at (2.8,1.8) {\tiny{$w^1$}};
\draw [arrow] (.5,1.3) to (2,1.3);
\node at (-.1,1.8) {$\wminuswin$};
\draw [arrow] (3,1.3) to (4.5,1.3);
\node at (5.2,1.8) {$\wminuswout$};
\draw [ultra thick] (1,0.7) rectangle (4,2.5);
\node at (2.5,2.8) {$\lambda_{W,w}$};
\draw [arrow, looseness=3.5] (3,1.8) to [out=30, in=150] (2,1.8);
\end{tikzpicture}
\end{center}
in $\wdzero\wminusww$.  Then $\lambda_{W,w}(M,f) = (M,f_{\setminus w}) \in \G_{W\setminus w}$ has structure maps
\[\nicexy{\bR \times \bR \ar[r]^-{\fin_{\setminus w}} \ar[dr]_-{\pi_1} & TM = \bR^2 \ar[d]^-{\pi_1} & \bR \ar[d]^-{\fout_{\setminus w}}\\
& \bR & \bR}\]
given by
\[\begin{split}
\fout_{\setminus w}(x) &= \fout(x)_{\setminus w^1} = (dx,e^x)_{\setminus w^1} = e^x;\\
\fin_{\setminus w}(x,y) &= \fin\bigl(x,(\fout(x)_{w^1},y)\bigr) = \fin(x,dx,y) = \bigl(x, (a+bd)x + cy\bigr).
\end{split}\]
\end{example}

The following observation ensures that $\G$ is a well-defined $\wdzero$-algebra, i.e., that it satisfies the generating axioms.

\begin{theorem}
\label{prop:ods-algebra-is-algebra}
The algebra of open dynamical systems $\G$ in Def. \ref{def:ods-algebra} is actually  a $\wdzero$-algebra in the sense of  Def. \ref{def:strict-algebra},\index{finite presentation for the algebra of open dynamical systems} hence also in the sense of Def. \ref{def2:operad-algebra} by Theorem \ref{thm:strict-wd-algebra}.
\end{theorem}

\begin{proof}
We must check the $8$ generating axioms corresponding to the strict elementary relations listed in Def. \ref{def:strict-algebra}.  All of them follow from a quick inspection of the definitions except for \eqref{wd-algebra-doubleloop}.  This generating axiom says that, in the setting of the elementary relation \eqref{move:c1} corresponding to a double-loop, the diagram
\begin{equation}
\label{ods-doubleloop}
\nicexy@C+.4cm{
\G_X \ar[r]^-{\lambda_{X,x^2}} \ar[d]_{\lambda_{X,x^1}} 
& \G_{X \setminus x^2} \ar[d]^-{\lambda_{X\setminus x^2, x^1}} \\
\G_{X \setminus x^1} \ar[r]^-{\lambda_{X \setminus x^1, x^2}} & \G_{X \setminus x}}
\end{equation}
is commutative.  

To prove \eqref{ods-doubleloop}, suppose $(M,f) \in \G_X$.  First define the element $(M,\fminusx) \in \G_{X \setminus x}$ with the maps $\fminusx = (\finminusx, \foutminusx)$,
\[\nicexy{
M_v \times \vxminusxin \ar[r]^-{\finminusx} \ar[dr]_-{\text{project}} & TM_v \ar[d]^-{\pi} & M_v \ar[d]^-{\foutminusx}\\
& M_v & \vxminusxout}\]
defined as follows.  Given $\um \in M_v$ and
\[
\uy \in \vxminusxin = \prod_{y \in \xin \setminus\{\xoneminus, \xtwominus\}} v(y)\]
we define
\[\begin{split}
\finminusx\bigl(\um,\uy\bigr) &= \fin\Bigl(\um; \bigl(\fout(\um)_{\{\xoneplus, \xtwoplus\}}, \uy\bigr) \Bigr) \in TM_v\\
\foutminusx(\um) &= \fout(\um)_{\setminus\{\xoneplus, \xtwoplus\}} \in \vxminusxout.
\end{split}\]
Recalling that $\fout(\um) \in \vxout$, the elements
\[
\fout(\um)_{\{\xoneplus, \xtwoplus\}} \in v(\xoneplus) \times v(\xtwoplus) = v(\xoneminus) \times v(\xtwominus)\]
and $ \fout(\um)_{\setminus\{\xoneplus, \xtwoplus\}}$ are as in  \eqref{xminusn}.

Now it follows from a direct inspection using the definition \eqref{ods-loop} that both composites in \eqref{ods-doubleloop}, when applied to $(M,f)$, yields $(M,\fminusx) \in \G_{\xminusx}$.  This proves the generating axiom \eqref{ods-doubleloop} for $\G$.
\end{proof}

\begin{remark}
\label{rk:ods-algebra-agree}
Our definition of the algebra of open dynamical systems $\G$ actually agrees with the one in \cite{vsl} (Def. 4.4 and 4.5).  To see this, note that among the four generating structure maps in Def. \ref{def:ods-algebra}:
\begin{itemize}
\item
$\epsilon$ \eqref{ods-empty}, $\tau_{X,Y}$ \eqref{ods-name-change}, and $\lambda_{X,x}$ \eqref{ods-loop} agree with \cite{vsl} (Def. 4.4);
\item
$\theta_{X,Y}$ \eqref{ods-2cell} agrees with \cite{vsl} (Def. 4.5).
\end{itemize}
Theorem \ref{thm:strict-wd-algebra} then implies that the two definitions of $\G$--namely, the one in \cite{vsl} and our Def. \ref{def:ods-algebra}--are equivalent.
\end{remark}

\section{Summary of Chapter \ref{ch06-wd-algebras}}

\begin{enumerate}
\item For an $S$-colored operad $\sO$, an $\sO$-algebra $A$ consists of a class $A_c$ for each $c \in S$ and a structure map
\[\nicexy{A_{c_1} \times \cdots \times A_{c_n} \ar[r]^-{\zeta} & A_d}\]
for each $\zeta \in \sO\dconecn$ that satisfies the associativity, unity, and equivariance axioms.
\item Each $\WD$-algebra can be described using eight generating structure maps that satisfy twenty-eight generating axioms.
\item The propagator algebra is a $\WD$-algebra.
\item Each $\wddot$-algebra can be described using seven generating structure maps that satisfy twenty-eight generating axioms.
\item The algebra of discrete systems is a $\wddot$-algebra.
\item Each $\wdzero$-algebra can be described using four generating structure maps that satisfy eight generating axioms.
\item The algebra of open dynamical systems is a $\wdzero$-algebra.
\end{enumerate}

\part{Undirected Wiring Diagrams}

The main purpose of this part is to describe the combinatorial structure of the operad $\uwd$ of undirected wiring diagrams.  The main result is a finite presentation theorem that describes the operad $\uwd$ in terms of $6$ operadic generators and $17$ generating relations.

The operad $\uwd$ of undirected wiring diagrams is recalled in Chapter \ref{ch07-undirected-wiring-diagrams}.  Operadic generators and generating relations for the operad $\uwd$ are presented in Chapter \ref{ch08-generating-uwd}.  Various decompositions of undirected wiring diagrams are given in Chapter \ref{ch09-decomp-uwd}.  Using the results in Chapters \ref{ch08-generating-uwd} and \ref{ch09-decomp-uwd}, the finite presentation theorem for the operad $\uwd$ is proved in Chapter \ref{ch10-stratified-uwd}.  In Chapter \ref{ch11-uwd-algebras} we prove the corresponding finite presentation theorem for $\uwd$-algebras and discuss the main example of the relational algebra.  This finite presentation theorem for algebras describes each $\uwd$-algebra in terms of finitely many generating structure maps and relations among these maps.  Also given in this Chapter is a partial verification of a conjecture of Spivak about the rigidity of the relational algebra.

\textbf{Reading Guide}.  The reader who already knows about pushouts of finite sets may skip Section \ref{sec:pushout}.  In Section \ref{sec:uwd-operad-structure}, where we define the operad structure on undirected wiring diagrams, the reader may wish to skip the proofs of Lemmas \ref{uwd-compi-horizontal-associative} and \ref{uwd-compi-vertical-associative} and just study the accompanying pictures.

Section \ref{sec:wasted-cable} is not technically needed in later sections.  However, it contains several illuminating examples about how wasted cables can be created from undirected wiring diagrams without wasted cables.  So we urge the reader \emph{not} to skip these examples.

The decompositions in Chapter \ref{ch09-decomp-uwd} are illustrated with a detailed example in Section \ref{sec:ex-decomp-uwd}.  The reader may read that section and skip the rest of the Chapter during the first reading.

In Section \ref{sec:uwd-elementary-eq}, after the initial definitions and examples, the reader may skip the proofs of Lemmas \ref{lemma:simplex-to-stratified-uwd} and \ref{lemma:stratified-connected-uwd} and go straight to Theorem \ref{thm:uwd-generator-relation}, the finite presentation theorem for undirected wiring diagrams.

\chapter{Undirected Wiring Diagrams}
\label{ch07-undirected-wiring-diagrams}

The purposes of this chapter are
\begin{enumerate}
\item
to recall the definition of an undirected wiring diagram (Def. \ref{def:uwd}) from \cite{spivak13};
\item
to give a proof that the collection of undirected wiring diagrams forms an operad (Theorem \ref{uwd-operad}).
\end{enumerate}
There is a subtlety regarding the definition and the operadic composition of  undirected wiring diagrams; see Remark \ref{rk:spivak-def-uwd}(4) and Example \ref{ex:uwd-jointly-surjective}.  Many more examples of undirected wiring diagrams and their operadic composition will be given in the next chapter.

Fix a class $S$ for this chapter.

\section{Defining Undirected Wiring Diagrams}

In this section, we recall the definition of an undirected wiring diagram.  Recall from Def. \ref{def:Fins} that an \emph{$S$-finite set} is a pair $(X,v)$ with $X$ a finite set and $v : X \to S$ a function, called the \emph{value assignment}.  Maps between $S$-finite sets are functions compatible with the value assignments.  The category of $S$-finite sets is denoted by $\Fins$.   As in earlier chapters, if there is no danger of confusion, then we write an $S$-finite set $(X,v)$ simply as $X$.  The number of elements in a finite set $T$ is denoted by $|T|$.

As in Section \ref{sec:wiring-diagrams} we first define undirected prewiring diagrams.  Undirected wiring diagrams are then defined as the appropriate equivalence classes.  

\begin{motivation}
Before we define an undirected wiring diagram precisely, let us first provide some motivation for it.  An undirected wiring diagram is a picture similar to this:
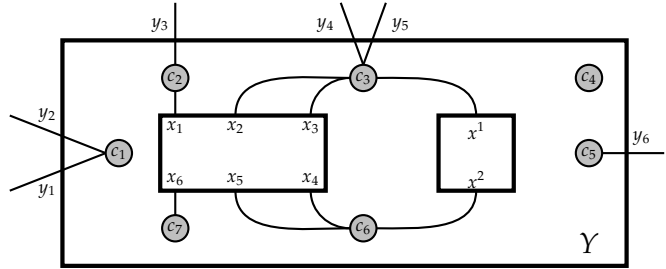
\begin{center}
\begin{tikzpicture}[scale=1]
\draw [ultra thick] (-1,0) rectangle (6.5,3);
\node at (-1.2,1) {\tiny{$y_1$}};
\node at (-1.2,2) {\tiny{$y_2$}};
\node at (.3,3.2) {\tiny{$y_3$}};
\node at (2.5,3.2) {\tiny{$y_4$}};
\node at (3.5,3.2) {\tiny{$y_5$}};
\node at (6.7,1.7) {\tiny{$y_6$}};
\node at (6,.3) {$Y$};
\draw [ultra thick] (.3,1) rectangle (2.5,2);
\node at (.5,1.85) {\tiny{$x_1$}};
\node at (1.3,1.85) {\tiny{$x_2$}};
\node at (2.3,1.85) {\tiny{$x_3$}};
\node at (2.3,1.15) {\tiny{$x_4$}};
\node at (1.3,1.15) {\tiny{$x_5$}};
\node at (.5,1.15) {\tiny{$x_6$}};
\draw [ultra thick] (4,1) rectangle (5,2);
\node at (4.5,1.15) {\tiny{$x^2$}};
\node at (4.5,1.8) {\tiny{$x^1$}};
\cable{(-.25,1.5)} \node at (-.25,1.5) {\tiny{$c_1$}};
\cable{(.5,2.5)} \node at (.5,2.5) {\tiny{$c_2$}};
\cable{(3,2.5)} \node at (3,2.5) {\tiny{$c_3$}};
\cable{(6,2.5)} \node at (6,2.5) {\tiny{$c_4$}};
\cable{(6,1.5)} \node at (6,1.5) {\tiny{$c_5$}};
\cable{(3,.5)} \node at (3,.5) {\tiny{$c_6$}};
\cable{(.5,.5)} \node at (.5,.5) {\tiny{$c_7$}};
\draw [thick] (-.425,1.5) -- (-1.7,1);
\draw [thick] (-.425,1.5) -- (-1.7,2);
\draw [thick] (.5,2) -- (.5,2.325);
\draw [thick] (.5,2.675) -- (.5,3.5);
\draw [thick] (1.3,2) to [out=90, in=180] (2.825,2.5);
\draw [thick] (2.3,2) to [out=90, in=180] (2.825,2.5);
\draw [thick] (4.5,2)  to [out=90, in=0] (3.175,2.5);
\draw [thick] (3,2.675) -- (2.7,3.5);
\draw [thick] (3,2.675) -- (3.3,3.5);
\draw [thick] (6.175,1.5) -- (7,1.5);
\draw [thick] (1.3,1) to [out=-90,in=180] (2.825,.5);
\draw [thick] (2.3,1) to [out=-90,in=180] (2.825,.5);
\draw [thick] (4.5,1) to [out=-90,in=0]  (3.175,.5);
\draw [thick] (.5,1) -- (.5,.675);
\end{tikzpicture}
\end{center}
There are two input boxes $X_1 = \{x_1,\ldots,x_6\}$ and $X_2 = \{x^1,x^2\}$, an output box $Y = \{y_1,\ldots,y_6\}$, and seven cables $\{c_1,\ldots,c_7\}$.  In general, there can be any finite numbers of input boxes and cables, and each input/output box is a finite set.  For each element $z$ in the input boxes and the output box, we need to specify a cable $c$ to which $z$ is connected.   For instance, $x_4$, $x_5$, and $x^2$ are connected to the cable $c_6$.  

Depending on the given undirected wiring diagrams, there are four possible kinds of cables.  First, a cable may only be connected to elements in the inside boxes, such as $c_6$ and $c_7$ in the picture above.  Second, a cable may only be connected to elements in the outside box, such as $c_1$ and $c_5$ above.  Third, a cable may be connected to elements in both the inside boxes and the outside box, such as $c_2$ and $c_3$ above.  Finally, there may be standalone cables that are not connected to anything in the inside and the outside boxes, such as $c_4$ above.  Such distinction among the set of cables will be important later when we discuss the finite presentation for the operad of undirected wiring diagrams.  We will come back to this picture in  Example \ref{ex:uwd-first-example} below.
\end{motivation}

The following definition is a slight generalization of \cite{spivak13} (Examples 2.1.7 and 4.1.1); see Remark \ref{rk:spivak-def-uwd}.  

\begin{definition}
\label{def:pre-uwd}
Suppose $S$ is a class.  An \emph{undirected $S$-prewiring diagram}\index{undirected prewiring diagram} is a tuple
\begin{equation}
\label{undirected-prewiring}
\varphi = \bigl(\uX, Y, C, f, g\bigr)
\end{equation}
consisting of the following data.
\begin{enumerate}
\item
$Y \in \Fins$, called the \emph{output box} of $\varphi$.  An element in $Y$ is called an \emph{output wire} for $\varphi$.
\item
$\uX = (X_1,\ldots,X_n)$ is a $\Fins$-profile for some $n \geq 0$ (Def. \ref{def:profile}); i.e., each $X_i \in \Fins$.  
\begin{itemize}
\item
We call $X_i$ the \emph{$i$th input box}\index{input box} of $\varphi$.  
\item
An element in each $X_i$ is called an \emph{input wire} for $\varphi$. 
\item
Denote by $X = \coprod_{i=1}^n X_i \in \Fins$ the coproduct.
\item
Each element in $X \amalg Y$ is called a \emph{wire}. 
\end{itemize}
\item
$C \in \Fins$, called the \emph{set of cables} of $\varphi$.  Each element in $C$ is called a \emph{cable}.\index{cable}
\item
$f$ and $g$ are maps in the diagram, called a \emph{cospan}\index{cospan}
\begin{equation}
\label{uwd-cospan}
\nicexy{X_1 \amalg \cdots \amalg X_n = X \ar[r]^-{f} & C & Y \ar[l]_-{g}}
\end{equation}
in $\Fins$.
\begin{itemize}
\item
$f$ is called the \emph{input soldering function}\index{input soldering function} and $g$ the \emph{output soldering function}.\index{output soldering function}
\item
If $f(x) = c$, then we say that \emph{$x$ is soldered to $c$} and that \emph{$c$ is soldered to $x$} via $f$.  If $g(y) = c$, then we say that \emph{$y$ is soldered to $c$} and that \emph{$c$ is soldered to $y$} via $g$.\index{soldered to}
\item
If $c \in C$ is a cable and if $m = |f^{-1}(c)|$ and $n = |g^{-1}(c)|$, then $c$ is called an \emph{$(m,n)$-cable}.\index{mncable@$(m,n)$-cable}
\item
A $(0,0)$-cable is also called a \emph{wasted cable}.\index{wasted cable}  In other words, a wasted cable is a cable that is in neither the image of $f$ nor the image of $g$.
\end{itemize}
\end{enumerate}
Given $Y$ and $\uX$, we will denote $\varphi$ as either the tuple $(C,f,g)$ or the cospan \eqref{uwd-cospan}.
\end{definition}

\begin{remark}
In Def. \ref{def:s-box} an $S$-box is an element in $\Fins \times \Fins$.  In the context of undirected (pre)wiring diagrams, the name \emph{box}\index{box} refers to an element in $\Fins$, such as the output box or one of the input boxes.  The context itself should make it clear what box means.
\end{remark}

The cables tell us how to wire the input wires and the output wires together.  So the names of the cables should not matter.  This is made precise in the following definition.

\begin{definition}
\label{def:uwd}
Suppose $\varphi = \left(\uX,Y,C,f,g\right)$ and $\varphi' = \left(\uX,Y,C',f',g'\right)$ are two undirected $S$-prewiring diagrams with the same output box $Y$ and input boxes $\uX$.
\begin{enumerate}
\item
An \emph{equivalence} $h : \varphi \to \varphi'$ is an isomorphism $h : C \to C' \in \Fins$ such that the diagram
\[
\nicexy{X \ar[r]^-{f} \ar[dr]_-{f'} & C \ar[d]^-{h}_-{\cong} & Y \ar[l]_-{g} \ar[dl]^-{g'}\\
& C'}\]
in $\Fins$ is commutative.
\item
We say that $\varphi$ and $\varphi'$ are \emph{equivalent} if there exists an equivalence $\varphi \to \varphi'$.
\item
An \emph{undirected $S$-wiring diagram}\index{undirected wiring diagram} is an equivalence class of undirected $S$-prewiring diagrams.  If $S$ is clear from the context, we will drop $S$ and just say \emph{undirected wiring diagram}.
\item
The class of undirected $S$-wiring diagrams with output box $Y$ and input boxes $\uX = (X_1, \ldots, X_n)$ is denoted by\label{notation:uwd-yux}
\begin{equation}
\label{uwdyux}
\uwd\yux \orspace \uwd\yxonexn.
\end{equation}
The class of all undirected $S$-wiring diagrams is denoted by \index{uwd@$\uwd$}\label{notation:uwd}$\uwd$.  If we want to emphasize the class $S$, we will write $\uwd^S$.
\end{enumerate}
\end{definition}

\begin{remark}
\label{rk:spivak-def-uwd}
Consider Def. \ref{def:pre-uwd} and \ref{def:uwd}.
\begin{enumerate}
\item
If $S = \{*\}$, a one-point set, then what we call an undirected $\{*\}$-wiring diagram is called a \emph{singly-typed wiring diagram}\index{singly-typed wiring diagram}  in \cite{spivak13} (Example 2.1.7).
\item
If $S = \set$, the collection of sets, then what we call an undirected $\set$-wiring diagram is called a \emph{typed wiring diagram}\index{typed wiring diagram}  in \cite{spivak13} (Example 4.1.1).
\item
Cospans \eqref{uwd-cospan} are also used in other work about processes and networks.  For example, cospans in a category, rather than just $\Fins$, are used in \cite{fong}.  That setting is then used in \cite{baez-fong,baez-fong-pollard} to study passive linear networks and Markov processes.
\item
In the book \cite{spivak14} p.464 (but not in \cite{spivak13} Example 2.1.7), Spivak's definition of an undirected $\{*\}$-wiring diagram is slightly different from ours.  More precisely, Spivak insisted that the maps $(f,g)$ in the cospan \eqref{uwd-cospan} be \emph{jointly surjective}, meaning that there are no wasted cables.  However, undirected wiring diagrams whose structure maps $(f,g)$ are jointly surjectivity are \emph{not} closed under the operad structure in $\uwd$, to be defined in Section \ref{sec:uwd-operad-structure}.  In Example \ref{ex:uwd-jointly-surjective} and Section \ref{sec:wasted-cable} we will illustrate that the operadic composition of undirected wiring diagrams without wasted cables can have wasted cables.  In other words, while individual undirected wiring diagrams may have no wasted cables, the operadic composition can actually create wasted cables.  So there is no such thing as the operad of undirected wiring diagrams without wasted cables. In \cite{fong} (Def. 3.1) Fong also used cospans, but did not insist that they be jointly surjective in any way.
\end{enumerate}
\end{remark}

\begin{convention}
\label{conv:undirected-prewiring}
To simplify the presentation, we usually suppress the difference between an undirected prewiring diagram and an undirected wiring diagram.  Each undirected wiring diagram $\varphi = \left(\uX,Y,C,f,g\right)$ has a unique representative in which:
\begin{itemize}
\item
each cable is an element in $S$;
\item
the value assignment $v : C \to S$ sends each cable to itself.  
\end{itemize}
Unless otherwise specified, we will always use this representative of an undirected wiring diagram.
\end{convention}

\begin{example}
\label{ex:uwd-first-example}
Suppose $S$ is any class.  Consider the undirected wiring diagram $\varphi \in \uwd\yxonextwo$ defined as follows.
\begin{itemize}
\item
The input boxes are $X_1 = \{x_1,x_2,x_3,x_4,x_5,x_6\}$ and $X_2 = \{x^1,x^2\} \in \Fins$.
\item
The output box is $Y = \{y_1,y_2,y_3,y_4,y_5,y_6\} \in \Fins$.
\item
The set of cables is $C = \{c_1,c_2,c_3,c_4,c_5,c_6,c_7\} \in \Fins$.
\end{itemize}
Their value assignments satisfied the following conditions:
\begin{itemize}
\item
$v(c_1) = v(y_1) = v(y_2) \in S$.
\item
$v(c_2) = v(x_1) = v(y_3) \in S$.
\item
$v(c_3) = v(x_2) = v(x_3) = v(x^1) = v(y_4) = v(y_5) \in S$.
\item
$v(c_4) \in S$ is arbitrary.
\item
$v(c_5) = v(y_6) \in S$.
\item
$v(c_6) = v(x_4) = v(x_5) = v(x^2) \in S$.
\item
$v(c_7) = v(x_6) \in S$.
\end{itemize}
The input and the output soldering functions
\[\nicexy{X_1 \amalg X_2 \ar[r]^-{f} & C & Y \ar[l]_-{g}}\]
are defined as follows:
\begin{itemize}
\item
$c_1 = g(y_1) = g(y_2)$ is a $(0,2)$-cable.
\item
$c_2 = f(x_1) = g(y_3)$ is a $(1,1)$-cable.
\item
$c_3 = f(x_2) = f(x_3) = f(x^1) = g(y_4) = g(y_5)$ is a $(3,2)$-cable.
\item
$c_4$ is a $(0,0)$-cable, i.e., a wasted cable.
\item
$c_5 = g(y_6)$ is a $(0,1)$-cable.
\item
$c_6 = f(x_4) = f(x_5) = f(x^2)$ is a $(3,0)$-cable.
\item
$c_7 = f(x_6)$ is a $(1,0)$-cable.
\end{itemize}
Graphically we represent this undirected wiring diagram $\varphi \in \uwd\yxonextwo$ as follows.
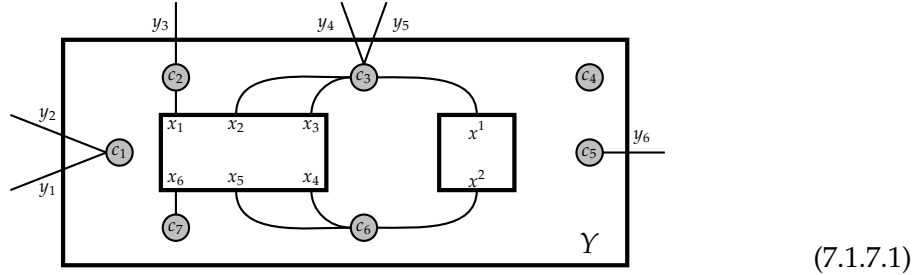
\begin{equation}
\label{uwd-first-picture}
\begin{tikzpicture}[scale=1]
\draw [ultra thick] (-1,0) rectangle (6.5,3);
\node at (-1.2,1) {\tiny{$y_1$}};
\node at (-1.2,2) {\tiny{$y_2$}};
\node at (.3,3.2) {\tiny{$y_3$}};
\node at (2.5,3.2) {\tiny{$y_4$}};
\node at (3.5,3.2) {\tiny{$y_5$}};
\node at (6.7,1.7) {\tiny{$y_6$}};
\node at (6,.3) {$Y$};
\draw [ultra thick] (.3,1) rectangle (2.5,2);
\node at (.5,1.85) {\tiny{$x_1$}};
\node at (1.3,1.85) {\tiny{$x_2$}};
\node at (2.3,1.85) {\tiny{$x_3$}};
\node at (2.3,1.15) {\tiny{$x_4$}};
\node at (1.3,1.15) {\tiny{$x_5$}};
\node at (.5,1.15) {\tiny{$x_6$}};
\draw [ultra thick] (4,1) rectangle (5,2);
\node at (4.5,1.15) {\tiny{$x^2$}};
\node at (4.5,1.8) {\tiny{$x^1$}};
\cable{(-.25,1.5)} \node at (-.25,1.5) {\tiny{$c_1$}};
\cable{(.5,2.5)} \node at (.5,2.5) {\tiny{$c_2$}};
\cable{(3,2.5)} \node at (3,2.5) {\tiny{$c_3$}};
\cable{(6,2.5)} \node at (6,2.5) {\tiny{$c_4$}};
\cable{(6,1.5)} \node at (6,1.5) {\tiny{$c_5$}};
\cable{(3,.5)} \node at (3,.5) {\tiny{$c_6$}};
\cable{(.5,.5)} \node at (.5,.5) {\tiny{$c_7$}};
\draw [thick] (-.425,1.5) -- (-1.7,1);
\draw [thick] (-.425,1.5) -- (-1.7,2);
\draw [thick] (.5,2) -- (.5,2.325);
\draw [thick] (.5,2.675) -- (.5,3.5);
\draw [thick] (1.3,2) to [out=90, in=180] (2.825,2.5);
\draw [thick] (2.3,2) to [out=90, in=180] (2.825,2.5);
\draw [thick] (4.5,2)  to [out=90, in=0] (3.175,2.5);
\draw [thick] (3,2.675) -- (2.7,3.5);
\draw [thick] (3,2.675) -- (3.3,3.5);
\draw [thick] (6.175,1.5) -- (7,1.5);
\draw [thick] (1.3,1) to [out=-90,in=180] (2.825,.5);
\draw [thick] (2.3,1) to [out=-90,in=180] (2.825,.5);
\draw [thick] (4.5,1) to [out=-90,in=0]  (3.175,.5);
\draw [thick] (.5,1) -- (.5,.675);
\end{tikzpicture}
\end{equation}
The input boxes $X_1 = \{x_1,\ldots,x_6\}$ and $X_2 = \{x^1,x^2\}$ are drawn as the smaller boxes inside.  The output box $Y = \{y_1,\ldots,y_6\}$ is drawn as the outer rectangle.  Each element within each box is drawn along the boundary, either just inside (as in $X_1$ and $X_2$) or just outside (as in $Y$).  Note that no orientation is attached to the sides of these squares and rectangles.  For example, we could have drawn $y_6 \in Y$ on the bottom side of the outer rectangle.  Each cable $c_i \in C$ is drawn as a small gray disk, which is not to be confused with a delay node in a wiring diagram (such as $d$ in \eqref{wd-first-example}).  The soldering functions $f$ and $g$ tell us how to connect the wires in $X = X_1 \amalg X_2$ and $Y$ to the cables.  We will revisit this example in Example \ref{ex:uwd-decomp} below.
\end{example}

\section{Pushouts}\label{sec:pushout}

The operadic composition on the collection of undirected wiring diagrams $\uwd$ (Def. \ref{def:uwd}) involves the basic categorical concept of a pushout, which we recall in this section.  The reader may consult \cite{awodey} (5.6) and \cite{maclane} (p.65-66) for more discussion of pushouts.  Roughly speaking, a pushout is a way of summing two objects with some identification.  We will only use the following definition when the category is $\Fins$ (Def. \ref{def:Fins}), so the reader may simply take the category $\C$ below to be $\Fins$ and objects and maps to be those in $\Fins$.

\begin{definition}
\label{def:pushout}
Suppose $\C$ is a category (Def. \ref{def:category}), and 
\begin{equation}
\label{pushout-diagram}
\nicexy{Y & X \ar[l]_-{f} \ar[r]^-{g} & Z}
\end{equation}
is a diagram in $\C$.  Then a \emph{pushout}\index{pushout} of this diagram is a tuple $(W,\alpha,\beta)$ consisting of an object $W \in \C$ and maps $\alpha : Y \to W$ and $\beta : Z \to W$ in $\C$ such that the following two conditions hold.
\begin{enumerate}
\item
The square
\[\nicexy{X \ar[d]_-{f} \ar[r]^-{g} & Z \ar[d]^-{\beta}\\ Y \ar[r]^-{\alpha} & W}\]
in $\C$ is commutative, i.e., $\alpha f = \beta g$.
\item
Suppose $(W',\alpha',\beta')$ is another such tuple, i.e., $\alpha' f = \beta' g$.  Then there exists a unique map $h : W \to W'$ in $\C$ such that the diagram
\begin{equation}
\label{pushout-universal-property}
\nicexy{X \ar[d]_-{f} \ar[r]^-{g} & Z \ar[d]^-{\beta} \ar@/^1pc/[ddr]^-{\beta'} &\\ 
Y \ar[r]^-{\alpha} \ar@/_1pc/[drr]_-{\alpha'} & W \ar[dr]|-{h} &\\
&& W'}
\end{equation}
in $\C$ is commutative, i.e., $\alpha' = h\alpha$ and $\beta' = h\beta$.
\end{enumerate} 
\end{definition}

If a pushout exists, then by definition it is unique up to unique isomorphisms.  In a general category, a pushout may not exist for a diagram of the form \eqref{pushout-diagram}.  Even if it exists, it may be difficult to describe.  Luckily, for $S$-finite sets (Def. \ref{def:Fins}), pushouts always exist and are easy to describe, as the following observation shows.

\begin{proposition}
\label{fins-pushouts-exist}
In the category $\Fins$ of $S$-finite sets, each diagram of the form \eqref{pushout-diagram} has a pushout given by the quotient
\begin{equation}
\label{pushout-formula}
W = \frac{Y \amalg Z}{\Bigl\{f(x) = g(x) : x \in X\Bigr\}}
\end{equation}
taken in $\Fins$.
\end{proposition}

\begin{proof}
The maps $\alpha : Y \to W$ and $\beta : Z \to W$ are the obvious maps, each being the inclusion into $Y \amalg Z$ followed by the quotient map to $W$.  Then the tuple $(W,\alpha,\beta)$ has the required universal property of a pushout in $\Fins$.
\end{proof}

\begin{example}
\label{ex:pushout-id}
A commutative square with opposite identity maps is a pushout square.  In other words, a pushout of the diagram
\[\nicexy{X & X \ar[l]_-{=} \ar[r]^-{g} & Z}\]
in any category is given by the commutative square
\[
\nicexy{
X \ar[d]_-{=} \ar[r]^-{g} & Z \ar[d]^-{=}\\
X \ar[r]^-{g} & Z}\]
as can be checked by a direct inspection.
\end{example}

\section{Operad Structure}
\label{sec:uwd-operad-structure}

Fix a class $S$.  In this section we define the $\Fins$-colored operad structure (Def. \ref{def:pseudo-operad}) on the collection of undirected wiring diagrams $\uwd$ (Def. \ref{def:uwd}).  When $S$ is either $\{*\}$ or the collection of sets, this operad structure on $\uwd$ was first introduced in \cite{spivak13} using the structure map $\gamma$ \eqref{operadic-composition}.

\begin{definition}[Equivariance in $\uwd$]
\label{def:uwd-equivariance}
Suppose $Y \in \Fins$, $\uX= (X_1,\ldots,X_n)$ is a $\Fins$-profile of length $n$, and $\sigma \in \Sigma_n$ is a permutation.  Define the function\index{equivariance in $\uwd$}
\begin{equation}
\label{uwd-right-action}
\nicexy{\uwd\yxonexn = \uwd\yux \ar[r]^-{\sigma}_-{\cong} 
& \uwd\yuxsigma = \uwd\smallbinom{Y}{X_{\sigma(1)}, \ldots, X_{\sigma(n)}}}
\end{equation}
by sending $\varphi = (C,f,g) \in \uwd\yux$ to $\varphi = (C,f,g) \in \uwd\yuxsigma$, using the fact that $\coprod_{i=1}^n X_i = \coprod_{i=1}^n X_{\sigma(i)}$.
\end{definition}

In other words, the equivariant structure in $\uwd$ simply relabels the input boxes.

Next we define the colored units in $\uwd$.  The $Y$-colored unit in $\uwd$, for $Y \in \Fins$, may be depicted as follows.
\begin{center}
\begin{tikzpicture}[scale=.5]
\draw [ultra thick] (-1,0) rectangle (3,2.5);
\node at (2.5,.4) {$Y$};
\draw [ultra thick] (1,.5) rectangle (2,2);
\node at (1.5,1.25) {$Y$};
\draw [thick] (1,.7) -- (-1.5,.7);
\draw [thick] (1,1.8) -- (-1.5,1.8);
\cable{(0,.7)} \node at (0,1.25) {{$\vdots$}}; \cable{(0,1.8)}
\end{tikzpicture}
\end{center}

\begin{definition}[Units in $\uwd$]
\label{def:uwd-units}
For each $Y \in \Fins$, the \emph{$Y$-colored unit} \index{units in $\uwd$} is defined as the undirected wiring diagram
\begin{equation}
\label{uwd-unit}
\tensorunit_Y = \Bigl(\nicexy{Y \ar[r]^-{=} & Y & Y \ar[l]_-{=}}\Bigr) \in \uwd\yy.
\end{equation}
\end{definition}

\begin{motivation}
Next we define the $\compi$-composition in $\uwd$.  The operadic composition $\varphi \compi \psi$ can be visualized in the following picture.
\begin{equation}
\label{uwd-compi-picture}
\begin{tikzpicture}[scale=.7]
\node at (6,5.3) {$\varphi$};
\draw [ultra thick] (0,0) rectangle (7,5);
\node at (3.5,2) {$\cdots$};
\draw [thick] (0,2.5) -- (-.5,2.5);
\draw [thick] (3.5,5) -- (3.5,5.5);
\draw [thick] (7,2.5) -- (7.5,2.5);
\draw [thick] (3.5,0) -- (3.5,-.5);
\draw [ultra thick] (1,1) rectangle (2,2);
\node at (1.5,1.5) {$X_1$};
\draw [thick] (1,1.5) -- (.5,1.5);
\draw [thick] (1.5,2) -- (1.5,2.5);
\draw [thick] (2,1.5) -- (2.5,1.5);
\draw [thick] (1.5,1) -- (1.5,.5);
\draw [ultra thick] (3,3) rectangle (4,4);
\node at (3.5,3.5) {$X_i$};
\draw [thick] (3,3.5) -- (2.5,3.5);
\draw [thick] (3.5,4) -- (3.5,4.5);
\draw [thick] (4,3.5) -- (4.5,3.5);
\draw [thick] (3.5,3) -- (3.5,2.5);
\draw [ultra thick] (5,1) rectangle (6,2);
\node at (5.5,1.5) {$X_n$};
\draw [thick] (5,1.5) -- (4.5,1.5);
\draw [thick] (5.5,2) -- (5.5,2.5);
\draw [thick] (6,1.5) -- (6.5,1.5);
\draw [thick] (5.5,1) -- (5.5,.5);
\draw [ultra thick] (1,6) rectangle (2,7);
\node at (1.5,6.5) {$\psi$};
\draw [thick] (1,6.5) -- (.5,6.5);
\draw [thick] (1.5,7) -- (1.5,7.5);
\draw [thick] (2,6.5) -- (2.5,6.5);
\draw [thick] (1.5,6) -- (1.5,5.5);
\draw [gray, ->, line width=3pt, semitransparent] (1.8,6) to [out=-90, in=180] (2.95,3.9);
\node at (2.2,5.3) {$\compi$};
\end{tikzpicture}
\end{equation}
Intuitively, to form the operadic composition $\varphi \compi \psi$ in $\uwd$, we replace the $i$th input box $X_i$ in $\varphi$ by the input boxes in $\psi$.  The set of cables in $\psi$ is added to the set of cables in $\varphi$, with appropriate identification from the input and the output soldering functions in $\varphi$ and $\psi$.
\end{motivation}

The following notation will be useful in the definition of the $\compi$-composition.

\begin{notation}
\label{fins-prof-notation}
Suppose $\uX = (X_1,\ldots,X_n)$ is a $\Fins$-profile.
\begin{enumerate}
\item
Write $X = X_1 \amalg \cdots \amalg X_n \in \Fins$.
\item
For integers $i$ and $j$, define
\begin{equation}
\label{x-sub-ij}
X_{[i,j]} = \begin{cases}
X_i \amalg \cdots \amalg X_j & \text{ if $1 \leq i \leq j \leq n$};\\
\varnothing & \text{ otherwise}.
\end{cases}\
\end{equation}
\end{enumerate}
\end{notation}

\begin{definition}[$\compi$-Composition in $\uwd$]
\label{def:compi-uwd}
Suppose:
\begin{itemize}
\item
$\varphi = \bigl\{\nicexy{X \ar[r]^-{\fphi} & \cphi & Y \ar[l]_-{\gphi}}\bigr\} \in \uwd\yux$ with $\uX = (X_1,\ldots,X_n)$, $n \geq 1$, and $1 \leq i \leq n$;
\item
$\psi = \bigl\{\nicexy{W \ar[r]^-{\fpsi} & \cpsi & X_i \ar[l]_-{\gpsi}}\bigr\} \in \uwd\xiuw$ with $|\uW| = m \geq 0$.
\end{itemize}
Define the undirected wiring diagram\index{compi-composition in UWD@$\compi$-composition in $\uwd$}
\[\phicompipsi = (C,f,g) \in \uwd\yxcompiw\]
as the cospan
\begin{equation}
\label{uwd-compi-cospan}
\nicexy@C+.3cm{
&& Y \ar[d]_-{\gphi} \ar@/^2pc/[dd]^-{g}\\
& X_1 \amalg \cdots \amalg X_n \ar[r]^-{\fphi} \ar[d]_-{(\Id, \gpsi, \Id)} \ar@{}[dr]|-{\mathrm{pushout}} & \cphi \ar[d]\\
X_{[1,i-1]} \amalg W \amalg X_{[i+1,n]} \ar[r]^-{(\Id, \fpsi, \Id)} \ar@/_2pc/[rr]_-{f}
& X_{[1,i-1]} \amalg \cpsi \amalg X_{[i+1,n]} \ar[r] & C}
\end{equation}
in $\Fins$.  Here:
\begin{itemize}
\item
The square is a pushout in $\Fins$, which always exists by Prop. \ref{fins-pushouts-exist}. 
\item
The $S$-finite sets $W$, $X_{[1,i-1]}$, and $X_{[i+1,n]}$ are as in Notation \ref{fins-prof-notation}.
\item
The $\Fins$-profile
\[
\xcompiw = \Bigl(X_1, \ldots, X_{i-1}, \uW, X_{i+1}, \ldots , X_n\Bigr)\]
is as in \eqref{compi-profile}.
\item
The map $f$ is the bottom horizontal composition, and the map $g$ is the right vertical composition.
\end{itemize}
\end{definition}

In the following observation, we describe the undirected wiring diagram $\phicompipsi$ more explicitly.

\begin{proposition}
\label{uwd-compi-explicit}
Consider the diagram \eqref{uwd-compi-cospan}.
\begin{enumerate}
\item
A choice of a pushout $C$ is the quotient
\begin{equation}
\label{uwd-compi-pushout-choice}
C = \frac{\cphi \amalg \cpsi}{\left\{\fphi(x) = \gpsi(x) : x \in X_i\right\}}
\end{equation}
in $\Fins$.   The following statements use this choice of $C$.
\item
The maps $\cphi \to C$ and $\cpsi \to C$ are the obvious maps, each being the inclusion into $\cphi \amalg \cpsi$ followed by the quotient map to $C$. 
\item
In the map $f$ and in the horizontal unnamed map, for $j \not= i$, the map $X_j \to C$ is the composition
\[
\nicexy{X_j \ar[r]^-{\fphi} & \cphi \ar[r] & C}.\]
\item
On $W$ the map $f$ is the composition
\[
\nicexy{W \ar[r]^-{\fpsi} & \cpsi \ar[r] & C}.\]
\end{enumerate}
\end{proposition}

\begin{proof}
A direct inspection shows that the first three statements indeed describe a pushout of the square in the diagram \eqref{uwd-compi-cospan}.  The last assertion follows from the definition of $f$ as the bottom horizontal composition.
\end{proof}

\begin{remark}
Consider Def. \ref{def:compi-uwd}.
\begin{enumerate}
\item
Pushouts are unique up to unique isomorphisms.  Since undirected wiring diagrams are defined as equivalence classes of undirected prewiring diagrams (Def. \ref{def:uwd}), the undirected wiring diagram $\phicompipsi$ is well-defined.
\item
In \cite{spivak13} $S$ was taken to be either the one-point set or the collection of sets.   The operadic composition in $\uwd$ was defined in terms of the operadic composition $\gamma$ \eqref{operadic-composition}.  By Prop. \ref{prop:operad-def-equiv} the two descriptions--i.e., the one in \cite{spivak13} and Def. \ref{def:compi-uwd}--are equivalent.
\end{enumerate}
\end{remark}

Recall from Def. \ref{def:pre-uwd} that a \emph{wasted cable} in an undirected wiring diagram is a cable that is not in the image of the input and the output soldering functions.

\begin{example}
\label{ex:uwd-jointly-surjective}
In this example, we observe that wasted cables can be created by the $\compi$-composition, even if the original undirected wiring diagrams have no wasted cables.  Suppose:
\begin{itemize}
\item
$S = \{*\}$ and $X = \{x_1,x_2\} \in \Fin$ is a two-element set.
\item
$\varphi = \left\{\nicexy{X \ar[r] & \{*\} & \varnothing \ar[l]}\right\} \in \uwd\emptyx$.
\item
$\psi = \left\{\nicexy{\varnothing \ar[r] & X & X \ar[l]_-{=}}\right\} \in \uwd\xempty$.
\end{itemize}
Note that neither $\varphi$ nor $\psi$ has a wasted cable.  On the other hand, the undirected wiring diagram $\varphi \compone \psi \in \uwd\emptyprof$ is the cospan
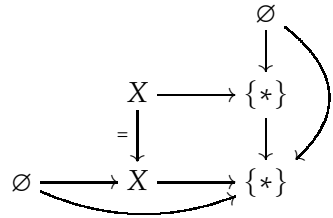
\[\nicexy@R-.3cm{
&& \varnothing \ar[d] \ar@/^2pc/[dd]\\
& X \ar[d]_-{=} \ar[r] & \{*\} \ar[d]\\
\varnothing \ar[r] \ar@/_1pc/[rr] & X \ar[r] & \{*\}}\]
in $\Fin$, in which the square is a pushout by Example \ref{ex:pushout-id}.  The following picture gives a visualization of this $\compone$-composition.
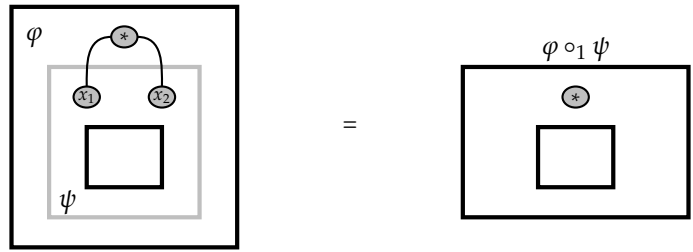
\begin{center}
\begin{tikzpicture}[xscale=1,yscale=.8]
\draw [ultra thick] (0,0) rectangle (3,4);
\draw [ultra thick, gray, semitransparent] (.5,.5) rectangle (2.5,3);
\draw [ultra thick] (1,1) rectangle (2,2);
\cable{(1,2.5)} \node at (1,2.5) {\tiny{$x_1$}};
\cable{(2,2.5)} \node at (2,2.5) {\tiny{$x_2$}};
\cable{(1.5,3.5)} \node at (1.5,3.5) {\tiny{$*$}};
\draw [thick] (1,2.675) to [out=90,in=180] (1.325,3.5);
\draw [thick] (2,2.675) to [out=90,in=0] (1.675,3.5);
\node at (.75,.75) {\small{$\psi$}};
\node at (.3,3.5) {\small{$\varphi$}};
\node at (4.5,2) {$=$};
\draw [ultra thick] (6,0.5) rectangle (9,3);
\draw [ultra thick] (7,1) rectangle (8,2);
\cable{(7.5,2.5)} \node at (7.5,2.5) {\tiny{$*$}};
\node at (7.5,3.3) {\small{$\phicomponepsi$}};
\end{tikzpicture}
\end{center}
In particular, the unique cable in $\phicomponepsi$ is a wasted cable.  This example illustrates that the $\compi$-composition of two undirected wiring diagrams without wasted cables may have wasted cables, which we mentioned in Remark \ref{rk:spivak-def-uwd}.  We will revisit this example in Section \ref{sec:wasted-cable} and Example \ref{ex:wasted-cable-simplex} below.
\end{example}
 
We now prove that the collection $\uwd$ of undirected wiring diagrams is a $\Fins$-colored operad in the sense of Def. \ref{def:pseudo-operad}.

\begin{lemma}
\label{uwd-compi-unity}
The $\compi$-composition in Def. \ref{def:compi-uwd} satisfies the left unity axiom \eqref{compi-left-unity}, the right unity axiom \eqref{compi-right-unity}, and the equivariance axiom \eqref{compi-eq}.\index{unity in $\uwd$}\index{equivariance in $\uwd$}
\end{lemma}

\begin{proof}
The equivariance axiom holds because the equivariant structure \eqref{def:uwd-equivariance} simply relabels the input boxes.  The unity axioms follow from the definitions of the colored units in $\uwd$ \eqref{uwd-unit} and Example \ref{ex:pushout-id}.
\end{proof}

\begin{motivation}
The horizontal associativity axiom in $\uwd$ may be visualized as follows.
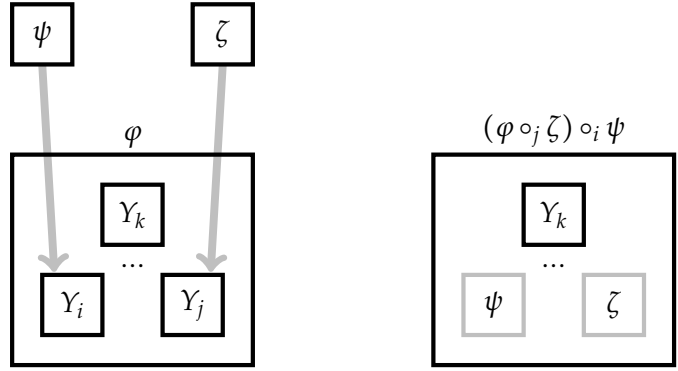
\begin{center}
\begin{tikzpicture}[scale=.8]
\draw [lightgray, ->, line width=3pt] (1.5,5) to (1.7,1.55);
\draw [lightgray, ->, line width=3pt] (4.5,5) to (4.3,1.55);
\draw [ultra thick] (1,5) rectangle (2,6);
\node at (1.5,5.5) {$\psi$};
\draw [ultra thick] (4,5) rectangle (5,6);
\node at (4.5,5.5) {$\zeta$};
\draw [ultra thick] (1,0) rectangle (5,3.5);
\node at (3,3.8) {$\varphi$};
\node at (3,1.6) {$\cdots$};
\draw [ultra thick] (1.5,.5) rectangle (2.5,1.5);
\node at (2,1) {$Y_i$};
\draw [ultra thick] (3.5,.5) rectangle (4.5,1.5);
\node at (4,1) {$Y_j$};
\draw [ultra thick] (2.5,2) rectangle (3.5,3);
\node at (3,2.5) {$Y_k$};
\draw [ultra thick] (8,0) rectangle (12,3.5);
\node at (10,3.9) {$(\varphi \compj \zeta) \compi \psi$};
\node at (10,1.6) {$\cdots$};
\draw [ultra thick, lightgray] (8.5,.5) rectangle (9.5,1.5);
\node at (9,1) {$\psi$};
\draw [ultra thick, lightgray]  (10.5,.5) rectangle (11.5,1.5);
\node at (11,1) {$\zeta$};
\draw [ultra thick] (9.5,2) rectangle (10.5,3);
\node at (10,2.5) {$Y_k$};
\end{tikzpicture}
\end{center}
To keep the picture simple, we omitted drawing the wires and the cables.  Note that this is basically the undirected version of the picture \eqref{wd-hass-picture}.
\end{motivation}

\begin{lemma}
\label{uwd-compi-horizontal-associative}
The $\compi$-composition in Def. \ref{def:compi-uwd} satisfies the horizontal associativity axiom \eqref{compi-associativity}.\index{horizontal associativity in $\uwd$}
\end{lemma}

\begin{proof}
Suppose:
\begin{itemize}
\item
$\varphi = \left(\cphi, \fphi, \gphi\right) \in \uwd\zuy$ with $|\uY| = n \geq 2$ and $1 \leq i < j \leq n$;
\item
$\psi = \left(\cpsi, \fpsi, \gpsi\right) \in \uwd\yiuw$ with $|\uW| = l$;
\item
$\zeta = \left(\czeta, \fzeta, \gzeta\right) \in \uwd\yjux$ with $|\uX| = m$.
\end{itemize}
We must show that
\begin{equation}
\label{uwd-horizontal-associativity}
\left( \varphi \compj \zeta\right) \compi \psi
= 
\left(\varphi \compi \psi\right) \compjonel \zeta \in \uwd\zycompjxcompiw.
\end{equation}

Consider the undirected wiring diagram
\[
\left(C,f,g\right) \in \uwd\zycompjxcompiw\]
given by the cospan
\[\nicexy{Y_{[1,i-1]} \amalg W \amalg Y_{[i+1,j-1]} \amalg X \amalg Y_{[j+1,n]} \ar[r]^-{f} & C & Z \ar[l]_-{g}}\]
in $\Fins$, where Notation \ref{fins-prof-notation} was used.  In this cospan:
\begin{itemize}
\item
The set of cables is  the quotient
\[
C = \frac{\cphi \amalg \cpsi \amalg \czeta}{\Bigl\{\fphi(y_i) = \gpsi(y_i),\, 
\fphi(y_j) = \gzeta(y_j) : y_i \in Y_i,\, y_j \in Y_j\Bigr\}.}\]
\item
The output soldering function $g$ is the composition $\nicexy{Z \ar[r]^-{\gphi} & \cphi \ar[r] & C}$.
\item
For the input soldering function $f$:
\begin{itemize}
\item
The restriction to $Y_k$ is the composition $\nicexy{Y_k \ar[r]^-{\fphi} & \cphi \ar[r] & C}$ for $k\not= i,j$.
\item
The restriction to $W$ is the composition $\nicexy{W \ar[r]^-{\fpsi} & \cpsi \ar[r] & C}$.
\item
The restriction to $X$ is the composition $\nicexy{X \ar[r]^-{\fzeta} & \czeta \ar[r] & C}$.
\end{itemize}
\end{itemize}
Using the description of $\compi$ in Prop. \ref{uwd-compi-explicit}, a direct inspection reveals that both sides of \eqref{uwd-horizontal-associativity} are equal to $(C,f,g)$.
\end{proof}

\begin{motivation}
The vertical associativity axiom in $\uwd$ may be visualized as follows.
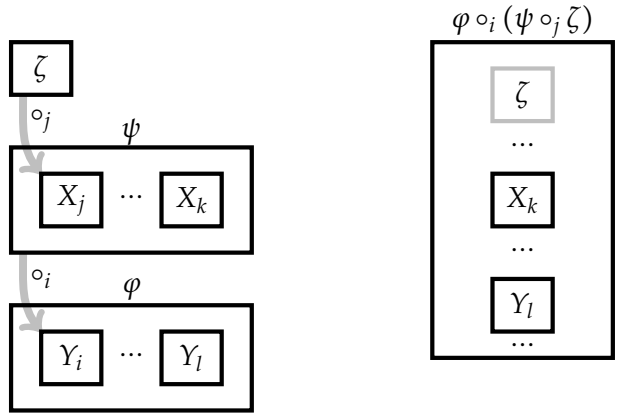
\begin{center}
\begin{tikzpicture}[xscale=.8,yscale=.7]
\draw [lightgray, ->, line width=3pt] (1.2,5) to [out=-90, in=125] (1.5,3.5);
\node at (1.5,4.5) {$\compj$};
\draw [lightgray, ->, line width=3pt] (1.2,2) to [out=-90, in=125]  (1.5,.5);
\node at (1.5,1.5) {$\compi$};
\draw [ultra thick] (1,5) rectangle (2,6);
\node at (1.5,5.5) {$\zeta$};
\draw [ultra thick] (1,2) rectangle (5,4);
\node at (3,4.3) {$\psi$};
\node at (3,3) {$\cdots$};
\draw [ultra thick] (1.5,2.5) rectangle (2.5,3.5);
\node at (2,3) {$X_j$};
\draw [ultra thick] (3.5,2.5) rectangle (4.5,3.5);
\node at (4,3) {$X_k$};
\draw [ultra thick] (1,-1) rectangle (5,1);
\node at (3,1.3) {$\varphi$};
\node at (3,0) {$\cdots$};
\draw [ultra thick] (1.5,-.5) rectangle (2.5,.5);
\node at (2,0) {$Y_i$};
\draw [ultra thick] (3.5,-.5) rectangle (4.5,0.5);
\node at (4,0) {$Y_l$};
\draw [ultra thick] (8,0) rectangle (11,6);
\node at (9.5,6.4) {$\varphi \compi (\psi \compj \zeta)$};
\node at (9.5,.2) {$\cdots$};
\node at (9.5,2) {$\cdots$};
\node at (9.5,4) {$\cdots$};
\draw [ultra thick, gray, semitransparent]  (9,4.5) rectangle (10,5.5);
\node at (9.5,5) {$\zeta$};
\draw [ultra thick] (9,2.5) rectangle (10,3.5);
\node at (9.5,3) {$X_k$};
\draw [ultra thick] (9,.5) rectangle (10,1.5);
\node at (9.5,1) {$Y_l$};
\end{tikzpicture}
\end{center}
As before, to keep the picture simple, we did not draw the wires and the cables.  Note that this is basically the undirected version of the picture \eqref{wd-vass-picture}.
\end{motivation}

\begin{lemma}
\label{uwd-compi-vertical-associative}
The $\compi$-composition in Def. \ref{def:compi-uwd} satisfies the vertical associativity axiom \eqref{compi-associativity-two}.\index{vertical associativity in $\uwd$}
\end{lemma}

\begin{proof}
Suppose:
\begin{itemize}
\item
$\varphi = \left(\cphi,\fphi,\gphi\right) \in \uwd\zuy$ with $|\uY| = n \geq 1$ and $1 \leq i \leq n$;
\item
$\psi = \left(\cpsi,\fpsi,\gpsi\right) \in \uwd\yiux$ with $|\uX| = m \geq 1$ and $1 \leq j \leq m$;
\item
$\zeta = \left(\czeta,\fzeta,\gzeta\right) \in \uwd\xjuw$ with $|\uW| = l$.
\end{itemize}
We must show that
\begin{equation}
\label{uwd-vertical-assoc}
\left(\varphi \compi \psi\right) \compionej \zeta 
=
\varphi \compi \left(\psi \compj \zeta\right)
\in \uwd\zycompixcompionejw.
\end{equation}

Consider the undirected wiring diagram
\[
\left(C,f,g\right) \in \uwd\zycompixcompionejw\]
given by the cospan
\[\nicexy{Y_{[1,i-1]} \amalg X_{[1,j-1]} \amalg W \amalg X_{[j+1,m]} \amalg Y_{[i+1,n]} \ar[r]^-{f} & C & Z \ar[l]_-{g}}\]
in $\Fins$, where Notation \ref{fins-prof-notation} was used.  In this cospan:
\begin{itemize}
\item
The set of cables is the quotient
\[
C = \frac{\cphi \amalg \cpsi \amalg \czeta}{\Bigl\{\fphi(y) = \gpsi(y),\, 
\fpsi(x) = \gzeta(x) : y \in Y_i,\, x \in X_j\Bigr\}.}\]
\item
The output soldering function $g$ is the composition $\nicexy{Z \ar[r]^-{\gphi} & \cphi \ar[r] & C}$.
\item
For the input soldering function $f$:
\begin{itemize}
\item
The restriction to $Y_l$ is the composition $\nicexy{Y_l \ar[r]^-{\fphi} & \cphi \ar[r] & C}$ for $l\not= i$.
\item
The restriction to $X_k$ is the composition $\nicexy{X_k \ar[r]^-{\fpsi} & \cpsi \ar[r] & C}$ for $k \not=j$.
\item
The restriction to $W$ is the composition $\nicexy{W \ar[r]^-{\fzeta} & \czeta \ar[r] & C}$.
\end{itemize}
\end{itemize}
Using the description of $\compi$ in Prop. \ref{uwd-compi-explicit}, a direct inspection reveals that both sides of \eqref{uwd-vertical-assoc} are equal to $(C,f,g)$.
\end{proof}

\begin{theorem}
\label{uwd-operad}
For any class $S$, when equipped with the structure in Def. \ref{def:uwd-equivariance}--\ref{def:compi-uwd}, $\uwd$ in Def. \ref{def:uwd} is a $\Fins$-colored operad, called the operad of undirected wiring diagrams. \index{UWD is an operad@$\uwd$ is an operad}
\end{theorem}

\begin{proof}
In view of Def. \ref{def:pseudo-operad}, this follows from Lemmas \ref{uwd-compi-unity}, \ref{uwd-compi-horizontal-associative}, and \ref{uwd-compi-vertical-associative}.
\end{proof}

\begin{example}
\label{ex:spivak-uwd-operads}
Consider Theorem \ref{uwd-operad}.
\begin{enumerate}
\item
If $S = \{*\}$, a one-point set, then our $\Fin$-colored operad $\uwd$ is called the \emph{operad of singly-typed wiring diagrams} in \cite{spivak13} (Example 2.1.7).
\item
If $S = \set$, the collection of sets, then our $\Fin_{\set}$-colored operad $\uwd$ is called the \emph{operad of typed wiring diagrams} in \cite{spivak13} (Example 4.1.1).
\end{enumerate}
\end{example}

We defined the operad $\uwd$ in terms of the $\compi$-compositions (Def. \ref{def:compi-uwd}).  The following observation expresses the operad $\uwd$ in terms of the operadic composition $\gamma$ \eqref{operadic-composition}.  In \cite{spivak13} the operad structure on undirected wiring diagrams was actually defined in terms of $\gamma$.

\begin{proposition}
\label{uwd-gamma}
Suppose:
\begin{itemize}
\item
$\varphi = (\cphi, \fphi, \gphi) \in \uwd\yux$ with $\uX = (X_1,\ldots,X_n)$ for some $n \geq 1$ and $X = X_1 \amalg \cdots \amalg X_n$.
\item
For each $1 \leq i \leq n$, $\psi_i = \left(C_i,f_i,g_i\right) \in \uwd\xiuwi$ with $\uW_i = \left(W_{i,1},\ldots,W_{i,k_i}\right)$ for some $k_i \geq 0$.
\item
$\uW =  (\uW_1,\ldots,\uW_n)$, $W_i = W_{i,1} \amalg \cdots \amalg W_{i,k_i}$, and $W = \coprod_{1\leq i \leq n} W_i$.
\end{itemize}
Then
\[\gamma\Bigl(\varphi; \psi_1, \ldots, \psi_n\Bigr) = (C,f,g) \in \uwd\yuw\]
is given by the cospan
\begin{equation}
\label{uwd-gamma-cospan}
\nicexy{&& Y \ar[d]_-{\gphi} \ar@/^2pc/[dd]^-{g}\\
& X \ar[r]^-{\fphi} \ar[d]_-{\amalg g_i} \ar@{}[dr]|-{\mathrm{pushout}} & \cphi \ar[d]\\
W \ar[r]^-{\amalg f_i} \ar@/_2pc/[rr]|-{f} & C_1 \amalg \cdots \amalg C_n \ar[r] & C}
\end{equation}
in $\Fins$, in which the square is a pushout.  In this diagram:
\begin{enumerate}
\item
$C$ is the quotient
\[C = \frac{\cphi \amalg C_1 \amalg \cdots \amalg C_n}{\Bigl\{\fphi(x) = g_i(x) : x \in X_i,\, 1 \leq i \leq n\Bigr\}}\]
in $\Fins$.
\item
The maps $\cphi \to C$ and $C_i \to C$ are the obvious maps, each being an  inclusion followed by a quotient map to $C$. 
\item
The restriction of $f$ to $W_i$ is the composition of $f_i : W_i \to C_i$ and $C_i \to C$.
\end{enumerate}
\end{proposition}

\begin{proof}
This follows from (i) the correspondence \eqref{gamma-in-comps} between $\gamma$ and the $\compi$-compositions and (ii) the description of $\compi$ given in Proposition \ref{uwd-compi-explicit}.
\end{proof}

\section{Summary of Chapter \ref{ch07-undirected-wiring-diagrams}}

\begin{enumerate}
\item An undirected $S$-wiring diagram has a finite number of input boxes, an output box, an $S$-finite set of cables, an input soldering function, and an output soldering function.
\item For each class $S$, the collection of $S$-wiring diagrams $\uwd$ is a $\Fins$-colored operad.
\end{enumerate}

\chapter{Generators and Relations}
\label{ch08-generating-uwd}

Fix a class $S$, and consider the $\Fins$-colored operad $\uwd$ of undirected wiring diagrams (Theorem \ref{uwd-operad}).  The purpose of this chapter is to describe a finite number of undirected wiring diagrams that we will later show to be sufficient to describe the entire operad $\uwd$.  One may also regard this chapter as consisting of a long list of examples of undirected wiring diagrams.

In Section \ref{sec:generating-uwd} we describe $6$ undirected wiring diagrams, called the \emph{generating undirected wiring diagrams}.  Later we will show that they generate the operad $\uwd$ of undirected wiring diagrams.  This means that every undirected wiring diagram can be obtained from finitely many generating undirected wiring diagrams via iterated operadic compositions.  For now one may think of the generating undirected wiring diagrams as examples of undirected wiring diagrams.


In Section \ref{sec:elementary-relations-uwd} we describe $17$ \emph{elementary relations} among the generating undirected wiring diagrams.   Later we will show that these elementary relations together with the operad associativity and unity axioms--\eqref{compi-associativity}, \eqref{compi-associativity-two}, \eqref{compi-left-unity}, and \eqref{compi-right-unity}--for the generating undirected wiring diagrams generate \emph{all} the relations in the operad $\uwd$ of undirected wiring diagrams.  In other words, suppose an arbitrary undirected wiring diagram can be built in two ways using the generating undirected wiring diagrams.  Then there exists a finite sequence of steps connecting them in which each step is given by one of the $17$ elementary relations or an operad associativity/unity axiom for the generating undirected wiring diagrams.   For now one may think of the elementary relations as examples of the operadic composition in the operad $\uwd$.

\section{Generating Undirected Wiring Diagrams}
\label{sec:generating-uwd}

Recall the definition of an undirected wiring diagram (Def. \ref{def:uwd}).  In this section, we introduce $6$ undirected wiring diagrams, called the generating undirected wiring diagrams.  They will be used in later chapters to give a finite presentation for the operad $\uwd$ of undirected wiring diagrams.  The undirected wiring diagrams in this section all have directed analogues in Section \ref{sec:generating-wd}.

The following undirected wiring diagram is an undirected analogue of the empty wiring diagram (Def. \ref{def:empty-wd}).

\begin{definition}
\label{def:uwd-gen-empty}
Define the \emph{empty cell}\index{empty cell}\index{epsilon@$\epsilon$}
\[\epsilon  = \Bigl(\nicexy{\varnothing \ar[r] & \varnothing & \varnothing \ar[l]}\Bigr) 
\in \uwd\emptynothing,\]
where $\varnothing$ is the empty $S$-finite set.  Note that the empty cell is a $0$-ary element in $\uwd$.
\end{definition}

Next we define the undirected wiring diagram
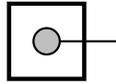
\begin{center}
\begin{tikzpicture}[scale=1]
\draw [ultra thick] (1,0.5) rectangle (2,1.5);
\cable{(1.5,1)}
\draw [thick] (1.675,1) -- (2.5,1);
\end{tikzpicture}
\end{center}
with no input boxes and whose unique cable is a $(0,1)$-cable.  This is an undirected analogue of a $1$-wasted wire (Def. \ref{def:wasted-wire-wd}).  To simplify the typography, we will often write $x$ for the one-point set $\{x\}$.

\begin{definition}
\label{def:uwd-gen-1output}
Suppose $* \in \Fins$ is a one-element $S$-finite set.  Define the \emph{$1$-output wire}\index{1-output wire@$1$-output wire}\index{omegasub@$\omega_*$}
\[\omega_* = \Bigl(\nicexy{\varnothing \ar[r] & \ast & \ast \ar[l]}\Bigr) \in \uwd\starnothing.\]
Note that a $1$-output wire is a $0$-ary element in $\uwd$.
\end{definition}

Next we define the undirected wiring diagram
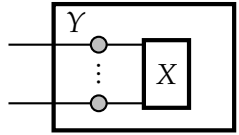
\begin{center}
\begin{tikzpicture}[scale=.6]
\draw [ultra thick] (-1,0) rectangle (3,2.8);
\node at (-.5,2.4) {$Y$};
\draw [ultra thick] (1,.5) rectangle (2,2);
\node at (1.5,1.25) {$X$};
\draw [thick] (1,.6) -- (-2,.6); \cable{(0,.6)}
\node at (0,1.25) {{$\vdots$}};
\draw [thick] (1,1.9) -- (-2,1.9); \cable{(0,1.9)}
\end{tikzpicture}
\end{center}
with $1$ input box and whose cables are all $(1,1)$-cables.  This is an undirected analogue of a name change (Def. \ref{def:name-change}).

\begin{definition}
\label{def:uwd-gen-namechange}
Suppose $f : X \to Y \in \Fins$ is a bijection.  Define the \emph{undirected name change}\index{undirected name change} \index{tauxy@$\tau_{X,Y}$}
\[\tau_f  = \Bigl(\nicexy{X \ar[r]^-{f}_-{\cong} & Y & Y \ar[l]_-{=}} \Bigr) \in \uwd\yx.\]
If the bijection $f$ is clear from the context, then we write $\tau_f$ as $\tau_{X,Y}$ or just $\tau$.  If there is no danger of confusion, then we will call $\tau_f$ a \emph{name change}.
\end{definition}

Next we define the undirected wiring diagram
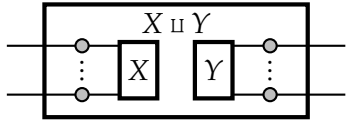
\begin{center}
\begin{tikzpicture}[scale=.5]
\draw [ultra thick] (-1,0) rectangle (6,3);
\node at (2.5,2.5) {$X \amalg Y$};
\draw [ultra thick] (1,.5) rectangle (2,2);
\node at (1.5,1.25) {$X$};
\draw [thick] (1,.6) -- (-2,.6); \cable{(0,.6)}
\node at (0,1.25) {{$\vdots$}};
\draw [thick] (1,1.9) -- (-2,1.9); \cable{(0,1.9)}
\draw [ultra thick] (3,.5) rectangle (4,2);
\node at (3.5,1.25) {$Y$};
\draw [thick] (4,.6) -- (7,.6); \cable{(5,.6)}
\node at (5,1.25) {{$\vdots$}};
\draw [thick] (4,1.9) -- (7,1.9); \cable{(5,1.9)}
\end{tikzpicture}
\end{center}
with $2$ input boxes and whose cables are all $(1,1)$-cables.  This is an undirected analogue of a $2$-cell (Def. \ref{def:theta-wd}).

\begin{definition}
\label{def:uwd-gen-2cell}
Suppose $X,Y \in \Fins$ and $X \amalg Y$ is their coproduct.  Define the \emph{undirected $2$-cell}\index{undirected 2-cell@undirected $2$-cell} \index{thetaxy@$\theta_{(X,Y)}$}
\[\theta_{(X,Y)} = \Bigl(\nicexy{X \amalg Y \ar[r]^-{=} & X \amalg Y & X \amalg Y \ar[l]_-{=}}\Bigr) \in \uwd\xplusyxy.\]
If there is no danger of confusion, then we will call it a \emph{$2$-cell}.
\end{definition}

Next we define the undirected wiring diagram
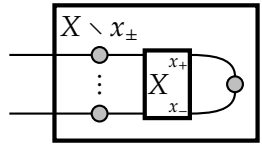
\begin{center}
\begin{tikzpicture}[scale=.6]
\draw [ultra thick] (-1,0) rectangle (3.5,3);
\node at (0,2.5) {$X \setminus x_{\pm}$};
\draw [ultra thick] (1,.5) rectangle (2,2);
\node at (1.75,.7) {\tiny{$x_-$}};
\node at (1.75,1.7) {\tiny{$x_+$}};
\node at (1.3,1.25) {$X$};
\draw [thick] (1,.6) -- (-2,.6); \cable{(0,.6)}
\node at (0,1.25) {{$\vdots$}};
\draw [thick] (1,1.9) -- (-2,1.9); \cable{(0,1.9)}
\cable{(3,1.25)}
\draw [thick] (2,1.9) to [out=0,in=90] (3,1.425);
\draw [thick] (2,.6) to [out=0,in=-90] (3,1.075);
\end{tikzpicture}
\end{center}
with a $(2,0)$-cable, all other cables being $(1,1)$-cables.  This is an undirected analogue of a $1$-loop (Def. \ref{def:loop-wd}).

\begin{definition}
\label{def:uwd-gen-loop}
Suppose:
\begin{itemize}
\item
$X \in \Fins$, and $x_+,x_- \in X$ are two distinct elements with $v(x_+) = v(x_-) \in S$.
\item
$\xminusplusminus \in \Fins$ is obtained from $X$ by removing $x_+$ and $x_-$.
\item
$X/(x_+=x_-) \in \Fins$ is the quotient of $X$ with $x_+$ and $x_-$ identified.
\end{itemize}
Define the \emph{loop}\label{uwd-loop}\index{loop} \index{lambdasubx@$\lambda_{(X,\xsubpm)}$}
\[\lambda_{(X,\xsubpm)} = \Bigl(\nicexy@C+.6cm{X \ar[r]^-{\mathrm{projection}} & \frac{X}{(x_+=x_-)} & \xminusplusminus \ar[l]_-{\mathrm{inclusion}}}\Bigr) \in \uwd\xminusplusminusx.\]
\end{definition}

Next we define the undirected wiring diagram
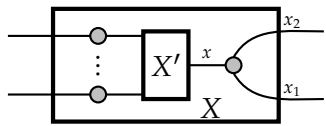
\begin{center}
\begin{tikzpicture}[scale=.6]
\draw [ultra thick] (-1,0) rectangle (4,2.5);
\node at (4.3,.7) {\tiny{$x_1$}};
\node at (4.3,2.2) {\tiny{$x_2$}};
\node at (2.5,.3) {$X$};
\draw [ultra thick] (1,.5) rectangle (2,2);
\node at (1.5,1.25) {$X'$};
\draw [thick] (1,.6) -- (-2,.6); \cable{(0,.6)}
\node at (0,1.25) {{$\vdots$}};
\draw [thick] (1,1.9) -- (-2,1.9); \cable{(0,1.9)}
\cable{(3,1.25)}
\draw [thick] (2,1.25) -- (2.825,1.25);
\node at (2.4,1.5) {\tiny{$x$}};
\draw [thick] (3,1.425) to [out=60,in=180] (5,2);
\draw [thick] (3,1.075) to [out=-60,in=180] (5,.5);
\end{tikzpicture}
\end{center}
with a $(1,2)$-cable, all other cables being $(1,1)$-cables.  This is an undirected analogue of an out-split (Def. \ref{def:out-split}).

\begin{definition}
\label{def:uwd-gen-split}
Suppose:
\begin{itemize}
\item
$X \in \Fins$, and $x_1, x_2$ are two distinct elements in $X$ with $v(x_1) = v(x_2) \in S$.
\item
$X' \in \Fins$, and $x \in X'$ such that $v(x) = v(x_1)$ and that $X' \setminus \{x\} = X \setminus \{x_1,x_2\}$.
\end{itemize}
Define the \emph{split}\label{uwd-split}\index{split} \index{sigmasup@$\sigma^{(X,x_1,x_2)}$} 
\[\sigma^{(X,x_1,x_2)} = \Bigl(\nicexy{X' \ar[r]^-{=} & X' & X \ar[l]}\Bigr) \in \uwd\xxprime\]
in which the output soldering function $X \to X'$ sends $x_1, x_2 \in X$ to $x \in X'$ and is the identity function on $X \setminus \{x_1,x_2\}$.
\end{definition}

\begin{definition}
\label{def:generating-uwd}
The $6$ undirected wiring diagrams in Def. \ref{def:uwd-gen-empty}--\ref{def:uwd-gen-split} will be referred to as \emph{generating undirected wiring diagrams}.\index{generating undirected wiring diagrams}  If the context is clear, we will simply call them \emph{generators}.\index{generating undirected wiring diagrams}\index{generators}
\end{definition}

\begin{remark}
\label{rk:uwd-generators-arity}
Among the generating undirected wiring diagrams:
\begin{enumerate}
\item
None has a wasted cable (Def. \ref{def:pre-uwd}).  As we will see in Section \ref{sec:wasted-cable}, wasted cables can be created by the generators.
\item
The empty cell $\epsilon$ (Def. \ref{def:uwd-gen-empty}) and a $1$-output wire $\omega_*$ (Def. \ref{def:uwd-gen-1output}) are $0$-ary elements in $\uwd$.
\item
A name change $\tau$ (Def. \ref{def:uwd-gen-namechange}), a loop $\lambda_{(X,\xsubpm)}$ (Def. \ref{def:uwd-gen-loop}), and a split $\sigma^{(X,x^1,x^2)}$ (Def. \ref{def:uwd-gen-split}) are unary elements in $\uwd$.
\item
A $2$-cell $\theta_{(X,Y)}$ (Def. \ref{def:uwd-gen-2cell}) is a binary element in $\uwd$.
\end{enumerate}
\end{remark}

\section{Elementary Relations}
\label{sec:elementary-relations-uwd}

The purpose of this section is to introduce $17$ elementary relations among the generating undirected wiring diagrams (Def. \ref{def:generating-uwd}).  Each elementary relation is proved using Prop. \ref{uwd-compi-explicit} and Example \ref{ex:pushout-id} and by a simple inspection of the relevant definitions of the generating undirected wiring diagrams and operadic compositions.  Each proof is  similar to Example \ref{ex:uwd-jointly-surjective} and the proofs of Lemma \ref{uwd-compi-horizontal-associative} and Lemma \ref{uwd-compi-vertical-associative}.  Therefore, we will omit the proofs, providing a picture instead in most cases.  Some, but not all, of the following relations have directed analogues in Section \ref{sec:elementary-relations}.

Recall the operadic composition in the $\Fins$-colored operad $\uwd$ (Def. \ref{def:compi-uwd}) and Notation \ref{comp-is-compone} for (iterated) $\compone$.  The first five relations are about name changes (Def. \ref{def:uwd-gen-namechange}).  The first one says that two consecutive name changes can be composed down into one name change.

\begin{proposition}
\label{prop:uwd-move-a1}
Suppose:
\begin{itemize}
\item
$f : X \to Y$ and $g : Y \to Z \in \Fins$ are bijections.  
\item
$\tau_f \in \uwd\yx$, $\tau_g \in \uwd\zy$, and $\tau_{gf} \in \uwd\zx$ are the corresponding name changes.
\end{itemize}
Then
\begin{equation}
\label{uwd-move-a1}
\tau_g \comp \tau_f = \tau_{gf} \in \uwd\zx.
\end{equation}
\end{proposition}

The next relation says that a name change of a $1$-output wire (Def. \ref{def:uwd-gen-1output}) is again a $1$-output wire.

\begin{proposition}
\label{prop:uwd-move-a2}
Suppose:
\begin{itemize}
\item
$X = \{x\}$ and $Y = \{y\}$ are two $1$-element $S$-finite sets with $v(x) = v(y) \in S$.
\item
$\omega_x \in \uwd\xnothing$ and $\omega_y \in \uwd\ynothing$ are the corresponding $1$-output wires.
\item
$\tau_{X,Y} \in \uwd\yx$ is the name change corresponding to the bijection $X \to Y \in \Fins$.
\end{itemize}
Then
\begin{equation}
\label{uwd-move-a2}
\tau_{X,Y} \comp \omega_x = \omega_y \in \uwd\ynothing.
\end{equation}
\end{proposition}

The next relation says that name changes inside a $2$-cell (Def. \ref{def:uwd-gen-2cell}) can be rewritten as a name change of a $2$-cell.

\begin{proposition}
\label{prop:uwd-move-a3}
Suppose:
\begin{itemize}
\item
$f_1 : X_1 \to Y_1$ and $f_2 : X_2 \to Y_2 \in \Fins$ are bijections.
\item
$f_1 \amalg f_2 : X_1 \amalg X_2 \to Y_1 \amalg Y_2 \in \Fins$ is their coproduct.
\item
$\tau_{f_1} \in \uwd\yonexone$, $\tau_{f_2} \in \uwd\ytwoxtwo$, and $\tau_{f_1 \amalg f_2} \in \uwd\yoneplusytwoxoneplusxtwo$ are the corresponding name changes.
\item
$\theta_{(X_1,X_2)} \in \uwd\xoneplusxtwoxonextwo$ and $\theta_{(Y_1,Y_2)} \in \uwd\yoneplusytwoyoneytwo$ are $2$-cells.
\end{itemize}
Then
\begin{equation}
\label{uwd-move-a3}
\Bigl(\theta_{(Y_1,Y_2)} \compone \tau_{f_1}\Bigr) \comptwo \tau_{f_2}
= \tau_{f_1 \amalg f_2} \comp \theta_{(X_1,X_2)}  \in \uwd\yoneplusytwoxonextwo.
\end{equation}
\end{proposition}

The next relation says that a name change inside a loop (Def. \ref{def:uwd-gen-loop}) can be rewritten as a name change of a loop.

\begin{proposition}
\label{prop:uwd-move-a4}
Suppose:
\begin{itemize}
\item
$X \in \Fins$, and $x_+,x_- \in X$ are two distinct elements with $v(x_+) = v(x_-) \in S$.
\item
$f : X \to Y \in \Fins$ is a bijection with $y_+ = f(x_+)$ and $y_- = f(x_-)$.
\item
$\xminusxpm \in \Fins$ and $\yminusypm \in \Fins$ are obtained from $X$ and $Y$ by removing the indicated elements.
\item
$f' : \xminusxpm \to \yminusypm$ is the corresponding bijection.
\item
$\tau_f \in \uwd\yx$ and $\tau_{f'} \in \uwd\yminusypmxminusxpm$ are name changes.
\item
$\lambda_{(X,\xpm)} \in \uwd\xminusxpmx$ and $\lambda_{(Y,\ypm)} \in \uwd\yminusypmy$ are loops.
\end{itemize}
Then
\begin{equation}
\label{uwd-move-a4}
\lambda_{(Y,\ypm)} \comp \tau_f = \tau_{f'} \comp \lambda_{(X,\xpm)} \in \uwd\yminusypmx.
\end{equation}
\end{proposition}

The next relation says that a name change inside a split (Def. \ref{def:uwd-gen-split}) can be rewritten as a name change of a split.

\begin{proposition}
\label{prop:uwd-move-a5}
Suppose:
\begin{itemize}
\item
$X \in \Fins$, and $x_1, x_2$ are two distinct elements in $X$ with $v(x_1) = v(x_2) \in S$.
\item
$X' \in \Fins$, and $x \in X'$ such that $v(x) = v(x_1)$ and that $X' \setminus \{x\} = X \setminus \{x_1,x_2\}$.
\item
$f : X \to Y \in \Fins$ is a bijection with $y_1 = f(x_1)$ and $y_2 = f(x_2)$.
\item
$Y' \in \Fins$, and $y \in Y'$ such that $v(y) = v(y_1) \in S$ and that $Y' \setminus \{y\} = Y \setminus \{y_1,y_2\}$.
\item
$f' : X' \to Y' \in \Fins$ is a bijection such that $f'(x) = y$ and that its restriction to $X' \setminus \{x\} = X \setminus \{x_1,x_2\}$ is equal to that of $f$.
\item
$\tau_f \in \uwd\yx$ and $\tau_{f'} \in \uwd\yprimexprime$ are name changes.
\item
$\sigma^{(X,x_1,x_2)} \in \uwd\xxprime$ and $\sigma^{(Y,y_1,y_2)} \in \uwd\yyprime$ are splits.
\end{itemize}
Then
\begin{equation}
\label{uwd-move-a5}
\sigma^{(Y,y_1,y_2)} \comp \tau_{f'} = \tau_f \comp \sigma^{(X,x_1,x_2)} \in \uwd\yxprime.
\end{equation}
\end{proposition}

The following two relations involve $1$-output wires in somewhat subtle ways.  The next relation says that the undirected wiring diagram
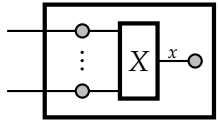
\begin{center}
\begin{tikzpicture}[scale=.5]
\draw [ultra thick] (-1,0) rectangle (3.5,3);
\draw [ultra thick] (1,.5) rectangle (2,2.5);
\node at (1.5,1.5) {$X$};
\cable{(0,.7)}
\node at (0,1.5) {{$\vdots$}};
\cable{(0,2.3)}
\cable{(3,1.5)}
\draw [thick] (1,.7) -- (.175,.7);
\draw [thick] (-.175,.7) -- (-2,.7);
\draw [thick] (1,2.3) -- (.175,2.3);
\draw [thick] (-.175,2.3) -- (-2,2.3);
\draw [thick] (2,1.5) -- (2.825,1.5);
\node at (2.4,1.7) {\tiny{$x$}};
\end{tikzpicture}
\end{center}
with a $(1,0)$-cable, all other cables being $(1,1)$-cables, can be obtained from the generators as either one of the following two (iterated) operadic compositions.
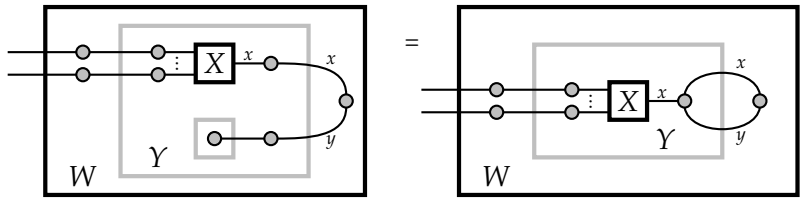
\begin{center}
\begin{tikzpicture}[scale=.5]
\draw [ultra thick] (1,3) rectangle (2,4);
\node at (1.5,3.5) {$X$};
\node at (.5,3.5) {\tiny{$\vdots$}};
\draw [ultra thick,lightgray] (1,1) rectangle (2,2);
\cable{(1.5,1.5)}
\cable{(0,3.2)} \cable{(0,3.8)}
\cable{(3,3.5)} \cable{(3,1.5)}
\draw [ultra thick, lightgray] (-1,.5) rectangle (4,4.5);
\node at (0,1) {$Y$};
\cable{(-2,3.2)} \cable{(-2,3.8)}
\cable{(5,2.5)}
\draw [thick] (1,3.8) -- (.175,3.8);
\draw [thick] (-.175,3.8) -- (-1.825,3.8);
\draw [thick] (-2.175,3.8) -- (-4,3.8);
\draw [thick] (1,3.2) -- (.175,3.2);
\draw [thick] (-.175,3.2) -- (-1.825,3.2);
\draw [thick] (-2.175,3.2) -- (-4,3.2);
\draw [thick] (2,3.5) -- (2.825,3.5);
\node at (2.4,3.7) {\tiny{$x$}};
\draw [thick] (3.175,3.5) to [out=0,in=90] (5,2.675);
\node at (4.6,3.6) {\tiny{$x$}};
\draw [thick] (5,2.325) to [out=-90,in=0] (3.175,1.5);
\node at (4.6,1.4) {\tiny{$y$}};
\draw [thick] (2.825,1.5) to (1.675,1.5);
\draw [ultra thick] (-3,0) rectangle (5.5,5);
\node at (-2,.5) {$W$};
\node at (6.75,4) {$=$};
\draw [ultra thick] (8,0) rectangle (17,5);
\node at (9,.5) {$W$};
\draw [ultra thick, lightgray] (10,1) rectangle (15,4);
\node at (13.5,1.5) {$Y$};
\draw [ultra thick] (12,2) rectangle (13,3);
\node at (12.5,2.5) {$X$};
\node at (11.5,2.5) {\tiny{$\vdots$}};
\cable{(11,2.2)} \cable{(11,2.8)} \cable{(14,2.5)}
\cable{(9,2.2)} \cable{(9,2.8)} \cable{(16,2.5)}
\draw [thick] (12,2.8) -- (11.175,2.8);
\draw [thick] (10.825,2.8) -- (9.175,2.8);
\draw [thick] (8.825,2.8) -- (7,2.8);
\draw [thick] (12,2.2) -- (11.175,2.2);
\draw [thick] (10.825,2.2) -- (9.175,2.2);
\draw [thick] (8.825,2.2) -- (7,2.2);
\draw [thick] (13,2.5) -- (13.825,2.5);
\node at (13.4,2.7) {\tiny{$x$}};
\draw [thick] (14,2.675) to [out=80,in=100] (16,2.675);
\node at (15.5,3.5) {\tiny{$x$}};
\draw [thick] (14,2.325) to [out=-80,in=-100] (16,2.325);
\node at (15.5,1.5) {\tiny{$y$}};
\end{tikzpicture}
\end{center}
As in Section \ref{sec:uwd-operad-structure}, the gray boxes here indicate an operadic composition.
\begin{itemize}
\item
On the left, a $1$-output wire $\omega_{y}$ is substituted into a $2$-cell $\theta_{(X,y)}$, which is then substituted into a loop $\lambda_{(Y,x,y)}$.
\item
On the right, a split $\sigma^{(Y,x,y)}$ is substituted into a loop $\lambda_{(Y,x,y)}$.
\end{itemize}

\begin{proposition}
\label{prop:uwd-move-b1}
Suppose:
\begin{itemize}
\item
$Y \in \Fins$, and $x,y \in Y$ are distinct elements with $v(x) = v(y) \in S$.
\item
$X = Y \setminus y \in \Fins$ is obtained from $Y$ by removing $y$.
\item
$\omega_{y} \in \uwd\smallynothing$ is the $1$-output wire for $y$.
\item
$\theta_{(X,y)} \in \uwd\yxy$ is a $2$-cell.
\item
$\lambda_{(Y,x,y)} \in \uwd\wy$ is a loop, where $W = Y \setminus \{x,y\} = X \setminus x$.
\item
$\sigma^{(Y,x,y)} \in \uwd\yx$ is a split.
\end{itemize}
Then
\begin{equation}
\label{uwd-move-b1}
\lambda_{(Y,x,y)}  \comp \Bigl(\theta_{(X,y)} \comptwo \omega_{y} \Bigr) =
\lambda_{(Y,x,y)}  \comp \sigma^{(Y,x,y)} \in \uwd\wx.
\end{equation}
\end{proposition}

The next relation says that the $X$-colored unit (Def. \ref{def:uwd-units}) can be obtained from the generators as the following iterated operadic composition.

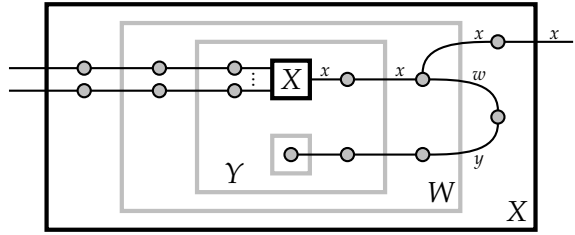
\begin{center}
\begin{tikzpicture}[scale=.5]
\draw [ultra thick] (-5,-.5) rectangle (8,5.5);
\node at (7.5,0) {$X$};
\draw [ultra thick, lightgray] (-3,0) rectangle (6,5);
\node at (5.5,.5) {$W$};
\draw [ultra thick, lightgray] (-1,.5) rectangle (4,4.5);
\node at (0,1) {$Y$};
\draw [ultra thick] (1,3) rectangle (2,4);
\node at (1.5,3.5) {$X$};
\node at (.5,3.5) {\tiny{$\vdots$}};
\draw [ultra thick,lightgray] (1,1) rectangle (2,2);
\cable{(1.5,1.5)}
\cable{(0,3.2)} \cable{(0,3.8)}
\cable{(3,3.5)} \cable{(3,1.5)}
\cable{(-2,3.2)} \cable{(-2,3.8)}
\cable{(5,3.5)} \cable{(5,1.5)}
\cable{(-4,3.2)} \cable{(-4,3.8)}
\cable{(7,4.5)} \cable{(7,2.5)}
\draw [thick] (1,3.8) -- (.175,3.8);
\draw [thick] (-.175,3.8) -- (-1.825,3.8);
\draw [thick] (-2.175,3.8) -- (-3.825,3.8);
\draw [thick] (-4.175,3.8) -- (-6,3.8);
\draw [thick] (1,3.2) -- (.175,3.2);
\draw [thick] (-.175,3.2) -- (-1.825,3.2);
\draw [thick] (-2.175,3.2) -- (-3.825,3.2);
\draw [thick] (-4.175,3.2) -- (-6,3.2);
\draw [thick] (2,3.5) -- (2.825,3.5);
\node at (2.4,3.7) {\tiny{$x$}};
\draw [thick] (3.175,3.5) -- (4.825,3.5);
\node at (4.4,3.7) {\tiny{$x$}};
\draw [thick] (5.175,3.5) to [out=0,in=90] (7,2.675);
\node at (6.5,3.6) {\tiny{$w$}};
\draw [thick] (7,2.325) to [out=-90,in=0] (5.175,1.5);
\node at (6.5,1.3) {\tiny{$y$}};
\draw [thick] (4.825,1.5) -- (3.175,1.5);
\draw [thick] (2.825,1.5) -- (1.675,1.5);
\draw [thick] (5,3.675) to [out=90,in=180] (6.825,4.5);
\node at (6.5,4.7) {\tiny{$x$}};
\draw [thick] (7.175,4.5) -- (9,4.5);
\node at (8.5,4.7) {\tiny{$x$}};
\end{tikzpicture}
\end{center}
More precisely, it says that the $X$-colored unit $\tensorunit_X$ can be obtained by substituting a $1$-output wire $\omega_y$ into a $2$-cell $\theta_{(X,y)}$, then into a split $\sigma^{(W,x,w)}$, and then into a loop $\lambda_{(W,w,y)}$.

\begin{proposition}
\label{prop:uwd-move-b2}
Suppose:
\begin{itemize}
\item
$W \in \Fins$, and $w,x,y$ are distinct elements in $W$ with $v(w) = v(x) = v(y) \in S$.
\item
$Y = W \setminus w \in \Fins$ is obtained from $W$ by removing $w$.
\item
$X = Y \setminus y \in \Fins$ is obtained from $Y$ by removing $y$.
\item
$\omega_{y} \in \uwd\smallynothing$ is the $1$-output wire for $y$.
\item
$\theta_{(X,y)} \in \uwd\yxy$ is a $2$-cell.
\item
$\sigma^{(W,x,w)} \in \uwd\wy$ is a split.
\item
$\lambda_{(W,w,y)} \in \uwd\xw$ is a loop.
\end{itemize}
Then
\begin{equation}
\label{uwd-move-b2}
\lambda_{(W,w,y)} \comp \sigma^{(W,x,w)} \comp 
\Bigl(\theta_{(X,y)} \comptwo \omega_{y} \Bigr) = \tensorunit_X \in \uwd\xx,
\end{equation}
in which $\tensorunit_X$ is the $X$-colored unit (Def. \ref{def:uwd-units}).
\end{proposition}

The next five relations are about $2$-cells (Def. \ref{def:uwd-gen-2cell}).  The following relation is the unity property of $2$-cells.

\begin{proposition}
\label{prop:uwd-move-c1}
Suppose:
\begin{itemize}
\item
$\theta_{(X,\varnothing)} \in \uwd\xxempty$ is a $2$-cell.
\item
$\epsilon \in \uwd\emptynothing$ is the empty cell (Def. \ref{def:uwd-gen-empty}).
\end{itemize}
Then
\begin{equation}
\label{uwd-move-c1}
\theta_{(X,\varnothing)} \comptwo \epsilon = \tensorunit_X \in \uwd\xx.
\end{equation}
\end{proposition}

The next relation is the associativity property of $2$-cells.  It gives two different ways to construct the following undirected wiring diagram using two $2$-cells.
\begin{center}
\begin{tikzpicture}[scale=.5]
\draw [ultra thick] (0,0) rectangle (7,4);
\draw [ultra thick] (1,1) rectangle (2,2);
\node at (1.5,1.5) {$X$};
\node at (1.5,2.5) {\tiny{$\cdots$}};
\draw [thick] (1.2,2) -- (1.2,4.5);
\draw [thick] (1.8,2) -- (1.8,4.5);
\draw [ultra thick] (3,1) rectangle (4,2);
\node at (3.5,1.5) {$Y$};
\node at (3.5,2.5) {\tiny{$\cdots$}};
\draw [thick] (3.2,2) -- (3.2,4.5);
\draw [thick] (3.8,2) -- (3.8,4.5);
\draw [ultra thick] (5,1) rectangle (6,2);
\node at (5.5,1.5) {$Z$};
\node at (5.5,2.5) {\tiny{$\cdots$}};
\draw [thick] (5.2,2) -- (5.2,4.5);
\draw [thick] (5.8,2) -- (5.8,4.5);
\cable{(1.2,3)} \cable{(1.8,3)}
\cable{(3.2,3)} \cable{(3.8,3)}
\cable{(5.2,3)} \cable{(5.8,3)}
\end{tikzpicture}
\end{center}

\begin{proposition}
\label{prop:uwd-move-c2}
Suppose:\index{associativity of undirected 2-cell@associativity of undirected $2$-cells}
\begin{itemize}
\item
$\theta_{(X \amalg Y, Z)} \in \uwd\xplusyz$ and $\theta_{(X,Y)} \in \uwd\xplusyxy$ are $2$-cells.
\item
$\theta_{(X, Y \amalg Z)} \in \uwd\xyplusz$ and $\theta_{(Y,Z)} \in \uwd\ypluszyz$ are $2$-cells.  
\end{itemize}
Then
\begin{equation}
\label{uwd-move-c2}
\theta_{(X \amalg Y, Z)} \compone \theta_{(X,Y)}
= \theta_{(X, Y \amalg Z)} \comptwo \theta_{(Y,Z)}
\in \uwd\xplusyplusz.
\end{equation}
\end{proposition}

The next relation is the commutativity property of $2$-cells.  It uses the equivariant structure \eqref{def:uwd-equivariance} in $\uwd$.

\begin{proposition}
\label{prop:uwd-move-c3}
Suppose:\index{commutativity of undirected 2-cell@commutativity of undirected $2$-cells}
\begin{itemize}
\item
$\theta_{X,Y} \in \uwd\xplusyxy$ is a $2$-cell.
\item
$(1~2) \in \Sigma_2$ is the non-trivial permutation.
\end{itemize}
Then
\begin{equation}
\label{uwd-move-c3}
\theta_{(X,Y)}(1~2) = \theta_{(Y,X)} \in \uwd\yplusxyx.
\end{equation}
\end{proposition}

The next relation is the commutativity property between a $2$-cell and a loop.  It gives two different ways to construct the following undirected wiring diagram.
\begin{center}
\begin{tikzpicture}[scale=.6]
\draw [ultra thick] (-1,0) rectangle (6,4);
\node at (.2,.4) {$X \amalg Y'$};
\draw [ultra thick] (1,1) rectangle (2,2);
\node at (1.5,1.5) {$X$};
\draw [thick] (1.2,2) -- (1.2,4.5);
\draw [thick] (1.8,2) -- (1.8,4.5);
\node at (1.5,2.5) {\tiny{$\cdots$}};
\cable{(1.2,3)} \cable{(1.8,3)}
\draw [ultra thick] (3,.5) rectangle (4,2);
\node at (3.5,1.25) {$Y$};
\node at (3.5,2.5) {\tiny{$\cdots$}};
\cable{(3.2,3)} \cable{(3.8,3)}
\draw [thick] (3.2,2) -- (3.2,2.825);
\draw [thick] (3.2,3.175) -- (3.2,4.5);
\draw [thick] (3.8,2) -- (3.8,2.825);
\draw [thick] (3.8,3.175) -- (3.8,4.5);
\cable{(5,1.25)} 
\draw [thick] (4,1.9) to [out=0,in=90] (5,1.425);
\node at (4.4,2.1) {\tiny{$y_+$}};
\draw [thick] (4,.6) to [out=0,in=-90] (5,1.075);
\node at (4.4,.3) {\tiny{$y_-$}};
\end{tikzpicture}
\end{center}

\begin{proposition}
\label{prop:uwd-move-c4}
Suppose:
\begin{itemize}
\item
$Y \in \Fins$, and $y_+, y_-$ are distinct elements in $Y$ with $v(y_+) = v(y_-) \in S$.
\item
$Y' = Y \setminus \ypm \in \Fins$ is obtained from $Y$ by removing $y_+$ and $y_-$.
\item
$\lambda_{(Y,\ypm)} \in \uwd\yprimey$ is a loop.
\item
$\theta_{(X,Y)} \in \uwd\xplusyxy$ and $\theta_{(X,Y')} \in \uwd\xplusyprimexyprime$ are $2$-cells for some $X \in \Fins$.
\item
$\lambda_{(X \amalg Y, \ypm)} \in \uwd\xplusyprimexplusy$ is a loop.
\end{itemize}
Then
\begin{equation}
\label{uwd-move-c4}
\theta_{(X,Y')} \comptwo \lambda_{(Y,\ypm)} 
= \lambda_{(X \amalg Y, \ypm)}  \comp \theta_{(X,Y)} \in \uwd\xplusyprimexy.
\end{equation}
\end{proposition}

The next relation is the commutativity between a $2$-cell and a split.  It gives two different ways to construct the following undirected wiring diagram.
\begin{center}
\begin{tikzpicture}[scale=.6]
\draw [ultra thick] (-1,0) rectangle (6,4);
\node at (.2,.4) {$X \amalg Y$};
\draw [ultra thick] (1,1) rectangle (2,2);
\node at (1.5,1.5) {$X$};
\node at (1.5,2.5) {\tiny{$\cdots$}};
\draw [thick] (1.2,2) -- (1.2,4.5);
\draw [thick] (1.8,2) -- (1.8,4.5);
\cable{(1.2,3)} \cable{(1.8,3)}
\draw [ultra thick] (3,1) rectangle (4,2);
\node at (3.5,1.5) {$Y'$};
\node at (3.5,2.5) {\tiny{$\cdots$}};
\cable{(3.2,3)} \cable{(3.8,3)}
\draw [thick] (3.2,2) -- (3.2,2.825);
\draw [thick] (3.2,3.175) -- (3.2,4.5);
\draw [thick] (3.8,2) -- (3.8,2.825);
\draw [thick] (3.8,3.175) -- (3.8,4.5);
\cable{(5,1.5)} 
\draw [thick] (4,1.5) -- (4.825,1.5);
\node at (4.4,1.7) {\tiny{$y$}};
\draw [thick] (5,1.675) to [out=90,in=180] (7,2.5);
\node at (6.4,2.7) {\tiny{$y_1$}};
\draw [thick] (5,1.325) to [out=-90,in=180] (7,.5);
\node at (6.4,.2) {\tiny{$y_2$}};
\end{tikzpicture}
\end{center}

\begin{proposition}
\label{prop:uwd-move-c5}
Suppose:
\begin{itemize}
\item
$Y \in \Fins$, and $y_1, y_2$ are distinct elements in $Y$ with $v(y_1) = v(y_2) \in S$.
\item
$Y' \in \Fins$, and $y \in Y'$ such that $v(y) = v(y_1)$ and that $Y' \setminus \{y\} = Y \setminus \{y_1,y_2\}$.
\item
$\sigma^{(Y,y_1,y_2)} \in \uwd\yyprime$ is a split.
\item
$\theta_{(X,Y)} \in \uwd\xplusyxy$ and $\theta_{(X,Y')} \in \uwd\xplusyprimexyprime$ are $2$-cells for some $X \in \Fins$.
\item
$\sigma^{(X \amalg Y, y_1, y_2)} \in \uwd\xplusyxplusyprime$ is a split.
\end{itemize}
Then
\begin{equation}
\label{uwd-move-c5}
\theta_{(X,Y)} \comptwo \sigma^{(Y,y_1,y_2)} 
= \sigma^{(X \amalg Y, y_1, y_2)}  \comp \theta_{(X,Y')} \in \uwd\xplusychoosexyprime.
\end{equation}
\end{proposition}

The following four relations are about splits.  The next relation is the commutativity property of splits.  It gives two different ways to construct the following undirected wiring diagram using two splits.

\begin{center}
\begin{tikzpicture}[scale=.6]
\draw [ultra thick] (1,0) rectangle (6,4);
\node at (1.5,3.5) {$X$};
\draw [ultra thick] (3,1) rectangle (4,2);
\node at (3.5,1.5) {$X'$};
\node at (3.5,2.5) {\tiny{$\cdots$}};
\cable{(3.2,3)} \cable{(3.8,3)}
\draw [thick] (3.2,2) -- (3.2,2.825);
\draw [thick] (3.2,3.175) -- (3.2,4.5);
\draw [thick] (3.8,2) -- (3.8,2.825);
\draw [thick] (3.8,3.175) -- (3.8,4.5);
\cable{(2,1.5)}
\draw [thick] (3,1.5) -- (2.175,1.5);
\node at (2.6,1.7) {\tiny{$z$}};
\draw [thick] (2,1.675) to [out=90,in=0] (0,2.5);
\node at (.6,2.7) {\tiny{$z_1$}};
\draw [thick] (2,1.325) to [out=-90,in=0] (0,.5);
\node at (.6,.2) {\tiny{$z_2$}};
\cable{(5,1.5)} 
\draw [thick] (4,1.5) -- (4.825,1.5);
\node at (4.4,1.7) {\tiny{$y$}};
\draw [thick] (5,1.675) to [out=90,in=180] (7,2.5);
\node at (6.4,2.7) {\tiny{$y_1$}};
\draw [thick] (5,1.325) to [out=-90,in=180] (7,.5);
\node at (6.4,.2) {\tiny{$y_2$}};
\end{tikzpicture}
\end{center}

\begin{proposition}
\label{prop:uwd-move-d1}
Suppose:\index{commutativity of splits}
\begin{itemize}
\item
$X \in \Fins$, and $y_1,y_2,z_1,z_2$ are distinct elements such that $v(y_1) = v(y_2)$ and $v(z_1) = v(z_2) \in S$.
\item
$X' \in \Fins$, and $y$ and $z$ are distinct elements in $X'$ such that
\begin{itemize}
\item
$v(y) = v(y_1)$ and $v(z) = v(z_1)$;
\item
$X' \setminus \{y,z\} = X \setminus \{y_1,y_2,z_1,z_2\}$.
\end{itemize}
\item
$Y = \bigl[X' \setminus \{y\}\bigr] \amalg \{y_1,y_2\}$ and 
$Z = \bigl[X' \setminus \{z\}\bigr] \amalg \{z_1,z_2\} \in\Fins$
\item
$\sigma^{(X,y_1,y_2)} \in \uwd\xz$ and $\sigma^{(Z,z_1,z_2)} \in \uwd\zxprime$ are splits.
\item
$\sigma^{(X,z_1,z_2)} \in \uwd\xchoosey$ and $\sigma^{(Y,y_1,y_2)} \in \uwd\yxprime$ are splits.
\end{itemize}
Then
\begin{equation}
\label{uwd-move-d1}
\sigma^{(X,y_1,y_2)} \comp \sigma^{(Z,z_1,z_2)}
=  \sigma^{(X,z_1,z_2)} \comp \sigma^{(Y,y_1,y_2)} \in \uwd\xxprime.
\end{equation}
\end{proposition}

The next relation is the associativity property of splits.  It gives two different ways to construct the following undirected wiring diagram using two splits.
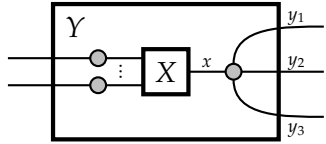
\begin{center}
\begin{tikzpicture}[scale=.6]
\draw [ultra thick] (-1,0) rectangle (4,3);
\node at (-.5,2.5) {$Y$};
\draw [ultra thick] (1,1) rectangle (2,2);
\node at (1.5,1.5) {$X$};
\node at (.5,1.5) {\tiny{$\vdots$}};
\cable{(0,1.2)} \cable{(0,1.8)}
\draw [thick] (1,1.8) -- (.175,1.8);
\draw [thick] (-.175,1.8) -- (-2,1.8);
\draw [thick] (1,1.2) -- (.175,1.2);
\draw [thick] (-.175,1.2) -- (-2,1.2);
\cable{(3,1.5)} 
\draw [thick] (2,1.5) -- (2.825,1.5);
\node at (2.4,1.7) {\tiny{$x$}};
\draw [thick] (3,1.675) to [out=90,in=180] (5,2.5);
\node at (4.4,2.7) {\tiny{$y_1$}};
\draw [thick] (3.175,1.5) -- (5,1.5);
\node at (4.4,1.7) {\tiny{$y_2$}};
\draw [thick] (3,1.325) to [out=-90,in=180] (5,.5);
\node at (4.4,.2) {\tiny{$y_3$}};
\end{tikzpicture}
\end{center}

\begin{proposition}
\label{prop:uwd-move-d2}
Suppose:\index{associativity of splits}
\begin{itemize}
\item
$Y \in \Fins$, and $y_1,y_2,y_3$ are distinct elements in $Y$ with $v(y_1) = v(y_2) = v(y_3) \in S$.
\item
$X \in \Fins$, and $x \in X$ such that 
\begin{itemize}
\item
$v(x) = v(y_1)$;
\item
$X \setminus \{x\} = Y \setminus \{y_1,y_2,y_3\}$.
\end{itemize}
\item
$Y_1 = Y/(y_1 = y_2) \in \Fins$ is the quotient of $Y$ with $y_1$ and $y_2$ identified, called $y_{12} \in Y_1$.
\item
$Y_2 = Y/(y_2 = y_3) \in \Fins$ is the quotient of $Y$ with $y_2$ and $y_3$ identified, called $y_{23} \in Y_2$.
\item
$\sigma^{(Y,y_1,y_2)} \in \uwd\yyone$ and $\sigma^{(Y_1,y_{12},y_3)} \in \uwd\yonex$ are splits.
\item
$\sigma^{(Y,y_2,y_3)} \in \uwd\yytwo$ and $\sigma^{(Y_2,y_1,y_{23})} \in \uwd\ytwox$ are splits.
\end{itemize}
Then
\begin{equation}
\label{uwd-move-d2}
\sigma^{(Y,y_1,y_2)} \comp \sigma^{(Y_1,y_{12},y_3)}
=  \sigma^{(Y,y_2,y_3)} \comp \sigma^{(Y_2,y_1,y_{23})} \in \uwd\yx.
\end{equation}
\end{proposition}

The next relation is the commutativity property between a split and a loop.  It gives two different ways to construct the following undirected wiring diagram using a split and a loop.
\begin{center}
\begin{tikzpicture}[scale=.6]
\draw [ultra thick] (1,0) rectangle (6,4);
\node at (1.5,3.5) {$Y$};
\draw [ultra thick] (3,.5) rectangle (4,2);
\node at (3.5,1.25) {$X$};
\node at (3.5,2.5) {\tiny{$\cdots$}};
\cable{(3.2,3)} \cable{(3.8,3)}
\draw [thick] (3.2,2) -- (3.2,2.825);
\draw [thick] (3.2,3.175) -- (3.2,4.5);
\draw [thick] (3.8,2) -- (3.8,2.825);
\draw [thick] (3.8,3.175) -- (3.8,4.5);
\cable{(5,1.25)} 
\draw [thick] (4,1.9) to [out=0,in=90] (5,1.425);
\node at (4.4,2.1) {\tiny{$x_+$}};
\draw [thick] (4,.6) to [out=0,in=-90] (5,1.075);
\node at (4.4,.3) {\tiny{$x_-$}};
\cable{(2,1.25)}
\draw [thick] (3,1.25) -- (2.175,1.25);
\node at (2.6,1.5) {\tiny{$x$}};
\draw [thick] (2,1.425) to [out=90,in=0] (0,2.25);
\node at (.5,2.5) {\tiny{$y_1$}};
\draw [thick] (2,1.075) to [out=-90,in=0] (0,.25);
\node at (.5,.5) {\tiny{$y_2$}};
\end{tikzpicture}
\end{center}

\begin{proposition}
\label{prop:uwd-move-d3}
Suppose:
\begin{itemize}
\item
$Y \in \Fins$, and $y_1$ and $y_2$ are distinct elements in $Y$ with $v(y_1)=v(y_2) \in S$.
\item
$X \in \Fins$, and $x,x_+,x_-$ are distinct elements in $X$ such that
\begin{itemize}
\item
$v(x) = v(y_1)$;
\item
$v(x_+) = v(x_-)$;
\item
$X \setminus \{x,x_+,x_-\} = Y \setminus \{y_1,y_2\}$.
\end{itemize}
\item
$Y' = \bigl[X \setminus \{x\}\bigr] \amalg \{y_1,y_2\}$ and 
$X' = X \setminus \{x_+,x_-\}$.
\item
$\sigma^{(Y',y_1,y_2)} \in \uwd\yprimex$ and $\sigma^{(Y,y_1,y_2)} \in \uwd\yxprime$ are splits.
\item
$\lambda_{(Y',\xsubpm)} \in \uwd\yyprime$ and $\lambda_{(X,\xsubpm)} \in \uwd\xprimex$ are loops.
\end{itemize}
Then
\begin{equation}
\label{uwd-move-d3}
\lambda_{(Y',\xsubpm)} \comp \sigma^{(Y',y_1,y_2)} 
=  \sigma^{(Y,y_1,y_2)} \comp \lambda_{(X,\xsubpm)} \in \uwd\yx.
\end{equation}
\end{proposition}

The next relation says that the undirected wiring diagram
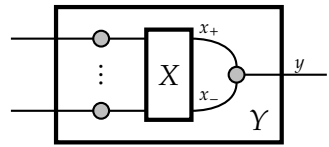
\begin{center}
\begin{tikzpicture}[scale=.6]
\draw [ultra thick] (-1,0) rectangle (4,3);
\node at (3.5,.5) {$Y$};
\draw [ultra thick] (1,.5) rectangle (2,2.5);
\node at (1.5,1.5) {$X$};
\cable{(0,.7)} \node at (0,1.5) {{$\vdots$}}; \cable{(0,2.3)}
\cable{(3,1.5)}
\draw [thick] (1,.7) -- (.175,.7);
\draw [thick] (-.175,.7) -- (-2,.7);
\draw [thick] (1,2.3) -- (.175,2.3);
\draw [thick] (-.175,2.3) -- (-2,2.3);
\draw [thick] (2,2.3) to [out=0,in=90] (3,1.675);
\node at (2.4,2.5) {\tiny{$x_+$}};
\draw [thick] (2,.7) to [out=0,in=-90] (3,1.325);
\node at (2.4,1) {\tiny{$x_-$}};
\draw [thick] (3.175,1.5) -- (5,1.5);
\node at (4.4,1.7) {\tiny{$y$}};
\end{tikzpicture}
\end{center}
can be obtained by substituting a split inside a loop as in the picture
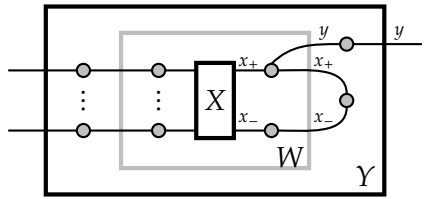
\begin{center}
\begin{tikzpicture}[scale=.5]
\draw [ultra thick] (-3,-1) rectangle (6,4);
\node at (5.5,-.5) {$Y$};
\draw [ultra thick, lightgray] (-1,-.3) rectangle (4,3.3);
\node at (-2,1.5) {{$\vdots$}};
\node at (3.5,0) {$W$};
\draw [ultra thick] (1,.5) rectangle (2,2.5);
\node at (1.5,1.5) {$X$};
\cable{(0,.7)} \cable{(-2,.7)}
\draw [thick] (1,.7) -- (.175,.7);
\draw [thick] (-.175,.7) -- (-1.825,.7);
\draw [thick] (-2.175,.7) -- (-4,.7);
\node at (0,1.5) {{$\vdots$}}; 
\cable{(0,2.3)} \cable{(-2,2.3)}
\draw [thick] (1,2.3) -- (.175,2.3);
\draw [thick] (-.175,2.3) -- (-1.825,2.3);
\draw [thick] (-2.175,2.3) -- (-4,2.3);
\cable {(3,2.3)} \cable{(5,1.5)} \cable{(3,.7)}  
\draw [thick] (2,2.3) -- (2.825,2.3); \node at (2.4,2.5) {\tiny{$x_+$}};
\draw [thick] (3.175,2.3) to [out=0,in=90] (5,1.675);
\node at (4.4,2.5) {\tiny{$x_+$}};
\draw [thick] (5,1.325) to [out=-90,in=0] (3.175,.7);
\node at (4.4,1) {\tiny{$x_-$}};
\draw [thick] (2.825,.7) -- (2,.7); \node at (2.4,1) {\tiny{$x_-$}};
\draw [thick] (3,2.475) to [out=45,in=180] (4.8,3);
\node at (4.4,3.3) {\tiny{$y$}};
\cable{(5,3)} \draw [thick] (5.175,3) -- (7,3); \node at (6.4,3.3) {\tiny{$y$}};
\end{tikzpicture}
\end{center}
or in the counterpart in which $x_+$ and $x_-$ are switched.  In \eqref{uwd-move-d4} below, this picture corresponds to the left side, and its counterpart corresponds to the right side.

\begin{proposition}
\label{prop:uwd-move-d4}
Suppose:
\begin{itemize}
\item
$X \in \Fins$, and $x_+$ and $x_-$ are distinct elements in $X$ with $v(x_+) = v(x_-) \in S$.
\item
$Y \in \Fins$, and $y \in Y$ such that $v(y) = v(x_+)$ and that $X \setminus \xpm = Y \setminus y$.
\item
$W = X \amalg y = Y \amalg \xpm \in \Fins$.
\item
$\sigma^{(W,y, x_+)} \in \uwd\wx$ and $\sigma^{(W,y, x_-)} \in \uwd\wx$ are splits.
\item
$\lambda_{(W,\xpm)} \in \uwd\yw$ is a loop.
\end{itemize}
Then
\begin{equation}
\label{uwd-move-d4}
\lambda_{(W,\xpm)} \comp \sigma^{(W,y, x_+)} 
= \lambda_{(W,\xpm)} \comp \sigma^{(W,y, x_-)} \in \uwd\yx.
\end{equation}
\end{proposition}

The final relation is the commutativity property of loops.  It gives two different ways to construct the following undirected wiring diagram using two loops.
\begin{center}
\begin{tikzpicture}[scale=.5]
\draw [ultra thick] (-2,0) rectangle (4,4);
\node at (3.5,3.5) {$Y$};
\draw [ultra thick] (0,.5) rectangle (2,2.5);
\node at (1,1.5) {$X$};
\draw [thick] (.5,2.5) -- (.5,4.5); \cable{(.5,3.25)}
\node at (1,2.8) {$\cdots$};
\draw [thick] (1.5,2.5) -- (1.5,4.5); \cable{(1.5,3.25)}
\draw [thick] (0,2.3) to [out=180,in=90] (-1,1.675);
\node at (-.4,2.5) {\tiny{$x_1$}};
\cable{(-1,1.5)}
\draw [thick] (-1,1.325) to [out=-90,in=180] (0,.7);
\node at (-.4,1) {\tiny{$x_2$}};
\draw [thick] (2,2.3) to [out=0,in=90] (3,1.675);
\node at (2.4,2.5) {\tiny{$x_3$}};
\cable{(3,1.5)}
\draw [thick] (2,.7) to [out=0,in=-90] (3,1.325);
\node at (2.4,1) {\tiny{$x_4$}};
\end{tikzpicture}
\end{center}

\begin{proposition}
\label{prop:uwd-move-e1}
Suppose:\index{commutativity of loops}
\begin{itemize}
\item
$X \in \Fins$, and $x_1,x_2,x_3,x_4$ are distinct elements in $X$ with $v(x_1) = v(x_2)$ and $v(x_3) = v(x_4) \in S$.
\item
$W = X \setminus \{x_1,x_2\}$, $Z = X \setminus \{x_3,x_4\}$, and $Y = X \setminus \{x_1,x_2,x_3,x_4\} \in \Fins$.
\item
$\lambda_{(W,x_3,x_4)} \in \uwd\yw$ and $\lambda_{(X,x_1,x_2)} \in \uwd\wx$  are loops.
\item
$\lambda_{(Z,x_1,x_2)} \in \uwd\yz$ and $\lambda_{(X,x_3,x_4)} \in \uwd\zx$ are loops.
\end{itemize}
Then
\begin{equation}
\label{uwd-move-e1}
\lambda_{(W,x_3,x_4)} \comp \lambda_{(X,x_1,x_2)}
= \lambda_{(Z,x_1,x_2)} \comp \lambda_{(X,x_3,x_4)} \in \uwd\yx.  
\end{equation}
\end{proposition}

\begin{definition}
\label{def:elementary-relation-uwd}
The $17$ relations \eqref{uwd-move-a1}--\eqref{uwd-move-e1} are called  \emph{elementary relations in $\uwd$}.  If there is no danger of confusion, we will call them \emph{elementary relations}.\index{elementary relations in UWD@elementary relations in $\uwd$}
\end{definition}

\section{Wasted Cables}
\label{sec:wasted-cable}

The purpose of this section is to consider several examples of how the generators in the operad $\uwd$ can create wasted cables (Def. \ref{def:pre-uwd}).  Example \ref{ex:y-wasted-two-ways} provides an illustration of some of the elementary relations in $\uwd$.  The examples in this section provide a good warm-up exercise for the discussion in Chapter \ref{ch10-stratified-uwd} about stratified presentations and elementary equivalences.  

Recall from Remark \ref{def:generating-uwd} that none of the generators has a wasted cable.

\begin{example}
\label{ex:wasted-cable-by-generators}
In the context of Example \ref{ex:uwd-jointly-surjective} with $X = \{x_1,x_2\}$:
\begin{enumerate}
\item
$\varphi = \left\{\nicexy{X \ar[r] & \ast & \varnothing \ar[l]}\right\} \in \uwd\emptyx$ is the loop $\lambda_{(X,x_1,x_2)}$.
\item
$\psi = \left\{\nicexy{\varnothing \ar[r] & X & X \ar[l]_-{=}}\right\} \in \uwd\xempty$ is the iterated operadic composition
\[\psi = \left[\Bigl(\theta_{(\varnothing, X)} \comptwo \theta_{(x_1,x_2)}\Bigr) \comptwo \omega_{x_1}\right] \comptwo \omega_{x_2}\]
involving two $2$-cells and two $1$-output wires.  
\end{enumerate}
So the composition
\[\varphi \comp \psi = \Bigl(\nicexy{\varnothing \ar[r] & \ast & \varnothing \ar[l]}\Bigr)\in \uwd\emptyprof,\] 
which is depicted as
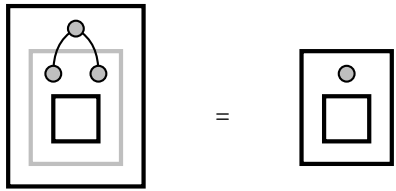
\begin{center}
\begin{tikzpicture}[scale=.6]
\draw [ultra thick] (-6,0) rectangle (-3,4);
\draw [ultra thick, lightgray] (-5.5,.5) rectangle (-3.5,3);
\draw [ultra thick] (-5,1) rectangle (-4,2);
\draw [thick] (-5,2.5) to [out=90,in=210] (-4.5,3.5);
\draw [thick] (-4.5,3.5) to [out=-30,in=90] (-4,2.5);
\cable{(-5,2.5)} \cable{(-4.5,3.5)} \cable{(-4,2.5)}
\node at (-1.25,1.5) {$=$};
\draw [ultra thick] (.5,.5) rectangle (2.5,3);
\draw [ultra thick] (1,1) rectangle (2,2);
\cable{(1.5,2.5)}
\end{tikzpicture}
\end{center}
and has one wasted cable, is the iterated operadic composition
\begin{equation}
\label{wasted-cable-simplex}
\lambda_{(X,x_1,x_2)} \comp \left[\Bigl(\left(\theta_{(\varnothing, X)} \comptwo \theta_{(x_1,x_2)}\right) \comptwo \omega_{x_1}\Bigr) \comptwo \omega_{x_2}\right]
\end{equation}
involving $5$ generators.  
\end{example}

\begin{example}
\label{ex:y-plus-wasted-cable}
As a variation of the previous example, consider any box $Y \in \Fins$ and the undirected wiring diagram with one wasted cable  
\[\zeta_Y = \Bigl(\nicexy@C+.4cm{Y \ar[r]^-{\mathrm{inclusion}} & Y \amalg \ast & Y \ar[l]_-{\mathrm{inclusion}}}\Bigr) \in \uwd\yy.\]
It is depicted as follows.
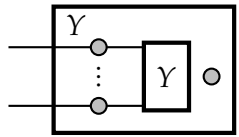
\begin{center}
\begin{tikzpicture}[scale=.6]
\draw [ultra thick] (-1,0) rectangle (3,2.8);
\node at (-.5,2.4) {$Y$};
\draw [ultra thick] (1,.5) rectangle (2,2);
\node at (1.5,1.25) {$Y$};
\draw [thick] (1,.6) -- (-2,.6); \cable{(0,.6)}
\node at (0,1.25) {{$\vdots$}};
\draw [thick] (1,1.9) -- (-2,1.9); \cable{(0,1.9)}
\cable{(2.5,1.25)}
\end{tikzpicture}
\end{center}
This undirected wiring diagram can be created by replacing the empty box $\varnothing$ in the $2$-cell $\theta_{(\varnothing, X)}$ by $Y$ and the loop $\lambda_{(X,x_1,x_2)}$ by the loop $\lambda_{(Y \amalg X,x_1,x_2)}$ in \eqref{wasted-cable-simplex} above.  The resulting operadic composition 
\begin{equation}
\label{wasted-cable-y-simplex}
\zeta_Y = 
\lambda_{(Y \amalg X,x_1,x_2)} \comp \left[\Bigl(\left(\theta_{(Y, X)} \comptwo \theta_{(x_1,x_2)}\right) \comptwo \omega_{x_1}\Bigr) \comptwo \omega_{x_2}\right] \in \uwd\yy
\end{equation}
involves $5$ generators: one loop, two $2$-cells, and two $1$-output wires.  It corresponds to the following picture.
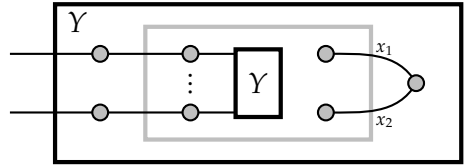
\begin{center}
\begin{tikzpicture}[scale=.6]
\draw [ultra thick] (-3,-.5) rectangle (6,3);
\node at (-2.5,2.6) {$Y$};
\draw [ultra thick, lightgray] (-1,0) rectangle (4,2.5);
\draw [ultra thick] (1,.5) rectangle (2,2);
\node at (1.5,1.25) {$Y$};
\draw [thick] (1,.6) -- (-4,.6); \cable{(0,.6)} \cable{(-2,.6)}
\node at (0,1.25) {{$\vdots$}};
\draw [thick] (1,1.9) -- (-4,1.9); \cable{(0,1.9)} \cable{(-2,1.9)}
\draw [thick] (3,1.9) to [out=0,in=120] (5,1.25);
\draw [thick] (5,1.25) to [out=-120,in=0] (3,.6);
\cable{(3,1.9)} \cable{(5,1.25)} \cable{(3,.6)}
\node at (4.3,2) {\tiny{$x_1$}};
\node at (4.3,.4) {\tiny{$x_2$}};
\end{tikzpicture}
\end{center}
The intermediate gray box is $Y \amalg X = Y \amalg \{x_1,x_2\}$.  Roughly speaking, the operadic composition \eqref{wasted-cable-y-simplex} says that a wasted cable can be created by applying a loop to two $1$-output wires.  Additional wasted cables can similarly be created using more $2$-cells, $1$-output wires, and loops.  
\end{example}

\begin{example}
\label{ex:y-plus-wasted-cable-b}
The undirected wiring diagram $\zeta_Y$ in Example \ref{ex:y-plus-wasted-cable} can also be created as in the following picture.
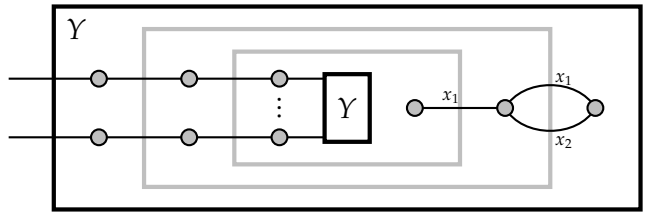
\begin{center}
\begin{tikzpicture}[scale=.6]
\draw [ultra thick] (-5,-1) rectangle (8,3.5);
\node at (-4.5,3) {$Y$};
\draw [ultra thick, lightgray] (-3,-.5) rectangle (6,3);
\draw [ultra thick, lightgray] (-1,0) rectangle (4,2.5);
\draw [ultra thick] (1,.5) rectangle (2,2);
\node at (1.5,1.25) {$Y$};
\draw [thick] (1,.6) -- (-6,.6); \cable{(0,.6)} \cable{(-2,.6)} \cable{(-4,.6)}
\node at (0,1.25) {{$\vdots$}};
\draw [thick] (1,1.9) -- (-6,1.9); \cable{(0,1.9)} \cable{(-2,1.9)} \cable{(-4,1.9)}
\draw [thick] (3,1.25) -- (5,1.25);
\draw [thick, bend left=60] (5,1.25) to (7,1.25);
\draw [thick, bend right=60] (5,1.25) to (7,1.25);
\cable{(3,1.25)} \cable{(5,1.25)} \cable{(7,1.25)}
\node at (3.8,1.5) {\tiny{$x_1$}};
\node at (6.3,1.9) {\tiny{$x_1$}};
\node at (6.3,.5) {\tiny{$x_2$}};
\end{tikzpicture}
\end{center}
The inner gray box is $Y \amalg x_1$, and the outer gray box is $Y \amalg X$. In terms of the generators, the above picture is realized as the operadic composition
\begin{equation}
\label{wasted-cable-y-simplex-b}
\zeta_Y = \Bigl[\lambda_{\left(Y \amalg X, x_1,x_2\right)} \comp \sigma^{\left(Y \amalg X, x_1,x_2\right)}\Bigr] \comp
\Bigl[\theta_{(Y,x_1)} \comptwo \omega_{x_1}\Bigr] \in \uwd\yy.
\end{equation}
It involves one loop, one split, one $2$-cell, and one $1$-output wire. Roughly speaking, the operadic composition \eqref{wasted-cable-y-simplex-b} says that a wasted cable can be created by applying a loop to a split that is attached to a $1$-output wire. 
\end{example}

\begin{example}
\label{ex:y-wasted-two-ways}
As an illustration of using the elementary relations in $\uwd$, recall the undirected wiring diagram $\zeta_Y \in \uwd\yy$ in Examples \ref{ex:y-plus-wasted-cable} and \ref{ex:y-plus-wasted-cable-b}.  It can be generated by the generators as either one of the two iterated operadic compositions \eqref{wasted-cable-y-simplex} and \eqref{wasted-cable-y-simplex-b}.  These two decompositions of $\zeta_Y$ are actually connected as follows.
\[\begin{split}
&\lambda_{(Y \amalg X,x_1,x_2)} \comp \left[\Bigl(\left(\theta_{(Y, X)} \comptwo \theta_{(x_1,x_2)}\right) \comptwo \omega_{x_1}\Bigr) \comptwo \omega_{x_2}\right]  \quad\text{by \eqref{wasted-cable-y-simplex}}\\
&= \lambda_{(Y \amalg X,x_1,x_2)} \comp \left[\Bigl(\left(\theta_{(Y \amalg x_1, x_2)} \compone \theta_{(Y,x_1)}\right) \comptwo \omega_{x_1}\Bigr) \comptwo \omega_{x_2}\right] \quad\text{by elem. rel. \eqref{uwd-move-c2}}\\
&= \lambda_{(Y \amalg X,x_1,x_2)} \comp \left[ \Bigl(\theta_{(Y \amalg x_1, x_2)} \compone \left(\theta_{(Y,x_1)} \comptwo \omega_{x_1}\right) \Bigr) \comptwo \omega_{x_2}\right] \quad\text{by vertical ass. \eqref{compi-associativity-two}}\\
&= \lambda_{(Y \amalg X,x_1,x_2)} \comp \left[ \Bigl(\theta_{(Y \amalg x_1, x_2)} \comptwo \omega_{x_2}\Bigr) \compone \Bigl(\theta_{(Y,x_1)} \comptwo \omega_{x_1}\Bigr)\right]\quad\text{by horizontal ass. \eqref{compi-associativity}}\\
&= \left[\lambda_{(Y \amalg X,x_1,x_2)} \comp \Bigl(\theta_{(Y \amalg x_1, x_2)} \comptwo \omega_{x_2}\Bigr)\right] \comp \Bigl[\theta_{(Y,x_1)} \comptwo \omega_{x_1}\Bigr] \quad\text{by vertical ass. \eqref{compi-associativity-two}}\\
&= \Bigl[\lambda_{(Y \amalg X,x_1,x_2)} \comp \sigma^{\left(Y \amalg X, x_1,x_2\right)}\Bigr] \comp \Bigl[\theta_{(Y,x_1)} \comptwo \omega_{x_1}\Bigr] \quad\text{by elem. rel. \eqref{uwd-move-b1}}
\end{split}\]
The last iterated operadic composition above is \eqref{wasted-cable-y-simplex-b}.  In other words, one can go from the decomposition \eqref{wasted-cable-y-simplex} of $\zeta_Y$ to \eqref{wasted-cable-y-simplex-b} using two elementary relations, the operad vertical associativity axiom twice, and the operad horizontal associativity axiom once.  In the terminology of Chapter \ref{ch10-stratified-uwd}, we say that \eqref{wasted-cable-y-simplex} and \eqref{wasted-cable-y-simplex-b} are stratified presentations (Def. \ref{def:stratified-presentation-uwd}) of $\zeta_Y$, and they are connected by a finite sequence of elementary equivalences (Def. \ref{def:equivalent-simplices-uwd}). 
\end{example}

\section{Summary of Chapter \ref{ch08-generating-uwd}}

\begin{enumerate}
\item There are six generating undirected wiring diagrams.
\item There are seventeen elementary relations in $\uwd$.
\item Wasted cables can arise from operadic composition of undirected wiring diagrams with no wasted cables.
\end{enumerate}

\chapter{Decomposition of Undirected Wiring Diagrams}
\label{ch09-decomp-uwd}

This chapter is the undirected analogue of Chapter \ref{ch04-decomposition}.  As part of the finite presentation theorem for the operad $\uwd$ of undirected wiring diagrams (Theorem \ref{uwd-operad}), in Theorem \ref{stratified-presentation-exists-uwd} we will observe that each undirected wiring diagram has a highly structured decomposition in terms of generators (Def. \ref{def:generating-uwd}), called a stratified presentation (Def. \ref{def:stratified-presentation-uwd}).  Stratified presentations are also needed to establish the second part of the finite presentation theorem for the operad $\uwd$ regarding relations (Theorem \ref{thm:uwd-generator-relation}).  The purpose of this chapter is to provide all the steps needed to establish the existence of a stratified presentation for each undirected wiring diagram.  We remind the reader about Notation \ref{comp-is-compone} for (iterated) operadic compositions.

Fix a class $S$, with respect to which the operad $\uwd$ of undirected wiring diagrams (Def. \ref{uwd-operad}) is defined.

\section{A Motivating Example}
\label{sec:ex-decomp-uwd}

Before we establish the desired decomposition of a general undirected wiring diagram, in this section we consider an elaborate example that will serve as a guide and motivation for the construction later in this chapter for the general case.  The point of this decomposition is to break the complexity of a general undirected wiring diagram into several stratified pieces, each of which is easy to understand and visualize.

The following notations regarding subsets of cables will be used frequently in this chapter.  Recall that an \emph{$(m,n)$-cable} is a cable to which exactly $m$ input wires and exactly $n$ output wires are soldered (Def. \ref{def:pre-uwd}).

\begin{notation}
\label{notation:cable-subsets}
Suppose $\psi = (\cpsi, \fpsi, \gpsi)$ is an undirected wiring diagram and $m, n \geq 0$.  Define:
\begin{itemize}
\item
$\cpsi^{(m,n)} \subseteq \cpsi$ as the subset of $(m,n)$-cables.\index{mncable@$(m,n)$-cable}  In particular, $\cpsizerozero$ is the set of wasted cables.
\item
$\cpsi^{(\geq m,n)} \subseteq \cpsi$\index{mncable2@$(\geq m, n)$-cable} as the subset of $(l,n)$-cables with $l \geq m$.
\item
$\cpsi^{(m,\geq n)} \subseteq \cpsi$\index{mncable3@$(m, \geq n)$-cable} as the subset of $(m,k)$-cables with $k \geq n$.
\item
$\cpsi^{(\geq m,\geq n)} \subseteq \cpsi$\index{mncable4@$(\geq m, \geq n)$-cable} as the subset of $(j,k)$-cables with $j \geq m$ and $k \geq n$.
\item
$\cpsithree \subseteq \cpsi$ as the subset of $(k,l)$-cables with $k,l \geq 1$ and $k+l \geq 3$.
\end{itemize}
A cable in $\cpsi^{(\geq m,n)}$ is called an \emph{$(\geq m,n)$-cable}, and similarly for cables in the other subsets defined above.
\end{notation}

As in the case of wiring diagrams (see Convention \ref{conv:ignore-name-change}), name changes (Def. \ref{def:uwd-gen-namechange}) are easy to deal with.  Therefore, in the following example, to keep the presentation simple, we will ignore name changes.

\begin{example}
\label{ex:uwd-decomp}
Consider $\varphi = \left(\cphi, \fphi, \gphi\right)\in \uwd\yxonextwo$ in \eqref{uwd-first-picture}, which is visualized as
\begin{center}
\begin{tikzpicture}[scale=1]
\draw [ultra thick] (-1,0) rectangle (6.5,3);
\node at (-1.2,1) {\tiny{$y_1$}};
\node at (-1.2,2) {\tiny{$y_2$}};
\node at (.3,3.2) {\tiny{$y_3$}};
\node at (2.5,3.2) {\tiny{$y_4$}};
\node at (3.5,3.2) {\tiny{$y_5$}};
\node at (6.7,1.7) {\tiny{$y_6$}};
\node at (6,.3) {$Y$};
\draw [ultra thick] (.3,1) rectangle (2.5,2);
\node at (.5,1.85) {\tiny{$x_1$}};
\node at (1.3,1.85) {\tiny{$x_2$}};
\node at (2.3,1.85) {\tiny{$x_3$}};
\node at (2.3,1.15) {\tiny{$x_4$}};
\node at (1.3,1.15) {\tiny{$x_5$}};
\node at (.5,1.15) {\tiny{$x_6$}};
\draw [ultra thick] (4,1) rectangle (5,2);
\node at (4.5,1.15) {\tiny{$x^2$}};
\node at (4.5,1.8) {\tiny{$x^1$}};
\cable{(-.25,1.5)} \node at (-.25,1.5) {\tiny{$c_1$}};
\cable{(.5,2.5)} \node at (.5,2.5) {\tiny{$c_2$}};
\cable{(3,2.5)} \node at (3,2.5) {\tiny{$c_3$}};
\cable{(6,2.5)} \node at (6,2.5) {\tiny{$c_4$}};
\cable{(6,1.5)} \node at (6,1.5) {\tiny{$c_5$}};
\cable{(3,.5)} \node at (3,.5) {\tiny{$c_6$}};
\cable{(.5,.5)} \node at (.5,.5) {\tiny{$c_7$}};
\draw [thick] (-.425,1.5) -- (-1.7,1);
\draw [thick] (-.425,1.5) -- (-1.7,2);
\draw [thick] (.5,2) -- (.5,2.325);
\draw [thick] (.5,2.675) -- (.5,3.5);
\draw [thick] (1.3,2) to [out=90, in=180] (2.825,2.5);
\draw [thick] (2.3,2) to [out=90, in=180] (2.825,2.5);
\draw [thick] (4.5,2)  to [out=90, in=0] (3.175,2.5);
\draw [thick] (3,2.675) -- (2.7,3.5);
\draw [thick] (3,2.675) -- (3.3,3.5);
\draw [thick] (6.175,1.5) -- (7,1.5);
\draw [thick] (1.3,1) to [out=-90,in=180] (2.825,.5);
\draw [thick] (2.3,1) to [out=-90,in=180] (2.825,.5);
\draw [thick] (4.5,1) to [out=-90,in=0]  (3.175,.5);
\draw [thick] (.5,1) -- (.5,.675);
\end{tikzpicture}
\end{center}
with $X_1 = \{x_1,\ldots, x_6\}$ and $X_2 = \{x^1,x^2\}$.  We can decompose it as 
\begin{equation}
\label{ex:varphi=varphi1varphi2}
\varphi = \phione \comp \phitwo
\end{equation}
as indicated in the following picture.
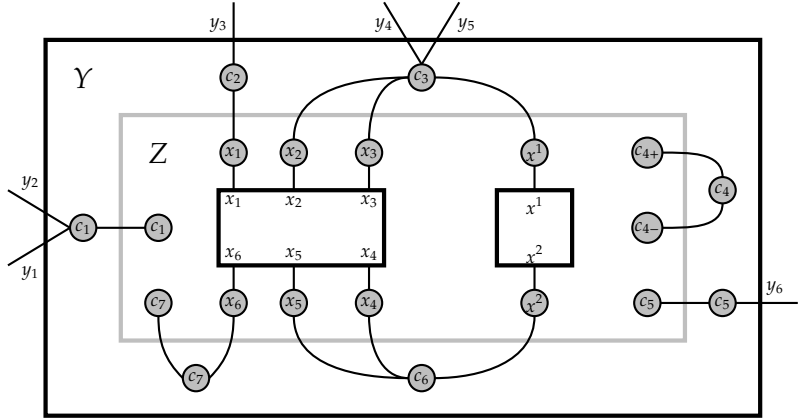
\begin{center}
\begin{tikzpicture}[scale=1]
\draw [ultra thick] (-2,-1) rectangle (7.5,4);
\node at (-1.5,3.5) {$Y$};
\node at (-2.2,.9) {\tiny{$y_1$}};
\node at (-2.2,2.1) {\tiny{$y_2$}};
\node at (.3,4.2) {\tiny{$y_3$}};
\node at (2.5,4.2) {\tiny{$y_4$}};
\node at (3.6,4.2) {\tiny{$y_5$}};
\node at (7.7,.7) {\tiny{$y_6$}};
\draw [ultra thick, lightgray] (-1,0) rectangle (6.5,3);
\node at (-.5,2.5) {$Z$};
\draw [ultra thick] (.3,1) rectangle (2.5,2);
\node at (.5,1.85) {\tiny{$x_1$}};
\node at (1.3,1.85) {\tiny{$x_2$}};
\node at (2.3,1.85) {\tiny{$x_3$}};
\node at (2.3,1.15) {\tiny{$x_4$}};
\node at (1.3,1.15) {\tiny{$x_5$}};
\node at (.5,1.15) {\tiny{$x_6$}};
\draw [ultra thick] (4,1) rectangle (5,2);
\node at (4.5,1.15) {\tiny{$x^2$}};
\node at (4.5,1.8) {\tiny{$x^1$}};
\cable{(-.5,1.5)} \node at (-.5,1.5) {\tiny{$c_1$}};
\cable{(.5,2.5)} \node at (.5,2.5) {\tiny{$x_1$}};
\cable{(1.3,2.5)} \node at (1.3,2.5) {\tiny{$x_2$}};
\cable{(2.3,2.5)} \node at (2.3,2.5) {\tiny{$x_3$}};
\cable{(4.5,2.5)} \node at (4.5,2.5) {\tiny{$x^1$}};
\midcable{(6,2.5)} \node at (6,2.5) {\tiny{$c_{4+}$}};
\midcable{(6,1.5)} \node at (6,1.5) {\tiny{$c_{4-}$}};
\cable{(6,.5)} \node at (6,.5) {\tiny{$c_5$}};
\cable{(4.5,.5)} \node at (4.5,.5) {\tiny{$x^2$}};
\cable{(2.3,.5)} \node at (2.3,.5) {\tiny{$x_4$}};
\cable{(1.3,.5)} \node at (1.3,.5) {\tiny{$x_5$}};
\cable{(.5,.5)} \node at (.5,.5) {\tiny{$x_6$}};
\cable{(-.5,.5)} \node at (-.5,.5) {\tiny{$c_7$}};
\cable{(-1.5,1.5)} \node at (-1.5,1.5) {\tiny{$c_1$}};
\cable{(.5,3.5)} \node at (.5,3.5) {\tiny{$c_2$}};
\cable{(3,3.5)} \node at (3,3.5) {\tiny{$c_3$}};
\cable{(7,2)} \node at (7,2) {\tiny{$c_4$}};
\cable{(7,.5)} \node at (7,.5) {\tiny{$c_5$}};
\cable{(3,-.5)} \node at (3,-.5) {\tiny{$c_6$}};
\cable{(0,-.5)} \node at (0,-.5) {\tiny{$c_7$}};
\draw [thick] (-.675,1.5) -- (-1.325,1.5);
\draw [thick] (-1.675,1.5) -- (-2.5,2);
\draw [thick] (-1.675,1.5) -- (-2.5,1);
\draw [thick] (.5,2) -- (.5,2.325);
\draw [thick] (.5,2.675) -- (.5,3.325);
\draw [thick] (.5,3.675) -- (.5,4.5);
\draw [thick] (1.3,2) -- (1.3,2.325);
\draw [thick] (1.3,2.675) to [out=90,in=180] (2.825,3.5);
\draw [thick] (2.3,2) -- (2.3,2.325);
\draw [thick] (2.3,2.675) to [out=90,in=180] (2.825,3.5);
\draw [thick] (4.5,2) -- (4.5,2.325);
\draw [thick] (4.5,2.675) to [out=90,in=0] (3.175,3.5);
\draw [thick] (3,3.675) to (2.5,4.5);
\draw [thick] (3,3.675) to (3.5,4.5);
\draw [thick] (6.2,2.5) to [out=0,in=90] (7,2.175);
\draw [thick] (7,1.825) to [out=-90,in=0] (6.2,1.5);
\draw [thick] (6.175,.5) -- (6.825,.5);
\draw [thick] (7.175,.5) -- (8,.5);
\draw [thick] (1.3,1) -- (1.3,.675);
\draw [thick] (1.3,.325) to [out=-90,in=180] (2.825,-.5);
\draw [thick] (2.3,1) -- (2.3,.675);
\draw [thick] (2.3,.325) to [out=-90,in=180] (2.825,-.5);
\draw [thick] (4.5,1) -- (4.5,.675);
\draw [thick] (4.5,.325) to [out=-90,in=0] (3.175,-.5);
\draw [thick] (.5,1) -- (.5,.675);
\draw [thick] (.5,.325) to [out=-90,in=45] (.175,-.5);
\draw [thick] (-.175,-.5) to [out=135,in=-90] (-.5,.325);
\end{tikzpicture}
\end{center}
As before the intermediate gray box $Z$ indicates that an operadic composition occurs along it.  The box $Z$ is defined as 
\[
Z = X \amalg \left\{c_{4+}, c_{4-}\right\} \amalg \{c_1,c_5\} \amalg \{c_7\} \in \Fins\]
in which:
\begin{itemize}
\item
$X = X_1 \amalg X_2$.
\item
$c_{4+}$ and $c_{4-}$ are two copies of the wasted cable $c_4$ in $\varphi$, so $\left\{c_{4+}, c_{4-}\right\} = \cphizerozero \amalg \cphizerozero$.
\item
$\{c_1,c_5\} = \cphizeroatleastone$.
\item
$\{c_7\} = \cphionezero$.
\end{itemize}
So we may also write $Z$ as
\begin{equation}
\label{ex:z-cables-decomp}
Z = X \amalg \cphizerozeropm \amalg \cphizeroatleastone \amalg \cphionezero
\end{equation}
in which $\cphizerozeropm = \cphizerozero \amalg \cphizerozero$ is the coproduct of two copies of the set of wasted cables $\cphizerozero$ in $\varphi$.

In the decomposition \eqref{ex:varphi=varphi1varphi2} of $\varphi$, the inside undirected wiring diagram is the cospan
\[\phitwo = \Bigl(\nicexy@C+.4cm{X \ar[r]^-{\mathrm{inclusion}} & Z & Z \ar[l]_-{=}}\Bigr) \in \uwd\zxonextwo.\]
Note that in $\phitwo$:
\begin{itemize}
\item
All the input wires--i.e., those in $X$--are soldered to $(1,1)$-cables.
\item
All other cables--i.e., those in $\cphizerozeropm \amalg \cphizeroatleastone \amalg \cphionezero$--are $(0,1)$-cables.
\item
There are no wasted cables.
\end{itemize}
As we will see later in \eqref{psitwo-2cell-1outputwire}, such an undirected wiring diagram can be decomposed into $2$-cells (Def. \ref{def:uwd-gen-2cell}) and $1$-output wires (Def. \ref{def:uwd-gen-1output}).  For example, this $\phitwo$ needs:
\begin{itemize}
\item
five $1$-output wires, exactly as many as the number of $(0,1)$-cables;
\item
six $2$-cells, where $6$ is the number of input boxes plus the number of $(0,1)$-cables minus $1$.
\end{itemize}

The outside undirected wiring diagram in the decomposition $\varphi = \phione \comp \phitwo$ is the cospan
\[\phione = \Bigl(\nicexy{Z \ar[r]^-{(\fphi, \iota)} & \cphi & Y \ar[l]_-{\gphi}} \Bigr) \in \uwd \yz.\]
Here:
\begin{itemize}
\item
$\fphi : X \to \cphi$ is the input soldering function of $\varphi$.
\item
$\iota : \cphizerozeropm \amalg \cphizeroatleastone \amalg \cphionezero \to \cphi$ is the inclusion map on each coproduct summand.
\item
$\gphi : Y \to \cphi$ is the output soldering function of $\varphi$.
\item
Every cable is an $(m,n)$-cable with $m \geq 1$ and $n \geq 0$.  In other words, every cable in $\phione$ is soldered to some input wires, so in particular there are no wasted cables in $\phione$.
\item
There are also no $(1,0)$-cables, but there are $(\geq 2,0)$-cables.
\end{itemize}
As we will see later, such an undirected wiring diagram can be decomposed into loops (Def. \ref{def:uwd-gen-loop}) and splits (Def. \ref{def:uwd-gen-split}).  In the case of $\phione$, which is the undirected wiring diagram
\begin{center}
\begin{tikzpicture}[scale=.8]
\draw [ultra thick] (-2,-1) rectangle (7.5,4);
\node at (-1.5,3.5) {$Y$};
\node at (-2.2,.9) {\tiny{$y_1$}};
\node at (-2.2,2.1) {\tiny{$y_2$}};
\node at (.3,4.2) {\tiny{$y_3$}};
\node at (2.5,4.2) {\tiny{$y_4$}};
\node at (3.6,4.2) {\tiny{$y_5$}};
\node at (7.8,.7) {\tiny{$y_6$}};
\draw [ultra thick] (-1,0) rectangle (6.5,3);
\node at (2.75,1.5) {$Z$};
\node at (.5,2.8) {\tiny{$x_1$}};
\node at (1.3,2.8) {\tiny{$x_2$}};
\node at (2.3,2.8) {\tiny{$x_3$}};
\node at (2.3,.15) {\tiny{$x_4$}};
\node at (1.3,.15) {\tiny{$x_5$}};
\node at (.5,.15) {\tiny{$x_6$}};
\node at (4.5,.2) {\tiny{$x^2$}};
\node at (4.5,2.75) {\tiny{$x^1$}};
\node at (-.8,1.5) {\tiny{$c_1$}};
\node at (6.25,2.5) {\tiny{$c_{4+}$}};
\node at (6.25,1.5) {\tiny{$c_{4-}$}};
\node at (6.3,.5) {\tiny{$c_5$}};
\node at (-.5,.2) {\tiny{$c_7$}};
\cable{(-1.5,1.5)} \node at (-1.5,1.5) {\tiny{$c_1$}};
\cable{(.5,3.5)} \node at (.5,3.5) {\tiny{$c_2$}};
\cable{(3,3.5)} \node at (3,3.5) {\tiny{$c_3$}};
\cable{(7,2)} \node at (7,2) {\tiny{$c_4$}};
\cable{(7,.5)} \node at (7,.5) {\tiny{$c_5$}};
\cable{(3,-.5)} \node at (3,-.5) {\tiny{$c_6$}};
\cable{(0,-.5)} \node at (0,-.5) {\tiny{$c_7$}};
\draw [thick] (-1,1.5) -- (-1.325,1.5);
\draw [thick] (-1.675,1.5) -- (-2.5,2);
\draw [thick] (-1.675,1.5) -- (-2.5,1);
\draw [thick] (.5,3) -- (.5,3.325);
\draw [thick] (.5,3.675) -- (.5,4.5);
\draw [thick] (1.3,3) to [out=60,in=180] (2.825,3.5);
\draw [thick] (2.3,3) to [out=80,in=180] (2.825,3.5);
\draw [thick] (4.5,3) to [out=120,in=0] (3.175,3.5);
\draw [thick] (3,3.675) to (2.5,4.5);
\draw [thick] (3,3.675) to (3.5,4.5);
\draw [thick] (6.5,2.5) to [out=0,in=90] (7,2.175);
\draw [thick] (7,1.825) to [out=-90,in=0] (6.5,1.5);
\draw [thick] (6.5,.5) -- (6.825,.5);
\draw [thick] (7.175,.5) -- (8,.5);
\draw [thick] (1.3,0) to [out=-60,in=180] (2.825,-.5);
\draw [thick] (2.3,0) to [out=-80,in=180] (2.825,-.5);
\draw [thick] (4.5,0) to [out=-120,in=0] (3.175,-.5);
\draw [thick] (.5,0) to [out=-90,in=15] (.175,-.5);
\draw [thick] (-.175,-.5) to [out=165,in=-90] (-.5,0);
\end{tikzpicture}
\end{center}
this further decomposition
\begin{equation}
\label{ex:varphi1=phi1phi2}
\phione = \phi_1 \comp \phi_2
\end{equation}
can be visualized as follows.
\begin{equation}
\label{phione-decomp-example-picture}
\begin{tikzpicture}[scale=.8]
\draw [ultra thick] (0,0) rectangle (6,3);
\node at (3,1.5) {$Z$};
\node at (.5,2.8) {\tiny{$x_1$}};
\node at (2.5,2.8) {\tiny{$x_2$}};
\node at (4,2.8) {\tiny{$x_3$}};
\node at (4,.15) {\tiny{$x_4$}};
\node at (2.5,.15) {\tiny{$x_5$}};
\node at (1.5,.15) {\tiny{$x_6$}};
\node at (5,.2) {\tiny{$x^2$}};
\node at (5,2.8) {\tiny{$x^1$}};
\node at (.2,1.5) {\tiny{$c_1$}};
\node at (5.75,2.5) {\tiny{$c_{4+}$}};
\node at (5.75,1.5) {\tiny{$c_{4-}$}};
\node at (5.8,.5) {\tiny{$c_5$}};
\node at (.5,.15) {\tiny{$c_7$}};
\draw [ultra thick, lightgray] (-1,-1) rectangle (7,4);
\node at (-.7,3.7) {$W$};
\draw [thick] (0,1.5) -- (-.5,1.5) to [bend right] (-1.5,2) -- (-3,2);
\draw [thick] (-.5,1.5) to [bend left] (-1.5,1) -- (-3,1);
\draw [thick] (.5,3) -- (.5,6);
\draw [thick] (2.5,3) -- (2.5,3.5) -- (1.5,4.5) -- (1.5,6);
\draw [thick] (2.5,3.5) -- (2,4.5) -- (2,6);
\draw [thick] (2.5,3.5) to (2.5,4.5) to [bend left=80] (5,3.5) to (5,3);
\draw [thick] (2.5,3.5) to [bend left] (3.25,4.5) to [bend left] (4,3.5) to (4,3);
\draw [thick] (6,2.5) to (6.5,2.5) to [bend left] (7.5,2) to [bend left] (6.5,1.5) to (6,1.5);
\draw [thick] (6,.5) -- (9,.5);
\draw [thick] (2.5,0) to (2.5,-.5) to [bend right] (3.25,-1.5) to [bend right] (4,-.5) to (4,0);
\draw [thick] (2.5,-.5) to (2.5,-1.5) to [bend right=80] (5,-.5) to (5,0);
\draw [thick] (.5,0) to (.5,-.5) to [bend right] (1,-1.5) to [bend right] (1.5,-.5) to (1.5,0);
\cable{(-.5,1.5)} 
\foreach \x in {.5,2.5,4,5} {\cable{(\x,3.5)}}
\cable{(6.5,2.5)} \cable{(6.5,1.5)} \cable{(6.5,.5)}
\foreach \x in {.5,1.5,2.5,4,5} {\cable{(\x,-.5)}}
\cable{(-1.5,1)} \cable{(-1.5,2)}
\foreach \x in {.5,1.5,2,2.5,3.25} {\cable{(\x,4.5)}}
\cable{(7.5,2)} \cable{(7.5,.5)}
\foreach \x in {1,2.5,3.25} {\cable{(\x,-1.5)}}
\draw [ultra thick] (-2,-2.5) rectangle (8,5.5);
\node at (-1.5,5) {$Y$};
\node at (-2.2,.8) {\tiny{$y_1$}};
\node at (-2.2,2.2) {\tiny{$y_2$}};
\node at (.3,5.7) {\tiny{$y_3$}};
\node at (1.3,5.7) {\tiny{$y_4$}};
\node at (2.2,5.7) {\tiny{$y_5$}};
\node at (8.3,.7) {\tiny{$y_6$}};
\end{tikzpicture}
\end{equation}

In this decomposition $\phione = \phi_1 \comp \phi_2$, the inner undirected wiring diagram is the cospan
\[\phi_2 = \Bigl(\nicexy{Z \ar[r]^-{=} & Z & W \ar[l]_-{g_{\phi_2}}} \Bigr) \in \uwd\wz\]
with $g_{\phi_2}$ surjective.  So every cable in $\phi_2$ is a $(1,n)$-cable for some $n \geq 1$ .  As we will see later, such an undirected wiring diagram is generated by splits (Def. \ref{def:uwd-gen-split}).  For example, this $\phi_2 $ is the iterated operadic composition of $5$ splits--one for the cable soldered to $c_1$, one for the cable soldered to $x_5$, and three for the cable soldered to $x_2$.

The outer undirected wiring diagram in the decomposition \eqref{ex:varphi1=phi1phi2} is the cospan
\[\phi_1 = \Bigl(\nicexy{W \ar[r] & C_{\phi_1} & Y \ar[l]} \Bigr) \in \uwd\yw\]
in which every cable is either a $(1,1)$-cable or a $(2,0)$-cable.  We will show later that such an undirected wiring diagram is generated by loops (Def. \ref{def:uwd-gen-loop}).  For example, this $\phi_1$ is the iterated operadic composition of $6$ loops, where $6$ is the number of $(2,0)$-cables in $\phi_1$.

In summary, we decompose $\varphi \in \uwd\yxonextwo$ as the iterated operadic composition
\begin{equation}
\label{sample-uwd-stratified}
\begin{split}
\varphi &= \phione \comp \phitwo = \phi_1 \comp \phi_2 \comp \phitwo\\
&= \Bigl(\underbrace{\lambda,\ldots,\lambda}_{6}, 
\underbrace{\sigma, \ldots, \sigma}_{5},
\underbrace{\theta,\ldots,\theta}_{6},
\underbrace{\omega,\ldots,\omega}_{5} \Bigr).
\end{split}
\end{equation}
Here $\lambda$, $\sigma$, $\theta$, and $\omega$ denote a loop, a split, a $2$-cell, and a $1$-output wire, respectively.  This decomposition in terms of the generators is called a stratified presentation (Def. \ref{def:stratified-presentation-uwd}).  In the next few sections, we will establish all the steps needed to obtain a stratified presentation for a general undirected wiring diagram.
\end{example}

\section{Factoring Undirected Wiring Diagrams}
\label{sec:decomp-uwd}

In this section, using Example \ref{ex:uwd-decomp} as a guide and motivation, we establish a decomposition of a general undirected wiring diagram into two simpler  undirected wiring diagrams (Theorem \ref{thm:uwd-first-decomp}).  This is the general version of the decomposition \eqref{ex:varphi=varphi1varphi2} above.  Each undirected wiring diagram in this decomposition will be decomposed further, eventually leading to the desired stratified presentation.

\begin{assumption}
\label{assumption:general-uwd}
Suppose
\begin{equation}
\label{psi-uwd}
\psi = \Bigl(\nicexy{X \ar[r]^-{\fpsi} & \cpsi & Y \ar[l]_-{\gpsi}}\Bigr) \in \uwd\yux
\end{equation}
is a general undirected wiring diagram with:
\begin{itemize}
\item
output box $Y \in \Fins$ and input boxes $\uX = \left(X_1, \ldots, X_N\right)$ for some $N \geq 0$;
\item
$X = X_1 \amalg \cdots \amalg X_N \in \Fins$.
\end{itemize}
\end{assumption}

Recall Notation \ref{notation:cable-subsets} for certain subsets of cables.  The undirected wiring diagrams $\psione$ and $\psitwo$ in the next definition are the general versions of $\phione$ and $\phitwo$ in the decomposition \eqref{ex:varphi=varphi1varphi2} above.

\begin{definition}
\label{def:psi1-psi2}
Suppose $\psi = \left(\cpsi,\fpsi,\gpsi\right) \in \uwd\yux$ is a general undirected wiring diagram as in \eqref{psi-uwd} with $\uX = (X_1,\ldots,X_N)$.
\begin{enumerate}
\item
Define
\begin{equation}
\label{uwd-psi-z}
Z = X \amalg \cpsizerozeropm \amalg \cpsizeroatleastone \amalg \cpsionezero \in \Fins
\end{equation}
in which
\[\cpsizerozeropm = \cpsizerozero \amalg \cpsizerozero\]
is the coproduct of two copies of the set of wasted cables $\cpsizerozero$ in $\psi$.
\item
Define the undirected wiring diagram 
\begin{equation}
\label{uwd-psi1}
\psione = \Bigl(\nicexy{Z \ar[r]^-{(\fpsi, \iota)} & \cpsi & Y \ar[l]_-{\gpsi}} \Bigr) \in \uwd \yz
\end{equation}
in which
\[\nicexy{\cpsizerozeropm \amalg \cpsizeroatleastone \amalg \cpsionezero \ar[r]^-{\iota} & \cpsi}\]
is the inclusion map on each coproduct summand.
\item
Define the undirected wiring diagram 
\begin{equation}
\label{uwd-psi2}
\psitwo = \Bigl(\nicexy@C+.4cm{X \ar[r]^-{\text{inclusion}} & Z & Z \ar[l]_-{=}}\Bigr) \in \uwd\zux.
\end{equation}
\end{enumerate}
\end{definition}

\begin{theorem}
\label{thm:uwd-first-decomp}
In the context of Def. \ref{def:psi1-psi2}, there is a decomposition \index{decomposition}
\begin{equation}
\label{psi-is-psi1psi2}
\psi = \psione \comp \psitwo \in \uwd\yux.
\end{equation}
\end{theorem}

\begin{proof}
By the definition of $\comp = \compone$ (Def. \ref{def:compi-uwd}), the operadic composition $\psione \comp \psitwo$ is given by the cospan
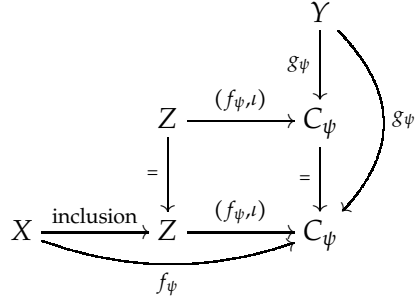
\[\nicexy@C+.4cm{
&& Y \ar[d]_-{\gpsi} \ar@/^2pc/[dd]^-{\gpsi}\\
& Z \ar[r]^-{(\fpsi,\iota)} \ar[d]_-{=} & \cpsi \ar[d]_-{=}\\
X \ar[r]^-{\text{inclusion}} \ar@/_1pc/[rr]_-{\fpsi} & Z \ar[r]^-{(\fpsi,\iota)} & \cpsi}\]
in $\Fins$.  The square is a pushout by Example \ref{ex:pushout-id}.  This cospan is equal to $\psi$.
\end{proof}

\begin{example}
\label{ex:uwd-empty-decomp}
If $\psi = \epsilon \in \uwd\emptynothing$ is the empty cell (Def. \ref{def:uwd-gen-empty}), then:
\begin{itemize}
\item
$\psi_1 = \tensorunit_{\varnothing}$, the $\varnothing$-colored unit (Def. \ref{def:uwd-units}) with $\varnothing \in \Fins$ the empty box;
\item
$\psitwo = \epsilon$.  
\end{itemize}
So in this case the decomposition \eqref{psi-is-psi1psi2} simply says $\epsilon = \tensorunit_{\varnothing} \comp \epsilon$.
\end{example}

\begin{remark}
\label{rk:psione-psitwo}
In the decomposition \eqref{psi-is-psi1psi2}, both $\psione$ and $\psitwo$ are simpler than $\psi$ for the following reasons.
\begin{enumerate}
\item
$\psione$ has the same set of cables $\cpsi$ and the same output soldering function $\gpsi$ as $\psi$.  Furthermore, its input soldering function $(\fpsi,\iota)$ includes the input soldering function $\fpsi$ of $\psi$.  However, every cable in $\psione$ is soldered to at least one input wire (i.e., $\cpsionezeroatleastzero = \varnothing$), whereas $\cpsizeroatleastzero$ may be non-empty.  In particular, $\psione$ has no wasted cables, even though $\psi$ may have some. Furthermore, $\psione$ has only one input box $Z$, while $\psi$ has $N \geq 0$ input boxes.
\item
$\psitwo$ has the same input boxes $\uX$ as $\psi$.  However, it is, in general, much simpler than $\psi$ and $\psione$ because its cables are either $(1,1)$-cables or $(0,1)$-cables. In particular, $\psitwo$ also has no wasted cables.
\item
Neither $\psione$ nor $\psitwo$ has any $(1,0)$-cables, even though $\psi$ may have some.
\end{enumerate}
\end{remark}

\section{The Inner Undirected Wiring Diagram}
\label{sec:psitwo-decomp}

The purpose of this section is to analyze the undirected wiring diagram $\psitwo$ in the decomposition \eqref{psi-is-psi1psi2}.  The undirected wiring diagram $\psione$ will be studied in the next few sections.  We begin with the following observation regarding iterated operadic compositions of $2$-cells (Def. \ref{def:uwd-gen-2cell}).

\begin{motivation}
The following result says that an undirected wiring diagram of the form
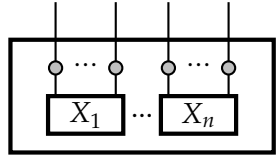
\begin{center}
\begin{tikzpicture}[scale=.5]
\draw [ultra thick] (-1,0) rectangle (6,3);
\node at (2.5,1) {$...$};
\draw [ultra thick] (0,.5) rectangle (2,1.5);
\node at (1,1) {$X_1$};
\draw [thick] (.2,1.5) -- (.2,4); \cable{(.2,2.25)}
\node at (1,2.25) {{$\cdots$}};
\draw [thick] (1.8,1.5) -- (1.8,4); \cable{(1.8,2.25)}
\draw [ultra thick] (3,.5) rectangle (5,1.5);
\node at (4,1) {$X_n$};
\draw [thick] (3.2,1.5) -- (3.2,4); \cable{(3.2,2.25)}
\node at (4,2.25) {{$\cdots$}};
\draw [thick] (4.8,1.5) -- (4.8,4); \cable{(4.8,2.25)}
\end{tikzpicture}
\end{center}
can be generated by $2$-cells.
\end{motivation}

\begin{proposition}
\label{uwd-iterated-2cells}
Suppose $n \geq 2$, $X_i \in \Fins$ for $1 \leq i \leq n$, and $X = \coprod_{i=1}^n X_i$.  Then the undirected wiring diagram
\begin{equation}
\label{iterated-2cells}
\Theta = \Bigl(\nicexy{X \ar[r]^-{=} & X & X \ar[l]_-{=}}\Bigr) \in \uwd\xxonexn
\end{equation}
is:
\begin{itemize}
\item
a $2$-cell if $n=2$;
\item
an iterated operadic composition
\[
\Theta = \Bigl(\bigl(\theta_1 \comptwo \theta_2\bigr) \compthree \cdots \Bigr) \comp_{n-1} \theta_{n-1}\]
with each $\theta_j$ a $2$-cell if $n \geq 3$.
\end{itemize}
\end{proposition}

\begin{proof}
This is proved by induction on $n \geq 2$.  The initial case simply says that $\Theta$ is the $2$-cell $\theta_{(X_1,X_2)}$ by Def. \ref{def:uwd-gen-2cell}.  

Suppose $n \geq 3$.  By the definition of $\comp_{n-1}$ (Def. \ref{def:compi-uwd}) and Example \ref{ex:pushout-id}, we may decompose $\Theta$ as
\[\Theta = \Theta_1 \comp_{n-1} \theta_{(X_{n-1},X_n)}\]
in which
\[\Theta_1 = \Bigl(\nicexy{X \ar[r]^-{=} & X & X \ar[l]_-{=}}\Bigr) \in \uwd\smallbinom{X}{X_1,\ldots,X_{n-2}, X_{n-1} \amalg X_n}\]
and
\[\theta_{(X_{n-1},X_n)} \in \uwd\xnminusonexn\] 
is a $2$-cell.  Since the induction hypothesis applies to $\Theta_1$, the proof is finished.
\end{proof}

\begin{example}
In the previous Proposition:
\begin{enumerate}
\item
If $n=3$, then $\Theta$ decomposes as
\[\Theta = \theta_{(X_1, X_2 \amalg X_3)} \comptwo \theta_{(X_2,X_3)}\]
into two $2$-cells.  
\item
If $n=4$, then $\Theta$ decomposes as
\[\Theta = \Bigl(\theta_{(X_1, X_2 \amalg X_3 \amalg X_4)} 
\comptwo \theta_{(X_2, X_3 \amalg X_4)}\Bigr)
\compthree \theta_{(X_3,X_4)}\]
into three $2$-cells.
\end{enumerate}
\end{example}

\begin{notation}
\label{uwd-psitwo-notation}
In the context of \eqref{uwd-psi-z}, write:
\begin{itemize}
\item
$\cpsi' = \cpsizerozeropm \amalg \cpsizeroatleastone \amalg \cpsionezero \in \Fins$, so $Z = X \amalg \cpsi'$.
\item
$p = \left|\cpsi'\right|$.
\end{itemize}
\end{notation}

The following observation covers the marginal cases for $\psitwo$.

\begin{lemma}
\label{uwd-psitwo-marginal}
For $\psitwo = \Bigl(\nicexy@C+.2cm{X \ar[r]^-{\mathrm{inclusion}} & Z & Z \ar[l]_-{=}}\Bigr) \in \uwd\zux$ in \eqref{uwd-psi2}:
\begin{enumerate}
\item
If $N = p = 0$, then $\psitwo$ is the empty cell $\epsilon$ (Def. \ref{def:uwd-gen-empty}).
\item
If $N=0$ and $p=1$, then $\psitwo$ is a $1$-output wire (Def. \ref{def:uwd-gen-1output}).
\item
If $N=1$ and $p=0$, then $\psitwo$ is the $X_1$-colored unit (Def. \ref{def:uwd-units}).
\end{enumerate}
\end{lemma}

\begin{proof}
Since $\uX = (X_1,\ldots,X_N)$ and $X = X_1 \amalg \cdots \amalg X_N$, all three statements follow immediately from the definition of $\psitwo$.
\end{proof}

The next observation covers the other cases for $\psitwo$.  Recall $\cpsi'$ in Notation \ref{uwd-psitwo-notation}.

\begin{motivation}
The following result says that an undirected wiring diagram of the form
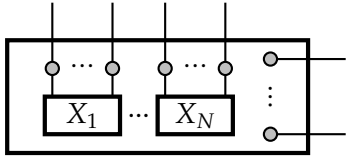
\begin{center}
\begin{tikzpicture}[scale=.5]
\draw [ultra thick] (-1,0) rectangle (7,3);
\node at (2.5,1) {$...$};
\draw [ultra thick] (0,.5) rectangle (2,1.5);
\node at (1,1) {$X_1$};
\draw [thick] (.2,1.5) -- (.2,4); \cable{(.2,2.25)}
\node at (1,2.25) {{$\cdots$}};
\draw [thick] (1.8,1.5) -- (1.8,4); \cable{(1.8,2.25)}
\draw [ultra thick] (3,.5) rectangle (5,1.5);
\node at (4,1) {$X_N$};
\draw [thick] (3.2,1.5) -- (3.2,4); \cable{(3.2,2.25)}
\node at (4,2.25) {{$\cdots$}};
\draw [thick] (4.8,1.5) -- (4.8,4); \cable{(4.8,2.25)}
\draw [thick] (6,2.5) -- (8,2.5); \cable{(6,2.5)}
\draw [thick] (6,.5) -- (8,.5); \cable{(6,.5)}
\node at (6,1.5) {$\vdots$};
\end{tikzpicture}
\end{center}
can be generated by $2$-cells and $1$-output wires.
\end{motivation}

\begin{proposition}
\label{uwd-psitwo-general}
For $\psitwo = \Bigl(\nicexy@C+.2cm{X \ar[r]^-{\mathrm{inclusion}} & Z & Z \ar[l]_-{=}}\Bigr) \in \uwd\zux$ in \eqref{uwd-psi2}, suppose:
\begin{itemize}
\item
$N, p \geq 1$, and $\cpsi' = \{c_1,\ldots,c_p\}$;
\item
$\omega_j =  \Bigl(\nicexy{\varnothing \ar[r] & c_j & c_j \ar[l]_-{=}}\Bigr)  \in \uwd\cjnothing$ is the $1$-output wire for $c_j$ (Def. \ref{def:uwd-gen-1output}) for $1 \leq j \leq p$.
\end{itemize}
Then there is a decomposition
\begin{equation}
\label{uwd-psi-two-general-decomp}
\psitwo = \Bigl(\bigl(\Theta \comp_{N+1} \omega_1\bigr) \cdots \Bigr) \comp_{N+1} \omega_p \in \uwd\zux
\end{equation}
in which every pair of parentheses starts on the left and
\[\Theta = \Bigl(\nicexy{Z \ar[r]^-{=} & Z & Z \ar[l]_-{=}}\Bigr) \in \uwd\zxonexnconecp.\]
\end{proposition}

\begin{proof}
The right side of \eqref{uwd-psi-two-general-decomp} is a well-defined element in $\zux$.  By the correspondence between the $\compi$-compositions and $\gamma$ \eqref{gamma-in-comps} in the operad $\uwd$, the right side of \eqref{uwd-psi-two-general-decomp} can be rewritten as
\[\psitwo' = \gamma\Bigl(\Theta; \tensorunit_{X_1},\ldots,\tensorunit_{X_N},\omega_1,\ldots,\omega_p\Bigr).\]
Since $Z = X \amalg \{c_1,\ldots,c_p\}$, by Prop. \ref{uwd-gamma} the cospan for $\psitwo'$ is 
\[\nicexy@C+.4cm{&& Z \ar[d]^-{=}\\
& Z \ar[d]_-{=} \ar[r]^-{=} & Z \ar[d]^-{=}\\
X \ar[r]^-{\mathrm{inclusion}} & Z \ar[r]^-{=} & Z}\]
in $\Fins$.  This is equal to the cospan that defines $\psitwo$.
\end{proof}

The following observation says that, if $N,p \geq 1$, then $\psitwo$ is generated by $2$-cells and $1$-output wires.

\begin{corollary}
\label{psitwo-is-2cells-1outputwires}
Suppose $\psitwo \in \uwd\zux$ in \eqref{uwd-psi2} has $N = |\uX|, p = \left|\cpsi'\right| \geq 1$.  Then there is a decomposition
\begin{equation}
\label{psitwo-2cell-1outputwire}
\psitwo = \left[\left[\Bigl(\bigl(\bigl(\theta_1 \comptwo \theta_2\bigr) \compthree \cdots \bigr) \comp_{N+p-1} \theta_{N+p-1}\Bigr) \comp_{N+1} \omega_1\right] \cdots \right] \comp_{N+1} \omega_p 
\end{equation}
with:
\begin{itemize}
\item
each $\theta_i$ a $2$-cell;
\item
each $\omega_j$ a $1$-output wire;
\item
each pair of parentheses starting on the left.
\end{itemize}
\end{corollary}

\begin{proof}
This is true by the decomposition \eqref{uwd-psi-two-general-decomp} above and Prop. \ref{uwd-iterated-2cells} with $n=N+p \geq 2$, applied to $\Theta$.
\end{proof}

\section{The Outer Undirected Wiring Diagram}
\label{sec:psione-decomp}

The purpose of this section is to establish a decomposition for the undirected wiring diagram $\psione$ \eqref{uwd-psi1} that appeared in \eqref{psi-is-psi1psi2}.  This is the general version of the decomposition \eqref{ex:varphi1=phi1phi2}, so the reader may wish to refer back there for specific examples of the constructions below.   Each of the constituent undirected wiring diagrams in this decomposition will be studied further in later sections.  The goal is to decompose $\psione$ into two undirected wiring diagrams in which the outer one, called $\phi_1$ below, is generated by loops (Def. \ref{def:uwd-gen-loop}), while the inner one, called $\phi_2$ below, is generated by splits (Def. \ref{def:uwd-gen-split}).

Recall Notation \ref{notation:cable-subsets} for certain subsets of cables.  Also recall from Remark \ref{rk:psione-psitwo} that $\psione \in \uwd\yz$ has neither  $(0, \geq 0)$-cables nor $(1,0)$-cables.  So $\psione$ satisfies the hypotheses of the next definition.

\begin{definition}
\label{def:psione-general}
Suppose $\varphi = \left(\cphi, \fphi, \gphi\right) \in \uwd\capba$ is an undirected wiring diagram with 
\begin{itemize}
\item one input box $A$ and 
\item $\cphizeroatleastzero = \varnothing = \cphionezero$.  
\end{itemize}
We will write $\fphi$ and $\gphi$ as $f$ and $g$, respectively.
\begin{enumerate}
\item
For each cable $c \in \cphiatleastthreezero \amalg \cphithree$, choose a wire $a_c \in \finv c \subseteq A$, where $\finv c = \finv(\{c\})$ is the set of $f$-preimages of $c$. 
\item
Define 
\begin{equation}
\label{w-for-phi}
W = B \amalg \finv\cphitwozero \amalg \coprod_{c \in \cphiatleastthreezero \amalg \cphithree} \Bigl[\finv c \setminus a_c \Bigr]_{\pm} \in \Fins
\end{equation}
in which
\[\Bigl[\finv c \setminus a_c \Bigr]_{\pm} = \Bigl[\finv c \setminus a_c \Bigr]_+ \amalg \Bigl[\finv c \setminus a_c \Bigr]_-\]
is the coproduct of two copies of $\Bigl[\finv c \setminus a_c \Bigr]$.  This $W$ is the general version of the $W$ in the example \eqref{phione-decomp-example-picture}.
\item
Define
\begin{equation}
\label{v-for-phi}
V = B \amalg \cphitwozero \amalg \coprod_{c \in \cphiatleastthreezero \amalg \cphithree} \Bigl[\finv c \setminus a_c \Bigr] \in \Fins.
\end{equation}
This $V$ is the general version of the set of cables between $W$ and $Y$ in the example \eqref{phione-decomp-example-picture}.
\item
Define
\begin{equation}
\label{phione-general}
\phi_1 = \Bigl(\nicexy@C+.4cm{W \ar[r]^-{f_1} & V  & B \ar[l]_-{g_1}^-{\text{inclusion}}}\Bigr) \in \uwd\bw 
\end{equation}
in which the restrictions of $f_1$ to the coproduct summands of $W$ are defined as follows.
\begin{itemize}
\item $f_1 : B \to B$ is the identity map.
\item $f_1 : \finv \cphitwozero \to \cphitwozero$ is the map $f$.
\item $f_1 : \Bigl[\finv c \setminus a_c \Bigr]_{\pm} \to \Bigl[\finv c \setminus a_c \Bigr]$ is the fold map for each $c \in \cphiatleastthreezero \amalg \cphithree$.  That is, the restriction of $f_1$ to each of $\Bigl[\finv c \setminus a_c \Bigr]_+$ and $\Bigl[\finv c \setminus a_c \Bigr]_-$ is the identity map.
\end{itemize}
This $\phi_1$ is the general version of that in the example \eqref{phione-decomp-example-picture}.
\item
Define
\begin{equation}
\label{phitwo-general}
\phi_2 = \Bigl(\nicexy{A \ar[r]^-{f_2}_-{=} & A & W \ar[l]_-{g_2}}\Bigr) \in \uwd\wa
\end{equation}
as follows.  We will use the equality
\[B = \ginv\cphioneatleastone \amalg \ginv\cphiatleasttwoatleastone\]
which is true because $\cphizeroatleastzero = \varnothing = \cphionezero$.  For
\begin{equation}
\label{w-coproduct-summands}
w \in W = \ginv\cphioneatleastone \amalg \ginv\cphiatleasttwoatleastone \amalg \finv\cphitwozero \amalg \coprod_{c \in \cphiatleastthreezero \amalg \cphithree} \Bigl[\finv c \setminus a_c \Bigr]_{\pm} \in \Fins
\end{equation}
define
\[g_2(w) = \begin{cases}
\finv g(w) \in A & \text{if $w \in \ginv\cphioneatleastone$};\\
a_{g(w)} \in \finv g(w) \subseteq A & \text{if $w \in \ginv\cphiatleasttwoatleastone$};\\
w & \text{if $w \in \finv\cphitwozero$ or $w \in \Bigl[\finv c \setminus a_c \Bigr]_+$};\\
a_c \in \finv c \subseteq A & \text{if $w \in \Bigl[\finv c \setminus a_c \Bigr]_-$}.\\ 
\end{cases}\]
In the first line of this definition, we used the fact that each cable in $\cphioneatleastone$ has a unique $f$-preimage in $A$.  In the second line, $a_?$ was defined earlier in the current definition, using the fact $\cphiatleasttwoatleastone \subseteq \cphithree$.   This $\phi_2$ is the general version of that in the example \eqref{phione-decomp-example-picture}.
\end{enumerate}
\end{definition}

\begin{remark}
\label{rk:phi1-phi2}
Consider the previous definition.
\begin{enumerate}
\item
The input soldering function $f_1$ of $\phi_1$ is surjective.  Furthermore, all the cables in $\phi_1$ are either $(1,1)$-cables (namely, those cables in $B \subseteq V$) or $(2,0)$-cables (namely, those in $V \setminus B$).
\item
The output soldering function $g_2$ of $\phi_2$ is surjective because of the assumption $\cphizeroatleastzero = \varnothing = \cphionezero$. 
\end{enumerate}
\end{remark}

\begin{theorem}
\label{thm:psione-decomp}
In the context of Def. \ref{def:psione-general}, there is a decomposition
\begin{equation}
\label{general-phi-is-phione-phitwo}
\varphi = \phi_1 \comp \phi_2.
\end{equation}
\end{theorem}

\begin{proof}
It suffices to check that the operadic composition $\phi_1 \comp \phi_2 \in \uwd\capba$ is given by the cospan
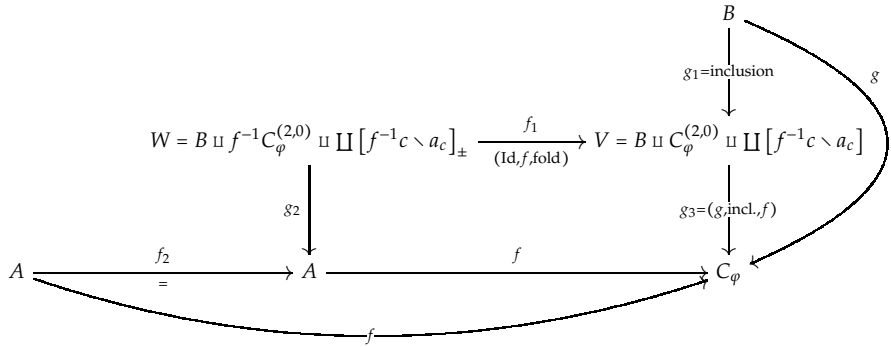
\begin{scriptsize}
\[\nicexy@C+.4cm@R+.3cm{
&& B \ar[d]|-{g_1=\text{inclusion}} \ar@/^5pc/[dd]^(.3){g}\\
& W = B \amalg \finv\cphitwozero \amalg \coprod\left[\finv c \setminus a_c\right]_{\pm} \ar[d]_-{g_2} \ar[r]^-{f_1}_-{(\Id,f,\text{fold})} 
& V = B \amalg \cphitwozero\amalg \coprod\left[\finv c \setminus a_c\right] \ar[d]|-{g_3=(g,\text{incl.},f)}\\
A \ar[r]^-{f_2}_-{=} \ar@/_2pc/[rr]|-{f}& A \ar[r]^-{f} & \cphi}\]
\end{scriptsize}
in $\Fins$.  Here the two coproducts $\coprod$ are both indexed by all $c \in \cphiatleastthreezero \amalg \cphithree$.  By the definition of $\comp = \compone$ \eqref{uwd-compi-cospan}, we just need to check that the rectangle is a pushout (Def. \ref{def:pushout}) in $\Fins$.  It follows from direct inspection of each coproduct summand of $W$ in \eqref{w-coproduct-summands} that the rectangle is commutative.

Next, suppose given a solid-arrow commutative diagram
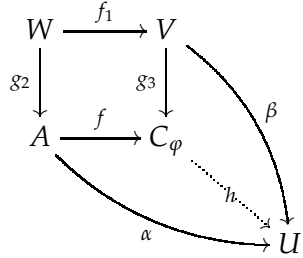
\[\nicexy{
W \ar[d]_-{g_2} \ar[r]^-{f_1} & V \ar[d]_-{g_3} \ar@/^1pc/[ddr]^-{\beta} &\\
A \ar@/_1pc/[drr]_-{\alpha} \ar[r]^-{f} & \cphi \ar@{.>}[dr]|-{h} &\\
&& U}\]
in $\Fins$.  We must show that there exists a unique map $h$ that makes the diagram commutative.  Recall that
\[\cphi = \cphioneone \amalg \cphitwozero \amalg \cphiatleastthreezero \amalg \cphithree\]
because $\cphizeroatleastzero = \varnothing = \cphionezero$.  Define $h : \cphi \to U$ as
\[h(c) = \begin{cases}
\alpha\finv(c) & \text{if $c \in \cphioneone$};\\
\beta(c) & \text{if $c \in \cphitwozero \subseteq V$};\\
\alpha(a_c) & \text{if $c \in \cphiatleastthreezero \amalg \cphithree$}. 
\end{cases}\]
One checks by direct inspection that (i) $hf = \alpha$ and $hg_3 = \beta$ and that (ii) $h$ is the only such map.
\end{proof}

As we mentioned just before Def. \ref{def:psione-general}, the decomposition \eqref{general-phi-is-phione-phitwo} applies to $\psione \in \uwd\yz$ defined in \eqref{uwd-psi1}.  In the next two sections, we will show that, up to name changes,  $\phi_1$ is generated by loops (Prop. \ref{prop:general-iterated-loops}), and $\phi_2$ is generated by splits (Prop. \ref{prop:general-iterated-splits}).

\section{Iterated Splits}
\label{sec:iterated-splits}

The purpose of this section is to show that $\phi_2$ \eqref{phitwo-general} is either a name change or is generated by splits.  First let us adopt the following convention, which is the undirected version of Convention \ref{conv:ignore-name-change}.

\begin{convention}
\label{conv:uwd-ignore-name-change}
Using the three elementary relations \eqref{uwd-move-a3}, \eqref{uwd-move-a4}, and \eqref{uwd-move-a5}, name changes (Def. \ref{def:uwd-gen-namechange}) can always be rewritten on the outside (i.e., left side) of an iterated operadic composition in $\uwd$.  Moreover, using the elementary relation \eqref{uwd-move-a1}, an iteration of name changes can be composed down into just one name change.  To simplify the presentation, in what follows these elementary relations regarding name changes are automatically applied wherever necessary.  With this in mind, in the sequel we will mostly not mention name changes.
\end{convention}

Recall from Remark \ref{rk:phi1-phi2} that in $\phi_2$ \eqref{phitwo-general}, the input soldering function is the identity function and the output soldering function is surjective.  So the following Proposition applies to $\phi_2$.

\begin{motivation}
The following result says that an undirected wiring diagram of the form
\begin{center}
\begin{tikzpicture}[scale=.6]
\draw [ultra thick] (-1,0) rectangle (5,3);
\draw [ultra thick] (1,.5) rectangle (3,2.5);
\node at (2,1.5) {$A$};
\node at (2.7,2.2) {\tiny{$a_1$}};
\node at (2.7,.8) {\tiny{$a_n$}};
\draw [thick] (1,.7) -- (-2,.7);
\draw [thick] (1,2.3) -- (-2,2.3);
\draw [thick] (3,2.2) -- (4,2.2);
\draw [thick] (4,2.2) to [out=45,in=180] (6,2.7);
\node at (5.5,2.2) {\tiny{$\vdots$}}; 
\draw [thick] (4,2.2) to [out=-45,in=180] (6,1.7);
\node at (6.5,2.3) {\tiny{$\ginv a_1$}};
\draw [thick] (3,.8) -- (4,.8);
\draw [thick] (4,.8) to [out=45,in=180] (6,1.3);
\draw [thick] (4,.8) to [out=-45,in=180] (6,.3);
\node at (5.5,.8) {\tiny{$\vdots$}}; 
\node at (6.5,.8) {\tiny{$\ginv a_n$}};
\cable{(0,.7)} \node at (0,1.5) {{$\vdots$}}; \cable{(0,2.3)}
\cable{(4,.8)} \node at (4,1.5) {{$\vdots$}};  \cable {(4,2.2)}
\end{tikzpicture}
\end{center}
can be generated by splits.
\end{motivation}

\begin{proposition}
\label{prop:general-iterated-splits}
Suppose $A, B \in \Fins$, and
\begin{equation}
\label{rho-iterated-splits}
\rho = \Bigl(\nicexy@C+.5cm{A \ar[r]^-{\frho}_-{=} & A & B \ar[l]_-{\grho}^-{\mathrm{surjective}}} \Bigr) \in \uwd\capba
\end{equation}
with $\frho = \Id_A$ and $\grho$ surjective.
\begin{enumerate}
\item
If $A = \varnothing$, then $\rho = \tensorunit_\varnothing\in \uwd\emptyprof$ (Def. \ref{def:uwd-units}).
\item
Suppose $A \not= \varnothing$.
\begin{enumerate}[(i)]
\item
If $\grho$ is a bijection, then $\rho$ is a name change $\tau_{A,B}$.
\item
Otherwise, $\rho$ is an iterated operadic composition of splits (Def. \ref{def:uwd-gen-split}).\index{iterated splits}
\end{enumerate}
\end{enumerate}
\end{proposition}

\begin{proof}
We will write $\frho$ and $\grho$ as $f$ and $g$, respectively.  If $A = \varnothing$, then $\rho$ is the cospan $(\nicexy{\varnothing \ar[r] & \varnothing & \varnothing \ar[l]})$, which is the $\varnothing$-colored unit in $\uwd$.   Suppose $A \not= \varnothing$.  If $g$ is a bijection, then by definition $g$ is the name change $\tau_{\ginv}$.  

So suppose $g$ is surjective but is not a bijection.  We must show that $\rho$ is an iterated operadic composition of splits.  The first step is to decompose $\rho$ in such a way that each constituent undirected wiring diagram creates one group of output wires $\ginv a_i$.  Decompose $A$ as $A = A_1 \amalg A_2 \in \Fins$ in which
\begin{itemize}
\item
$A_1 = \Bigl\{a \in A : \left|\ginv a \right| = 1\Bigr\}$;
\item
$A_2 = \Bigl\{a \in A : \left|\ginv a \right| \geq 2\Bigr\} = \{a_1,\ldots,a_n\}$.
\end{itemize}
By assumption $A_2 \not= \varnothing$.  To decompose $\rho$ we will use the following intermediate boxes.  For each $1 \leq i \leq n+1$, define
\[
D_i = A_1 \amalg \underbrace{\coprod_{1 \leq k < i} \ginv a_k}_{\text{$\varnothing$ if $i=1$}} \amalg \underbrace{\{a_i, \ldots, a_n\}}_{\text{$\varnothing$ if $i=n+1$}} \in \Fins.\]
Note that $D_1 = A$ and $D_{n+1} \cong B$.

For $1 \leq i \leq n$ define
\begin{equation}
\label{rho-sub-i}
\rho_i = \Bigl(\nicexy{D_i \ar[r]^-{f_i}_-{=} & D_i & D_{i+1} \ar[l]_-{g_i}}\Bigr) \in \uwd\diplusonedi
\end{equation}
in which, for $d \in D_{i+1}$,
\[
g_i(d) = \begin{cases}
a_i & \text{if $d \in \ginv a_i$},\\
d & \text{otherwise}.\end{cases}\]
A direct inspection using Example \ref{ex:pushout-id} and Prop. \ref{uwd-compi-explicit} shows that, up to a name change, there is a decomposition
\begin{equation}
\label{rho-is-rho1-rhon}
\rho = \rho_n \comp \cdots \comp \rho_1 \in \uwd\capba.
\end{equation}
When $n=2$, this decomposition is depicted in the following picture.
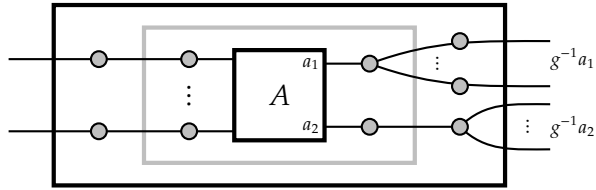
\begin{center}
\begin{tikzpicture}[scale=.6]
\draw [ultra thick] (-3,-.5) rectangle (7,3.5);
\draw [ultra thick, lightgray] (-1,0) rectangle (5,3);
\draw [ultra thick] (1,.5) rectangle (3,2.5);
\node at (2,1.5) {$A$};
\node at (2.7,2.2) {\tiny{$a_1$}};
\node at (2.7,.8) {\tiny{$a_2$}};
\draw [thick] (1,.7) -- (-4,.7);
\draw [thick] (1,2.3) -- (-4,2.3);
\draw [thick] (3,2.2) -- (4,2.2);
\draw [thick] (4,2.2) to [out=20,in=180] (8,2.7);
\node at (5.5,2.2) {\tiny{$\vdots$}}; 
\draw [thick] (4,2.2) to [out=-20,in=180] (8,1.7);
\node at (8.5,2.3) {\tiny{$\ginv a_1$}};
\draw [thick] (3,.8) -- (6,.8);
\draw [thick] (6,.8) to [out=45,in=180] (8,1.3);
\draw [thick] (6,.8) to [out=-45,in=180] (8,.3);
\node at (7.5,.8) {\tiny{$\vdots$}}; 
\node at (8.5,.8) {\tiny{$\ginv a_2$}};
\cable{(0,.7)} \node at (0,1.5) {{$\vdots$}}; \cable{(0,2.3)}
\cable{(-2,.7)} \cable{(-2,2.3)}
\cable{(4,2.2)} \cable{(6,2.7)} \cable{(6,1.7)}
\cable{(4,.8)} \cable{(6,.8)} 
\end{tikzpicture}
\end{center}
To finish the proof, it suffices to show that each $\rho_i$ is an iterated operadic composition of splits.  We will prove this assertion in the next result.
\end{proof}

\begin{motivation}
The following result says that an undirected wiring diagram of the form
\begin{center}
\begin{tikzpicture}[scale=.6]
\draw [ultra thick] (-1,0) rectangle (5,3);
\draw [ultra thick] (1,.5) rectangle (3,2.5);
\node at (2.7,1.5) {\tiny{$b$}};
\draw [thick] (1,.7) -- (-2,.7);
\draw [thick] (1,2.3) -- (-2,2.3);
\draw [thick] (3,1.5) -- (4,1.5);
\draw [thick] (4,1.5) to [out=45,in=180] (6,2.5);
\node at (5.5,1.5) {$\vdots$}; 
\draw [thick] (4,1.5) to [out=-45,in=180] (6,.5);
\node at (6.3,2.5) {\tiny{$b_1$}};
\node at (6.3,.5) {\tiny{$b_p$}};
\cable{(0,.7)} \node at (0,1.5) {{$\vdots$}}; \cable{(0,2.3)}
\cable{(4,1.5)}
\end{tikzpicture}
\end{center}
can be generated by splits.
\end{motivation}

\begin{proposition}
\label{prop:iterated-splits}
Suppose $D \in \Fins$, $D \not\ni b \in S$, and
\begin{equation}
\label{iterated-splits-pi}
\pi = \Bigl(\nicexy{D \amalg b \ar[r]^-{=} & D \amalg b & D \amalg \ginv b \ar[l]_-{g}}\Bigr) \in \uwd\dplusginvbdplusb
\end{equation}
such that
\begin{itemize}
\item
$g|_D = \Id_D$ and
\item
$p = \left|\ginv b\right| \geq 2$.  
\end{itemize}
Then $\pi$ is an iterated operadic composition of $p-1$ splits.
\end{proposition}

\begin{proof}
Write $\ginv b = \{b_1,\dots,b_p\}$.  If $p=2$ then $\pi$ is the split $\sigma^{\left(D \amalg \{b_1,b_2\}, b_1,b_2\right)}$ by definition.  

Suppose $p \geq 3$.  Then there is a decomposition of $\pi$ into $p-1$ splits as
\begin{equation}
\label{pi-decomposes-splits}
\pi = \sigma^{\left(D \amalg \{b_1,\ldots,b_p\}, b_{p-1}, b_p\right)} \comp \cdots \comp 
\sigma^{\left(D \amalg \{b_1, b_{[2,p]}\}, b_1, b_{[2,p]}\right)}.
\end{equation}
Here each $b_{[j,p]}$ means the wires $b_l$ for $j \leq l \leq p$ are identified into one element.  Starting from the right, the $j$th split in the above decomposition of $\pi$, namely
\[
\sigma^{\left(D \amalg \left\{b_1,\ldots,b_j, b_{[j+1,p]}\right\}, b_j, b_{[j+1,p]}\right)} \in \uwd\smallbinom{D \amalg \{b_1,\ldots,b_j, b_{[j+1,p]}\}}{D \amalg \{b_1,\ldots,b_{j-1}, b_{[j,p]}\}}, \]
is the cospan
\[\nicexy{
& D \amalg \bigl\{b_1,\ldots,b_j, b_{[j+1,p]}\bigr\} \ar[d]\\
D \amalg \bigl\{b_1,\ldots,b_{j-1}, b_{[j,p]}\bigr\} \ar[r]^-{=} 
& D \amalg \bigl\{b_1,\ldots,b_{j-1}, b_{[j,p]}\bigr\}}\]
in $\Fins$.  Here the output soldering function sends $b_j$ and $b_{[j+1,p]}$ to $b_{[j,p]}$ and is the identity function everywhere else.
\end{proof}

\begin{example}
For each $1 \leq i \leq n$ the undirected wiring diagram $\rho_i \in \uwd\diplusonedi$ in \eqref{rho-sub-i} is of the form $\pi$ \eqref{iterated-splits-pi} with 
\[D = A_1 \amalg \coprod_{1 \leq k < i} \ginv a_k \amalg \{a_{i+1}, \ldots, a_n\}
\andspace b = a_i.\]
Therefore, $\rho$ in \eqref{rho-is-rho1-rhon} is an iterated operadic composition of splits.
\end{example}

\begin{example}
In \eqref{pi-decomposes-splits} above:
\begin{enumerate}
\item
If $p=3$, then $\pi$ decomposes into two splits as
\[\pi = \sigma^{\left(D \amalg \{b_1,b_2,b_3\}, b_2, b_3\right)} \comp
\sigma^{\left(D \amalg \{b_1, b_{[2,3]}\}, b_1, b_{[2,3]}\right)}.\]
\item
If $p=4$, then $\pi$ decomposes into three splits as
\[\pi = \sigma^{\left(D \amalg \{b_1,b_2,b_3,b_4\}, b_3, b_4\right)} \comp
\sigma^{\left(D \amalg \{b_1,b_2,b_{[3,4]}\}, b_2, b_{[3,4]}\right)} \comp
\sigma^{\left(D \amalg \{b_1, b_{[2,4]}\}, b_1, b_{[2,4]}\right)}.\]
\end{enumerate}
\end{example}

\section{Iterated Loops}
\label{sec:iterated-loops}

The purpose of this section is to show that $\phi_1$ \eqref{phione-general} is either a name change or is generated by loops (Def. \ref{def:uwd-gen-loop}).  Recall Convention \ref{conv:uwd-ignore-name-change} regarding name changes.  Also recall from Remark \ref{rk:phi1-phi2} that in $\phi_1$ \eqref{phione-general} each cable is either a $(1,1)$-cable or a $(2,0)$-cable.  Therefore, the following result applies to $\phi_1$.

\begin{motivation}
The following result says that an undirected wiring diagram of the form
\begin{center}
\begin{tikzpicture}[scale=.6]
\draw [ultra thick] (-1,0) rectangle (7,3.6);
\draw [ultra thick] (1,.5) rectangle (6,2.5);
\draw [thick] (1,.7) -- (-2,.7); \cable{(0,.7)} 
\node at (0,1.5) {{$\vdots$}}; 
\draw [thick] (1,2.3) -- (-2,2.3); \cable{(0,2.3)}
\draw [thick] (2,2.5) to [out=90, in=180] (2.5,3);
\draw [thick, out=0, in=90] (2.5,3) to (3,2.5);
\cable{(2.5,3)}
\node at (3.5,3) {$\cdots$};
\draw [thick] (4,2.5) to [out=90, in=180] (4.5,3);
\draw [thick, out=0, in=90] (4.5,3) to (5,2.5);
\cable{(4.5,3)}
\end{tikzpicture}
\end{center}
can be generated by loops.
\end{motivation}

\begin{proposition}
\label{prop:general-iterated-loops}
Suppose $A,B,C \in \Fins$ and
\[
\xi = \Bigl(\nicexy@C+.5cm{A \ar[r]^-{f} & C & B \ar[l]_-{g}^-{\mathrm{inclusion}}} \Bigr) \in \uwd\capba\]
in which each cable is either a $(1,1)$-cable or a $(2,0)$-cable.  Suppose $\xi$ has $q$ $(2,0)$-cables.  Then:
\begin{enumerate}
\item
$\xi$ is a name change if $q=0$.
\item
$\xi$ is the iterated operadic composition of $q$ loops (Def. \ref{def:uwd-gen-loop}) if $q \geq 1$.\index{iterated loops}
\end{enumerate}
\end{proposition}

\begin{proof}
Since $\xi$ only has $(1,1)$-cables and $(2,0)$-cables, up to a name change we may write it as the cospan
\[
\xi = \Bigl(\nicexy@C+.5cm{B \amalg T_{\pm} \ar[r]^-{f}_-{(\Id, \mathrm{fold})} & B \amalg T & B \ar[l]_-{g}^-{\mathrm{inclusion}}} \Bigr)\]
with:
\begin{itemize}
\item
$T$ the set of $(2,0)$-cables in $\xi$;
\item
$\finv T = T_{\pm} = T_+ \amalg T_-$ the coproduct of two copies of $T$.
\end{itemize}
If $q = 0$ (i.e., $T=\varnothing$), then $\xi$ is the $B$-colored unit (Def. \ref{def:uwd-units}).  

If $q = 1$ with $T = \{t\}$ and $T_{\pm} = \{t_{\pm}\}$, then $\xi$ is the loop $\lambda_{(B \amalg \{t_{\pm}\}, t_{\pm})}$ by definition.

Suppose $q \geq 2$.  We may write $T = \{t^1,\ldots,t^q\}$ and $T_{\pm} = \left\{t^1_{\pm}, \ldots, t^q_{\pm}\right\}$.  The picture
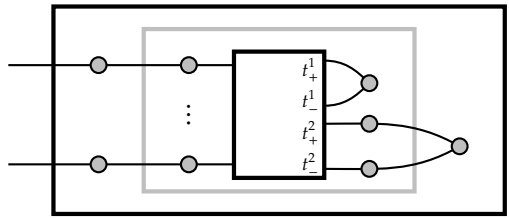
\begin{center}
\begin{tikzpicture}[scale=.6]
\draw [ultra thick] (-3,-1) rectangle (7,3.6);
\draw [ultra thick, lightgray] (-1,-.5) rectangle (5,3.1);
\draw [ultra thick] (1,-.2) rectangle (3,2.6);
\node at (2.7,2.2) {\tiny{$t^1_+$}};
\node at (2.7,1.5) {\tiny{$t^1_-$}};
\node at (2.7,.8) {\tiny{$t^2_+$}};
\node at (2.7,.1) {\tiny{$t^2_-$}};
\draw [thick] (1,.1) -- (-4,.1);
\draw [thick] (1,2.3) -- (-4,2.3);
\draw [thick] (3,2.4) to [out=0,in=135] (4,1.9);
\draw [thick] (4,1.9) to [out=225,in=0] (3,1.4);
\draw [thick] (3,1) to [out=0,in=150] (6,.5);
\draw [thick] (6,.5) to [out=210,in=0] (3,0);
\cable{(0,.1)} \node at (0,1.2) {{$\vdots$}}; \cable{(0,2.3)}
\cable{(-2,.1)} \cable{(-2,2.3)}
\cable{(4,1.9)} 
\cable{(4,0)} \cable{(4,1)} \cable{(6,.5)}
\end{tikzpicture}
\end{center}
depicts a decomposition of $\xi$ into two loops when $q = 2$.  A direct inspection shows that there is a decomposition of $\xi$ into $q$ loops as
\begin{equation}
\label{xi-decomposes-loops}
\xi = \lambda_{\left(B \amalg t^q_{\pm}, t^q_{\pm}\right)} \comp \cdots \comp
\lambda_{\left(B \amalg \{t^1_{\pm},\ldots,t^q_{\pm}\}, t^1_{\pm}\right)}.
\end{equation}
Starting from the right, the $j$th loop in the above decomposition of $\xi$, namely
\[
\lambda_{\left(B \amalg \{t^j_{\pm}, \ldots, t^q_{\pm}\}, t^j_{\pm}\right)} \in \uwd \smallbinom{B \amalg \{t^{j+1}_{\pm}, \ldots, t^q_{\pm}\}}{B \amalg \{t^{j}_{\pm}, \ldots, t^q_{\pm}\}},\]
is the cospan
\[
\nicexy@C+.5cm{& B \amalg \left\{t^{j+1}_{\pm}, \ldots, t^q_{\pm}\right\} \ar[d]^-{\mathrm{inclusion}}\\
B \amalg \left\{t^{j}_{\pm}, \ldots, t^q_{\pm}\right\} \ar[r]^-{t^j_{\pm} \mapsto t^j} & 
B \amalg \left\{t^{j+1}_{\pm}, \ldots, t^q_{\pm}\right\} \amalg t^j}\]
in $\Fins$.  Here the input soldering function sends $t^j_{\pm}$ to $t_j$ and is the identity function everywhere else.
\end{proof}

\begin{example}
In \eqref{xi-decomposes-loops} above:
\begin{enumerate}
\item
If $q=2$, then $\xi$ decomposes into two loops as
\[\xi = \lambda_{\left(B \amalg t^2_{\pm}, t^2_{\pm}\right)} \comp 
\lambda_{\left(B \amalg \{t^1_{\pm}, t^2_{\pm}\}, t^1_{\pm}\right)}. \]
\item
If $q=3$, then $\xi$ decomposes into three loops as
\[\xi = \lambda_{\left(B \amalg t^3_{\pm}, t^3_{\pm}\right)} \comp 
\lambda_{\left(B \amalg \{t^2_{\pm}, t^3_{\pm}\}, t^2_{\pm}\right)} \comp
\lambda_{\left(B \amalg \{t^1_{\pm}, t^2_{\pm}, t^3_{\pm}\}, t^1_{\pm}\right)}. \]
\end{enumerate}
\end{example}

\section{Summary of Chapter \ref{ch09-decomp-uwd}}

Every undirected wiring diagram $\varphi$ has a decomposition
\[\varphi = \phi_1 \comp \phi_2 \comp \varphi_2\]
in which:
\begin{itemize}
\item $\phi_1$ is generated by loops;
\item $\phi_2$ is generated by splits;
\item $\varphi_2$ is generated by $2$-cells and $1$-output wires.
\end{itemize}

\chapter{Finite Presentation for Undirected Wiring Diagrams}
\label{ch10-stratified-uwd}

Fix a class $S$, with respect to which the $\Fins$-colored operad $\uwd$ of undirected wiring diagrams is defined (Theorem \ref{uwd-operad}).  The main purpose of this chapter is to establish a finite presentation for the operad $\uwd$; see Theorem \ref{thm:uwd-generator-relation}.  This means the following two statements.
\begin{enumerate}
\item
The $6$ generating undirected wiring diagrams (Def. \ref{def:generating-uwd}) generate the operad $\uwd$.  This means that every undirected wiring diagram can be expressed as a finite iterated operadic composition involving only the $6$ generators.  See Theorem \ref{stratified-presentation-exists-uwd}.  
\item
If an undirected wiring diagram can be operadically generated by the generators in two different ways, then there exists a finite sequence of \emph{elementary equivalences} (Def. \ref{def:equivalent-simplices-uwd}) from the first iterated operadic composition to the other one.  See Theorem \ref{thm:uwd-generator-relation}.  An elementary equivalence is induced by either an elementary relation in $\uwd$ (Def. \ref{def:elementary-relation-uwd}) or an operad associativity/unity axiom for the generators.
\end{enumerate}
This finite presentation theorem for $\uwd$ is the undirected analogue of the finite presentation theorem for $\WD$ (Theorem \ref{thm:wd-generator-relation}).  As in the directed case, this result leads to a finite presentation theorem for $\uwd$-algebras, which we will discuss in Chapter \ref{ch11-uwd-algebras}.

We will continue to use Notation \ref{comp-is-compone} for (iterated) operadic compositions.

\section{Stratified Presentation}
\label{sec:stratified-uwd}

In this section, we define a \emph{stratified presentation} in $\uwd$ and show that every undirected wiring diagram has a stratified presentation (Theorem \ref{stratified-presentation-exists-uwd}).  The following definition is the undirected analogue of Def. \ref{def:simplex}.

\begin{motivation}
A simplex below is a finite parenthesized word whose alphabets are generating undirected wiring diagrams, in which each pair of parentheses has a well defined associated $\compi$-composition.  In particular, a simplex has a well defined operadic composition.  As we have seen in Chapter \ref{ch08-generating-uwd}, it is often possible to express an undirected wiring diagram as an operadic composition of generating undirected wiring diagrams in multiple ways.  In other words, an undirected wiring diagram can have many different simplex presentations.  We now start to develop the necessary language to say precisely that any two such simplex presentations of the same undirected wiring diagram are equivalent in some way.
\end{motivation}

\begin{definition}
\label{def:simplex-uwd}
Suppose $n \geq 1$.  An \emph{$n$-simplex} \index{simplex in UWD@simplex in $\uwd$} $\Psi$ and its \emph{composition} \index{composition of a simplex in UWD@composition of a simplex in $\uwd$} \label{notation:simplex-composition-uwd}$|\Psi| \in \uwd$ are defined inductively as follows.
\begin{enumerate}
\item
A \emph{$1$-simplex} is a generator (Def. \ref{def:generating-uwd}) $\psi$.  Its \emph{composition} $|\psi|$ is defined as $\psi$ itself.
\item
Suppose $n \geq 2$ and that $k$-simplices for $1 \leq k \leq n-1$ and their compositions in $\uwd$  are already defined.  An \emph{$n$-simplex} in $\uwd$ is a tuple $\Psi = \left(\upsi, i, \uphi\right)$ consisting of
\begin{itemize}
\item
an integer $i \geq 1$,
\item
a $p$-simplex $\upsi$ for some $p \geq 1$, and
\item
a $q$-simplex $\uphi$ for some $q \geq 1$
\end{itemize}
such that:
\begin{enumerate}[(i)]
\item
$p+q=n$;
\item
the operadic composition
\begin{equation}
\label{length-n-simplex-uwd}
|\Psi| \defn \left|\upsi\right| \compi \left|\uphi\right|
\end{equation}
is defined in $\uwd$ (Def. \ref{def:compi-uwd}).  
\end{enumerate}
The undirected wiring diagram $|\Psi|$ in \eqref{length-n-simplex-uwd} is the \emph{composition} of $\Psi$.
\end{enumerate}
A \emph{simplex} in $\uwd$ is an $m$-simplex in $\uwd$ for some $m \geq 1$.  We say that a simplex $\Psi$ is a \emph{presentation} \index{presentation} of the undirected wiring diagram $|\Psi|$.
\end{definition}

\begin{remark}
\label{rk:uwd-simplex-notation}
To simplify the presentation, as in Notation \ref{notation:simplex}, we will sometimes use either
\begin{itemize}
\item
the right side of \eqref{length-n-simplex-uwd} or
\item
even just the list of generators $(\psi_1,\ldots,\psi_n)$ in a simplex in the order in which they appear in \eqref{length-n-simplex-uwd}
\end{itemize} 
to denote a simplex in $\uwd$.
\end{remark}

\begin{example}
Elementary relations in $\uwd$ (Def. \ref{def:elementary-relation-uwd}) provide a large source of simplices in $\uwd$.  In fact, each side, either left or right, of each elementary relation in $\uwd$ is a simplex.  For example:
\begin{enumerate}
\item
The elementary relation \eqref{uwd-move-b1} says that the $3$-simplex
\[\lambda_{(Y,x,y)}  \comp \Bigl(\theta_{(X,y)} \comptwo \omega_{y} \Bigr)\]
and the $2$-simplex
\[\lambda_{(Y,x,y)}  \comp \sigma^{(Y,x,y)}\]
have the same composition.
\item
The elementary relation \eqref{uwd-move-b2} says that the $4$-simplex
\[\lambda_{(W,w,y)} \comp \sigma^{(W,x,w)} \comp 
\Bigl(\theta_{(X,y)} \comptwo \omega_{y} \Bigr) \]
and the $1$-simplex $(\tensorunit_X)$ have the same composition.  In other words, the former has composition $\tensorunit_X$.
\end{enumerate}
\end{example}

\begin{example}
\label{ex:wasted-cable-simplex}
In \eqref{wasted-cable-simplex} we considered the $5$-simplex
\[\lambda_{(X,x_1,x_2)} \comp \left[\Bigl(\left(\theta_{(\varnothing, X)} \comptwo \theta_{(x_1,x_2)}\right) \comptwo \omega_{x_1}\Bigr) \comptwo \omega_{x_2}\right],\]
whose composition has a wasted cable.
\end{example}

\begin{example}
\label{ex:wasted-cable-y-simplex}
In \eqref{wasted-cable-y-simplex} we considered the $5$-simplex
\[\lambda_{(Y \amalg X,x_1,x_2)} \comp \left[\Bigl(\left(\theta_{(Y, X)} \comptwo \theta_{(x_1,x_2)}\right) \comptwo \omega_{x_1}\Bigr) \comptwo \omega_{x_2}\right],\]
whose composition also has a wasted cable.
\end{example}

\begin{motivation}
If we think of a simplex as a parenthesized word whose alphabets are generating undirected wiring diagrams, then the stratified simplex in the next definition is a word where the same alphabets must occur in a consecutive string.  For example, all the loops must occur together as a string $(\lambda_1,\ldots, \lambda_n)$.  Furthermore, we can even insist that these strings for different types of generating undirected wiring diagrams occur in a specific order, with name changes and loops at the top, followed by splits, $2$-cells, and $1$-output wires at the bottom.
\end{motivation}

\begin{definition}
\label{def:stratified-presentation-uwd}
A \emph{stratified simplex} \index{stratified simplex in UWD@stratified simplex in $\uwd$} in $\uwd$ is a simplex in $\uwd$ (Def. \ref{def:simplex-uwd}) of one of the following two forms:
\begin{enumerate}
\item
$\left(\epsilon\right)$, where $\epsilon \in \uwd\emptynothing$ is the empty cell (Def. \ref{def:uwd-gen-empty}).
\item
$\bigl(\tau, \ulambda, \usigma, \utheta, \uomega\bigr)$, where:
\begin{itemize}
\item
$\tau$ is a name change (Def. \ref{def:uwd-gen-namechange});
\item
$\ulambda$ is a possibly empty string of loops (Def. \ref{def:uwd-gen-loop});
\item
$\usigma$ is a possibly empty string of splits (Def. \ref{def:uwd-gen-split});
\item
$\utheta$ is a possibly empty string of $2$-cells (Def. \ref{def:uwd-gen-2cell});
\item
$\uomega$ is a possibly empty string of $1$-output wires (Def. \ref{def:uwd-gen-1output}).
\end{itemize}
\end{enumerate}
We call these stratified simplices of type (1) and of type (2), respectively.  If $\Psi$ is a stratified simplex in $\uwd$, then we call it a \emph{stratified presentation} \index{stratified presentation} of the undirected wiring diagram $|\Psi|$.
\end{definition}

\begin{remark}
\label{rk:stratified-disjoint-uwd}
The composition of a stratified simplex of type (2) cannot be the empty cell $\epsilon$.  Indeed, if the composition of a stratified simplex of type (2) is in $\uwd\emptynothing$, then $1$-output wires must be involved.  So the composition must have at least one cable, and it cannot be the empty cell.
\end{remark}

\begin{example}
In \eqref{wasted-cable-y-simplex} and \eqref{wasted-cable-y-simplex-b} we gave two stratified presentations of the undirected wiring diagram $\zeta_Y \in \uwd\yy$ in Example \ref{ex:y-plus-wasted-cable}.
\end{example}

\begin{example}
In \eqref{sample-uwd-stratified} we gave a stratified presentation of the undirected wiring diagram $\varphi \in \uwd\yxonextwo$.
\end{example}

We now observe that the generators generate the operad $\uwd$ of undirected wiring diagrams in a highly structured way.

\begin{theorem}
\label{stratified-presentation-exists-uwd}
Every undirected wiring diagram admits a stratified presentation.\index{stratified presentations exist in UWD@stratified presentations exist in $\uwd$}
\end{theorem}

\begin{proof}
Suppose $\psi$ is an undirected wiring diagram.  If $\psi$ is the empty cell $\epsilon$ (Def. \ref{def:uwd-gen-empty}), then $(\epsilon)$ is a stratified presentation of $\psi$.  So let us now assume that $\psi$ is not the empty cell.  We will show that it admits a stratified presentation of type (2).  We remind the reader of Convention \ref{conv:uwd-ignore-name-change}.

Combining \eqref{psi-is-psi1psi2} and \eqref{general-phi-is-phione-phitwo}, there is a decomposition 
\[\psi = \phi_1 \comp \phi_2 \comp \psitwo.\]
In this decomposition:
\begin{itemize}
\item
The undirected wiring diagram $\psi_2$ either is a colored unit or has a stratified presentation $(\utheta, \uomega)$ by Lemma \ref{uwd-psitwo-marginal} and Corollary \ref{psitwo-is-2cells-1outputwires}.  
\item
The undirected wiring diagram $\phi_2$ either is a name change or has a stratified presentation $(\usigma)$ by Prop. \ref{prop:general-iterated-splits}.
\item
The undirected wiring diagram $\phi_1$ either is a name change or has a stratified presentation $(\ulambda)$ by Prop. \ref{prop:general-iterated-loops}.
\end{itemize}
By the above three statements, $\psi$ has a stratified presentation of type (2).
\end{proof}

\section{Elementary Equivalences}
\label{sec:uwd-elementary-eq}

The purpose of this section is to establish the second part of the finite presentation theorem for the operad $\uwd$ of undirected wiring diagrams.  Recall the $17$ elementary relations in $\uwd$ (Def. \ref{def:elementary-relation-uwd}).  In what follows, we will regard each operad associativity or unity axiom as an equality.  The following definition is the undirected analogue of Def. \ref{def:equivalent-simplices}.

\begin{definition}
\label{def:equivalent-simplices-uwd}
Suppose $\Psi$ is an $n$-simplex in $\uwd$ as in Def. \ref{def:simplex-uwd}.
\begin{enumerate}
\item
A \emph{subsimplex} \index{subsimplex in UWD@subsimplex in $\uwd$} of $\Psi$ is a simplex in $\uwd$ defined inductively as follows.
\begin{itemize}
\item
If $\Psi$ is a $1$-simplex, then a \emph{subsimplex} of $\Psi$ is $\Psi$ itself.
\item
Suppose $n \geq 2$ and $\Psi =  \left(\upsi, i, \uphi\right)$ for some $i \geq 1$, $p$-simplex $\upsi$, and $q$-simplex $\uphi$ with $p+q=n$.  Then a \emph{subsimplex} of $\Psi$ is 
\begin{itemize}
\item
a subsimplex of $\upsi$, 
\item
a subsimplex of $\uphi$, or
\item
$\Psi$ itself.
\end{itemize}
\end{itemize}
If $\Psi'$ is a subsimplex of $\Psi$, then we write $\Psi' \subseteq \Psi$.
\item
An \emph{elementary subsimplex} \index{elementary subsimplex in UWD@elementary subsimplex in $\uwd$} $\Psi'$ of $\Psi$ is a subsimplex of one of two forms: 
\begin{enumerate}[(i)]
\item
$\Psi'$ is one side (either left or right) of a specified elementary relation in $\uwd$ (Def. \ref{def:elementary-relation-uwd}).
\item
$\Psi'$ is one side (either left or right) of a specified operad associativity or unity axiom--\eqref{compi-associativity}, \eqref{compi-associativity-two}, \eqref{compi-left-unity}, or \eqref{compi-right-unity}--involving only the generators in $\uwd$ (Def. \ref{def:generating-uwd}). 
\end{enumerate}
\item
Suppose $\Phi$ is another simplex in $\uwd$.  Then $\Psi$ and $\Phi$ are said to be \emph{equivalent} \index{equivalent simplices in UWD@equivalent simplices in $\uwd$} if their compositions are equal; i.e., $|\Psi| = |\Phi| \in \uwd$.
\item
Suppose:
\begin{itemize}
\item
$\Psi' \subseteq \Psi$ is an elementary subsimplex corresponding to one side of $R$, which is either an elementary relation or an operad associativity/unity axiom for the generators in $\uwd$.
\item
$\Psi"$ is the simplex given by the other side of $R$.
\item
$\Psi^1$ is the simplex obtained from $\Psi$ by replacing the subsimplex $\Psi'$ by $\Psi"$.  
\end{itemize}
We say that $\Psi$ and $\Psi^1$ are \emph{elementarily equivalent}.\index{elementarily equivalent in UWD@elementarily equivalent in $\uwd$}  Note that elementarily equivalent simplices are also equivalent.
\item
If $\Psi$ and $\Phi$ are elementarily equivalent, we write $\Psi \sim \Phi$ and call this an \emph{elementary equivalence}.\index{elementary equivalence in UWD@elementary equivalence in $\uwd$}
\item
Suppose $\Psi_0,\ldots,\Psi_r$ are simplices in $\uwd$ for some $r \geq 1$ and that there exist elementary equivalences
\[\Psi_0 \sim \Psi_1 \sim \cdots \sim \Psi_r.\]
Then we say that $\Psi_0$ and $\Psi_r$ are \emph{connected by a finite sequence of elementary equivalences}.  Note that in this case $\Psi_0$ and $\Psi_r$ are equivalent.
\end{enumerate}
\end{definition}

\begin{example}
By definition, for each elementary relation in $\uwd$ (Def. \ref{def:elementary-relation-uwd}), the simplices given by its two sides are elementarily equivalent.  For example:
\begin{enumerate}
\item
By the elementary relation \eqref{uwd-move-b1}, replacing an elementary subsimplex of the form
\[\lambda_{(Y,x,y)}  \comp \Bigl(\theta_{(X,y)} \comptwo \omega_{y} \Bigr)\]
by one of the form
\[\lambda_{(Y,x,y)}  \comp \sigma^{(Y,x,y)}\]
yields an elementary equivalence.
\item
By the elementary relation \eqref{uwd-move-b2}, replacing an elementary subsimplex of the form
\[\lambda_{(W,w,y)} \comp \sigma^{(W,x,w)} \comp 
\Bigl(\theta_{(X,y)} \comptwo \omega_{y} \Bigr)\]
by the colored unit $\tensorunit_X$ yields an elementary equivalence.
\end{enumerate}
\end{example}

\begin{example}
In Example \ref{ex:y-wasted-two-ways}, we observed that for the undirected wiring diagram $\zeta_Y \in \uwd\yy$ in Example \ref{ex:y-plus-wasted-cable}, the stratified presentations \eqref{wasted-cable-y-simplex} and \eqref{wasted-cable-y-simplex-b} are connected by a finite sequence of elementary equivalences. 
\end{example}

\begin{convention}
\label{conv:operad-axiom-move-uwd}
As in Convention \ref{conv:operad-axiom-move}, to simplify the presentation, elementary equivalences corresponding to an operad associativity/unity axiom--\eqref{compi-associativity}, \eqref{compi-associativity-two}, \eqref{compi-left-unity}, or \eqref{compi-right-unity}--for the generators in $\uwd$ (Def. \ref{def:generating-uwd}) will often be applied tacitly wherever necessary.
\end{convention}

The goal of this section is to show that any two equivalent simplices in $\uwd$ are connected by a finite sequence of elementary equivalences.   The first step is to show that every simplex in $\uwd$ is connected to a stratified simplex (Def. \ref{def:stratified-presentation-uwd}) in the following sense.  The following observation is the undirected analogue of Lemma \ref{lemma:simplex-to-stratified}.

\begin{lemma}
\label{lemma:simplex-to-stratified-uwd}
Every simplex in $\uwd$ is either a stratified simplex (Def. \ref{def:stratified-presentation-uwd}) or is connected to an equivalent stratified simplex by a finite sequence of elementary equivalences (Def. \ref{def:equivalent-simplices-uwd}).
\end{lemma}

\begin{proof}
Suppose $\Psi$ is a simplex in $\uwd$ that is not a stratified simplex.  Using the three elementary relations \eqref{uwd-move-a3}--\eqref{uwd-move-a5}, first we  move all the name changes (Def. \ref{def:name-change}) in $\Psi$, if there are any, to the left.  Then we use the elementary relation \eqref{uwd-move-a1} repeatedly to compose them down into one name change.  Therefore, after a finite sequence of elementary equivalences, we may assume that there is at most one name change in $\Psi$, which is the left-most entry.  If there are further elementary equivalences later that create name changes, we will perform the same procedure without explicitly mentioning it.

The empty cell $\epsilon \in \uwd\emptynothing$ (Def. \ref{def:uwd-gen-empty}) and a $1$-output wire $\omega_* \in \uwd\starnothing$ (Def. \ref{def:uwd-gen-1output}) have no input boxes, so no operadic composition of the forms $\epsilon \compi -$ or  $\omega_* \compi -$ can be defined.  Therefore, after a finite sequence of elementary equivalences corresponding to the operad horizontal associativity axiom  \eqref{compi-associativity}, we may assume that $\Psi$ has the form
\[\Bigl(\tau,\Psi^1,\uepsilon, \uomega\Bigr).\]
Here:
\begin{itemize}
\item
$\tau$ is a name change;
\item
all the $1$-output wires $\uomega$ are at the right-most entries;
\item
all the empty cells $\uepsilon$ are just to their left;
\item
$\Psi^1$ is either empty or is a subsimplex involving $2$-cells  (Def. \ref{def:uwd-gen-2cell}), loops (Def. \ref{def:uwd-gen-loop}), and splits (Def. \ref{def:uwd-gen-split}).
\end{itemize}

Next we use the elementary relations \eqref{uwd-move-c2}--\eqref{uwd-move-c5} to move all the $2$-cells in $\Psi$ to just the left of $\uepsilon$.  Then we use the elementary relations \eqref{uwd-move-d3} and \eqref{uwd-move-e1} to move all the remaining loops to just the right of the name change $\tau$.  So after a finite sequence of elementary equivalences, we may assume that the simplex $\Psi$ has the form
\begin{equation}
\label{almost-stratified-uwd}
\bigl(\tau, \ulambda, \usigma, \utheta, \uepsilon, \uomega\bigr).
\end{equation}
If the string $\uepsilon$ of empty cells is empty, then we are done because this is now a stratified simplex of type (2).

So suppose the string $\uepsilon$ of empty cells in \eqref{almost-stratified-uwd} is non-empty.  Using finitely many elementary equivalences corresponding to the elementary relations \eqref{uwd-move-c1}--\eqref{uwd-move-c3}, we may cancel all the unnecessary empty cells in \eqref{almost-stratified-uwd}.  If there are no empty cells left after the cancellation, then we have a stratified simplex of type (2).  

Suppose that, after the cancellation in the previous paragraph, the resulting string $\uepsilon$ of empty cells is still non-empty.  Then it must contain a single empty cell $\epsilon$, and there are no $2$-cells $\utheta$ and no $1$-output wires $\uomega$ in the resulting simplex $\Psi$.  Since the output box of $\epsilon$ is the empty box, the current simplex $\Psi$ cannot have any loops $\ulambda$ or splits $\usigma$.  Therefore, in this case the simplex \eqref{almost-stratified-uwd} has the form
\[\bigl(\tau, \epsilon\bigr).\]
Since the output box of $\epsilon$ is the empty box $\varnothing$, this name change $\tau$ must be the colored unit $\tensorunit_{\varnothing}$.  So by the operad left unity axiom \eqref{compi-left-unity}, the simplex $\bigl(\tensorunit_{\varnothing}, \epsilon\bigr)$ is elementarily equivalent to $(\epsilon)$, which is a stratified simplex of type (1).
\end{proof}

The next step is to show that equivalent stratified simplices are connected.  The following observation is the undirected analogue of Lemmas \ref{lemma:stratified-type1} and \ref{lemma:stratified-type2}.

\begin{lemma}
\label{lemma:stratified-connected-uwd}
Any two equivalent stratified simplices (Def. \ref{def:stratified-presentation-uwd}) in $\uwd$ are either equal or are connected by a finite sequence of elementary equivalences (Def. \ref{def:equivalent-simplices-uwd}).
\end{lemma}

\begin{proof}
Suppose $\Psi^1$ and $\Psi^2$ are equivalent stratified simplices in $\uwd$.  By Remark \ref{rk:stratified-disjoint-uwd} $\Psi^1$ and $\Psi^2$ are both of type (1) or both of type (2).  If they are both of type (1), then they are both equal to $(\epsilon)$ by definition.  So let us now assume that $\Psi^1$ and $\Psi^2$ are distinct but equivalent stratified simplices of type (2) in $\uwd$.  The rest of the proof is similar to that of Lemma \ref{lemma:stratified-type2} and consists of a series of reductions.  

Write $\psi \in \uwd\yux$ for the common composition of $\Psi^1$ and $\Psi^2$.  Using elementary equivalences corresponding to
\begin{itemize}
\item
the operad unity axioms \eqref{compi-left-unity} and \eqref{compi-right-unity}, 
\item
the elementary relations \eqref{uwd-move-b2} and \eqref{uwd-move-c1} regarding colored units, and 
\item
other elementary relations in $\uwd$ that move the generators around the stratified simplices, 
\end{itemize}
we may assume that there are no unnecessary generators in these stratified simplices.  Here \emph{unnecessary} refers to either a colored unit or generators whose (iterated) operadic composition is a colored unit.  If $\psi$ itself is a name change $\tau$, then at this stage both $\Psi^1$ and $\Psi^2$ are equal to the $1$-simplex $(\tau)$.  So let us now assume that $\psi$ is not a name change.

The name change $\tau^1$ in $\Psi^1$ has output box $Y$, and the same is true for the name change $\tau^2$ in $\Psi^2$.  We may actually assume that the input boxes of $\tau^1$ and of $\tau^2$ are also equal, provided that we change the output boxes of other generators in $\Psi^2$ accordingly by a name change if necessary.  Such changes correspond to elementary equivalences coming from the operad unity and associativity axioms and the elementary relation \eqref{uwd-move-a1}.  Using finitely many elementary equivalences, we may therefore assume that $\tau^1$ is equal to $\tau^2$.  So we may as well assume that there are no name changes in the two stratified simplices $\Psi^i$.  At this stage, each stratified simplex $\Psi^i$ has the form
\[\left(\ulambda^i, \usigma^i, \utheta^i, \uomega^i\right).\]

Using finitely many elementary equivalences corresponding to the elementary relation \eqref{uwd-move-b1} and other elementary relations that move the generators around, we may assume that in the simplices $\Psi^1$ and $\Psi^2$ all the $(1,0)$-cables in $\psi$ are created as in the left side of \eqref{uwd-move-b1}, i.e., as $\bigl(\lambda, \theta, \omega\bigr)$ rather than as $\bigl(\lambda, \sigma\bigr)$.  In plain language, this means that every $(1,0)$-cable in $\psi$ is created in both $\Psi^i$ by applying a loop to a $1$-output wire and some other wire, rather than by applying a loop to a split.

Similarly, using Example \ref{ex:y-wasted-two-ways} and finitely many elementary equivalences, we may further assume that every wasted cable in $\psi$ is created in both simplices $\Psi^i$ as $\bigl(\lambda, \theta,\theta, \omega, \omega\bigr)$ as in  \eqref{wasted-cable-y-simplex}, rather than as $\bigl(\lambda, \sigma, \theta, \omega\bigr)$ as in \eqref{wasted-cable-y-simplex-b}.  In plain language, this means that every wasted cable in $\psi$ is created in both $\Psi^i$ by applying a loop to two $1$-output wires, as in the first picture after \eqref{wasted-cable-y-simplex}.

Recall that each $1$-output wire $\omega$ is a $0$-ary element in $\uwd$, while loops and splits are unary.  Therefore, by the operad associativity axioms \eqref{compi-associativity} and \eqref{compi-associativity-two}, we may assume that in each simplex $\Psi^i$, the right portion $\bigl(\utheta^i, \uomega^i\bigr)$ is a subsimplex.  By the above reductions, at this stage the two simplices $\bigl(\utheta^1, \uomega^1\bigr)$ and $\bigl(\utheta^2, \uomega^2\bigr)$ are uniquely determined by $\psi$ and have the same composition.  In fact, each simplex $\bigl(\utheta^i, \uomega^i\bigr)$ has composition $\psitwo$ \eqref{uwd-psi2}.  Using finitely many elementary equivalences corresponding to the elementary relations \eqref{uwd-move-c2} and \eqref{uwd-move-c3} and the operad associativity axioms, we may now assume that the simplices $\bigl(\utheta^1, \uomega^1\bigr)$ and $\bigl(\utheta^2, \uomega^2\bigr)$ are equal.

Similarly, at this stage the simplices $\bigl(\ulambda^1,\usigma^1\bigr)$ and $\bigl(\ulambda^2,\usigma^2\bigr)$ both have composition $\psione$ \eqref{uwd-psi1}.  Using finitely many elementary equivalences corresponding to the elementary relations \eqref{uwd-move-d1}--\eqref{uwd-move-e1} and the operad vertical associativity axiom \eqref{compi-associativity-two}, we may now assume that the simplices $\bigl(\ulambda^1,\usigma^1\bigr)$ and $\bigl(\ulambda^2,\usigma^2\bigr)$  are equal.
\end{proof}

We are now ready for the finite presentation theorem for undirected wiring diagrams.  It describes the operad $\uwd$ of undirected wiring diagrams (Theorem \ref{uwd-operad}) in terms of finitely many generators and finitely many relations.

\begin{theorem}
\label{thm:uwd-generator-relation}
Consider the operad $\uwd$ of undirected wiring diagrams.\index{finite presentation for uwd@finite presentation for $\uwd$}
\begin{enumerate}
\item
Every undirected wiring diagram can be obtained from finitely many generators (Def. \ref{def:generating-uwd}) via iterated operadic compositions (Def. \ref{def:pseudo-operad}).
\item
Any two equivalent simplices in $\uwd$ are either equal or are connected by a finite sequence of elementary equivalences  (Def. \ref{def:equivalent-simplices-uwd}).
\end{enumerate}
\end{theorem}

\begin{proof}
The first statement is a special case of Theorem \ref{stratified-presentation-exists-uwd}.  The second statement is a combination of Lemma \ref{lemma:simplex-to-stratified-uwd} twice and Lemma \ref{lemma:stratified-connected-uwd}.
\end{proof}

\section{Summary of Chapter \ref{ch10-stratified-uwd}}

\begin{enumerate}
\item A simplex in $\uwd$ is a finite non-empty parenthesized word of generating undirected wiring diagrams in which each pair of parentheses is equipped with an operadic $\compi$-composition.
\item A stratified simplex in $\uwd$ is a simplex of one of the following two forms.
\begin{itemize}
\item $\left(\epsilon\right)$
\item
$\bigl(\tau, \ulambda, \usigma, \utheta, \uomega\bigr)$
\end{itemize}
\item Every undirected wiring diagram has a stratified presentation.
\item Two simplices in $\uwd$ are elementarily equivalent if one can be obtained from the other by replacing a subsimplex $\Psi'$ by an equivalent simplex $\Psi''$ such that $|\Psi'| = |\Psi''|$ is either one of the seventeen elementary relations in $\uwd$ or an operad associativity/unity axiom involving only the six generating undirected wiring diagrams.
\item Any two simplex presentations of a given undirected wiring diagram are connected by a finite sequence of elementary equivalences.
\end{enumerate}

\chapter{Algebras of Undirected Wiring Diagrams}
\label{ch11-uwd-algebras}

The main purpose of this chapter is to provide a finite presentation theorem for algebras over the operad $\uwd$ of undirected wiring diagrams (Theorem \ref{uwd-operad}).  As in the case of wiring diagrams (Theorem \ref{thm:wd-algebra}), this finite presentation for $\uwd$-algebras is a consequence of the finite presentation for the operad $\uwd$ (Theorem \ref{thm:uwd-generator-relation}).  This finite presentation allows us to reduce the understanding of a $\uwd$-algebra to just a few basic structure maps and a small number of easy-to-check axioms.  We will illustrate this point further with the relational algebra of a set and the typed relational algebra.

In Section \ref{sec:uwd-algebra-generators} we first define a $\uwd$-algebra in terms of $6$ generating structure maps and $17$ generating axioms corresponding to the generators (Def. \ref{def:generating-uwd}) and the elementary relations (Def. \ref{def:elementary-relation-uwd}) in $\uwd$.  Then we observe that this finite presentation for a $\uwd$-algebra is in fact equivalent to the general definition (Def. \ref{def2:operad-algebra}) of a $\uwd$-algebra (Theorem \ref{thm:uwd-algebra}).   This is an application of the finite presentation theorem for the operad $\uwd$ (Theorem \ref{thm:uwd-generator-relation}).

In Section \ref{sec:relational-algebra} we provide, for each set $A$, a finite presentation for the $\uwd$-algebra called the relational algebra of $A$ (Theorem \ref{relational-algebra-is-algebra}).  In its original form, the relational algebra was the main algebra example in \cite{spivak13}.  In \cite{spivak13} Spivak pointed out that the relational algebra of a set and its variant called the typed relational algebra, to be discussed in Section \ref{sec:relational-typed}, have applications in digital circuits, machine learning, and database theory. 

In Section \ref{sec:rigidity-relational} we prove a rigidity property of the relational algebra of a set (Theorem \ref{thm:relational-rigidity}).  It says that a given map of sets is a bijection precisely when it induces a map between the two relational algebras.  The motivation for this rigidity property comes from \cite{spivak13} Section 3, where several examples suggest that there are very few interesting maps out of the relational algebra of a set.  In fact, in \cite{spivak13} Conjecture 3.1.6, Spivak conjectured that the relational algebra of any set is quotient-free.  Although our rigidity result does not prove Spivak's conjecture in its full generality, it adds further evidence that the conjecture should be true.

In Section \ref{sec:relational-typed} we consider a generalization of the relational algebra of a set, called the typed relational algebra.  We observe that, similar to the relational algebra of a set, the typed relational algebra has a finite presentation (Theorem \ref{typed-relational-algebra}).  In its original form, the typed relational algebra was first defined in \cite{spivak13} Section 4.

\section{Finite Presentation for Algebras}
\label{sec:uwd-algebra-generators}

The main purpose of this section is to prove a finite presentation theorem for $\uwd$-algebras.  We first define a $\uwd$-algebra in terms of a finite number of generators and relations.  Immediately afterwards we will show that this definition agrees with the general definition of an operad algebra (Def. \ref{def2:operad-algebra}) when the operad is $\uwd$.

Fix a class $S$, with respect to which the $\Fins$-colored operad $\uwd$ of undirected wiring diagrams is defined (Theorem \ref{uwd-operad}).

\begin{definition}
\label{def:uwd-algebra-biased}
A \emph{$\uwd$-algebra} \index{uwd-algebra@$\uwd$-algebra} $A$ consists of the following data.  For each $X \in \Fins$, $A$ is equipped with a set $A_X$ called the \emph{$X$-colored entry} of $A$.  Moreover, it is equipped with the following $6$ \emph{generating structure maps} corresponding to the generators in $\uwd$ (Def. \ref{def:generating-uwd}).
\begin{enumerate}
\item
Corresponding to the empty cell $\epsilon \in \uwd\emptynothing$ (Def. \ref{def:uwd-gen-empty}), it has a structure map
\begin{equation}
\label{uwd-structuremap-empty}
\nicexy{\ast \ar[r]^-{\epsilon} & A_{\varnothing}},
\end{equation}
i.e., a chosen element in $A_{\varnothing}$.
\item
Corresponding to each $1$-output wire $\omega_* \in \uwd\starnothing$ (Def. \ref{def:uwd-gen-1output}), it has a structure map
\begin{equation}
\label{uwd-structuremap-1outout}
\nicexy{\ast \ar[r]^-{\omega_*} & A_{\ast}},
\end{equation}
i.e., a chosen element in $A_{\ast}$.
\item
Corresponding to each name change $\tau_f \in \uwd\yx$ (Def. \ref{def:uwd-gen-namechange}), it has a structure map
\begin{equation}
\label{uwd-structuremap-namechange}
\nicexy{A_X \ar[r]^-{\tau_f} & A_Y}
\end{equation}
that is, furthermore, the identity map if $f$ is the identity map on $X$.
\item
Corresponding to each $2$-cell $\theta_{(X,Y)} \in \uwd\xplusyxy$ (Def. \ref{def:uwd-gen-2cell}), it has a structure map
\begin{equation}
\label{uwd-structuremap-2cell}
\nicexy{A_X \times A_Y \ar[r]^-{\theta_{(X,Y)}} & A_{X \amalg Y}}. 
\end{equation}
\item
Corresponding to each loop $\lambda_{(X,\xsubpm)} \in \uwd\xminusplusminusx$ (Def. \ref{def:uwd-gen-loop}), it has a structure map
\begin{equation}
\label{uwd-structuremap-loop}
\nicexy{A_X \ar[r]^-{\lambda_{(X,\xsubpm)}} & A_{X \setminus \xsubpm}}.
\end{equation}
\item
Corresponding to each split $\sigma^{(X,x_1,x_2)} \in \uwd\xxprime$ (Def. \ref{def:uwd-gen-split}), it has a structure map
\begin{equation}
\label{uwd-structuremap-split}
\nicexy@C+.5cm{A_{X'} \ar[r]^-{\sigma^{(X,x_1,x_2)}} & A_X}.
\end{equation}
\end{enumerate}
The following $17$ diagrams, called the \emph{generating axioms}, which correspond to the elementary relations in $\uwd$ (Def. \ref{def:elementary-relation-uwd}), are required to be commutative.
\begin{enumerate}
\item
In the setting of \eqref{uwd-move-a1}, the diagram
\[\nicexy{A_X \ar[r]^-{\tau_f} \ar[dr]_-{\tau_{gf}} & A_Y \ar[d]^-{\tau_g}\\ & A_Z}\]
is commutative.
\item
In the setting of \eqref{uwd-move-a2}, the diagram
\[\nicexy{\ast \ar[r]^-{\omega_x} \ar[dr]_-{\omega_y} & A_x \ar[d]^-{\tau_{x,y}}\\
& A_y}\]
is commutative.
\item
In the setting of \eqref{uwd-move-a3}, the diagram
\[\nicexy@C+.5cm{A_{X_1} \times A_{X_2} \ar[d]_-{(\tau_{f_1}, \tau_{f_2})} \ar[r]^-{\theta_{(X_1,X_2)}} & A_{X_1 \amalg X_2} \ar[d]^-{\tau_{f_1\amalg f_2}}\\
A_{Y_1} \times A_{Y_2} \ar[r]^-{\theta_{(Y_1,Y_2)}} & A_{Y_1 \amalg Y_2}}\]
is commutative.
\item
In the setting of \eqref{uwd-move-a4}, the diagram
\[\nicexy{A_X \ar[d]_-{\tau_f} \ar[r]^-{\lambda_{(X,\xsubpm)}} & A_{X \setminus \xsubpm} \ar[d]^-{\tau_{f'}}\\
A_Y \ar[r]^-{\lambda_{(Y,\ypm)}} & A_{Y \setminus \ypm}}\]
is commutative.
\item
In the setting of \eqref{uwd-move-a5}, the diagram
\[\nicexy@C+.5cm{A_{X'} \ar[d]_-{\tau_{f'}} \ar[r]^-{\sigma^{(X,x_1,x_2)}} & A_X \ar[d]^-{\tau_f}\\
A_{Y'} \ar[r]^-{\sigma^{(Y,y_1,y_2)}} & A_Y}\]
is commutative.
\item
In the setting of \eqref{uwd-move-b1}, the diagram
\begin{equation}
\label{uwd-algebra-move-b1}
\nicexy{A_X \times \ast \cong A_X \ar[d]_-{(\Id, \omega_y)} \ar[rr]^-{\sigma^{(Y,x,y)}} && A_Y \ar[d]^-{\lambda_{(Y,x,y)}}\\
A_X \times A_y \ar[r]^-{\theta_{(X,y)}} & A_Y \ar[r]^-{\lambda_{(Y,x,y)}} & A_W}
\end{equation}
is commutative.
\item
In the setting of \eqref{uwd-move-b2}, the diagram
\[\nicexy{A_X \cong A_X \times \ast \ar[d]_-{\Id} \ar[r]^-{(\Id,\omega_y)} & A_X \times A_y \ar[r]^-{\theta_{(X,y)}} & A_Y \ar[d]^-{\sigma^{(W,x,w)}}\\
A_X && A_W \ar[ll]_-{\lambda_{(W,w,y)}}}\]
is commutative.
\item
In the setting of \eqref{uwd-move-c1}, the diagram
\[\nicexy{A_X \ar[d]_-{\Id} \ar[r]^-{\cong} & A_X \times \ast \ar[d]^-{(\Id, \epsilon)}\\
A_X & A_X \times A_\varnothing \ar[l]_-{\theta_{(X,\varnothing)}}}\]
is commutative.
\item
In the setting of \eqref{uwd-move-c2}, the diagram
\begin{equation}
\label{uwd-algebra-2cell-associativity}
\nicexy@C+.5cm{A_X \times A_Y \times A_Z \ar[r]^-{(\Id,\theta_{(Y,Z)})} \ar[d]_-{(\theta_{(X,Y)}, \Id)} & A_X \times A_{Y \amalg Z} \ar[d]^-{\theta_{(X, Y \amalg Z)}}\\
A_{X\amalg Y} \times A_Z \ar[r]^-{\theta_{(X \amalg Y, Z)}} & A_{X \amalg Y \amalg Z}}
\end{equation}
is commutative.
\item
In the setting of \eqref{uwd-move-c3}, the diagram
\begin{equation}
\label{uwd-algebra-equivariance}
\nicexy@C+.4cm{A_Y \times A_X \ar[r]^-{\mathrm{permute}} \ar[d]_-{\theta_{(Y,X)}} & A_X \times A_Y \ar[d]^-{\theta_{(X,Y)}}\\
A_{Y \amalg X} \ar[r]^-{=} & A_{X \amalg Y}}
\end{equation}
is commutative.
\item
In the setting of \eqref{uwd-move-c4}, the diagram
\[\nicexy{A_X \times A_Y \ar[d]_-{\left(\Id, \lambda_{(Y,\ypm)}\right)} \ar[r]^-{\theta_{(X,Y)}} & A_{X \amalg Y} \ar[d]^-{\lambda_{(X \amalg Y, \ypm)}}\\
A_X \times A_{Y'} \ar[r]^-{\theta_{(X,Y')}} & A_{X \amalg Y'}}\]
is commutative.
\item
In the setting of \eqref{uwd-move-c5}, the diagram
\[\nicexy{A_X \times A_{Y'} \ar[d]_-{\left(\Id, \sigma^{(Y,y_1,y_2)}\right)} \ar[r]^-{\theta_{(X,Y')}} & A_{X \amalg Y'} \ar[d]^-{\sigma^{(X\amalg Y, y_1, y_2)}}\\
A_X \times A_Y \ar[r]^-{\theta_{(X,Y)}} & A_{X \amalg Y}}\]
is commutative.
\item
In the setting of \eqref{uwd-move-d1}, the diagram
\[\nicexy@C+.5cm{A_{X'} \ar[d]_-{\sigma^{(Z,z_1,z_2)}} \ar[r]^-{\sigma^{(Y,y_1,y_2)}} & A_Y \ar[d]^-{\sigma^{(X,z_1,z_2)}}\\
A_Z \ar[r]^-{\sigma^{(X,y_1,y_2)}} & A_X}\]
is commutative.
\item
In the setting of \eqref{uwd-move-d2}, the diagram
\[\nicexy@C+.5cm{A_X \ar[d]_-{\sigma^{(Y_1,y_{12},y_3)}} \ar[r]^-{\sigma^{(Y_2,y_1,y_{23})}} & A_{Y_2} \ar[d]^-{\sigma^{(Y,y_2,y_3)}}\\
A_{Y_1} \ar[r]^-{\sigma^{(Y,y_1,y_2)}} & A_Y}\]
is commutative.
\item
In the setting of \eqref{uwd-move-d3}, the diagram
\[\nicexy@C+.5cm{A_X \ar[d]_-{\sigma^{(Y',y_1,y_2)}} \ar[r]^-{\lambda_{(X,\xsubpm)}} & A_{X'} \ar[d]^-{\sigma^{(Y,y_1,y_2)}}\\
A_{Y'} \ar[r]^-{\lambda_{(Y',\xsubpm)}} & A_Y}\]
is commutative.
\item
In the setting of \eqref{uwd-move-d4}, the diagram
\[\nicexy@C+.5cm{A_X \ar[d]_-{\sigma^{(W,y,x_+)}} \ar[r]^-{\sigma^{(W,y,x_-)}} & A_W \ar[d]^-{\lambda_{(W,\xsubpm)}}\\
A_W \ar[r]^-{\lambda_{(W,\xsubpm)}} & A_Y}\]
is commutative.
\item
In the setting of \eqref{uwd-move-e1}, the diagram
\[\nicexy@C+.5cm{A_X \ar[d]_-{\lambda_{(X,x_1,x_2)}} \ar[r]^-{\lambda_{(X,x_3,x_4)}} & A_Z \ar[d]^-{\lambda_{(Z,x_1,x_2)}}\\
A_W \ar[r]^-{\lambda_{(W,x_3,x_4)}} & A_Y}\]
is commutative.
\end{enumerate}
This finishes the definition of a $\uwd$-algebra.
\end{definition}

At this moment we have two definitions of a $\uwd$-algebra.  
\begin{itemize}
\item
On the one hand, in Def. \ref{def2:operad-algebra} with $\O = \uwd$, a $\uwd$-algebra has a structure map $\mu_{\zeta}$ \eqref{operad-structure-map} for each undirected wiring diagram $\zeta \in \uwd$.  This structure map satisfies the associativity axiom \eqref{operad-algebra-associativity2} for a general operadic composition in $\uwd$, together with the unity and the equivariance axioms in Def. \ref{def1:operad-algebra}.  
\item
On the other hand, in Def. \ref{def:uwd-algebra-biased} a $\uwd$-algebra has $6$ generating structure maps and satisfies $17$ generating axioms.  
\end{itemize}
We now observe that these two definitions are equivalent, so $\uwd$-algebras indeed have a finite presentation as in Def. \ref{def:uwd-algebra-biased}.

\begin{theorem}
\label{thm:uwd-algebra}
For the operad  $\uwd$ of undirected wiring diagrams (Theorem \ref{uwd-operad}), Def. \ref{def2:operad-algebra} with $\O=\uwd$ and Def. \ref{def:uwd-algebra-biased} of a $\uwd$-algebra are equivalent.\index{finite presentation for uwd algebra@finite presentation for $\uwd$-algebras}
\end{theorem}

\begin{proof}
This proof proceeds as in the proof of Theorem \ref{thm:wd-algebra}.  First suppose $(A,\mu)$ is a $\uwd$-algebra in the sense of Def. \ref{def2:operad-algebra}.  To see that it is also a $\uwd$-algebra in the sense of Def. \ref{def:uwd-algebra-biased}, first note that the structure map $\mu_?$ \eqref{operad-structure-map} gives the $6$ generating  structure maps \eqref{uwd-structuremap-empty}--\eqref{uwd-structuremap-split}.  Moreover, the generating structure map $\mu_{\tensorunit_X}$ \eqref{uwd-structuremap-namechange} is the identity map by the unity axiom \eqref{operad-algebra-unity}.

The generating axiom \eqref{uwd-algebra-equivariance} is a special case of the equivariance diagram \eqref{operad-algebra-eq}, so it is commutative.  Each of the other $16$ generating axioms corresponds to an elementary relation in $\uwd$ that describes two different ways to construct the same undirected wiring diagram as an iterated operadic composition of generators in $\uwd$.  Each such generating axiom asserts that the two corresponding compositions of generating structure maps--defined using the composition in the diagram \eqref{operad-algebra-associativity2}--are equal.  The associativity axiom \eqref{operad-algebra-associativity2} of $(A,\mu)$ applied twice guarantees that two such compositions are indeed equal.

Conversely, suppose $A$ is a $\uwd$-algebra in the sense of Def. \ref{def:uwd-algebra-biased}, so it has $6$ generating structure maps that satisfy $17$ generating axioms.  We must show that it is a $\uwd$-algebra in the sense of Def. \ref{def2:operad-algebra}  For an undirected wiring diagram $\psi \in \uwd$ with a presentation $\Psi$ (Def. \ref{def:simplex-uwd}), we define its structure map $\mu_{\psi}$ \eqref{operad-structure-map} inductively as follows.
\begin{enumerate}
\item
If $\Psi$ is a $1$-simplex, then $\Psi = (\psi)$, and $\psi$ is a generator in $\uwd$ by definition.  In this case, we define $\mu_{\psi}$ as the corresponding generating structure map  \eqref{uwd-structuremap-empty}--\eqref{uwd-structuremap-split} of $A$.
\item
Inductively, suppose $\Psi$ is an $n$-simplex for some $n \geq 2$, so $\Psi = (\Phi, i, \Theta)$ for some $i \geq 1$, $p$-simplex $\Phi$, and $q$-simplex $\Theta$ with $p+q = n$.  Since $1 \leq p,q < n$, by the induction hypothesis, the structure maps $\mu_{|\Phi|}$ and $\mu_{|\Theta|}$ are already defined.  Then we define the structure map
\begin{equation}
\label{mu-compi-algebra-uwd}
\mu_{\psi} = \mu_{|\Phi|} \compi \mu_{|\Theta|}
\end{equation}
as in Notation \ref{notation:algebra-compi}.
\end{enumerate}

By Theorem \ref{stratified-presentation-exists-uwd} every undirected wiring diagram has a stratified presentation, hence a presentation.  To see that the structure map $\mu_{\psi}$ as above is well-defined, we need to show that the map $\mu_{\psi}$ is independent of the choice of a presentation $\Psi$.  Any two presentations of an undirected wiring diagram are by definition equivalent simplices.  By Theorem \ref{thm:uwd-generator-relation}(2) ($=$ the relations part of the finite presentation theorem for $\uwd$), any two equivalent simplices in $\uwd$ are either equal or are connected by a finite sequence of elementary equivalences.  Therefore, it suffices to show that every elementary equivalence in $\uwd$ yields a commutative diagram involving the generating structure maps of $A$, where $\compi$ is interpreted as in Notation \ref{notation:algebra-compi}.  Recall from Def. \ref{def:equivalent-simplices-uwd} that an elementary equivalence comes from either an elementary relation in $\uwd$ or an operad associativity/unity axiom for the generators in $\uwd$.

It follows from a direct inspection that the operad associativity and unity axioms--\eqref{compi-associativity}, \eqref{compi-associativity-two}, \eqref{compi-left-unity}, and \eqref{compi-right-unity}--for the generators yield commutative diagrams involving the generating structure maps of $A$.  In fact, the diagrams involving the horizontal and the vertical associativity axioms \eqref{compi-associativity} and \eqref{compi-associativity-two} are commutative because composition of functions is associative.  The diagrams for the two unity axioms \eqref{compi-left-unity} and \eqref{compi-right-unity} are commutative because the generating structure map for a colored unit \eqref{uwd-structuremap-namechange} is required to be the identity map.

By definition each of the $17$ generating axioms of $A$ corresponds to an elementary relation (Def. \ref{def:elementary-relation-uwd}) and is a commutative diagram.   Therefore, the structure map $\mu_{\psi}$ for each wiring diagram $\psi$ is well-defined.

It remains to check that the structure map $\mu$ satisfies the required unity, equivariance, and associativity axioms.  The unity axiom \eqref{operad-algebra-unity} holds because it is part of the assumption on the generating structure map corresponding to a name change \eqref{uwd-structuremap-namechange}.

The associativity axiom \eqref{operad-algebra-associativity2} holds because the structure map $\mu_{\psi}$ is defined above \eqref{mu-compi-algebra-uwd} by requiring that the  diagram  \eqref{operad-algebra-associativity2} be commutative.

For the equivariance axiom \eqref{operad-algebra-eq}, first note that it is enough to check it when the undirected wiring diagram in questioned is an iterated operadic composition of $2$-cells.  This is because $2$-cells are the only binary generators in $\uwd$ (Remark \ref{rk:uwd-generators-arity}).  All other generators are either $0$-ary or unary, for which equivariance is trivial.

So now suppose $\zeta$ in the equivariance axiom \eqref{operad-algebra-eq} is an iterated operadic composition of $2$-cells.  If $\zeta$ is a $2$-cell and the permutation $\sigma$ is the transposition $(1~2) \in \Sigma_2$, then  the equivariance axiom \eqref{operad-algebra-eq} is true by the generating axiom \eqref{uwd-algebra-equivariance}.  The general case now follows from this special case using:
\begin{itemize}
\item
the generating axiom \eqref{uwd-algebra-2cell-associativity} corresponding to the associativity property of $2$-cells \eqref{uwd-move-c2};
\item
the operad associativity axioms \eqref{compi-associativity} and \eqref{compi-associativity-two} when applied to $2$-cells;
\item 
the fact that the transpositions $(i, i+1)$ for $1 \leq i \leq n-1$ generate the symmetric group $\Sigma_n$.
\end{itemize}
So  $(A,\mu)$ is a $\uwd$-algebra in the sense of Def. \ref{def2:operad-algebra}.
\end{proof}

\begin{remark}
\label{rk:finite-presentation-algebra}
The proofs of the finite presentation theorems \ref{thm:uwd-algebra} and \ref{thm:wd-algebra} for $\uwd$-algebras and $\WD$-algebras are almost identical.  In fact, it is not difficult to formulate and prove a more general result that has both of these finite presentation theorems as special cases.  Such a result would say that, if an operad $\O$ has a finite presentation (i.e., specific finite sets of generators and generating relations expressed in terms of simplices and elementary equivalences similar to Def. \ref{def:simplex-uwd} and \ref{def:equivalent-simplices-uwd}), then $\O$-algebras have a corresponding finite presentation.  We purposely chose not to present the material this way in order to avoid the higher level of abstraction that is unnecessary for actual applications of (undirected) wiring diagrams.  Although the context is slightly different, the formulation and proof of such a finite presentation theorem for $\O$-algebras can be extracted from the Strong Biased Definition Theorem in \cite{jy2} page 193.
\end{remark}

\section{Finite Presentation for the Relational Algebra}
\label{sec:relational-algebra}

The purpose of this section is to provide an illustration of Theorem \ref{thm:uwd-algebra}, the finite presentation theorem of $\uwd$-algebras, using the relational algebra.  This algebra was originally introduced as the main algebra example in \cite{spivak13} using Def. \ref{def1:operad-algebra}.  We will describe the relational algebra in terms of $6$ generating structure maps and $17$ generating axioms.  First we need some notations.

\begin{definition}
\label{def:relational-construction}
Suppose $A$ and $X$ are sets.
\begin{enumerate}
\item
$A^X = \set(X,A)$ is the set of functions $X \to A$.  
\item
An element $u \in A^X$ is called an \emph{$X$-vector in $A$}, and $u(x) \in A$ is called the \emph{$x$-entry of $u$} for $x \in X$.
\item
$\powerset(X) = \bigl\{T \subseteq X\bigr\}$ is the set of subsets of $X$, called the \emph{power set of $X$}.\index{power set}
\item
$\rela(X) = \powerset\left(A^X\right)$ is the power set of the set $A^X$ of $X$-vectors in $A$.
\end{enumerate}
\end{definition}

\begin{example}
\label{ex:rela-functoriality}
Suppose $A, B, X, Y$ are sets and $*$ is a one-element set.
\begin{enumerate}
\item
Since $A^{\varnothing} = *$, it follows that
\begin{equation}
\label{rel-empty}
\rela(\varnothing) = \powerset(*) = \{\varnothing, *\}.
\end{equation}
\item
Since $A^* = A$, it follows that
\begin{equation}
\label{rel-star}
\rela(*) = \powerset(A), 
\end{equation}
the power set of $A$.
\item
There is a canonical bijection
\begin{equation}
\label{a-to-xplusy-bijection}
A^X \times A^Y \cong A^{X \amalg Y}.
\end{equation}
\item
Suppose $h : X \to Y$ is a function.
\begin{enumerate}[(i)]
\item
There is an induced map $h^* : A^Y \to A^X$ sending each map $Y \to A$ to the composition of $X \to Y \to A$; i.e., $h^*$ is pre-composition with $h$.
\item
Likewise, there is an induced map
\begin{equation}
\label{rel-contravariant}
\nicexy{\rela(Y) \ar[r]^-{h^*} & \rela(X)}
\end{equation}
given by pre-composition with $h$.
\end{enumerate}
\item
Suppose $p : A \to B$ is a function.
\begin{enumerate}[(i)]
\item
There is an induced map $p_* : A^X \to B^X$ sending each map $X \to A$ to the composition $X \to A \to B$; i.e., $p_*$ is post-composition with $p$.  
\item
Likewise, there is an induced map
\begin{equation}
\label{relational-induced-map}
\nicexy{\rela(X) \ar[r]^-{p_*} & \relb(X)} 
\end{equation}
given by post-composition with $p$.
\end{enumerate}
\end{enumerate}
\end{example}

\begin{assumption}
Throughout this section, $S$ is the one-element set.   So $\Fins = \Fin$ is the collection of finite sets, and $\uwd$ is a $\Fin$-colored operad (Theorem \ref{uwd-operad}).  
\end{assumption}

We now define the relational algebra in terms of finitely many generating structure maps and generating axioms.

\begin{definition}
\label{def:relational-algebra-biased}
Suppose $A$ is a set.  The \emph{relational algebra of $A$},\index{relational algebra of a set} denoted $\rela$, is the $\uwd$-algebra in the sense of Def. \ref{def:uwd-algebra-biased} defined as follows.  For each finite set $X$, the $X$-colored entry is
\[\rela(X) = \powerset\left(A^X\right)\]
as in Def. \ref{def:relational-construction}.  Its $6$ generating structure maps are defined as follows.
\begin{enumerate}
\item
For the empty cell $\epsilon \in \uwd\emptynothing$ (Def. \ref{def:uwd-gen-empty}), the chosen element in $\rela(\varnothing) = \{\varnothing, *\}$ \eqref{rel-empty} is $*$.
\item
For a $1$-output wire $\omega_* \in \uwd\starnothing$ (Def. \ref{def:uwd-gen-1output}), the chosen element in $\rela(*) = \powerset(A)$ \eqref{rel-star} is $A$.
\item
For a bijection $f : X \to Y \in \Fin$ and the name change $\tau_f \in \uwd\yx$ (Def. \ref{def:uwd-gen-namechange}), the generating structure map is the pre-composition map \eqref{rel-contravariant}
\[\nicexy@C+1cm{\rela(X) \ar[r]^-{\tau_f ~=~ (\finv)^*} & \rela(Y)}.\]
In other words, each $X$-vector in $A$ is reindexed as a $Y$-vector in $A$ using the bijection $f$.
\item
For a $2$-cell $\theta_{(X,Y)} \in \uwd\xplusyxy$ (Def. \ref{def:uwd-gen-2cell}), the generating structure map is defined using the bijection \eqref{a-to-xplusy-bijection} as follows.
\[\nicexy{\rela(X) \times \rela(Y) \ar[d]^-{\theta_{(X,Y)}} & \bigl(U \subseteq A^X, V \subseteq A^Y\bigr) \ar@{|->}[d] \\
\rela(X \amalg Y) & \bigl(U \times V \subseteq A^X \times A^Y \cong A^{X \amalg Y} \bigr)}\]
In other words, concatenate every given $X$-vector with every given $Y$-vector to form an $(X \amalg Y)$-vector.
\item
For a loop $\lambda_{(X,\xsubpm)} \in \uwd\xminusplusminusx$ (Def. \ref{def:uwd-gen-loop}), the generating structure map is defined as
\begin{equation}
\label{relational-structure-loop}
\nicexy{\rela(X) \ar[d]^-{\lambda_{(X,\xsubpm)}} & \bigl(U \subseteq A^X\bigr) \ar@{|->}[d]\\
\rela\left(X \setminus \xsubpm\right) & \Bigl\{ u_{\setminus \xsubpm} : u \in U,\, u(x_+) = u(x_-) \Bigr\} \subseteq A^{X \setminus \xsubpm}}
\end{equation}
in which $u_{\setminus \xsubpm}$ is the composition
\[\nicexy@C+.5cm{X \setminus \xsubpm \ar[r]^-{\mathrm{inclusion}} & X \ar[r]^-{u} & A}.\] 
In other words, for a given subset of $X$-vectors in $A$, look for those whose $x_+$-entry and $x_-$-entry are equal, and then delete these two entries to form a subset of $(X \setminus \xsubpm)$-vectors.
\item
For a split $\sigma^{(X,x_1,x_2)} \in \uwd\xxprime$ (Def. \ref{def:uwd-gen-split}), denote by $p : X \to X'$ the projection map that sends both $x_1,x_2 \in X$ to $x \in X'$ and is the identity function everywhere else.  The generating structure map is defined as the pre-composition map \eqref{rel-contravariant}
\[\nicexy@C+1.3cm{\rela(X') \ar[r]^-{\sigma^{(X,x_1,x_2)} ~=~ p^*} & \rela(X)}.\]
In other words, for a given subset of $X'$-vectors in $A$, use the $x$-entry for both the $x_1$-entry and the $x_2$-entry to form a subset of $X$-vectors.
\end{enumerate}
\end{definition}

The following observation is the finite presentation theorem for the relational algebra of a set.

\begin{theorem}
\label{relational-algebra-is-algebra}
For each set $A$, the relational algebra $\rela$ in Def. \ref{def:relational-algebra-biased} is actually a $\uwd$-algebra in the sense of Def. \ref{def:uwd-algebra-biased}, hence also in the sense of Def. \ref{def2:operad-algebra} by Theorem \ref{thm:uwd-algebra}.\index{finite presentation for the relational algebra}
\end{theorem}

\begin{proof}
We need to check that $\rela$ satisfies the $17$ generating axioms in Def. \ref{def:uwd-algebra-biased}.  Due to the simplicity of the $6$ generating structure maps, all 
the generating axioms can be checked by a direct inspection.  For example, the generating axiom \eqref{uwd-algebra-move-b1}, which corresponds to the elementary relation \eqref{uwd-move-b1}, is the assertion that the diagram
\[\nicexy{\rela(X) \times \ast \cong \rela(X) \ar[d]_-{(\Id, \omega_y)} \ar[rr]^-{\sigma^{(Y,x,y)}} && \rela(Y) \ar[d]^-{\lambda_{(Y,x,y)}}\\
\rela(X) \times \rela(y) \ar[r]^-{\theta_{(X,y)}} & \rela(Y) \ar[r]^-{\lambda_{(Y,x,y)}} & \rela(W)}\]
is commutative, where $X = Y \setminus y$ and $W = Y \setminus \{x,y\} = X \setminus x$.  For each element
\[\left(U \subseteq A^X\right) \in \rela(X),\]
one can check that both compositions in the above diagram send $U$ to the element
\[\Bigl\{\minushspace\nicexy@C+.3cm{W \ar[r]^-{\mathrm{inclusion}} & X \ar[r]^-{u} & A} : u \in U \Bigr\} \subseteq A^W\]
in $\rela(W)$.  The other $16$ generating axioms are checked similarly.
\end{proof}

\begin{remark}
\label{rk:relational-algebra-agree}
To see that our relational algebra $\rela$ in Def. \ref{def:relational-algebra-biased} agrees with the one in \cite{spivak13} Example 2.2.10, note that the latter is based on Def. \ref{def1:operad-algebra}, which is equivalent to Def. \ref{def2:operad-algebra}.  A direct inspection of \cite{spivak13} Eq. (11) reveals that Spivak's structure map of $\rela$, when applied to the $6$ generators in $\uwd$ (Def. \ref{def:generating-uwd}), reduces to our  generating structure maps in  Def. \ref{def:relational-algebra-biased}.  Theorem \ref{relational-algebra-is-algebra} then guarantees that the two definitions of the relational algebra $\rela$--i.e., our Def. \ref{def:relational-algebra-biased}  and \cite{spivak13} Example 2.2.10--are equivalent.
\end{remark}

\section{Spivak's Conjecture: Rigidity of the Relational Algebra}
\label{sec:rigidity-relational}

The purpose of this section is to partially verify a conjecture in \cite{spivak13} (Conjecture 3.1.6) that states that the relational algebra $\rela$ in Def. \ref{def:relational-algebra-biased} is quotient-free.  To state this conjecture, we first need the definition of a map between operad algebras.  Recall from Section \ref{sec:operad-algebras} the definition of an operad algebra.

\begin{definition}
\label{def:operad-algebra-map}
Suppose $\O$ is an $S$-colored operad, and $(A,\mu^A)$ and $(B,\mu^B)$ are $\O$-algebras.  A \emph{map of $\O$-algebras}\index{map of operad algebras} $f : A \to B$ consists of a collection of maps
\[\Bigl\{\hspace{-.1cm}\nicexy{A_c \ar[r]^-{f_c} & B_c} :  c \in S\Bigr\}\]
that is compatible with the structure maps in the following sense.
\begin{description}
\item[Compatibility with Structure Maps]
For each $d \in S$, $\uc = (c_1,\ldots,c_n) \in \profs$, and $\zeta \in \O\duc$, the diagram
\begin{equation}
\label{algebra-map-compatibility}
\nicexy{\prod\limits_{i=1}^n A_{c_i} \ar[d]_-{\mu^A_{\zeta}} \ar[r]^-{\prod f_{c_i}} & \prod\limits_{i=1}^n B_{c_i} \ar[d]^-{\mu^B_{\zeta}}\\
A_d \ar[r]^-{f_d} & B_d}
\end{equation}
is commutative.
\end{description}
Furthermore, we call $f$ an \emph{isomorphism}\index{isomorphism of operad algebras} if there exists a map of $\O$-algebras $g : B \to A$ such that $gf = \Id_A$ and $fg = \Id_B$.  If such a map $g$ exists, then it is necessarily unique.
\end{definition}

The following conjecture regarding the relational algebra is \cite{spivak13} Conjecture 3.1.6.  Since we are talking about the relational algebra, here $S$ is a one-element set, and $\uwd$ is a $\Fin$-colored operad.

\begin{conjecture}
\label{conj:relational-quotient-free}
Suppose:\index{Spivak's Conjecture}
\begin{itemize}
\item
$A$ is a set, and $\rela$ is the relational algebra of $A$ in Def. \ref{def:relational-algebra-biased}.
\item
$f : \rela \to B$ is a map of $\uwd$-algebras.
\end{itemize}
Then at least one of the following two statements holds.
\begin{enumerate}
\item
$f$ is an isomorphism.
\item
$B_X$ is a one-element set for each finite set $X$.
\end{enumerate}
\end{conjecture}

Roughly speaking this conjecture states that there are no interesting maps out of any relational algebra.  In the following observation, we will verify Conjecture \ref{conj:relational-quotient-free} in the special case when the map $\rela \to B$ is induced by a map of sets out of $A$.

\begin{theorem}
\label{thm:relational-rigidity}
Suppose $f : A \to B$ is a map of sets.  Then the following statements are equivalent.
\begin{enumerate}
\item
$f$ is a bijection of sets.
\item
The post-composition maps $f_* : \rela(X) \to \relb(X)$ \eqref{relational-induced-map}, with $X$ running through all the finite sets, form an isomorphism of $\uwd$-algebras.
\item
The post-composition maps $f_* : \rela(X) \to \relb(X)$, with $X$ running through all the finite sets, form a map of $\uwd$-algebras.
\end{enumerate}
\end{theorem}

\begin{proof}
The implications $(1) \Longrightarrow (2) \Longrightarrow (3)$ are both immediate from the definitions.  It remains to check the implication $(3) \Longrightarrow (1)$.  So let us now assume that $f_* : \rela \to \relb$ is a map of $\uwd$-algebras.  This means that the diagram
\begin{equation}
\label{relational-algebra-map-compatibility}
\nicexy{\prod\limits_{i=1}^n \rela(X_i) \ar[d]_-{\mu^A_{\zeta}} \ar[r]^-{\prod f_{X_i}} & \prod\limits_{i=1}^n \relb(X_i) \ar[d]^-{\mu^B_{\zeta}}\\
\rela(Y) \ar[r]^-{f_Y} & \relb(Y)}
\end{equation}
is commutative for each $\zeta \in \uwd\yux$ with $\uX = (X_1,\ldots,X_n)$.  We will consider two special cases, which will show that $f$ is surjective and injective.
\begin{enumerate}
\item
To show that $f$ is surjective, consider a $1$-output wire $\omega_* \in \uwd\starnothing$ (Def. \ref{def:uwd-gen-1output}).  In this case, the commutative diagram \eqref{relational-algebra-map-compatibility} becomes the diagram
\[\nicexy{\ast \ar[d]_-{\omega_*^A} \ar[r]^-{=} & \ast \ar[d]^-{\omega_*^B}\\ \rela(*) = \powerset(A) \ar[r]^-{f_*} & \powerset(B) = \relb(*)}\]
in which
\[\omega_*^A(*) = A \in \powerset(A) \andspace \omega_*^B(*) = B \in \powerset(B)\]
by Def. \ref{def:relational-algebra-biased}.  The bottom horizontal map $f_*$ is post-composition with $f$, so it sends each subset $U \subseteq A$ to its image $f(U) \subseteq B$.  Therefore, the commutativity of the above diagram forces
\[f(A) = f_*(A) = B,\]
so $f$ is surjective.
\item
To show that $f$ is injective, we argue by contradiction.  So suppose there exist distinct elements $a_+, a_- \in A$ such that $f(a_+) = f(a_-) \in B$.  We will show that this assumption leads to a contradiction.   Consider the loop $\lambda_{(A,\apm)} \in \uwd\aminusplusminusa$ (Def. \ref{def:uwd-gen-loop}).
In this case, the commutative diagram \eqref{relational-algebra-map-compatibility} becomes the diagram
\begin{equation}
\label{relational-rigid-loop}
\nicexy{\rela(A) \ar[d]_-{\lambda_{(A,\apm)}} \ar[r]^-{f_*} & \relb(A) \ar[d]^-{\lambda_{(A,\apm)}}\\
\rela(A\setminus \apm) \ar[r]^-{f_*} & \relb(A \setminus \apm)}
\end{equation}
in which both horizontal maps $f_*$ are post-composition with $f$.  Consider the element
\[U = \left\{\nicexy{A \ar[r]^-{=} & A}\right\} \in \rela(A) = \powerset\left(A^A\right),\]
which is the single-element set consisting of the identity map of $A$.  We will show that the two compositions in \eqref{relational-rigid-loop} do not agree at $U$.

On the one hand, since $a_+ \not= a_-$, applying the left vertical map to $U$ yields
\[\lambda_{(A,\apm)}(U) = \varnothing \subseteq A^{A \setminus \apm}\]
by the definition of the generating structure map for a loop \eqref{relational-structure-loop}.  Applying the bottom horizontal map $f_*$ in \eqref{relational-rigid-loop} yields
\[f_*\lambda_{(A,\apm)}(U) = \varnothing \subseteq B^{A \setminus \apm}.\]
On the other hand, applying the top horizontal map $f_*$ to $U$ yields
\[f_*(U) = \Bigl\{\nicexy{A \ar[r]^-{f} & B}\Bigr\} \in \relb(A) = \powerset\left(B^A\right),\]
which is the single-element set consisting of the map $f$.  Since $f(a_+) = f(a_-)$, applying the right vertical map now yields the single-element set
\[\lambda_{(A,\apm)} f_*(U) = \Bigl\{\nicexy@C+.4cm{A \setminus \apm \ar[r]^-{\mathrm{inclusion}} & A \ar[r]^-{f} & B}\Bigr\} \subseteq B^{A \setminus \apm}\]
by \eqref{relational-structure-loop}.  Since this is not empty, the two compositions in the commutative diagram \eqref{relational-rigid-loop} do not agree at the element $U$.  This is a contradiction, so the map $f$ is injective.
\end{enumerate}
\end{proof}

\begin{example}
Let us provide a simple illustration of the non-commutativity of the diagram \eqref{relational-rigid-loop} for the non-injective function $f : A = \{0,1\} \to \{*\} = B$ that sends both $0$ and $1$ to $*$.  Consider the singleton $U = \{\Id_A\} \in \powerset(A^A) = \rela(A)$.  Since $0 \not= 1$, we have
\[f_*\lambda_{(A,0,1)}U = f_*(\varnothing) = \varnothing \in \powerset(B^{\varnothing}) = \relb(\varnothing).\]
On the other hand, $f_*U$ is the singleton $\{f\} \in \powerset(B^A) = \relb(A)$ with $f(0) = f(1)$.  So
\[\lambda_{(A,0,1)}f_*U = \lambda_{(A,0,1)}\{f\} = \{\varnothing \to B\} \in \powerset(B^\varnothing) = \relb(\varnothing),\]
which contains one element.  So $\lambda_{(A,0,1)}f_* \not= f_*\lambda_{(A,0,1)}$, and the diagram \eqref{relational-rigid-loop} in this case is not commutative.  In particular, this non-injective function $f$ does \emph{not} induce a map of $\uwd$-algebras $\rela \to \relb$.
\end{example}

\section{Finite Presentation for the Typed Relational Algebra}
\label{sec:relational-typed}

The relational algebra $\rela$ in Def. \ref{def:relational-algebra-biased} has a fixed set $A$ as the set of potential values in each coordinate in an $X$-vector (Def. \ref{def:relational-construction}).  There is a more general version of the relational algebra, called the typed relational algebra, in which each coordinate in an $X$-vector has its own set of potential values.  The typed relational algebra was first introduced in \cite{spivak13} (Section 4) using Def. \ref{def1:operad-algebra}.  In this section, we observe that the typed  relational algebra also has a finite presentation, similar to the one for the relational algebra in Theorem \ref{relational-algebra-is-algebra}.

\begin{assumption}
Throughout this section, $S$ is the collection of sets, so $\uwd$ is the $\Fin_{\set}$-colored operad of undirected wiring diagrams (Theorem \ref{uwd-operad}).
\end{assumption}

\begin{definition}
Suppose $(X,v) \in \Finset$ (Def. \ref{def:Fins}), so $X$ is a finite set and $v : X \to \set$ assigns to each $x \in X$ a set $v(x)$.  
\begin{enumerate}
\item
Define the set
\begin{equation}
\label{x-sub-v-finset}
X_v = \prod_{x\in X} v(x)
\end{equation}
in which an empty product, for the case $X = \varnothing$, means the one-point set.
\item
An element in $X_v$ is also called an \emph{$X$-vector}.
\item
For $x \in X$, the $x$-entry of an $X$-vector $u$ is denoted by $u_x$.
\end{enumerate}
As before, we will omit writing $v$ in $(X,v)$ if there is no danger of confusion.
\end{definition}

\begin{example}
Suppose $A$ is a set and $(X,v) \in \Finset$ such that $v(x) = A$ for all $x \in X$.  Then 
\[X_v = \prod_{x \in X} A = A^X\]
as in Def. \ref{def:relational-construction}
\end{example}

\begin{example}
Each map $f : (X,v) \to (Y,v) \in \Finset$ induces a map
\begin{equation}
\label{f-sub-v-map}
\nicexy{Y_v = \prod\limits_{y \in Y} v(y) \ar[r]^-{f_v} & \prod\limits_{x \in X} v(x) = X_v}.
\end{equation}
For $(a_y)_{y \in Y} \in Y_v$ with each $a_y \in v(y)$, its image is defined as
\[f_v\Bigl((a_y)_{y\in Y}\Bigr) = \left(a_{f(x)}\right)_{x \in X}\]
using the fact that $a_{f(x)} \in v(f(x)) = v(x)$.  If $A$ is a set and $v(y) = v(x) = A$ for all $y \in Y$ and $x \in X$, then $f_v$ is the pre-composition map in Example \ref{ex:rela-functoriality}.  In what follows, a map induced by the map $f_v$ will often be denoted by $f_v$ as well.
\end{example}

We now define the typed relational algebra in terms of finitely many generating structure maps and generating axioms.

\begin{definition}
\label{def:typed-relational-algebra}
The \emph{typed relational algebra}\index{typed relational algebra} $\relation$ is the $\uwd$-algebra in the sense of Def. \ref{def:uwd-algebra-biased} defined as follows.  For each $(X,v) \in \Finset$, the $X$-colored entry is
\[\relation(X) = \powerset\left(X_v\right)\]
with $X_v$ as in \eqref{x-sub-v-finset} and $\powerset(-)$ the power set (Def. \ref{def:relational-construction}).  Its $6$ generating structure maps are defined as follows.
\begin{enumerate}
\item
For the empty cell $\epsilon \in \uwd\emptynothing$ (Def. \ref{def:uwd-gen-empty}), the chosen element in 
\[\relation(\varnothing) = \powerset(\varnothing_v) = \powerset(*) = \{\varnothing, *\}\]
is $*$.
\item
For a $1$-output wire $\omega_x \in \uwd\smallxnothing$ (Def. \ref{def:uwd-gen-1output}), the chosen element in $\relation(x) = \powerset(v(x))$ is $v(x)$.
\item
For a bijection $f : X \to Y \in \Finset$ and the name change $\tau_f \in \uwd\yx$ (Def. \ref{def:uwd-gen-namechange}), the generating structure map is the bijection
\[\nicexy@C+.5cm{\relation(X) = \powerset\left(X_v\right) \ar[r]^-{\tau_f ~=~ \finv_v} & \powerset\left(Y_v\right) = \relation(Y)}\]
induced by $\finv_v$ as in \eqref{f-sub-v-map}.
\item
For a $2$-cell $\theta_{(X,Y)} \in \uwd\xplusyxy$ (Def. \ref{def:uwd-gen-2cell}), the generating structure map is defined as follows.
\[\nicexy{\relation(X) \times \relation(Y) \ar[d]^-{\theta_{(X,Y)}} & \bigl(U \subseteq X_v, V \subseteq Y_v\bigr) \ar@{|->}[d] \\
\relation(X \amalg Y) & \Bigl(U \times V \subseteq X_v \times Y_v \cong (X \amalg Y)_v \Bigr)}\]
\item
For a loop $\lambda_{(X,\xsubpm)} \in \uwd\xminusplusminusx$ (Def. \ref{def:uwd-gen-loop}), the generating structure map is defined as
\begin{equation}
\label{typed-relational-structure-loop}
\nicexy{\relation(X) \ar[d]^-{\lambda_{(X,\xsubpm)}} & \bigl(U \subseteq X_v\bigr) \ar@{|->}[d]\\
\relation\left(X \setminus \xsubpm\right) & \Bigl\{ u_{\setminus \xsubpm} : u \in U,\, u_{x_+} = u_{x_-} \Bigr\} \subseteq (X \setminus \xsubpm)_v}
\end{equation}
in which $u_{\setminus \xsubpm}$ is obtained from $u$ by deleting the $\xsubpm$-entries.  In other words, if $\iota : (X \setminus \xsubpm) \to X$ is the inclusion map, then
\[u_{\setminus \xsubpm} = \iota_v(u),\]
provided $u_{x_+} = u_{x_-}$, where $\iota_v$ is as in \eqref{f-sub-v-map}.
\item
For a split $\sigma^{(X,x_1,x_2)} \in \uwd\xxprime$ (Def. \ref{def:uwd-gen-split}), denote by $p : X \to X'$ the projection map that sends both $x_1,x_2 \in X$ to $x \in X'$ and is the identity function everywhere else.  The generating structure map is the map
\[\nicexy@C+1.5cm{\relation(X') \ar[r]^-{\sigma^{(X,x_1,x_2)} ~=~ p_v} & \relation(X)}\]
induced by $p_v$ in \eqref{f-sub-v-map}.  In other words, $\sigma^{(X,x_1,x_2)}$ 
takes each $X'$-vector $u$ in any given subset of $X'$-vectors to the $X$-vector whose $x_1$-entry and $x_2$-entry are both equal to the $x$-entry of $u$, and all other entries remain the same.
\end{enumerate}
\end{definition}

\begin{example}
\label{ex:rela-reduces-to-rela}
Suppose $A$ is a set and $(X,v) \in \Finset$ such that $v(x) = A$ for all $x \in X$.  Then
\[\relation(X) = \powerset(X_v) = \powerset\left(A^X\right) = \rela(X),\]
the $X$-colored entry of the relational algebra of $A$ (Def. \ref{def:relational-algebra-biased}).  Furthermore, if in Def. \ref{def:typed-relational-algebra} all the value assignments $v$ take the constant value $A$, then the $6$ generating structure maps of $\relation$ reduce to those of $\rela$ in Def. \ref{def:relational-algebra-biased}.
\end{example}

The following observation is the finite presentation theorem for the typed relational algebra.

\begin{theorem}
\label{typed-relational-algebra}
The typed relational algebra $\relation$ in Def. \ref{def:typed-relational-algebra} is actually a $\uwd$-algebra in the sense of Def. \ref{def:uwd-algebra-biased}, hence also in the sense of Def. \ref{def2:operad-algebra} by Theorem \ref{thm:uwd-algebra}.\index{finite presentation for the typed relational algebra}
\end{theorem}

\begin{proof}
As in the proof of Theorem \ref{relational-algebra-is-algebra} for the relational algebra of $A$, we just need to check that $\relation$ satisfies the $17$ generating axioms in Def. \ref{def:uwd-algebra-biased}.  Due to the simplicity of the $6$ generating structure maps, all the generating axioms can be checked by a direct inspection.   For example, the generating axiom \eqref{uwd-algebra-move-b1}, which corresponds to the elementary relation \eqref{uwd-move-b1}, is the assertion that the diagram
\[\nicexy{\relation(X) \times \ast \cong \relation(X) \ar[d]_-{(\Id, \omega_y)} \ar[rr]^-{\sigma^{(Y,x,y)}} && \relation(Y) \ar[d]^-{\lambda_{(Y,x,y)}}\\
\relation(X) \times \relation(y) \ar[r]^-{\theta_{(X,y)}} & \relation(Y) \ar[r]^-{\lambda_{(Y,x,y)}} & \relation(W)}\]
is commutative, where $X = Y \setminus y$ and $W = Y \setminus \{x,y\} = X \setminus x$.  For each element
\[\Bigl(U \subseteq X_v = \prod v(x)\Bigr) \in \relation(X),\]
one can check that both compositions in the above diagram send $U$ to the element
\[\Bigl\{ \iota_v(u) : u \in U \Bigr\} \subseteq W_v\]
in $\relation(W)$, where $\iota : W \to X$ is the inclusion map and $\iota_v$ is as in \eqref{f-sub-v-map}.  The other $16$ generating axioms are checked similarly.
\end{proof}

\begin{remark}
\label{rk:typed-relational-algebra-agree}
To see that our relational algebra $\relation$ in Def. \ref{def:typed-relational-algebra} agrees with the one in \cite{spivak13} Section 4, note that the latter is based on Def. \ref{def1:operad-algebra}, which is equivalent to Def. \ref{def2:operad-algebra}.  A direct inspection of \cite{spivak13} Lemma 4.1.2 reveals that Spivak's structure map of $\relation$, when applied to the $6$ generators in $\uwd$ (Def. \ref{def:generating-uwd}), reduces to our  generating structure maps in  Def. \ref{def:typed-relational-algebra}.  Theorem \ref{typed-relational-algebra} then guarantees that the two definitions of the relational algebra $\relation$--i.e., our Def. \ref{def:typed-relational-algebra}  and \cite{spivak13} Section 4--are equivalent.
\end{remark}

\section{Summary of Chapter \ref{ch11-uwd-algebras}}

\begin{enumerate}
\item Each $\uwd$-algebra can be described using six generating structure maps and seventeen generating axioms.
\item The (typed) relational algebra is a $\uwd$-algebra.
\item Spivak's Conjecture holds when restricted to relational algebras.
\end{enumerate}

\part{Maps Between Operads of Wiring Diagrams}
\label{part:maps-operad}

So far we have considered four operads constructed from wiring diagrams and undirected wiring diagrams:
\begin{enumerate}
\item the $\boxs$-colored operad of wiring diagrams $\WD$ (Theorem \ref{wd-operad});
\item the $\boxs$-colored operad of normal wiring diagrams $\wddot$ (Prop. \ref{prop:without-dn-operads});
\item the $\boxs$-colored operad of strict wiring diagrams $\wdzero$ (Prop. \ref{prop:strict-operads});
\item the $\Fins$-colored operad of undirected wiring diagrams $\uwd$ (Theorem \ref{uwd-operad}). 
\end{enumerate}
The purpose of this part is to study maps between these operads.  We will show that there is a commutative diagram
\[\nicexy{\wdzero \ar[dr]_-{\chizero} \ar[r] & \wddot \ar[d]^-{\chi} \ar[r] & \WD \ar[dl]^-{\rho}\\
& \uwd &}\]
of operad maps, in which the horizontal maps are operad inclusions.  All the operad maps in this diagram except $\rho$ are discussed in Chapter \ref{ch12-maps}.  The operad map $\rho$ is discussed in Chapter \ref{ch13-wd-uwd}.  The existence of such operad maps implies that:
\begin{itemize}
\item Every $\uwd$-algebra induces a $\WD$-algebra along $\rho$.
\item Every $\WD$-algebra restricts to a $\wddot$-algebra.
\item Every $\wddot$-algebra restricts to a $\wdzero$-algebra.
\end{itemize}

For each of the three operad maps that end at $\uwd$, we will compute precisely the image.  An undirected wiring diagram is in the image of
\begin{itemize}
\item $\chizero$ if and only if its cables are either $(1,1)$-cables or $(2,0)$-cables;
\item $\chi$ if and only if it has no $(0,0)$-cables and no $(0, \geq 2)$-cables.
\end{itemize}
Furthermore, the operad map $\rho : \WD \to \uwd$ is surjective, so every undirected wiring diagram is the $\rho$-image of some wiring diagram.  Delay nodes play a crucial role in the surjectivity of the operad map $\rho$.

\textbf{Reading Guide}.  The reader who already knows about operad maps may skip most of Section \ref{sec:operad-maps} and go straight to Propositions \ref{prop:strict-normal-wd} and \ref{prop:uwd-vary-s}.  Instead of the proof of Theorem \ref{wddot-uwd-operad-map} about the existence of the operad map $\chi$, the reader may wish to concentrate on the motivating Example \ref{ex:wddot-to-uwd-motivation}.  Likewise, the proofs of Theorems \ref{thm:chi-image} and \ref{chizero-wdzero-uwd} about the images of the operad maps $\chi$ and $\chizero$ may be skipped during the first reading.

Instead of the proof of Theorem \ref{wd-uwd-operad-map} on the existence of the operad map $\rho$, the reader may wish to concentrate on the motivating Example \ref{ex:rho-motivation}.  Likewise, before reading the proof of Theorem \ref{thm:wd-uwd-surjective} on the surjectivity of the map $\rho$, the reader may wish to first read the illustrating Example \ref{ex:rho-preimage}.

\chapter{Map from Normal to Undirected Wiring Diagrams}
\label{ch12-maps}

This chapter has three main purposes.
\begin{enumerate}
\item We show that there exists an operad map $\chi : \wddot \to \uwd$  from the operad $\wddot$ of normal wiring diagrams to the operad $\uwd$ of undirected wiring diagrams.  See Theorem \ref{wddot-uwd-operad-map}.   Recall that a normal wiring diagram is a wiring diagram with no delay nodes.  Intuitively, the map $\chi$ is given by forgetting directions.  
\item We compute precisely the image of the operad map $\chi : \wddot \to \uwd$ in Theorem \ref{thm:chi-image}.   An undirected wiring diagram is in the image of $\chi$ if and only if it has no wasted cables and no $(0, \geq 2)$-cables.
\item We consider the restriction $\chizero : \wdzero \to \uwd$ of the operad map $\chi : \wddot \to \uwd$ to the operad $\wdzero$ of strict wiring diagrams and compute precisely its image in Theorem \ref{chizero-wdzero-uwd}.   An undirected wiring diagram is in the image of $\chizero$ if and only if its cables are either $(1,1)$-cables or $(2,0)$-cables. 
\end{enumerate}
In Theorem \ref{wd-uwd-operad-map} we will extend the operad map $\chi : \wddot \to \uwd$ to an operad map $\rho : \WD \to \uwd$ defined for all wiring diagrams.  Furthermore, we will show in Theorem \ref{thm:wd-uwd-surjective} that the operad map $\rho : \WD \to \uwd$ is surjective.

In Section \ref{sec:operad-maps} we define a map of operads and observe that there are inclusions of operads $\wdzero \to \wddot \to \WD$.  Furthermore, if the underlying class $S$ changes, then there are corresponding maps of operads for each of the four operads of (undirected) wiring diagrams.

Section \ref{sec:normal-undirected} contains the first main result Theorem \ref{wddot-uwd-operad-map} of this chapter.  It says that there is a map of operads $\chi : \wddot \to \uwd$ defined by forgetting directions.

In Section \ref{sec:examples-normal-undirected} we provide a series of examples to illustrate the operad map $\chi$.

Section \ref{sec:image-normal-undirected} contains the second main result Theorem \ref{thm:chi-image} of this chapter.  This result identifies precisely the image of the operad map $\chi : \wddot \to \uwd$ as consisting of the undirected wiring diagrams with no wasted cables and no $(0, \geq 2)$-cables.

In Section \ref{sec:strict-to-undirected} we consider the restriction of the operad map $\chi$ to the operad $\wdzero$ of strict wiring diagrams.  Recall that a strict wiring diagram is a wiring diagram with no delay nodes and whose supplier assignment is a bijection.  In Theorem \ref{chizero-wdzero-uwd} we will show that the image of the operad map $\wdzero \to \uwd$ consists of precisely the undirected wiring diagrams with only $(1,1)$-cables and $(2,0)$-cables.

\section{Operad Maps}
\label{sec:operad-maps}

In this section, we define an operad map and record some obvious maps among the various operads of (undirected) wiring diagrams (Prop. \ref{prop:strict-normal-wd} and \ref{prop:uwd-vary-s}).  Recall from Def. \ref{def:pseudo-operad} the definition of an $S$-colored operad.

\begin{definition}
\label{def:operad-map}
Suppose $\O$ is an $S$-colored operad and $\P$ is a $T$-colored operad.  A \emph{map of operads},\index{map of operads} also called an \emph{operad map},\index{operad map} $f : \O \to \P$ consists of a pair of maps $(f_0,f_1)$ as follows.
\begin{enumerate}
\item
$f_0 : S \to T$, called the \emph{color map}.\index{color map}
\item
For each $d \in S$ and $\uc = (c_1,\ldots,c_n) \in \profs$ with $n \geq 0$, it has a map, called an \emph{entry map},\index{entry map}
\[\nicexy{\O\duc \ar[r]^-{f_1} & \P\fdfuc}\]
in which $fd = f_0 d \in T$ and $f\uc = (f_0 c_1,\ldots,f_0 c_n) \in \proft$.
\end{enumerate}
We will usually write both $f_0$ and $f_1$ as $f$.  These maps are required to preserve the operad structure in the sense that the following three conditions hold.
\begin{description}
\item[Preservation of Equivariance]
For each $\duc \in \profs \times S$ as above and permutation $\sigma \in \Sigma_{n}$, the diagram
\begin{equation}
\label{operad-map-eq}
\nicexy{\O\duc \ar[d]_-{\sigma}^-{\cong} \ar[r]^-{f} & \P\fdfuc \ar[d]^-{\sigma}_-{\cong}\\
\O\dcsigma \ar[r]^-{f} & \P\fdfcsigma}
\end{equation}
is commutative, in which $f\uc\sigma = \left(fc_{\sigma(1)}, \ldots, fc_{\sigma(n)}\right)$.
\item[Preservation of Colored Units]
For each $c \in S$, 
\begin{equation}
\label{operad-map-unit}
f\tensorunit_c = \tensorunit_{fc}
\end{equation}
in which $\tensorunit_c \in \O\ccsingle$ is the $c$-colored unit in $\O$ and $\tensorunit_{fc} \in \P\fcsingle$ is the $(fc)$-colored unit in $\P$.
\item[Preservation of Operadic Composition]
For each $\duc \in \profs \times S$ with $|\uc| \geq 1$, $\ub \in \profs$, and $1 \leq i \leq |\uc|$, the diagram
\begin{equation}
\label{operad-map-comp}
\nicexy{\O\duc \times \O\ciub \ar[d]_-{\compi} \ar[r]^-{(f,f)} & \P\fdfuc \times \P\fcifub \ar[d]^-{\compi}\\
\O\dccompib \ar[r]^-{f} & \P\fdfccompib}
\end{equation}
is commutative.
\end{description}
\end{definition}

\begin{definition}
\label{def:operad-inclusion}
Suppose $\O$ is an $S$-colored operad and $\P$ is a $T$-colored operad.  
\begin{enumerate}
\item
A map of operads $f : \O \to \P$ is called an \emph{operad inclusion}\index{operad inclusion} if the color map $f_0 : S \to T$ and all the entry maps $f_1$ are inclusions. 
\item
If there exists an operad inclusion $\O \to \P$, then we call $\O$ a \emph{sub-operad of $\P$}.\index{sub-operad}
\end{enumerate}
\end{definition}

\begin{example}
\label{ex:terminal-operad}
With the one-point set $*$ as the set of color, there is an operad $\T$ in which every entry $\T\smallbinom{*}{*,\ldots,*} = *$.  Then for each colored operad $\O$, there exists a unique map of operads $\O \to \T$.  In category theoretical terms, $\T$ is a \emph{terminal object} in the category of all colored operads.\index{terminal operad}
\end{example}

\begin{example}
\label{ex:initial-operad}
There is an operad $\I$ with the empty set as its set of color, and $\I$ has no entries because it has no colors.  Then for each colored operad $\O$, there exists a unique operad inclusion $\I \to \O$.  In category theoretical terms, $\I$ is an \emph{initial object} in the category of all colored operads.\index{initial operad}
\end{example}

\begin{example}
\label{ex:compose-operad-maps}
Suppose $f : \O \to \P$ and $g : \P \to \Q$ are maps of operads.  Composing the color maps and the entry maps, there is a composition map of operads $gf : \O \to \Q$.\index{operad map composition}
\end{example}

\begin{example}
\label{ex:operad-map-algebra}
Suppose $f : \O \to \P$ is an operad map as in Def. \ref{def:operad-map} and $A = \{A_t\}_{t \in T}$ is a $\P$-algebra as in Def. \ref{def1:operad-algebra} or Def. \ref{def2:operad-algebra}.  Then there is an induced $\O$-algebra $A^f$ defined by the following data.\index{induced operad algebra} \index{operad algebra induced along a map}
\begin{enumerate}
\item For each $c \in S$, $A^f$ is equipped with the $c$-colored entry $A^f_c = A_{fc}$.
\item For each $\dconecn \in \profs \times S$ and $\zeta \in \O\dconecn$, $A^f$ is equipped with the structure map
\[\nicexy@C+1cm{A^f_{c_1} \times \cdots \times A^f_{c_n} = A_{fc_1} \times \cdots \times A_{fc_n} \ar[r]^-{\mu^f_{\zeta}\,=\,\mu_{f\zeta}} & A_{fd} = A^f_d}.\]
Here $\mu_?$ is the $\P$-algebra structure map of $A$, and $f\zeta \in \P\fdfconefcn$.  Since an operad map is assumed to preserve all the operad structure, a direct inspection reveals that $A^f$ is indeed an $\O$-algebra.
\end{enumerate}
We say that the $\O$-algebra $A^f$ is \emph{induced along $f$}.
\end{example}

In the next two observations, we record some obvious operad maps among the various operads of (undirected) wiring diagrams.

Recall the $\boxs$-colored operad $\WD$ of wiring diagrams (Theorem \ref{wd-operad}), the $\boxs$-colored operad $\wddot$ of normal wiring diagrams (Prop. \ref{prop:without-dn-operads}), and the $\boxs$-colored operad $\wdzero$ of strict wiring diagrams (Prop. \ref{prop:strict-operads}).  Remember that a normal wiring diagram is a wiring diagram without delay nodes, and a strict wiring diagram is a normal wiring diagram whose supplier assignment is a bijection.  Also recall that $\WD^S$ means the operad of $S$-wiring diagrams, and the symbol $S$ is suppressed from the notation $\WD^S$ unless we need to emphasize it.  Similar remarks apply to $\wdzero^S$ and $\wddot^S$.

\begin{proposition}
\label{prop:strict-normal-wd}
Given a map $f : S \to T$ of classes, there exists an induced commutative diagram
\begin{equation}
\label{wd-maps-vary-s}
\nicexy{\wdzero^S \ar[r] \ar[d]_-{f_*} & \wddot^S \ar[r] \ar[d]_-{f_*} & \WD^S \ar[d]_-{f_*}\\
\wdzero^T \ar[r] & \wddot^T \ar[r] & \WD^T}
\end{equation}
of maps of operads in which:
\begin{itemize}
\item
the horizontal maps are operad inclusions;
\item
the vertical maps are induced by $f$ on value assignments.
\end{itemize}
\end{proposition}

\begin{proof}
In each row of the diagram \eqref{wd-maps-vary-s}, the maps on colors (i.e., either $\boxs$ or $\boxt$) are the identity map.  For a fixed class, a strict wiring diagram is by definition also a normal wiring diagram, which by definition is also a wiring diagram.  In each of $\wdzero$ and $\wddot$, the operad structure--i.e., the equivariant structure, the colored units, and the operad composition--is defined as that in $\WD$ (Def. \ref{wd-equivariance}--\ref{def:compi-wd}).  So the horizontal entrywise inclusions actually define operad inclusions.  

For the vertical maps in the diagram \eqref{wd-maps-vary-s}, first note that $f$ induces maps $\Fins \to \Fint$ and $\boxs \to \boxt$ that are the identity map on the underlying finite sets.  On value assignments, these maps are post-composition with $f$.  A direct inspection reveals that, using  these two maps, every (strict/normal) $S$-wiring diagram is sent to a (strict/normal) $T$-wiring diagram.  Moreover, all the operad structure (Def. \ref{wd-equivariance}--\ref{def:compi-wd}) is preserved by these maps.  So the vertical entrywise defined maps in the diagram \eqref{wd-maps-vary-s} are maps of  operads.  The commutativity of the diagram is immediate from the definitions of the operad maps.
\end{proof}

\begin{example}
By Example \ref{ex:operad-map-algebra} and Prop. \ref{prop:strict-normal-wd}, every $\WD$-algebra (Def. \ref{def:wd-algebra}) restricts to a $\wddot$-algebra (Def. \ref{def:normal-algebra}), and every $\wddot$-algebra restricts to a $\wdzero$-algebra (Def. \ref{def:strict-algebra}).  For example, the propagator algebra (Def. \ref{def:propagator-algebra}), which is a $\WD$-algebra, restricts to a $\wddot$-algebra and also to a $\wdzero$-algebra.
\end{example}

Recall that $\uwd^S$ is the $\Fins$-colored operad of undirected $S$-wiring diagrams (Theorem \ref{uwd-operad}), and the symbol $S$ in $\uwds$ is dropped unless we need to emphasize $S$.  Essentially the same as in Prop. \ref{prop:strict-normal-wd}, we have the following operad maps for undirected wiring diagrams.

\begin{proposition}
\label{prop:uwd-vary-s}
Given a map $f : S \to T$ of classes, there exists an induced map of operads
\[\nicexy{\uwd^S \ar[r] & \uwd^T}.\]
\end{proposition}

\begin{example}
Suppose $S = *$, a one-point set, and $T = \set$, the collection of sets.  Then the inclusion $* \to \set$ induces an operad inclusion $\uwd^* \to \uwd^{\set}$.  In \cite{spivak13} Example 2.1.7, $\uwd^*$ is denoted by $\mathcal{S}$ and is called the operad of singly-typed wiring diagrams.    In \cite{spivak13} Example 4.1.1, $\uwdset$ is denoted by $\mathcal{T}$ and is called the operad of typed wiring diagrams.  
\end{example}

\section{Normal to Undirected Wiring Diagrams}
\label{sec:normal-undirected}

Fix a class $S$ for the rest of this chapter, with respect to which the operad $\wddot$ of normal wiring diagrams (Prop. \ref{prop:without-dn-operads}) and the operad $\uwd$ of undirected wiring diagrams (Theorem \ref{uwd-operad}) are defined.  Recall that a normal wiring diagram is a wiring diagram without delay nodes. The purpose of this section is to construct a map of operads $\wddot \to \uwd$ given by forgetting directions (Theorem \ref{wddot-uwd-operad-map}).  The existence of such a map of operads was hinted at in the discussion in \cite{rupel-spivak} Section 4.1.

\begin{example}
\label{ex:wddot-to-uwd-motivation}
To motivate the definition of the map of operads $\wddot \to \uwd$ to be defined below, consider the following normal wiring diagram.
\begin{center}
\begin{tikzpicture}[scale=0.6]
\draw [ultra thick] (1,0.2) rectangle (6,4.8);
\node at (3.5,5.3) {$\varphi \in \wddot\yxonextwo$};
\draw [ultra thick] (3,3) rectangle (4,4.1);
\node at (3.5,3.5) {$X_1$};
\draw [ultra thick] ((3,.5) rectangle (4,1.5);
\node at (3.5,1) {$X_2$};
\draw [thick] (0,3.5) -- (1,3.5);
\node at (0.5,3.8) {\tiny{$y_1$}};
\draw [arrow] (0,1) -- (.95,1);
\node at (0.5,1.3) {\tiny{$y_2$}};
\draw [arrow] (6,3.8) -- (7,3.8);
\node at (6.5,4.1) {\tiny{$y^1$}};
\draw [arrow] (6,2) -- (7,2);
\node at (6.5,2.3) {\tiny{$y^2$}};
\draw [arrow] (6,1) -- (7,1);
\node at (6.5,.6) {\tiny{$y^3$}};
\draw [arrow] (1,3.5) -- (2.95,3.5);
\draw [arrow, thick] (1.5,3.5) to [out=0, in=180] (2.95,1);
\draw [thick] (4,1) -- (6,1);
\draw [thick] (4.5,1) to [out=0, in=180] (6,2);
\draw [arrow, looseness=1.4] (4.5,1) to [out=0, in=180] (2.95,3.2);
\draw [arrow] (4,3.2) -- (5,3.2);
\draw [thick] (4,3.8) -- (6,3.8);
\draw [arrow, looseness=5] (4,3.8) to [out=30, in=150] (2.95,3.8);
\end{tikzpicture}
\end{center}
Here
\[\xin_1 = \left\{\smallxin_{11}, \smallxin_{12}, \smallxin_{13}\right\},\quad
\xout_1 = \left\{\smallxout_{11}, \smallxout_{12}\right\},\quad
\xin_2 = \left\{\smallxin_2\right\}, \andspace \xout_2 = \left\{\smallxout_2\right\}.\]
The supplier assignment of $\varphi$ (Def. \ref{def:wiring-diagram}), $s : \dmphi \to \supplyphi$, is given by
\begin{itemize}
\item $y_1 = s\left(\smallxin_{12}\right) = s\left(\smallxin_2\right)$;
\item $\smallxout_{11} = s\left(\smallxin_{11}\right) = s(y^1)$;
\item $\smallxout_2 = s(y^2) = s(y^3) = s\left(\smallxin_{13}\right)$.
\end{itemize}
Note that $y_2 \in \yin$ is an external wasted wire, and $\smallxout_{12} \in \xout_1$ is an internal wasted wire.

A natural way to make $\varphi$ into an undirected wiring diagram is to forget the directions of all the arrows.  For instance, we send 
\[Y = \left(\yin, \yout\right) = \Bigl(\{y_1,y_2\}, \{y^1,y^2,y^3\}\Bigr) \in \boxs\]
to
\[\ybar = \yin \amalg \yout = \left\{y_1,y_2,y^1,y^2,y^3 \right\} \in \Fins,\]
and similarly we send $X_i \in \boxs$ to $\xbar_i = \xin_i \amalg \xout_i \in \Fins$ for $i=1,2$.  So we have
\[\xbar_1 = \left\{\smallxin_{11}, \smallxin_{12}, \smallxin_{13}, \smallxout_{11}, \smallxout_{12}\right\} \andspace
\xbar_2 = \left\{\smallxin_2, \smallxout_2\right\}.\]

Inserting cables at appropriate places, we obtain the following undirected wiring diagram.
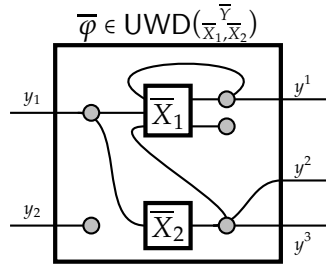
\begin{center}
\begin{tikzpicture}[scale=0.6]
\draw [ultra thick] (1,0.2) rectangle (6,5);
\node at (3.5,5.5) {$\varphibar \in \uwd\yxonextwobar$};
\draw [ultra thick] (3,3) rectangle (4,4.1);
\node at (3.5,3.5) {$\xbar_1$};
\draw [ultra thick] ((3,.5) rectangle (4,1.5);
\node at (3.5,1) {$\xbar_2$};
\draw [thick] (0,3.5) -- (1,3.5);
\node at (0.5,3.8) {\tiny{$y_1$}};
\draw [thick] (0,1) -- (1.8,1);
\node at (0.5,1.3) {\tiny{$y_2$}};
\draw [thick] (6,3.8) -- (7,3.8);
\node at (6.5,4.1) {\tiny{$y^1$}};
\draw [thick] (6,2) -- (7,2);
\node at (6.5,2.3) {\tiny{$y^2$}};
\draw [thick] (6,1) -- (7,1);
\node at (6.5,.6) {\tiny{$y^3$}};
\draw [thick] (1,3.5) -- (3,3.5);
\draw [thick] (1.5,3.5) to [out=0, in=180] (3,1);
\draw [thick] (4,1) -- (6,1);
\draw [thick] (4.5,1) to [out=0, in=180] (6,2);
\draw [thick, looseness=1.4] (4.5,1) to [out=0, in=180] (3,3.2);
\draw [thick] (4,3.2) -- (4.8,3.2);
\draw [thick] (4,3.8) -- (6,3.8);
\draw [thick, looseness=3] (4.8,3.8) to [out=30, in=150] (3,3.8);
\cable{(1.8,3.5)} \cable{(1.8,1)}
\cable{(4.8,3.8)} \cable{(4.8,3.2)} \cable{(4.8,1)} 
\end{tikzpicture}
\end{center}
Call the two cables on the left, from top to bottom, $c_1$ and $c_2$ and the three cables on the right, also from top to bottom, $c_3$, $c_4$, and $c_5$.  Then on  the left side $c_1$ is a $(2,1)$-cable, and $c_2$ is a $(0,1)$-cable.  On the right side,  $c_3$, $c_4$, and $c_5$ are a $(2,1)$-cable, a $(1,0)$-cable, and a $(2,2)$-cable (Def. \ref{def:pre-uwd}), respectively.

Note that the set of cables $\{c_1,\ldots,c_5\}$ in $\varphibar$ is in canonical bijection with the set of supply wires in $\varphi$ (Def. \ref{def:wiring-diagram}), namely
\[\supplyphi = \yin \amalg \xout_1 \amalg \xout_2 \in \Fins.\] 
With this identification, the input and output soldering functions of $\varphibar$ are completely determined by the identity map on $\supplyphi$ and the supplier assignment of $\varphi$.  We will make this precise in \eqref{psibar-cospan} below.
\end{example}

With the previous example as motivation, we now define the map of operads $\wddot \to \uwd$.  Recall the definitions of normal wiring diagrams (Def. \ref{def:wiring-diagram} and \ref{def:wd-without-dn}) and of undirected wiring diagrams (Def. \ref{def:pre-uwd} and \ref{def:uwd}).

\begin{definition}
\label{def:wddot-uwd-map}
Fix a class $S$.
\begin{enumerate}
\item
Define the map $\chi_0 : \boxs \to \Fins$ by
\begin{equation}
\label{wddot-to-uwd-colors}
\chi_0 Y = \ybar = \yin \amalg \yout \in \Fins
\end{equation}
for each $Y = (\yin,\yout) \in \boxs$.
\item
For each $\yxonexn \in \profboxsboxs$ with $n \geq 0$, define the map
\begin{equation}
\label{wddot-to-uwd-entries}
\nicexy{\wddot\yxonexn \ar[r]^-{\chi_1} & \uwd\yxonexnbar}
\end{equation}
as follows.  For $\psi \in \wddot\yxonexn$ (so $\dnpsi = \varnothing$), its image
\[\chi_1 \psi = \psibar \in \uwd\yxonexnbar\]
is the cospan
\begin{equation}
\label{psibar-cospan}
\nicexy@C+1.3cm{& \ybar = \yin \amalg \yout \ar[d]^-{\Id_{\yin} \amalg \spsi|_{\yout}}\\
\coprod\limits_{i=1}^n \xbar_i = \xin \amalg \xout \ar[r]^-{\left(\spsi|_{\xin},\Id_{\xout}\right)} & \supplypsi = \yin \amalg \xout}
\end{equation}
in $\Fins$.  Here:
\begin{itemize}
\item
$\xout = \coprod_{i=1}^n \xout_i$ and $\xin = \coprod_{i=1}^n \xin_i$.
\item
$\supplypsi$ is the set of supply wires of $\psi$.
\item
The map
\[\nicexy{\dmpsi = \yout \amalg \xin\ar[r]^-{\spsi} & \yin \amalg \xout = \supplypsi}\]
is the supplier assignment for $\psi$ (Def. \ref{def:wiring-diagram}).
\end{itemize}
\end{enumerate}
\end{definition}

\begin{remark}
Consider the output soldering function of $\psibar$ in \eqref{psibar-cospan}.  Due to the non-instantaneity requirement \eqref{non-instant}, the restriction of the supplier assignment $\spsi$ to $\yout$ is a map $\yout \to \xout$.
\end{remark}

\begin{theorem}
\label{wddot-uwd-operad-map}
The maps $\chi_0$ \eqref{wddot-to-uwd-colors} and $\chi_1$ \eqref{wddot-to-uwd-entries} define a map of operads
\begin{equation}
\label{wddot-uwd-chi}
\nicexy{\wddot \ar[r]^-{\chi} & \uwd}.
\end{equation}
\end{theorem}

\begin{proof}
As before we will write both $\chi_0$ and$\chi_1$ as $\chi$.  We must check that $\chi$ preserves the operad structure in the sense of Def. \ref{def:operad-map}.  In both $\wddot$ \eqref{wd-permutation} and $\uwd$ \eqref{uwd-right-action}, the equivariant structure is given by permuting the labels of the input boxes.  So $\chi$ preserves equivariance in the sense of \eqref{operad-map-eq}.  Likewise, it follows immediately from the definitions of the colored units in $\wddot$ \eqref{wd-unit} and $\uwd$ \eqref{uwd-unit} that they are preserved by $\chi$ in the sense of \eqref{operad-map-unit}.

To check that $\chi$ preserves operadic composition in the sense of \eqref{operad-map-comp}, suppose $\varphi \in \wddot\yxonexn$ with $n \geq 1$, $1 \leq i \leq n$, and $\psi \in \wddot \xiwonewm$ with $m \geq 0$.  We must show that
\begin{equation}
\label{chi-preserves-operad-comp}
\chi\Bigl(\phicompipsi\Bigr) = \bigl(\chi \varphi\bigr) \compi \bigl(\chi\psi\bigr) \in \uwd\yzbar
\end{equation}
in which 
\[\zbar = \overline{\left(\uX \compi \uW\right)} = \Bigl(\xbar_1,\ldots, \xbar_{i-1}, \wbar_1,\ldots, \wbar_m, \xbar_{i+1}, \ldots, \xbar_n \Bigr) \in \proffins\]
as in \eqref{compi-profile}, $\uX = (X_1,\ldots,X_n)$, and $\uW = (W_1,\ldots,W_m)$.

To prove \eqref{chi-preserves-operad-comp}, on the one hand, by Def. \ref{def:compi-wd} $\phicompipsi \in \wddot \yxcompiw$ has supplier assignment
\[\nicexy{
\dmphicompipsi = \yout \amalg \coprod\limits_{j\not= i} \xin_j \amalg
\win \ar[r]^-{\sphicompipsi} &
\yin \amalg \coprod\limits_{j\not= i} \xout_j \amalg \wout = \supplyphicompipsi}\]
that is given by $\sphi$, $\spsi\sphi$, $\spsi$, $\sphi\spsi$, or $\spsi\sphi\spsi$ according to \eqref{compi-supply1} and \eqref{compi-supply2}.  Here
\[\win = \coprod_{k=1}^m \win_k \andspace \wout = \coprod_{k=1}^m \wout_k \in \Fins.\]
So by \eqref{psibar-cospan}  $\chi\left(\phicompipsi\right) \in \uwd\yzbar$ is the cospan
\begin{equation}
\label{chi-phicompipsi-cospan}
\nicexy@C+.5cm{& \ybar = \yin \amalg \yout \ar[d]^-{\Id \amalg \sphicompipsi}\\
\coprod\limits_{j\not=i} \xbar_j \amalg \coprod\limits_k \wbar_k \ar[r]^-{\left(\sphicompipsi, \Id\right)} & \supplyphicompipsi = \yin \amalg \coprod\limits_{j\not=i} \xout_j \amalg \wout}
\end{equation}
in $\Fins$.  Here:
\begin{itemize}
\item
The input soldering function is made up of
\begin{itemize}
\item
the identity map on $\coprod_{j\not=i} \xout_j \amalg \wout$;
\item
the supplier assignment $\sphicompipsi : \coprod_{j\not=i} \xin_j \amalg \win \to \supplyphicompipsi$.
\end{itemize}
\item
The output soldering function is the coproduct of
\begin{itemize}
\item
the identity map on $\yin$;
\item
the supplier assignment $\sphicompipsi : \yout \to  \coprod_{j\not=i} \xout_j \amalg \wout$.
\end{itemize}
\end{itemize}

On the other hand, by \eqref{psibar-cospan}  and Def. \ref{def:compi-uwd}, the $\compi$-composition
\[\bigl(\chi \varphi\bigr) \compi \bigl(\chi\psi\bigr) \in \uwd\yzbar\]
is the cospan
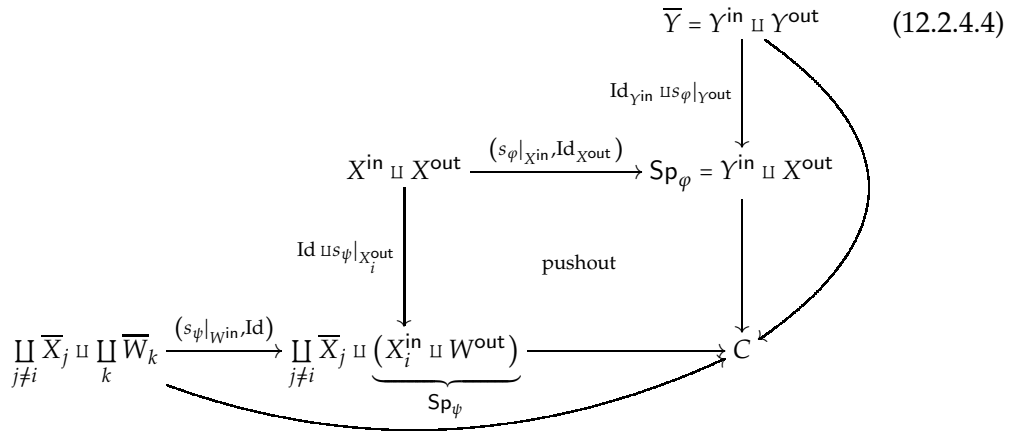
\begin{equation}
\label{chi-phicompipsi-cospan2}
\begin{small}
\nicexy@C+.5cm@R+.5cm{&& \ybar = \yin \amalg \yout \ar[d]_-{\Id_{\yin} \amalg \sphi|_{\yout}} \ar@/^4pc/[dd]\\
& \xin \amalg \xout \ar[d]_-{\Id \amalg \spsi\vert_{\xout_i}} \ar[r]^-{\left(\sphi|_{\xin}, \Id_{\xout}\right)} \pushout & \supplyphi = \yin \amalg \xout \ar[d]\\
\coprod\limits_{j\not=i} \xbar_j  \amalg \coprod\limits_k \wbar_k \ar[r]^-{\left(\spsi|_{\win}, \Id\right)} \ar@/_2.5pc/[rr] & \coprod\limits_{j\not=i} \xbar_j  \amalg \underbrace{\left(\xin_i \amalg \wout\right)}_{\supplypsi} \ar[r] & C}
\end{small}
\end{equation}
\bigskip
in $\Fins$, in which the square is defined as a pushout.  In this diagram:
\begin{itemize}
\item
In the middle vertical map, the restriction to 
\begin{itemize}
\item
$\xout_i$ is the supplier assignment $\spsi : \xout_i \to \wout$;
\item
all other coproduct summands is the identity map.
\end{itemize}
\item
In the bottom left horizontal map, the restriction to
\begin{itemize}
\item
$\win$ is the supplier assignment $\spsi : \win \to \supplypsi$;
\item
all other coproduct summands is the identity map.
\end{itemize}
\end{itemize}

Therefore, to check the condition \eqref{chi-preserves-operad-comp}, it suffices to show that the cospans \eqref{chi-phicompipsi-cospan} and \eqref{chi-phicompipsi-cospan2} are the same.  So it is enough to show that the square
\begin{equation}
\label{chi-operadmap-pushout}
\nicexy@R+1cm{\xin \amalg \coprod\limits_{j\not=i} \xout_j \amalg \xout_i = \xin \amalg \xout \ar[d]_-{\left(\Id_{\xin \amalg \coprod_{j\not=i} \xout_j}\right) \amalg \left(\spsi\vert_{\xout_i}\right)} \ar[r]^-{\left(\sphi\vert_{\xin}, \Id_{\xout}\right)} 
& \supplyphi = \yin \amalg \xout \ar[d]_-{\left(\Id_{\yin \amalg \coprod_{j\not=i} \xout_j}\right) \amalg \left(\spsi\vert_{\xout_i}\right)}\\
\xin \amalg \coprod\limits_{j\not=i} \xout_j \amalg \wout = \coprod\limits_{j\not=i} \xbar_j  \amalg \left(\xin_i \amalg \wout\right) \ar[r]^-{h} & \supplyphicompipsi = \yin \amalg \coprod\limits_{j\not=i} \xout_j \amalg \wout}
\end{equation}
in $\Fins$ is a pushout (Def. \ref{def:pushout}).  Here the restriction of the map $h$ to
\begin{itemize}
\item
$\coprod\limits_{j\not=i} \xout_j \amalg \wout$ is the identity map;
\item
$\xin$ is the composition
\[\nicexy@C+.6cm{\xin \ar[r]^-{\sphi} & \yin \amalg \coprod\limits_{j\not=i} \xout_j \amalg \xout_i \ar[r]^-{\Id \amalg \spsi\vert_{\xout_i}} & \yin \amalg \coprod\limits_{j\not=i} \xout_j \amalg \wout}.\]
\end{itemize}
To see that the square \eqref{chi-operadmap-pushout} is commutative, observe that both compositions when restricted to
\begin{itemize}
\item
$\xin$ is $\Bigl(\Id \amalg \spsi\vert_{\xout_i}\Bigr) \comp \sphi$;
\item
$\coprod\limits_{j\not=i} \xout_j$ is the identity map;
\item
$\xout_i$ is $\spsi$.
\end{itemize}

It remains to check the condition \eqref{pushout-universal-property} of a pushout.  So suppose given a solid-arrow commutative diagram
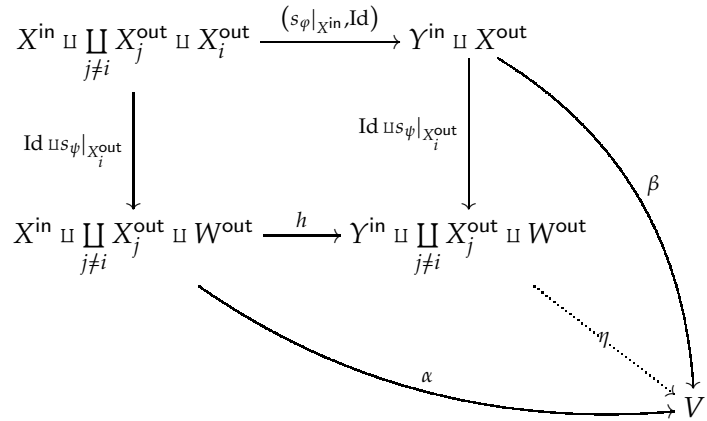
\[\nicexy@R+.7cm{\xin \amalg \coprod\limits_{j\not=i} \xout_j \amalg \xout_i \ar[d]_-{\Id \amalg \spsi\vert_{\xout_i}} \ar[r]^-{\left(\sphi\vert_{\xin}, \Id\right)} & \yin \amalg \xout \ar[d]_-{\Id \amalg \spsi\vert_{\xout_i}} \ar@/^2pc/[ddr]^-{\beta} & \\
\xin \amalg \coprod\limits_{j\not=i} \xout_j \amalg \wout \ar[r]^-{h} \ar@/_2pc/[drr]^-{\alpha} & \yin \amalg \coprod\limits_{j\not=i} \xout_j \amalg \wout \ar@{.>}[dr]|-{\eta}& \\ & & V}\]
in $\Fins$.  We must show that there exists a unique dotted map $\eta$ that makes the diagram commutative.  By a direct inspection the only possible candidate for $\eta$ is given by the restrictions
\[\nicexy{\yin \amalg \coprod\limits_{j\not=i} \xout_j \ar[r]^-{\eta ~=~ \beta} & V} \andspace \nicexy{\wout \ar[r]^-{\eta ~=~ \alpha} &  V}.\]
With this definition of $\eta$, it remains to check  the equalities
\begin{equation}
\label{wddot-uwd-check-equalities}
\eta \left(\Id \amalg \spsi\vert_{\xout_i}\right) = \beta \andspace \eta h = \alpha.
\end{equation}
Both of these equalities can be checked by a direct inspection.  This finishes the proof that the square \eqref{chi-operadmap-pushout} is a pushout and, therefore, that $\chi$ preserves operadic composition \eqref{chi-preserves-operad-comp}.
\end{proof}

\begin{example}
By Example \ref{ex:operad-map-algebra} and Theorem \ref{wddot-uwd-operad-map}, every $\uwd$-algebra (Def. \ref{def:uwd-algebra-biased}) induces a $\wddot$-algebra (Def. \ref{def:normal-algebra}) along the operad map $\chi : \wddot \to \uwd$.  For example:
\begin{itemize}
\item The relational algebra of a set $A$ (Def. \ref{def:relational-algebra-biased} with $S = *$) induces a $\wddot$-algebra along the operad map $\chi$.
\item The typed relational algebra (Def. \ref{def:typed-relational-algebra} with $S = \set$) also induces a $\wddot$-algebra along the operad map $\chi$.
\end{itemize}
\end{example}

\section{Examples of the Operad Map}
\label{sec:examples-normal-undirected}

The purpose of this section is to provide concrete examples of the map of operads $\chi : \wddot \to \uwd$ in Theorem \ref{wddot-uwd-operad-map}.  Recall that the map $\chi$ was defined in Def. \ref{def:wddot-uwd-map}.  Similar to Section  \ref{sec:elementary-relations-uwd}, all the assertions in this section are checked by a direct inspection of the normal and undirected wiring diagrams involved.  So we will omit the proofs.

First we consider the normal generating wiring diagrams (Def. \ref{def:generating-without-dn}) that are sent by $\chi$ to generators in $\uwd$ (Def. \ref{def:generating-uwd}).

\begin{example}
\label{wddot-uwd-empty}
For the empty wiring diagram $\epsilon \in \wddot \emptynothing$ (Def. \ref{def:empty-wd}), the image
\[\chi\epsilon \in \uwd\emptynothing\]
is the empty cell (Def. \ref{def:uwd-gen-empty}). 
\end{example}

\begin{example}
\label{wddot-uwd-namechange}
For a name change $\tau_f \in \wddot\yx$ (Def. \ref{def:name-change}), the image 
\[\chi\tau_f \in \uwd\yxbar\]
is the name change $\tau_{\fbar}$ (Def. \ref{def:uwd-gen-namechange}) corresponding to the bijection
\[\nicexy@C+1cm{\xbar = \xin \amalg \xout \ar[r]^-{\fin \amalg (\fout)^{-1}} & \yin \amalg \yout = \ybar} \in \Fins\]
induced by $f$.
\end{example}

\begin{example}
\label{wddot-uwd-2cell}
For a $2$-cell $\theta_{X,Y} \in \wddot\xplusyxy$ (Def. \ref{def:theta-wd}), the image 
\[\chi\theta_{X,Y} \in \uwd\xplusyxybar\]
is the $2$-cell $\theta_{(\xbar,\ybar)}$ (Def. \ref{def:uwd-gen-2cell}).
\end{example}

\begin{example}
\label{wddot-uwd-loop}
For a $1$-loop $\lambda_{X,x} \in \wddot\xminusxx$ (Def. \ref{def:loop-wd}), the image
\[\chi\lambda_{X,x} \in \uwd\xminusxxbar\]
is the loop $\lambda_{(\xbar, \xsubpm)}$ (Def. \ref{def:uwd-gen-loop}).
\end{example}

\begin{example}
\label{wddot-uwd-outsplit}
For an out-split $\sigma^{Y,y_1,y_2} \in \wddot\yx$ (Def. \ref{def:out-split}), the image
\[\chi\sigma^{Y,y_1,y_2} \in \uwd\yxbar\]
is the split $\sigma^{(\ybar,y_1,y_2)}$ (Def. \ref{def:uwd-gen-split}).
\end{example}

Next we consider an in-split and a $1$-wasted wire.  They are \emph{not} sent by the map of operads $\chi$ to generators in $\uwd$.  So we will express their $\chi$-images as operadic compositions of the generators in $\uwd$.  By Theorem \ref{stratified-presentation-exists-uwd} this is always possible.

\begin{example}
\label{wddot-uwd-insplit}
For an in-split $\sigma_{X,x_1,x_2} \in \wddot\yx$ (Def. \ref{def:insplit-wd}) with $Y = \frac{X}{(x_1 = x_2)}$, suppose
\[Z = \ybar \amalg \bigl\{x_1^+, x_1^-\bigr\} \in \Fins\]
in which $v(x_1^+) = v(x_1^-) = v(x_1) \in S$.  Identify
\[\frac{Z}{(x_{12} = x_1^+)} = \frac{\xbar \amalg \left\{x_1^+, x_1^-\right\}}{(x_1 = x_2 = x_1^+)} = \xbar\]
via $x_1^+ \mapsto x_1$ and $x_1^- \mapsto x_2$.  Then we have
\begin{equation}
\label{chi-insplit}
\chi\sigma_{X,x_1,x_2} = \lambda_{\left(Z, x_1^{\pm}\right)} \comp \sigma^{\left(Z, x_{12}, x_1^+\right)} \in \uwd\yxbar
\end{equation}
in which:
\begin{itemize}
\item $\lambda_{\left(Z, x_1^{\pm}\right)} \in \uwd\ybarz$ is a loop (Def. \ref{def:uwd-gen-loop});
\item $\sigma^{\left(Z, x_{12}, x_1^+\right)} \in \uwd\zxbar$ is a split (Def. \ref{def:uwd-gen-split}).
\end{itemize}
Observe that the right side of \eqref{chi-insplit} also appeared in the elementary relation \eqref{uwd-move-d4} in $\uwd$.  The equality \eqref{chi-insplit} may be visualized as the following picture.
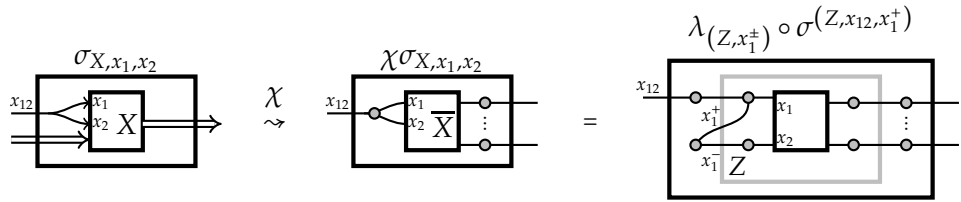
\begin{center}
\begin{tikzpicture}[scale=.7]
\draw [ultra thick] (2,1) rectangle (3,2.1);
\node at (2.7,1.5) {$X$};
\node at (2.2,1.9) {\tiny{$x_1$}};
\node at (2.2,1.5) {\tiny{$x_2$}};
\draw [implies] (.5,1.2) to (2,1.2);
\node at (.7,1.9) {\tiny{$x_{12}$}};
\draw [thick] (.5,1.7) -- (1.2,1.7);
\draw [arrow] (1.2,1.7) to [out=0, in=180] (2,1.9);
\draw [arrow] (1.2,1.7) to [out=0, in=180] (2,1.5);
\draw [implies] (3,1.5) to (4.5,1.5);
\draw [ultra thick] (1,0.7) rectangle (4,2.4);
\node at (2.5,2.7) {$\sigma_{X, x_1, x_2}$};
\node at (5.5,1.5) {$\leadsto$};
\node at (5.5,2) {$\chi$};
\draw [ultra thick] (8,1) rectangle (9,2.1);
\node at (8.7,1.5) {$\xbar$};
\node at (8.2,1.9) {\tiny{$x_1$}};
\node at (8.2,1.5) {\tiny{$x_2$}};
\draw [thick] (6.5,1.7) -- (7.2,1.7);
\node at (6.7,1.9) {\tiny{$x_{12}$}};
\draw [thick] (7.2,1.7) to [out=0, in=180] (8,1.9);
\draw [thick] (7.2,1.7) to [out=0, in=180] (8,1.5);
\draw [thick] (9,1.1) to (10.5,1.1);
\node at (9.5,1.5) {\tiny{$\vdots$}};
\draw [thick] (9,1.9) to (10.5,1.9);
\smallcable{(7.4,1.7)} \smallcable{(9.5,1.9)} \smallcable{(9.5,1.1)}
\draw [ultra thick] (7,0.7) rectangle (10,2.4);
\node at (8.5,2.7) {$\chi\sigma_{X, x_1, x_2}$};
\node at (11.5,1.5) {$=$};
\draw [ultra thick] (15,1) rectangle (16,2.1);
\node at (15.2,1.8) {\tiny{$x_1$}};
\node at (15.2,1.2) {\tiny{$x_2$}};
\draw [ultra thick, lightgray] (14,0.4) rectangle (17,2.4);
\node at (14.3,.7) {\small{$Z$}};
\draw [ultra thick] (13,.1) rectangle (18,2.7);
\node at (15.5,3.3) {$\lambda_{\left(Z, x_1^{\pm}\right)} \comp \sigma^{\left(Z, x_{12}, x_1^+\right)}$};
\draw [thick] (15,2) to (12.5,2);
\node at (12.6,2.2) {\tiny{$x_{12}$}};
\draw [thick] (14.5,1.9) to [out=-90,in=60] (13.5,1.1);
\draw [thick] (13.5,1.1) to (15,1.1);
\node at (13.8,1.65) {\tiny{$x_{1}^+$}};
\node at (13.8,.8) {\tiny{$x_{1}^-$}};
\draw [thick] (16,1.1) to (18.5,1.1);
\node at (17.5,1.5) {\tiny{$\vdots$}};
\draw [thick] (16,1.9) to (18.5,1.9);
\smallcable{(13.5,2)} \smallcable{(13.5,1.1)}
\smallcable{(14.5,2)} \smallcable{(14.5,1.1)}
\smallcable{(16.5,1.9)} \smallcable{(17.5,1.9)} 
\smallcable{(16.5,1.1)} \smallcable{(17.5,1.1)} 
\end{tikzpicture}
\end{center}
On the right side, the intermediate gray box is $Z$.  In $\chi\sigma_{X, x_1, x_2}$ we drew all of $\xbar \setminus \{x_1,x_2\}$ on the right side of the box to make the picture easier to read.  It has a $(2,1)$-cable, and all other cables are $(1,1)$-cables.
\end{example}

Next we consider the $\chi$-image of a $1$-wasted wire.

\begin{example}
\label{wddot-uwd-1wasted}
For a $1$-wasted wire $\omega_{Y,y} \in \wddot\yx$ (Def. \ref{def:wasted-wire-wd}) with $X = Y \setminus y$, we have
\begin{equation}
\label{chi-1wasted}
\chi\omega_{Y,y} = \theta_{(\xbar, y)} \comptwo \omega_y \in \uwd\yxbar
\end{equation}
in which:
\begin{itemize}
\item $\theta_{(\xbar, y)} \in \uwd\xbarplusyxbary$ is a $2$-cell (Def. \ref{def:uwd-gen-2cell});
\item $\omega_y \in \uwd\smallynothing$ is a $1$-output wire (Def. \ref{def:uwd-gen-1output}).
\end{itemize}
Observe that the right side of \eqref{chi-1wasted} also appeared in the elementary relations \eqref{uwd-move-b1} and \eqref{uwd-move-b2} and the example \eqref{wasted-cable-y-simplex-b}.  The equality \eqref{chi-1wasted} may be visualized as the following picture.
\begin{equation}
\label{chi-of-1wasted-wire}
\begin{tikzpicture}[scale=.8]
\draw [ultra thick] (2,1) rectangle (3,2);
\node at (2.5,1.5) {$X$};
\draw [implies] (.5,1.3) to (2,1.3);
\draw [arrow] (.5,1.7) -- (1,1.7);
\node at (.7,1.9) {\tiny{$y$}};
\draw [implies] (3,1.5) to (4.5,1.5);
\draw [ultra thick] (1,0.7) rectangle (4,2.3);
\node at (2.5,2.6) {$\omega_{Y,y}$};
\node at (5.5,1.5) {$\leadsto$};
\node at (5.5,2) {$\chi$};
\draw [ultra thick] (8,1) rectangle (9,2);
\node at (8.5,1.5) {$\xbar$};
\draw [thick] (6.5,1.5) -- (7.5,1.5);
\node at (6.7,1.7) {\tiny{$y$}};
\draw [thick] (9,1.9) to (10.5,1.9);
\node at (9.5,1.5) {\tiny{$\vdots$}};
\draw [thick] (9,1.1) to (10.5,1.1);
\smallcable{(7.5,1.5)} \smallcable{(9.5,1.9)} \smallcable{(9.5,1.1)}
\draw [ultra thick] (7,0.7) rectangle (10,2.3);
\node at (8.5,2.6) {$\chi\omega_{Y,y}$};
\end{tikzpicture}
\end{equation}
In $\chi\omega_{Y,y}$ we drew all of $\xbar$ on the right side to make the picture easier to read.  It has a $(0,1)$-cable, and all other cables are $(1,1)$-cables.
\end{example}

\begin{example}
\label{wddot-uwd-1output}
Each $1$-output wire is in the image of $\chi$.  Indeed, for a $1$-output wire $\omega_y \in \uwd\smallynothing$ (Def. \ref{def:uwd-gen-1output}), we have
\begin{equation}
\label{1output-chi-image}
\begin{split}
\omega_y &= \left[\theta_{(\varnothing,y)} \comptwo \omega_y\right] \comp \epsilon\\
&= \bigl(\chi\omega_{y,y}\bigr) \comp \bigl(\chi \epsilon\bigr) \quad \text{by Examples \ref{wddot-uwd-1wasted} and \ref{wddot-uwd-empty}}\\
&= \chi\bigl(\omega_{y,y} \comp \epsilon\bigr).
\end{split}
\end{equation}
Here:
\begin{itemize}
\item $\theta_{(\varnothing,y)} \in \uwd\smallyemptyy$ is a $2$-cell (Def. \ref{def:uwd-gen-2cell}).
\item $\epsilon \in \uwd\emptynothing$ is the empty cell (Def. \ref{def:uwd-gen-empty}).
\item $\omega_{y,y} \in \wddot\smallyempty$ is the $1$-wasted wire $\omega_{Y,y}$ (Def. \ref{def:wasted-wire-wd}) with $\yout = \varnothing$ and $\yin = y$.
\item
$\epsilon \in \wddot\emptynothing$ is the empty wiring diagram (Def. \ref{def:empty-wd}).
\end{itemize}
Note that in the first two lines of \eqref{1output-chi-image}, the operadic compositions are in the operad $\uwd$.  On the other hand, in the last line of \eqref{1output-chi-image} the operadic composition is in the operad $\wddot$.
\end{example}

\begin{example}
\label{uwd-gen-chi-image}
By Examples \ref{wddot-uwd-empty}--\ref{wddot-uwd-outsplit} and \ref{wddot-uwd-1output}, all $6$ types of generators in $\uwd$ (Def. \ref{def:generating-uwd}) are in the image of the operad map $\chi : \wddot \to \uwd$.  However, one must be careful that this does \emph{not} imply that every undirected wiring diagram is in the image of $\chi$.  We will make precise the image of the operad map $\chi$ in Theorem \ref{thm:chi-image} below.
\end{example}

Next we consider a $1$-internal wasted wire, which by Prop. \ref{prop:internal-wasted-wire} can be generated by a $1$-loop and a $1$-wasted wire.  

\begin{example}
\label{wddot-uwd-1internalwasted}
For a $1$-internal wasted wire $\omega^{X,x} \in \wddot\yx$ (Def. \ref{def:internal-wasted-wire}) with $Y = X \setminus x$, we have
\begin{equation}
\label{chi-1internal}
\chi\omega^{X,x} = \lambda_{\left(Z, \xsubpm\right)} \comp \sigma^{\left(Z, \xsubpm\right)} \in \uwd\yxbar
\end{equation}
in which:
\begin{itemize}
\item $Z = \ybar \amalg \xsubpm \in \Fins$ with $v(x_+) = v(x_-) = v(x) \in S$;
\item $\lambda_{\left(Z, \xsubpm\right)} \in \uwd\ybarz$ is a loop (Def. \ref{def:uwd-gen-loop});
\item $\sigma^{\left(Z, \xsubpm\right)} \in \uwd\zxbar$ is a split (Def. \ref{def:uwd-gen-split}).
\end{itemize}
Observe that the right side of \eqref{chi-1internal} also appeared in the elementary relation \eqref{uwd-move-b1}.   The equality \eqref{chi-1internal} may be visualized as the following picture.
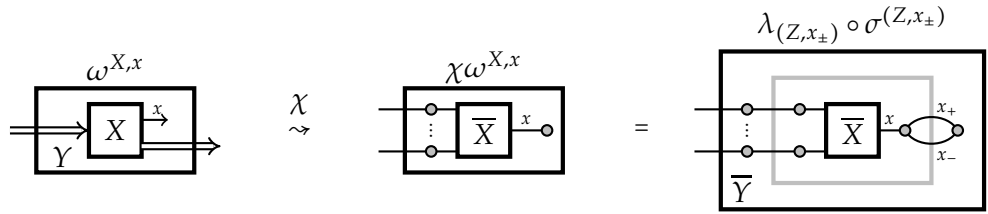
\begin{center}
\begin{tikzpicture}[scale=.7]
\draw [ultra thick] (2,1) rectangle (3,2);
\node at (2.5,1.5) {$X$};
\draw [implies] (.5,1.5) to (2,1.5);
\node at (1.5,1) {$Y$};
\draw [arrow] (3,1.7) -- (3.5,1.7);
\node at (3.3,1.9) {\tiny{$x$}};
\draw [implies] (3,1.2) to (4.5,1.2);
\draw [ultra thick] (1,0.7) rectangle (4,2.3);
\node at (2.5,2.7) {$\omega^{X,x}$};
\node at (6,1.5) {$\leadsto$};
\node at (6,2) {$\chi$};
\draw [ultra thick] (8,0.7) rectangle (11,2.3);
\node at (9.5,2.7) {$\chi\omega^{X,x}$};
\draw [ultra thick] (9,1) rectangle (10,2);
\node at (9.5,1.5) {$\xbar$};
\draw [thick] (9,1.1) -- (7.5,1.1);
\draw [thick] (9,1.9) -- (7.5,1.9);
\draw [thick] (10,1.5) -- (10.7,1.5);
\node at (10.3,1.7) {\tiny{$x$}};
\smallcable{(8.5,1.1)} \node at (8.5,1.5) {\tiny{$\vdots$}}; \smallcable{(8.5,1.9)}
\smallcable{(10.7,1.5)}
\node at (12.5,1.5) {$=$};
\draw [ultra thick] (14,0) rectangle (19,3);
\node at (14.4,.4) {$\ybar$};
\node at (16.5,3.5) {$\lambda_{\left(Z, \xsubpm\right)} \comp \sigma^{\left(Z, \xsubpm\right)}$};
\draw [ultra thick, lightgray] (15,.5) rectangle (18,2.5);
\draw [ultra thick] (16,1) rectangle (17,2);
\node at (16.5,1.5) {$\xbar$};
\draw [thick] (16,1.9) to (13.5,1.9);
\draw [thick] (16,1.1) to (13.5,1.1);
\node at (14.5,1.5) {\tiny{$\vdots$}};
\draw [thick] (17,1.5) to (17.5,1.5);
\node at (17.2,1.7) {\tiny{$x$}};
\draw [thick, bend left=60] (17.5,1.5) to (18.5,1.5);
\draw [thick, bend right=60] (17.5,1.5) to (18.5,1.5);
\node at (18.3,1.9) {\tiny{$x_+$}};
\node at (18.3,1) {\tiny{$x_-$}};
\smallcable{(14.5,1.9)} \smallcable{(14.5,1.1)}
\smallcable{(15.5,1.9)} \smallcable{(15.5,1.1)}
\smallcable{(17.5,1.5)} \smallcable{(18.5,1.5)}
\end{tikzpicture}
\end{center}
On the right side, the gray box is $Z$.  In $\chi\omega^{X,x}$, all of $\xbar \setminus x = \ybar$ is drawn on the left side.  It has a $(1,0)$-cable, and all other cables are $(1,1)$-cables.
\end{example}

\begin{example}
Consider the wiring diagram $\pi \in \wddot\yx$ in Example \ref{ex:factoring-pi}.  Then $\chi\pi \in \uwd\yxbar$ is the right side of the following picture.
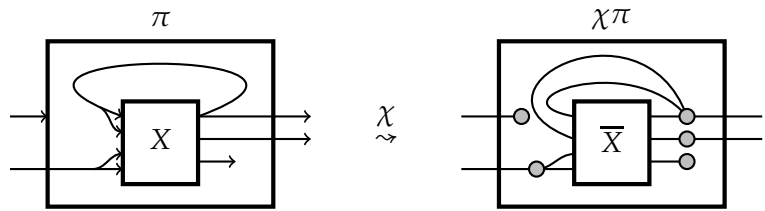
\begin{center}
\begin{tikzpicture}
\draw [ultra thick] (1,1.9) rectangle (2,3);
\node at (1.5,2.5) {$X$};
\draw [arrow, looseness=7] (2,2.8) to [out=20, in=160] (1,2.8);
\draw [arrow] (.7,2.92) to [out=-30, in=180] (1,2.6);
\draw [arrow] (2,2.8) to (3.5,2.8);
\draw [arrow] (2,2.5) to (3.5,2.5);
\draw [arrow] (2,2.2) to (2.5,2.2);
\draw [ultra thick] (0,1.6) rectangle (3,3.8);
\node at (1.5,4.1) {$\pi$};
\draw [arrow] (-.5,2.8) to (0,2.8);
\draw [arrow] (-.5,2.1) to (1,2.1);
\draw [arrow] (.6,2.1) to [out=0,in=180] (1,2.3);
\node at (4.5,2.5) {$\leadsto$};
\node at (4.5,2.8) {$\chi$};
\draw [ultra thick] (7,1.9) rectangle (8,3);
\node at (7.5,2.5) {$\xbar$};
\draw [thick, looseness=3] (8.5,2.8) to [out=110, in=160] (7,2.5);
\draw [thick, looseness=2] (8.5,2.8) to [out=130, in=170] (7,2.8);
\draw [thick] (8,2.8) to (9.5,2.8);
\draw [thick] (8,2.5) to (9.5,2.5);
\draw [thick] (8,2.2) to (8.5,2.2);
\draw [ultra thick] (6,1.6) rectangle (9,3.8);
\node at (7.5,4.1) {$\chi\pi$};
\draw [thick] (5.5,2.8) to (6.3,2.8);
\draw [thick] (5.5,2.1) to (7,2.1);
\draw [thick] (6.5,2.1) to [out=0,in=180] (7,2.3);
\smallcable{(6.3,2.8)} \smallcable{(6.5,2.1)}
\smallcable{(8.5,2.8)} \smallcable{(8.5,2.5)} \smallcable{(8.5,2.2)}
\end{tikzpicture}
\end{center}
So $\chi\pi$ has a $(0,1)$-cable, a $(2,1)$-cable, a $(3,1)$-cable, a $(1,1)$-cable, and a $(1,0)$-cable.
\end{example}

\begin{example}
Consider the wiring diagram $\pi_2 \in \wddot\zx$ in Example \ref{ex:factor-pitwo}.  Then $\chi\pi_2 \in \uwd\zbarxbar$ is the right side of the following picture.
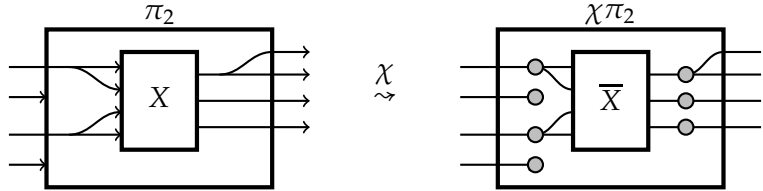
\begin{center}
\begin{tikzpicture}[xscale=1]
\draw [ultra thick] (1,1.8) rectangle (2,3.1);
\node at (1.5,2.5) {$X$};
\draw [arrow] (-.5,2.9) to (1,2.9);
\draw [arrow] (-.5,2) to (1,2);
\draw [arrow] (2,2.8) to (3.5,2.8);
\draw [arrow] (2,2.45) to (3.5,2.45);
\draw [arrow] (2,2.1) to (3.5,2.1);
\draw [ultra thick] (0,1.3) rectangle (3,3.4);
\node at (1.5,3.6) {$\pi_2$};
\draw [arrow] (-.5,2.5) to (0,2.5);
\draw [arrow] (-.5,1.6) to (0,1.6);
\draw [thick] (2.3,2.8) to [out=0, in=180] (3,3.1);
\draw [arrow] (3,3.1) to (3.5,3.1);
\draw [arrow] (.3,2.9) to [out=0, in=180] (1,2.6);
\draw [arrow] (.3,2) to [out=0, in=180] (1,2.3);
\node at (4.5,2.5) {$\leadsto$};
\node at (4.5,2.8) {$\chi$};
\draw [ultra thick] (7,1.8) rectangle (8,3.1);
\node at (7.5,2.5) {$\xbar$};
\draw [thick] (5.5,2.9) to (7,2.9);
\draw [thick] (5.5,2) to (7,2);
\draw [thick] (8,2.8) to (9.5,2.8);
\draw [thick] (8,2.45) to (9.5,2.45);
\draw [thick] (8,2.1) to (9.5,2.1);
\draw [ultra thick] (6,1.3) rectangle (9,3.4);
\node at (7.5,3.6) {$\chi\pi_2$};
\draw [arrow] (5.5,2.5) to (6.5,2.5);
\draw [arrow] (5.5,1.6) to (6.5,1.6);
\draw [thick] (8.5,2.8) to [out=0, in=180] (9,3.1);
\draw [thick] (9,3.1) to (9.5,3.1);
\draw [thick] (6.5,2.9) to [out=0, in=180] (7,2.6);
\draw [thick] (6.5,2) to [out=0, in=180] (7,2.3);
\smallcable{(6.5,2.9)} \smallcable{(6.5,2.5)} \smallcable{(6.5,2)} \smallcable{(6.5,1.6)}
\smallcable{(8.5,2.8)} \smallcable{(8.5,2.45)} \smallcable{(8.5,2.1)}
\end{tikzpicture}
\end{center}
So $\chi\pi_2$ has two $(2,1)$-cables, two $(0,1)$-cables, a $(1,2)$-cable, and two $(1,1)$-cables.
\end{example}

\section{Image of the Operad Map}
\label{sec:image-normal-undirected}

The purpose of this section is to give an explicit description of the image of the map of operads $\chi : \wddot \to \uwd$ in Theorem \ref{wddot-uwd-operad-map}.  Recall that the map $\chi$ was defined in Def. \ref{def:wddot-uwd-map}.  Also recall the notations and terminology in Notation \ref{notation:cable-subsets} regarding subsets of cables.

\begin{theorem}
\label{thm:chi-image}
Consider the operad map $\chi : \wddot \to \uwd$ in Theorem \ref{wddot-uwd-operad-map}.  Then:
\begin{enumerate}
\item
The color map $\chi_0 : \boxs \to \Fins$ \eqref{wddot-to-uwd-colors} is surjective.
\item
The image of the entry map $\chi_1 : \wddot \to \uwd$ \eqref{wddot-to-uwd-entries} consists of precisely the undirected wiring diagrams with no wasted cables and no $(0,\geq 2)$-cables.
\end{enumerate}
\end{theorem}

\begin{proof}
The color map $\chi_0$ is surjective because, for each $X \in \Fins$, we have $(\varnothing, X) \in \boxs$ and $\chi(\varnothing,X) = X$.

For the second assertion, we will prove the required inclusions in both directions.  First suppose $\psi \in \wddot\yxonexn$ for some $n \geq 0$.  Recall that $\chi\psi = \psibar \in \uwd\yxonexnbar$ is the cospan \eqref{psibar-cospan}
\[\nicexy@C+1.3cm{& \ybar = \yin \amalg \yout \ar[d]^-{\Id_{\yin} \amalg \spsi|_{\yout}}\\
\coprod\limits_{i=1}^n \xbar_i = \xin \amalg \xout \ar[r]^-{\left(\spsi|_{\xin},\Id_{\xout}\right)} & \supplypsi = \yin \amalg \xout}\]
in $\Fins$.  To see that $\psibar$ has no wasted cables (i.e., $(0,0)$-cables) and no $(0,\geq 2)$-cables, suppose $c \in \supplypsi$ is \emph{not} in the image of the input soldering function $\left(\spsi|_{\xin},\Id_{\xout}\right)$.  We must show that $c$ is a $(0,1)$-cable in $\psibar$.  Since $c$ is not in the image of $\Id_{\xout}$, we have $c \in \yin$.  By the non-instantaneity requirement \eqref{non-instant}, we also have $c \not\in \spsi(\yout)$.  Therefore, $c$ is in the image of the output soldering function $\Id_{\yin} \amalg \spsi|_{\yout}$ exactly once, so $c$ is a $(0,1)$-cable in $\psibar$.

To improve readability, the other half of the second assertion--i.e., that every undirected wiring diagram with no wasted cables and no $(0,\geq 2)$-cables is in the image of the operad map $\chi$--will be proved in Proposition \ref{prop:chi-image} below.
\end{proof}

First we consider the special case when there are no $(0,\geq 0)$-cables.

\begin{definition}
\label{def:chi-image-input-surjective}
Suppose
\begin{equation}
\label{wddot-uwd-phi-cospan}
\varphi = \Bigl(\nicexy{\coprod\limits_{j=1}^N U_j \ar[r]^-{\fphi} & \cphi & V \ar[l]_-{\gphi}}\Bigr) \in \uwd\vuoneun
\end{equation}
with $N \geq 0$ and $\cphizeroatleastzero = \varnothing$. 
\begin{enumerate}
\item
For each cable $c \in \cphi$, pick a wire
\begin{equation}
\label{usubc-fphi-preimage}
u_c \in \fphiinv(c) \subseteq \coprod_{j=1}^N U_j
\end{equation}
in the $\fphi$-preimage of $c$.  This is possible because the assumption $\cphizeroatleastzero = \varnothing$ means exactly that $\fphi$ is surjective.  We will use the canonical bijection
\begin{equation}
\label{xout-is-cphi}
\nicexy{\Bigl\{u_c : c \in \cphi\Bigr\} \ar[r]^-{\fphi}_-{\cong} & \cphi} \in \Fins
\end{equation}
below.
\item
For each $1 \leq j \leq N$ define a box $X_j \in \boxs$ as
\[\xout_j = \Bigl\{u_c \in U_j : c \in \cphi\Bigr\} \subseteq U_j \andspace 
\xin_j = U_j \setminus \xout_j.\]
Note that we have $\xbar_j = U_j$,
\[\xout = \coprod_{j=1}^N \xout_j = \bigl\{u_c : c \in \cphi\bigr\} \cong \cphi, \andspace
\xin = \coprod_{j=1}^N U_j \setminus \bigl\{u_c : c \in \cphi\bigr\}.\]
\item
Define $Y = (\varnothing, V) \in \boxs$, so $\ybar = V$.
\item
Using the bijection \eqref{xout-is-cphi}, define $\psi \in \wddot\yxonexbign$ whose supplier assignment
\begin{equation}
\label{wddot-uwd-psi-supplier}
\nicexy{\dmpsi = \yout \amalg \xin = V \amalg \left[\coprod\limits_{j=1}^N U_j \setminus \bigl\{u_c : c \in \cphi\bigr\}\right] \ar[d]_-{\spsi}^-{=\, (\gphi, \fphi)}\\
\supplypsi = \yin \amalg \xout \cong \cphi}
\end{equation}
is 
\begin{itemize}
\item $\gphi$ when restricted to $V$;
\item $\fphi$ when restricted to $\coprod_j U_j \setminus \{u_c\}$. 
\end{itemize}
\end{enumerate}
\end{definition}

\begin{lemma}
\label{lem1:chi-image}
In the context of Def. \ref{def:chi-image-input-surjective}, we have
\[\chi\psi = \varphi \in \uwd\yxonexbignbar = \uwd\vuoneun.\]
\end{lemma}

\begin{proof}
By definition $\chi\psi$ is the cospan \eqref{psibar-cospan}
\[\nicexy@C+.5cm{& V = \yin \amalg \yout \ar[d]^-{\Id_{\varnothing} \amalg \spsi}\\
\coprod\limits_{j=1}^N U_j = \xin \amalg \xout \ar[r]^-{(\spsi, \Id)} 
& \supplypsi = \yin \amalg \xout \cong \cphi}\]
in $\Fins$.  Using the bijection \eqref{xout-is-cphi} and the definition of $\spsi$ \eqref{wddot-uwd-psi-supplier}, it follows that this cospan is equivalent to the given cospan \eqref{wddot-uwd-phi-cospan} in the sense of Def. \ref{def:uwd}.   So they define the same undirected wiring diagram, i.e., $\chi\psi = \varphi$.
\end{proof}

\begin{example}
Consider $\varphi$ in \eqref{wddot-uwd-phi-cospan} with $N=0$.  Then $\coprod_j U_j = \varnothing$.  Since $\fphi$ is sujective, it follows that $\cphi = V = \varnothing$.  So $\varphi$ is the empty cell $\epsilon \in \uwd\emptynothing$ (Def. \ref{def:uwd-gen-empty}).  The construction \eqref{wddot-uwd-psi-supplier} above yields $\psi = \epsilon \in \wddot\emptynothing$, the empty wiring diagram.  So the conclusion $\chi\epsilon = \epsilon$ agrees with Example \ref{wddot-uwd-empty}.  
\end{example}

\begin{example}
\label{ex:phi-liftedto-psi-prelim}
Suppose $\varphi$ with $\cphizeroatleastzero = \varnothing$ is the following undirected wiring diagram.
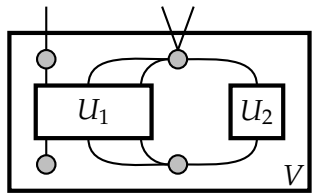
\begin{center}
\begin{tikzpicture}[scale=.7]
\draw [ultra thick] (-.2,0) rectangle (5.5,3);
\node at (5.2,.3) {$V$};
\draw [ultra thick] (.3,1) rectangle (2.5,2);
\node at (1.4,1.5) {$U_1$};
\draw [ultra thick] (4,1) rectangle (5,2);
\node at (4.5,1.5) {$U_2$};
\draw [thick] (.5,2) -- (.5,3.5);
\draw [thick] (1.3,2) to [out=90, in=180] (2.825,2.5);
\draw [thick] (2.3,2) to [out=90, in=180] (2.825,2.5);
\draw [thick] (4.5,2)  to [out=90, in=0] (3.175,2.5);
\draw [thick] (3,2.675) -- (2.7,3.5);
\draw [thick] (3,2.675) -- (3.3,3.5);
\draw [thick] (1.3,1) to [out=-90,in=180] (2.825,.5);
\draw [thick] (2.3,1) to [out=-90,in=180] (2.825,.5);
\draw [thick] (4.5,1) to [out=-90,in=0]  (3.175,.5);
\draw [thick] (.5,1) -- (.5,.675);
\cable{(.5,2.5)} \cable{(3,2.5)} \cable{(3,.5)} \cable{(.5,.5)} 
\end{tikzpicture}
\end{center}
Then one choice of a $\chi$-preimage $\psi \in \wddot$, as constructed in \eqref{wddot-uwd-psi-supplier}, is the following normal wiring diagram.
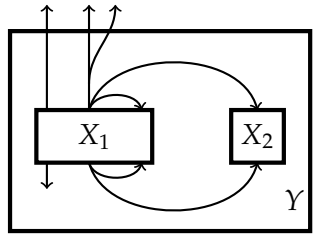
\begin{center}
\begin{tikzpicture}[scale=.7]
\draw [ultra thick] (-.2,-.3) rectangle (5.5,3.5);
\node at (5.2,.3) {$Y$};
\draw [ultra thick] (.3,1) rectangle (2.5,2);
\node at (1.4,1.5) {$X_1$};
\draw [ultra thick] (4,1) rectangle (5,2);
\node at (4.5,1.5) {$X_2$};
\draw [arrow] (.5,2) -- (.5,4);
\draw [arrow] (1.3,2) to [out=75, in=105] (4.5,2);
\draw [arrow] (1.3,2) to [out=75, in=105] (2.3,2);
\draw [arrow] (1.3,2) to (1.3,4);
\draw [arrow] (1.3,2.5) to [out=90,in=-90] (1.8,4);
\draw [arrow] (1.3,1) to [out=-75,in=-105] (4.5,1);
\draw [arrow] (1.3,1) to [out=-75,in=-105] (2.3,1);
\draw [arrow] (.5,1) -- (.5,.5);
\end{tikzpicture}
\end{center}
So $X_1$ has $2$ inputs and $4$ outputs, and $X_2$ has $2$ inputs and no outputs.  Note that according to Convention \ref{conv:box} we should draw inputs on the left and outputs on the right.  However, we drew $\psi$ to resemble $\varphi$ to make the construction easier to understand.
\end{example}

Next we consider the general case where there may be $(0,1)$-cables in $\varphi$.  The idea is to decompose $\varphi$ as $\varphi_1 \compone \varphi_0$ such that the following statements hold.
\begin{itemize}
\item $\varphi_0$ satisfies the hypotheses of Lemma \ref{lem1:chi-image}, so it has no $(0, \geq 0)$-cables.  
\item $\varphi_1$ contains all the $(0,1)$-cables in $\varphi$; its other cables are all $(1,1)$-cables.
\item
Each of $\varphi_0$ and $\varphi_1$ can be lifted back to $\wddot$ in such a way that the lifted wiring diagrams are operadically composable in $\wddot$.
\end{itemize}
The fact that $\chi : \wddot \to \uwd$ is an operad map will then show that $\varphi$ has a $\chi$-preimage.

\begin{definition}
\label{def:chi-image-general-psi}
Suppose
\begin{equation}
\label{wddot-uwd-phi-cospan-general}
\varphi = \Bigl(\nicexy{\coprod\limits_{j=1}^N U_j \ar[r]^-{\fphi} & \cphi & V \ar[l]_-{\gphi}}\Bigr) \in \uwd\vuoneun
\end{equation}
with $N \geq 0$ and $\cphizerozero = \cphizeroatleasttwo = \varnothing$.  By assumption there is a decomposition
\[\cphi = \im(\fphi) \amalg \cphizeroone\]
in which $\im(\fphi)$ is the image of $\fphi$.
\begin{itemize}
\item Define 
\begin{equation}
\label{vzero-vone}
V_0 = \gphiinv\bigl(\im(\fphi)\bigr) \andspace V_1 = \gphiinv\bigl(\cphizeroone\bigr).\end{equation}
So $V = V_0 \amalg V_1$, and there is a bijection $\gphi : V_1 \cong \cphizeroone$.
\item Define
\begin{equation}
\label{chi-image-phi0}
\varphi_0 = \Bigl(\nicexy{\coprod\limits_{j=1}^N U_j \ar[r]^-{\fphi} & \im(\fphi) & V_0 \ar[l]_-{\gphi}}\Bigr) \in \uwd\vzerouoneun,
\end{equation}
in which the input soldering function is surjective, i.e., $\cphizerozeroatleastzero = \varnothing$.
\item Define
\begin{equation}
\label{chi-image-phi1}
\varphi_1 = \Bigl(\nicexy@C+.4cm{V_0 \ar[r]^-{\text{inclusion}} & V_0 \amalg V_1 = V & V \ar[l]_-{=}} \Bigr) \in \uwd\vvzero,
\end{equation}
which has only $1$ input box $V_0$.
\item Define $Y_0 = (\varnothing, V_0) \in \boxs$, so $\ybar_0 = V_0 \in \Fins$.
\item Define $Y = (V_1, V_0) \in \boxs$, so $\ybar = V \in \Fins$.
\end{itemize}
\end{definition}

First we observe that the two undirected wiring diagrams in the previous definition give a decomposition of $\varphi$.

\begin{lemma}
\label{chi-wddot-uwd-phi-decomp}
In the context of Def. \ref{def:chi-image-general-psi}, there is a decomposition
\[\varphi = \varphi_1 \comp \varphi_0.\]
\end{lemma}

\begin{proof}
Since $\varphi$ is the cospan
\[\nicexy@C+.4cm{&& V \ar[d]^-{=} \ar@/^3pc/[dd]^-{\gphi}\\
& V_0 \ar[d]_-{\gphi} \ar[r]^-{\text{inclusion}} & V_0 \amalg V_1 \ar[d]^-{\gphi}\\
\coprod\limits_{j=1}^N U_j \ar[r]^-{\fphi} \ar@/_2pc/[rr]|-{\fphi} & \im(\fphi) \ar[r]^-{\text{inclusion}} & \cphi = \im(\fphi) \amalg \cphizeroone}\]
in $\Fins$, by the definition of $\comp = \compone$ in $\uwd$ (Def. \ref{def:compi-uwd}) it is enough to check that the square is a pushout.  Since  $\gphi : V_1 \cong \cphizeroone$ is a bijection, a direct inspection reveals that the square is a pushout.
\end{proof}

\begin{lemma}
\label{phi1-chi-preimage}
For $\varphi_1 \in \uwd\vvzero$ in \eqref{chi-image-phi1}, there exists $\psi_1 \in \wddot \yyzero$ such that
\[\chi\psi_1 = \varphi_1.\]
\end{lemma}

\begin{proof}
Define $\psi_1 \in \wddot \yyzero$ whose supplier assignment
\begin{equation}
\label{psi1-chi-wddot-uwd}
\nicexy{\dm_{\psi_1} = \yout \amalg \yin_0 = V_0 \ar[d]_-{s_{\psi_1}}^-{=\, \text{inclusion}}\\
\supply_{\psi_1} = \yin \amalg \yout_0 = V_1 \amalg V_0 = V}
\end{equation}
is the inclusion map.  Then it follows from the definition of $\chi$ \eqref{psibar-cospan} that $\chi\psi_1 = \varphi_1$.
\end{proof}

\begin{remark}
By Lemma \ref{lemma:betaone} $\psi_1$ in \eqref{psi1-chi-wddot-uwd} is an iterated operadic composition of $|V_1|$ $1$-wasted wires (Def. \ref{def:wasted-wire-wd}).  Since $V_1 \cong \gphiinv\left(\cphizeroone\right)$ \eqref{vzero-vone}, this means that the $(0,1)$-cables in $\varphi$ are lifted to external wasted wires in $\psi_1$.
\end{remark}

\begin{proposition}
\label{prop:chi-image}
Every undirected wiring diagram $\varphi$ with $\cphizerozero = \cphizeroatleasttwo = \varnothing$ is in the image of the operad map $\chi : \wddot \to \uwd$.
\end{proposition}

\begin{proof}
Suppose $\varphi \in \uwd\vuoneun$ with $\cphizerozero = \cphizeroatleasttwo = \varnothing$.  By Lemma \ref{chi-wddot-uwd-phi-decomp} there is a decomposition
\begin{equation}
\label{chi-phi-is-phi1phi0}
\varphi = \varphi_1 \comp \varphi_0
\end{equation}
with $\varphi_0 \in \uwd\vzerouoneun$ and $\varphi_1 \in \uwd\vvzero$ as in Def. \ref{def:chi-image-general-psi}.  Moreover, $\varphi_0$ satisfies the hypotheses of Lemma \ref{lem1:chi-image} (i.e., that its input soldering function is surjective).  So there exists $\psi_0 \in \wddot\yzeroxonexbign$ such that 
\begin{equation}
\label{chipsi0-is-phi0}
\chi\psi_0 = \varphi_0.
\end{equation}
With $\psi_1 \in \wddot\yyzero$ as in Lemma \ref{phi1-chi-preimage}, we have
\[\begin{split}
\chi\bigl(\psi_1 \comp \psi_0\bigr) 
&= \bigl(\chi\psi_1\bigr) \comp \bigl(\chi\psi_0\bigr) \quad\text{by Theorem \ref{wddot-uwd-operad-map}}\\
&= \varphi_1 \comp \varphi_0 \quad\text{by Lemma \ref{phi1-chi-preimage} and \eqref{chipsi0-is-phi0}}\\
&= \varphi \quad\text{by \eqref{chi-phi-is-phi1phi0}}.
\end{split}\]
This proves that $\varphi$ is in the image of $\chi$.
\end{proof}

Proposition \ref{prop:chi-image} finishes the proof of Theorem \ref{thm:chi-image}.

\begin{example}
\label{ex:phi-liftedto-psi}
This is an extension of Example \ref{ex:phi-liftedto-psi-prelim}.  Suppose $\varphi$ with $\cphizerozero = \cphizeroatleasttwo = \varnothing$ is the following undirected wiring diagram.
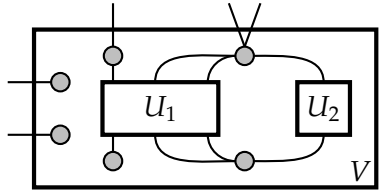
\begin{center}
\begin{tikzpicture}[scale=.7]
\draw [ultra thick] (-1,0) rectangle (5.5,3);
\node at (5.2,.3) {$V$};
\draw [ultra thick] (.3,1) rectangle (2.5,2);
\node at (1.4,1.5) {$U_1$};
\draw [ultra thick] (4,1) rectangle (5,2);
\node at (4.5,1.5) {$U_2$};
\draw [thick] (.5,2) -- (.5,3.5);
\draw [thick] (1.3,2) to [out=90, in=180] (2.825,2.5);
\draw [thick] (2.3,2) to [out=90, in=180] (2.825,2.5);
\draw [thick] (4.5,2)  to [out=90, in=0] (3.175,2.5);
\draw [thick] (3,2.675) -- (2.7,3.5);
\draw [thick] (3,2.675) -- (3.3,3.5);
\draw [thick] (1.3,1) to [out=-90,in=180] (2.825,.5);
\draw [thick] (2.3,1) to [out=-90,in=180] (2.825,.5);
\draw [thick] (4.5,1) to [out=-90,in=0]  (3.175,.5);
\draw [thick] (.5,1) -- (.5,.675);
\draw [thick] (-.5,2) to (-1.5,2);
\draw [thick] (-.5,1) to (-1.5,1);
\cable{(-.5,2)} \cable{(-.5,1)}
\cable{(.5,2.5)} \cable{(3,2.5)} \cable{(3,.5)} \cable{(.5,.5)} 
\end{tikzpicture}
\end{center}
Then one choice of a $\chi$-preimage $\psi \in \wddot$, as constructed in Prop. \ref{prop:chi-image}, is the following normal wiring diagram.
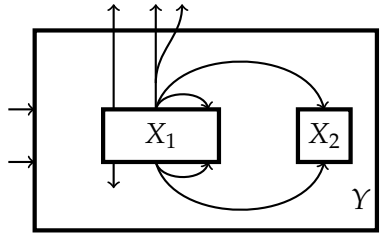
\begin{center}
\begin{tikzpicture}[scale=.7]
\draw [ultra thick] (-1,-.3) rectangle (5.5,3.5);
\node at (5.2,.3) {$Y$};
\draw [ultra thick] (.3,1) rectangle (2.5,2);
\node at (1.4,1.5) {$X_1$};
\draw [ultra thick] (4,1) rectangle (5,2);
\node at (4.5,1.5) {$X_2$};
\draw [arrow] (.5,2) -- (.5,4);
\draw [arrow] (1.3,2) to [out=75, in=105] (4.5,2);
\draw [arrow] (1.3,2) to [out=75, in=105] (2.3,2);
\draw [arrow] (1.3,2) to (1.3,4);
\draw [arrow] (1.3,2.5) to [out=90,in=-90] (1.8,4);
\draw [arrow] (1.3,1) to [out=-75,in=-105] (4.5,1);
\draw [arrow] (1.3,1) to [out=-75,in=-105] (2.3,1);
\draw [arrow] (.5,1) -- (.5,.5);
\draw [arrow] (-1.5,2) to (-1,2);
\draw [arrow] (-1.5,1) to (-1,1);
\end{tikzpicture}
\end{center}
Note that the two $(0,1)$-cables in $\varphi$ are lifted to external wasted wires in $\psi$.
\end{example}

\section{Map from Strict to Undirected Wiring Diagrams}
\label{sec:strict-to-undirected}

For a fixed class $S$, recall the $\boxs$-colored operad of strict wiring diagrams $\wdzero$ (Prop. \ref{prop:strict-operads}).  As pointed out in Example \ref{ex:initial-operad}, we can compose the operad map $\chi : \wddot \to \uwd$ in Theorem \ref{wddot-uwd-operad-map} with the operad inclusion $\wdzero \to \wddot$ in Prop. \ref{prop:strict-normal-wd} to obtain an operad map
\begin{equation}
\label{chizero}
\nicexy{\wdzero \ar[r] \ar@/^2pc/[rr]^-{\chizero} & \wddot \ar[r]^-{\chi} & \uwd}.
\end{equation}
The purpose of this section is to identify precisely the image of this operad map.

\begin{theorem}
\label{chizero-wdzero-uwd}
The image of the operad map $\chizero : \wdzero \to \uwd$ consists of precisely the undirected wiring diagrams whose cables are either $(1,1)$-cables or $(2,0)$-cables.
\end{theorem}

\begin{proof}
To make the argument easier to read, we will prove the two required inclusions in Lemmas \ref{lem1:wdzero-uwd} and \ref{lem2:wdzero-uwd} below.
\end{proof}

\begin{lemma}
\label{lem1:wdzero-uwd}
The image of each strict wiring diagram under the operad map $\chizero : \wdzero \to \uwd$ has only $(1,1)$-cables and $(2,0)$-cables.
\end{lemma}

\begin{proof}
Suppose $\psi \in \wdzero \yxonexn$, so it has no delay nodes and its supplier assignment
\[\nicexy{\dmpsi = \yout \amalg \xin \ar[r]^-{\spsi}_-{\cong} & \yin \amalg \xout = \supplypsi}\]
is a bijection, where $\xin = \coprod_{i=1}^n \xin_i$ and $\xout = \coprod_{i=1}^n \xout_i$.  Since $\spsi\left(\yout\right) \subseteq \xout$ by the non-instantaneity requirement \eqref{non-instant}, there is a decomposition 
\[\xout = \xout_+ \amalg \xout_-\]
such that there are bijections
\[\nicexy{\yout \ar[r]^-{\spsi}_-{\cong} & \xout_- = \im\left(\spsi|_{\yout}\right)} \andspace
\nicexy{\xin \ar[r]^-{\spsi}_-{\cong} & \yin \amalg \xout_+}\]
in $\Fins$.  By definition \eqref{psibar-cospan}, $\chizero\psi \in \uwd\yxonexnbar$ is the following cospan.
\[\nicexy@C+1.3cm{& \ybar = \yin \amalg \yout \ar[d]^-{\Id_{\yin} \amalg \spsi|_{\yout}}\\
\coprod\limits_{i=1}^n \xbar_i = \xin \amalg \xout_+ \amalg \xout_- \ar[r]^-{\left(\spsi|_{\xin},\Id_{\xout}\right)} & \supplypsi = \yin \amalg \xout_+ \amalg \xout_-}\]
Observe that:
\begin{itemize}
\item $\xout_+ \subseteq \supplypsi$ consists of only $(2,0)$-cables in $\chizero\psi$;
\item $\yin \amalg \xout_- \subseteq \supplypsi$ consists of only $(1,1)$-cables in $\chizero\psi$.
\end{itemize}
Since there are no other cables, this finishes the proof.
\end{proof}

\begin{lemma}
\label{lem2:wdzero-uwd}
If $\varphi \in \uwd$ has only $(1,1)$-cables and $(2,0)$-cables, then it is in the image of the operad map $\chizero : \wdzero \to \uwd$.
\end{lemma}

\begin{proof}
Suppose
\[\varphi = \Bigl(\nicexy{U = \coprod\limits_{j=1}^N U_j \ar[r]^-{\fphi} & \cphi & V \ar[l]_-{\gphi}}\Bigr) \in \uwd\vuoneun\]
has only $(1,1)$-cables and $(2,0)$-cables, i.e., $\cphi = \cphioneone \amalg \cphitwozero$.  To construct a $\chizero$-preimage of $\varphi$, first note that  there is a decomposition
\[U = \uone \amalg \utwoplus \amalg \utwominus\]
such that the following statements hold.
\begin{itemize}
\item $\uone = \Bigl\{u \in U : \fphi u \in \cphioneone\Bigr\}$, so there are bijections
\begin{equation}
\label{fphi-gphi-bijections}
\nicexy{\uone \ar[r]^-{\fphi}_-{\cong} & \cphioneone & V \ar[l]_-{\gphi}^-{\cong}}.
\end{equation}
\item $\utwoplus = \Bigl\{u_c \in \fphiinv(c) : c \in \cphitwozero\Bigr\}$ with $u_c$ as in \eqref{usubc-fphi-preimage}.
\item For each $c \in \cphitwozero$, there exist unique $u_+ = u_c \in \utwoplus$ and $u_- \in \utwominus$ such that $\fphi\left(u_{\pm}\right) = c$.  The correspondence $u_+ \leftrightarrow u_-$ defines a bijection $\utwoplus \cong \utwominus$.
\end{itemize}

In \eqref{wddot-uwd-psi-supplier} we already defined $\psi \in \wddot\yxonexbign$ with supplier assignment
\[\nicexy{\dmpsi = \yout \amalg \xin = V \amalg U \setminus \bigl\{u_c : c \in \cphi\bigr\} \ar[d]_-{\spsi}^-{=\, (\gphi, \fphi)}\\
\supplypsi = \yin \amalg \xout \cong \cphi = \cphioneone \amalg \cphitwozero}\]
such that $\chi\psi = \varphi$ by Lemma \ref{lem1:chi-image}.  So it is enough to show that $\psi$ is a strict wiring diagram.  Since $\psi \in \wddot$ has no delay nodes, it suffices to show that its supplier assignment $\spsi$ is a bijection.  First note that the map $\gphi : V \to \cphioneone$ \eqref{fphi-gphi-bijections} is a bijection.  

It remains to show that the map
\[\nicexy{U \setminus \bigl\{u_c : c \in \cphi\bigr\} \ar[r]^-{\fphi} & \cphitwozero}\]
is also a bijection.  We have
\[\begin{split}
U \setminus \bigl\{u_c : c \in \cphi\bigr\} 
&= \uone \amalg \utwoplus \amalg \utwominus \setminus \Bigl\{u_c : c \in \cphioneone \amalg \cphitwozero\Bigr\} \\
&= \utwoplus \amalg \utwominus \setminus \Bigl\{u_c : c \in \cphitwozero\Bigr\}\\
&= \utwominus. 
\end{split}\]
Since $\fphi : \utwominus \to \cphitwozero$ is a bijection, the proof is complete.
\end{proof}

The proof of Theorem \ref{chizero-wdzero-uwd} is complete.

\begin{example}
The following generators in $\uwd$ (Def. \ref{def:generating-uwd}) are in the image of the operad map $\chizero : \wdzero \to \uwd$.
\begin{itemize}
\item  the empty cell $\epsilon \in \uwd\emptynothing$ (Def. \ref{def:uwd-gen-empty});
\item a name change $\tau_f \in \uwd\yx$ (Def. \ref{def:uwd-gen-namechange});
\item a $2$-cell $\theta_{X,Y} \in \uwd\xplusyxy$ (Def. \ref{def:uwd-gen-2cell}); 
\item a $1$-loop $\lambda_{X,x} \in \uwd\xminusxx$ (Def. \ref{def:uwd-gen-loop}).
\end{itemize}
In fact, by Examples \ref{wddot-uwd-empty}--\ref{wddot-uwd-loop}, these generators are $\chizero$-images of strict generating wiring diagrams (Def. \ref{def:generating-strict}).  On the other hand,
\begin{itemize}
\item a $1$-output wire $\omega_* \in \uwd\starnothing$ (Def. \ref{def:uwd-gen-1output}) and
\item a split $\sigma^{(X,x_1,x_2)} \in \uwd\xxprime$ (Def. \ref{def:uwd-gen-split})
\end{itemize}
are not in the image of $\chizero$.
\end{example}

\begin{example}
In the following picture, the strict wiring diagram on the left is sent by $\chizero : \wdzero \to \uwd$ to the undirected wiring diagram on the right.
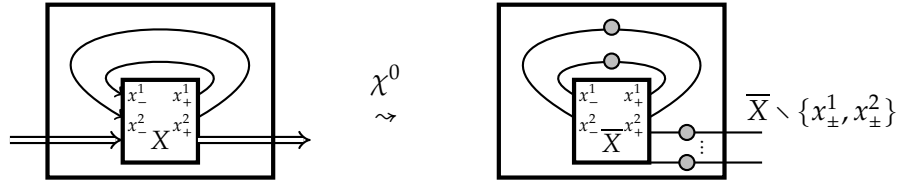
\begin{center}
\begin{tikzpicture}
\draw [ultra thick] (2,.9) rectangle (3,2);
\node at (2.5,1.2) {\small{$X$}};
\node at (2.2,1.8) {\tiny{$x^1_-$}};
\node at (2.8,1.8) {\tiny{$x^1_+$}};
\node at (2.2,1.4) {\tiny{$x^2_-$}};
\node at (2.8,1.4) {\tiny{$x^2_+$}};
\draw [implies] (.5,1.2) to (2,1.2);
\draw [implies] (3,1.2) to (4.5,1.2);
\draw [arrow, looseness=3] (3,1.8) to [out=30, in=150] (2,1.8);
\draw [arrow, looseness=8] (3,1.5) to [out=30, in=150] (2,1.5);
\draw [ultra thick] (1,0.7) rectangle (4,3);
\node at (5.5,1.5) {$\leadsto$};
\node at (5.5,2) {$\chizero$};
\draw [ultra thick] (8,.9) rectangle (9,2);
\node at (8.5,1.2) {\small{$\xbar$}};
\node at (8.2,1.8) {\tiny{$x^1_-$}};
\node at (8.8,1.8) {\tiny{$x^1_+$}};
\node at (8.2,1.4) {\tiny{$x^2_-$}};
\node at (8.8,1.4) {\tiny{$x^2_+$}};
\draw [thick] (9,1.3) to (10.5,1.3);
\draw [thick] (9,.9) to (10.5,.9);
\node at (9.7,1.1) {\tiny{$\vdots$}};
\draw [thick, looseness=3] (9,1.8) to [out=30, in=150] (8,1.8);
\draw [thick, looseness=8] (9,1.5) to [out=30, in=150] (8,1.5);
\draw [ultra thick] (7,0.7) rectangle (10,3);
\smallcable{(8.5,2.25)} \smallcable{(8.5,2.7)}
\smallcable{(9.5,1.3)} \smallcable{(9.5,.9)}
\node at (11.3,1.6) {$\xbar \setminus\{\xsubpm^1, \xsubpm^2\}$};
\end{tikzpicture}
\end{center}
On the right, there are two $(2,0)$-cables, and the other cables are $(1,1)$-cables.
\end{example}

\begin{example}
In the following picture, the strict wiring diagram on the left is sent by $\chizero : \wdzero \to \uwd$ to the undirected wiring diagram on the right.
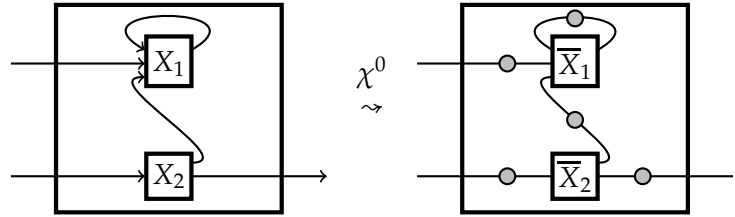
\begin{center}
\begin{tikzpicture}[scale=0.6]
\draw [ultra thick] (1,0.2) rectangle (6,4.8);
\draw [ultra thick] (3,3) rectangle (4,4.1);
\node at (3.5,3.5) {$X_1$};
\draw [ultra thick] ((3,.5) rectangle (4,1.5);
\node at (3.5,1) {$X_2$};
\draw [arrow] (0,3.5) -- (2.97,3.5);
\draw [arrow] (0,1) -- (2.97,1);
\draw [arrow] (4,1) -- (7,1);
\draw [arrow, looseness=1.4] (4,1.3) to [out=0, in=180] (2.97,3.2);
\draw [arrow, looseness=5] (4,3.8) to [out=30, in=150] (2.97,3.8);
\node at (8,2.5) {$\leadsto$};
\node at (8,3.3) {$\chizero$};
\draw [ultra thick] (10,0.2) rectangle (15,4.8);
\draw [ultra thick] (12,3) rectangle (13,4.1);
\node at (12.5,3.5) {$\xbar_1$};
\draw [ultra thick] (12,.5) rectangle (13,1.5);
\node at (12.5,1) {$\xbar_2$};
\draw [thick] (9,3.5) -- (12,3.5);
\draw [thick] (9,1) -- (12,1);
\draw [thick] (13,1) -- (16,1);
\draw [thick, looseness=1.4] (13,1.3) to [out=0, in=180] (12,3.2);
\draw [thick, looseness=5] (13,3.8) to [out=30, in=150] (12,3.8);
\cable{(11,3.5)} \cable{(11,1)} \cable{(12.5,4.5)}
\cable{(12.5,2.25)} \cable{(14,1)}
\end{tikzpicture}
\end{center}
On the right, there are two $(2,0)$-cables and three $(1,1)$-cables.
\end{example}

\section{Summary of Chapter \ref{ch12-maps}}

\begin{enumerate}
\item There is an operad map $\chi : \wddot \to \uwd$ given by forgetting directions whose image consists of precisely the undirected wiring diagrams with no wasted cables and no $(0, \geq 2)$-cables.
\item The restriction of $\chi$ to $\wdzero$ is an operad map $\chizero : \wdzero \to \uwd$.  Its image consists of precisely the undirected wiring diagrams whose cables are either $(1,1)$-cables or $(2,0)$-cables.
\end{enumerate}

\chapter{Map from Wiring Diagrams to Undirected Wiring Diagrams}
\label{ch13-wd-uwd}

This chapter has two main purposes.
\begin{enumerate}
\item We extend the operad map $\chi : \wddot \to \uwd$ in Theorem \ref{wddot-uwd-operad-map}, defined for normal wiring diagrams (i.e., those without delay nodes), to an operad map $\rho : \WD \to \uwd$ that is defined for all wiring diagrams.  See Theorem \ref{wd-uwd-operad-map}.
\item Furthermore, we will show that the operad map $\rho : \WD \to \uwd$ is surjective; see Theorem \ref{thm:wd-uwd-surjective}.  In other words, every undirected wiring diagram is the $\rho$-image of some wiring diagram.  
\end{enumerate}
We remind the reader that the image of the operad map $\chi : \wddot \to \uwd$ was identified in Theorem \ref{thm:chi-image}.  It consists of precisely those undirected wiring diagrams with no wasted cables and no $(0, \geq 2)$-cables.

At first glance, the existence of the operad map $\rho : \WD \to \uwd$ is not obvious because a general wiring diagram has delay nodes, but undirected wiring diagrams have no obvious analogues of delay nodes.  In fact, for this reason Rupel and Spivak (\cite{rupel-spivak} 4.1) expressed doubt that there exists an operad map from $\WD$ to $\uwd$.  Our main results in this chapter, Theorems \ref{wd-uwd-operad-map} and \ref{thm:wd-uwd-surjective}, show that not only  is there an operad map $\rho : \WD \to \uwd$, but also it is surjective.  As we will see in \eqref{psi-delay-nodes}, delay nodes play a critical role in realizing wasted cables and $(0,\geq 2)$-cables in undirected wiring diagrams.

In Section \ref{sec:wd-undirected} we prove that an operad map $\rho : \WD \to \uwd$ exists and that it extends the existing operad map $\chi : \wddot \to \uwd$.  The construction of the operad map $\rho$ is motivated in Example \ref{ex:rho-motivation}, where we discuss how delay nodes should be sent to undirected wiring diagrams.

In Section \ref{sec:wd-uwd-map-example} we provide a series of examples, all containing delay nodes, to further illustrate the operad map $\rho : \WD \to \uwd$.

In Section \ref{sec:wd-undirected-onto} we prove that the operad map $\rho : \WD \to \uwd$ is surjective.  This section ends with Example \ref{ex:rho-preimage}, which provides a detailed illustration of how to lift an undirected wiring diagram back to a wiring diagram.

\section{Wiring Diagrams to Undirected Wiring Diagrams}
\label{sec:wd-undirected}

The purpose of this section is to construct an operad map $\rho : \WD \to \uwd$ that extends the operad map $\chi : \wddot \to \uwd$ in Theorem \ref{wddot-uwd-operad-map}.  Since normal wiring diagrams are wiring diagrams with no delay nodes, to construct the operad map $\rho$, we need to decide how to map the delay nodes to undirected wiring diagrams.

\begin{example}
\label{ex:rho-motivation}
Before we define the desired operad map $\rho : \WD \to \uwd$, let us consider a motivating example that explains what happens to delay nodes.  In the following picture, the wiring diagram $\varphi \in \WD\ynothing$ on the left is sent to the undirected wiring diagram $\varphibar \in \uwd\ybarnothing$ on the right.  
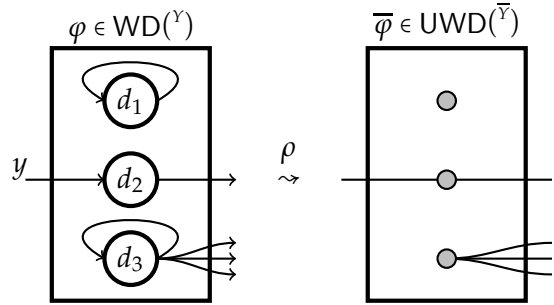
\begin{center}
\begin{tikzpicture}[scale=0.7]
\draw [ultra thick] (2,-.8) rectangle (5,4);
\node at (3.5,4.4) {$\varphi  \in \WD\ynothing$};
\draw [ultra thick] (3.5,3) circle [radius=0.5];
\node at (3.5,3) {$d_1$};
\draw [arrow, looseness=5] (4,3) to [out=30, in=150] (3,3);
\draw [ultra thick] (3.5,1.5) circle [radius=0.5];
\node at (3.5,1.5) {$d_2$};
\draw [arrow] (1.5,1.5) to (3,1.5);
\node at (1.4,1.7) {$y$};
\draw [arrow] (4,1.5) to (5.5,1.5);
\draw [ultra thick] (3.5,0) circle [radius=0.5];
\draw [arrow, looseness=5] (4,0) to [out=30, in=150] (3,0);
\node at (3.5,0) {$d_3$};
\draw [arrow] (4,0) to [out=0, in=180] (5.5,.3);
\draw [arrow] (4,0) to (5.5,0);
\draw [arrow] (4,0) to [out=0, in=180] (5.5,-.3);
\node at (6.5,1.5) {$\leadsto$};
\node at (6.5,2) {$\rho$};
\draw [ultra thick] (8,-.8) rectangle (11,4);
\node at (9.5,4.5) {$\varphibar  \in \uwd\ybarnothing$};
\draw [thick] (7.5,1.5) to (11.5,1.5);
\draw [thick] (9.5,0) to [out=0, in=180] (11.5,.3);
\draw [thick] (9.5,0) to (11.5,0);
\draw [thick] (9.5,0) to [out=0, in=180] (11.5,-.3);
\cable{(9.5,3)} \cable{(9.5,1.5)} \cable{(9.5,0)}
\end{tikzpicture}
\end{center}
In $\varphi$ there are $3$ delay nodes and no input boxes.  In $\varphibar$ there are $3$ cables and no input boxes.  As in the operad map $\chi : \wddot \to \uwd$, every supply wire $\{y,d_1,d_2,d_3\}$ in $\varphi$ yields a cable in $\varphibar$.  However, since a delay node is both a demand wire and a supply wire, we need to identify the cables corresponding to a delay node $d$ and its supply wire $s(d)$.
\begin{itemize}
\item
In $\varphi$ the top delay node $d_1$ supplies only itself, so the identification is trivial.  It yields a wasted cable in $\varphibar$.  
\item
The middle delay node $d_2$ is supplied by the global input $y$, so their cables are identified, yielding a $(0,2)$-cable in $\varphibar$.  
\item
The bottom delay node $d_3$ supplies itself and three global outputs, so the identification is trivial.  Its cable is a $(0,3)$-cable in $\varphibar$.
\end{itemize}
For a general wiring diagram, this identification is defined in \eqref{rhopsi-cospan} below.

Observe that in $\varphibar$, there are a wasted cable, a $(0,2)$-cable, and a $(0,3)$-cable, none of which is possible in the image of $\chi : \wddot \to \uwd$ by Theorem \ref{thm:chi-image}.  In fact, this example suggests that the desired operad map $\rho : \WD \to \uwd$ is surjective because wasted cables and $(0,\geq 2)$-cables are now realizable by carefully chosen delay nodes.  We will prove in Theorem \ref{thm:wd-uwd-surjective} that this is indeed the case.
\end{example}

We now define the operad map $\rho : \WD \to \uwd$ that extends the operad map $\chi : \wddot \to \uwd$ in Theorem \ref{wddot-uwd-operad-map}.  Recall the color map $\chi_0 : \boxs \to \Fins$ in \eqref{wddot-to-uwd-colors} with
\[\chi_0 Y = \ybar = \yin \amalg \yout \in \Fins\]
for each $Y = (\yin, \yout) \in \boxs$.

\begin{definition}
\label{def:wd-uwd-map}
Fix a class $S$.  For each $\yxonexn \in \profboxsboxs$ with $n \geq 0$, define the map
\begin{equation}
\label{wd-to-uwd-entries}
\nicexy{\WD\yxonexn \ar[r]^-{\rho_1} & \uwd\yxonexnbar}
\end{equation}
as follows.  For $\psi \in \WD\yxonexn$, its image
\[\rho_1 \psi = \psibar \in \uwd\yxonexnbar\]
is the cospan
\begin{equation}
\label{rhopsi-cospan}
\nicexy@C+1cm{&& \ybar = \yin \amalg \yout \ar[d]_-{\Id_{\yin} \amalg \spsi|_{\yout}} \ar@<4ex>@/^4pc/[dd]\\
&& \supplypsi = \yin \amalg \xout \amalg \dnpsi \ar[d]_-{\text{quotient}}\\
\xin \amalg \xout \ar[r]^-{\left(\spsi|_{\xin},\Id_{\xout}\right)} \ar@/_1.5pc/[rr]
& \supplypsi \ar[r]^-{\text{quotient}} & C_{\rho\psi} = \dfrac{\supplypsi}{\Bigl(d \,=\, \spsi d \,:\, d \in \dnpsi\Bigr)}}
\end{equation}
in $\Fins$.  Here:
\begin{itemize}
\item
$\xout = \coprod_{i=1}^n \xout_i$ and $\xin = \coprod_{i=1}^n \xin_i$.
\item
The map
\[\nicexy{\dmpsi = \yout \amalg \xin \amalg \dnpsi \ar[r]^-{\spsi} & \yin \amalg \xout \amalg \dnpsi = \supplypsi}\]
is the supplier assignment for $\psi$ (Def. \ref{def:wiring-diagram}).
\item
$\spsi|_{\yout} : \yout \to \xout \amalg \dnpsi$ by the non-instantaneity requirement \eqref{non-instant}.
\item
$C_{\rho\psi}$ in the lower right corner is the quotient set obtained from $\supplypsi$ by identifying $d$ and $\spsi d$ for each delay node $d \in \dnpsi$.
\end{itemize}
\end{definition}

\begin{remark}
\label{rk:rho=chi-in-wddot}
In Def. \ref{def:wd-uwd-map} suppose $\psi \in \wddot\yxonexn$, i.e., $\dnpsi = \varnothing$.  Then
\[C_{\rho\psi} = \supplypsi = \yin \amalg \xout,\]
so $\rho_1\psi$ in \eqref{rhopsi-cospan} is equal to $\chi_1\psi$ in \eqref{psibar-cospan}.  In other words, when applied to normal wiring diagrams, the entry maps $\rho_1$ and $\chi_1$ are the same.  So Def. \ref{def:wd-uwd-map} is indeed an extension of Def. \ref{def:wddot-uwd-map} to all wiring diagrams.
\end{remark}

\begin{theorem}
\label{wd-uwd-operad-map}
The maps $\rho_0 = \chi_0$ \eqref{wddot-to-uwd-colors} and $\rho_1$ \eqref{wd-to-uwd-entries} define a map of operads
\begin{equation}
\label{wd-uwd-rho}
\nicexy{\WD \ar[r]^-{\rho} & \uwd}.
\end{equation}
\end{theorem}

\begin{proof}
This proof is similar to that of Theorem \ref{wddot-uwd-operad-map}.  The difference here is that we now need to take into account the delay nodes.

We will write both $\rho_0$ and $\rho_1$ as $\rho$.  We must check that $\rho$ preserves the operad structure in the sense of Def. \ref{def:operad-map}.  In both $\WD$ \eqref{wd-permutation} and $\uwd$ \eqref{uwd-right-action}, the equivariant structure is given by permuting the labels of the input boxes.  So $\rho$ preserves equivariance in the sense of \eqref{operad-map-eq}.  Likewise, it follows immediately from the definitions of the colored units in $\WD$ \eqref{wd-unit} and $\uwd$ \eqref{uwd-unit} that they are preserved by $\rho$ in the sense of \eqref{operad-map-unit}.

To check that $\rho$ preserves operadic composition in the sense of \eqref{operad-map-comp}, suppose $\varphi \in \WD\yxonexn$ with $n \geq 1$, $1 \leq i \leq n$, and $\psi \in \WD \xiwonewm$ with $m \geq 0$.  We must show that
\begin{equation}
\label{rho-preserves-operad-comp}
\rho\Bigl(\phicompipsi\Bigr) = \bigl(\rho \varphi\bigr) \compi \bigl(\rho \psi\bigr) \in \uwd\yzbar
\end{equation}
in which 
\[\zbar = \overline{\left(\uX \compi \uW\right)} = \Bigl(\xbar_1,\ldots, \xbar_{i-1}, \wbar_1,\ldots, \wbar_m, \xbar_{i+1}, \ldots, \xbar_n \Bigr) \in \proffins\]
as in \eqref{compi-profile}, $\uX = (X_1,\ldots,X_n)$, and $\uW = (W_1,\ldots,W_m)$.

To prove \eqref{rho-preserves-operad-comp}, on the one hand, by Def. \ref{def:compi-wd} $\phicompipsi \in \WD \yxcompiw$ has supplier assignment
\[\nicexy{\dmphicompipsi = \yout \amalg \coprod\limits_{j\not= i} \xin_j \amalg
\win \amalg \dnphi \amalg \dnpsi \ar[d]_-{\sphicompipsi}\\
\supplyphicompipsi = \yin \amalg \coprod\limits_{j\not= i} \xout_j \amalg \wout \amalg \dnphi \amalg \dnpsi}\]
that is given by $\sphi$, $\spsi\sphi$, $\spsi$, $\sphi\spsi$, or $\spsi\sphi\spsi$ according to \eqref{compi-supply1} and \eqref{compi-supply2}.  Here
\[\win = \coprod_{k=1}^m \win_k \andspace \wout = \coprod_{k=1}^m \wout_k \in \Fins.\]
So by \eqref{rhopsi-cospan}  $\rho\left(\phicompipsi\right) \in \uwd\yzbar$ is the cospan
\begin{equation}
\label{rho-phicompipsi-cospan}
\nicexy@C+1.5cm{& \ybar = \yin \amalg \yout \ar[d]^-{\left(\Id_{\yin}, \sphicompipsi|_{\yout}\right)}\\
\zin \amalg \zout \ar[r]^-{\left(\sphicompipsi|_{\zin}, \Id_{\zout}\right)} & C_{\rho(\phicompipsi)} = \dfrac{\supplyphicompipsi = \yin \amalg \zout \amalg \dnphi \amalg \dnpsi}{\Bigl(d\,=\, \sphicompipsi d \,:\, d \in \dnphi \amalg \dnpsi\Bigr)}}
\end{equation}
in $\Fins$.  In \eqref{rho-phicompipsi-cospan}:
\begin{itemize}
\item
$\uZ = \uX \compi \uW \in \boxs$, so
\begin{equation}
\label{zin-zout-xcompiw}
\zin = \coprod\limits_{j\not=i} \xin_j \amalg \win \andspace 
\zout = \coprod\limits_{j\not=i} \xout_j \amalg \wout \in \Fins.
\end{equation} 
\item
The input soldering function is induced by
\begin{itemize}
\item
the identity map on $\zout$;
\item
the supplier assignment $\sphicompipsi : \zin \to \supplyphicompipsi$.
\end{itemize}
\item
The output soldering function is induced by
\begin{itemize}
\item
the identity map on $\yin$;
\item
the supplier assignment $\sphicompipsi : \yout \to \zout \amalg \dnphi \amalg \dnpsi$.
\end{itemize}
\end{itemize}
In \eqref{rho-phicompipsi-cospan} and in what follows, to simplify the notation, we use the same symbol to denote a map and a map induced by it.

On the other hand, by \eqref{rhopsi-cospan}  and Def. \ref{def:compi-uwd}, the $\compi$-composition
\[\bigl(\rho \varphi\bigr) \compi \bigl(\rho \psi\bigr) \in \uwd\yzbar\]
is the cospan
\begin{equation}
\label{rho-phicompipsi-cospan2}
\begin{small}
\nicexy@C+.7cm@R+.5cm{&& \ybar = \yin \amalg \yout \ar[d]_(.4){\left(\Id_{\yin}, \sphi|_{\yout}\right)} \ar@/^2pc/[dd]\\
& \xin \amalg \xout \pushout \ar[d]_(.4){\left(\Id, \spsi\vert_{\xout_i}\right)} \ar[r]^-{\left(\sphi|_{\xin}, \Id_{\xout}\right)} & C_{\rho\varphi} \ar[d]\\
\coprod\limits_{j\not=i} \xbar_j  \amalg \win \amalg \wout \ar[r]^-{\left(\spsi|_{\win}, \Id\right)} \ar@/_2pc/[rr] & \coprod\limits_{j\not=i} \xbar_j  \amalg C_{\rho\psi} \ar[r] & C}
\end{small}
\end{equation}
\bigskip
in $\Fins$, in which the square is defined as a pushout.  In \eqref{rho-phicompipsi-cospan2}:
\begin{itemize}
\item
$C_{\rho\varphi}$ and $C_{\rho\psi}$ are the sets of cables in $\rho\varphi$ and $\rho\psi$, i.e.,
\[C_{\rho\varphi} = \frac{\supplyphi = \yin \amalg \xout \amalg \dnphi}{\Bigl(d \,=\, \sphi d \,:\, d \in \dnphi\Bigr)}
\andspace C_{\rho\psi} = \frac{\supplypsi = \xin_i \amalg \wout \amalg \dnpsi}{\Bigl(d \,=\, \spsi d \,:\, d \in \dnpsi\Bigr)}.\]
\item
In the middle vertical map, the restriction to 
\begin{itemize}
\item
$\xout_i$ is the composition $\nicexy{\xout_i \ar[r]^-{\spsi} & \wout \amalg \dnpsi \subseteq \supplypsi \ar[r] & C_{\rho\psi}}$;
\item
all other coproduct summands is induced by the identity map.
\end{itemize}
\item
In the bottom left horizontal map, the restriction to
\begin{itemize}
\item
$\win$ is the composition $\nicexy{\win \ar[r]^-{\spsi} & \supplypsi \ar[r] & C_{\rho\psi}}$;
\item
all other coproduct summands is induced by the identity map.
\end{itemize}
\item
In the middle right horizontal map, the restriction to
\begin{itemize}
\item $\xin$ is the composition $\nicexy{\xin \ar[r]^-{\sphi} & \supplyphi \ar[r] & C_{\rho\varphi}}$;
\item $\xout$ is induced by the identity map.
\end{itemize}
\end{itemize}

Therefore, to check the condition \eqref{rho-preserves-operad-comp}, it suffices to show that the cospans \eqref{rho-phicompipsi-cospan} and \eqref{rho-phicompipsi-cospan2} are the same.  So it is enough to show that the square
\begin{equation}
\label{rho-operadmap-pushout}
\nicexy@C+1cm@R+.5cm{\xin \amalg \coprod\limits_{j\not=i} \xout_j \amalg \xout_i = \xin \amalg \xout \ar[d]_-{\left(\Id, \spsi\vert_{\xout_i}\right)} \ar[r]^-{\left(\sphi\vert_{\xin}, \Id_{\xout}\right)} 
& C_{\rho\varphi} \ar[d]_-{\left(\Id, \spsi\vert_{\xout_i}\right)}\\
\coprod\limits_{j\not=i} \xin_j \amalg \coprod\limits_{j\not=i} \xout_j \amalg C_{\rho\psi} = \coprod\limits_{j\not=i} \xbar_j  \amalg C_{\rho\psi} \ar[r]^-{h} & C_{\rho(\phicompipsi)}}
\end{equation}
in $\Fins$ is a pushout (Def. \ref{def:pushout}).  In \eqref{rho-operadmap-pushout} the map $h$ is induced by
\begin{itemize}
\item
the composition $\nicexy@C+.5cm{\xin \ar[r]^-{\sphi} & \supplyphi \ar[r] & C_{\rho\varphi} \ar[r]^-{\left(\Id,\spsi\vert_{\xout_i}\right)} & C_{\rho(\phicompipsi)}}$;
\item
the identity map of $\coprod_{j\not=i} \xout_j  \amalg \wout \amalg \dnpsi$.
\end{itemize}
A direct inspection reveals that the right vertical map and the map $h$ in \eqref{rho-operadmap-pushout} are both well-defined.

To see that the square \eqref{rho-operadmap-pushout} is commutative, observe that both compositions when restricted to
\begin{itemize}
\item
$\xin$ is induced by $\Bigl(\Id, \spsi\vert_{\xout_i}\Bigr) \comp \left(\sphi|_{\xin}\right)$;
\item
$\coprod\limits_{j\not=i} \xout_j$ is induced by the identity map;
\item
$\xout_i$ is induced by $\spsi|_{\xout_i}$.
\end{itemize}

It remains to check the condition \eqref{pushout-universal-property} of a pushout.  So suppose given a solid-arrow commutative diagram
\[\nicexy@C+.5cm@R+.5cm{\xin \amalg \coprod\limits_{j\not=i} \xout_j \amalg \xout_i \ar[d]_-{\left(\Id, \spsi\vert_{\xout_i}\right)} \ar[r]^-{\left(\sphi\vert_{\xin}, \Id\right)} & C_{\rho\varphi} \ar[d]_-{\left(\Id, \spsi\vert_{\xout_i}\right)} \ar@/^1pc/[ddr]^-{\beta} & \\
\coprod\limits_{j\not=i} \xin_j \amalg \coprod\limits_{j\not=i} \xout_j \amalg C_{\rho\psi} \ar[r]^-{h} \ar@/_1pc/[drr]^-{\alpha} & C_{\rho(\phicompipsi)} \ar@{.>}[dr]|-{\eta}& \\ & & V}\]
in $\Fins$.  We must show that there exists a unique dotted map $\eta$ that makes the diagram commutative.  Recall that
\[C_{\rho\varphi} = \frac{\supplyphi = \yin \amalg \xout \amalg \dnphi}{\Bigl(d \,=\, \sphi d \,:\, d \in \dnphi\Bigr)}, \quad
 C_{\rho\psi} = \frac{\supplypsi = \xin_i \amalg \wout \amalg \dnpsi}{\Bigl(d \,=\, \spsi d \,:\, d \in \dnpsi\Bigr)},\]
and
\[C_{\rho(\phicompipsi)} = \frac{\supplyphicompipsi = \yin \amalg  \coprod\limits_{j\not=i} \xout_j \amalg \wout \amalg \dnphi \amalg \dnpsi}{\Bigl(d\,=\, \sphicompipsi d \,:\, d \in \dnphi \amalg \dnpsi\Bigr)}.\]
A direct inspection reveals that the only possible candidate for $\eta$ is the map induced by the compositions
\begin{equation}
\label{rho-operad-map-eta}
\begin{split}
&\nicexy{\yin \amalg \coprod\limits_{j\not=i} \xout_j \amalg \dnphi \subseteq \supplyphi \ar[r] & C_{\rho\varphi} \ar[r]^-{\beta} & V};\\
&\nicexy{\wout \amalg \dnpsi \subseteq \supplypsi \ar[r] & C_{\rho\psi} \ar[r]^-{\alpha} &  V}.
\end{split}
\end{equation}
One can check that these definitions indeed yield a well-defined map $\eta$.  So with $\eta$ defined as in \eqref{rho-operad-map-eta}, it remains to check  the equalities
\begin{equation}
\label{etah=alpha-etastuff=beta}
\eta \left(\Id, \spsi\vert_{\xout_i}\right) = \beta \andspace \eta h = \alpha.
\end{equation}
Both of these equalities can be checked by a direct inspection.  This finishes the proof that the square \eqref{rho-operadmap-pushout} is a pushout and, therefore, that $\rho$ preserves operadic composition \eqref{chi-preserves-operad-comp}.
\end{proof}

\begin{example}
By Example \ref{ex:operad-map-algebra} and Theorem \ref{wd-uwd-operad-map}, every $\uwd$-algebra (Def. \ref{def:uwd-algebra-biased}) induces a $\WD$-algebra (Def. \ref{def:wd-algebra}) along the operad map $\rho : \WD \to \uwd$.  For example:
\begin{itemize}
\item The relational algebra of a set $A$ (Def. \ref{def:relational-algebra-biased} with $S = *$) induces a $\WD$-algebra along the operad map $\rho$.
\item The typed relational algebra (Def. \ref{def:typed-relational-algebra} with $S = \set$) also induces a $\WD$-algebra along the operad map $\rho$.
\end{itemize}
\end{example}

\section{Examples of the Operad Map}
\label{sec:wd-uwd-map-example}

In this section we provide examples of the operad map $\rho : \WD \to \uwd$ in Theorem \ref{wd-uwd-operad-map}.

Recall from Remark \ref{rk:rho=chi-in-wddot} that the operad map $\rho : \WD \to \uwd$ extends the operad map $\chi : \wddot \to \uwd$.  Therefore, the next example and Examples \ref{wddot-uwd-empty}--\ref{wddot-uwd-1wasted} give a complete description of the $\rho$-images of all the generating wiring diagrams (Def. \ref{def:generating-wiring-diagrams}).

\begin{example}
For an element $d \in S$, the $1$-delay node $\delta_d \in \WD\dnothing$ (Def. \ref{def:one-dn}) is sent by $\rho$ to the undirected wiring diagram on the right.
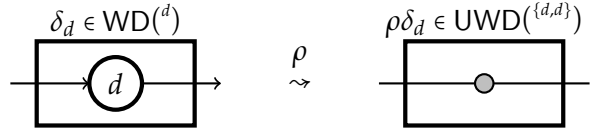
\begin{center}
\begin{tikzpicture}[scale=0.7]
\draw [ultra thick] (0,1.2) rectangle (3,2.8);
\node at (1.5,3.2) {$\delta_d \in \WD\dnothing$};
\draw [ultra thick] (1.5,2) circle [radius=0.5];
\node at (1.5,2) {$d$};
\draw [arrow] (-.5,2) -- (1,2);
\draw [arrow] (2,2) -- (3.5,2);
\node at (5,2) {$\leadsto$};
\node at (5,2.5) {$\rho$};
\draw [ultra thick] (7,1.2) rectangle (10,2.8);
\node at (8.5,3.2) {$\rho\delta_d \in \uwd\ddnothing$};
\draw [thick] (6.5,2) to (10.5,2);
\cable{(8.5,2)}
\end{tikzpicture}
\end{center}
In $\delta_d$ the delay node is supplied by the unique global input, so in $\rho\delta_d$ there is only one cable.  In $\rho\delta_d$ the output box $\{d,d\} \in \Fins$ is the $S$-finite set with two copies of the element $d \in S$.  There are no input boxes in $\rho\delta_d$, and the only cable in it is a $(0,2)$-cable.  
\end{example}

\begin{example}
In the following picture, the wiring diagram $\varphi \in \WD\yx$ with one delay node $d$ is sent by $\rho$ to the undirected wiring diagram on the right.
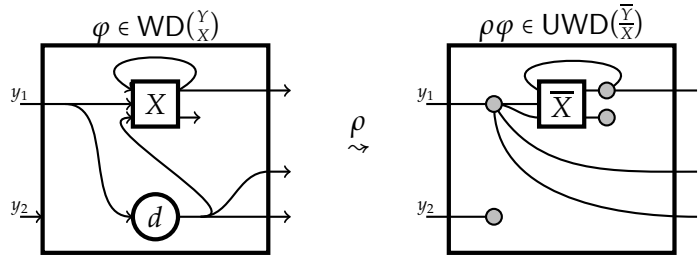
\begin{center}
\begin{tikzpicture}[scale=0.6]
\draw [ultra thick] (1,0.2) rectangle (6,4.8);
\node at (3.5,5.2) {$\varphi \in \WD\yx$};
\draw [ultra thick] (3,3) rectangle (4,4);
\node at (3.5,3.5) {$X$};
\draw [ultra thick] (3.5,1) circle [radius=0.5];
\node at (3.5,1) {$d$};
\draw [thick] (0.5,3.5) -- (1,3.5);
\node at (0.5,3.7) {\tiny{$y_1$}};
\draw [arrow] (0.5,1) -- (1,1);
\node at (0.5,1.2) {\tiny{$y_2$}};
\draw [arrow] (6,3.8) -- (6.5,3.8);
\draw [arrow] (6,2) -- (6.5,2);
\draw [arrow] (6,1) -- (6.5,1);
\draw [arrow] (1,3.5) -- (3,3.5);
\draw [arrow, thick] (1.5,3.5) to [out=0, in=180] (3,1);
\draw [thick] (4,1) -- (6,1);
\draw [thick] (4.5,1) to [out=0, in=180] (6,2);
\draw [arrow, looseness=1.2] (4.5,1) to [out=0, in=180] (3,3.2);
\draw [arrow] (4,3.2) -- (4.5,3.2);
\draw [thick] (4,3.8) -- (6,3.8);
\draw [arrow, looseness=5] (4,3.8) to [out=30, in=150] (3,3.8);
\node at (8,2.5) {$\leadsto$};
\node at (8,3) {$\rho$};;
\draw [ultra thick] (10,0.2) rectangle (15,4.8);
\node at (12.5,5.3) {$\rho\varphi \in \uwd\yxbar$};
\node at (9.5,3.7) {\tiny{$y_1$}};
\node at (9.5,1.2) {\tiny{$y_2$}};
\draw [ultra thick] (12,3) rectangle (13,4);
\node at (12.5,3.5) {$\xbar$};
\draw [thick] (13,3.8) -- (15.5,3.8);
\draw [thick, looseness=3] (13.5,3.8) to [out=30,in=150] (12,3.8);
\draw [thick] (9.5,1) -- (11,1);
\draw [thick] (13,3.2) -- (13.5,3.2);
\draw [thick] (9.5,3.5) -- (12,3.5);
\draw [thick] (11,3.5) to [out=0,in=180] (12,3.2);
\draw [thick] (11,3.5) to [out=-60,in=180] (15.5,2);
\draw [thick] (11,3.5) to [out=-90,in=180] (15.5,1);
\cable{(11,3.5)} \cable{(11,1)}
\cable{(13.5,3.8)} \cable{(13.5,3.2)}
\end{tikzpicture}
\end{center}
Indeed, the set of supply wires in $\varphi$ is 
\[\supplyphi = \yin \amalg \dnphi \amalg \xout = \bigl\{y_1,y_2,d \bigr\} \amalg \xout.\]
Since the delay node $d$ is supplied by the global input $y_1$, by definition \eqref{rhopsi-cospan} the set of cables of $\rho\varphi$ is
\[C_{\rho\varphi} = \frac{\bigl\{y_1,y_2,d \bigr\} \amalg \xout}{\Bigl(d = \sphi d = y_1\Bigr)} = \bigl\{y_1,y_2\bigr\} \amalg \xout.\]
Therefore, in $\rho\varphi$ the cable represented by $y_1$ is a $(2,3)$-cable.  It is soldered to: $y_1$, the input of $X$ supplied by $y_1$ in $\varphi$, and the two global output wires and the input wire of $X$ supplied by $d$ in $\varphi$.  The cable represented by $y_2$ is a $(0,1)$-cable.  The other two cables are a $(1,0)$-cable and a $(2,1)$-cable.
\end{example}

\begin{example}
In the following picture, the wiring diagram $\varphi \in \WD\yx$ with one delay node $d$ is sent by $\rho$ to the undirected wiring diagram on the right.
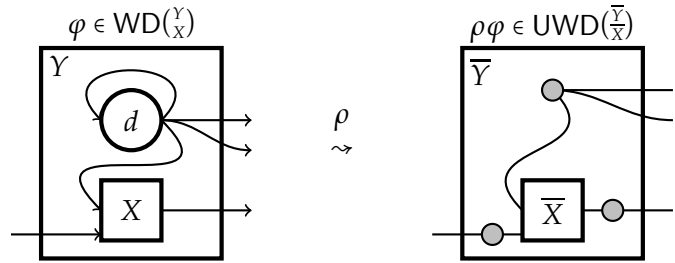
\begin{center}
\begin{tikzpicture}[scale=.8]
\draw [ultra thick] (1.5,2.5) circle [radius=0.5];
\node at (1.5,2.5) {$d$};
\draw [arrow, looseness=4] (2,2.5) to [out=45, in=135] (1,2.5);
\draw [arrow] (2,2.5) to (3.5,2.5);
\draw [arrow] (2,2.5) to [out=0,in=180] (3.5,2);
\draw [arrow, looseness=3] (2,2.5) to [out=-45, in=135] (1,1);
\draw [ultra thick] (1,.5) rectangle (2,1.5);
\draw [arrow] (-.5,.6) to (1,.6);
\draw [arrow] (2,1) to (3.5,1);
\node at (1.5,1) {$X$};
\draw [ultra thick] (0,0.2) rectangle (3,3.7);
\node at (1.5,4.1) {$\varphi \in \WD\yx$};
\node at (0.3,3.4) {$Y$};
\node at (5,2) {$\leadsto$};
\node at (5,2.5) {$\rho$};
\draw [ultra thick] (7,0.2) rectangle (10,3.7);
\node at (8.5,4.1) {$\rho\varphi \in \uwd\yxbar$};
\node at (7.3,3.3) {$\ybar$};
\draw [thick] (8.5,3) to (10.5,3);
\draw [thick] (8.5,3) to [out=0,in=180](10.5,2.5);
\draw [thick, looseness=2] (8.5,3) to [out=-45, in=135]  (8,1);
\cable{(8.5,3)}
\draw [ultra thick] (8,.5) rectangle (9,1.5);
\draw [thick] (6.5,.6) to (8,.6);
\draw [thick] (9,1) to (10.5,1);
\node at (8.5,1) {$\xbar$};
\cable{(7.5,.6)} \cable{(9.5,1)}
\end{tikzpicture}
\end{center}
In $\varphi$ the delay node is supplied by itself, so the set of cables in $\rho\varphi$ is
\[C_{\rho\varphi} = \supplyphi = \yin \amalg \xout \amalg \{d\}.\]
In $\rho\varphi$ the cable corresponding to $d$ is a $(1,2)$-cable.  The other two cables are both $(1,1)$-cables.
\end{example}

\begin{example}
In the following picture, the wiring diagram $\varphi \in \WD\yx$ with one delay node $d$ is sent by $\rho$ to the undirected wiring diagram on the right.
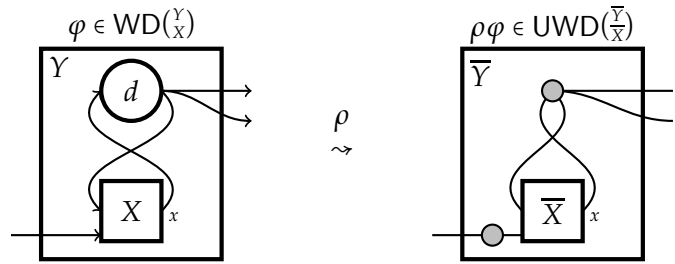
\begin{center}
\begin{tikzpicture}[scale=.8]
\draw [ultra thick] (1.5,3) circle [radius=0.5];
\node at (1.5,3) {$d$};
\draw [arrow] (2,3) to (3.5,3);
\draw [arrow] (2,3) to [out=0,in=180] (3.5,2.5);
\draw [arrow, looseness=1.5] (2,3) to [out=-45, in=135] (1,1);
\draw [ultra thick] (1,.5) rectangle (2,1.5);
\node at (2.2,.9) {\tiny{$x$}};
\draw [arrow] (-.5,.6) to (1,.6);
\draw [arrow, looseness=1.5] (2,1) to [out=45, in=225] (1,3);
\node at (1.5,1) {$X$};
\draw [ultra thick] (0,0.2) rectangle (3,3.7);
\node at (1.5,4.1) {$\varphi \in \WD\yx$};
\node at (0.3,3.4) {$Y$};
\node at (5,2) {$\leadsto$};
\node at (5,2.5) {$\rho$};
\draw [ultra thick] (7,0.2) rectangle (10,3.7);
\node at (8.5,4.1) {$\rho\varphi \in \uwd\yxbar$};
\node at (7.3,3.3) {$\ybar$};
\draw [thick] (8.5,3) to (10.5,3);
\draw [thick] (8.5,3) to [out=0,in=180] (10.5,2.5);
\draw [thick, looseness=1.5] (8.5,3) to [out=-45, in=135] (8,1);
\draw [thick, looseness=1.5] (8.5,3) to [out=225, in=45]  (9,1);
\cable{(8.5,3)}
\draw [ultra thick] (8,.5) rectangle (9,1.5);
\node at (9.2,.9) {\tiny{$x$}};
\draw [thick] (6.5,.6) to (8,.6);
\node at (8.5,1) {$\xbar$};
\cable{(7.5,.6)}
\end{tikzpicture}
\end{center}
In $\varphi$ the delay node $d$ is supplied by the unique output wire $x$ of $X$, so their corresponding cables are identified in $\rho\varphi$.  This cable is a $(2,2)$-cable.  The other cable is a $(1,1)$-cable.
\end{example}

\section{Surjectivity of the Operad Map}
\label{sec:wd-undirected-onto}

The reader is reminded of Notation \ref{notation:cable-subsets} regarding subsets of cables.  The purpose of this section is to show that the operad map $\rho : \WD \to \uwd$ in Theorem \ref{wd-uwd-operad-map} is surjective.  Our strategy is similar to the proof of Lemma \ref{lem1:chi-image}, except that here the input soldering function may not be surjective.  Cables not in the image of the input soldering function are $(0,\geq 0)$-cables.  Wasted cables (i.e., $(0,0)$-cables) and $(0, \geq 2)$-cables are realized using delay nodes, similar to the delay nodes $d_1$ and $d_3$ in Example \ref{ex:rho-motivation}.  Moreover, $(0,1)$-cables are realized using external wasted wires, similar to $y$ in the picture \eqref{chi-of-1wasted-wire}.

Given an undirected wiring diagram, we now define a wiring diagram that will be shown to be a $\rho$-preimage.  Below we will use the map $\chi_0 = \rho_0 : \boxs \to \Fins$ \eqref{wddot-to-uwd-colors}, usually denoted by $\chi_0 Y = \ybar = \yin \amalg \yout$.  The following definition is the general version of Def. \ref{def:chi-image-input-surjective} in the sense that now we do not require the input soldering function to be surjective.  A detailed example of the following definition will be given in Example \ref{ex:rho-preimage}.

\begin{definition}
\label{def:rho-preimage-phi}
Suppose
\begin{equation}
\label{general-uwd-varphi-uv}
\varphi = \Bigl(\nicexy{U = \coprod\limits_{j=1}^N U_i \ar[r]^-{\fphi} & \cphi & V \ar[l]_-{\gphi}}\Bigr) \in \uwd\vuoneun
\end{equation}
for some $N \geq 0$.  We will use the equalities
\begin{equation}
\label{rho-preimage-equalities}
\begin{split}
\cphiatleastoneatleastzero &= \im(\fphi)\\
\cphi &= \cphizeroone \amalg \cphiatleastoneatleastzero \amalg \cphizerozero \amalg \cphizeroatleasttwo\\
V &= \gphiinv \cphizeroone \amalg \gphiinv \cphiatleastoneatleastzero \amalg \gphiinv \cphizeroatleasttwo\\
\end{split}
\end{equation}
below.
\begin{enumerate}
\item For each $c \in \cphiatleastoneatleastzero$, pick a preimage
\[ u_c \in \fphiinv(c) \subseteq U.\]
We will use the bijection
\begin{equation}
\label{fphi-on-usubc}
\nicexy{\{u_c\} = \Bigl\{u_c : c \in \cphiatleastoneatleastzero\Bigr\} \ar[r]^-{\fphi}_-{\cong} & \cphiatleastoneatleastzero}
\end{equation}
and its inverse below.
\item
For each $1 \leq j \leq N$ define a box $X_j = (\xin_j, \xout_j) \in \boxs$ as
\begin{equation}
\label{uj-is-xjbar}
\xout_j = \Bigl\{u_c \in U_j : c \in \cphiatleastoneatleastzero\Bigr\} \subseteq U_j \andspace 
\xin_j = U_j \setminus \xout_j.
\end{equation}
Note that we have $\xbar_j = \xin_j \amalg \xout_j = U_j$;
\[\begin{split}
\xout &= \coprod_{j=1}^N \xout_j = \Bigl\{u_c : c \in \cphiatleastoneatleastzero\Bigr\} \cong \cphiatleastoneatleastzero;\\
\xin &= \coprod_{j=1}^N \xin_j = U \setminus \Bigl\{u_c : c \in \cphiatleastoneatleastzero\Bigr\}.
\end{split}\]
\item Define a box $Y = (\yin, \yout) \in \boxs$ as
\begin{equation}
\label{v-is-ybar}
\begin{split}
\yin &= \gphiinv \cphizeroone;\\
\yout &= \gphiinv \cphiatleastoneatleastzero \amalg \gphiinv \cphizeroatleasttwo.
\end{split}
\end{equation}
Note that $\ybar = \yin \amalg \yout = V$.
\item Define $\psi \in \WD\yxonexbign$ with delay nodes 
\begin{equation}
\label{psi-delay-nodes}
\dnpsi = \cphizerozero \amalg \cphizeroatleasttwo.
\end{equation}
Recall the sets 
\[\dmpsi = \yout \amalg \xin \amalg \dnpsi \andspace \supplypsi = \yin \amalg \xout \amalg \dnpsi\]
of demand wires and of supply wires.  The supplier assignment for $\psi$
\begin{equation}
\label{wd-uwd-psi-supplier}
\nicexy{\dmpsi = \gphiinv\cphiatleastoneatleastzero \amalg \gphiinv \cphizeroatleasttwo \amalg \bigl[U \setminus \{u_c\}\bigr] \amalg \cphizerozero \amalg \cphizeroatleasttwo 
\ar[d]_-{\spsi}\\
\supplypsi = \gphiinv \cphizeroone \amalg \{u_c\} \amalg \cphizerozero \amalg \cphizeroatleasttwo}
\end{equation}
is defined by the restrictions:
\[\begin{split}
& \nicexy{\gphiinv \cphiatleastoneatleastzero \amalg \bigl[U \setminus \{u_c\}\bigr] \ar[r]^-{(\gphi, \fphi)} & \cphiatleastoneatleastzero \ar[r]^-{\fphiinv}_-{\cong} & \{u_c\}};\\
& \nicexy{\gphiinv \cphizeroatleasttwo \ar[r]^-{\gphi} & \cphizeroatleasttwo};\\
& \nicexy{\cphizerozero \amalg \cphizeroatleasttwo \ar[r]^-{=} & \cphizerozero \amalg \cphizeroatleasttwo}.
\end{split}\]
Here $\fphiinv$ is the inverse of the bijection \eqref{fphi-on-usubc}.
\end{enumerate}
\end{definition}

\begin{remark}
Consider Def. \ref{def:rho-preimage-phi}.
\begin{enumerate}
\item If $\cphizeroatleastzero = \varnothing$ (i.e., $\fphi$ is surjective), then Def. \ref{def:rho-preimage-phi} reduces to Def. \ref{def:chi-image-input-surjective}.
\item The definition \eqref{uj-is-xjbar} of each box $X_j$ means that:
\begin{itemize}
\item For each cable $c \in \cphiatleastoneatleastzero = \im(\fphi)$, one wire in $U = \coprod_{j=1}^N U_j$, namely $u_c$, soldered to $c$ is made into an output wire in $\psi$.
\item All other wires in $U$ are made into input wires in $\psi$.
\end{itemize}
\item The definition \eqref{v-is-ybar} of the box $Y$ means that:
\begin{itemize}
\item
Elements in $V$ that are soldered to $(0,1)$-cables in $\varphi$ are made into global input wires in $\psi$.
\item
All other elements in $V$ are made into global output wires in $\psi$.
\end{itemize}
\item The supplier assignment $\spsi$ in \eqref{wd-uwd-psi-supplier} satisfies the non-instantaneity requirement \eqref{non-instant} because $\yin = \gphiinv \cphizeroone$ is disjoint from the image of $\spsi$, which is $\{u_c\} \amalg \cphizerozero \amalg \cphizeroatleasttwo$.  In fact, $\yin = \gphiinv \cphizeroone$ is exactly the set of external wasted wires in $\psi$.
\item Each delay node \eqref{psi-delay-nodes} is supplied by itself, i.e., $d = \spsi d$ for each $d \in \dnpsi$.
\end{enumerate}
\end{remark}

\begin{theorem}
\label{thm:wd-uwd-surjective}
Consider the operad map $\rho : \WD \to \uwd$ in Theorem \ref{wd-uwd-operad-map}.
\begin{enumerate}
\item
$\rho$ is surjective on colors.
\item
In the context of Def. \ref{def:rho-preimage-phi}, we have that $\rho\psi = \varphi$.  In particular, the operad map $\rho$ is surjective on entries as well.
\end{enumerate}
\end{theorem}

\begin{proof}
The color map of $\rho$ is $\rho_0 = \chi_0 : \boxs \to \Fins$ \eqref{wddot-to-uwd-colors}, which is surjective by Theorem \ref{thm:chi-image}(1).

For the second assertion, the undirected wiring diagram
\[\rho\psi \in \uwd\yxonexnbar = \uwd\vuoneun\]
is by definition the cospan in \eqref{rhopsi-cospan}.  Since every delay node in $\psi$ \eqref{psi-delay-nodes} is supplied by itself, the set of cables in $\rho\psi$ is $C_{\rho\psi} = \supplypsi$, the set of supply wires in $\psi$.  Furthermore, there is a bijection
\begin{equation}
\label{crhopsi-cphi}
\nicexy@R+.5cm{C_{\rho\psi} = \supplypsi = \gphiinv \cphizeroone \amalg \{u_c\} \amalg \cphizerozero \amalg \cphizeroatleasttwo \ar[d]_-{\gphi \amalg \fphi \amalg \Id}^-{\cong}\\
\cphi = \cphizeroone \amalg \cphiatleastoneatleastzero \amalg \cphizerozero \amalg \cphizeroatleasttwo}
\end{equation}
in which $\fphi$ is the bijection \eqref{fphi-on-usubc}.

By the definition of $\spsi$ \eqref{wd-uwd-psi-supplier}, $\rho\psi \in \uwd\vuoneun$ is the cospan \eqref{rhopsi-cospan}
\[\nicexy@R+.5cm@C+1cm{
& V = \ybar = \overbrace{\gphiinv \cphizeroone}^{\yin} \amalg \overbrace{\gphiinv \cphiatleastoneatleastzero \amalg \gphiinv \cphizeroatleasttwo}^{\yout} \ar[d]_-{\Id_{\yin} \amalg \spsi|_{\yout}}^-{=\, \Id \amalg \fphiinv\gphi \amalg \gphi}\\
U = \underbrace{\bigl[U \setminus \{u_c\}\bigr]}_{\xin} \amalg \underbrace{\{u_c\}}_{\xout} \ar[r]^-{\left(\spsi|_{\xin}, \Id_{\xout}\right)}_-{=\, \left(\fphiinv \fphi, \Id\right)} & 
C_{\rho\psi} = \gphiinv \cphizeroone \amalg \{u_c\} \amalg \cphizerozero \amalg \cphizeroatleasttwo}\]
in which $\fphiinv$ is the inverse of the bijection \eqref{fphi-on-usubc}.  Combining this cospan for $\rho\psi$ with the bijection \eqref{crhopsi-cphi}, there is a commutative diagram
\[\nicexy@R+.5cm@C+.8cm{
& V = \gphiinv \cphizeroone \amalg \gphiinv \cphiatleastoneatleastzero \amalg \gphiinv \cphizeroatleasttwo \ar[d]_-{\Id \amalg \fphiinv\gphi \amalg \gphi} \ar@/^2pc/@<2.7cm>[dd]^-{\gphi}\\
& C_{\rho\psi} = \gphiinv \cphizeroone \amalg \{u_c\} \amalg \cphizerozero \amalg \cphizeroatleasttwo \ar[d]_-{\gphi \amalg \fphi \amalg \Id}^-{\cong}\\
U = \bigl[U \setminus \{u_c\}\bigr] \amalg \{u_c\} \ar@(u,ul)[ur]^-{\left(\fphiinv \fphi, \Id\right)} \ar[r]^-{\fphi} 
& \cphi = \cphizeroone \amalg \cphiatleastoneatleastzero \amalg \cphizerozero \amalg \cphizeroatleasttwo}\]
in $\Fins$.  In this diagram, the outer cospan is $\varphi$ \eqref{general-uwd-varphi-uv}.  Therefore, by Def. \ref{def:uwd} we have proved $\varphi = \rho\psi$.
\end{proof}

\begin{example}
\label{ex:rho-preimage}
This is an illustration of Def. \ref{def:rho-preimage-phi} and Theorem \ref{thm:wd-uwd-surjective}.  Consider the undirected wiring diagram $\varphi \in \uwd\vuoneutwo$ in \eqref{uwd-first-picture}, depicted as
\begin{center}
\begin{tikzpicture}[scale=1]
\draw [ultra thick] (-1,0) rectangle (6.5,3);
\node at (-1.2,1) {\tiny{$v_1$}};
\node at (-1.2,2) {\tiny{$v_2$}};
\node at (.3,3.2) {\tiny{$v_3$}};
\node at (2.5,3.2) {\tiny{$v_4$}};
\node at (3.5,3.2) {\tiny{$v_5$}};
\node at (6.7,1.7) {\tiny{$v_6$}};
\node at (6,.3) {$V$};
\draw [ultra thick] (.3,1) rectangle (2.5,2);
\node at (.5,1.85) {\tiny{$u_1$}};
\node at (1.3,1.85) {\tiny{$u_2$}};
\node at (2.3,1.85) {\tiny{$u_3$}};
\node at (2.3,1.15) {\tiny{$u_4$}};
\node at (1.3,1.15) {\tiny{$u_5$}};
\node at (.5,1.15) {\tiny{$u_6$}};
\draw [ultra thick] (4,1) rectangle (5,2);
\node at (4.5,1.2) {\tiny{$u^2$}};
\node at (4.5,1.8) {\tiny{$u^1$}};
\cable{(-.25,1.5)} \node at (-.25,1.5) {\tiny{$c_1$}};
\cable{(.5,2.5)} \node at (.5,2.5) {\tiny{$c_2$}};
\cable{(3,2.5)} \node at (3,2.5) {\tiny{$c_3$}};
\cable{(6,2.5)} \node at (6,2.5) {\tiny{$c_4$}};
\cable{(6,1.5)} \node at (6,1.5) {\tiny{$c_5$}};
\cable{(3,.5)} \node at (3,.5) {\tiny{$c_6$}};
\cable{(.5,.5)} \node at (.5,.5) {\tiny{$c_7$}};
\draw [thick] (-.425,1.5) -- (-1.7,1);
\draw [thick] (-.425,1.5) -- (-1.7,2);
\draw [thick] (.5,2) -- (.5,2.325);
\draw [thick] (.5,2.675) -- (.5,3.5);
\draw [thick] (1.3,2) to [out=90, in=180] (2.825,2.5);
\draw [thick] (2.3,2) to [out=90, in=180] (2.825,2.5);
\draw [thick] (4.5,2)  to [out=90, in=0] (3.175,2.5);
\draw [thick] (3,2.675) -- (2.7,3.5);
\draw [thick] (3,2.675) -- (3.3,3.5);
\draw [thick] (6.175,1.5) -- (7,1.5);
\draw [thick] (1.3,1) to [out=-90,in=180] (2.825,.5);
\draw [thick] (2.3,1) to [out=-90,in=180] (2.825,.5);
\draw [thick] (4.5,1) to [out=-90,in=0]  (3.175,.5);
\draw [thick] (.5,1) -- (.5,.675);
\end{tikzpicture}
\end{center}
with $V = \{v_1,\ldots,v_6\}$, $U_1 = \{u_1,\ldots,u_6\}$,  $U_2 = \{u^1,u^2\}$, and $\cphi = \{c_1,\ldots,c_7\}$.

Following Def. \ref{def:rho-preimage-phi} first note that we have the subsets
\[\begin{split}
\cphizerozero &= \{c_4\}, \quad \cphizeroone = \{c_5\}, \quad
\cphizeroatleasttwo = \{c_1\}, \quad\text{and}\\
\cphiatleastoneatleastzero &= \im(\fphi) = \{c_2,c_3,c_6,c_7\}.
\end{split}\]
Next, for each cable $c \in \cphiatleastoneatleastzero$, we are supposed to choose an $\fphi$-preimage $u_c \in U = U_1 \amalg U_2$.  We may choose, for example,
\[u_{c_2} = u_1,\quad u_{c_3} = u_2, \quad u_{c_6} = u_5, \andspace u_{c_7} = u_6, \]
all in $U_1$.  With such choices, the boxes $X_1$ and $X_2 \in \boxs$ \eqref{uj-is-xjbar} are
\[\begin{split}
X_1 &= (\xin_1,\xout_1) = \Bigl(\{u_3,u_4\}, \{u_1,u_2,u_5,u_6\}\Bigr);\\
X_2 &= (\xin_2,\xout_2) = \Bigl(\{u^1,u^2\}, \varnothing\Bigr).
\end{split}\]
The box $Y \in \boxs$ \eqref{v-is-ybar} is
\[\begin{split}
(\yin,\yout) &= \Bigl(\gphiinv \cphizeroone, \gphiinv \cphiatleastoneatleastzero \amalg \gphiinv \cphizeroatleasttwo\Bigr)\\
&= \Bigl(\{v_6\}, \{v_1,v_2,v_3,v_4,v_5\}\Bigr).
\end{split}\]

The set of delay nodes of $\psi \in \WD\yxonextwo$ \eqref{psi-delay-nodes} is
\[\dnpsi = \cphizerozero \amalg \cphizeroatleasttwo = \{c_1, c_4\}.\]
The supplier assignment for $\psi$ \eqref{wd-uwd-psi-supplier} is the function
\[\begin{small}
\nicexy{\dmpsi = \yout \amalg \xin \amalg \dnpsi = \{v_1,v_2,v_3,v_4,v_5\} \amalg \{u_3,u_4, u^1, u^2\} \amalg \{c_1, c_4\} \ar[d]_-{\spsi}\\
\supplypsi = \yin \amalg \xout \amalg \dnpsi = \{v_6\} \amalg  \{u_1,u_2,u_5,u_6\} \amalg \{c_1, c_4\}}\end{small}\]
given by
\[\begin{split}
\spsi(v_1) &= \spsi(v_2) = c_1, \quad \spsi(v_3) = u_1, \quad \spsi(v_4) = \spsi(v_5) = u_2,\\
\spsi(u_3) &= \spsi(u^1) = u_2, \quad \spsi(u_4) = \spsi(u^2) = u_5,\\
\spsi(c_1) &= c_1, \andspace \spsi(c_4) = c_4.
\end{split}\]
The supply wire $v_6$ is an external wasted wire, and $u_6$ is an internal wasted wire in $\psi$.  We may draw $\psi \in \WD\yxonextwo$ as follows.
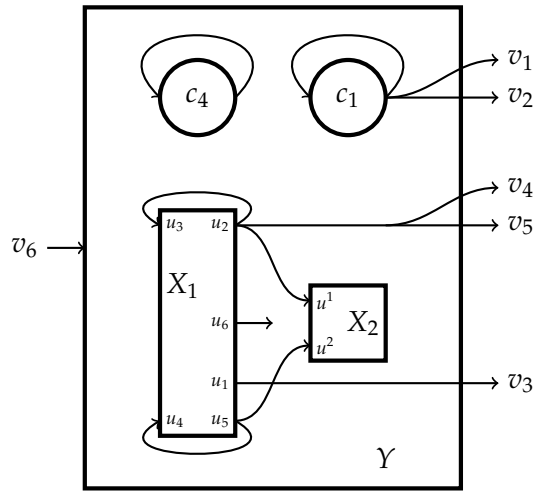
\begin{center}
\begin{tikzpicture}[scale=1]
\draw [ultra thick] (0,0.3) rectangle (5,6.7);
\node at (4,.7) {$Y$};
\draw [arrow, thick] (-.5,3.5) to (0,3.5);
\node at (-.8,3.5) {$v_6$}; 
\draw [ultra thick] (1.5,5.5) circle [radius=.5];
\node at (1.5,5.5) {$c_4$};
\draw [arrow, looseness=4] (2,5.5) to [out=45, in=135] (1,5.5);
\draw [ultra thick] (3.5,5.5) circle [radius=.5];
\node at (3.5,5.5) {$c_1$};
\draw [arrow, looseness=4] (4,5.5) to [out=45, in=135] (3,5.5);
\draw [arrow, thick] (4,5.5) to (5.5,5.5);
\node at (5.8,5.5) {$v_2$};
\draw [arrow, thick] (4,5.5) to [out=0,in=180] (5.5,6);
\node at (5.8,6) {$v_1$};
\draw [ultra thick] (1,1) rectangle (2,4);
\node at (1.3,3) {$X_1$};
\node at (1.8,3.8) {\tiny{$u_2$}};
\node at (1.2,3.8) {\tiny{$u_3$}};
\draw [arrow, looseness=3] (2,3.8) to [out=30, in=150] (1,3.8);
\draw [arrow, thick] (2,3.8) to [out=0,in=180] (3,2.8);
\node at (3.2,2.8) {\tiny{$u^1$}};
\draw [arrow, thick] (2,3.8) to (5.5,3.8);
\node at (5.8,3.8) {$v_5$};
\draw [arrow, thick] (4,3.8) to [out=0,in=180] (5.5,4.3);
\node at (5.8,4.3) {$v_4$};
\draw [arrow, thick] (2,2.5) to (2.5,2.5);
\node at (1.8,2.5) {\tiny{$u_6$}};
\node at (1.8,1.7) {\tiny{$u_1$}};
\draw [arrow, thick] (2,1.7) to (5.5,1.7);
\node at (5.8,1.7) {$v_3$};
\node at (1.8,1.2) {\tiny{$u_5$}};
\node at (1.2,1.2) {\tiny{$u_4$}};
\draw [arrow, looseness=3] (2,1.2) to [out=-30, in=210] (1,1.2);
\draw [arrow, thick] (2,1.2) to [out=0,in=180] (3,2.2);
\node at (3.2,2.2) {\tiny{$u^2$}};
\draw [ultra thick] (3,2) rectangle (4,3);
\node at (3.7,2.5) {$X_2$};
\end{tikzpicture}
\end{center}
By Theorem \ref{thm:wd-uwd-surjective}(2) or a direct inspection, the map $\rho : \WD \to \uwd$ sends $\psi$ to $\varphi$.

Note that $\psi \in \WD\yxonextwo$ is certainly not the only $\rho$-preimage of $\varphi \in \uwd\vuoneutwo$.  For example, the wiring diagram
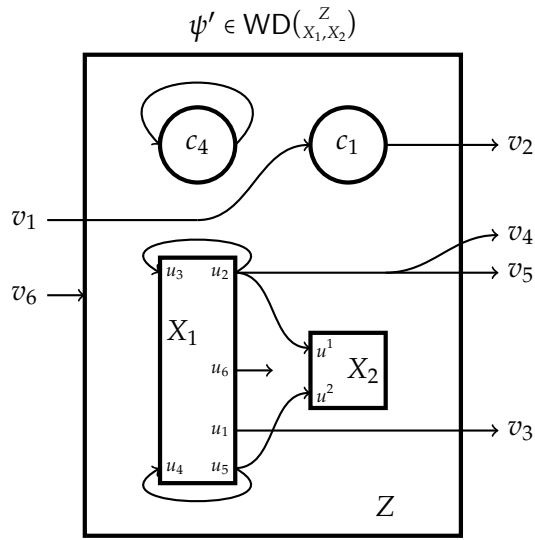
\begin{center}
\begin{tikzpicture}[scale=1]
\draw [ultra thick] (0,0.3) rectangle (5,6.7);
\node at (4,.7) {$Z$};
\node at (2.5, 7.1) {$\psi' \in \WD\zxonextwo$};
\draw [arrow, thick] (-.5,3.5) to (0,3.5);
\node at (-.8,3.5) {$v_6$}; 
\draw [ultra thick] (1.5,5.5) circle [radius=.5];
\node at (1.5,5.5) {$c_4$};
\draw [arrow, looseness=4] (2,5.5) to [out=45, in=135] (1,5.5);
\draw [ultra thick] (3.5,5.5) circle [radius=.5];
\node at (3.5,5.5) {$c_1$};
\node at (-.8,4.5) {$v_1$};
\draw [thick] (-.5,4.5) to (1.5,4.5);
\draw [arrow, thick] (1.5,4.5) to [out=0, in=180] (3,5.5);
\node at (5.8,5.5) {$v_2$};
\draw [arrow, thick] (4,5.5) to (5.5,5.5);
\draw [ultra thick] (1,1) rectangle (2,4);
\node at (1.3,3) {$X_1$};
\node at (1.8,3.8) {\tiny{$u_2$}};
\node at (1.2,3.8) {\tiny{$u_3$}};
\draw [arrow, looseness=3] (2,3.8) to [out=30, in=150] (1,3.8);
\draw [arrow, thick] (2,3.8) to [out=0,in=180] (3,2.8);
\node at (3.2,2.8) {\tiny{$u^1$}};
\draw [arrow, thick] (2,3.8) to (5.5,3.8);
\node at (5.8,3.8) {$v_5$};
\draw [arrow, thick] (4,3.8) to [out=0,in=180] (5.5,4.3);
\node at (5.8,4.3) {$v_4$};
\draw [arrow, thick] (2,2.5) to (2.5,2.5);
\node at (1.8,2.5) {\tiny{$u_6$}};
\node at (1.8,1.7) {\tiny{$u_1$}};
\draw [arrow, thick] (2,1.7) to (5.5,1.7);
\node at (5.8,1.7) {$v_3$};
\node at (1.8,1.2) {\tiny{$u_5$}};
\node at (1.2,1.2) {\tiny{$u_4$}};
\draw [arrow, looseness=3] (2,1.2) to [out=-30, in=210] (1,1.2);
\draw [arrow, thick] (2,1.2) to [out=0,in=180] (3,2.2);
\node at (3.2,2.2) {\tiny{$u^2$}};
\draw [ultra thick] (3,2) rectangle (4,3);
\node at (3.7,2.5) {$X_2$};
\end{tikzpicture}
\end{center}
also satisfies $\rho\psi' = \varphi$.  Here the output box is
\[Z = (\zin, \zout) = \Bigl(\{v_1,v_6\}, \{v_2,v_3,v_4,v_5\}\Bigr) \in \boxs,\]
and $s_{\psi'}(c_1) = v_1$.  Everything else is the same as in $\psi$.
\end{example}

\section{Summary of Chapter \ref{ch13-wd-uwd}}

\begin{enumerate}
\item The operad map $\chi : \wddot \to \uwd$ extends to an operad map $\rho : \WD \to \uwd$ that sends each delay node to a cable.
\item The operad map $\rho$ is surjective.
\end{enumerate}

\chapter{Problems}
\label{ch-problems}

This final chapter contains some problems from the earlier chapters about operads and (undirected) wiring diagrams.

\begin{problem}
Give a detailed proof of Prop. \ref{prop:operad-def-equiv}, which states that the two definitions of a colored operad--one in terms of May's $\gamma$ (Def. \ref{def:colored-operad}) and the other in terms of the $\compi$-compositions (Def. \ref{def:pseudo-operad})--are equivalent.  The equivalence of the unity axioms is rather easy to check.  However, to prove the equivalence of the associativity axioms and the equivariance axioms in the two definitions, a fair amount of bookkeeping and notations are needed.
\end{problem}

\begin{problem}
For the collection $\WD$ of $S$-wiring diagrams (Def. \ref{wd-equivalence}), write down its structure map $\gamma$ \eqref{operadic-composition} and prove that:
\begin{enumerate}
\item
$\WD$ is a $\boxs$-colored operad in the sense of Def. \ref{def:colored-operad}.
\item
This structure map $\gamma$ corresponds to the $\compi$-compositions in Def. \ref{def:compi-wd} in the sense of Prop. \ref{prop:operad-def-equiv}.
\end{enumerate}
In \cite{rupel-spivak} the operad $\WD$ was in fact defined in terms of the structure map $\gamma$.
\end{problem}

\begin{problem}
Consider the $28$ elementary relations in Section \ref{sec:elementary-relations}.
\begin{enumerate}
\item
Give a detailed proof for each elementary relation.  These proofs are similar to those for Lemma \ref{compi-horizontal-associative}, Lemma \ref{compi-vertical-associative}, and Prop. \ref{prop:internal-wasted-wire}. 
\item
For each elementary relation, draw a picture that depicts the operadic compositions, similar to those just before Prop. \ref{prop:move:b1} and Prop. \ref{prop:move:b3}, if one was not given. 
\end{enumerate}
\end{problem}

\begin{problem}
Write down the proof for Lemma \ref{lemma:loop-element}.
\end{problem}

\begin{problem}
In Example \ref{ex:factoring-pi}:
\begin{enumerate}
\item
Write down precisely the wiring diagrams $\pi$, $\pione$, and $\pitwo$, including their supplier assignments.  
\item
Check carefully that there is indeed a decomposition $\pi = \pione \comp \pitwo$.
\end{enumerate}
\end{problem}

\begin{problem}
In Example \ref{ex:factor-pitwo}:
\begin{enumerate}
\item
Write down precisely the wiring diagrams $\beta_1$, $\beta_2$, and $\beta_3$, including their supplier assignments.  
\item
Check carefully that there is indeed a decomposition $\pitwo = \beta_1 \comp \beta_2 \comp \beta_3$.
\end{enumerate}
\end{problem}

\begin{problem}
In Remark \ref{rk:stratified-disjoint} it was stated that stratified simplices of type (1) and of type (2) are mutually exclusive.  Write down a detailed proof for this claim.
\end{problem}

\begin{problem}
For the wiring diagram in \eqref{wd-first-example}, without using Theorem \ref{stratified-presentation-exists}, prove directly that it has a stratified presentation.
\end{problem}

\begin{problem}
Check carefully the proof of Theorem \ref{thm:without-dn-coherence}, which is the finite presentation theorem for the operad $\wddot$ of normal wiring diagrams.
\end{problem}

\begin{problem}
Give a direct proof of Theorem \ref{thm:strict-wd-coherence}--the finite presentation theorem for the operad $\wdzero$ of strict wiring diagrams--without referencing the proofs of Theorem \ref{stratified-presentation-exists},  Lemma \ref{lemma:simplex-to-stratified}, Lemma \ref{lemma:stratified-type1}, and Lemma \ref{lemma:stratified-type2}.
\end{problem}

\begin{problem}
Give a detailed proof that  Def. \ref{def1:operad-algebra} and Def. \ref{def2:operad-algebra} of an operad algebra are indeed equivalent.  The reader may consult \cite{yau-operad} (Chapter 16) for more information about operad algebras.
\end{problem}

\begin{problem}
In Def. \ref{def:wd-algebra}, check that the $28$ generating axioms in fact correspond to the $28$ elementary relations in the sense of the associativity diagram \eqref{operad-algebra-associativity2}.
\end{problem}

\begin{problem}
In the proof of Theorem \ref{prop:propagator-algebra-is-algebra}--the finite presentation theorem for the propagator algebra--we checked the generating axioms  \eqref{wd-algebra-doubleloop} and \eqref{wd-algebra-loopelement} that are the least obvious.  Give detailed proofs for the other $26$ generating axioms for the propagator algebra.
\end{problem}

\begin{problem}
In Remark \ref{rk:propagator-agree} we pointed out that the structure map of the propagator algebra in \cite{rupel-spivak}, when applied to the generating wiring diagrams (section \ref{sec:generating-wd}), reduces to our $8$ generating structure maps in  Def. \ref{def:propagator-algebra}.  Check carefully that this is indeed the case.
\end{problem}

\begin{problem}
Check carefully the proofs of Theorems \ref{thm:normal-wd-algebra} and \ref{thm:strict-wd-algebra}, the finite presentation theorems for $\wddot$-algebras and $\wdzero$-algebras.
\end{problem}

\begin{problem}
In the proof of Theorem \ref{thm:algebra-ds}--the finite presentation theorem for the algebra of discrete systems--one generating axiom was written down in detail.  Check the other $27$ generating axioms carefully.
\end{problem}

\begin{problem}
In the proof of Theorem \ref{prop:ods-algebra-is-algebra}--the finite presentation theorem for the algebra of open dynamical systems--we checked the generating axiom \eqref{wd-algebra-doubleloop} corresponding to a double-loop.  Give detailed proofs for the other $7$ generating axioms for the algebra of open dynamical systems.
\end{problem}

\begin{problem}
In Remark \ref{rk:ods-algebra-agree}  we pointed out that the structure map of the algebra of open dynamical systems in \cite{vsl}, when applied to the strict generating wiring diagrams--namely, $\epsilon$, $\tau_{X,Y}$, $\theta_{X,Y}$, and $\lambda_{X,x}$--reduces to our $4$ generating structure maps in Def. \ref{def:ods-algebra}.  Check carefully that this is indeed the case.
\end{problem}

\begin{problem}
Check carefully the proof of Propositions \ref{uwd-compi-explicit} and \ref{uwd-gamma}.
\end{problem}

\begin{problem}
Check that each elementary relation in Section \ref{sec:elementary-relations-uwd} is indeed an equality in $\uwd$.
\end{problem}

\begin{problem}
In Example \ref{ex:wasted-cable-by-generators}, check that the iterated operadic composition \eqref{wasted-cable-simplex} is the intended undirected wiring diagram $\varphi \comp \psi$.
\end{problem}

\begin{problem}
In Example \ref{ex:y-plus-wasted-cable}, check that the iterated operadic composition \eqref{wasted-cable-y-simplex} is actually equal to $\zeta_Y$.
\end{problem}

\begin{problem}
In Example \ref{ex:y-plus-wasted-cable-b}, check that the iterated operadic composition \eqref{wasted-cable-y-simplex-b} is actually equal to $\zeta_Y$.
\end{problem}

\begin{problem}
Following the hint in Remark \ref{rk:finite-presentation-algebra}, formulate and prove a finite presentation theorem for a colored operad with given finite sets of generators and relations.
\end{problem}

\begin{problem}
In the proofs of Theorems \ref{relational-algebra-is-algebra} and \ref{typed-relational-algebra}--the finite presentation theorems for the (typed) relational algebra--we checked one of the generating axioms.  Give detailed proofs for the other $16$ generating axioms.
\end{problem}

\begin{problem}
Prove or disprove Spivak's Conjecture \ref{conj:relational-quotient-free} regarding the quotient freeness of the relational algebra.  Then send me an email and tell me how you do it.
\end{problem}

\begin{problem}
In Example \ref{ex:operad-map-algebra}, check that $A^f$ is indeed an $\O$-algebra.
\end{problem}

\begin{problem}
In the proof of Theorem \ref{wddot-uwd-operad-map}, check that it is actually sufficient to prove that the diagram \eqref{chi-operadmap-pushout} is a pushout in $\Fins$.  Then check the equalities \eqref{wddot-uwd-check-equalities}.
\end{problem}

\begin{problem}
In the proof of Theorem \ref{wd-uwd-operad-map}:
\begin{enumerate}
\item Check that it is actually sufficient to prove that the diagram \eqref{rho-operadmap-pushout} is a pushout in $\Fins$.
\item In the diagram \eqref{rho-operadmap-pushout}, check that the right vertical map and the map $h$ are indeed well-defined.
\item Check that the definitions \eqref{rho-operad-map-eta} actually yield a well-defined map $\eta$.
\item Check the equalities \eqref{etah=alpha-etastuff=beta}.
\end{enumerate}
\end{problem}

\chapter{Further Reading}
\label{ch:further-reading}

Listed below are some references for categories, operads, props, and their applications.  Each topic is a huge subject, so this list is not meant to be complete. It represents only a small sample of the existing literature.  The reader is encouraged to consult these books and articles and the references therein.

\section{Category Theory}

These are references for basic category theory, of which \cite{maclane} is the most advanced and \cite{awodey,leinster,riehl} are more basic.  The basic concepts of categories, functors, and natural transformations were all introduced in the founding article \cite{em}.  The paper \cite{jsv} introduced \index{traced monoidal category} traced monoidal categories, which appear in many recent applications of category theory.  The other books all have a view toward applications in the sciences.

\cite{awodey} S. Awodey, Category Theory, 2nd. ed., Oxford Logic Guides 52, Oxford Univ. Press, Oxford, 2010.

\cite{barr-wells} M. Barr and C. Wells, Category Theory for Computing Science, Prentice-Hall , 1990.

\cite{em} S. Eilenberg and S. MacLane, General theory of natural equivalences, Trans. Amer. Math. Soc. 58 (1945), 231-294.

\cite{jsv} A. Joyal, R. Street, and D. Verity, Traced monoidal categories, Mathematical
Proceedings of the Cambridge Philosophical Society 119 (1996), 447-468.

\cite{ls} F.W. Lawvere and S.H. Schanuel, Conceptual Mathematics: A first introduction to categories, 2nd ed., Cambridge, 2009.

\cite{leinster} T. Leinster, Basic Category Theory, Cambridge Studies in Adv. Math. 143, Cambridge Univ. Press, Cambridge, 2014.

\cite{maclane} S. Mac Lane, Categories for the working mathematician, Grad. Texts in Math. 5, 2nd ed., Springer-Verlag, New York, 1998.

\cite{pierce} B.C. Pierce, Basic Category Theory for Computer Scientists, MIT Press,1991.

\cite{riehl} E. Riehl, Category Theory in Context, Dover, New York, 2016.

\cite{spivak14} D.I. Spivak, Category Theory for the Sciences, MIT Press, 2014.

\cite{walters} R.F.C. Walters, Categories and Computer Science, Cambridge, 1991.

\section{Operads}

These are references for operads, originally defined by Lambek \cite{lambek} without symmetric group action and called \index{multicategory} \emph{multicategories}.  The name operad was coined by May in \cite{may72}.  Many applications of operads in mathematics and physics are discussed in \cite{mss}.  The book \cite{lv} is an in-depth study of algebraic operads, and \cite{yau-operad} is a basic introduction to operads in a symmetric monoidal category.

\cite{kelly05} G.M. Kelly, On the operads of J.P. May, Reprints in Theory Appl. Categ. 13 (2005), 1-13. 

\cite{lambek} J. Lambek, Deductive systems and categories. II. Standard constructions and
closed categories, in: 1969 Category Theory, Homology Theory and their
Applications, I (Battelle Inst. Conf., Seattle, Wash., 1968, Vol. 1)
p.76-122, Springer, Berlin, 1969.

\cite{mss} M. Markl, S. Shnider, and J. Stasheff, Operads in Algebra, Topology and Physics, Math. Surveys and Monographs 96, Amer. Math. Soc., Providence, 2002.

\cite{may72} J.P. May, The geometry of iterated loop spaces, Lecture Notes in Math. 271, Springer-Verlag, New York, 1972.

\cite{may97} J.P. May, Definitions: operads, algebras and modules, Contemp. Math. 202, p.1--7, 1997.

\cite{lv} J.-L. Loday and B. Vallette, Algebraic Operads, Grund. der math. Wiss. 346, Springer-Verlag, New York, 2012.

\cite{yau-operad} D. Yau, Colored Operads, Graduate Studies in Math. 170, Amer. Math. Soc., Providence, RI, 2016.

\section{Props}

While an operad models operations with multiple inputs and one output, a \index{prop} prop--short for \emph{product} and \emph{permutation}--models operations with multiple inputs and multiple outputs.  A typical example of a prop is the collection of functions $\Map(X^m,X^n)$ for a set $X$ with $m,n \geq 0$.  Using the kind of pictures in Motivation \ref{mot:operad}, a function $f : X^m \to X^n$ may be depicted as follows.
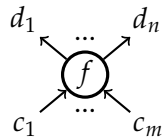
\begin{center}
\begin{tikzpicture}
\matrix[row sep=1cm,column sep=1cm] {
\node [plain,label=above:$...$,label=below:$...$] (p) {$f$}; \\};
\draw [outputleg] (p) to node[above left=.1cm]{$d_1$} +(-.6cm,.5cm);
\draw [outputleg] (p) to node[above right=.1cm]{$d_n$} +(.6cm,.5cm);
\draw [inputleg] (p) to node[below left=.1cm]{$c_1$} +(-.6cm,-.5cm);
\draw [inputleg] (p) to node[below right=.1cm]{$c_m$} +(.6cm,-.5cm);
\end{tikzpicture}
\end{center}
Props were originally defined by Adams and Mac Lane \cite{maclane65}.  Variations of props include \index{wheeled prop} wheeled props, in which there are contraction operations, modeling maps
\[\xi^i_j : \Map(X^m,X^n) \to \Map(X^{m-1},X^{n-1})\]
with $1 \leq i \leq n$ and $1 \leq j \leq m$.  If $f$ is represented as in the previous picture, then its contraction $\xi^i_jf$ may be represented as follows.
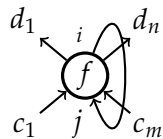
\begin{center}
\begin{tikzpicture}
\matrix[row sep=1cm,column sep=1cm] {
\node [plain] (p) {$f$}; \\};
\draw [outputleg] (p) to node[above left=.1cm]{$d_1$} +(-.6cm,.5cm);
\draw [outputleg] (p) to node[above right=.1cm]{$d_n$} +(.6cm,.5cm);
\draw [inputleg] (p) to node[below left=.1cm]{$c_1$} +(-.6cm,-.5cm);
\draw [inputleg] (p) to node[below right=.1cm]{$c_m$} +(.6cm,-.5cm);
\draw [arrow, out=70, in=-70, looseness=7] (p) to node[at start]{\footnotesize{$i$}} 
node[at end]{$j$} (p);
\end{tikzpicture}
\end{center}
The articles \cite{markl08,vallette} provide surveys of these compositional structures.  The book \cite{jy2} is a comprehensive foundation of this subject.

\cite{maclane65} S. Mac Lane, Categorical algebra, Bull. Amer. Math. Soc. 71 (1965), 40-106.

\cite{markl08} M. Markl, Operads and PROPs, Handbook of Algebra 5, p.87-140, Elsevier, 2008.

\cite{mms} M. Markl, S. Merkulov, and S. Shadrin, Wheeled PROPs, graph complexes and the master equation, J. Pure Appl. Algebra 213 (2009), 496-535.

\cite{vallette} B. Vallette, Algebra $+$ Homotopy $=$ Operad, preprint, arXiv:1202.3245.

\cite{jy2} D. Yau and M.W. Johnson, A Foundation for PROPs, Algebras, and Modules, Math. Surveys and Monographs 203, Amer. Math. Soc., Providence, RI, 2015

\section{Applications of Compositional Structures}

This is a short list of applications of categories, operads, and props in the sciences, including dynamical systems, computer science, engineering, network theory, linguistics, biology, neuroscience, and machine learning.

\cite{asudeh} A. Asudeh, Monads: Some Linguistic Applications, available at\\ \texttt{http://users.ox.ac.uk/$\sim$cpgl0036/handouts/asudeh-se-lfg13.pdf}.

\cite{baez-erbele} J.C. Baez and J. Erbele, Categories in control, arXiv:1405.6881.

\cite{baez-fong} J.C. Baez and B. Fong, A compositional framework for passing linear networks, arXiv:1504.05625.

\cite{baez-fong-pollard} J.C. Baez, B. Fong, and B.S. Pollard, A compositional framework for Markov processes, J. Math. Phys. 57, No. 3, 033301, 30 p. (2016).

\cite{bo} J. Baez and N. Otter, Operads and the tree of life, available at\\  \texttt{http://math.ucr.edu/home/baez/tree$\underline{\hspace{.2cm}}$of$\underline{\hspace{.2cm}}$life/tree$\underline{\hspace{.2cm}}$of$\underline{\hspace{.2cm}}$life.pdf}.

\cite{baez-stay} J.C. Baez and M. Stay, Physics, Topology, Logic and Computation: A Rosetta Stone, in: New Structures for Physics, ed. Bob Coecke, Lecture Notes in Physics vol. 813, Springer, Berlin, 2011, pp. 95-174.

\cite{bsz} F. Bonchi, P. Soboci\'{n}ski, and F. Zanasi, A categorical semantics of signal flow graphs, pp. 435-450, Lecture Notes in Comp. Sci. 8704, Springer, 2014.

\cite{bp} R. Brown and T. Porter, Category Theory and Higher Dimensional Algebra: potential descriptive tools in neuroscience, in: Proceedings of the International Conference on Theoretical Neurobiology, Delhi, February 2003, edited by Nandini Singh, National Brain Research Centre, Conference Proceedings 1 (2003) 80-92.

\cite{coecke} B. Coecke, Quantum picturalism, Contemporary Physics 51 (2010), 59-83.

\cite{hsg} C. Heunen, M. Sadrzadeh, and E. Grefenstette, ed., Quantum Physics and Linguistics: A Compositional, Diagrammatic Discourse, Oxford, 2013.

\cite{ior} O. Iordache, Self-Evolvable Systems: Machine Learning in Social Media, Springer-Verlag, Berlin, 2012.

\cite{mendez} M.A. M\'{e}ndez, Set Operads in Combinatorics and Computer Science, Springer, New York, 2015.

\cite{rupel-spivak} D. Rupel and D.I. Spivak, The operad of temporal wiring diagrams: formalizing a graphical language for discrete-time processes, arXiv:1307.6894.

\cite{spivak13} D.I. Spivak, The operad of wiring diagrams: formalizing a graphical language for databases, recursion, and plug-and-play circuits, arXiv:1305.0297.

\cite{spivak15} D.I. Spivak, Nesting of dynamic systems and mode-dependent networks, arXiv:1502.07380.

\cite{spivak15b} D.I. Spivak, The steady states of coupled dynamical systems
compose according to matrix arithmetic, arXiv:1512.00802.

\cite{ssr15} D.I. Spivak, P. Schultz, and D. Rupel, String diagrams for traced and compact categories are oriented $1$-cobordisms, arXiv:1508.01069.

\cite{vsl} D. Vagner, D.I. Spivak, and E. Lerman, Algebras of open dynamical systems on the operad of wiring diagrams, Theory Appl. Categ. 30 (2015), 1793-1822.

\appendix

\backmatter

\bibliographystyle{amsalpha}

\printindex

\end{document}